\newcommand{\xBc}{\langle}
\newcommand{\xBe}{\rangle}

\newcommand{\xbD}{\Delta}
\newcommand{\xbF}{\Phi}
\newcommand{\xbG}{\Gamma}

\newcommand{\xbL}{\Lambda}
\newcommand{\xbO}{\Omega}
\newcommand{\xbP}{\Pi}

\newcommand{\xbS}{\Sigma}

\newcommand{\xba}{\alpha}
\newcommand{\xbb}{\beta}

\newcommand{\xbd}{\delta}
\newcommand{\xbe}{\in}
\newcommand{\xbf}{\phi}
\newcommand{\xbg}{\gamma}
\newcommand{\xbh}{\eta}

\newcommand{\xbj}{\vartheta}
\newcommand{\xbk}{\kappa}
\newcommand{\xbl}{\lambda}
\newcommand{\xbm}{\mu}
\newcommand{\xbn}{\nu}
\newcommand{\xbo}{\omega}
\newcommand{\xbp}{\pi}
\newcommand{\xbq}{\psi}
\newcommand{\xbr}{\rho}
\newcommand{\xbs}{\sigma}
\newcommand{\xbt}{\tau}

\newcommand{\xCB}{A}

\newcommand{\xCI}{{\Big(}}
\newcommand{\xCJ}{{\Big)}}
\newcommand{\xCK}{\times}
\newcommand{\xCL}{\pm}

\newcommand{\xCN}{\neg}
\newcommand{\xCQ}{\emptyset}

\newcommand{\xCc}{<}
\newcommand{\xCd}{\approx}
\newcommand{\xCe}{>}
\newcommand{\xCf}{\hspace{0.1em}}

\newcommand{\xCq}{\sim}

\newcommand{\xcA}{\forall}
\newcommand{\xcB}{\stackrel{\subset}{\neq}}
\newcommand{\xcC}{\not\subseteq}

\newcommand{\xcE}{\exists}

\newcommand{\xcH}{\not\Rightarrow}
\newcommand{\xcI}{\not\Leftarrow}
\newcommand{\xcJ}{\not\Leftrightarrow}

\newcommand{\xcL}{\not\vdash}
\newcommand{\xcM}{\not\models}
\newcommand{\xcN}{\hspace{0.2em}\not\sim\hspace{-0.9em}\mid\hspace{0.8em}}
\newcommand{\xcO}{\bigvee}
\newcommand{\xcP}{\not\rightarrow}

\newcommand{\xcS}{\bigcap}
\newcommand{\xcT}{\bot}
\newcommand{\xcU}{\bigwedge}
\newcommand{\xcV}{\bigcup}

\newcommand{\xcX}{\Box}

\newcommand{\xca}{\infty}
\newcommand{\xcb}{\subset}
\newcommand{\xcc}{\subseteq}
\newcommand{\xcd}{\supseteq}
\newcommand{\xce}{\not\in}
\newcommand{\xcf}{\supset}
\newcommand{\xcg}{\geq}
\newcommand{\xch}{\Rightarrow}
\newcommand{\xci}{\Leftarrow}
\newcommand{\xcj}{\Leftrightarrow}
\newcommand{\xck}{\leq}
\newcommand{\xcl}{\vdash}
\newcommand{\xcm}{\models}
\newcommand{\xcn}{\hspace{0.2em}\sim\hspace{-0.9em}\mid\hspace{0.58em}}

\newcommand{\xco}{\vee}
\newcommand{\xcp}{\rightarrow}
\newcommand{\xcq}{\leftarrow}
\newcommand{\xcr}{\leftrightarrow}
\newcommand{\xcs}{\cap}
\newcommand{\xcu}{\wedge}
\newcommand{\xcv}{\cup}
\newcommand{\xcx}{\Diamond}

\newcommand{\xcz}{\Box}

\newcommand{\xDC}{\hspace{2em}}
\newcommand{\xDH}{\item }

\newcommand{\xDM}{\circ}
\newcommand{\xDN}{\ominus}

\newcommand{\xDc}{\ll}
\newcommand{\xDd}{\equiv}

\newcommand{\xdA}{\mbox{\boldmath$A$}}
\newcommand{\xdB}{\mbox{\boldmath$B$}}
\newcommand{\xdC}{\mbox{\boldmath$C$}}
\newcommand{\xdD}{\mbox{\boldmath$D$}}

\newcommand{\xdL}{\mbox{\boldmath$L$}}

\newcommand{\xdO}{\mbox{\boldmath$O$}}
\newcommand{\xdP}{\mbox{\boldmath$P$}}

\newcommand{\xdR}{\Re}

\newcommand{\xda}{{\cal A}}

\newcommand{\xdc}{{\cal C}}
\newcommand{\xdd}{{\cal D}}
\newcommand{\xde}{{\cal E}}
\newcommand{\xdf}{{\cal F}}

\newcommand{\xdi}{{\cal I}}

\newcommand{\xdl}{{\cal L}}
\newcommand{\xdm}{{\cal M}}
\newcommand{\xdn}{{\cal N}}
\newcommand{\xdo}{{\cal O}}
\newcommand{\xdp}{{\cal P}}

\newcommand{\xds}{{\cal S}}

\newcommand{\xdu}{{\cal U}}

\newcommand{\xdw}{{\cal W}}
\newcommand{\xdx}{{\cal X}}
\newcommand{\xdy}{{\cal Y}}
\newcommand{\xdz}{{\cal Z}}

\newcommand{\xEC}{\spadesuit}
\newcommand{\xEH}{ & }
\newcommand{\xEI}{\begin{itemize}}
\newcommand{\xEJ}{\end{itemize}}
\newcommand{\xEP}{ \\ }

\newcommand{\xEc}{\not<}
\newcommand{\xEd}{\neq}
\newcommand{\xEh}{\begin{enumerate}}
\newcommand{\xEj}{\end{enumerate}}

\newcommand{\xeA}{\nabla}
\newcommand{\xeB}{\not\prec}
\newcommand{\xeC}{\not\preceq}

\newcommand{\xeY}{\triangle}

\newcommand{\xeb}{\prec}
\newcommand{\xec}{\preceq}
\newcommand{\xed}{\succeq}
\newcommand{\xee}{\succ}

\newcommand{\xej}{\lhd}

\newcommand{\xem}{\rhd}
\newcommand{\xen}{\blacktriangleleft}

\newcommand{\xeq}{\blacktriangleright}
\newcommand{\xer}{\sqsubset}
\newcommand{\xes}{\sqsubseteq}

\newcommand{\xex}{\lceil}

\newcommand{\xFB}{\cdots}
\newcommand{\xFO}{\parallel}

\newcommand{\xfA}{\mid}
\newcommand{\xfB}{\uparrow}

\newcommand{\xfo}{\hookrightarrow}

\newcommand{\Xl}{\ldots}

\newcommand{\ol}{\overline}
\newcommand{\ul}{\underline}

\newcommand{\wt}{\overbrace}

\newcommand{\xssc}{\scriptsize}

\newcommand{\bl}{\begin{lemma} \rm}
\newcommand{\el}{\end{lemma}}
\newcommand{\br}{\begin{remark} \rm}
\newcommand{\er}{\end{remark}}
\newcommand{\be}{\begin{example} \rm}
\newcommand{\ee}{\end{example}}
\newcommand{\bco}{\begin{corollary} \rm}
\newcommand{\eco}{\end{corollary}}
\newcommand{\bc}{\begin{claim} \rm}
\newcommand{\ec}{\end{claim}}
\newcommand{\bfa}{\begin{fact} \rm}
\newcommand{\efa}{\end{fact}}
\newcommand{\bp}{\begin{proposition} \rm}
\newcommand{\ep}{\end{proposition}}
\newcommand{\bd}{\begin{definition} \rm}
\newcommand{\ed}{\end{definition}}
\newcommand{\bcs}{\begin{construction} \rm}
\newcommand{\ecs}{\end{construction}}
\newcommand{\bcd}{\begin{condition} \rm}
\newcommand{\ecd}{\end{condition}}
\newcommand{\bt}{\begin{theorem} \rm}
\newcommand{\et}{\end{theorem}}
\newcommand{\bn}{\begin{notation} \rm}
\newcommand{\en}{\end{notation}}
\newcommand{\bfi}{\begin{bild} \rm}
\newcommand{\efi}{\end{bild}}
\newcommand{\bsta}{\begin{statement} \rm}
\newcommand{\esta}{\end{statement}}
\newcommand{\bcom}{\begin{comment} \rm}
\newcommand{\ecom}{\end{comment}}
\newcommand{\bdia}{\begin{diagram} \rm}
\newcommand{\edia}{\end{diagram}}

\newcommand{\bfc}{\begin{figure}[htb] \begin{center}}
\newcommand{\efc}{\end{center} \end{figure}}

\documentclass{book}

\usepackage{amssymb,latexsym,epic,eepic,rotating}

\title{
Logical tools for handling change in agent-based systems
}

\author{Dov M Gabbay
\thanks{
Dov.Gabbay@kcl.ac.uk, www.dcs.kcl.ac.uk/staff/dg
} \\
King's College, London
\thanks{
Department of Computer Science, King's College London, Strand,
London WC2R 2LS, UK
} \\ \\
Karl Schlechta
\thanks{
ks@cmi.univ-mrs.fr, karl.schlechta@web.de, http://www.cmi.univ-mrs.fr/ $\sim$ ks
} \\
Laboratoire d'Informatique Fondamentale de Marseille
\thanks{
UMR 6166, CNRS and Universit\'{e} de Provence,
Address: CMI, 39, rue Joliot-Curie, F-13453 Marseille Cedex 13, France
}
}

\begin{document}

\maketitle

\newtheorem{lemma}{Lemma}[section]
\newtheorem{theorem}[lemma]{Theorem}
\newtheorem{proposition}[lemma]{Proposition}
\newtheorem{corollary}[lemma]{Corollary}
\newtheorem{claim}[lemma]{Claim}
\newtheorem{fact}[lemma]{Fact}
\newtheorem{remark}[lemma]{Remark}
\newtheorem{definition}{Definition}[section]
\newtheorem{construction}{Construction}[section]
\newtheorem{condition}{Condition}[section]
\newtheorem{example}{Example}[section]
\newtheorem{notation}{Notation}[section]
\newtheorem{bild}{Figure}[section]
\newtheorem{comment}{Comment}[section]
\newtheorem{statement}{Statement}[section]
\newtheorem{diagram}{Diagram}[section]

\renewcommand{\labelenumi}
  {(\arabic{enumi})}
\renewcommand{\labelenumii}
  {(\arabic{enumi}.\arabic{enumii})}
\renewcommand{\labelenumiii}
  {(\arabic{enumi}.\arabic{enumii}.\arabic{enumiii})}
\renewcommand{\labelenumiv}
  {(\arabic{enumi}.\arabic{enumii}.\arabic{enumiii}.\arabic{enumiv})}

\setcounter{secnumdepth}{4}
\setcounter{tocdepth}{4}

\tableofcontents

\pagenumbering{arabic}

\setcounter{page}{7}

\chapter{Introduction and Motivation}

% \include{1-1-int}

% ******* BEGIN LATEX SOURCE FILE 1-3-icr.tex *******
%
% Uebers. aus Karltex File: 1-3-icr.m
%
%

Throughout, unless said otherwise, we will work in propositional logic.
\section{
Program
}
\index{Motivation IBRS}

The human agent in his daily activity has to deal with many situations
involving change. Chief among them are the following

 \xEh

 \xDH Common sense reasoning from available data. This involves
predication
of what unavailable data is supposed to be (nonmonotonic deduction)
but it is a defeasible prediction, geared towards immediate change.
This is formally known as nonmonotonic reasoning and is studied by
the nonmonotonic community.

 \xDH Belief revision, studied by a very large community. The agent
is unhappy with the totality of his beliefs which he finds internally
unacceptable (usually logically inconsistent but not necessarily so)
and needs to change/revise it.

 \xDH Receiving and updating his data, studied by the update community.

 \xDH Making morally correct decisions, studied by the deontic logic
community.

 \xDH Dealing with hypothetical and counterfactual situations. This
is studied by a large community of philosophers and AI researchers.

 \xDH Considering temporal future possibilities, this is covered by
modal and temporal logic.

 \xDH Dealing with properties that persist through time in the near
future and with reasoning that is constructive. This is covered by
intuitionistic logic.

 \xEj

All the above types of reasoning exist in the human mind and are used
continuously and coherently every hour of the day. The formal modelling
of these types is done by diverse communities which are largely distinct
with no significant communication or cooperation. The formal models
they use are very similar and arise from a more general theory, what
we might call:

``Reasoning with information bearing binary relations''.

\section{
Short overview of the different logics
}
\label{Section Logic-Overview}

We will discuss the semantics of the propositional logic situation only.

In all cases except the last two (i.e. Inheritance and Argumentation), the
semantics consist of a set of classical models for the
underlying language, with an additional structure, usually
a binary relation (sometimes relative to a point of origin).
This additional structure is not unique, and the result of the reasoning
based
on this additional structure will largely depend on the specific choice of
this
structure. The laws which are usually provided (as axioms or rationality
postulates) are those which hold for any such additional structure.
\subsection{
Nonmonotonic logics
}

Nonmonotonic logics (NML) were created to deal with principled reasoning
about
``normal'' situation. Thus, ``normal'' birds will (be able to) fly, but there
are
many exceptions, like penguins, roasted chickens, etc., and it is usually
difficult to enumerate all exceptions, so they will be treated in bulk as
``abnormal'' birds.

The standard example is - as we began to describe already - that ``normal''
birds
will (be able to) fly, that there are exceptions, like penguins, that
``normal''
penguins will not fly, but that there might be exceptions to the
exceptions,
that some abnormal penguin might be able to fly - due to a jet pack on its
back,
some genetic tampering, etc. Then, if we know that some animal is a bird,
call
it ``Tweety'' as usual, and if we want to keep it as a pet, we should make
sure
that its cage has a roof, as it might fly away otherwise. If, however, we
know
that Tweety is not only a bird, but also a penguin, then its cage need not
have
a roof.

Note that this reasoning is nonmonotonic: From the fact ``Tweety is a
bird'', we
conclude that it will (normally) fly, but from the facts that ``Tweety is a
bird''
and ``Tweety is a penguin'', we will not conclude that it will (normally)
fly any
more, we will even conclude the contrary, that it will (normally) not fly.

We can also see here a general principle at work: more specific
information
(Tweety is a penguin) and its consequences (Tweety will not fly) will
usually be
considered more reliable than the more general information (Tweety is a
bird)
and its consequences (Tweety will fly). Then, NML can also be considered
as a
principled treatment of information of different quality or reliability.
The
classical information is the best one, and the conjecture that the case at
hand
is a normal one is less reliable.

Note that normality is absolute here in the following sense: normal birds
will
be normal with respect to all ``normal'' properties of birds, i.e. they will
fly,
lay eggs, build nests, etc. In this treatment, there are no birds normal
with
respect to flying, but not laying eggs, etc.

It is sometimes useful to introduce a generalized quantifier $ \xeA.$
In a first order (FOL) setting $ \xeA x \xbf (x)$ will mean that $ \xbf
(x)$ holds almost
everywhere, in a propositional setting $ \xeA \xbf $ will mean that in
almost all
models $ \xbf $ holds. Of course, this ``almost everywhere'' or
``almost all'' has to be made precise, e.g. by a filter over the
FOL universe, or the set of all propositional models.

Inheritance systems will be discussed separately below.

 \xEI
 \xDH Formal semantics by preferential systems

The semantics for preferential logics are preferential structures, a set
of
classical models with an arbitrary binary relation. This relation need not
be
transitive, nor does it need to have any other of the usual properties. If
$m \xeb m',$
then $m$ is considered more normal (or less abnormal) than $m'.$ $m$ is
said to be
minimal in a set of models $M$ iff there is no $m' \xbe M,$ $m' \xeb m$ -
a word of warning:
there might be $m' \xeb m,$ but $m' \xce M!$

This defines a semantic consequence relation as follows: we say $ \xbf
\xcn \xbq $ iff
$ \xbq $ holds in all minimal models of $ \xbf.$

As a model $m$ might be minimal in $M( \xbf )$ - the set of models of $
\xbf $ - but not
minimal in $M( \xbq ),$ where $ \xcm \xbf \xcp \xbq $ classically, this
consequence relation $ \xcn $ is
nonmonotonic. Non-flying penguins are normal $(=minimally$ abnormal)
penguins, but
all non-flying birds are abnormal birds.

Minimal models of $ \xbf $ need not
exist, even if $ \xbf $ is consistent - there might be cycles or infinite
descending chains. We will write $M( \xbf )$ for the set of $ \xbf
-$models,
and $ \xbm ( \xbf )$ or $ \xbm (M( \xbf ))$ for the set of minimal models
of $ \xbf.$
If there is some set $X$ and some $x' \xbe X$ s.t. $x' \xeb x,$ we say
that
$x' $ minimzes $x,$ likewise that $X$ minimizes $x.$ We will be more
precise in Chapter \ref{Chapter Pref} (page \pageref{Chapter Pref}).

One can impose various restrictions on $ \xeb,$ they will
sometimes change the resulting logic.
The most important one is perhaps rankedness: If $m$ and $m' $ are
$ \xeb -$incomparable, then for all $m'' $ $m'' \xeb m$ iff $m'' \xeb m' $
and also $m \xeb m'' $ iff $m' \xeb m''.$ We can
interpret the fact that $m$ and $m' $ are $ \xeb -$incomparable by putting
them at the same
distance from some imaginary point of maximal normality. Thus, if $m$ is
closer
to this point than $m'' $ is, then so will be $m',$ and if $m$ is farther
away from
this point than $m'' $ is, then so will be $m'.$ (The also very important
condition,
smoothness, is more complicated, so the reader is referred to
Chapter \ref{Chapter Pref} (page \pageref{Chapter Pref})  for discussion.

Preferential structures are presented and discussed
in Chapter \ref{Chapter Pref} (page \pageref{Chapter Pref}).

 \xEJ
\subsection{
Theory revision
}

The problem of Theory Revision is to ``integrate'' some new information $
\xbf $ into an
old body of knowledge $K$ such that the result is consistent, even if $K$
together
with $ \xbf $ (i.e. the union $K \xcv \{ \xbf \})$ is inconsistent. (We
will assume that $K$ and $ \xbf $
are consistent separately.)

The best examined approach was first published in  \cite{AGM85},
and is know for the intials of its authors as the AGM approach.
The formal presentation of this approach (and more) is in
Chapter \ref{Chapter TR} (page \pageref{Chapter TR}).

This problem is well known in juridical thinking, where a new law might be
inconsistent with the old set of laws, and the task is to ``throw away''
enough,
but not too many, of the old laws, so we can incorporate the new law into
the
old system in a consistent way.

We can take up the example for NML, and modify it slightly. Suppose our
background theory $K$ is that birds fly, in the form: Blackbirds fly,
ravens fly,
penguins fly, robins fly,  \Xl., and that the new information is that
penguins
$don' t$ fly. Then, of course, the minimal change to the old theory is to
delete
the information that penguins fly, and replace it with the new
information.

Often, however, the situation is not so simple. $K$ might be that $ \xbq $
holds, and so
does $ \xbq \xcp \xbr.$ The new information might be that $ \xCN \xbr $
holds. The radical - and
usually excessive - modification will be to delete all information from
$K,$ and
just take the new information. More careful modifications will be to
delete
either $ \xbq $ or $ \xbq \xcp \xbr,$ but not both. But there is a
decision problem here: which
of the two do we throw out? Logic alone cannot tell us, and we will need
more
information to take this decision.

 \xEI
 \xDH Formal semantics

In many cases, revising $K$ by $ \xbf $ is required to contain $ \xbf,$
thus, if $*$ denotes
the revision operation, then $K* \xbf \xcl \xbf $ (classically). Dropping
this requirement
does not change the underlying mechanism enormously, we will uphold it.

Speaking semantically, $K* \xbf $ will then be defined by some subset of
$M( \xbf ).$ If we
choose all of $ \xbf,$ then any influence of $K$ is forgotten. A good way
to capture
this influence seems to choose those models of $ \xbf,$ which are closest
to the
$K-$models, in some way, and with respect to some distance $d.$ We thus
choose those
$ \xbf -$models $m$ such that there is $n \xbe M(K)$ with $d(n,m)$ minimal
among all $d(n',m' ),$
$n' \xbe M(K),$ $m' \xbe M( \xbf ).$ (We assume again that the minimal
distance exists, i.e.
that there are no infinite descending chains of distances, without any
minimum.)
Of course, the choice of the distance is left open, and will influence the
outcome. For instance, choosing as $d$ the trivial distance, i.e.
$d(x,y)=1$ iff
$x \xEd y,$ and 0 otherwise, will give us just $ \xbf $ - if $K$ is
inconsistent with $ \xbf.$

This semantic approach corresponds well to the classical, syntactic AGM
revision
approach in the following sense: When we fix $K,$ this
semantics corresponds exactly to the AGM postulates (which leave $K$
fixed). When
we allow $K$ to change, we can also treat iterated revision, i.e.
something like
$(K* \xbf )* \xbq,$ thus go beyond the AGM approach (but pay the price of
arbitrarily
long axioms). This semantics leaves the order (or distance) untouched, and
is
thus fundamentally different from e.g. Spohn's Ordinal Conditional
Functions.

 \xEJ
\subsection{
Theory update
}

Theory Update is the problem of ``guessing'' the results of in some way
optimal
developments.

Consider the following situation: There is a house, at time 0, the light
is on,
and so is the deep freezer. At time 1, the light is off. Problem: Is the
deep
freezer still on? The probably correct answer depends on circumstances.
Suppose
in situation $ \xCf A,$ there is someone in the house, and weather
conditions are
normal. In situation $B,$ there is no one in the house, and there is a
very heavy
thunderstorm going on. Then, in situation $ \xCf A,$ we will conjecture
that the
$person(s)$ in the house have switched the light off, but left the deep
freezer
on. In situation $B,$ we might conjecture a general power failure, and
that the
deep freezer is now off, too.

We can describe the states at time 0 and 1 by a triple: light on/off,
freezer
on/off, power failure yes/no.

In situation $ \xCf A,$ we will consider the development
(light on, freezer on, no power failure) to
(light off, freezer on, no power failure)
as the most likely (or normal) one.

In situation $B$ we will consider the development
(light on, freezer on, no power failure) to
(light off, freezer off, power failure)
as the most likely (or normal) one.

Often, we will assume a general principle of inertia: things stay the way
as
they are, unless they are forced to change. Thus, when the power failure
is
repaired, freezer and light will go on again.

 \xEI
 \xDH Formal semantics

In the general case, we will consider a set of fixed length sequences of
classical models, say $m= \xBc m_{0},m_{1}, \Xl,m_{n} \xBe,$ which represent
develoments considered
possible. Among this set, we have some relation $ \xeb,$ which is
supposed to single
out the most natural, or probable ones. We then look at some coordinate,
say $i,$
and try to find the most probable situation at coordinate $i.$ For
example, we
have a set $S$ of sequences, and look at the theory defined by the
information at
coordinate $i$ of the most probable sequences of $S:$ $Th(\{m_{i}:m \xbe
\xbm (S)\})$ - where
$ \xbm (S)$ is the set of the most probable sequences of $S,$ and $Th(X)$
is the set
of formulas which hold in all models $x \xbe X.$

Looking back at our above intuitive example, $S$ will be the set of
sequences
consisting of

$ \xBc (l,f,-p),(l,f,p) \xBe,$

$ \xBc (l,f,-p),(l,f,-p) \xBe,$

$ \xBc (l,f,-p),(l,-f,p) \xBe,$

$ \xBc (l,f,-p),(l,-f,-p) \xBe,$

$ \xBc (l,f,-p),(-l,f,p) \xBe,$

$ \xBc (l,f,-p),(-l,f,-p) \xBe,$

$ \xBc (l,f,-p),(-l,-f,p) \xBe,$

$ \xBc (l,f,-p),(-l,-f,-p) \xBe,$

where ``l'' stands for ``light on'', ``f'' for ``freezer on'', ``p'' for ``power
failure''
etc. The ``best'' sequence in situation A will be $ \xBc (l,f,-p),(-l,f,-p) \xBe
,$
and in
situation $B$ $ \xBc (l,f,-p),(-l,-f,p) \xBe.$ Thus, in situation A, the result
is
defined by
-l, $f,$ -p, etc. - the theory of the second coordinate.

Thus, again, the choice of the actual distance has an enormous influence
on the
outsome.

 \xEJ
\subsection{
Deontic logic
}

Deontic logic treats (among other things) the moral acceptability of
situations
or acts.

For instance, when driving a car, you should not cause accidents and hurt
someone. So, in all ``good'' driving situations, there are no accidents and
no
victims. Yet, accidents unfortunately happen. And if you have caused an
accident, you should stop and help the possibly injured. Thus, in the
``morally
bad'' situations where you have caused an accident, the morally best
situations
are those where you help the victims, if there are any.

The parallel to above example for NML is obvious, and, as a matter of
fact, the
first preferential semantics was given for deontic, and not for
nonmonotonic
logics - see Section \ref{Section Deon-Semantik} (page \pageref{Section
Deon-Semantik}).

There is, however, an important difference to be made. Preferential
structures for NML describe what $ \xCf holds$ in the normally best
models, those for
deontic logic what holds in ``morally'' best models. But obligations
are not supposed to say what holds in the morally best worlds, but
are supposed to distinguish in some way the ``good'' from
the ``bad'' models. This problem is discussed in extenso
in Section \ref{Section Deon-Semantik} (page \pageref{Section Deon-Semantik}).

 \xEI
 \xDH Formal semantics

As said already, preferential structures as defined above for NML were
given as
a semantics for deontic logics, before NML came into existence.

A word of warning: Here, the morally optimal models describe ``good''
situations,
and not directly actions to take. This is already obvious by the law of
weakening, which holds for all such structures: If $ \xbf $ holds in all
minimal
models, and $ \xcl \xbf \xcp \xbq $ (classically), then so does $ \xbq.$
But if one should be kind,
then it does not follow that one should be kind or kill one's grandmother.
Of
course, we can turn this reasoning into advice for action: act the way
that the
outcome of your actions assures you to be in a morally good situation.

 \xEJ
\subsection{
Counterfactual conditionals
}

A counterfactual conditional states an implication, where the antecedent
is (at
the moment, at least) wrong. ``If it were to rain, he would open his
umbrella.''
This is comprehensible, the person has an umbrella, and if it were to
start to
rain now, he would open the umbrella. If, however, the rain would fall in
the midst of a hurricane, then opening the umbrella would only lead to its
destruction. Thus, if, at the moment the sentence was uttered, there was
no
hurricane, and no hurricane announced, then the speaker was referring to a
situation which was different from the present situation only in so far as
it is
raining, or, in other words, minimally different from the actual
situation, but
with rain falling. If, however, there was a hurricane in sight at the
moment of
uttering the sentence, we might doubt the speakers good sense, and point
the
problem out to him/her. We see here again a reasoning about minimal
change, or
normal situations.

 \xEI
 \xDH Formal semantics

Stalnaker and Lewis first gave a minimal distance semantics in the
following
way:

If we are in the actual situation $m,$ then $ \xbf > \xbq $ (read: if $
\xbf $ were the case, then
$ \xbq $ would also hold) holds in $m,$ iff in all $ \xbf -$models which
are closest to $m,$ $ \xbq $
also holds. Thus, there might well be $ \xbf -$models where $ \xbq $
fails, but these are
not among the $ \xbf -$models closest to $m.$ The distance will, of
course, express the
difference between the situation $m$ and the models considered. Thus, in
the first
scenario, situations where it rains and there is no extreme wind condition
are
closer to the original one than those where a hurricane blows.

In the original approach, distances from each possible actual situation
are
completely independent. It can, however, be shown that we can achieve the
same
results with one uniform distance over the whole structure, see
 \cite{SM94}.

 \xEJ
\subsection{
Modal logic
}

Modal logic reasons about the possible or necessary. If we are in the
midwest of
the US, and it is a hurricane season, then the beautiful sunny weather
might
turn into a hurricane over the next hours. Thus, the weather need not
necessarily stay the way it is, but it might become a very difficult
situation.
Note that we reason here not about what is likely, or normal, but about
what is
considered possible at all. We are not concerned only about what might
happen in
time $t+1,$ but about what might happen in some (foreseeable, reasonable)
future -
and not about what will be the case at the end of the developments
considered
possible either. Just everything which might be the case some time in the
near
future.

``Necessary'' and ``possible'' are dual: if $ \xbf $ is necessarily the case,
this means
that it will always hold in all situations evolving from the actual
situation,
and if $ \xbf $ is possibly the case, this means that it is not necessary
that $ \xCN \xbf $
holds, i.e. there is at least some situation into which the present can
evolve
and where $ \xbf $ holds.

 \xEI
 \xDH Formal semantics

Kripke gave a semantics for Modal Logic by possible worlds, i.e. a set of
classical models, with an additional binary relation, expressing
accessibility.

If $m$ is in relation $R$ with $n,$ $ \xCf mRn,$ then $m$ can possibly
become $n,$ is a
possibility seen from $m,$ or whatever one might want to say. Again, $R$
can be any
binary relation. The necessity operator is essentially a universal
quantifier,
$ \xcX \xbf $ holds in $m$ iff $ \xbf $ holds in all $n$ accessible via
$R$ from $m.$ Likewise, the
possibility operator is an existential quantifier, $ \xcx \xbf $ holds in
$m$ iff there is
at least one $n$ accessible from $m$ where $ \xbf $ holds.

Again, it is interesting to impose additional postulates on the relation
$R,$ like
reflexivity, transitivity, etc.

 \xEJ
\subsection{
Intuitionistic logic
}

Intuitionistic Logic is (leaving philosophy aside) reasoning about
performed
constructions and proofs in mathematics, or development of (certain)
knowledge.
We may have a conjecture - or, simply,
any statement -, and a proof for it, a proof for its contrary, or neither.
Proofs are supposed to be correct, so what is considered a proof will stay
one
forever. Knowledge can only be won, but not lost. If we have neither a
proof
nor a refutation, then we might one day have one or the other, or we might
stay
ignorant forever.

 \xEI
 \xDH Formal semantics

Intuitionistic Logic can also be given a semantics in the style of the one
for
Modal Logics. There are two, equivalent, variants. The one closer to Modal
Logic
interprets intuitionistic statements (in the above sense that a
construction or
proof has been performed) as preceeded by the necessity quantifier. Thus,
it is
possible that in $m$ neither $ \xcX \xbf $ nor $ \xcX \xCN \xbf $ hold, as
we have neither a proof for
$ \xbf,$ nor for its contrary, and might well find one for one or the
other in some
future, possible, situation. Progressing along $R$ may, if the relation is
transitive, only make more statements of the form $ \xcX \xbf $ true, as
we quantify then
over less situations.

 \xEJ
\subsection{
Inheritance systems
}

Inheritance systems or diagrams are directed acyclic graphs with
two types of arrows, positive and negative ones. Roughly, nodes stand for
sets of objects, like birds, penguins, etc., or properties like
``able to fly''.
A positive arrow $a \xcp b$ stands for
`` (almost) all $x \xbe a$ are also in $b$ '' - so it admits exceptions.
A negative arrow $a \xcp b$ stands for
`` (almost) all $x \xbe a$ are $ \xCf not$ in $b$ '' - so it also admits
exceptions.
Negation is thus very strong.
The problem is to find the valid paths (concatenations of arrows)
in a given diagram, considering contradictions and
specificity. See Chapter \ref{Chapter Inheritance} (page \pageref{Chapter
Inheritance})
for a deeper explanaton and formal definitions.

 \xEI
 \xDH Formal semantics

 \xEJ
\subsection{
Argumentation theory
}

???

 \xEI
 \xDH Formal semantics
 \xEJ
???
\subsection{
A summarizing table for the semantics
}

The following table summarizes important properties of the semantics of
some
of the logics discussed.

``cpm'' stands for ``classical propositional model'', ``CPM'' for a set of such,

``T.'' means: number of types of nodes/arrows,

``Counterfactuals'': counterfactual conditionals,

$'' \xbL '' $ means: Limit version in the absence of minimal/closest or
best points,

``collect'' means: all points in relation are collected,

``Relativization'' means: some reasonable ways to make the logic
dependent on a fixed model as starting point. (This is evident for
modal logic, less evident for some other logics.)

\begin{turn}{90}

{\tiny

\begin{tabular}{|l|c|l|c|l|l|l|l|l|}

\hline

\multicolumn{9}{|c|}{Meaning, Properties, Use} \\

\hline

Logic & \multicolumn{2}{c|}{Nodes} & \multicolumn{4}{c|}{Arrows (Relations)} &
Remarks & Relativization \\
\cline{2-7}

 & T. & Meaning & T. & Meaning & Properties & Use & & \\

\hline

Modal & 1 & cpm & 1 & poss. development? & transitive & collect &
& What is reachable from one point \\
 & & & & & reflexive,... & (reachable) &
 & Summary over all points (monotonic)  \\
 & & & & & & nodes & & \\

\hline

Intuitionistic & 1 & cpm & 1 & epist. change: & transitive & collect &
& What is reachable from one point  \\
 & & & & progress of knowl. & & & & Summary over all points (monotonic)  \\

\hline

Preferential & 1 & cpm & 1 & epist. change: & ranked & best nodes &
& What is optimal below one point  \\
 & & & & conjecture & smooth,... & or $\xbL$ &
 & Summary over all points (monotonic)  \\

\hline

Deontic & 1 & cpm & 1 & moral (value) & ranked... & best/$\xbL$ &
& What is optimal below one point  \\
 & & & & change & & & & Summary over all points (monotonic)  \\

\hline

Counterfactuals & 1 & cpm & 1 & change in world to & ranked & individ. best/ &
& What is closest to one point  \\
 & & & & make prereq. real & (distance) & $\xbL$ &
 & Summary over all points (monotonic) \\

\hline

full TR & 1 & cpm & 1 & epist. change: & ranked & global best/ & global best
& $A\xfA B$ is all of B not   \\
 & & & & incorp. new info & (distance) & $\xbL$ & needs
 & more distant from A than m is from A \\
 & & & & & & & distance for & If m is closer to A than B is, then $\xCQ$  \\
 & & & & & & & comparison & \\

\hline

fixed & 2 & CPM & 1 & epist. change: & ranked? & global best/ &
& Same as full TR \\
K TR & & cpm & & incorp. new info & (AGM & $\xbL$ & & \\
 & & & & & tradition) & & & \\

\hline

Update & 1 & CPM & 1 & change of world, & ranked? & individ. best/ &
& What is optimal in all threads through \\
 & & & & evolving sit. & distance? & $\xbL$ &
 & a point satisfying start condition \\
 & & & & & & & & Summary over all points (monotonic) \\

\hline

Inheritance & 1 & CPM & 2 & normally, $\xbf$'s & acyclic? & collect
& comparison & - \\
 & & & & are $\xbq$'s (or $\neg \xbq$'s) & & & of paths & \\

\hline

Argumentation & 1 & CPM & 2 & ``soft'' implication & acyclic? & collect &
& Same as block \\

\hline

Causation & 1 & CPM & 2 & ``normally causes'' & acyclic? & collect & & Same as
block \\
 & & & & ``normally inhibits'' & & & & \\

\hline

Block & 1 &? & 2 & enable/inhibit & acyclic? & collect & inhib. & What is
reachable and not blocked \\
 & & & & & & & win & Not monotonic, not simple summary \\

\hline

Typed & 2 &?, & 2 & inhib. go from & acyclic? & collect & inhib. & Same as block
\\
block & & causation & & causation to & & & block  & \\
 & & nodes & & causation & & & specific & \\
 & & & & & & & causes & \\

\hline

\end{tabular}

}

\end{turn}

\section{
A discussion of concepts
}
\label{Section Discussion-Concepts}

The aim of this Section is to describe the concepts and roles of models
and
operators in various propositional logics.

\bn

$\hspace{0.01em}$

% (+++ Orig. No.:  Notation: +++)

\label{Notation:}

We will use $ \xEC $ for the global universal modal quantifier: $ \xEC
\xbf $ holds in a model
iff $ \xbf $ holds everywhere - it is the dual of consistency.

$ \xcX $ and $ \xcx $ are the usual universal and existential modal
quantifiers, recall
that $ \xeA $ is some normality quantifier, see e.g.
Chapter \ref{Chapter Size} (page \pageref{Chapter Size}).
% Basic semantic entities, truth values, and operators
% Basic semantic entities, truth values, and operators
% ----------------------------------------------------

\subsection{
Basic semantic entities, truth values, and operators
}
\label{Section 1.2.2.1}

% The levels of language and semantics
% The levels of language and semantics
% ------------------------------------
\subsubsection{
The levels of language and semantics
}

\en

We have several levels:

(1) the language and the truth values

(2) the basic semantical entities, e.g. classical models, maximal
consistent
sets of modal formulas, etc.

(3) abstract or algebraic semantics, which describe the interpretation of
the
operators of the language in set-theoretic terms, like the interpretation
of $ \xcu $ by $ \xcs,$ $ \xeA $ ( ``the normal cases of'' ) by $ \xbm $
(choice of minimal
elements),
etc. These semantics do not indicate any mechanism which generates these
abstract operators.

(4) structural semantics, like Kripke structures, preferential structures,
which
give such mechanisms, and generate the abstract behaviour of the
operators.
They are or should be the intuitive basis of the whole enterprise.

(In analogue, we have a structural, an abstract, and a logical limit - see
Section \ref{Section Limit} (page \pageref{Section Limit}).)
% Language and truth values
% Language and truth values
% -------------------------
\subsubsection{
Language and truth values
}

A language has
 \xEI
 \xDH variable parts, like propositional variables

 \xDH constant parts, like
 \xEI
 \xDH operators, e.g. $ \xcu,$ $ \xeA,$ $ \xcX,$ etc, and
 \xDH relations like a consequence relation $ \xcn,$ etc.

 \xEJ
 \xEJ

Operators may have a
 \xEI
 \xDH unique interpretation, like $ \xcu,$ which is always interpreted by
$ \xcs,$
 \xDH or only restrictions on the interpretation, like $ \xeA,$ $ \xcX $
 \xEJ

Operators and relations may be
 \xEI
 \xDH nested, like $ \xcu,$ $ \xeA,$
 \xDH or only flat (perhaps on the top level), like $ \xcn $
 \xEJ

The truth values are part of the overall framework. For the moment, we
will
tacitly assume that there is only TRUE and FALSE. We will see in (1.3.1)
that
this restriction is unimportant for our purposes.
% Basic semantical entities
% Basic semantical entities
% -------------------------
\subsubsection{
Basic semantical entities
}

The language speaks about the basic semantic entities.
Note that the language will usually NOT speak about the relation of a
Kripke or
preferential structure, etc., only about the resulting function, resp. the
operator which is interpreted by the function, as part of formulas - we do
not
speak directly about operators, but only as part of formulas.

For the same language, there may be different semantic entities.
The semantic entities are (perhaps consistent, perhaps complete wrt. the
logic)
sets of formulae of the language. They are descriptions - in the language
- of
situations. They are NOT objects (in FOL we have names for objects), nor
situations, but only descriptions, even if it may help to consider them as
such
objects, (but unreal ones, just as mannequins in shop windows are unreal -
they
are only there to exhibit the garments.)

An example for different semantic entities is intuitionist logics, where
we may
take
 \xEI
 \xDH knowledge states, which may be incomplete (this forces the relation
in
Kripke structures to be monotonic), or
 \xDH classical models, where $ \xcX $ codes knowledge, and,
automatically, its
growth.
 \xEJ
(Their equivalence is a mathematical result, the former approach is
perhaps
the philosophically better one, the second one easier to handle, as
intuitionistic formulas are distinct by the preceding $ \xcX.)$

The entities need not contain all formulas with all operators, perhaps
they are
only sets of propositional variables, with no operators at all, perhaps
they
are all consistent sets of formulas of some sublanguage. For classical
logic,
we can take as basic entities either just sets of propositional variables,
or
maximal consistent sets of formulas.

In the case of maximal consistent formula sets in the full language,
perhaps the
simplest way to find all semantic entities is:

 \xEI
 \xDH Any set of formulas is a candidate for a semantic entity.

 \xDH If we want them to be complete, eliminate those who are not.

 \xDH Eliminate all which are contradictory under the operators and
the logic which governs their behaviour - e.g. $p,$ $q,$ and $ \xCN p \xco
\xCN q$
together cannot hold in classical models.
 \xEJ

In this approach, the language determines all situations, the logic those
which
are possible. We thus have a clean distinction between the work of the
language
and that of the logic.

For preferential reasoning (with some relation $ \xcn $ outside
the ``core language'' ),
we may again take all classical models - however defined - and introduce
the
interpretation of $ \xcn $ in the algebraic or abstract superstructure
(see below),
but we may also consider a normality operator $ \xeA $ directly in the
language,
and all consistent sets of such formulas (see e.g.  \cite{SGMRT00}).
Our picture is large enough to admit both possibilities.

In modal logic, we may again consider classical models, or, as is mostly
done,
(maximal consistent) sets of formulas of the full language.

The choice of these entities is a philosophical decision, not dictated by
the
language, but it has some consequences. - See below (nested preferential
operators) for details.
We call these entities models or basic models.
% Several truth values in basic semantic entities
% Several truth values in basic semantic entities
% -----------------------------------------------
\subsubsection{
Several truth values in basic semantic entities
}

When we consider sets of formulas as basic models, we assume implicitly
two
truth values: everything which is in the set, has truth value TRUE, the
rest
truth value FALSE (or undecided - context will tell). Of course, we can
instead consider (partial) functions from the set of formulas into any set
of truth values - this does not change the overall approach. E.g., in an
information state, we might have been informed that $ \xbf $ holds with
some
reliability $r,$ in another one with reliability $r',$ etc., so $r,$ $r'
$ etc. may be
truth values, and even pairs $\{r,r' \}$ when we have been informed with
reliability
$r$ that $ \xbf,$ and with reliability $r' $ that $ \xCN \xbf.$ Whatever
is reasonable in the
situation considered should be admitted as truth value.
% Algebraic and structural semantics
% Algebraic and structural semantics
% ----------------------------------
\subsection{
Algebraic and structural semantics
}
\label{Section 1.2.2.2}

We make now a major conceptual distinction, between an ``algebraic'' and a
``structural'' semantics, which can best be illustrated by an example.

Consider nonmonotonic logics as discussed above. In preferential
structures, we
only consider
the minimal elements, say $ \xbm (X),$ if $X$ is a set of models.
Abstractly, we thus
have a choice function $ \xbm,$ defined on the power set of the model
set, and $ \xbm $ has
certain properties, e.g. $ \xbm (X) \xcc X.$ More important is the
following property:
$X \xcc Y$ $ \xcp $ $ \xbm (Y) \xcs X \xcc \xbm (X).$ (The proof is
trivial: suppose there were $x \xbe \xbm (Y) \xcs X,$
$x \xce \xbm (X).$ Then there must be $x' \xeb x,$ $x' \xbe X \xcc Y,$ but
then $x$ cannot be minimal in $Y.)$

Thus, all preferential structures generate $ \xbm $ functions with certain
properties,
and once we have a complete list, we can show that any arbitrary model
choice
function with these properties can be generated by an appropriate
preferential
structure.

Note that we do not need here the fact that we have a relation between
models,
just any relation on an arbitrary set suffices. It seems natural to call
the complete list of properties of such $ \xbm -$functions an algebraic
semantics,
forgetting that the function itself was created by a preferential
structure,
which is the structural semantics.

This distinction is very helpful, it not
only incites us to separate the two semantics conceptually, but also to
split
completeness proof in two parts: One part, where we show correspondence
between
the logical side and the algebraic semantics, and a second one, where we
show
the correspondence between the algebraic and the structural semantics. The
latter part will usually be more difficult, but any result obtained here
is
independent from logics itself, and can thus often be re-used in other
logical
contexts. On the other hand, there are often some subtle problems for the
correspondence between the logics and the algebraic semantics (see
definability
preservation, in particular the discussion in  \cite{Sch04}),
which we can then more clearly isolate, identify, and solve.
% Abstract or algebraic semantics
% Abstract or algebraic semantics
% -------------------------------
\subsubsection{
Abstract or algebraic semantics
}

$\hspace{0.01em}$

In all cases, we see that the structural semantics define a set operator,
and
thus an algebraic semantics:
 \xEI
 \xDH in nonmonotonic logics (and Deontic Logic), the function chooses the
minimal
(morally best) models, a subset, $ \xbm (X) \xcc X$
 \xDH in (distance based) Theory Revision, we have abinary operator, say $
\xfA $
which chooses the $ \xbf -$models closest to the set of $K-$models: $M(K)
\xfA M( \xbf )$
 \xDH in Theory Update, the operator chooses the $i-$th coordinate of all
best sequences
 \xDH in the Logic of Counterfactual Conditionals, whave again a binary
operator
$m \xfA M( \xbf )$ which chooses the $ \xbf -$models closest to $m,$ or,
when we consider
a whole set $X$ of models as starting points $X \xfA M( \xbf )= \xcV \{m
\xfA M( \xbf ):m \xbe X\}.$
 \xDH in Modal and Intuitionistic Logic, seen from some model $m,$ we
choose a subset
of all the models (thus not a subset of a more restricted model set),
those
which can be reached from $m.$
 \xEJ

Thus, in each case, the structure ``sends'' us to another model set, and
this
expresses the change from the original situation to the ``most plausible'',
``best'', ``possible'' etc. situations. It seems natural to call all such
logics
``generalized modal logics'', as they all use the idea of a model choice
function.

(Note again that we have neglected here the possibility that there are no
best
or closest models (or sequences), but only ever better ones.)

Abstract semantics are interpretations of the operators of the language
(all,
flat, top level or not) by functions (or relations in the case of $ \xcn
),$ which
assign to sets of models sets of models, $ \xdo: \xdp ( \xdm ) \xcp \xdp
( \xdm )$ - $ \xdp $ the power set
operator, and $ \xdm $ the set of basic models -, or binary functions for
binary
operators, etc.

These functions are determined or restricted by the laws for the
corresponding
operators. E.g., in classical, preferential, or modal logic, $ \xcu $ is
interpreted by $ \xcs,$ etc.; in preferential logic $ \xeA $ by
$ \xbm;$ in modal logic, we interpret $ \xcX,$ etc.

Operators may be truth-functional or not. $ \xCN $ is truth-functional. It
suffices
to know the truth value of $ \xbf $ at some point, to know that of $ \xCN
\xbf $ at the same
point. $ \xcX $ is not truth-functional: $ \xbf $ and $ \xbq $ may hold,
and $ \xcX \xbf,$ but not $ \xcX \xbq,$
all at the same point $(=$ base model), we have to look at the full
picture,
not only at some model.

We consider first those operators, which have a unique possible
interpretation, like $ \xcu,$ which is interpreted by $ \xcs,$ $ \xCN $
by $ \xdC,$ the set
theoretic complement, etc.
$ \xeA $ (standing for ``most'', ``the important'', etc.) e.g.
has only restrictions to its interpretation, like $ \xbm (X) \xcc X,$ etc.
Given a set of
models without additional structure, we do not know its exact form, we
know it only once we have fixed the additional structure (the relation in
this
case).

If the models contain already the operator, the function will respect it,
i.e.
we cannot have $ \xbf $ and $ \xCN \xbf $ in the same model, as $ \xCN $
is interpreted by $ \xdC.$ Thus,
the functions can, at least in some cases, control consistency.

If, e.g. the models contain $ \xcu,$ then we have two ways to evaluate $
\xbf \xcu \xbq:$ we
can first evaluate $ \xbf,$ then $ \xbq,$ and use the function for $
\xcu $ to evaluate $ \xbf \xcu \xbq.$
Alternatively, we can look directly at the model for $ \xbf \xcu \xbq $ -
provided we
considered the full language in constructing the models.

As we can apply one function to the result of the other, we can evaluate
complicated formulas, using the functions on the set of models.
Consequently,
if $ \xcn $ or $ \xeA $ is evaluated by $ \xbm,$ we can consider $ \xbm (
\xbm (X))$ etc., thus, the
machinery for the flat case gives immediately an interpretation for nested
formulas, too - whether we looked for it, or not.

As far as we see, our picture covers the usual presentations of classical
logic, preferential, intuionist, and modal logic, but also of linear logic
(where we have more structure on the set of basic models, a monoid, with a
distinct set $ \xcT,$ plus some topology for! and? - see below), and
quantum logic
a la Birkhoff/von Neumann.

We can introduce new truth-functional operators into the language as
follows:
Suppose we have
a distinct truth value TRUE, then we may define $ \xdo_{X}( \xbf )=TRUE$
iff the truth-value
of $ \xbf $ is an element of $X.$ This might sometimes be helpful. Making
the truth
value explicit as element of the object language may facilitate the
construction
of an accompanying proof system - experience will tell whether this is the
case.
In this view, $ \xCN $ has now a double meaning in the classical
situation: it is an
operator for the truth value ``false'', and an operator on the model set,
and
corresponds to the complement. ``Is true'' is the identical truth functional
operator, $is-true( \xbf )$ and $ \xbf $ have the same truth value.

If the operators have a unique interpretation, this might be all there is
to
say in this abstract framework.
(This does not mean that it is impossible to introduce new operators which
are
independent from any additional structure, and based only on the set of
models
for the basic language. We can, for instance, introduce a ``CON'' operator,
saying
that $ \xbf $ is consistent, and $CON( \xbf )$ will hold everywhere iff $
\xbf $ is consistent,
i.e. holds in at least one model. Or, for a more bizarre example, a 3
operator,
which says that $ \xbf $ has at least 3 models (which is then dependent on
the
language). We can also provide exactly one additional structure, e.g. in
the
following way: Introduce a ranked order between models as follows: At the
bottom, put the single model which makes all propositional variables true,
on
the next level those which make exactly one propositional variable true,
then
two, etc., with the model making all false on top. So there is room to
play, if
one can find many useful examples is another question.)

If the operator has no unique interpretation (like $ \xeA,$ $ \xcX,$
etc., which are only
restricted, e.g. by $ \xEC ( \xbf \xcp \xbq )$ $ \xcp $ $ \xEC ( \xeA \xbq
\xcu \xbf \xcp \xeA \xbf ),$ the situation seems more
complicated, and is discussed below
in Section \ref{Section Restricted-Operators} (page \pageref{Section
Restricted-Operators}).

It is sometimes useful to consider the abstract semantics as a (somehow
coherent) system
of filters. For instance, in preferential structures, $ \xbm (X) \xcc X$
can be seen as
the basis of a principal filter. Thus, $ \xbf \xcn \xbq $ iff $ \xbq $
holds in all minimal models
of $ \xbf,$ iff there is a ``big'' subset of $M( \xbf )$ where $ \xbq $
holds, recalling that a
filter is an abstraction of size - sets in the filter are big, their
complements
small, and the other sets have medium size. Thus, the ``normal'' elements
form the
smallest big subset. Rules like $X \xcc Y$ $ \xcp $ $ \xbm (Y) \xcs X \xcc
\xbm (X)$ form the coherence
between the individual filters, we cannot choose them totally
independently.
Particularly for preferential structures, the reasoning with small and big
subsets can be made very precise and intuitively appealing, and we will
come
back to this point later.
We can also introduce a generalized quantifier, say $ \xeA,$ with the
same meaning,
i.e. $ \xbf \xcn \xbq $ iff $ \xeA ( \xbf ). \xbq,$ i.e. ``almost
everywhere'', or ``in the important cases''
where $ \xbf $ holds, so will $ \xbq.$ This is then the syntactic
analogue of the
semantical filter system.
These aspects are discussed in detail in
Chapter \ref{Chapter Size} (page \pageref{Chapter Size}).
% Structural semantics
% Structural semantics
% --------------------
\subsubsection{
Structural semantics
}

Structural semantics generate the abstract or algebraic semantics, i.e.
the
behaviour of the functions or relations (and of the operators in the
language
when we work with ``rich'' basic models). Preferences between models
generate
corresponding $ \xbm -$functions, relations in Kripke structures generate
the
functions corresponding to $ \xcX -$operators, etc.

Ideally, structural semantics capture the essence of what we want to
reason and speak about (beyond classical logic), they come, or should
come,
first. Next, we try to see
the fundamental ingredients and laws of such structures, code them in an
algebraic semantics and
the language, i.e. extract the functions and operators, and their laws. In
a
backward movement, we make the roles of the operators (or relations)
precise
(should they be nested or not?, etc.), and define the basic models and the
algebraic operators. This may result in minor modifications of the
structural
semantics (like introduction of copies), but should still be close to the
point of outset. In this view, the construction of a logic is a
back-and-forth
movement.
% Restricted operators and relations
% Restricted operators and relations
% ----------------------------------
\subsection{
Restricted operators and relations
}
\label{Section Restricted-Operators}

We discuss only operators, relations seem to be similar. The discussion
applies
as well to abstract as to structural semantics.

An operator, which is only restricted in his behaviour, but not fully
defined,
has to be interpreted by a unique function. Thus, the interpretation will
be
more definite than the operator. It seems that the problem has no
universal
solution.

(1) If there is tacitly a ``best choice'', it seems natural to make this
choice.
At the same time, such a best choice may also serve to code our ignorance,
without enumerating all possible cases among which we do not know how to
decide.

For instance, in reasoning about normality (preferential structures), the
interpretation which makes nothing more normal than explicitly required -
corresponding to a set-theoretically minimal relation - seems a natural
choice. This will NOT always give the same result as a disjunction over
all
possible interpretations: e.g., if the operator is in the language, and we
have
finitely many possibilities, we can express them by ``or'', and this need
not
be the same as considering the unique minimal solution. (Note that this
will
usually force us to consider ``copies'' in preferential structures - see
below.)

(2) We can take all possible interpretations, and consider them
separately, and take as result only those facts which hold in all
possibilities.
Bookkeeping seems difficult, especially when we have nested operators,
which
have all to be interpreted in the various ways. In a second step, we can
unite
all possibilities in one grand picture (a universal structure, as it
permits
to find exactly all consequences in one construction, and not in several
ones
as is done for classical logic), essentially by a disjoint union - this
was
done (more or less) in the authors'  \cite{SGMRT00} for preferential
structures.

(3) We can work with basic models already in the full language, and
capture
the different possibilities already on the basic level. The interpretation
of
the operator will then be on the big set of models for the full language,
which serve essentially as bookkeeping device for the different
possibilities
of interpretation - again an universal structure. This is done in the
usual completeness proofs for modal logic.
% Copies in preferential models
% Copies in preferential models
% -----------------------------
\subsection{
Copies in preferential models
}
\label{Section 1.2.2.4}

Copies in preferential structures (variant (2) in
Section \ref{Section Restricted-Operators} (page \pageref{Section
Restricted-Operators}) )
thus seem to serve to
construct universal structures, or code our ignorance, i.e. we know that
$x$ is
minimized by $X,$ but we do not know by which element of $X,$ they are in
this view
artificial. But they have an intuitive justification, too:
They allow minimization by sets of other elements only.
We may consider an element $m$ only abnormal in the
presence of several other elements together. E.g., considering penguins,
nightingales, woodpeckers, ravens, they all have some exceptional
qualities,
so we may perhaps not consider a penguin more abnormal than a woodpecker,
etc.,
but seen all these birds together, the penguin stands out as the most
abnormal
one. But we cannot code minimization by a set, without minimization by its
elements, without the use of copies. Copies will then code the different
aspects of abnormality.
% Further remarks on universality of representation proofs
% Further remarks on universality of representation proofs
% --------------------------------------------------------
\subsection{
Further remarks on universality of representation proofs
}
\label{Section 1.2.2.5}

There is a fundamental difference between considering all possible
ramifications
and coding ignorance. For instance, if we know that $\{a,b,c\}$ is
minimized by
$\{b,c\},$ we can create two structures, one, where a is minimized by $b,$
the other,
where a is minimized by $c.$ These are all possible ramifications (if
$ \xbm (\{a\}) \xEd \xCQ ).$ Or, we can code our ignorance with copies of
a, as is done in our
completeness
constructions. $((\{a,b\} \xcn b)$ or $(\{a,c\} \xcn c))$ is different
from $\{a,b,c\} \xcn \{b,c\},$ and
if the language is sufficient, we can express this. In a ``directly
ignorant''
structure, none of the two disjoints hold, so the disjunction will fail.

Our proofs try to express ignorance directly. Note that a representing
structure
can be optimal in several ways: (1) optimal, or universal, as it expresses
exactly the logic, (2) optimal, as it has exactly the required properties,
but
not more. For instance, a smooth structure can be optimal, as it expresses
exactly the logic it is supposed to code, or optimal, as it preserves
exactly
smoothness, there is no room to move left. Usually, one
seeks only the first variant. If both variants disagree, then structure
and
model do not coincide exactly, the structure still has more space to move
than
necessary for representation.

In our constructions, all possibilities are coded into the choice
functions
(the indices), but as they are not directly visible to the language, they
are
only seen together, so there is no way to analyze the ``or'' of the
different
possibilities separately.
% %n in the object language?
% %n in the object language?
% --------------------------
\subsection{
$\xcn$ in the object language?
}
\label{Section 1.2.2.6}

It is tempting to try and put a consequence relation $ \xcn $ into the
object
language by creating a new modal operator $ \xeA $ (expressing ``most'' or
so), with the
translation $ \xbf \xcn \xbq $ iff $ \xcl \xeA \xbf \xcp \xbq.$

We examine now this possibility and the resulting consequences.

We suppose that $ \xcp $ will be interpreted in the usual way, i.e. by the
subset
relation. The aim is then to define the interpretation of $ \xeA $ s.t. $
\xbf \xcn \xbq $ iff
$M( \xeA \xbf ) \xcc M( \xbq )$ - $M( \xbf )$ the set of models of $ \xbf
.$

It need not be the case - but it will be desirable - that $ \xeA $ is
insensitive
to logical equivalence. $ \xeA \xbf $ and $ \xeA \xbf ' $ may well be
interpreted by different
model sets, even if $ \xbf $ and $ \xbf ' $ are interpreted by the same
model sets. Thus,
we need not left logical equivalence - whatever the basic logic is. On the
right hand side, as we define $ \xcn $ via $ \xcp,$ and $ \xcp $ via
model subsets, we will
have closure under semantic consequence (even infinite closure, if this
makes
a difference). It seems obvious, that this is the only property a relation
$ \xcn $
has to have to be translatable into object language via ``classical'' $ \xcp
,$ or,
better, subsets of models.

Note that the interpretation of $ \xeA \xbf $ can be equivalent to a
formula, a theory,
or to a set of models logically equivalent to a formula or theory, even if
some
models are missing (lack of definability preservation). Likewise, we may
also
consider $ \xeA T$ for a full theory $T.$ The standard solution will, of
course, be
$M( \xeA \xbf )$ $:=$ $ \xcS \{M( \xbq ): \xbf \xcn \xbq \}.$
% Possibilities and problems of external %n vs. internal #A
% Possibilities and problems of external %n vs. internal #A
% ---------------------------------------------------------
\subsubsection{
Possibilities and problems of external $\xcn$ vs. internal $\xeA$
}

We first enumerate a small number of various differences, before we turn
to
a few deeper ones.

 \xEI

 \xDH
The translation of $ \xcn $ into the object language - whenever possible -
introduces
contraposition into the logic, see also  \cite{Sch04} for a
discussion.

 \xDH
It may also help to clarify questions like validity of the deduction
theorem
(we see immediately that one half is monotony), cut, situations like $
\xeA \xbf \xcp \xbq $
and $ \xeA \xbq \xcp \xbf,$ etc.

 \xDH
It is more usual to consider full theories outside the language, even if,
a
priori, there is no problem to define a formula $ \xeA T$ using a theory
$T$ in the
inductive step. Thus, the external variant can be more expressive in one
respect.

 \xDH
At least in principle, having $ \xeA $ inside the language makes it more
amenable
to relativize it to different viewpoints, as in modal logic. It seems at
least
more usual to write $m \xcm \xeA \xbf \xcp \xbq $ than to say that in $m,$
$ \xbf \xcn \xbq $ holds.

 \xEJ
% Disjunction in preferential and modal structures
% Disjunction in preferential and modal structures
% ------------------------------------------------
\subsubsection{
Disjunction in preferential and modal structures
}

In a modal Kripke model, $ \xcX \xbf \xco \xcX \xbq $ may hold everywhere,
but neither $ \xcX \xbf $ nor
$ \xcX \xbq $ may hold everywhere. This is possible, as formulas are
evaluated locally,
and at one point $m,$ we may make $ \xcX \xbf $ hold, at another $m' $ $
\xcX \xbq.$

This is not the case for the globally evaluated modal operator $ \xEC.$
Then, in $ \xCf one$ structure, if $ \xEC ( \xbf ) \xco \xEC ( \xbq )$
holds, either $ \xEC ( \xbf )$ or $ \xEC ( \xbq )$
holds, but a structure where $ \xEC ( \xbf )$ (or $ \xEC ( \xbq ))$ holds,
has more information than
a structure in which $ \xEC ( \xbf ) \xco \xEC ( \xbq )$ holds.
Consequently, one Kripke structure
is not universal any more for $ \xEC $ (and $ \xco ),$ we need several
such structures to
represent disjunctions of $ \xEC.$

The same is the case for preferential structures. If we put $ \xcn $ into
the language
as $ \xeA,$ to represent disjunctions, we may need several preferential
structures:
$( \xeA \xbf \xcp \xbq ) \xco ( \xeA \xbf \xcp \xbq ' )$ will only be
representable as $ \xeA \xbf \xcp \xbq $ or as $ \xeA \xbf \xcp \xbq ',$
but in one structure, only one of them may hold. Thus, again, $ \xCf one$
structure
will say more than the original formula $( \xeA \xbf \xcp \xbq ) \xco (
\xeA \xbf \xcp \xbq ' ),$ and to express
this formula faithfully, we will need again two structures. Thus, putting
$ \xcn $
into the object language destroys the universality of preferential
structures
by its richer language.

(Remark: rational monotony is, semantically, also universally quantified,
$ \xba \xcn \xbg $ $ \xcp $ $ \xba \xcu \xbb \xcn \xbg $ or $ \xCf
everywhere$ $ \xba \xcn \xCN \xbb.$
% Iterated #A in preferential structures
% Iterated #A in preferential structures
% --------------------------------------
\subsubsection{
Iterated $\xeA$ in preferential structures
}

Once we have $ \xeA $ in the object language, we can form $ \xeA \xeA \xbf
$ etc.

When we consider preferential structures in the standard way, it is
obvious that
$ \xeA \xeA \xbf = \xeA \xbf,$ etc. will hold.

But, even in a preferential context, it is not obvious that this has to
hold,
it suffices to interpret the relation slightly differently. Instead of
setting
$ \xbm (X)$ the set of minimal elements of $X,$ it suffices to define $
\xbm (X)$ the set of
non-worst elements of $X,$ i.e. everything except the upper layer. (One of
the
authors
once had a similar discussion with M.Magidor.) (Note that, in
this interpretation, for instance $ \xeA \xeA \xbf $ may seem to be the
same as $ \xeA \xbf,$ but
$ \xeA \xeA \xeA \xbf $ clearly is different from $ \xeA \xeA \xbf:$ In
going from $ \xeA \xbf $ to $ \xeA \xeA \xbf,$ we
loose some models, but not enough to be visible by logics - a problem of
definability preservation. The loss becomes visible in the next step.)

But, we can have a similar result even in the usual interpretation of
``normal''
by ``best'':

Consider for any $X \xcc \xbo $ the logic $ \xcn_{X}$ defined by the
axioms
$\{A_{i}:i \xbe X\},$ where $A_{i}:= \xeA^{i+1} \xbf \xcr \xeA^{i} \xbf,$
and $ \xeA^{i} \xbf $ is, of course, $i$ many $ \xeA ' $s, followed by $
\xbf,$ etc., plus the usual
axioms
for preferential logics. This defines $2^{ \xbo }$ many logics, and we
show that they
are all different. The semantics we will give show at the same time that
the
usual axioms for preferential logics do not entail $ \xeA \xeA \xbf \xcr
\xeA \xbf.$

For simplicity, we first show that the $ \xbo $ many logics defined
by the axioms $B_{i}$ $=$ $ \xeA^{i+1} \xbf \xcr \xeA^{i} \xbf $ for
arbitrary $ \xbf $ are different.
We consider sequences of $ \xbo $ many preferential structures over some
infinite
language, s.t. we choose exactly at place $i$ the same structure twice,
and all
the other times different structures. Let $S_{i}$ be the structure which
minimizes
$ \xCN p_{i}$ to $p_{i},$ i.e. every $p_{i}-model$ $m$ is smaller than its
opposite $m',$ which is like
$m,$ only $m' \xcm \xCN p_{i}.$ It is obvious, that the associated $ \xbm
-$functions $ \xbm_{i}$ will
give different result on many sets $X$ (if a model $x \xbe X$ is minimized
at all in $X,$
it will be minimized in different ways). Consider now e.g. a formula $
\xeA \xeA \xbf.$
We start to evaluate at $S_{0},$ evaluating the leftmost $ \xeA $ by $
\xbm_{0}.$ The second $ \xeA $
will be evaluated at $S_{1},$ by $ \xbm_{1}.$ If, for instance, $ \xbf $
is a tautology, we
eliminate in the first step the $ \xCN p_{0}-$models, in the second step
the $ \xCN p_{1}-models,$
so $ \xeA \xeA $ is not equivalent to $ \xeA.$ If, instead of taking at
position $i+1$ $S_{i+1},$
we just take $S_{i}$ again, then axiom $B_{i}$ will hold, but not the
other $B_{j},$ $i \xEd j.$
Thus, the logics $B_{i}$ are really different.

In the first case, for the $A_{i},$ we repeat $S_{j}$ for all $j \xbe X,$
instead of taking
$S_{j+1},$ and start evaluation again at $S_{0}.$ Again, the different
sequences of
structures will distinguish the different logics, and we are done.

Yet the standard interpretation of $ \xeA $ in one single preferential
structure does
not allow to distinguish between all these logics, as sequences of $ \xeA
' $s will
always be collapsed to one $ \xeA.$ As long as a (preferential)
consequence relation
satisfies this, it will hold in the standard interpretation, if not, it
will
fail.

This is another example which shows that putting $ \xcn $ as $ \xeA $ in
the object language
can increase substantially the expressiveness of the logics.

Compare the situation to modal logic and Kripke semantics. In Kripke
semantics, any evaluation of $ \xcX $ takes us to different points in the
structure,
and, if the relation is not transitive, what is beyond the first step is
not
visible from the origin. In preferential structures, this ``hiding'' is not
possible, by the very definition of $ \xbm,$ we have to follow any chain
as far as
it goes.

(One of the authors had this idea when looking at the very interesting
article
``Mathematical modal logic: a view of its evolution'' by Robert Goldblatt,
 \cite{Gol03},
and begun to mull a bit over the short remark it contains on the continuum
many
modal logics indistinguishable by Kripke structures.)
% Various considerations on abstract semantics
% Various considerations on abstract semantics
% --------------------------------------------
\subsection{
Various considerations on abstract semantics
}
\label{Section 1.2.2.7}

We can think of preferential structures as elements attacking each
other: if $x \xeb y,$ then $x$ attacks $y,$ so $y$ cannot be minimal any
more. Thus the
non-monotonicity.

In topology e.g., things are different. Consider the open interval
$(0,1),$ and
its closure $[0,1].$ Sequences converging to 0 or 1 ``defend'' 0 or 1 -
thus, the
more elements there are, the bigger the closure will be, monotony. The
open core
has a similar property.
The same is true for e.g. Kripke structures for modal logic: If $y$ is in
relation
$R$ with $x,$ then any time $y$ is there, $x$ can be reached: $y$ defends
$x.$

A neuron may have inhibitory (attackers) and excitatory (defending)
inputs.
In preferential models, we may need many attackers to destroy minimality
of one
model, provided this model occurs in several copies.

Of course, in a neuron system, there are usually many attackers and many
defenders, so we have here a rather complicated system.

Abstractly, both defense and attack are combined in Gabbay's reactive
diagrams,
see  \cite{Gab04},
and Section \ref{Section IBRS} (page \pageref{Section IBRS}).

Now, back to our
% Set functions
% Set functions
% -------------
\subsubsection{
Set functions
}

We can make a list of possible formal properties of such functions, which
might
include (U is the universe or base set):

$f(X)=X,$ $f(X) \xcc X,$ $X \xcc f(X)$

$f(X)=f(f(X)),$ $f(X) \xcc f(f(X)),$ $f(f(X)) \xcc f(X)$

$X \xEd \xCQ $ $ \xcp $ $f(X) \xEd \xCQ $

$X \xEd U$ $ \xcp $ $f(X) \xcs X= \xCQ $

$X \xEd \xCQ $ $ \xcp $ $f(X) \xcs X \xEd \xCQ $

$f( \xdC (X))= \xdC (f(X))$

$f(X \xcv Y)=f(X) \xcv f(Y)$

$f(X \xcv Y)=f(X)$ or $f(Y)$ or $f(X) \xcv f(Y)$ (true in ranked
structures)

$X \xcc Y$ $ \xcp $ $f(Y) \xcs X \xcc f(X)$ (attack, basic law of
preferential structures)

$X \xcc Y$ $ \xcp $ $f(X) \xcc f(Y)$ (defense, e.g. Modal Logic)

$f(X) \xcc Y \xcc X$ $ \xcp $ $f(X)=f(Y)$ (smoothness), holds also for the
open core of a set
in topology

$X \xcc Y \xcc f(X)$ $ \xcp $ $f(X)=f(Y)$ counterpart, holds for
topological closure

(Note that the last two properties will also hold in all other situations
where
one chooses the biggest subset or smallest superset from a set of
candidates.)

$f( \xcV \xdx )= \xcV \{f(X):X \xbe \xdx \}$

etc.

General, distance based revision is a two argument set function (in the
AGM approach, it has only one argument):

$M(K) \xfA M( \xbf ) \xcc M( \xbf )$

This is non-monotonic in both arguments, as the result depends more on the
``shape'' of both model sets than their size.

Counterfactuals (and update) are also two argument set functions:

$M( \xbf ) \xfB M( \xbq )$ is defined as the set of $ \xbq -$models
closest to some individual
$ \xbf -$model,
here the function is monotonic in the first argument, and non-monotonic in
the
second - we collect for all $m \xbe M( \xbf )$ the closest $ \xbq
-$models.
% A comparison with Reiter defaults
% A comparison with Reiter defaults
% ---------------------------------
\subsection{
A comparison with Reiter defaults
}
\label{Section 1.2.2.8}

The meaning of Reiter defaults differs from that of preferential
structures
in a number of aspects.

(1) The simple (Reiter) default ``normally $ \xbf $ '' does not only mean
that in normal
cases $ \xbf $ holds, but also that, if $ \xbq $ holds, then normally also
$ \xbf \xcu \xbq $ holds.
It thus inherits ``normally $ \xbf $ '' down on subsets.

(2) Of course, this is itself a default rule, as we might have that for
$ \xbq -$cases, normally $ \xCN \xbf $ holds. But this is a meta-default.

(3) Defaults can also be concatenated, if normally $ \xbf $ holds, and
normally, if
$ \xbf $ holds, then also $ \xbq $ holds, we conclude that normally $ \xbf
\xcu \xbq $ holds. Again,
this is a default rule.

Thus, Reiter defaults give us (at least) three levels of certainty:
classical
information, the information directly expressed by defaults, and the
information
concluded by the usual treatment of defaults.
\section{
IBRS
}
\label{Section IBRS}

\subsection{
Definition and comments
}
\index{Definition IBRS}

\bd

$\hspace{0.01em}$

% (+++ Orig. No.:  Definition IBRS +++)

\label{Definition IBRS}

 \xEh
 \xDH An information bearing binary relation frame IBR, has the form
$(S, \xdR ),$ where $S$ is a non-empty set and $ \xdR $ is a subset of
$S,$ where $S$ is defined by induction as follows:

 \xEh

 \xDH $S_{0}=S$

 \xDH $S_{n+1}$ $=$ $S_{n} \xcv (S_{n} \xCK S_{n}).$

 \xDH $S$ $=$ $ \xcV \{S_{n}:n \xbe \xbo \}$

 \xEj

We call elements from $S$ points or nodes, and elements from $ \xdR $
arrows.
Given $(S, \xdR ),$ we also set $ \xdP ((S, \xdR )):=S,$ and $ \xdA ((S,
\xdR )):= \xdR.$

If $ \xba $ is an arrow, the origin and destination of $ \xba $ are
defined
as usual, and we write $ \xba:x \xcp y$ when $x$ is the origin, and $y$
the destination of the arrow $ \xba.$ We also write $o( \xba )$ and $d(
\xba )$ for
the origin and destination of $ \xba.$

 \xDH Let $Q$ be a set of atoms, and $ \xdL $ be a set of labels (usually
$\{0,1\}$ or $[0,1]).$ An information assignment $h$ on $(S, \xdR )$
is a function $h:Q \xCK \xdR \xcp \xdL.$

 \xDH An information bearing system IBRS, has the form
$(S, \xdR,h,Q, \xdL ),$ where $S,$ $ \xdR,$ $h,$ $Q,$ $ \xdL $ are as
above.

 \xEj

See Diagram \ref{Diagram IBRS} (page \pageref{Diagram IBRS})  for an
illustration.

\vspace{10mm}

% \begin{diagram}[ht]
% \bfc
\begin{diagram}

\centering
\setlength{\unitlength}{0.00083333in}
{\renewcommand{\dashlinestretch}{30}
\begin{picture}(4961,5004)(0,0)
\path(1511,1583)(611,3683)
\blacken\path(685.845,3584.520)(611.000,3683.000)(630.696,3560.885)(672.451,3539.613)(685.845,3584.520)
\path(1511,1583)(2411,3683)
\blacken\path(2391.304,3560.885)(2411.000,3683.000)(2336.155,3584.520)(2349.549,3539.613)(2391.304,3560.885)
\path(3311,1583)(4361,4133)
\blacken\path(4343.050,4010.616)(4361.000,4133.000)(4287.570,4033.461)(4301.603,3988.750)(4343.050,4010.616)
\path(3316,1574)(2416,3674)
\blacken\path(2490.845,3575.520)(2416.000,3674.000)(2435.696,3551.885)(2477.451,3530.613)(2490.845,3575.520)
\path(986,2783)(2621,2783)
\blacken\path(2501.000,2753.000)(2621.000,2783.000)(2501.000,2813.000)(2465.000,2783.000)(2501.000,2753.000)
\path(2486,2783)(2786,2783)
\blacken\path(2666.000,2753.000)(2786.000,2783.000)(2666.000,2813.000)(2630.000,2783.000)(2666.000,2753.000)
\path(3311,1583)(2051,2368)
\blacken\path(2168.714,2330.008)(2051.000,2368.000)(2136.987,2279.083)(2183.406,2285.509)(2168.714,2330.008)
\path(2166,2288)(1906,2458)
\blacken\path(2022.854,2417.439)(1906.000,2458.000)(1990.019,2367.221)(2036.567,2372.629)(2022.854,2417.439)

\put(1511,1358) {{\xssc $a$}}
\put(3311,1358) {{\xssc $d$}}
\put(3311,1058)  {{\xssc $(p,q)=(1,0)$}}
\put(1511,1058)  {{\xssc $(p,q)=(0,0)$}}
\put(2411,3758){{\xssc $c$}}
\put(4361,4433){{\xssc $(p,q)=(1,1)$}}
\put(4361,4208){{\xssc $e$}}
\put(2411,3983){{\xssc $(p,q)=(0,1)$}}
\put(611,3983) {{\xssc $(p,q)=(0,1)$}}
\put(611,3758) {{\xssc $b$}}
\put(1211,2883){{\xssc $(p,q)=(1,1)$}}
\put(260,2333) {{\xssc $(p,q)=(1,1)$}}
\put(2261,1583) {{\xssc $(p,q)=(1,1)$}}
\put(1286,3233){{\xssc $(p,q)=(1,1)$}}
\put(2711,3083){{\xssc $(p,q)=(1,1)$}}
\put(3836,2633){{\xssc $(p,q)=(1,1)$}}

\put(300,700)
{{\rm\bf
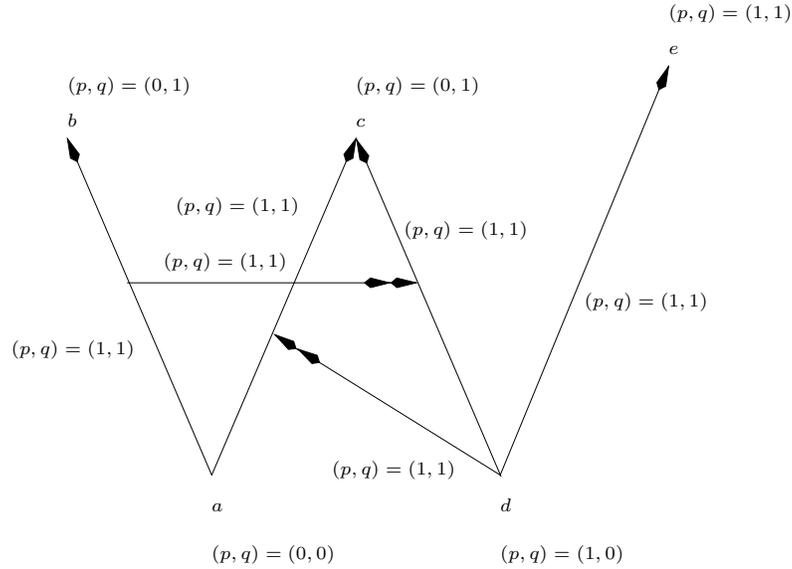
A simple example of an information bearing system.
}}

\end{picture}
}
\label{Diagram IBRS}
\index{Diagram IBRS}

\end{diagram}

We have here:
\[\begin{array}{l}
S =\{a,b,c,d,e\}.\\
\xdR = S \cup \{(a,b), (a,c), (d,c), (d,e)\} \cup \{((a,b), (d,c)),
(d,(a,c))\}.\\
Q = \{p,q\}
\end{array}
\]
The values of $h$ for $p$ and $q$ are as indicated in the figure. For
example $h(p,(d,(a,c))) =1$.

\vspace{4mm}

\index{Comment IBRS}

\ed

\bcom

$\hspace{0.01em}$

% (+++ Orig. No.:  Comment +++)

\label{Comment}

\label{Comment IBRS}

The elements in Figure Diagram \ref{Diagram IBRS} (page \pageref{Diagram IBRS}) 
can be interpreted in
many ways,
depending on the area of application.

 \xEh

 \xDH The points in $S$ can be interpreted as possible worlds, or
as nodes in an argumentation network or nodes in a neural net or
states, etc.

 \xDH The direct arrows from nodes to nodes can be interpreted as
accessibility relation, attack or support arrows in an argumentation
networks, connection in a neural nets, a preferential ordering in
a nonmonotonic model, etc.

 \xDH The labels on the nodes and arrows can be interpreted as fuzzy
values in the accessibility relation or weights in the neural net
or strength of arguments and their attack in argumentation nets, or
distances in a counterfactual model, etc.

 \xDH The double arrows can be interpreted as feedback loops to nodes
or to connections, or as reactive links changing the system which are
activated
as we pass between the nodes.

 \xEj
\index{IBRS as abstraction}
\label{IBRS as abstraction}

\ecom

Thus, IBRS can be used as a source of information for various logics
based on the atoms in $Q.$ We now illustrate by listing several such
logics.
\subsubsection*{Modal Logic}

One can consider the figure as giving rise to two modal logic models.
One with actual world a and one with $d,$ these being the two minimal
points of the relation. Consider a language with $ \xcX q.$ how
do we evaluate $a \xcm \xcX q?$

The modal logic will have to give an algorithm for calculating the
values.

Say we choose algorithm $ \xda_{1}$ for $a \xcm \xcX q,$ namely:

[
$ \xda_{1}(a, \xcX q)=1$
]
iff for all $x \xbe S$ such that $a=x$ or $(a,x) \xbe \xdR $ we have
$h(q,x)=1.$

According to $ \xda_{1}$ we get that $ \xcX q$ is false at a. $ \xda_{1}$
gives rise to a $T-$modal logic. Note that the reflexivity is not
anchored at the relation $ \xdR $ of the network but in the algorithm
$ \xda_{1}$ in the way we evaluate. We say $(S, \xdR, \Xl.) \xcm \xcX $
$q$ iff $ \xcX q$ holds in all minimal points of $(S, \xdR ).$

For orderings without minimal points we may choose a subset of
distinguished
points.
\subsubsection*{Nonmonotonic Deduction}

We can ask whether $p \xcn q$ according to algorithm $ \xda_{2}$ defined
below. $ \xda_{2}$ says that $p \xcn q$ holds iff $q$ holds in all
minimal models of $p.$ Let us check the value of $ \xda_{2}$ in this
case:

Let $S_{p}=\{s \xbe S \xfA h(p,s)=1\}.$ Thus $S_{p}=\{d,e\}.$

The minimal points of $S_{p}$ are $\{d\}.$ Since $h(q,d)=0,$ we have that
$p \xcN q.$

Note that in the cases of modal logic and
nonmonotonic logic we ignored the arrows $(d,(a,c))$ (i.e. the double
arrow from $d$ to the arc $(a,c))$ and the $h$ values to arcs. These
values do not play a part in the traditional modal or nonmonotonic
logic. They do play a part in other logics. The attentive reader may
already suspect that we have her an opportunity for generalisation
of say nonmonotonic logic, by giving a role to arc annotations.
\subsubsection*{Argumentation Nets}

Here the nodes of $S$ are interpreted as arguments. The atoms $\{p,q\}$
can be interpreted as types of arguments and the arrows e.g. $(a,b) \xbe
\xdR $
as indicating that the argument a is attacking the argument $b.$

So, for example, let
\begin{quote}

a $=$ we must win votes.

$b$ $=$ death sentence for murderers.

$c$ $=$ We must allow abortion for teenagers

$d$ $=$ Bible forbids taking of life.

$q$ $=$ the argument is a social argument

$p$ $=$ the argument is a religious argument.

$(d,(a,c))$ $=$ there should be no connection between winning votes
and
abortion.

$((a,b),(d,c))$ $=$ If we attack the death sentence in order to win
votes then we must stress (attack) that there should be no connection
between religion (Bible) and social issues.
\end{quote}

Thus we have according to this model that supporting abortion can
lose votes. The argument for abortion is a social one and the argument
from the Bible against it is a religious one.

We can extract information from this IBRS using two algorithms. The
modal logic one can check whether for example every social argument
is attacked by a religious argument. The answer is no, since the social
argument $b$ is attacked only by a which is not a religious argument.

We can also use algorithm $ \xda_{3}$ (following Dung) to extract the
winning arguments of this system. The arguments a and $d$ are winning
since they are not attacked. $d$ attacks the connection between a
and $c$ (i.e. stops a attacking $c).$

The attack of a on $b$ is successful and so $b$ is out. However the
arc $(a,b)$ attacks the arc $(d,c).$ So $c$ is not attacked at all
as both arcs leading into it are successfully eliminated. So $c$
is in. $e$ is out because it is attacked by $d.$

So the winning arguments are $\{a,c,d\}$

In this model we ignore the annotations on arcs. To be consistent
in our mathematics we need to say that $h$ is a partial function on
$ \xdR.$ The best way is to give more specific definition on IBRS to
make it suitable for each logic.

See also  \cite{Gab08b} and  \cite{BGW05}.
\subsubsection*{Counterfactuals}

The traditional semantics for counterfactuals involves closeness of
worlds. The clauses $y \xcm p \xfo q,$ where $ \xfo $
is a counterfactual implication is that $q$ holds in all worlds $y' $
``near enough'' to $y$ in which $p$ holds. So if we interpret the
annotation on arcs as distances then we can define ``near'' as distance
$ \xck $ 2, we get: $a \xcm p \xfo q$ iff in all worlds
of $p-$distance $ \xck 2$ if $p$ holds so does $q.$ Note that the distance
depends on $p.$

In this case we get that $a \xcm p \xfo q$ holds. The
distance function can also use the arrows from arcs to arcs, etc.
There are many opportunities for generalisation in our IBRS set up.
\subsubsection*{Intuitionistic Persistence}

We can get an intuitionistic Kripke model out of this IBRS by letting,
for $t,s \xbe S,t \xbr_{0}s$ iff $t=s$ or $[tRs \xcu \xcA q \xbe Q(h(q,t)
\xck h(q,s))].$
We get that

[
$r_{0}=\{(y,y) \xfA y \xbe S\} \xcv \{(a,b),(a,c),(d,e)\}.$
]

Let $ \xbr $ be the transitive closure of $ \xbr_{0}.$ Algorithm $
\xda_{4}$
evaluates $p \xch q$ in this model, where $ \xch $ is intuitionistic
implication.

$ \xda_{4}:$ $p \xch q$ holds at the IBRS iff $p \xch q$ holds
intuitionistically
at every $ \xbr -$minimal point $of(S, \xbr ).$
\subsection{The power of IBRS}
\label{Section 1.2}

We show now how a number of logics fit into our general picture of IBRS.

 \xEh

 \xDH Nonmonotonic logics in the form of preferential logics:

There are only arrows from nodes to nodes, and they are unlabelled.
The nodes are classical models, and as such all propositional variables
of the base language are given a value from $\{0,1\}.$

The structure is used as described above, i.e. the $R-$minimal models of a
formula or theory are considered.

 \xDH Theory Revision

In the full case, i.e. where the left hand side can change, nodes are
again
classical models, arrows exist only between nodes, and express by their
label
the distance between nodes. Thus, there is just one (dummy) $p,$ and a
real
value as label. In the AGM situation, where the left hand side is fixed,
nodes are classical models (on the right) or sets thereof (on the left),
arrows go from sets of models to models, and express again distance from a
set
(the $K-$models in AGM notation) to a model (of the new formula $ \xbf ).$

The structure is used by considering the closest $ \xbf -$models.

The framework is sufficiently general to express revision also
differently:
Nodes are pairs of classical models, and arrows express that in pair
$(a,b)$
the distance from a to $b$ is smaller than the distance in the pair $(a'
,b' ).$

 \xDH Theory Update

As developments of length 2 can be expressed by a binary relation and the
distance associated, we can - at least in the simple case - proceed
analogously
to the first revision situation. It seems, however, more natural to
consider
as nodes threads of developments, i.e. sequences of classical models, as
arrows
comparisons between such threads, i.e. unlabelled simple arrows only,
expressing
that one thread is more natural or likely than another.

The evaluation is then by considering the ``best'' threads under above
comparison,
and taking a projection on the desired coordinate (i.e. classical model).
The
result is then the theory defined by these projections.

 \xDH Deontic Logic

Just as for preferential logics.

 \xDH The Logic of Counterfactual Conditionals

Again, we can compare pairs (with same left element) as above, or,
alternatively, compare single models with respect to distance from a fixed
other model. This would give arrows with indices, which stand for this
other model.

Evaluation will then be as usual, taking the closest $ \xbf -$models, and
examining
whether $ \xbq $ holds in them.

 \xDH Modal Logic

Nodes are classical models, and thus have the usual labels, arrows are
unlabelled, and only between nodes, and express reachability.

For evaluation, starting from some point, we collect all reachable other
models,
perhaps adding the point of departure.

 \xDH Intuitionistic Logic

Just as for modal logic.

 \xDH Inheritance Systems

Nodes are properties (or sets of models), arrows come in two flavours,
positive
and negative, and exist between nodes only.

The evaluation is relatively complicated, and the subject of ongoing
discussion.

 \xDH Argumentation Theory

There is no unique description of an argumentation system as an IBRS.
For instance, an inheritance system is an argumentation system, so we can
describe such a system as detailed above. But an argument can also be a
deontic statement, as we saw in the first part of this introduction, and
a deontic statement can be described as an IBRS itself. Thus, a node can
be,
under finer granularity, itself an IBRS. Labels can describe the type of
argument (social, etc.) or its validity, etc.

 \xEj
\subsection{
Abstract semantics for IBRS and its engineering realization
}
\label{Section 1.2.2}
\label{Section Reac-Sem}
\index{Section Reac-Sem}

\subsubsection{
Introduction
}

% {\LARGE karl-search= Start Reac-Sem-Intro }

\label{Section Reac-Sem-Intro}
\index{Section Reac-Sem-Intro}

(1) Nodes and arrows

As we may have counterarguments not only against nodes, but also against
arrows, they must be treated basically the same way, i.e. in some way
there has
to be a positive, but also a negative influence on both. So arrows cannot
just
be concatenation between the contents of nodes, or so.

We will differentiate between nodes and arrows by labelling arrows in
addition
with a time delay. We see nodes as situations, where the output is
computed
instantenously from the input, whereas arrows describe some ``force'' or
``mechanism'' which may need some time to ``compute'' the result from the
input.

Consequently, if $ \xba $ is an arrow, and $ \xbb $ an arrow pointing to $
\xba,$ then it
should point to the input of $ \xba,$ i.e. before the time lapse.
Conversely,
any arrow originating in $ \xba $ should originate after the time lapse.

Apart this distinction, we will treat nodes and arrows the same way, so
the
following discussion will apply to both - which we call just ``objects''.

(2) Defeasibility

The general idea is to code each object, say $X,$ by $I(X):U(X) \xcp
C(X):$ If $I(X)$
holds then, unless $U(X)$ holds, consequence $C(X)$ will hold. (We adopted
Reiter's
notation for defaults, as IBRS have common points with the former.)

The situation is slightly more complicated, as there can be several
counterarguments, so $U(X)$ really is an ``or''. Likewise, there can be
several
supporting arguments, so $I(X)$ also is an ``or''.

A counterargument must not always be an argument against a specific
supporting
argument, but it can be. Thus, we should admit both possibilties. As we
can use
arrows to arrows, the second case is easy to treat (as is the dual, a
supporting
argument can be against a specific counterargument). How do we treat the
case
of unspecific pro- and counterarguments? Probably the easiest way is to
adopt
Dung's idea: an object is in, if it has at least one support, and no
counterargument - see  \cite{Dun95}.
Of course, other possibilities may be adopted, counting, use
of labels, etc., but we just consider the simple case here.

(3) Labels

In the general case, objects stand for some kind of defeasible
transmission.
We may in some cases see labels as restricting this transmission to
certain
values. For instance, if the label is $p=1$ and $q=0,$ then the $p-$part
may be
transmitted and the $q-$part not.

Thus, a transmission with a label can sometimes be considered as a family
of
transmissions, which ones are active is indicated by the label.

\be

$\hspace{0.01em}$

% (+++ Orig. No.:  Example 2.1 +++)

\label{Example 2.1}

In fuzzy Kripke models, labels are elements of $[0,1].$ $p=0.5$ as label
for a node
$m' $ which stands for a fuzzy model means that the value of $p$ is 0.5.
$p=0.5$ as
label for an arrow from $m$ to $m' $ means that $p$ is transmitted with
value 0.5.
Thus, when we look from $m$ to $m',$ we see $p$ with value
$0.5*0.5=0.25.$ So, we have
$ \xcx p$ with value 0.25 at $m$ - if, e.g., $m,m' $ are the only models.

\ee

(4) Putting things together

If an arrow leaves an object, the object's output will be connected to the
(only) positive input of the arrow. (An arrow has no negative inputs from
objects it leaves.) If a positive arrow enters an object, it is connected
to
one of the positive inputs of the object, analogously for negative arrows
and
inputs.

When labels are present, they are transmitted through some operation.

% karl-search= End Reac-Sem-Intro
\vspace{7mm}

% *************************************

\vspace{7mm}

%  2.2  Formal definition
%  2.2  Formal definition
% %
% =======================

\subsubsection{
Formal definition
}

% {\LARGE karl-search= Start Reac-Sem-Def }

\label{Section Reac-Sem-Def}
\index{Section Reac-Sem-Def}

\bd

$\hspace{0.01em}$

% (+++ Orig. No.:  Definition 2.1 +++)

\label{Definition 2.1}

In the most general case, objects of IBRS have the form:
$( \xBc I_{1},L_{1} \xBe, \Xl, \xBc I_{n},L_{n} \xBe ):( \xBc U_{1},L'_{1} \xBe
, \Xl
, \xBc U_{n},L'_{n} \xBe ),$ where the $L_{i},L'_{i}$ are labels
and the $I_{i},U_{i}$ might be just truth values, but can also be more
complicated,
a (possibly infinite) sequence of some values. Connected objects have,
of course, to have corresponding such sequences. In addition, the object
$X$ has a criterion for each input, whether it is valid or not (in the
simple
case,
this will just be the truth value $'' true'' ).$ If there is at least one
positive
valid input $I_{i},$ and no valid negative input $U_{i},$ then the output
$C(X)$ and its
label are calculated on the basis of the valid inputs and their labels. If
the
object is an arrow, this will take some time, $t,$ otherwise, this is
instantaneous.

\ed

Evaluating a diagram

An evaluation is relative to a fixed input, i.e. some objects will be
given
certain values, and the diagram is left to calculate the others. It may
well
be that it oscillates, i.e. shows a cyclic behaviour. This may be true for
a
subset of the diagram, or the whole diagram. If it is restricted to an
unimportant part, we might neglect this. Whether it oscillates or not can
also depend on the time delays of the arrows (see Example \ref{Example 2.2}
(page \pageref{Example 2.2}) ).

We therefore define for a diagram $ \xbD $

$ \xba \xcn_{ \xbD } \xbb $ iff

(a) $ \xba $ is a (perhaps partial) input - where the other values are set
``not valid''

(b) $ \xbb $ is a (perhaps partial) output

(c) after some time, $ \xbb $ is stable, i.e. all still possible
oscillations
do not affect $ \xbb $

(d) the other possible input values do not matter, i.e. whatever the
input,
the result is the same.

In the cases examined here more closely, all input values will be defined.

% karl-search= End Reac-Sem-Def
\vspace{7mm}

% *************************************

\vspace{7mm}

%  2.3  A circuit semantics for simple IBRS without labels
%  2.3  A circuit semantics for simple IBRS without labels
% %
% ========================================================

\subsubsection{
A circuit semantics for simple IBRS without labels
}

% {\LARGE karl-search= Start Reac-Sem-Circuit }

\label{Section Reac-Sem-Circuit}
\index{Section Reac-Sem-Circuit}

It is standard to implement the usual logical connectives by electronic
circuits. These components are called gates. Circuits with feedback
sometimes
show undesirable behaviour when the initial conditions are not specified.
(When
we switch a circuit on, the outputs of the individual gates can have
arbitrary
values.) The technical
realization of these initial values shows the way to treat defaults. The
initial
values are set via resistors (in the order of 1 $k \xbO )$ between the
point in the
circuit we want to intialize and the desired tension (say 0 Volt for
false,
5 Volt for true). They are called pull-down or pull-up resistors (for
default
0 or 5 Volt). When a ``real'' result comes in, it will override the tension
applied via the resistor.

Closer inspection reveals that we have here a 3 level default situation:
The initial value will be the weakest, which can be overridden by any
``real''
signal, but a positive argument can be overridden by a negative one. Thus,
the biggest resistor will be for the initialization, the smaller one for
the
supporting arguments, and the negative arguments have full power.

Technical details will be left to the experts.

We give now an example which shows that the delays of the arrows can
matter.
In one situation, a stable state is reached, in another, the circuit
begins to
oscillate.

\be

$\hspace{0.01em}$

% (+++ Orig. No.:  Example 2.2 +++)

\label{Example 2.2}

(In engineering terms, this is a variant of a JK flip-flop with $R*S=0,$ a
circuit
with feedback.)

We have 8 measuring points.

$In1,In2$ are the overall input, $Out1,Out2$ the overall output,
$A1,A2,A3,A4$ are auxiliary internal points. All points can be true or
false.

The logical structure is as follows:

A1 $=$ $In1 \xcu Out1,$ A2 $=$ $In2 \xcu Out2,$

A3 $=$ $A1 \xco Out2,$ A4 $=$ $A2 \xco Out1,$

Out1 $=$ $ \xCN A3,$ Out2 $=$ $ \xCN A4.$

Thus, the circuit is symmetrical, with In1 corresponding to In2, A1 to A2,
A3 to A4, Out1 to Out2.

The input is held constant. See Diagram \ref{Diagram Gate-Sem} (page
\pageref{Diagram Gate-Sem}).

\vspace{10mm}

\begin{diagram}

\label{Diagram Gate-Sem}
\index{Diagram Gate-Sem}

\centering
\setlength{\unitlength}{1mm}
{\renewcommand{\dashlinestretch}{30}
\begin{picture}(150,170)(0,0)

% Upper part

\put(15,130){\arc{10}{-1.57}{1.57}}
\path(15,125)(15,135)
\put(16.3,129.3){\xssc{$\xcu$}}
\put(22,132){\xssc{$A1$}}
\path(13,132)(15,133)(13,134)
\path(13,126)(15,127)(13,128)

\put(45,133){\arc{10}{-1.57}{1.57}}
\path(45,128)(45,138)
\put(46.3,132.3){\xssc{$\xco$}}
\put(52,135){\xssc{$A3$}}
\path(43,135)(45,136)(43,137)
\path(43,129)(45,130)(43,131)

\path(75,128)(75,138)(83,133)(75,128)
\put(76.3,132.3){\xssc{$\xCN$}}
\path(73,132)(75,133)(73,134)

\path(0,127)(15,127)
\path(20,130)(45,130)
\path(50,133)(75,133)
\path(83,133)(108,133)
\path(106,132)(108,133)(106,134)
\put(93,133){\circle*{1}}
\put(101,133){\circle*{1}}

\put(0,124){\xssc{$In1$}}
\put(110,132){\xssc{$Out1$}}

\path(15,133)(5,133)(5,150)(101,150)(101,133)
\path(45,136)(35,136)(35,143)(85,143)(93,77)

% Lower part

\put(15,80){\arc{10}{-1.57}{1.57}}
\path(15,75)(15,85)
\put(16.3,79.3){\xssc{$\xcu$}}
\put(22,76){\xssc{$A2$}}
\path(13,82)(15,83)(13,84)
\path(13,76)(15,77)(13,78)

\put(45,77){\arc{10}{-1.57}{1.57}}
\path(45,72)(45,82)
\put(46.3,76.3){\xssc{$\xco$}}
\put(52,73){\xssc{$A4$}}
\path(43,79)(45,80)(43,81)
\path(43,73)(45,74)(43,75)

\path(75,72)(75,82)(83,77)(75,72)
\put(76.3,76.3){\xssc{$\xCN$}}
\path(73,78)(75,77)(73,76)

\path(0,83)(15,83)
\path(20,80)(45,80)
\path(50,77)(75,77)
\path(83,77)(108,77)
\path(106,76)(108,77)(106,78)
\put(93,77){\circle*{1}}
\put(101,77){\circle*{1}}

\put(0,85.5){\xssc{$In2$}}
\put(110,76){\xssc{$Out2$}}

\path(15,77)(5,77)(5,60)(101,60)(101,77)
\path(45,74)(35,74)(35,67)(85,67)(93,133)

\put(30,20) {{\rm\bf Gate Semantics}}

\end{picture}
}

\end{diagram}

\vspace{4mm}

\ee

We suppose that the output of the individual gates is present $n$ time
slices
after the input was present. $n$ will in the first circuit be equal to 1
for all
gates, in the second circuit equal to 1 for all but the AND gates, which
will
take 2 time slices. Thus, in both cases, e.g. Out1 at time $t$ will be the
negation of A3 at time $t-1.$ In the first case, A1 at time $t$ will be
the
conjunction of In1 and Out1 at time $t-1,$ and in the second case the
conjunction
of In1 and Out1 at time $t-2.$

We initialize In1 as true, all others as false. (The initial value of A3
and A4
does not matter, the behaviour is essentially the same for all such
values.)

The first circuit will oscillate with a period of 4, the second circuit
will go
to a stable state.

We have the following transition tables (time slice shown at left):

Circuit 1, $delay=1$ everywhere:
%     In1 In2 A1 A2 A3 A4 Out1 Out2 &;
%     In1 In2 A1 A2 A3 A4 Out1 Out2 &;
% 1:  T   F   F  F  F  F  F    F, &;
% 1:  T   F   F  F  F  F  F    F, &;
% 2:  T   F   F  F  F  F  T    T, &;
% 2:  T   F   F  F  F  F  T    T, &;
% 3:  T   F   T  F  T  T  T    T, &;
% 3:  T   F   T  F  T  T  T    T, &;
% 4:  T   F   T  F  T  T  F    F, &;
% 4:  T   F   T  F  T  T  F    F, &;
% 5:  T   F   F  F  T  F  F    F     - oscillation starts &;
% 5:  T   F   F  F  T  F  F    F     - oscillation starts &;
% 6:  T   F   F  F  F  F  F    T, &;
% 6:  T   F   F  F  F  F  F    T, &;
% 7:  T   F   F  F  T  F  T    T, &;
% 7:  T   F   F  F  T  F  T    T, &;
% 8:  T   F   T  F  T  T  F    T, &;
% 8:  T   F   T  F  T  T  F    T, &;
% 9:  T   F   F  F  T  F  F    F     - back to start of oscillation.
% 9:  T   F   F  F  T  F  F    F     - back to start of oscillation.

\begin{tabular}{lccccccccl}
   & In1 & In2 & A1 & A2 & A3 & A4 & Out1 & Out2 & \\
1: &   T &   F &  F &  F &  F &  F &    F &    F & \\
2: &   T &   F &  F &  F &  F &  F &    T &    T & \\
3: &   T &   F &  T &  F &  T &  T &    T &    T & \\
4: &   T &   F &  T &  F &  T &  T &    F &    F & \\
5: &   T &   F &  F &  F &  T &  F &    F &    F & oscillation starts \\
6: &   T &   F &  F &  F &  F &  F &    F &    T & \\
7: &   T &   F &  F &  F &  T &  F &    T &    T & \\
8: &   T &   F &  T &  F &  T &  T &    F &    T & \\
9: &   T &   F &  F &  F &  T &  F &    F &    F &
back to start of oscillation \\
\end{tabular}

Circuit 2, $delay=1$ everywhere, except for AND with $delay=2:$

(Thus, A1 and A2 are held at their intial value up to time 2, then they
are
calculated using the values of time $t-2.)$
%     In1 In2 A1 A2 A3 A4 Out1 Out2 &;
%     In1 In2 A1 A2 A3 A4 Out1 Out2 &;
% 1:  T   F   F  F  F  F  F    F, &;
% 1:  T   F   F  F  F  F  F    F, &;
% 2:  T   F   F  F  F  F  T    T, &;
% 2:  T   F   F  F  F  F  T    T, &;
% 3:  T   F   F  F  T  T  T    T, &;
% 3:  T   F   F  F  T  T  T    T, &;
% 4:  T   F   T  F  T  T  F    F, &;
% 4:  T   F   T  F  T  T  F    F, &;
% 5:  T   F   T  F  T  F  F    F, &;
% 5:  T   F   T  F  T  F  F    F, &;
% 6:  T   F   F  F  T  F  F    T     - stable state is reached &;
% 6:  T   F   F  F  T  F  F    T     - stable state is reached &;
% 7:  T   F   F  F  T  F  F    T.
% 7:  T   F   F  F  T  F  F    T.

\begin{tabular}{lccccccccl}
   & In1 & In2 & A1 & A2 & A3 & A4 & Out1 & Out2 & \\
1: &   T &   F &  F &  F &  F &  F &    F &    F & \\
2: &   T &   F &  F &  F &  F &  F &    T &    T & \\
3: &   T &   F &  F &  F &  T &  T &    T &    T & \\
4: &   T &   F &  T &  F &  T &  T &    F &    F & \\
5: &   T &   F &  T &  F &  T &  F &    F &    F & \\
6: &   T &   F &  F &  F &  T &  F &    F &    T & stable state reached \\
7: &   T &   F &  F &  F &  T &  F &    F &    T & \\
8: &   T &   F &  F &  F &  T &  F &    F &    T & \\
\end{tabular}

Note that state 6 of circuit 2 is also stable in circuit 1, but it is
never
reached in that circuit.

% karl-search= End Reac-Sem-Circuit
\vspace{7mm}

% *************************************

\vspace{7mm}

\chapter{
Basic definitions and results
}
\label{Section Basic-definitions-and-results}

\section{
Algebraic definitions
}
\index{Notation FOL-Tilde}

\bn

$\hspace{0.01em}$

% (+++ Orig. No.:  Notation D-FOL-Tilde +++)

\label{Notation D-FOL-Tilde}

\index{FOL}
\index{NML}
We use sometimes FOL as abbreviation for first order logic, and NML for
nonmonotonic logic.
To avoid Latex complications in bigger expressions, we replace
$\widetilde{xxxxx}$ by $\wt{xxxxx}$.
\index{Definition Alg-Base}

\en

\bd

$\hspace{0.01em}$

% (+++ Orig. No.:  Definition Alg-Base +++)

\label{Definition Alg-Base}

We use $ \xdp $ to denote the power set operator,
$ \xbP \{X_{i}:i \xbe I\}$ $:=$ $\{g:$ $g:I \xcp \xcV \{X_{i}:i \xbe I\},$
$ \xcA i \xbe I.g(i) \xbe X_{i}\}$ is the general cartesian
product, $card(X)$ shall denote the cardinality of $X,$ and $V$ the
set-theoretic
universe we work in - the class of all sets. Given a set of pairs $ \xdx
,$ and a
set $X,$ we denote by $ \xdx \xex X:=\{ \xBc x,i \xBe  \xbe \xdx:x \xbe X\}.$
When
the context is clear, we
will sometime simply write $X$ for $ \xdx \xex X.$ (The intended use is
for preferential
structures, where $x$ will be a point (intention: a classical
propositional
model), and $i$ an index, permitting copies of logically identical
points.)

$A \xcc B$ will denote that $ \xCf A$ is a subset of $B$ or equal to $B,$
and $A \xcb B$ that $ \xCf A$ is
a proper subset of $B,$ likewise for $A \xcd B$ and $A \xcf B.$

Given some fixed set $U$ we work in, and $X \xcc U,$ then $ \xdC (X):=U-X$
.

If $ \xdy \xcc \xdp (X)$ for some
$X,$ we say that $ \xdy $ satisfies

$( \xcs )$ iff it is closed under finite intersections,

$( \xcS )$ iff it is closed under arbitrary intersections,

$( \xcv )$ iff it is closed under finite unions,

$( \xcV )$ iff it is closed under arbitrary unions,

$( \xdC )$ iff it is closed under complementation,

$ \xCf (-)$ iff it is closed under set difference.

We will sometimes write $A=B \xFO C$ for: $A=B,$ or $A=C,$ or $A=B \xcv
C.$

We make ample and tacit use of the Axiom of Choice.
\index{Definition Rel-Base}

\ed

\bd

$\hspace{0.01em}$

% (+++ Orig. No.:  Definition Rel-Base +++)

\label{Definition Rel-Base}

$ \xeb^{*}$ will denote the transitive closure of the relation $ \xeb.$
If a relation $<,$
$ \xeb,$ or similar is given, $a \xcT b$ will express that a and $b$ are
$<-$ (or $ \xeb -)$
incomparable - context will tell. Given any relation $<,$ $ \xck $ will
stand for
$<$ or $=,$ conversely, given $ \xck,$ $<$ will stand for $ \xck,$ but
not $=,$ similarly
for $ \xeb $ etc.
\index{Definition Tree-Base}

\ed

\bd

$\hspace{0.01em}$

% (+++ Orig. No.:  Definition Tree-Base +++)

\label{Definition Tree-Base}

A child (or successor) of an element $x$ in a tree $t$ will be a direct
child in $t.$
A child of a child, etc. will be called an indirect child. Trees will be
supposed to grow downwards, so the root is the top element.
\index{Definition Seq-Base}

\ed

\bd

$\hspace{0.01em}$

% (+++ Orig. No.:  Definition Seq-Base +++)

\label{Definition Seq-Base}

A subsequence $ \xbs_{i}:i \xbe I \xcc \xbm $ of a sequence $ \xbs_{i}:i
\xbe \xbm $ is called cofinal, iff
for all $i \xbe \xbm $ there is $i' \xbe I$ $i \xck i'.$

Given two sequences $ \xbs_{i}$ and $ \xbt_{i}$ of the same length, then
their Hamming distance
is the quantity of $i$ where they differ.
\index{Definition Def-Clos}

\ed

\bd

$\hspace{0.01em}$

% (+++ Orig. No.:  Definition Def-Clos +++)

\label{Definition Def-Clos}

Let $ \xdy \xcc \xdp (Z)$ be given and closed under arbitrary
intersections.

(1) For $A \xcc Z,$ let $ \wt{A}$ $:=$ $ \xcS \{X \xbe \xdy:A \xcc X\}.$

(2) For $B \xbe \xdy,$ we call $A \xcc B$ a small subset of $B$ iff there
is no
$X \xbe \xdy $ such that $B-A \xcc X \xcb B.$

(Context will disambiguate from other uses of
$'' small''.)$

Intuitively, $Z$ is the set of all models for $ \xdl,$ $ \xdy $ is $
\xdD_{ \xdl }$, and $ \wt{A}=M(Th(A)),$
this is the intended application - $Th(A)$ is the set of formulas which
hold in all $a \xbe A,$ and $M(Th(A))$ is the set of models of $Th(A).$
Note that then $ \wt{ \xCQ }= \xCQ.$
\index{Fact Def-Clos}

\ed

\bfa

$\hspace{0.01em}$

% (+++ Orig. No.:  Fact Def-Clos +++)

\label{Fact Def-Clos}

(1) If $ \xdy \xcc \xdp (Z)$ is closed under arbitrary intersections and
finite unions,
$Z \xbe \xdy,$ $X,Y \xcc Z,$ then the following hold:

$(Cl \xcv )$ $ \wt{X \xcv Y}$ $=$ $ \wt{X} \xcv \wt{Y}$

$(Cl \xcs )$ $ \wt{X \xcs Y} \xcc \wt{X} \xcs \wt{Y},$ but usually not
conversely,

$ \xCf (Cl-)$ $ \wt{A}- \wt{B} \xcc \wt{A-B},$

$(Cl=)$ $X=Y$ $ \xch $ $ \wt{X}= \wt{Y},$ but not conversely,

$(Cl \xcc 1)$ $ \wt{X} \xcc Y$ $ \xch $ $X \xcc Y,$ but not conversely,

$(Cl \xcc 2)$ $X \xcc \wt{Y}$ $ \xch $ $ \wt{X} \xcc \wt{Y}.$

(2) If, in addition, $X \xbe \xdy $ and $ \xdC X:=Z-X \xbe \xdy,$ then
the following two properties
hold, too:

$(Cl \xcs +)$ $ \wt{A} \xcs X= \wt{A \xcs X},$

$(Cl-+)$ $ \wt{A}-X= \wt{A-X}.$

(3) In the intended application, i.e. $ \wt{A}=M(Th(A)),$ the following
hold:

(3.1) $Th(X)$ $=$ $Th( \wt{X}),$

(3.2) Even if $A= \wt{A},$ $B= \wt{B},$ it is not necessarily true that $
\wt{A-B} \xcc \wt{A}- \wt{B}.$
\index{Fact Def-Clos Proof}

\efa

\subparagraph{
Proof
}

$\hspace{0.01em}$

% (+++ Orig.:  Proof +++)

$(Cl=),$ $(Cl \xcc 1),$ $(Cl \xcc 2),$ (3.1) are trivial.

$(Cl \xcv )$ Let $ \xdy (U):=\{X \xbe \xdy:U \xcc X\}.$ If $A \xbe \xdy
(X \xcv Y),$ then $A \xbe \xdy (X)$ and $A \xbe \xdy (Y),$ so
$ \wt{X \xcv Y}$ $ \xcd $ $ \wt{X} \xcv \wt{Y}.$
If $A \xbe \xdy (X)$ and $B \xbe \xdy (Y),$ then $A \xcv B \xbe \xdy (X
\xcv Y),$ so $ \wt{X \xcv Y}$ $ \xcc $ $ \wt{X} \xcv \wt{Y}.$

$(Cl \xcs )$ Let $X',Y' \xbe \xdy,$ $X \xcc X',$ $Y \xcc Y',$ then $X
\xcs Y \xcc X' \xcs Y',$ so $ \wt{X \xcs Y} \xcc \wt{X} \xcs \wt{Y}.$
For the converse, set $X:=M_{ \xdl }-\{m\},$ $Y:=\{m\}$ in Example \ref{Example
Not-Def} (page \pageref{Example Not-Def}).
$(M_{ \xdl }$ is the set of all models of the language $ \xdl.)$

$ \xCf (Cl-)$ Let $A-B \xcc X \xbe \xdy,$ $B \xcc Y \xbe \xdy,$ so $A
\xcc X \xcv Y \xbe \xdy.$ Let $x \xce \wt{B}$ $ \xch $ $ \xcE Y \xbe \xdy
(B \xcc Y,$ $x \xce Y),$
$x \xce \wt{A-B}$ $ \xch $ $ \xcE X \xbe \xdy (A-B \xcc X,$ $x \xce X),$
so $x \xce X \xcv Y,$ $A \xcc X \xcv Y,$ so $x \xce \wt{A}.$ Thus, $x \xce
\wt{B},$ $x \xce \wt{A-B}$ $ \xch $
$x \xce \wt{A},$ or $x \xbe \wt{A}- \wt{B}$ $ \xch $ $x \xbe \wt{A-B}.$

$(Cl \xcs +)$ $ \wt{A} \xcs X \xcd \wt{A \xcs X}$ by $(Cl \xcs ).$
For `` $ \xcc $ '': Let $A \xcs X \xcc A' \xbe \xdy,$ then by closure
under $( \xcv ),$
$A \xcc A' \xcv \xdC X \xbe \xdy,$ $(A' \xcv \xdC X) \xcs X \xcc A'.$ So
$ \wt{A} \xcs X \xcc \wt{A \xcs X}.$

$(Cl-+)$ $ \wt{A-X}= \wt{A \xcs \xdC X}= \wt{A} \xcs \xdC X= \wt{A}-X$ by
$(Cl \xcs +).$

(3.2) Set $A:=M_{ \xdl },$ $B:=\{m\}$ for $m \xbe M_{ \xdl }$ arbitrary, $
\xdl $ infinite. So $A= \wt{A},$ $B= \wt{B},$ but
$ \wt{A-B}=A \xEd A-B.$

$ \xcz $
\\[3ex]
\section{
Basic logical definitions
}
\index{Definition Log-Base}

\bd

$\hspace{0.01em}$

% (+++ Orig. No.:  Definition Log-Base +++)

\label{Definition Log-Base}

We work here in a classical propositional language $ \xdl,$ a theory $T$
will be an
arbitrary set of formulas. Formulas will often be named $ \xbf,$ $ \xbq
,$ etc., theories
$T,$ $S,$ etc.

$v( \xdl )$ will be the set of propositional variables of $ \xdl.$

$M_{ \xdl }$ will be the set of (classical) models for $ \xdl,$ $M(T)$ or
$M_{T}$
is the set of models of $T,$ likewise $M( \xbf )$ for a formula $ \xbf.$

$ \xdD_{ \xdl }:=\{M(T):$ $T$ a theory in $ \xdl \},$ the set of definable
model sets.

Note that, in classical propositional logic, $ \xCQ,M_{ \xdl } \xbe
\xdD_{ \xdl },$ $ \xdD_{ \xdl }$ contains
singletons, is closed under arbitrary intersections and finite unions.

An operation $f: \xdy \xcp \xdp (M_{ \xdl })$ for $ \xdy \xcc \xdp (M_{
\xdl })$ is called definability
preserving, $ \xCf (dp)$ or $( \xbm dp)$ in short, iff for all $X \xbe
\xdD_{ \xdl } \xcs \xdy $ $f(X) \xbe \xdD_{ \xdl }.$

We will also use $( \xbm dp)$ for binary functions $f: \xdy \xCK \xdy \xcp
\xdp (M_{ \xdl })$ - as needed
for theory revision - with the obvious meaning.

$ \xcl $ will be classical derivability, and

$ \ol{T}:=\{ \xbf:T \xcl \xbf \},$ the closure of $T$ under $ \xcl.$

$Con(.)$ will stand for classical consistency, so $Con( \xbf )$ will mean
that
$ \xbf $ is clasical consistent, likewise for $Con(T).$ $Con(T,T' )$ will
stand for
$Con(T \xcv T' ),$ etc.

Given a consequence relation $ \xcn,$ we define

$ \ol{ \ol{T} }:=\{ \xbf:T \xcn \xbf \}.$

(There is no fear of confusion with $ \ol{T},$ as it just is not useful to
close
twice under classical logic.)

$T \xco T':=\{ \xbf \xco \xbf ': \xbf \xbe T, \xbf ' \xbe T' \}.$

If $X \xcc M_{ \xdl },$ then $Th(X):=\{ \xbf:X \xcm \xbf \},$ likewise
for $Th(m),$ $m \xbe M_{ \xdl }.$ $( \xcm $ will
usually be classical validity.)
\index{Fact Log-Base}

\ed

We recollect and note:

\bfa

$\hspace{0.01em}$

% (+++ Orig. No.:  Fact Log-Base +++)

\label{Fact Log-Base}

Let $ \xdl $ be a fixed propositional language, $ \xdD_{ \xdl } \xcc X,$ $
\xbm:X \xcp \xdp (M_{ \xdl }),$ for a $ \xdl -$theory $T$
$ \ol{ \ol{T} }:=Th( \xbm (M_{T})),$ let $T,$ $T' $ be arbitrary theories,
then:

(1) $ \xbm (M_{T}) \xcc M_{ \ol{ \ol{T} }}$,

(2) $M_{T} \xcv M_{T' }=M_{T \xco T' }$ and $M_{T \xcv T' }=M_{T} \xcs
M_{T' }$,

(3) $ \xbm (M_{T})= \xCQ $ $ \xcj $ $ \xcT \xbe \ol{ \ol{T} }$.

If $ \xbm $ is definability preserving or $ \xbm (M_{T})$ is finite, then
the following also hold:

(4) $ \xbm (M_{T})=M_{ \ol{ \ol{T} }}$,

(5) $T' \xcl \ol{ \ol{T} }$ $ \xcj $ $M_{T' } \xcc \xbm (M_{T}),$

(6) $ \xbm (M_{T})=M_{T' }$ $ \xcj $ $ \ol{T' }= \ol{ \ol{T} }.$
$ \xcz $
\\[3ex]
\index{Fact Th-Union}

\efa

\bfa

$\hspace{0.01em}$

% (+++ Orig. No.:  Fact Th-Union +++)

\label{Fact Th-Union}

Let $A,B \xcc M_{ \xdl }.$

Then $Th(A \xcv B)$ $=$ $Th(A) \xcs Th(B).$
\index{Fact Th-Union Proof}

\efa

\subparagraph{
Proof
}

$\hspace{0.01em}$

% (+++ Orig.:  Proof +++)

$ \xbf \xbe Th(A \xcv B)$ $ \xcj $ $A \xcv B \xcm \xbf $ $ \xcj $ $A \xcm
\xbf $ and $B \xcm \xbf $ $ \xcj $ $ \xbf \xbe Th(A)$ and $ \xbf \xbe
Th(B).$

$ \xcz $
\\[3ex]
\index{Fact Log-Form}

\bfa

$\hspace{0.01em}$

% (+++ Orig. No.:  Fact Log-Form +++)

\label{Fact Log-Form}

Let $X \xcc M_{ \xdl },$ $ \xbf, \xbq $ formulas.

(1) $X \xcs M( \xbf ) \xcm \xbq $ iff $X \xcm \xbf \xcp \xbq.$

(2) $X \xcs M( \xbf ) \xcm \xbq $ iff $M(Th(X)) \xcs M( \xbf ) \xcm \xbq
.$

(3) $Th(X \xcs M( \xbf ))= \ol{Th(X) \xcv \{ \xbf \}}$

(4) $X \xcs M( \xbf )= \xCQ $ $ \xcj $ $M(Th(X)) \xcs M( \xbf )= \xCQ $

(5) $Th(M(T) \xcs M(T' ))= \ol{T \xcv T' }.$
\index{Fact Log-Form Proof}

\efa

\subparagraph{
Proof
}

$\hspace{0.01em}$

% (+++ Orig.:  Proof +++)

(1) `` $ \xch $ '': $X=(X \xcs M( \xbf )) \xcv (X \xcs M( \xCN \xbf )).$ In
both parts holds $ \xCN \xbf \xco \xbq,$ so
$X \xcm \xbf \xcp \xbq.$ `` $ \xci $ '': Trivial.

(2) $X \xcs M( \xbf ) \xcm \xbq $ (by (1)) iff $X \xcm \xbf \xcp \xbq $
iff $M(Th(X)) \xcm \xbf \xcp \xbq $ iff (again by (1))
$M(Th(X)) \xcs M( \xbf ) \xcm \xbq.$

(3) $ \xbq \xbe Th(X \xcs M( \xbf ))$ $ \xcj $ $X \xcs M( \xbf ) \xcm \xbq
$ $ \xcj_{(2)}$ $M(Th(X) \xcv \{ \xbf \})=M(Th(X)) \xcs M( \xbf ) \xcm
\xbq $ $ \xcj $
$Th(X) \xcv \{ \xbf \} \xcl \xbq.$

(4) $X \xcs M( \xbf )= \xCQ $ $ \xcj $ $X \xcm \xCN \xbf $ $ \xcj $
$M(Th(X)) \xcm \xCN \xbf $ $ \xcj $ $M(Th(X)) \xcs M( \xbf )= \xCQ.$

(5) $M(T) \xcs M(T' )=M(T \xcv T' ).$

$ \xcz $
\\[3ex]
\index{Fact Dp-Base}

\bfa

$\hspace{0.01em}$

% (+++ Orig. No.:  Fact Dp-Base +++)

\label{Fact Dp-Base}

If $X=M(T),$ then $M(Th(X))=X.$
\index{Fact Dp-Base Proof}

\efa

\subparagraph{
Proof
}

$\hspace{0.01em}$

% (+++ Orig.:  Proof +++)

$X \xcc M(Th(X))$ is trivial. $Th(M(T))= \ol{T}$ is trivial by classical
soundness and
completeness. So $M(Th(M(T))=M( \ol{T})=M(T)=X.$ $ \xcz $
\\[3ex]
\index{Example Not-Def}

\be

$\hspace{0.01em}$

% (+++ Orig. No.:  Example Not-Def +++)

\label{Example Not-Def}

If $v( \xdl )$ is infinite, and $m$ any model for $ \xdl,$ then $M:=M_{
\xdl }-\{m\}$ is not definable
by any theory $T.$ (Proof: Suppose it were, and let $ \xbf $ hold in $M,$
but not in $m,$ so in $m$ $ \xCN \xbf $ holds, but as $ \xbf $ is finite,
there is a model $m' $ in
$M$ which coincides on all propositional variables of $ \xbf $ with $m,$
so in $m' $ $ \xCN \xbf $
holds, too, a contradiction.) Thus, in the infinite case, $ \xdp (M_{ \xdl
}) \xEd \xdD_{ \xdl }.$

(There is also a simple cardinality argument, which shows that almost no
model sets are definable, but it is not constructive and thus less
instructive
than above argument. We give it nonetheless: Let $ \xbk:=card(v( \xdl
)).$ Then
there are $ \xbk $ many formulas, so $2^{ \xbk }$ many theories, and thus
$2^{ \xbk }$ many
definable model sets. But there are $2^{ \xbk }$ many models, so $(2^{
\xbk })^{ \xbk }$ many model
sets.)

$ \xcz $
\\[3ex]
\section{
Basic definitions and results for nonmonotonic logics
}
\index{Definition Log-Cond-Ref-Size}
\label{Definition Log-Cond-Ref-Size}

\ee

The numbers in the first column ``Correspondence''
refer to Proposition \ref{Proposition Alg-Log} (page \pageref{Proposition
Alg-Log}),
published as Proposition 21 in  \cite{GS08c},
those in the second column ``Correspondence''
to Proposition \ref{Proposition Ref-Class-Mu-neu} (page \pageref{Proposition
Ref-Class-Mu-neu}).

\begin{turn}{90}

% {\tiny

{\xssc

\begin{tabular}{|c|c|c|c|c|c|}

\hline

\multicolumn{2}{|c|}{Logical rule}
\xEH
Correspondence
\xEH
Model set
\xEH
Correspondence
\xEH
Size Rules
\xEP

\hline

\multicolumn{6}{|c|}{Basics}
\xEP

\hline

$(SC)$ Supraclassicality
\xEH
$(SC)$
\xEH
$\xch$ (4.1)
\xEH
$( \xbm \xcc )$
\xEH
trivial
\xEH
$(Opt)$
\xEP

\cline{3-3}

$ \xbf \xcl \xbq $ $ \xch $ $ \xbf \xcn \xbq $
\xEH
$ \ol{T} \xcc \ol{ \ol{T} }$
\xEH
$\xci$ (4.2)
\xEH
$f(X) \xcc X$
\xEH
\xEH
\xEP

\cline{1-1}

$(REF)$ Reflexivity
\xEH
\xEH
\xEH
\xEH
\xEH
\xEP

$ T \xcv \{\xba\} \xcn \xba $
\xEH
\xEH
\xEH
\xEH
\xEH
\xEP

\hline

$(LLE)$
\xEH
$(LLE)$
\xEH
\xEH
\xEH
\xEH
\xEP

Left Logical Equivalence
\xEH
\xEH
\xEH
\xEH
\xEH
\xEP

$ \xcl \xbf \xcr \xbf ',  \xbf \xcn \xbq   \xch $
\xEH
$ \ol{T}= \ol{T' }  \xch   \ol{\ol{T}} = \ol{\ol{T'}}$
\xEH
\xEH
\xEH
\xEH
\xEP

$ \xbf ' \xcn \xbq $
\xEH
\xEH
\xEH
\xEH
\xEH
\xEP

\hline

$(RW)$ Right Weakening
\xEH
$(RW)$
\xEH
\xEH
\xEH
trivial
\xEH
$(iM)$
\xEP

$ \xbf \xcn \xbq,  \xcl \xbq \xcp \xbq '   \xch $
\xEH
$ T \xcn \xbq,  \xcl \xbq \xcp \xbq '   \xch $
\xEH
\xEH
\xEH
\xEH
\xEP

$ \xbf \xcn \xbq ' $
\xEH
$T \xcn \xbq ' $
\xEH
\xEH
\xEH
\xEH
\xEP

\hline

$(wOR)$
\xEH
$(wOR)$
\xEH
$\xch$ (3.1)
\xEH
$( \xbm wOR)$
\xEH
$\xci$ (1.1)
\xEH
$(eM\xdi)$
\xEP

\cline{3-3}
\cline{5-5}

$ \xbf \xcn \xbq,$ $ \xbf ' \xcl \xbq $ $ \xch $
\xEH
$ \ol{ \ol{T} } \xcs \ol{T' }$ $ \xcc $ $ \ol{ \ol{T \xco T' } }$
\xEH
$\xci$ (3.2)
\xEH
$f(X \xcv Y) \xcc f(X) \xcv Y$
\xEH
$\xch$ (1.2)
\xEH
\xEP

$ \xbf \xco \xbf ' \xcn \xbq $
\xEH
\xEH
\xEH
\xEH
\xEH
\xEP

\hline

$(disjOR)$
\xEH
$(disjOR)$
\xEH
$\xch$ (2.1)
\xEH
$( \xbm disjOR)$
\xEH
$\xci$ (4.1)
\xEH
$(I\xcv disj)$
\xEP

\cline{3-3}
\cline{5-5}

$ \xbf \xcl \xCN \xbf ',$ $ \xbf \xcn \xbq,$
\xEH
$\xCN Con(T \xcv T') \xch$
\xEH
$\xci$ (2.2)
\xEH
$X \xcs Y= \xCQ $ $ \xch $
\xEH
$\xch$ (4.2)
\xEH
\xEP

$ \xbf ' \xcn \xbq $ $ \xch $ $ \xbf \xco \xbf ' \xcn \xbq $
\xEH
$ \ol{ \ol{T} } \xcs \ol{ \ol{T' } } \xcc \ol{ \ol{T \xco T' } }$
\xEH
\xEH
$f(X \xcv Y) \xcc f(X) \xcv f(Y)$
\xEH
\xEH
\xEP

\hline

$(CP)$
\xEH
$(CP)$
\xEH
$\xch$ (5.1)
\xEH
$( \xbm \xCQ )$
\xEH
trivial
\xEH
$(I_1)$
\xEP

\cline{3-3}

Consistency Preservation
\xEH
\xEH
$\xci$ (5.2)
\xEH
\xEH
\xEH
\xEP

$ \xbf \xcn \xcT $ $ \xch $ $ \xbf \xcl \xcT $
\xEH
$T \xcn \xcT $ $ \xch $ $T \xcl \xcT $
\xEH
\xEH
$f(X)= \xCQ $ $ \xch $ $X= \xCQ $
\xEH
\xEH
\xEP

\hline

\xEH
\xEH
\xEH
$( \xbm \xCQ fin)$
\xEH
\xEH
$(I_1)$
\xEP

\xEH
\xEH
\xEH
$X \xEd \xCQ $ $ \xch $ $f(X) \xEd \xCQ $
\xEH
\xEH
\xEP

\xEH
\xEH
\xEH
for finite $X$
\xEH
\xEH
\xEP

\hline

\xEH
$(AND_1)$
\xEH
\xEH
\xEH
\xEH
$(I_2)$
\xEP

\xEH
$\xba\xcn\xbb \xch \xba\xcN\xCN\xbb$
\xEH
\xEH
\xEH
\xEH
\xEP

\hline

\xEH
$(AND_n)$
\xEH
\xEH
\xEH
\xEH
$(I_n)$
\xEP

\xEH
$\xba\xcn\xbb_1, \ldots, \xba\xcn\xbb_{n-1} \xch $
\xEH
\xEH
\xEH
\xEH
\xEP

\xEH
$\xba\xcN(\xCN\xbb_1 \xco \ldots \xco \xCN\xbb_{n-1})$
\xEH
\xEH
\xEH
\xEH
\xEP

\hline

$(AND)$
\xEH
$(AND)$
\xEH
\xEH
\xEH
trivial
\xEH
$(I_\xbo)$
\xEP

$ \xbf \xcn \xbq,  \xbf \xcn \xbq '   \xch $
\xEH
$ T \xcn \xbq, T \xcn \xbq '   \xch $
\xEH
\xEH
\xEH
\xEH
\xEP

$ \xbf \xcn \xbq \xcu \xbq ' $
\xEH
$ T \xcn \xbq \xcu \xbq ' $
\xEH
\xEH
\xEH
\xEH
\xEP

\hline

$(CCL)$ Classical Closure
\xEH
$(CCL)$
\xEH
\xEH
\xEH
trivial
\xEH
$(iM)+(I_\xbo)$
\xEP

\xEH
$ \ol{ \ol{T} }$ classically closed
\xEH
\xEH
\xEH
\xEH
\xEP

\hline

$(OR)$
\xEH
$(OR)$
\xEH
$\xch$ (1.1)
\xEH
$( \xbm OR)$
\xEH
$\xci$ (2.1)
\xEH
$(eM\xdi)+(I_\xbo)$
\xEP

\cline{3-3}
\cline{5-5}

$ \xbf \xcn \xbq,  \xbf ' \xcn \xbq   \xch $
\xEH
$ \ol{\ol{T}} \xcs \ol{\ol{T'}} \xcc \ol{\ol{T \xco T'}} $
\xEH
$\xci$ (1.2)
\xEH
$f(X \xcv Y) \xcc f(X) \xcv f(Y)$
\xEH
$\xch$ (2.2)
\xEH
\xEP

$ \xbf \xco \xbf ' \xcn \xbq $
\xEH
\xEH
\xEH
\xEH
\xEH
\xEP

\hline

\xEH
$(PR)$
\xEH
$\xch$ (6.1)
\xEH
$( \xbm PR)$
\xEH
$\xci$ (3.1)
\xEH
$(eM\xdi)+(I_\xbo)$
\xEP

\cline{3-3}
\cline{5-5}

$ \ol{ \ol{ \xbf \xcu \xbf ' } }$ $ \xcc $ $ \ol{ \ol{ \ol{ \xbf } } \xcv
\{ \xbf ' \}}$
\xEH
$ \ol{ \ol{T \xcv T' } }$ $ \xcc $ $ \ol{ \ol{ \ol{T} } \xcv T' }$
\xEH
$\xci (\xbm dp)+(\xbm\xcc)$ (6.2)
\xEH
$X \xcc Y$ $ \xch $
\xEH
$\xch$ (3.2)
\xEH
\xEP

\cline{3-3}

\xEH
\xEH
$\xcI$ without $(\xbm dp)$ (6.3)
\xEH
$f(Y) \xcs X \xcc f(X)$
\xEH
\xEH
\xEP

\cline{3-3}

\xEH
\xEH
$\xci (\xbm\xcc)$ (6.4)
\xEH
\xEH
\xEH
\xEP

\xEH
\xEH
$T'$ a formula
\xEH
\xEH
\xEH
\xEP

\cline{3-4}

\xEH
\xEH
$\xci$ (6.5)
\xEH
$(\xbm PR ')$
\xEH
\xEH
\xEP

\xEH
\xEH
$T'$ a formula
\xEH
$f(X) \xcs Y \xcc f(X \xcs Y)$
\xEH
\xEH
\xEP

\hline

$(CUT)$
\xEH
$(CUT)$
\xEH
$\xch$ (7.1)
\xEH
$ (\xbm CUT) $
\xEH
$\xci$ (8.1)
\xEH
$(eM\xdi)+(I_\xbo)$
\xEP

\cline{3-3}
\cline{5-5}

$ T  \xcn \xba; T \xcv \{ \xba\} \xcn \xbb \xch $
\xEH
$T \xcc \ol{T' } \xcc \ol{ \ol{T} }  \xch $
\xEH
$\xci$ (7.2)
\xEH
$f(X) \xcc Y \xcc X  \xch $
\xEH
$\xcH$ (8.2)
\xEH
\xEP

$ T  \xcn \xbb $
\xEH
$ \ol{ \ol{T'} } \xcc \ol{ \ol{T} }$
\xEH
\xEH
$f(X) \xcc f(Y)$
\xEH
\xEH
\xEP

\hline

% \end{tabular*}
\end{tabular}

}

\end{turn}

\begin{turn}{90}

{\xssc

\begin{tabular}{|c|c|c|c|c|c|}

\hline

\multicolumn{2}{|c|}{Logical rule}
\xEH
Correspondence
\xEH
Model set
\xEH
Correspondence
\xEH
Size-Rule
\xEP

\hline

\multicolumn{6}{|c|}{Cumulativity}
\xEP

\hline

$(wCM)$
\xEH
\xEH
\xEH
\xEH
trivial
\xEH
$(eM\xdf)$
\xEP

$\xba\xcn\xbb, \xba'\xcl\xba, \xba\xcu\xbb\xcl\xba' \xch \xba'\xcn\xbb$
\xEH
\xEH
\xEH
\xEH
\xEH
\xEP

\hline

$(CM_2)$
\xEH
\xEH
\xEH
\xEH
\xEH
$(I_2)$
\xEP

$\xba\xcn\xbb, \xba\xcn\xbb' \xch \xba\xcu\xbb\xcL\xCN\xbb'$
\xEH
\xEH
\xEH
\xEH
\xEH
\xEP

\hline

$(CM_n)$
\xEH
\xEH
\xEH
\xEH
\xEH
$(I_n)$
\xEP

$\xba\xcn\xbb_1, \ldots, \xba\xcn\xbb_n \xch $
\xEH
\xEH
\xEH
\xEH
\xEH
\xEP

$\xba \xcu \xbb_1 \xcu \ldots \xcu \xbb_{n-1} \xcL\xCN\xbb_n$
\xEH
\xEH
\xEH
\xEH
\xEH
\xEP

\hline

$(CM)$ Cautious Monotony
\xEH
$(CM)$
\xEH
$\xch$ (8.1)
\xEH
$ (\xbm CM) $
\xEH
$\xci$ (5.1)
\xEH
$(\xdm^+_\xbo)(4)$
\xEP

\cline{3-3}
\cline{5-5}

$ \xbf \xcn \xbq,  \xbf \xcn \xbq '   \xch $
\xEH
$T \xcc \ol{T' } \xcc \ol{ \ol{T} }  \xch $
\xEH
$\xci$ (8.2)
\xEH
$f(X) \xcc Y \xcc X  \xch $
\xEH
$\xch$ (5.2)
\xEH
\xEP

$ \xbf \xcu \xbq \xcn \xbq ' $
\xEH
$ \ol{ \ol{T} } \xcc \ol{ \ol{T' } }$
\xEH
\xEH
$f(Y) \xcc f(X)$
\xEH
\xEH
\xEP

\cline{1-1}

\cline{3-4}

or $(ResM)$ Restricted Monotony
\xEH
\xEH
$\xch$ (9.1)
\xEH
$(\xbm ResM)$
\xEH
\xEH
\xEP

\cline{3-3}

$ T  \xcn \xba, \xbb \xch T \xcv \{\xba\} \xcn \xbb $
\xEH
\xEH
$\xci$ (9.2)
\xEH
$ f(X) \xcc A \xcs B \xch f(X \xcs A) \xcc B $
\xEH
\xEH
\xEP

\hline

$(CUM)$ Cumulativity
\xEH
$(CUM)$
\xEH
$\xch$ (11.1)
\xEH
$( \xbm CUM)$
\xEH
$\xci$ (9.1)
\xEH
$(eM\xdi)+(I_\xbo)+(\xdm^{+}_{\xbo})(4)$
\xEP

\cline{3-3}
\cline{5-5}

$ \xbf \xcn \xbq   \xch $
\xEH
$T \xcc \ol{T' } \xcc \ol{ \ol{T} }  \xch $
\xEH
$\xci$ (11.2)
\xEH
$f(X) \xcc Y \xcc X  \xch $
\xEH
$\xcH$ (9.2)
\xEH
\xEP

$( \xbf \xcn \xbq '   \xcj   \xbf \xcu \xbq \xcn \xbq ' )$
\xEH
$ \ol{ \ol{T} }= \ol{ \ol{T' } }$
\xEH
\xEH
$f(Y)=f(X)$
\xEH
\xEH
\xEP

\hline

\xEH
$ (\xcc \xcd) $
\xEH
$\xch$ (10.1)
\xEH
$ (\xbm \xcc \xcd) $
\xEH
$\xci$ (10.1)
\xEH
$(eM\xdi)+(I_\xbo)+(eM\xdf)$
\xEP

\cline{3-3}
\cline{5-5}

\xEH
$T \xcc \ol{\ol{T'}}, T' \xcc \ol{\ol{T}} \xch $
\xEH
$\xci$ (10.2)
\xEH
$ f(X) \xcc Y, f(Y) \xcc X \xch $
\xEH
$\xcH$ (10.2)
\xEH
\xEP

\xEH
$ \ol{\ol{T'}} = \ol{\ol{T}}$
\xEH
\xEH
$ f(X)=f(Y) $
\xEH
\xEH
\xEP

\hline

\multicolumn{6}{|c|}{Rationality}
\xEP

\hline

$(RatM)$ Rational Monotony
\xEH
$(RatM)$
\xEH
$\xch$ (12.1)
\xEH
$( \xbm RatM)$
\xEH
$\xci$ (6.1)
\xEH
$(\xdm^{++})$
\xEP

\cline{3-3}
\cline{5-5}

$ \xbf \xcn \xbq,  \xbf \xcN \xCN \xbq '   \xch $
\xEH
$Con(T \xcv \ol{\ol{T'}})$, $T \xcl T'$ $ \xch $
\xEH
$\xci$ $(\xbm dp)$ (12.2)
\xEH
$X \xcc Y, X \xcs f(Y) \xEd \xCQ   \xch $
\xEH
$\xch$ (6.2)
\xEH
\xEP

\cline{3-3}

$ \xbf \xcu \xbq ' \xcn \xbq $
\xEH
$ \ol{\ol{T}} \xcd \ol{\ol{\ol{T'}} \xcv T} $
\xEH
$\xcI$ without $(\xbm dp)$ (12.3)
\xEH
$f(X) \xcc f(Y) \xcs X$
\xEH
\xEH
\xEP

\cline{3-3}

\xEH
\xEH
$\xci$ $T$ a formula (12.4)
\xEH
\xEH
\xEH
\xEP

\hline

\xEH
$(RatM=)$
\xEH
$\xch$ (13.1)
\xEH
$( \xbm =)$
\xEH
\xEH
\xEP

\cline{3-3}

\xEH
$Con(T \xcv \ol{\ol{T'}})$, $T \xcl T'$ $ \xch $
\xEH
$\xci$ $(\xbm dp)$ (13.2)
\xEH
$X \xcc Y, X \xcs f(Y) \xEd \xCQ   \xch $
\xEH
\xEH
\xEP

\cline{3-3}

\xEH
$ \ol{\ol{T}} = \ol{\ol{\ol{T'}} \xcv T} $
\xEH
$\xcI$ without $(\xbm dp)$ (13.3)
\xEH
$f(X) = f(Y) \xcs X$
\xEH
\xEH
\xEP

\cline{3-3}

\xEH
\xEH
$\xci$ $T$ a formula (13.4)
\xEH
\xEH
\xEH
\xEP

\hline

\xEH
$(Log=' )$
\xEH
$\xch$ (14.1)
\xEH
$( \xbm =' )$
\xEH
\xEH
\xEP

\cline{3-3}

\xEH
$Con( \ol{ \ol{T' } } \xcv T)$ $ \xch $
\xEH
$\xci$ $(\xbm dp)$ (14.2)
\xEH
$f(Y) \xcs X \xEd \xCQ $ $ \xch $
\xEH
\xEH
\xEP

\cline{3-3}

\xEH
$ \ol{ \ol{T \xcv T' } }= \ol{ \ol{ \ol{T' } } \xcv T}$
\xEH
$\xcI$ without $(\xbm dp)$ (14.3)
\xEH
$f(Y \xcs X)=f(Y) \xcs X$
\xEH
\xEH
\xEP

\cline{3-3}

\xEH
\xEH
$\xci$ $T$ a formula (14.4)
\xEH
\xEH
\xEH
\xEP

\hline

\xEH
$(Log \xFO )$
\xEH
$\xch$ (15.1)
\xEH
$( \xbm \xFO )$
\xEH
\xEH
\xEP

\cline{3-3}

\xEH
$ \ol{ \ol{T \xco T' } }$ is one of
\xEH
$\xci$ (15.2)
\xEH
$f(X \xcv Y)$ is one of
\xEH
\xEH
\xEP

\xEH
$\ol{\ol{T}},$ or $\ol{\ol{T'}},$ or $\ol{\ol{T}} \xcs \ol{\ol{T'}}$ (by (CCL))
\xEH
\xEH
$f(X),$ $f(Y)$ or $f(X) \xcv f(Y)$
\xEH
\xEH
\xEP

\hline

\xEH
$(Log \xcv )$
\xEH
$\xch$ $(\xbm\xcc)+(\xbm=)$ (16.1)
\xEH
$( \xbm \xcv )$
\xEH
\xEH
\xEP

\cline{3-3}

\xEH
$Con( \ol{ \ol{T' } } \xcv T),$ $ \xCN Con( \ol{ \ol{T' } }
\xcv \ol{ \ol{T} })$ $ \xch $
\xEH
$\xci$ $(\xbm dp)$ (16.2)
\xEH
$f(Y) \xcs (X-f(X)) \xEd \xCQ $ $ \xch $
\xEH
\xEH
\xEP

\cline{3-3}

\xEH
$ \xCN Con( \ol{ \ol{T \xco T' } } \xcv T' )$
\xEH
$\xcI$ without $(\xbm dp)$ (16.3)
\xEH
$f(X \xcv Y) \xcs Y= \xCQ$
\xEH
\xEH
\xEP

\hline

\xEH
$(Log \xcv ' )$
\xEH
$\xch$ $(\xbm\xcc)+(\xbm=)$ (17.1)
\xEH
$( \xbm \xcv ' )$
\xEH
\xEH
\xEP

\cline{3-3}

\xEH
$Con( \ol{ \ol{T' } } \xcv T),$ $ \xCN Con( \ol{ \ol{T' }
} \xcv \ol{ \ol{T} })$ $ \xch $
\xEH
$\xci$ $(\xbm dp)$ (17.2)
\xEH
$f(Y) \xcs (X-f(X)) \xEd \xCQ $ $ \xch $
\xEH
\xEH
\xEP

\cline{3-3}

\xEH
$ \ol{ \ol{T \xco T' } }= \ol{ \ol{T} }$
\xEH
$\xcI$ without $(\xbm dp)$ (17.3)
\xEH
$f(X \xcv Y)=f(X)$
\xEH
\xEH
\xEP

\hline

\xEH
\xEH
\xEH
$( \xbm \xbe )$
\xEH
\xEH
\xEP

\xEH
\xEH
\xEH
$a \xbe X-f(X)$ $ \xch $
\xEH
\xEH
\xEP

\xEH
\xEH
\xEH
$ \xcE b \xbe X.a \xce f(\{a,b\})$
\xEH
\xEH
\xEP

\hline

\end{tabular}

}

\end{turn}

The proof of the following Fact - together with the subsequent examples -
requires some knowledge of preferential structures, which will be
introduced ``officially'' only
in Chapter \ref{Chapter Pref} (page \pageref{Chapter Pref}).
We chose to give those results here, so the reader will have
immediately a global picture, and can come back later, if desired,
and read the proofs and (counter)examples.
\index{Fact Mu-Base}

\bfa

$\hspace{0.01em}$

% (+++ Orig. No.:  Fact Mu-Base +++)

\label{Fact Mu-Base}

The following table is to be read as follows: If the left hand side holds
for
some
function $f: \xdy \xcp \xdp (U),$ and the auxiliary properties noted in
the middle also
hold for $f$ or $ \xdy,$ then the right hand side will hold, too - and
conversely.

{\small

% \begin{tabular*}{12.8cm}{|c|c|c|c|}
\begin{tabular}{|c|c|c|c|}

\hline

\multicolumn{4}{|c|}{Basics} \xEP

\hline

(1.1)
\xEH
$(\xbm PR)$
\xEH
$\xch$ $(\xcs)+(\xbm \xcc)$
\xEH
$(\xbm PR')$
\xEP

\cline{1-1}

\cline{3-3}

(1.2)
\xEH
\xEH
$\xci$
\xEH
\xEP

\hline

(2.1)
\xEH
$(\xbm PR)$
\xEH
$\xch$ $(\xbm \xcc)$
\xEH
$(\xbm OR)$
\xEP

\cline{1-1}

\cline{3-3}

(2.2)
\xEH
\xEH
$\xci$ $(\xbm \xcc)$ + $(-)$
\xEH
\xEP

\cline{1-1}

\cline{3-4}

(2.3)
\xEH
\xEH
$\xch$ $(\xbm \xcc)$
\xEH
$(\xbm wOR)$
\xEP

\cline{1-1}

\cline{3-3}

(2.4)
\xEH
\xEH
$\xci$ $(\xbm \xcc)$ + $(-)$
\xEH
\xEP

\hline

(3)
\xEH
$(\xbm PR)$
\xEH
$\xch$
\xEH
$( \xbm CUT)$
\xEP

\hline

(4)
\xEH
$(\xbm \xcc )+(\xbm \xcc \xcd )+(\xbm CUM)+$
\xEH
$\xcH$
\xEH
$( \xbm PR)$
\xEP

\xEH
$(\xbm RatM)+(\xcs )$
\xEH
\xEH
\xEP

\hline

\multicolumn{4}{|c|}{Cumulativity} \xEP

\hline

(5.1)
\xEH
$(\xbm CM)$
\xEH
$\xch$ $(\xcs)+(\xbm \xcc)$
\xEH
$(\xbm ResM)$
\xEP

\cline{1-1}

\cline{3-3}

(5.2)
\xEH
\xEH
$\xci$ (infin.)
\xEH
\xEP

\hline

(6)
\xEH
$(\xbm CM)+(\xbm CUT)$
\xEH
$\xcj$
\xEH
$(\xbm CUM)$
\xEP

\hline

(7)
\xEH
$( \xbm \xcc )+( \xbm \xcc \xcd )$
\xEH
$\xch$
\xEH
$( \xbm CUM)$
\xEP

\hline

(8)
\xEH
$( \xbm \xcc )+( \xbm CUM)+( \xcs )$
\xEH
$\xch$
\xEH
$( \xbm \xcc \xcd )$
\xEP

\hline

(9)
\xEH
$( \xbm \xcc )+( \xbm CUM)$
\xEH
$\xcH$
\xEH
$( \xbm \xcc \xcd )$
\xEP

\hline

\multicolumn{4}{|c|}{Rationality} \xEP

\hline

(10)
\xEH
$( \xbm RatM )+( \xbm PR )$
\xEH
$\xch$
\xEH
$( \xbm =)$
\xEP

\hline

(11)
\xEH
$( \xbm =)$
\xEH
$ \xch $
\xEH
$( \xbm PR)+(\xbm RatM)$
\xEP

\hline

(12.1)
\xEH
$( \xbm =)$
\xEH
$ \xch $ $(\xcs)+( \xbm \xcc )$
\xEH
$( \xbm =' )$
\xEP
\cline{1-1}
\cline{3-3}
(12.2)
\xEH
\xEH
$ \xci $
\xEH
\xEP

\hline

(13)
\xEH
$( \xbm \xcc )+( \xbm =)$
\xEH
$ \xch $ $(\xcv)$
\xEH
$( \xbm \xcv )$
\xEP

\hline

(14)
\xEH
$( \xbm \xcc )+( \xbm \xCQ )+( \xbm =)$
\xEH
$ \xch $ $(\xcv)$
\xEH
$( \xbm \xFO ),$ $( \xbm \xcv ' ),$ $( \xbm CUM)$
\xEP

\hline

(15)
\xEH
$( \xbm \xcc )+( \xbm \xFO )$
\xEH
$ \xch $ $(-)$ of $\xdy$
\xEH
$( \xbm =)$
\xEP

\hline

(16)
\xEH
$( \xbm \xFO )+( \xbm \xbe )+( \xbm PR)+$
\xEH
$ \xch $ $(\xcv)$ + $\xdy$ contains singletons
\xEH
$( \xbm =)$
\xEP
\xEH
$( \xbm \xcc )$
\xEH
\xEH
\xEP

\hline

(17)
\xEH
$( \xbm CUM)+( \xbm =)$
\xEH
$ \xch $ $(\xcv)$ + $\xdy$ contains singletons
\xEH
$( \xbm \xbe )$
\xEP

\hline

(18)
\xEH
$( \xbm CUM)+( \xbm =)+( \xbm \xcc )$
\xEH
$ \xch $ $(\xcv)$
\xEH
$( \xbm \xFO )$
\xEP

\hline

(19)
\xEH
$( \xbm PR)+( \xbm CUM)+( \xbm \xFO )$
\xEH
$ \xch $ sufficient, e.g. true in $\xdD_{\xdl}$
\xEH
$( \xbm =)$.
\xEP

\hline

(20)
\xEH
$( \xbm \xcc )+( \xbm PR)+( \xbm =)$
\xEH
$ \xcH $
\xEH
$( \xbm \xFO )$
\xEP

\hline

(21)
\xEH
$( \xbm \xcc )+( \xbm PR)+( \xbm \xFO )$
\xEH
$ \xcH $ (without $(-)$)
\xEH
$( \xbm =)$
\xEP

\hline

(22)
\xEH
$( \xbm \xcc )+( \xbm PR)+( \xbm \xFO )+$
\xEH
$ \xcH $
\xEH
$( \xbm \xbe )$
\xEP
\xEH
$( \xbm =)+( \xbm \xcv )$
\xEH
\xEH
(thus not representability
\xEP
\xEH
\xEH
\xEH
by ranked structures)
\xEP

\hline

\end{tabular}

}

\index{Fact Mu-Base Proof}

\efa

\subparagraph{
Proof
}

$\hspace{0.01em}$

% (+++ Orig.:  Proof +++)

All sets are to be in $ \xdy.$

(1.1) $( \xbm PR)+( \xcs )+( \xbm \xcc )$ $ \xch $ $( \xbm PR' ):$

By $X \xcs Y \xcc X$ and $( \xbm PR),$ $f(X) \xcs X \xcs Y \xcc f(X \xcs
Y).$ By $( \xbm \xcc )$ $f(X) \xcs Y=f(X) \xcs X \xcs Y.$

(1.2) $( \xbm PR' ) \xch ( \xbm PR):$

Let $X \xcc Y,$ so $X=X \xcs Y,$ so by $( \xbm PR' )$ $f(Y) \xcs X \xcc
f(X \xcs Y)=f(X).$

(2.1) $( \xbm PR)+( \xbm \xcc )$ $ \xch $ $( \xbm OR):$

$f(X \xcv Y) \xcc X \xcv Y$ by $( \xbm \xcc ),$ so $f(X \xcv Y)$ $=$ $(f(X
\xcv Y) \xcs X) \xcv (f(X \xcv Y) \xcs Y)$ $ \xcc $ $f(X) \xcv f(Y).$

(2.2) $( \xbm OR)$ $+$ $( \xbm \xcc )$ $+$ $ \xCf (-)$ $ \xch $ $( \xbm
PR):$

Let $X \xcc Y,$ $X':=Y-X$. $f(Y) \xcc f(X) \xcv f(X' )$ by $( \xbm OR),$
so $f(Y) \xcs X$ $ \xcc $
$(f(X) \xcs X) \xcv (f(X' ) \xcs X)$ $=_{( \xbm \xcc )}$ $f(X) \xcv \xCQ $
$=$ $f(X).$

(2.3) $( \xbm PR)+( \xbm \xcc )$ $ \xch $ $( \xbm wOR):$

Trivial by (2.1).

(2.4) $( \xbm wOR)$ $+$ $( \xbm \xcc )$ $+$ $ \xCf (-)$ $ \xch $ $( \xbm
PR):$

Let $X \xcc Y,$ $X':=Y-X$. $f(Y) \xcc f(X) \xcv X' $ by $( \xbm wOR),$
so $f(Y) \xcs X$ $ \xcc $
$(f(X) \xcs X) \xcv (X' \xcs X)$ $=_{( \xbm \xcc )}$ $f(X) \xcv \xCQ $ $=$
$f(X).$

(3) $( \xbm PR)$ $ \xch $ $( \xbm CUT):$

$f(X) \xcc Y \xcc X$ $ \xch $ $f(X) \xcc f(X) \xcs Y \xcc f(Y)$ by $( \xbm
PR).$

(4) $( \xbm \xcc )+( \xbm \xcc \xcd )+( \xbm CUM)+( \xbm RatM)+( \xcs )$ $
\xcH $ $( \xbm PR):$

This is shown in Example \ref{Example Need-Pr} (page \pageref{Example Need-Pr})
.

(5.1) $( \xbm CM)+( \xcs )+( \xbm \xcc )$ $ \xch $ $( \xbm ResM):$

Let $f(X) \xcc A \xcs B,$ so $f(X) \xcc A,$ so by $( \xbm \xcc )$ $f(X)
\xcc A \xcs X \xcc X,$
so by $( \xbm CM)$ $f(A \xcs X) \xcc f(X) \xcc B.$

(5.2) $( \xbm ResM) \xch ( \xbm CM):$

We consider here the infinitary version, where all sets can be model sets
of infinite theories.
Let $f(X) \xcc Y \xcc X,$ so $f(X) \xcc Y \xcs f(X),$ so by $( \xbm ResM)$
$f(Y)=f(X \xcs Y) \xcc f(X).$

(6) $( \xbm CM)+( \xbm CUT)$ $ \xcj $ $( \xbm CUM):$

Trivial.

(7) $( \xbm \xcc )+( \xbm \xcc \xcd )$ $ \xch $ $( \xbm CUM):$

Suppose $f(D) \xcc E \xcc D.$ So by $( \xbm \xcc )$ $f(E) \xcc E \xcc D,$
so by $( \xbm \xcc \xcd )$ $f(D)=f(E).$

(8) $( \xbm \xcc )+( \xbm CUM)+( \xcs )$ $ \xch $ $( \xbm \xcc \xcd ):$

Let $f(D) \xcc E,$ $f(E) \xcc D,$ so by $( \xbm \xcc )$ $f(D) \xcc D \xcs
E \xcc D,$ $f(E) \xcc D \xcs E \xcc E.$ As $f(D \xcs E)$
is defined, so $f(D)=f(D \xcs E)=f(E)$ by $( \xbm CUM).$

(9) $( \xbm \xcc )+( \xbm CUM)$ $ \xcH $ $( \xbm \xcc \xcd ):$

This is shown in Example \ref{Example Mu-Cum-Cd} (page \pageref{Example
Mu-Cum-Cd}).

(10) $( \xbm RatM)+( \xbm PR)$ $ \xch $ $( \xbm =):$

Trivial.

(11) $( \xbm =)$ entails $( \xbm PR)$ and $( \xbm RatM):$

Trivial.

(12.1) $( \xbm =) \xch ( \xbm =' ):$

Let $f(Y) \xcs X \xEd \xCQ,$ we have to show $f(X \xcs Y)=f(Y) \xcs X.$
By $( \xbm \xcc )$ $f(Y) \xcc Y,$ so $f(Y) \xcs X=f(Y) \xcs (X \xcs Y),$
so by $( \xbm =)$ $f(Y) \xcs X$ $=$
$f(Y) \xcs (X \xcs Y)$ $=$ $f(X \xcs Y).$

(12.2) $( \xbm =' ) \xch ( \xbm =):$

Let $X \xcc Y,$ $f(Y) \xcs X \xEd \xCQ,$ then $f(X)=f(Y \xcs X)=f(Y) \xcs
X.$

(13) $( \xbm \xcc ),$ $( \xbm =)$ $ \xch $ $( \xbm \xcv ):$

If not, $f(X \xcv Y) \xcs Y \xEd \xCQ,$ but $f(Y) \xcs (X-f(X)) \xEd \xCQ
.$ By (11), $( \xbm PR)$ holds,
so $f(X \xcv Y) \xcs X \xcc f(X),$ so $ \xCQ $ $ \xEd $ $f(Y) \xcs
(X-f(X))$ $ \xcc $ $f(Y) \xcs (X-f(X \xcv Y)),$ so
$f(Y)-f(X \xcv Y) \xEd \xCQ,$ so by $( \xbm \xcc )$ $f(Y) \xcc Y$ and
$f(Y) \xEd f(X \xcv Y) \xcs Y.$
But by $( \xbm =)$ $f(Y)=f(X \xcv Y) \xcs Y,$ a contradiction.

(14)

$( \xbm \xcc ),$ $( \xbm \xCQ ),$ $( \xbm =)$ $ \xch $ $( \xbm \xFO ):$

If $X$ or $Y$ or both are empty, then this is trivial.
Assume then $X \xcv Y \xEd \xCQ,$ so by $( \xbm \xCQ )$ $f(X \xcv Y) \xEd
\xCQ.$
By $( \xbm \xcc )$ $f(X \xcv Y) \xcc X \xcv Y,$ so $f(X \xcv Y) \xcs X=
\xCQ $ and $f(X \xcv Y) \xcs Y= \xCQ $ together are
impossible.
Case 1, $f(X \xcv Y) \xcs X \xEd \xCQ $ and $f(X \xcv Y) \xcs Y \xEd \xCQ
:$ By $( \xbm =)$ $f(X \xcv Y) \xcs X=f(X)$ and
$f(X \xcv Y) \xcs Y=f(Y),$ so by $( \xbm \xcc )$ $f(X \xcv Y)=f(X) \xcv
f(Y).$
Case 2, $f(X \xcv Y) \xcs X \xEd \xCQ $ and $f(X \xcv Y) \xcs Y= \xCQ:$
So by $( \xbm =)$ $f(X \xcv Y)=f(X \xcv Y) \xcs X=f(X).$
Case 3, $f(X \xcv Y) \xcs X= \xCQ $ and $f(X \xcv Y) \xcs Y \xEd \xCQ:$
Symmetrical.

$( \xbm \xcc ),$ $( \xbm \xCQ ),$ $( \xbm =)$ $ \xch $ $( \xbm \xcv ' ):$

Let $f(Y) \xcs (X-f(X)) \xEd \xCQ.$
If $X \xcv Y= \xCQ,$ then $f(X \xcv Y)=f(X)= \xCQ $ by $( \xbm \xcc ).$
So suppose $X \xcv Y \xEd \xCQ.$ By
(13), $f(X \xcv Y) \xcs Y= \xCQ,$ so $f(X \xcv Y) \xcc X$ by $( \xbm \xcc
).$ By $( \xbm \xCQ ),$ $f(X \xcv Y) \xEd \xCQ,$ so
$f(X \xcv Y) \xcs X \xEd \xCQ,$ and $f(X \xcv Y)=f(X)$ by $( \xbm =).$

$( \xbm \xcc ),$ $( \xbm \xCQ ),$ $( \xbm =)$ $ \xch $ $( \xbm CUM):$

Let $f(Y) \xcc X \xcc Y.$
If $Y= \xCQ,$ this is trivial by $( \xbm \xcc ).$ If $Y \xEd \xCQ,$ then
by $( \xbm \xCQ )$ - which is
crucial here - $f(Y) \xEd \xCQ,$ so by $f(Y) \xcc X$ $f(Y) \xcs X \xEd
\xCQ,$ so by $( \xbm =)$
$f(Y)=f(Y) \xcs X=f(X).$

(15) $( \xbm \xcc )+( \xbm \xFO )$ $ \xch $ $( \xbm =):$

Let $X \xcc Y,$ $X \xcs f(Y) \xEd \xCQ,$ and consider $Y=X \xcv (Y-$X).
Then $f(Y)=f(X) \xFO f(Y-$X). As
$f(Y) \xcs X \xEd \xCQ,$ $f(Y)=f(Y-$X) is impossible. Otherwise,
$f(X)=f(Y) \xcs X,$ and we are
done.

(16) $( \xbm \xFO )+( \xbm \xbe )+( \xbm PR)+( \xbm \xcc )$ $ \xch $ $(
\xbm =):$

Suppose $X \xcc Y,$ $x \xbe f(Y) \xcs X,$ we have to show $f(Y) \xcs
X=f(X).$ `` $ \xcc $ '' is trivial
by $( \xbm PR).$ `` $ \xcd $ '': Assume $a \xce f(Y)$ (by $( \xbm \xcc )),$
but $a \xbe f(X).$ By $( \xbm \xbe )$ $ \xcE b \xbe Y.a \xce f(\{a,b\}).$
As $a \xbe f(X),$ by $( \xbm PR),$ $a \xbe f(\{a,x\}).$ By $( \xbm \xFO
),$ $f(\{a,b,x\})$ $=$ $f(\{a,x\}) \xFO f(\{b\}).$
As $a \xce f(\{a,b,x\}),$ $f(\{a,b,x\})$ $=$ $f(\{b\}),$ so $x \xce
f(\{a,b,x\}),$ contradicting $( \xbm PR),$
as $a,b,x \xbe Y.$

(17) $( \xbm CUM)+( \xbm =)$ $ \xch $ $( \xbm \xbe ):$

Let $a \xbe X-f(X).$ If $f(X)= \xCQ,$ then $f(\{a\})= \xCQ $ by $( \xbm
CUM).$ If not: Let
$b \xbe f(X),$ then $a \xce f(\{a,b\})$ by $( \xbm =).$

(18) $( \xbm CUM)+( \xbm =)+( \xbm \xcc )$ $ \xch $ $( \xbm \xFO ):$

By $( \xbm CUM),$ $f(X \xcv Y) \xcc X \xcc X \xcv Y$ $ \xch $ $f(X)=f(X
\xcv Y),$ and $f(X \xcv Y) \xcc Y \xcc X \xcv Y$ $ \xch $
$f(Y)=f(X \xcv Y).$ Thus, if $( \xbm \xFO )$ were to fail, $f(X \xcv Y)
\xcC X,$ $f(X \xcv Y) \xcC Y,$ but then
by $( \xbm \xcc )$ $f(X \xcv Y) \xcs X \xEd \xCQ,$ so $f(X)=f(X \xcv Y)
\xcs X,$ and $f(X \xcv Y) \xcs Y \xEd \xCQ,$ so
$f(Y)=f(X \xcv Y) \xcs Y$ by $( \xbm =).$ Thus, $f(X \xcv Y)$ $=$ $(f(X
\xcv Y) \xcs X) \xcv (f(X \xcv Y) \xcs Y)$ $=$
$f(X) \xcv f(Y).$

(19) $( \xbm PR)+( \xbm CUM)+( \xbm \xFO )$ $ \xch $ $( \xbm =):$

Suppose $( \xbm =)$ does not hold. So, by $( \xbm PR),$ there are $X,Y,y$
s.t. $X \xcc Y,$ $X \xcs f(Y) \xEd \xCQ,$
$y \xbe Y-f(Y),$ $y \xbe f(X).$ Let $a \xbe X \xcs f(Y).$ If $f(Y)=\{a\},$
then by $( \xbm CUM)$ $f(Y)=f(X),$ so
there must be $b \xbe f(Y),$ $b \xEd a.$ Take now $Y',$ $Y'' $ s.t. $Y=Y'
\xcv Y'',$ $a \xbe Y',$ $a \xce Y'',$ $b \xbe Y'',$
$b \xce Y',$ $y \xbe Y' \xcs Y''.$ Assume now $( \xbm \xFO )$ to hold,
we show a contradiction.
If $y \xce f(Y'' ),$ then by $( \xbm PR)$ $y \xce f(Y'' \xcv \{a\}).$ But
$f(Y'' \xcv \{a\})$ $=$ $f(Y'' ) \xFO f(\{a,y\}),$
so $f(Y'' \xcv \{a\})=f(Y'' ),$ contradicting $a \xbe f(Y).$ If $y \xbe
f(Y'' ),$ then by $f(Y)$ $=$
$f(Y' ) \xFO f(Y'' ),$ $f(Y)=f(Y' ),$ $contradiction$ as $b \xce f(Y' ).$

(20) $( \xbm \xcc )+( \xbm PR)+( \xbm =)$ $ \xcH $ $( \xbm \xFO ):$

See Example \ref{Example Mu-Barbar} (page \pageref{Example Mu-Barbar}).

(21) $( \xbm \xcc )+( \xbm PR)+( \xbm \xFO )$ $ \xcH $ $( \xbm =):$

See Example \ref{Example Mu-Equal} (page \pageref{Example Mu-Equal}).

(22) $( \xbm \xcc )+( \xbm PR)+( \xbm \xFO )+( \xbm =)+( \xbm \xcv )$ $
\xcH $ $( \xbm \xbe ):$

See Example \ref{Example Mu-Epsilon} (page \pageref{Example Mu-Epsilon}).

Thus, by Fact \ref{Fact Rank-Hold} (page \pageref{Fact Rank-Hold}), the
conditions do not assure
representability by ranked structures.

$ \xcz $
\\[3ex]
\index{Remark RatM=}

\br

$\hspace{0.01em}$

% (+++ Orig. No.:  Remark RatM= +++)

\label{Remark RatM=}

Note that $( \xbm =' )$ is very close to $ \xCf (RatM):$ $ \xCf (RatM)$
says:
$ \xba \xcn \xbb,$ $ \xba \xcN \xCN \xbg $ $ \xch $ $ \xba \xcu \xbg \xcn
\xbb.$ Or, $f(A) \xcc B,$ $f(A) \xcs C \xEd \xCQ $ $ \xch $
$f(A \xcs C) \xcc B$ for all $A,B,C.$ This is not quite, but almost: $f(A
\xcs C) \xcc f(A) \xcs C$
(it depends how many $B$ there are, if $f(A)$ is some such $B,$ the fit is
perfect).
\index{Example Mu-Cum-Cd}

\er

\be

$\hspace{0.01em}$

% (+++ Orig. No.:  Example Mu-Cum-Cd +++)

\label{Example Mu-Cum-Cd}

We show here $( \xbm \xcc )+( \xbm CUM)$ $ \xcH $ $( \xbm \xcc \xcd ).$

Consider $X:=\{a,b,c\},$ $Y:=\{a,b,d\},$ $f(X):=\{a\},$ $f(Y):=\{a,b\},$ $
\xdy:=\{X,Y\}.$
(If $f(\{a,b\})$ were defined, we would have $f(X)=f(\{a,b\})=f(Y),$
$contradiction.)$

Obviously, $( \xbm \xcc )$ and $( \xbm CUM)$ hold, but not $( \xbm \xcc
\xcd ).$

$ \xcz $
\\[3ex]
\index{Example Need-Pr}

\ee

\be

$\hspace{0.01em}$

% (+++ Orig. No.:  Example Need-Pr +++)

\label{Example Need-Pr}

We show here $( \xbm \xcc )+( \xbm \xcc \xcd )+( \xbm CUM)+( \xbm RatM)+(
\xcs )$ $ \xcH $ $( \xbm PR).$

Let $U:=\{a,b,c\}.$ Let $ \xdy = \xdp (U).$ So $( \xcs )$ is trivially
satisfied.
Set $f(X):=X$ for all $X \xcc U$ except for $f(\{a,b\})=\{b\}.$ Obviously,
this cannot be represented by a preferential structure and $( \xbm PR)$ is
false
for $U$ and $\{a,b\}.$ But it satisfies $( \xbm \xcc ),$ $( \xbm CUM),$ $(
\xbm RatM).$ $( \xbm \xcc )$ is trivial.
$( \xbm CUM):$ Let $f(X) \xcc Y \xcc X.$ If $f(X)=X,$ we are done.
Consider $f(\{a,b\})=\{b\}.$ If
$\{b\} \xcc Y \xcc \{a,b\},$ then $f(Y)=\{b\},$ so we are done again. It
is shown in
Fact \ref{Fact Mu-Base} (page \pageref{Fact Mu-Base}), (8) that $( \xbm \xcc
\xcd )$ follows.
$( \xbm RatM):$ Suppose $X \xcc Y,$ $X \xcs f(Y) \xEd \xCQ,$ we have to
show $f(X) \xcc f(Y) \xcs X.$ If $f(Y)=Y,$ the result holds by $X \xcc Y,$
so it does if $X=Y.$
The only remaining case is $Y=\{a,b\},$ $X=\{b\},$ and the result holds
again.

$ \xcz $
\\[3ex]
\index{Example Mu-Barbar}

\ee

\be

$\hspace{0.01em}$

% (+++ Orig. No.:  Example Mu-Barbar +++)

\label{Example Mu-Barbar}

The example shows that $( \xbm \xcc )+( \xbm PR)+( \xbm =)$ $ \xcH $ $(
\xbm \xFO ).$

Consider the following structure without transitivity:
$U:=\{a,b,c,d\},$ $c$ and $d$ have $ \xbo $ many copies in descending
order $c_{1} \xed c_{2}$  \Xl., etc.
$a,b$ have one single copy each. $a \xed b,$ $a \xed d_{1},$ $b \xed a,$
$b \xed c_{1}.$
$( \xbm \xFO )$ does not hold: $f(U)= \xCQ,$ but $f(\{a,c\})=\{a\},$
$f(\{b,d\})=\{b\}.$
$( \xbm PR)$ holds as in all preferential structures.
$( \xbm =)$ holds: If it were to fail, then for some $A \xcc B,$ $f(B)
\xcs A \xEd \xCQ,$ so $f(B) \xEd \xCQ.$
But the only possible cases for $B$ are now: $(a \xbe B,$ $b,d \xce B)$ or
$(b \xbe B,$ $a,c \xce B).$
Thus, $B$ can be $\{a\},$ $\{a,c\},$ $\{b\},$ $\{b,d\}$ with $f(B)=$
$\{a\},$ $\{a\},$ $\{b\},$ $\{b\}.$
If $A=B,$ then the result will hold trivially. Moreover, $ \xCf A$ has to
be $ \xEd \xCQ.$
So the remaining cases of $B$ where it might fail are $B=$ $\{a,c\}$ and
$\{b,d\},$ and
by $f(B) \xcs A \xEd \xCQ,$ the only cases of $ \xCf A$ where it might
fail, are $A=$ $\{a\}$ or $\{b\}$
respectively.
So the only cases remaining are: $B=\{a,c\},$ $A=\{a\}$ and $B=\{b,d\},$
$A=\{b\}.$
In the first case, $f(A)=f(B)=\{a\},$ in the second $f(A)=f(B)=\{b\},$ but
$( \xbm =)$
holds in both.

$ \xcz $
\\[3ex]
\index{Example Mu-Equal}

\ee

\be

$\hspace{0.01em}$

% (+++ Orig. No.:  Example Mu-Equal +++)

\label{Example Mu-Equal}

The example shows that $( \xbm \xcc )+( \xbm PR)+( \xbm \xFO )$ $ \xcH $
$( \xbm =).$

Work in the set of theory definable model sets of an infinite
propositional
language. Note that this is not closed under set difference, and closure
properties will play a crucial role in the argumentation.
Let $U:=\{y,a,x_{i< \xbo }\},$ where $x_{i} \xcp a$ in the standard
topology. For the order,
arrange s.t. $y$ is minimized by any set iff this set contains a cofinal
subsequence of
the $x_{i},$ this can be done by the standard construction. Moreover, let
the $x_{i}$
all kill themselves, i.e. with $ \xbo $ many copies $x^{1}_{i} \xed
x^{2}_{i} \xed $  \Xl. There are no other
elements in the relation. Note that if $a \xce \xbm (X),$ then $a \xce X,$
and $X$ cannot contain
a cofinal subsequence of the $x_{i},$ as $X$ is closed in the standard
topology.
(A short argument: suppose $X$ contains such a subsequence, but $a \xce
X.$ Then the
theory of a $Th(a)$ is inconsistent with $Th(X),$ so already a finite
subset of
$Th(a)$ is inconsistent with $Th(X),$ but such a finite subset will
finally hold
in a cofinal sequence converging to a.)
Likewise, if $y \xbe \xbm (X),$ then $X$ cannot contain a cofinal
subsequence of the $x_{i}.$

Obviously, $( \xbm \xcc )$ and $( \xbm PR)$ hold, but $( \xbm =)$ does not
hold: Set $B:=U,$ $A:=\{a,y\}.$
Then $ \xbm (B)=\{a\},$ $ \xbm (A)=\{a,y\},$ contradicting $( \xbm =).$

It remains to show that $( \xbm \xFO )$ holds.

$ \xbm (X)$ can only be $ \xCQ,$ $\{a\},$ $\{y\},$ $\{a,y\}.$ As $ \xbm
(A \xcv B) \xcc \xbm (A) \xcv \xbm (B)$ by $( \xbm PR),$

Case 1, $ \xbm (A \xcv B)=\{a,y\}$ is settled.

Note that if $y \xbe X- \xbm (X),$ then $X$ will contain a cofinal
subsequence, and thus
$a \xbe \xbm (X).$

Case 2: $ \xbm (A \xcv B)=\{a\}.$

Case 2.1: $ \xbm (A)=\{a\}$ - we are done.

Case 2.2: $ \xbm (A)=\{y\}:$ $ \xCf A$ does not contain $ \xCf a,$ nor a
cofinal subsequence.
If $ \xbm (B)= \xCQ,$ then $a \xce B,$ so $a \xce A \xcv B,$ a
contradiction.
If $ \xbm (B)=\{a\},$ we are done.
If $y \xbe \xbm (B),$ then $y \xbe B,$ but $B$ does not contain a cofinal
subsequence, so
$A \xcv B$ does not either, so $y \xbe \xbm (A \xcv B),$ $contradiction.$

Case 2.3: $ \xbm (A)= \xCQ:$ $ \xCf A$ cannot contain a cofinal
subsequence.
If $ \xbm (B)=\{a\},$ we are done.
$a \xbe \xbm (B)$ does have to hold, so $ \xbm (B)=\{a,y\}$ is the only
remaining possibility.
But then $B$ does not contain a cofinal subsequence, and neither does $A
\xcv B,$ so
$y \xbe \xbm (A \xcv B),$ $contradiction.$

Case 2.4: $ \xbm (A)=\{a,y\}:$ $ \xCf A$ does not contain a cofinal
subsequence.
If $ \xbm (B)=\{a\},$ we are done.
If $ \xbm (B)= \xCQ,$ $B$ does not contain a cofinal subsequence (as $a
\xce B),$ so neither
does $A \xcv B,$ so $y \xbe \xbm (A \xcv B),$ $contradiction.$
If $y \xbe \xbm (B),$ $B$ does not contain a cofinal subsequence, and we
are done again.

Case 3: $ \xbm (A \xcv B)=\{y\}:$
To obtain a contradiction, we need $a \xbe \xbm (A)$ or $a \xbe \xbm (B).$
But in both cases
$a \xbe \xbm (A \xcv B).$

Case 4: $ \xbm (A \xcv B)= \xCQ:$
Thus, $A \xcv B$ contains no cofinal subsequence. If, e.g. $y \xbe \xbm
(A),$ then $y \xbe \xbm (A \xcv B),$
if $a \xbe \xbm (A),$ then $a \xbe \xbm (A \xcv B),$ so $ \xbm (A)= \xCQ
.$

$ \xcz $
\\[3ex]
\index{Example Mu-Epsilon}

\ee

\be

$\hspace{0.01em}$

% (+++ Orig. No.:  Example Mu-Epsilon +++)

\label{Example Mu-Epsilon}

The example show that $( \xbm \xcc )+( \xbm PR)+( \xbm \xFO )+( \xbm =)+(
\xbm \xcv )$ $ \xcH $ $( \xbm \xbe ).$

Let $U:=\{y,x_{i< \xbo }\},$ $x_{i}$ a sequence, each $x_{i}$ kills
itself, $x^{1}_{i} \xed x^{2}_{i} \xed  \Xl $
and $y$ is killed by all cofinal subsequences of the $x_{i}.$ Then for any
$X \xcc U$
$ \xbm (X)= \xCQ $ or $ \xbm (X)=\{y\}.$

$( \xbm \xcc )$ and $( \xbm PR)$ hold obviously.

$( \xbm \xFO ):$ Let $A \xcv B$ be given. If $y \xce X,$ then for all $Y
\xcc X$ $ \xbm (Y)= \xCQ.$
So, if $y \xce A \xcv B,$ we are done. If $y \xbe A \xcv B,$ if $ \xbm (A
\xcv B)= \xCQ,$ one of $A,B$ must contain
a cofinal sequence, it will have $ \xbm = \xCQ.$ If not, then $ \xbm (A
\xcv B)=\{y\},$ and this will
also hold for the one $y$ is in.

$( \xbm =):$ Let $A \xcc B,$ $ \xbm (B) \xcs A \xEd \xCQ,$ show $ \xbm
(A)= \xbm (B) \xcs A.$ But now $ \xbm (B)=\{y\},$ $y \xbe A,$
so $B$ does not contain a cofinal subsequence, neither does A, so $ \xbm
(A)=\{y\}.$

$( \xbm \xcv ):$ $(A- \xbm (A)) \xcs \xbm (A' ) \xEd \xCQ,$ so $ \xbm (A'
)=\{y\},$ so $ \xbm (A \xcv A' )= \xCQ,$ as $y \xbe A- \xbm (A).$

But $( \xbm \xbe )$ does not hold: $y \xbe U- \xbm (U),$ but there is no
$x$ s.t. $y \xce \xbm (\{x,y\}).$

$ \xcz $
\\[3ex]

\ee

We turn to interdependencies of the different $ \xbm -$conditions.
Again, we will sometimes use preferential structures in our arguments.
\index{Fact Mwor}

\bfa

$\hspace{0.01em}$

% (+++ Orig. No.:  Fact Mwor +++)

\label{Fact Mwor}

$( \xbm wOR)+( \xbm \xcc )$ $ \xch $ $f(X \xcv Y) \xcc f(X) \xcv f(Y) \xcv
(X \xcs Y)$
\index{Fact Mwor Proof}

\efa

\subparagraph{
Proof
}

$\hspace{0.01em}$

% (+++ Orig.:  Proof +++)

$f(X \xcv Y) \xcc f(X) \xcv Y,$ $f(X \xcv Y) \xcc X \xcv f(Y),$ so $f(X
\xcv Y)$ $ \xcc $ $(f(X) \xcv Y) \xcs (X \xcv f(Y))$ $=$
$f(X) \xcv f(Y) \xcv (X \xcs Y)$ $ \xcz $
\\[3ex]
\index{Proposition Alg-Log}

\bp

$\hspace{0.01em}$

% (+++ Orig. No.:  Proposition Alg-Log +++)

\label{Proposition Alg-Log}

The following table is to be read as follows:

Let a logic $ \xcn $ satisfy $ \xCf (LLE)$ and $ \xCf (CCL),$ and define a
function $f: \xdD_{ \xdl } \xcp \xdD_{ \xdl }$
by $f(M(T)):=M( \ol{ \ol{T} }).$ Then $f$ is well defined, satisfies $(
\xbm dp),$ and $ \ol{ \ol{T} }=Th(f(M(T))).$

If $ \xcn $ satisfies a rule in the left hand side,
then - provided the additional properties noted in the middle for $ \xch $
hold, too -
$f$ will satisfy the property in the right hand side.

Conversely, if $f: \xdy \xcp \xdp (M_{ \xdl })$ is a function, with $
\xdD_{ \xdl } \xcc \xdy,$ and we define a logic
$ \xcn $ by $ \ol{ \ol{T} }:=Th(f(M(T))),$ then $ \xcn $ satisfies $ \xCf
(LLE)$ and $ \xCf (CCL).$
If $f$ satisfies $( \xbm dp),$ then $f(M(T))=M( \ol{ \ol{T} }).$

If $f$ satisfies a property in the right hand side,
then - provided the additional properties noted in the middle for $ \xci $
hold, too -
$ \xcn $ will satisfy the property in the left hand side.

If ``formula'' is noted in the table, this means that, if one of the
theories
(the one named the same way in Definition \ref{Definition Log-Cond-Ref-Size}
(page \pageref{Definition Log-Cond-Ref-Size}) )
is equivalent to a formula, we do not need $( \xbm dp).$

{\small

% \begin{tabular*}{10.0cm}{|c|c|c|c|}
\begin{tabular}{|c|c|c|c|}

\hline

\multicolumn{4}{|c|}{Basics} \xEP

\hline

(1.1) \xEH $(OR)$ \xEH $\xch$ \xEH $(\xbm OR)$ \xEP

\cline{1-1}

\cline{3-3}

(1.2) \xEH \xEH $\xci$ \xEH \xEP

\hline

(2.1) \xEH $(disjOR)$ \xEH $\xch$ \xEH $(\xbm disjOR)$ \xEP

\cline{1-1}

\cline{3-3}

(2.2) \xEH \xEH $\xci$ \xEH \xEP

\hline

(3.1) \xEH $(wOR)$ \xEH $\xch$ \xEH $(\xbm wOR)$ \xEP

\cline{1-1}

\cline{3-3}

(3.2) \xEH \xEH $\xci$ \xEH \xEP

\hline

(4.1) \xEH $(SC)$ \xEH $\xch$ \xEH $(\xbm \xcc)$ \xEP

\cline{1-1}

\cline{3-3}

(4.2) \xEH \xEH $\xci$ \xEH \xEP

\hline

(5.1) \xEH $(CP)$ \xEH $\xch$ \xEH $(\xbm \xCQ)$ \xEP

\cline{1-1}

\cline{3-3}

(5.2) \xEH \xEH $\xci$ \xEH \xEP

\hline

(6.1) \xEH $(PR)$ \xEH $\xch$ \xEH $(\xbm PR)$ \xEP

\cline{1-1}

\cline{3-3}

(6.2) \xEH \xEH $\xci$ $(\xbm dp)+(\xbm \xcc)$ \xEH \xEP

\cline{1-1}

\cline{3-3}

(6.3) \xEH \xEH $\xcI$ without $(\xbm dp)$ \xEH \xEP

\cline{1-1}

\cline{3-3}

(6.4) \xEH \xEH $\xci$ $(\xbm \xcc)$ \xEH \xEP

\xEH \xEH $T'$ a formula \xEH \xEP

\hline

(6.5) \xEH $(PR)$ \xEH $\xci$ \xEH $(\xbm PR')$ \xEP

\xEH \xEH $T'$ a formula \xEH \xEP

\hline

(7.1) \xEH $(CUT)$ \xEH $\xch$ \xEH $(\xbm CUT)$ \xEP

\cline{1-1}

\cline{3-3}

(7.2) \xEH \xEH $\xci$ \xEH \xEP

\hline

\multicolumn{4}{|c|}{Cumulativity} \xEP

\hline

(8.1) \xEH $(CM)$ \xEH $\xch$ \xEH $(\xbm CM)$ \xEP

\cline{1-1}

\cline{3-3}

(8.2) \xEH \xEH $\xci$ \xEH \xEP

\hline

(9.1) \xEH $(ResM)$ \xEH $\xch$ \xEH $(\xbm ResM)$ \xEP

\cline{1-1}

\cline{3-3}

(9.2) \xEH \xEH $\xci$ \xEH \xEP

\hline

(10.1) \xEH $(\xcc \xcd)$ \xEH $\xch$ \xEH $(\xbm \xcc \xcd)$ \xEP

\cline{1-1}

\cline{3-3}

(10.2) \xEH \xEH $\xci$ \xEH \xEP

\hline

(11.1) \xEH $(CUM)$ \xEH $\xch$ \xEH $(\xbm CUM)$ \xEP

\cline{1-1}

\cline{3-3}

(11.2) \xEH \xEH $\xci$ \xEH \xEP

\hline

\multicolumn{4}{|c|}{Rationality} \xEP

\hline

(12.1) \xEH $(RatM)$ \xEH $\xch$ \xEH $(\xbm RatM)$ \xEP

\cline{1-1}

\cline{3-3}

(12.2) \xEH \xEH $\xci$ $(\xbm dp)$ \xEH \xEP

\cline{1-1}

\cline{3-3}

(12.3) \xEH \xEH $\xcI$ without $(\xbm dp)$ \xEH \xEP

\cline{1-1}

\cline{3-3}

(12.4) \xEH \xEH $\xci$ \xEH \xEP

\xEH \xEH $T$ a formula \xEH \xEP

\hline

(13.1) \xEH $(RatM=)$ \xEH $\xch$ \xEH $(\xbm =)$ \xEP

\cline{1-1}

\cline{3-3}

(13.2) \xEH \xEH $\xci$ $(\xbm dp)$ \xEH \xEP

\cline{1-1}

\cline{3-3}

(13.3) \xEH \xEH $\xcI$ without $(\xbm dp)$ \xEH \xEP

\cline{1-1}

\cline{3-3}

(13.4) \xEH \xEH $\xci$ \xEH \xEP

\xEH \xEH $T$ a formula \xEH \xEP

\hline

(14.1) \xEH $(Log = ')$ \xEH $\xch$ \xEH $(\xbm = ')$ \xEP

\cline{1-1}
\cline{3-3}

(14.2) \xEH \xEH $\xci$ $(\xbm dp)$ \xEH \xEP

\cline{1-1}
\cline{3-3}

(14.3) \xEH \xEH $\xcI$ without $(\xbm dp)$ \xEH \xEP

\cline{1-1}
\cline{3-3}

(14.4) \xEH \xEH $\xci$ $T$ a formula \xEH \xEP

\hline

(15.1) \xEH $(Log \xFO )$ \xEH $\xch$ \xEH $(\xbm \xFO )$ \xEP

\cline{1-1}
\cline{3-3}

(15.2) \xEH \xEH $\xci$ \xEH \xEP

\hline

(16.1)
\xEH
$(Log \xcv )$
\xEH
$\xch$ $(\xbm \xcc)+(\xbm =)$
\xEH
$(\xbm \xcv )$
\xEP

\cline{1-1}
\cline{3-3}

(16.2) \xEH \xEH $\xci$ $(\xbm dp)$ \xEH \xEP

\cline{1-1}
\cline{3-3}

(16.3) \xEH \xEH $\xcI$ without $(\xbm dp)$ \xEH \xEP

\hline

(17.1)
\xEH
$(Log \xcv ')$
\xEH
$\xch$ $(\xbm \xcc)+(\xbm =)$
\xEH
$(\xbm \xcv ')$
\xEP

\cline{1-1}
\cline{3-3}

(17.2) \xEH \xEH $\xci$ $(\xbm dp)$ \xEH \xEP

\cline{1-1}
\cline{3-3}

(17.3) \xEH \xEH $\xcI$ without $(\xbm dp)$ \xEH \xEP

\hline

\end{tabular}

}

\index{Proposition Alg-Log Proof}

\ep

\subparagraph{
Proof
}

$\hspace{0.01em}$

% (+++ Orig.:  Proof +++)

Set $f(T):=f(M(T)),$ note that $f(T \xcv T' ):=f(M(T \xcv T' ))=f(M(T)
\xcs M(T' )).$

We show first the general framework.

Let $ \xcn $ satisfy $ \xCf (LLE)$ and $ \xCf (CCL).$ Let $f: \xdD_{ \xdl
} \xcp \xdD_{ \xdl }$ be defined by $f(M(T)):=M( \ol{ \ol{T} }).$
If $M(T)=M(T' ),$ then $ \ol{T}= \ol{T' },$ so by $ \xCf (LLE)$ $ \ol{
\ol{T} }= \ol{ \ol{T' } },$ so $f(M(T))=f(M(T' )),$ so $f$ is
well defined and satisfies $( \xbm dp).$ By $ \xCf (CCL)$ $Th(M( \ol{
\ol{T} }))= \ol{ \ol{T} }.$

Let $f$ be given, and $ \xcn $ be defined by $ \ol{ \ol{T}
}:=Th(f(M(T))).$ Obviously, $ \xcn $ satisfies
$ \xCf (LLE)$ and $ \xCf (CCL)$ (and thus $ \xCf (RW)).$ If $f$ satisfies
$( \xbm dp),$ then $f(M(T))=M(T' )$ for
some $T',$ and $f(M(T))=M(Th(f(M(T))))=M( \ol{ \ol{T} })$ by Fact \ref{Fact
Dp-Base} (page \pageref{Fact Dp-Base}). (We will use
Fact \ref{Fact Dp-Base} (page \pageref{Fact Dp-Base})  now without further
mentioning.)

Next we show the following fact:

(a) If $f$ satisfies $( \xbm dp),$ or $T' $ is equivalent to a formula,
then
$Th(f(T) \xcs M(T' ))= \ol{ \ol{ \ol{T} } \xcv T' }.$

Case 1, $f$ satisfies $( \xbm dp).$ $Th(f(M(T)) \xcs M(T' ))$ $=$ $Th(M(
\ol{ \ol{T} }) \xcs M(T' )$ $=$ $ \ol{ \ol{ \ol{T} } \xcv T' }$
by Fact \ref{Fact Log-Form} (page \pageref{Fact Log-Form})  (5).

Case 2, $T' $ is equivalent to $ \xbf '.$ $Th(f(M(T)) \xcs M( \xbf ' ))$
$=$ $ \ol{Th(f(M(T))) \xcv \{ \xbf ' \}}$ $=$
$ \ol{ \ol{ \ol{T} } \xcv \{ \xbf ' \}}$ by Fact \ref{Fact Log-Form} (page
\pageref{Fact Log-Form})
(3).

We now prove the individual properties.

(1.1) $ \xCf (OR)$ $ \xch $ $( \xbm OR)$

Let $X=M(T),$ $Y=M(T' ).$ $f(X \xcv Y)$ $=$ $f(M(T) \xcv M(T' ))$ $=$
$f(M(T \xco T' ))$ $:=$ $M( \ol{ \ol{T \xco T' } })$ $ \xcc_{(OR)}$
$M( \ol{ \ol{T} } \xcs \ol{ \ol{T' } })$ $=_{(CCL)}$ $M( \ol{ \ol{T} })
\xcv M( \ol{ \ol{T' } })$ $=:$ $f(X) \xcv f(Y).$

(1.2) $( \xbm OR)$ $ \xch $ $ \xCf (OR)$

$ \ol{ \ol{T \xco T' } }$ $:=$ $Th(f(M(T \xco T' )))$ $=$ $Th(f(M(T) \xcv
M(T' )))$ $ \xcd_{( \xbm OR)}$ $Th(f(M(T)) \xcv f(M(T' )))$ $=$
(by Fact \ref{Fact Th-Union} (page \pageref{Fact Th-Union}) ) $Th(f(M(T))) \xcs
Th(f(M(T' )))$ $=:$ $
\ol{ \ol{T} } \xcs \ol{ \ol{T' } }.$

(2) By $ \xCN Con(T,T' ) \xcj M(T) \xcs M(T' )= \xCQ,$ we can use
directly the proofs for 1.

(3.1) $ \xCf (wOR)$ $ \xch $ $( \xbm wOR)$

Let $X=M(T),$ $Y=M(T' ).$ $f(X \xcv Y)$ $=$ $f(M(T) \xcv M(T' ))$ $=$
$f(M(T \xco T' ))$ $:=$ $M( \ol{ \ol{T \xco T' } })$ $ \xcc_{(wOR)}$
$M( \ol{ \ol{T} } \xcs \ol{T' })$ $=_{(CCL)}$ $M( \ol{ \ol{T} }) \xcv M(
\ol{T' })$ $=:$ $f(X) \xcv Y.$

(3.2) $( \xbm wOR)$ $ \xch $ $ \xCf (wOR)$

$ \ol{ \ol{T \xco T' } }$ $:=$ $Th(f(M(T \xco T' )))$ $=$ $Th(f(M(T) \xcv
M(T' )))$ $ \xcd_{( \xbm wOR)}$ $Th(f(M(T)) \xcv M(T' ))$ $=$
(by Fact \ref{Fact Th-Union} (page \pageref{Fact Th-Union}) ) $Th(f(M(T))) \xcs
Th(M(T' ))$ $=:$ $
\ol{ \ol{T} } \xcs \ol{T' }.$

(4.1) $ \xCf (SC)$ $ \xch $ $( \xbm \xcc )$

Trivial.

(4.2) $( \xbm \xcc )$ $ \xch $ $ \xCf (SC)$

Trivial.

(5.1) $ \xCf (CP)$ $ \xch $ $( \xbm \xCQ )$

Trivial.

(5.2) $( \xbm \xCQ )$ $ \xch $ $ \xCf (CP)$

Trivial.

(6.1) $ \xCf (PR)$ $ \xch $ $( \xbm PR):$

Suppose $X:=M(T),$ $Y:=M(T' ),$ $X \xcc Y,$ we have to show $f(Y) \xcs X
\xcc f(X).$
By prerequisite, $ \ol{T' } \xcc \ol{T},$ so $ \ol{T \xcv T' }= \ol{T},$
so $ \ol{ \ol{T \xcv T' } }= \ol{ \ol{T} }$ by $ \xCf (LLE).$ By $ \xCf
(PR)$ $ \ol{ \ol{T \xcv T' } } \xcc \ol{ \ol{ \ol{T'
} } \xcv T},$
so $f(Y) \xcs X=f(T' ) \xcs M(T)=M( \ol{ \ol{T' } } \xcv T) \xcc M( \ol{
\ol{T \xcv T' } })=M( \ol{ \ol{T} })=f(X).$

(6.2) $( \xbm PR)+( \xbm dp)+( \xbm \xcc )$ $ \xch $ $ \xCf (PR):$

$f(T) \xcs M(T' )=_{( \xbm \xcc )}f(T) \xcs M(T) \xcs M(T' )=f(T) \xcs M(T
\xcv T' ) \xcc_{( \xbm PR)}f(T \xcv T' ),$ so
$ \ol{ \ol{T \xcv T' } }=Th(f(T \xcv T' )) \xcc Th(f(T) \xcs M(T' ))= \ol{
\ol{ \ol{T} } \xcv T' }$ by (a) above and $( \xbm dp).$

(6.3) $( \xbm PR)$ $ \xcH $ $ \xCf (PR)$ without $( \xbm dp):$

$( \xbm PR)$ holds in all preferential structures
(see Definition \ref{Definition Pref-Str} (page \pageref{Definition Pref-Str}) )
by Fact \ref{Fact Pref-Sound} (page \pageref{Fact Pref-Sound}).
Example \ref{Example Pref-Dp} (page \pageref{Example Pref-Dp})  shows that $
\xCf (DP)$ may fail in the
resulting logic.

(6.4) $( \xbm PR)+( \xbm \xcc )$ $ \xch $ $ \xCf (PR)$ if $T' $ is
classically equivalent to a formula:

It was shown in the proof of (6.2) that $f(T) \xcs M( \xbf ' ) \xcc f(T
\xcv \{ \xbf ' \}),$ so
$ \ol{ \ol{T \xcv \{ \xbf ' \}} }=Th(f(T \xcv \{ \xbf ' \})) \xcc Th(f(T)
\xcs M( \xbf ' ))= \ol{ \ol{ \ol{T} } \xcv \{ \xbf ' \}}$ by (a) above.

(6.5) $( \xbm PR' )$ $ \xch $ $ \xCf (PR),$ if $T' $ is classically
equivalent to a formula:

$f(M(T)) \xcs M( \xbf ' )$ $ \xcc_{( \xbm PR' )}$ $f(M(T) \xcs M( \xbf '
))$ $=$ $f(M(T \xcv \{ \xbf ' \})).$ So again
$ \ol{ \ol{T \xcv \{ \xbf ' \}} }=Th(f(T \xcv \{ \xbf ' \})) \xcc Th(f(T)
\xcs M( \xbf ' ))= \ol{ \ol{ \ol{T} } \xcv \{ \xbf ' \}}$ by (a) above.

(7.1) $ \xCf (CUT)$ $ \xch $ $( \xbm CUT)$

So let $X=M(T),$ $Y=M(T' ),$ and $f(T):=M( \ol{ \ol{T} }) \xcc M(T' ) \xcc
M(T)$ $ \xch $ $ \ol{T} \xcc \ol{T' } \xcc \ol{ \ol{T} }=_{ \xCf (LLE)}
\ol{ \ol{( \ol{T})} }$ $ \xch $
(by $ \xCf (CUT))$
$ \ol{ \ol{T} }= \ol{ \ol{( \ol{T})} } \xcd \ol{ \ol{( \ol{T' })} }= \ol{
\ol{T' } }$ $ \xch $ $f(T)=M( \ol{ \ol{T} }) \xcc M( \ol{ \ol{T' } })=f(T'
),$ $thus$ $f(X) \xcc f(Y).$

(7.2) $( \xbm CUT)$ $ \xch $ $ \xCf (CUT)$

Let $T$ $ \xcc $ $ \ol{T' }$ $ \xcc $ $ \ol{ \ol{T} }.$ Thus $f(T) \xcc M(
\ol{ \ol{T} })$ $ \xcc $ $M(T' )$ $ \xcc $ $M(T),$ so by $( \xbm CUT)$
$f(T) \xcc f(T' ),$
so $ \ol{ \ol{T} }$ $=$ $Th(f(T))$ $ \xcd $ $Th(f(T' ))$ $=$ $ \ol{ \ol{T'
} }.$

(8.1) $ \xCf (CM)$ $ \xch $ $( \xbm CM)$

So let $X=M(T),$ $Y=M(T' ),$ and $f(T):=M( \ol{ \ol{T} }) \xcc M(T' ) \xcc
M(T)$ $ \xch $ $ \ol{T} \xcc \ol{T' } \xcc \ol{ \ol{T} }=_{ \xCf (LLE)}
\ol{ \ol{( \ol{T})} }$ $ \xch $
(by $ \xCf (LLE),$ $ \xCf (CM))$
$ \ol{ \ol{T} }= \ol{ \ol{( \ol{T})} } \xcc \ol{ \ol{( \ol{T' })} }= \ol{
\ol{T' } }$ $ \xch $ $f(T)=M( \ol{ \ol{T} }) \xcd M( \ol{ \ol{T' } })=f(T'
),$ $thus$ $f(X) \xcd f(Y).$

(8.2) $( \xbm CM)$ $ \xch $ $ \xCf (CM)$

Let $T$ $ \xcc $ $ \ol{T' }$ $ \xcc $ $ \ol{ \ol{T} }.$ Thus by $( \xbm
CM)$ and $f(T) \xcc M( \ol{ \ol{T} })$ $ \xcc $ $M(T' )$ $ \xcc $ $M(T),$
so $f(T) \xcd f(T' )$
by $( \xbm CM),$ so $ \ol{ \ol{T} }$ $=$ $Th(f(T))$ $ \xcc $ $Th(f(T' ))$
$=$ $ \ol{ \ol{T' } }.$

(9.1) $ \xCf (ResM)$ $ \xch $ $( \xbm ResM)$

Let $f(X):=M( \ol{ \ol{ \xbD } }),$ $A:=M( \xba ),$ $B:=M( \xbb ).$ So
$f(X) \xcc A \xcs B$ $ \xch $ $ \xbD \xcn \xba, \xbb $ $ \xch_{ \xCf
(ResM)}$
$ \xbD, \xba \xcn \xbb $ $ \xch $ $M( \ol{ \ol{ \xbD, \xba } }) \xcc M(
\xbb )$ $ \xch $ $f(X \xcs A) \xcc B.$

(9.2) $( \xbm ResM)$ $ \xch $ $ \xCf (ResM)$

Let $f(X):=M( \ol{ \ol{ \xbD } }),$ $A:=M( \xba ),$ $B:=M( \xbb ).$ So $
\xbD \xcn \xba, \xbb $ $ \xch $ $f(X) \xcc A \xcs B$ $ \xch_{( \xbm
ResM)}$
$f(X \xcs A) \xcc B$ $ \xch $ $ \xbD, \xba \xcn \xbb.$

(10.1) $( \xcc \xcd )$ $ \xch $ $( \xbm \xcc \xcd )$

Let $f(T) \xcc M(T' ),$ $f(T' ) \xcc M(T).$
So $Th(M(T' )) \xcc Th(f(T)),$ $Th(M(T)) \xcc Th(f(T' )),$ so $T' \xcc
\ol{T' } \xcc \ol{ \ol{T} },$ $T \xcc \ol{T} \xcc \ol{ \ol{T' } },$
so by $( \xcc \xcd )$ $ \ol{ \ol{T} }= \ol{ \ol{T' } },$ so $f(T):=M( \ol{
\ol{T} })=M( \ol{ \ol{T' } })=:f(T' ).$

(10.2) $( \xbm \xcc \xcd )$ $ \xch $ $( \xcc \xcd )$

Let $T \xcc \ol{ \ol{T' } }$ and $T' \xcc \ol{ \ol{T} }.$ So by $ \xCf
(CCL)$ $Th(M(T))= \ol{T} \xcc \ol{ \ol{T' } }=Th(f(T' )).$
But $Th(M(T)) \xcc Th(X) \xch X \xcc M(T):$ $X \xcc M(Th(X)) \xcc
M(Th(M(T)))=M(T).$
So $f(T' ) \xcc M(T),$ likewise $f(T) \xcc M(T' ),$ so by $( \xbm \xcc
\xcd )$ $f(T)=f(T' ),$ so $ \ol{ \ol{T} }= \ol{ \ol{T' } }.$

(11.1) $ \xCf (CUM)$ $ \xch $ $( \xbm CUM):$

So let $X=M(T),$ $Y=M(T' ),$ and $f(T):=M( \ol{ \ol{T} }) \xcc M(T' ) \xcc
M(T)$ $ \xch $ $ \ol{T} \xcc \ol{T' } \xcc \ol{ \ol{T} }=_{ \xCf (LLE)}
\ol{ \ol{( \ol{T})} }$ $ \xch $
$ \ol{ \ol{T} }= \ol{ \ol{( \ol{T})} }= \ol{ \ol{( \ol{T' })} }= \ol{
\ol{T' } }$ $ \xch $ $f(T)=M( \ol{ \ol{T} })=M( \ol{ \ol{T' } })=f(T' ),$
$thus$ $f(X)=f(Y).$

(11.2) $( \xbm CUM)$ $ \xch $ $ \xCf (CUM)$:

Let $T$ $ \xcc $ $ \ol{T' }$ $ \xcc $ $ \ol{ \ol{T} }.$ Thus by $( \xbm
CUM)$ and $f(T) \xcc M( \ol{ \ol{T} })$ $ \xcc $ $M(T' )$ $ \xcc $ $M(T),$
so $f(T)=f(T' ),$
so $ \ol{ \ol{T} }$ $=$ $Th(f(T))$ $=$ $Th(f(T' ))$ $=$ $ \ol{ \ol{T' }
}.$

(12.1) $ \xCf (RatM)$ $ \xch $ $( \xbm RatM)$

Let $X=M(T),$ $Y=M(T' ),$ and $X \xcc Y,$ $X \xcs f(Y) \xEd \xCQ,$ so $T
\xcl T' $ and $M(T) \xcs f(M(T' )) \xEd \xCQ,$
so $Con(T, \ol{ \ol{T' } }),$ so $ \ol{ \ol{ \ol{T' } } \xcv T} \xcc \ol{
\ol{T} }$ by $ \xCf (RatM),$ so $f(X)=f(M(T))=M( \ol{ \ol{T} }) \xcc M(
\ol{ \ol{T' } } \xcv T)=$
$M( \ol{ \ol{T' } }) \xcs M(T)=f(Y) \xcs X.$

(12.2) $( \xbm RatM)+( \xbm dp)$ $ \xch $ $ \xCf (RatM):$

Let $X=M(T),$ $Y=M(T' ),$ $T \xcl T',$ $Con(T, \ol{ \ol{T' } }),$ so $X
\xcc Y$ and by $( \xbm dp)$ $X \xcs f(Y) \xEd \xCQ,$
so by $( \xbm RatM)$ $f(X) \xcc f(Y) \xcs X,$ so
$ \ol{ \ol{T} }= \ol{ \ol{T \xcv T' } }=Th(f(T \xcv T' )) \xcd Th(f(T' )
\xcs M(T))= \ol{ \ol{ \ol{T' } } \xcv T}$ by (a) above and $( \xbm dp).$

(12.3) $( \xbm RatM)$ $ \xcH $ $ \xCf (RatM)$ without $( \xbm dp):$

$( \xbm RatM)$ holds in all ranked preferential structures
(see Definition \ref{Definition Rank-Rel} (page \pageref{Definition Rank-Rel}) )
by Fact \ref{Fact Rank-Hold} (page \pageref{Fact Rank-Hold}). Example
\ref{Example Rank-Dp} (page \pageref{Example Rank-Dp})  (2)
shows
that $ \xCf (RatM)$ may fail in the resulting logic.

(12.4) $( \xbm RatM)$ $ \xch $ $ \xCf (RatM)$ if $T$ is classically
equivalent to a formula:

$ \xbf \xcl T' $ $ \xch $ $M( \xbf ) \xcc M(T' ).$ $Con( \xbf, \ol{
\ol{T' } })$ $ \xcj $ $M( \ol{ \ol{T' } }) \xcs M( \xbf ) \xEd \xCQ $ $
\xcj $ $f(T' ) \xcs M( \xbf ) \xEd \xCQ $
by Fact \ref{Fact Log-Form} (page \pageref{Fact Log-Form})  (4). Thus $f(M( \xbf
)) \xcc f(M(T' ))
\xcs M( \xbf )$ by $( \xbm RatM).$
Thus by (a) above $ \ol{ \ol{ \ol{T' } } \xcv \{ \xbf \}} \xcc \ol{ \ol{
\xbf } }.$

(13.1) $ \xCf (RatM=)$ $ \xch $ $( \xbm =)$

Let $X=M(T),$ $Y=M(T' ),$ and $X \xcc Y,$ $X \xcs f(Y) \xEd \xCQ,$ so $T
\xcl T' $ and $M(T) \xcs f(M(T' )) \xEd \xCQ,$
so $Con(T, \ol{ \ol{T' } }),$ so $ \ol{ \ol{ \ol{T' } } \xcv T}= \ol{
\ol{T} }$ by $ \xCf (RatM=),$ so $f(X)=f(M(T))=M( \ol{ \ol{T} })=M( \ol{
\ol{T' } } \xcv T)=$
$M( \ol{ \ol{T' } }) \xcs M(T)=f(Y) \xcs X.$

(13.2) $( \xbm =)+( \xbm dp)$ $ \xch $ $ \xCf (RatM=)$

Let $X=M(T),$ $Y=M(T' ),$ $T \xcl T',$ $Con(T, \ol{ \ol{T' } }),$ so $X
\xcc Y$ and by $( \xbm dp)$ $X \xcs f(Y) \xEd \xCQ,$
so by $( \xbm =)$ $f(X)=f(Y) \xcs X.$ So $ \ol{ \ol{ \ol{T' } } \xcv T}=
\ol{ \ol{T} }$ (a) above and $( \xbm dp).$

(13.3) $( \xbm =)$ $ \xcH $ $ \xCf (RatM=)$ without $( \xbm dp):$

$( \xbm =)$ holds in all ranked preferential structures
(see Definition \ref{Definition Rank-Rel} (page \pageref{Definition Rank-Rel}) )
by Fact \ref{Fact Rank-Hold} (page \pageref{Fact Rank-Hold}). Example
\ref{Example Rank-Dp} (page \pageref{Example Rank-Dp})  (1)
shows
that $ \xCf (RatM=)$ may fail in the resulting logic.

(13.4) $( \xbm =)$ $ \xch $ $ \xCf (RatM=)$ if $T$ is classically
equivalent to a formula:

The proof is almost identical to the one for (12.4). Again, the
prerequisites of $( \xbm =)$ are satisfied, so $f(M( \xbf ))=f(M(T' ))
\xcs M( \xbf ).$
Thus, $ \ol{ \ol{ \ol{T' } } \xcv \{ \xbf \}}= \ol{ \ol{ \xbf } }$ by (a)
above.

Of the last four, we show (14), (15), (17), the proof for (16) is similar
to the
one for (17).

(14.1) $(Log=' )$ $ \xch $ $( \xbm =' ):$

$f(M(T' )) \xcs M(T) \xEd \xCQ $ $ \xch $ $Con( \ol{ \ol{T' } } \xcv T)$ $
\xch_{(Log=' )}$ $ \ol{ \ol{T \xcv T' } }= \ol{ \ol{ \ol{T' } } \xcv T}$ $
\xch $
$f(M(T \xcv T' ))=f(M(T' )) \xcs M(T).$

(14.2) $( \xbm =' )+( \xbm dp)$ $ \xch $ $(Log=' ):$

$Con( \ol{ \ol{T' } } \xcv T)$ $ \xch_{( \xbm dp)}$ $f(M(T' )) \xcs M(T)
\xEd \xCQ $ $ \xch $ $f(M(T' \xcv T))=f(M(T' ) \xcs M(T))$ $=_{( \xbm ='
)}$
$f(M(T' )) \xcs M(T),$ so $ \ol{ \ol{T' \xcv T} }$ $=$ $ \ol{ \ol{ \ol{T'
} } \xcv T}$ by (a) above and $( \xbm dp).$

(14.3) $( \xbm =' )$ $ \xcH $ $(Log=' )$ without $( \xbm dp):$

By Fact \ref{Fact Rank-Hold} (page \pageref{Fact Rank-Hold})  $( \xbm =' )$
holds in ranked
structures.
Consider Example \ref{Example Rank-Dp} (page \pageref{Example Rank-Dp})  (2).
There, $Con(T, \ol{ \ol{T' }
}),$ $T=T \xcv T',$ and
it was shown that $ \ol{ \ol{ \ol{T' } } \xcv T}$ $ \xcC $ $ \ol{ \ol{T}
}$ $=$ $ \ol{ \ol{T \xcv T' } }$

(14.4) $( \xbm =' )$ $ \xch $ $(Log=' )$ if $T$ is classically equivalent
to a formula:

$Con( \ol{ \ol{T' } } \xcv \{ \xbf \})$ $ \xch $ $ \xCQ \xEd M( \ol{
\ol{T' } }) \xcs M( \xbf )$ $ \xch $ $f(T' ) \xcs M( \xbf ) \xEd \xCQ $ by
Fact \ref{Fact Log-Form} (page \pageref{Fact Log-Form})  (4). So $f(M(T' \xcv \{
\xbf \}))$ $=$
$f(M(T' ) \xcs M( \xbf ))$ $=$
$f(M(T' )) \xcs M( \xbf )$ by $( \xbm =' ),$ so $ \ol{ \ol{T' \xcv \{ \xbf
\}} }= \ol{ \ol{ \ol{T' } } \xcv \{ \xbf \}}$ by (a) above.

(15.1) $(Log \xFO )$ $ \xch $ $( \xbm \xFO ):$

Trivial.

(15.2) $( \xbm \xFO )$ $ \xch $ $(Log \xFO ):$

Trivial.

(16) $(Log \xcv )$ $ \xcj $ $( \xbm \xcv ):$ Analogous to the proof of
(17).

(17.1) $(Log \xcv ' )+( \xbm \xcc )+( \xbm =)$ $ \xch $ $( \xbm \xcv ' ):$

$f(M(T' )) \xcs (M(T)-f(M(T))) \xEd \xCQ $ $ \xch $ (by $( \xbm \xcc ),$
$( \xbm =),$
Fact \ref{Fact Rank-Auxil} (page \pageref{Fact Rank-Auxil}) ) $f(M(T' )) \xcs
M(T) \xEd \xCQ,$ $f(M(T'
)) \xcs f(M(T))= \xCQ $ $ \xch $
$Con( \ol{ \ol{T' } },T),$ $ \xCN Con( \ol{ \ol{T' } }, \ol{ \ol{T} })$ $
\xch $ $ \ol{ \ol{T \xco T' } }= \ol{ \ol{T} }$ $ \xch $ $f(M(T))=f(M(T
\xco T' ))=f(M(T) \xcv M(T' )).$

(17.2) $( \xbm \xcv ' )+( \xbm dp)$ $ \xch $ $(Log \xcv ' ):$

$Con( \ol{ \ol{T' } } \xcv T),$ $ \xCN Con( \ol{ \ol{T' } } \xcv \ol{
\ol{T} })$ $ \xch_{( \xbm dp)}$ $f(T' ) \xcs M(T) \xEd \xCQ,$ $f(T' )
\xcs f(T)= \xCQ $ $ \xch $
$f(M(T' )) \xcs (M(T)-f(M(T))) \xEd \xCQ $ $ \xch $ $f(M(T))$ $=$ $f(M(T)
\xcv M(T' ))$ $=$ $f(M(T \xco T' )).$
So $ \ol{ \ol{T} }= \ol{ \ol{T \xco T' } }.$

(17.3) and (16.3) are solved by Example \ref{Example Rank-Dp} (page
\pageref{Example Rank-Dp})  (3).

$ \xcz $
\\[3ex]
\index{Example Rank-Dp}

\be

$\hspace{0.01em}$

% (+++ Orig. No.:  Example Rank-Dp +++)

\label{Example Rank-Dp}

(1) $( \xbm =)$ without $( \xbm dp)$ does not imply $(RatM=):$

Take $\{p_{i}:i \xbe \xbo \}$ and put $m:=m_{ \xcU p_{i}},$ the model
which makes all $p_{i}$ true, in the top
layer, all the other in the bottom layer. Let $m' \xEd m,$ $T':= \xCQ,$
$T:=Th(m,m' ).$ Then
Then $ \ol{ \ol{T' } }=T',$ so $Con( \ol{ \ol{T' } },T),$ $ \ol{ \ol{T}
}=Th(m' ),$ $ \ol{ \ol{ \ol{T' } } \xcv T}=T.$

So $(RatM=)$ fails, but $( \xbm =)$ holds in all ranked structures.

(2) $( \xbm RatM)$ without $( \xbm dp)$ does not imply (RatM):

Take $\{p_{i}:i \xbe \xbo \}$ and let $m:=m_{ \xcU p_{i}},$ the model
which makes all $p_{i}$ true.

Let $X:=M( \xCN p_{0}) \xcv \{m\}$ be the top layer, put the rest of $M_{
\xdl }$ in the bottom layer.
Let $Y:=M_{ \xdl }.$ The structure is ranked, as shown in Fact \ref{Fact
Rank-Hold} (page \pageref{Fact Rank-Hold}),
$( \xbm RatM)$ holds.

Let $T':= \xCQ,$ $T:=Th(X).$ We have to show that $Con(T, \ol{ \ol{T' }
}),$ $T \xcl T',$ but
$ \ol{ \ol{ \ol{T' } } \xcv T} \xcC \ol{ \ol{T} }.$ $ \ol{ \ol{T' } }$ $=$
$Th(M(p_{0})-\{m\})$ $=$ $ \ol{p_{0}}.$ $T$ $=$ $ \ol{\{ \xCN p_{0}\} \xco
Th(m)},$ $ \ol{ \ol{T} }=T.$ $So$ $Con(T, \ol{
\ol{T' } }).$
$M( \ol{ \ol{T' } })=M(p_{0}),$ $M(T)=X,$ $M( \ol{ \ol{T' } } \xcv T)=M(
\ol{ \ol{T' } }) \xcs M(T)=\{m\},$ $m \xcm p_{1},$ so $p_{1} \xbe \ol{
\ol{ \ol{T' } } \xcv T},$ but
$X \xcM p_{1}.$

(3) This example shows that we need $( \xbm dp)$ to go from $( \xbm \xcv
)$ to
$(Log \xcv )$ and from $( \xbm \xcv ' )$ to $(Log \xcv ' ).$

Let $v( \xdl ):=\{p,q\} \xcv \{p_{i}:i< \xbo \}.$ Let $m$ make all
variables true.

Put all models of $ \xCN p,$ and $m,$ in the upper layer, all other models
in the
lower layer. This is ranked, so by Fact \ref{Fact Rank-Hold} (page \pageref{Fact
Rank-Hold})  $( \xbm
\xcv )$
and $( \xbm \xcv ' )$ hold.
Set $X:=M( \xCN q) \xcv \{m\},$ $X':=M(q)-\{m\},$ $T:=Th(X)= \xCN q \xco
Th(m),$ $T':=Th(X' )= \ol{q}.$
Then $ \ol{ \ol{T} }= \ol{p \xcu \xCN q},$ $ \ol{ \ol{T' } }= \ol{p \xcu
q}.$ We have $Con( \ol{ \ol{T' } },T),$ $ \xCN Con( \ol{ \ol{T' } }, \ol{
\ol{T} }).$
But $ \ol{ \ol{T \xco T' } }= \ol{p} \xEd \ol{ \ol{T} }= \ol{p \xcu \xCN
q}$ and $Con( \ol{ \ol{T \xco T' } },T' ),$ so $(Log \xcv )$ and $(Log
\xcv ' )$ fail.

$ \xcz $
\\[3ex]
\index{Fact Cut-Pr}

\ee

\bfa

$\hspace{0.01em}$

% (+++ Orig. No.:  Fact Cut-Pr +++)

\label{Fact Cut-Pr}

$ \xCf (CUT)$ $ \xcH $ $ \xCf (PR)$
\index{Fact Cut-Pr Proof}

\efa

\subparagraph{
Proof
}

$\hspace{0.01em}$

% (+++ Orig.:  Proof +++)

We give two proofs:

(1) If $ \xCf (CUT)$ $ \xch $ $ \xCf (PR),$ then by $( \xbm PR)$ $ \xch $
(by Fact \ref{Fact Mu-Base} (page \pageref{Fact Mu-Base})  (3))
$( \xbm CUT)$ $ \xch $ (by Proposition \ref{Proposition Alg-Log} (page
\pageref{Proposition Alg-Log})  (7.2) $
\xCf (CUT)$ $ \xch $ $ \xCf (PR)$
we would have a proof of $( \xbm PR)$ $ \xch $ $ \xCf (PR)$ without $(
\xbm dp),$ which is impossible,
as shown by Example \ref{Example Pref-Dp} (page \pageref{Example Pref-Dp}).

(2) Reconsider Example \ref{Example Need-Pr} (page \pageref{Example Need-Pr}),
and say $a \xcm p \xcu q,$
$b \xcm p \xcu \xCN q,$ $c \xcm \xCN p \xcu q.$
It is shown there that $( \xbm CUM)$ holds, so $( \xbm CUT)$ holds, so by
Proposition \ref{Proposition Alg-Log} (page \pageref{Proposition Alg-Log}) 
(7.2) $ \xCf (CUT)$ holds, if we
define $ \ol{ \ol{T} }:=Th(f(M(T)).$
Set $T:=\{p \xco ( \xCN p \xcu q)\},$ $T':=\{p\},$ then
$ \ol{ \ol{T \xcv T' } }= \ol{ \ol{T' } }= \ol{\{p \xcu \xCN q\}},$ $ \ol{
\ol{T} }= \ol{T},$ $ \ol{T \xcv T' }= \ol{T' }= \ol{\{p\}},$ $so$ $ \xCf
(PR)$ $fails.$

$ \xcz $
\\[3ex]
\chapter{
Abstract semantics by size
}
\label{Chapter Size}
\section{
The first order setting
}

We first introduce a generalized quantifier in a first order setting,
as this is very natural, and prepares the more abstract discussion
to come.
\index{Definition Nabla}

\bd

$\hspace{0.01em}$

% (+++ Orig. No.:  Definition Nabla +++)

\label{Definition Nabla}

Augment the language of first order logic by the new quantifier:
If $ \xbf $ and $ \xbq $ are formulas, then so are $ \xeA x \xbf (x),$ $
\xeA x \xbf (x): \xbq (x),$
for any variable $x.$ The
:-versions are the restricted variants.
We call any formula of $ \xdl,$ possibly containing $ \xeA $ a $ \xeA -
\xdl -$formula.
\index{Definition N-Model}

\ed

\bd

$\hspace{0.01em}$

% (+++ Orig. No.:  Definition N-Model +++)

\label{Definition N-Model}

$( \xdn -$Model)

Let $ \xdl $ be a first order language, and $M$ be a $ \xdl -$structure.
Let $ \xdn (M)$ be
a weak filter, or $ \xdn -$system - $ \xdn $ for normal - over $M.$
Define $ \xBc M, \xdn (M) \xBe $ $ \xcm $ $ \xbf $ for any $ \xeA - \xdl
-$formula
inductively as usual, with
one additional induction step:

$ \xBc M, \xdn (M) \xBe $ $ \xcm $ $ \xeA x \xbf (x)$ iff there is $A \xbe \xdn
(M)$
s.t. $ \xcA a \xbe A$ $( \xBc M, \xdn (M) \xBe $ $ \xcm $ $ \xbf [a]).$
\index{Definition NablaAxioms}

\ed

\bd

$\hspace{0.01em}$

% (+++ Orig. No.:  Definition NablaAxioms +++)

\label{Definition NablaAxioms}

Let any axiomatization of predicate calculus be given. Augment this with
the axiom schemata

(1) $ \xeA x \xbf (x)$ $ \xcu $ $ \xcA x( \xbf (x) \xcp \xbq (x))$ $ \xch
$ $ \xeA x \xbq (x),$

(2) $ \xeA x \xbf (x)$ $ \xch $ $ \xCN \xeA x \xCN \xbf (x),$

(3) $ \xcA x \xbf (x)$ $ \xch $ $ \xeA x \xbf (x)$ and $ \xeA x \xbf (x)$
$ \xch $ $ \xcE x \xbf (x),$

(4) $ \xeA x \xbf (x)$ $ \xcr $ $ \xeA y \xbf (y)$ if $x$ does not occur
free in $ \xbf (y)$ and $y$ does not
occur free in $ \xbf (x).$

(for all $ \xbf,$ $ \xbq )$.
\index{Proposition NablaRepr}

\ed

\bp

$\hspace{0.01em}$

% (+++ Orig. No.:  Proposition NablaRepr +++)

\label{Proposition NablaRepr}

The axioms given in Definition \ref{Definition NablaAxioms} (page
\pageref{Definition NablaAxioms})
are sound and complete for the semantics of Definition \ref{Definition N-Model}
(page \pageref{Definition N-Model})

See  \cite{Sch95-1} or  \cite{Sch04}.
\index{Definition Nabla-System}

\ep

\bd

$\hspace{0.01em}$

% (+++ Orig. No.:  Definition Nabla-System +++)

\label{Definition Nabla-System}

Call $ \xdn^{+}(M)= \xBc  \xdn (N):N \xcc M \xBe $ a $ \xdn^{+}-system$ or
system of
weak filters over $M$ iff
for each $N \xcc M$ $ \xdn (N)$ is a weak filter or $ \xdn -$system over
$N.$
(It suffices to consider the definable subsets of $M.)$
\index{Definition N-Model-System}

\ed

\bd

$\hspace{0.01em}$

% (+++ Orig. No.:  Definition N-Model-System +++)

\label{Definition N-Model-System}

Let $ \xdl $ be a first order language, and $M$ a $ \xdl -$structure. Let
$ \xdn^{+}(M)$ be
a $ \xdn^{+}-system$ over $M.$

Define $ \xBc M, \xdn^{+}(M) \xBe $ $ \xcm $ $ \xbf $ for any formula
inductively as
usual, with
the additional induction steps:

1. $ \xBc M, \xdn^{+}(M) \xBe $ $ \xcm $ $ \xeA x \xbf (x)$ iff there is $A \xbe
\xdn (M)$ s.t. $ \xcA a \xbe A$ $( \xBc M, \xdn^{+}(M) \xBe $ $ \xcm $ $ \xbf
[a]),$

2. $ \xBc M, \xdn^{+}(M) \xBe $ $ \xcm $ $ \xeA x \xbf (x): \xbq (x)$ iff there
is
$A \xbe \xdn (\{x: \xBc M, \xdn^{+}(M) \xBe  \xcm \xbf (x)\})$ s.t.
$ \xcA a \xbe A$ $( \xBc M, \xdn^{+}(M) \xBe $ $ \xcm $ $ \xbq [a]).$
\index{Definition NablaAxioms-System}

\ed

\bd

$\hspace{0.01em}$

% (+++ Orig. No.:  Definition NablaAxioms-System +++)

\label{Definition NablaAxioms-System}

Extend the logic of first order predicate calculus by adding the axiom
schemata

(1) a. $ \xeA x \xbf (x)$ $ \xcj $ $ \xeA x(x=x): \xbf (x),$
$b.$ $ \xcA x( \xbs (x) \xcr \xbt (x))$ $ \xcu $ $ \xeA x \xbs (x): \xbf
(x)$ $ \xch $ $ \xeA x \xbt (x): \xbf (x),$

(2) $ \xeA x \xbf (x): \xbq (x)$ $ \xcu $ $ \xcA x( \xbf (x) \xcu \xbq (x)
\xcp \xbj (x))$ $ \xch $ $ \xeA x \xbf (x): \xbj (x),$

(3) $ \xcE x \xbf (x)$ $ \xcu $ $ \xeA x \xbf (x): \xbq (x)$ $ \xch $ $
\xCN \xeA x \xbf (x): \xCN \xbq (x),$

(4) $ \xcA x( \xbf (x) \xcp \xbq (x))$ $ \xch $ $ \xeA x \xbf (x): \xbq
(x)$
and $ \xeA x \xbf (x): \xbq (x)$ $ \xcp $ $[ \xcE x \xbf (x)$ $ \xcp $ $
\xcE x( \xbf (x) \xcu \xbq (x))],$

(5) $ \xeA x \xbf (x): \xbq (x)$ $ \xcr $ $ \xeA y \xbf (y): \xbq (y)$
(under the usual caveat for substitution).

(for all $ \xbf,$ $ \xbq,$ $ \xbj,$ $ \xbs,$ $ \xbt )$.
\index{Proposition NablaRepr-System}

\ed

\bp

$\hspace{0.01em}$

% (+++ Orig. No.:  Proposition NablaRepr-System +++)

\label{Proposition NablaRepr-System}

The axioms of Definition \ref{Definition NablaAxioms-System} (page
\pageref{Definition NablaAxioms-System})  are
sound and complete for the $ \xdn^{+}-semantics$
of $ \xeA $ as defined in Definition \ref{Definition N-Model-System} (page
\pageref{Definition N-Model-System}).

See  \cite{Sch95-1} or  \cite{Sch04}.
\section{
General size semantics
}
\subsection{
Introduction
}
\subsubsection{
Context
}

\ep

We show how one can develop a multitude of rules for nonmonotonic
logics from a very small
set of principles about reasoning with size. The notion of size gives an
algebraic semantics to nonmonotonic logics, in the sense that $ \xba $
implies $ \xbb $
iff the set of cases where $ \xba \xcu \xCN \xbb $ holds is a small subset
of all $ \xba -$cases.
In a similar way, e.g. Heyting algebras are an algebraic semantics for
intuitionistic logic.

In our understanding, algebraic semantics describe
the abstract properties corresponding model sets have. Structural
semantics, on the other hand, give intuitive concepts like accessibility
or preference, from which properties of model sets, and thus algebraic
semantics, originate.

Varying properties of structural semantics (e.g. transitivity, etc.)
result in
varying properties of algebraic semantics, and thus of logical rules.
We consider operations directly on the algebraic semantics and their
logical
consequences, and we see that simple manipulations of the size concept
result in most rules of nonmonotonic logics. Even more, we show how
to generate new rules from those manipulations. The result is one big
table, which, in a much more modest scale, can be seen as a
``periodic table'' of the ``elements'' of nonmonotonic logic.
Some simple underlying principles allow to generate them all.

Historical remarks: The first time that abstract size was related to
nonmonotonic logics was, to our knowledge,
in the second author's  \cite{Sch90} and  \cite{Sch95-1}, and,
independently, in  \cite{BB94}.
More detailed remarks can e.g. be found in  \cite{GS08c}. But,
again to our knowledge, connections are elaborated systematically and in
fine
detail here for the first time.
\subsubsection{
Overview
}

The main part of this Section is the big table
in Section \ref{Section Main-Table} (page \pageref{Section Main-Table}).
It shows connections and how to develop a multitude of logical rules
known from nonmonotonic logics by combining a small number of principles
about size. We use them as building blocks to construct the rules from.

These principles are some basic and very natural postulates,
$ \xCf (Opt),$ $ \xCf (iM),$ $(eM \xdi ),$ $(eM \xdf ),$ and a continuum
of power of the notion of
``small'', or, dually, ``big'', from $(1*s)$ to $(< \xbo *s).$
From these, we can develop the rest except, essentially, Rational
Monotony,
and thus an infinity of different rules.

This is a conceptual Section, and it does not contain any more difficult
formal results. The interest lies, in our opinion, in the simplicity,
paucity,
and naturalness of the basic building blocks.
We hope that this schema brings more and deeper order into the rich fauna
of nonmonotonic and related logics.
\subsection{
Main table
}
\label{Section Table}
\subsubsection{
Notation
}

 \xEh

 \xDH

$ \xdi (X) \xcc \xdp (X)$ and $ \xdf (X) \xcc \xdp (X)$ are dual abstract
notions of size, $ \xdi (X)$ is the
set of ``small'' subsets of $X,$ $ \xdf (X)$ the set of ``big'' subsets of
$X.$ They are
dual in the sense that $A \xbe \xdi (X) \xcj X-A \xbe \xdf (X).$ `` $ \xdi
$ '' evokes ``ideal'',
`` $ \xdf $ '' evokes ``filter'' though the full strength of both is reached
only
in $(< \xbo *s).$ ``s'' evokes ``small'', and `` $(x*s)$ '' stands for
`` $x$ small sets together are still not everything''.

 \xDH

If $A \xcc X$ is neither in $ \xdi (X),$ nor in $ \xdf (X),$ we say it has
medium size, and
we define $ \xdm (X):= \xdp (X)-( \xdi (X) \xcv \xdf (X)).$ $
\xdm^{+}(X):= \xdp (X)- \xdi (X)$ is the set of subsets
which are not small.

 \xDH

$ \xeA x \xbf $ is a generalized first order quantifier, it is read
``almost all $x$ have property $ \xbf $ ''. $ \xeA x( \xbf: \xbq )$ is the
relativized version, read:
``almost all $x$ with property $ \xbf $ have also property $ \xbq $ ''. To
keep the table
simple, we write mostly only the non-relativized versions.
Formally, we have $ \xeA x \xbf: \xcj \{x: \xbf (x)\} \xbe \xdf (U)$
where $U$ is the universe, and
$ \xeA x( \xbf: \xbq ): \xcj \{x:( \xbf \xcu \xbq )(x)\} \xbe \xdf (\{x:
\xbf (x)\}).$
Soundness and completeness results on $ \xeA $ can be found in
 \cite{Sch95-1}.

 \xDH

Analogously, for propositional logic, we define:

$ \xba \xcn \xbb $ $: \xcj $ $M( \xba \xcu \xbb ) \xbe \xdf (M( \xba )),$

where $M( \xbf )$ is the set of models of $ \xbf.$

 \xDH

In preferential structures, $ \xbm (X) \xcc X$ is the set of minimal
elements of $X.$
This generates a principal filter by $ \xdf (X):=\{A \xcc X: \xbm (X) \xcc
A\}.$ Corresponding
properties about $ \xbm $ are not listed systematically.

 \xDH

The usual rules $ \xCf (AND)$ etc. are named here $(AND_{ \xbo }),$ as
they are in a
natural ascending line of similar rules, based on strengthening of the
filter/ideal properties.

 \xEj
\subsubsection{
The groupes of rules
}

The rules are divided into 5 groups:

 \xEh

 \xDH $ \xCf (Opt),$ which says that ``All'' is optimal - i.e. when there
are no
exceptions, then a soft rule $ \xcn $ holds.

 \xDH 3 monotony rules:

 \xEh
 \xDH $ \xCf (iM)$ is inner monotony, a subset of a small set is small,
 \xDH $(eM \xdi )$ external monotony for ideals: enlarging the base set
keeps small
sets small,
 \xDH $(eM \xdf )$ external monotony for filters: a big subset stays big
when the base
set shrinks.
 \xEj

These three rules are very natural if ``size'' is anything coherent over
change
of base sets. In particular, they can be seen as weakening.

 \xDH $( \xCd )$ keeps proportions, it is here mainly to point the
possibility out.

 \xDH a group of rules $x*s,$ which say how many small sets will not yet
add to
the base set.
 \xDH Rational monotony, which can best be understood as robustness of $
\xdm^{+},$
see $( \xdm^{++})(3).$

 \xEj
\paragraph{
Regularities
}

 \xEh

 \xDH

The group of rules $(x*s)$ use ascending strength of $ \xdi / \xdf.$

 \xDH

The column $( \xdm^{+})$ contains interesting algebraic properties. In
particular,
they show a strengthening from $(3*s)$ up to Rationality. They are not
necessarily
equivalent to the corresponding $(I_{x})$ rules, not even in the presence
of the basic rules. The examples show that care has to be taken when
considering
the different variants.

 \xDH

Adding the somewhat superflous $(CM_{2}),$ we have increasing cautious
monotony from $ \xCf (wCM)$ to full $(CM_{ \xbo }).$

 \xDH

We have increasing ``or'' from $ \xCf (wOR)$ to full $(OR_{ \xbo }).$

 \xDH

The line $(2*s)$ is only there because there seems to be no $(
\xdm^{+}_{2}),$ otherwise
we could begin $(n*s)$ at $n=2.$

 \xEj
\subsubsection{
Direct correspondences
}

Several correspondences are trivial and are mentioned now.
Somewhat less obvious (in)dependencies are given in
Section \ref{Section Coherent-Systems} (page \pageref{Section Coherent-Systems})
.
Finally, the connections with the $ \xbm -$rules are given in
Section \ref{Section Principal} (page \pageref{Section Principal}).
In those rules, $(I_{ \xbo })$ is implicit, as they are about principal
filters.
Still, the $ \xbm -$rules are written in the main table in their
intuitively
adequate place.

 \xEh

 \xDH

The columns ``Ideal'' and ``Filter'' are mutually dual, when both entries are
defined.

 \xDH

The correspondence between the ideal/filter column and the $ \xeA -$column
is obvious,
the latter is added only for completeness' sake, and to point out the
trivial translation to first order logic.

 \xDH

The ideal/filter and the AND-column correspond directly.

 \xDH

We can construct logical rules from the $ \xdm^{+}-column$ by direct
correspondence,
e.g. for $( \xdm^{+}_{ \xbo }),$ (1):

Set $Y:=M( \xbg ),$ $X:=M( \xbg \xcu \xbb ),$ $A:=M( \xbg \xcu \xbb \xcu
\xba ).$

 \xEI
 \xDH $X \xbe \xdm^{+}(Y)$ will become $ \xbg \xcN \xCN \xbb $
 \xDH $A \xbe \xdf (X)$ will become $ \xbg \xcu \xbb \xcn \xba $
 \xDH $A \xbe \xdm^{+}(Y)$ will become $ \xbg \xcN \xCN ( \xba \xcu \xbb
).$
 \xEJ

so we obtain $ \xbg \xcN \xCN \xbb,$ $ \xbg \xcu \xbb \xcn \xba $ $ \xch
$ $ \xbg \xcN \xCN ( \xba \xcu \xbb ).$

We did not want to make the table too complicated, so
such rules are not listed in the table.

 \xDH

Various direct correspondences:

 \xEI

 \xDH In the line $ \xCf (Opt),$ the filter/ideal entry corresponds to $
\xCf (SC),$
 \xDH in the line $ \xCf (iM),$ the filter/ideal entry corresponds to $
\xCf (RW),$
 \xDH in the line $(eM \xdi ),$ the ideal entry corresponds to $(PR' )$
and $ \xCf (wOR),$
 \xDH in the line $(eM \xdf ),$ the filter entry corresponds to $ \xCf
(wCM),$
 \xDH in the line $( \xCd ),$ the filter/ideal entry corresponds to $ \xCf
(disjOR),$
 \xDH in the line $(1*s),$ the filter/ideal entry corresponds to $ \xCf
(CP),$
 \xDH in the line $(2*s),$ the filter/ideal entry corresponds to
$(CM_{2})=(OR_{2}).$

 \xEJ

 \xDH

Note that one can, e.g., write $(AND_{2})$ in two flavours:

 \xEI
 \xDH $ \xba \xcn \xbb,$ $ \xba \xcn \xbb ' $ $ \xch $ $ \xba \xcL \xCN
\xbb \xco \xCN \xbb ',$ or
 \xDH $ \xba \xcn \xbb $ $ \xch $ $ \xba \xcN \xCN \xbb $
 \xEJ

(which is $(CM_{2})=(OR_{2}).)$

For reasons of simplicity, we mention only one.

 \xEj
\subsubsection{
Rational Monotony
}

$ \xCf (RatM)$ does not fit into adding small sets. We have exhausted the
combination
of small sets by $(< \xbo *s),$ unless we go to languages with infinitary
formulas.

The next idea would be to add medium size sets. But, by definition,
$2*medium$
can be all. Adding small and medium sets would not help either: Suppose we
have a rule $medium+n*small \xEd all.$ Taking the complement of the first
medium
set, which is again medium, we have the rule $2*n*small \xEd all.$ So we
do not
see any meaningful new internal rule. i.e. without changing the base set.

Probably, $ \xCf (RatM)$ has more to do with independence: by default, all
``normalities'' are independent, and intersecting with another formula
preserves normality.
\subsubsection{
Summary
}

We can obtain all rules except $ \xCf (RatM)$ and $( \xCd )$ from $ \xCf
(Opt),$ the monotony
rules - $ \xCf (iM),$ $(eM \xdi ),$ $(eM \xdf )$ -, and $(x*s)$ with
increasing $x.$
\subsubsection{
Main table
}
\label{Section Main-Table}
\label{Definition Gen-Filter}
\newpage

\begin{turn}{90}

{\tiny
% {\scriptsize
% {\footnotesize

\begin{tabular}{|c|c@{:}c|c|c|c|c|c|c|}

\hline

\xEH
% \begin{turn}{270}
``Ideal''
% \end{turn}
\xEH
``Filter''
\xEH
$ \xdm^+ $
\xEH
$ \xeA $
\xEH
various rules
\xEH
AND
\xEH
OR
\xEH
Caut./Rat.Mon.
\xEP

\hline
\hline

\multicolumn{9}{|c|}{Optimal proportion} \xEP

\hline

$(Opt)$
\xEH
$ \xCQ \xbe \xdi (X)$
\xEH
$X \xbe \xdf (X)$
\xEH
\xEH
$ \xcA x \xbf \xcp \xeA x \xbf$
\xEH
$(SC)$
\xEH
\xEH
\xEH
\xEP

\xEH
\xEH
\xEH
\xEH
\xEH
$ \xba \xcl \xbb \xch \xba \xcn \xbb $
\xEH
\xEH
\xEH
\xEP

\hline
\hline

\multicolumn{9}{|c|}
{Monotony (Improving proportions). $(iM)$: internal monotony,
$(eM \xdi )$: external monotony for ideals,
$(eM \xdf )$: external monotony for filters}
\xEP

\hline

$(iM)$
\xEH
$A \xcc B \xbe \xdi (X)$ $ \xch $
\xEH
$A \xbe \xdf (X)$, $A \xcc B \xcc X$
\xEH
\xEH
$\xeA x \xbf \xcu \xcA x (\xbf \xcp \xbf')$
\xEH
$(RW)$
\xEH
\xEH
\xEH
\xEP

\xEH
$A \xbe \xdi (X)$
\xEH
$ \xch $ $B \xbe \xdf (X)$
\xEH
\xEH
$ \xcp $ $ \xeA x \xbf'$
\xEH
$ \xba \xcn \xbb, \xbb \xcl \xbb' \xch $
\xEH
\xEH
\xEH
\xEP

\xEH
\xEH
\xEH
\xEH
\xEH
$ \xba \xcn \xbb' $
\xEH
\xEH
\xEH
\xEP

\hline

$(eM \xdi )$
\xEH
$X \xcc Y \xch$
\xEH
\xEH
\xEH
$\xeA x (\xbf: \xbq) \xcu$
\xEH
$(PR')$
\xEH
\xEH
$(wOR)$
\xEH
\xEP

\xEH
$\xdi (X) \xcc \xdi (Y)$
\xEH
\xEH
\xEH
$\xcA x (\xbf' \xcp \xbq) \xcp$
\xEH
$\xba \xcn \xbb, \xba \xcl \xba',$
\xEH
\xEH
$ \xba \xcn \xbb,$ $ \xba ' \xcl \xbb $ $ \xch $
\xEH
\xEP

\xEH
\xEH
\xEH
\xEH
$\xeA x (\xbf \xco \xbf': \xbq)$
\xEH
$\xba' \xcu \xCN \xba \xcl \xbb \xch$
\xEH
\xEH
$ \xba \xco \xba ' \xcn \xbb $
\xEH
\xEP

\xEH
\xEH
\xEH
\xEH
\xEH
$\xba' \xcn \xbb$
\xEH
\xEH
$(\xbm wOR)$
\xEH
\xEP

\xEH
\xEH
\xEH
\xEH
\xEH
$(\xbm PR)$
\xEH
\xEH
$\xbm(X \xcv Y) \xcc \xbm(X) \xcv Y$
\xEH
\xEP

\xEH
\xEH
\xEH
\xEH
\xEH
$X \xcc Y \xch$
\xEH
\xEH
\xEH
\xEP

\xEH
\xEH
\xEH
\xEH
\xEH
$\xbm(Y) \xcs X \xcc \xbm(X)$
\xEH
\xEH
\xEH
\xEP

\hline

$(eM \xdf )$
\xEH
\xEH
$X \xcc Y \xch$
\xEH
\xEH
$\xeA x (\xbf: \xbq) \xcu$
\xEH
\xEH
\xEH
\xEH
$(wCM)$
\xEP

\xEH
\xEH
$\xdf (Y) \xcs \xdp (X) \xcc \xdf (X)$
\xEH
\xEH
$\xcA x (\xbq \xcu \xbf \xcp \xbf') \xcp$
\xEH
\xEH
\xEH
\xEH
$\xba \xcn \xbb, \xba' \xcl \xba,$
\xEP

\xEH
\xEH
\xEH
\xEH
$\xeA x (\xbf \xcu \xbf': \xbq)$
\xEH
\xEH
\xEH
\xEH
$\xba \xcu \xbb \xcl \xba' \xch$
\xEP

\xEH
\xEH
\xEH
\xEH
\xEH
\xEH
\xEH
\xEH
$\xba' \xcn \xbb$
\xEP

\hline
\hline

\multicolumn{9}{|c|}{Keeping proportions} \xEP

\hline

$(\xCd)$
\xEH
$(\xdi \xcv disj)$
\xEH
$(\xdf \xcv disj)$
\xEH
\xEH
$\xeA x(\xbf: \xbq) \xcu$
\xEH
\xEH
\xEH
$(disjOR)$
\xEH
\xEP

\xEH
$A \xbe \xdi (X),$ $B \xbe \xdi (Y),$
\xEH
$A \xbe \xdf (X),$ $B \xbe \xdf (Y),$
\xEH
\xEH
$\xeA x(\xbf': \xbq) \xcu$
\xEH
\xEH
\xEH
$ \xbf \xcn \xbq,$ $ \xbf ' \xcn \xbq ' $
\xEH
\xEP

\xEH
$X \xcs Y= \xCQ $ $ \xch $
\xEH
$X \xcs Y= \xCQ $ $ \xch $
\xEH
\xEH
$\xCN \xcE x(\xbf \xcu \xbf') \xcp$
\xEH
\xEH
\xEH
$ \xbf \xcl \xCN \xbf ',$ $ \xch $
\xEH
\xEP

\xEH
$A \xcv B \xbe \xdi (X \xcv Y)$
\xEH
$A \xcv B \xbe \xdf (X \xcv Y)$
\xEH
\xEH
$\xeA x(\xbf \xco \xbf': \xbq)$
\xEH
\xEH
\xEH
$ \xbf \xco \xbf ' \xcn \xbq \xco \xbq ' $
\xEH
\xEP

\xEH
\xEH
\xEH
\xEH
\xEH
\xEH
\xEH
$(\xbm disjOR)$
\xEH
\xEP

\xEH
\xEH
\xEH
\xEH
\xEH
\xEH
\xEH
$X \xcs Y = \xCQ \xch$
\xEH
\xEP

\xEH
\xEH
\xEH
\xEH
\xEH
\xEH
\xEH
$\xbm(X \xcv Y) \xcc \xbm(X) \xcv \xbm(Y)$
\xEH
\xEP

\hline
\hline

\multicolumn{9}{|c|}{Robustness of proportions: $n*small \xEd All$} \xEP

\hline

$(1*s)$
\xEH
$(\xdi_1)$
\xEH
$(\xdf_1)$
\xEH
\xEH
$(\xeA_1)$
\xEH
$(CP)$
\xEH
$(AND_{1})$
\xEH
\xEH
\xEP

\xEH
$X \xce \xdi (X)$
\xEH
$ \xCQ \xce \xdf (X)$
\xEH
\xEH
$ \xeA x \xbf \xcp \xcE x \xbf $
\xEH
$\xbf\xcn\xcT \xch \xbf\xcl\xcT$
\xEH
$ \xba \xcn \xbb $ $ \xch $ $ \xba \xcL \xCN \xbb $
\xEH
\xEH
\xEP

\hline

$(2*s)$
\xEH
$(\xdi_2)$
\xEH
$(\xdf_2)$
\xEH
\xEH
$(\xeA_2)$
\xEH
\xEH
$(AND_{2})$
\xEH
$(OR_{2})$
\xEH
$(CM_{2})$
\xEP

\xEH
$A,B \xbe \xdi (X) \xch $
\xEH
$A,B \xbe \xdf (X) \xch $
\xEH
\xEH
$ \xeA x \xbf \xcu \xeA x \xbq $
\xEH
\xEH
$ \xba \xcn \xbb,$ $ \xba \xcn \xbb ' $ $ \xch $
\xEH
$ \xba \xcn \xbb \xch \xba \xcN \xCN \xbb $
\xEH
$ \xba \xcn \xbb \xch \xba \xcN \xCN \xbb $
\xEP

\xEH
$ A \xcv B \xEd X$
\xEH
$A \xcs B \xEd \xCQ $
\xEH
\xEH
$ \xcp $ $ \xcE x( \xbf \xcu \xbq )$
\xEH
\xEH
$ \xba \xcL \xCN \xbb \xco \xCN \xbb ' $
\xEH
\xEH
\xEP

\hline

%  $(3*s)$
%  \xEH
%  $(\xdi_3)$
%  \xEH
%  $(\xdf_3)$
%  \xEH
%  $( \xdm^{+}_{3})$
%  \xEH
%  $(\xeA_3)$
%  \xEH
%  \xEH
%  $(AND_{3})$
%  \xEH
%  $(OR_{3})$
%  \xEH
%  $(CM_{3})$
%  \xEP
%
%  \xEH
%  $A,B,C \xbe \xdi (X) \xch $
%  \xEH
%  $A,B,C \xbe \xdf (X) \xch $
%  \xEH
%  $A \xbe \xdf (X),$ $X \xbe \xdf (Y)$
%  \xEH
%  $ \xeA x \xbf \xcu \xeA x \xbq \xcu \xeA x \xbs $
%  \xEH
%  \xEH
%  $ \xba \xcn \xbb,$ $ \xba \xcn \xbb ',$ $ \xba \xcn \xbb '' $
%  \xEH
%  $ \xba \xcn \xbb,$ $ \xba ' \xcn \xbb $ $ \xch $
%  \xEH
%  $ \xba \xcn \xbb,$ $ \xba \xcn \xbb ' $ $ \xch $
%  \xEP
%
%  \xEH
%  $ A \xcv B \xcv C \xEd X$
%  \xEH
%  $ A \xcs B \xcs C \xEd \xCQ $
%  \xEH
%   $ \xch $ $A \xbe \xdm^{+}(Y)$
%  \xEH
%  $ \xcp \xcE x( \xbf \xcu \xbq \xcu \xbs )$
%  \xEH
%  \xEH
%   $ \xch $ $ \xba \xcL \xCN \xbb \xco \xCN \xbb ' \xco \xCN \xbb '' $
%  \xEH
%  $ \xba \xco \xba ' \xcN \xCN \xbb $
%  \xEH
%  $ \xba \xcu \xbb \xcN \xCN \xbb ' $
%  \xEP
%
%  \hline

$(n*s)$
\xEH
$(\xdi_n)$
\xEH
$(\xdf_n)$
\xEH
$( \xdm^{+}_{n})$
\xEH
$(\xeA_n)$
\xEH
\xEH
$(AND_{n})$
\xEH
$(OR_{n})$
\xEH
$(CM_{n})$
\xEP

$(n \xcg 3)$
\xEH
$A_{1},.,A_{n} \xbe \xdi (X) $
\xEH
$A_{1},.,A_{n} \xbe \xdi (X) $
\xEH
$X_{1} \xbe \xdf (X_{2}),., $
\xEH
$ \xeA x \xbf_{1} \xcu.  \xcu \xeA x \xbf_{n} $
\xEH
\xEH
$ \xba \xcn \xbb_{1},., \xba \xcn \xbb_{n}$ $ \xch $
\xEH
$ \xba_{1} \xcn \xbb,., \xba_{n-1} \xcn \xbb $
\xEH
$ \xba \xcn \xbb_{1},., \xba \xcn \xbb_{n-1}$
\xEP

\xEH
$ \xch $
\xEH
$ \xch $
\xEH
$ X_{n-1} \xbe \xdf (X_{n})$ $ \xch $
\xEH
$ \xcp $
\xEH
\xEH
$ \xba \xcL \xCN \xbb_{1} \xco.  \xco \xCN \xbb_{n}$
\xEH
$ \xch $
\xEH
$ \xch $
\xEP

\xEH
$ A_{1} \xcv.  \xcv A_{n} \xEd X $
\xEH
$A_{1} \xcs.  \xcs A_{n} \xEd \xCQ$
\xEH
$X_{1} \xbe \xdm^{+}(X_{n})$
\xEH
$ \xcE x (\xbf_{1} \xcu.  \xcu \xbf_{n}) $
\xEH
\xEH
\xEH
$ \xba_{1} \xco.  \xco \xba_{n-1} \xcN \xCN \xbb $
\xEH
$ \xba \xcu \xbb_1 \xcu.  \xcu \xbb_{n-2} \xcN \xCN \xbb_{n-1}$
\xEP

\hline

$(< \xbo*s)$
\xEH
$(\xdi_\xbo)$
\xEH
$(\xdf_\xbo)$
\xEH
$( \xdm^{+}_{ \xbo })$
\xEH
$(\xeA_{\xbo})$
\xEH
\xEH
$(AND_{ \xbo })$
\xEH
$(OR_{ \xbo })$
\xEH
$(CM_{ \xbo })$
\xEP

\xEH
$A,B \xbe \xdi (X) \xch $
\xEH
$A,B \xbe \xdf (X) \xch $
\xEH
(1)
\xEH
$ \xeA x \xbf \xcu \xeA x \xbq \xcp $
\xEH
\xEH
$ \xba \xcn \xbb,$ $ \xba \xcn \xbb ' $ $ \xch $
\xEH
$ \xba \xcn \xbb,$ $ \xba ' \xcn \xbb $ $ \xch $
\xEH
$ \xba \xcn \xbb,$ $ \xba \xcn \xbb ' $ $ \xch $
\xEP

\xEH
$ A \xcv B \xbe \xdi (X)$
\xEH
$ A \xcs B \xbe \xdf (X)$
\xEH
$A \xbe \xdf (X),$ $X \xbe \xdm^{+}(Y)$
\xEH
$ \xeA x( \xbf \xcu \xbq )$
\xEH
\xEH
$ \xba \xcn \xbb \xcu \xbb ' $
\xEH
$ \xba \xco \xba ' \xcn \xbb $
\xEH
$ \xba \xcu \xbb \xcn \xbb ' $
\xEP

\xEH
\xEH
\xEH
$ \xch $ $A \xbe \xdm^{+}(Y)$
\xEH
\xEH
\xEH
\xEH
$(\xbm OR)$
\xEH
$(\xbm CM)$
\xEP

\xEH
\xEH
\xEH
(2)
\xEH
\xEH
\xEH
\xEH
$\xbm(X \xcv Y) \xcc \xbm(X) \xcv \xbm(Y)$
\xEH
$\xbm(X) \xcc Y \xcc X \xch$
\xEP

\xEH
\xEH
\xEH
$A \xbe \xdm^{+}(X),$ $X \xbe \xdf (Y)$
\xEH
\xEH
\xEH
\xEH
\xEH
$\xbm(Y) \xcc \xbm(X)$
\xEP

\xEH
\xEH
\xEH
$ \xch $ $A \xbe \xdm^{+}(Y)$
\xEH
\xEH
\xEH
\xEH
\xEH
\xEP

\xEH
\xEH
\xEH
(3)
\xEH
\xEH
\xEH
\xEH
\xEH
\xEP

\xEH
\xEH
\xEH
$A \xbe \xdf (X),$ $X \xbe \xdf (Y)$
\xEH
\xEH
\xEH
\xEH
\xEH
\xEP

\xEH
\xEH
\xEH
$ \xch $ $A \xbe \xdf (Y)$
\xEH
\xEH
\xEH
\xEH
\xEH
\xEP

\xEH
\xEH
\xEH
(4)
\xEH
\xEH
\xEH
\xEH
\xEH
\xEP

\xEH
\xEH
\xEH
$A,B \xbe \xdi (X)$ $ \xch $
\xEH
\xEH
\xEH
\xEH
\xEH
\xEP

\xEH
\xEH
\xEH
$A-B \xbe \xdi (X-$B)
\xEH
\xEH
\xEH
\xEH
\xEH
\xEP

\hline
\hline

\multicolumn{9}{|c|}{Robustness of $\xdm^+$} \xEP

\hline

$(\xdm^{++})$
\xEH
\xEH
\xEH
$(\xdm^{++})$
\xEH
\xEH
\xEH
\xEH
\xEH
$(RatM)$
\xEP

\xEH
\xEH
\xEH
(1)
\xEH
\xEH
\xEH
\xEH
\xEH
$ \xbf \xcn \xbq,  \xbf \xcN \xCN \xbq '   \xch $
\xEP

\xEH
\xEH
\xEH
$A \xbe \xdi (X),$ $B \xce \xdf (X)$
\xEH
\xEH
\xEH
\xEH
\xEH
$ \xbf \xcu \xbq ' \xcn \xbq $
\xEP

\xEH
\xEH
\xEH
$ \xch $ $A-B \xbe \xdi (X-B)$
\xEH
\xEH
\xEH
\xEH
\xEH
$(\xbm RatM)$
\xEP

\xEH
\xEH
\xEH
(2)
\xEH
\xEH
\xEH
\xEH
\xEH
$X \xcc Y,$
\xEP

\xEH
\xEH
\xEH
$A \xbe \xdf (X), B \xce \xdf (X)$
\xEH
\xEH
\xEH
\xEH
\xEH
$X \xcs \xbm(Y) \xEd \xCQ \xch$
\xEP

\xEH
\xEH
\xEH
$ \xch $ $A-B \xbe \xdf (X-$B)
\xEH
\xEH
\xEH
\xEH
\xEH
$\xbm(X) \xcc \xbm(Y) \xcs X$
\xEP

\xEH
\xEH
\xEH
(3)
\xEH
\xEH
\xEH
\xEH
\xEH
\xEP

\xEH
\xEH
\xEH
$A \xbe \xdm^+ (X),$
\xEH
\xEH
\xEH
\xEH
\xEH
\xEP

\xEH
\xEH
\xEH
$X \xbe \xdm^+ (Y)$
\xEH
\xEH
\xEH
\xEH
\xEH
\xEP

\xEH
\xEH
\xEH
$ \xch $ $A \xbe \xdm^+ (Y)$
\xEH
\xEH
\xEH
\xEH
\xEH
\xEP

\hline

\hline

\end{tabular}

}

\end{turn}

\newpage
\index{Remark Ref-Class-Short}

\br

$\hspace{0.01em}$

% (+++ Orig. No.:  Remark Ref-Class-Short +++)

\label{Remark Ref-Class-Short}

There is, however, an important conceptual distinction to make here.
Filters
express ``size'' in an abstract way, in the context of
nonmonotonic logics, $ \xba \xcn \xbb $ iff the set of
$ \xba \xcu \xCN \xbb $ is small in $ \xba.$ But here, we were
interested in ``small'' changes in the reference set $X$ (or $ \xba $ in our
example). So
we have two quite different uses of ``size'', one for
nonmonotonic logics, abstractly expressed by
a filter, the other for coherence conditions. It is possible, but not
necessary,
to consider both essentially the same notions. But we should not forget
that we
have two conceptually different uses of size here.
\subsection{
Coherent systems
}
\label{Section Coherent-Systems}
\subsubsection{
Definition and basic facts
}

\er

Note that whenever we work with model sets, the rule

$ \xCf (LLE),$ left logical equivalence, $ \xcl \xba \xcr \xba ' $ $ \xch
$ $( \xba \xcn \xbb $ $ \xcj $ $ \xba ' \xcn \xbb )$

will hold. We will not mention this any further.

\bd

$\hspace{0.01em}$

% (+++ Orig. No.:  Definition CoherentSystem +++)

\label{Definition CoherentSystem}

A coherent system of sizes, $ \xdc \xds,$ consists of a universe $U,$ $
\xCQ \xce \xdy \xcc \xdp (U),$ and
for all $X \xbe \xdy $ a system $ \xdi (X) \xcc \xdp (X)$ (dually $ \xdf
(X),$ i.e. $A \xbe \xdf (X) \xcj X-A \xbe \xdi (X)).$
$ \xdy $ may satisfy certain closure properties like closure under $ \xcv
,$ $ \xcs,$
complementation, etc. We will mention this when needed, and not obvious.

We say that $ \xdc \xds $ satisfies a certain property iff all $X,Y \xbe
\xdy $ satisfy this
property.

$ \xdc \xds $ is called basic or level 1 iff it satisfies $ \xCf (Opt),$ $
\xCf (iM),$ $(eM \xdi ),$
$(eM \xdf ),$ $(1*s).$

$ \xdc \xds $ is level $x$ iff it satisfies $ \xCf (Opt),$ $ \xCf (iM),$
$(eM \xdi ),$ $(eM \xdf ),$ $(x*s).$

\ed

\bfa

$\hspace{0.01em}$

% (+++ Orig. No.:  Fact 1-element +++)

\label{Fact 1-element}

Note that, if for any $Y$ $ \xdi (Y)$ consists only of subsets of at most
1 element,
then $(eM \xdf )$ is trivially satisfied for $Y$ and its subsets by $ \xCf
(Opt).$ $ \xcz $
\\[3ex]

\efa

\bfa

$\hspace{0.01em}$

% (+++ Orig. No.:  Fact Not-2*s +++)

\label{Fact Not-2*s}

Let a $ \xdc \xds $ be given s.t. $ \xdy = \xdp (U).$ If $X \xbe \xdy $
satisfies $( \xdm^{++}),$ but not
$(< \xbo *s),$ then there is $Y \xbe \xdy $ which does not satisfy
$(2*s).$

\efa

\subparagraph{
Proof
}

$\hspace{0.01em}$

% (+++ Orig.:  Proof +++)

We work with version (1) of $( \xdm^{++}),$ we will see in
Fact \ref{Fact M-plus-plus} (page \pageref{Fact M-plus-plus})  that all three
versions are equivalent.

As $X$ does not satisfy $(< \xbo *s),$ there are $A,B \xbe \xdi (X)$ s.t.
$A \xcv B \xbe \xdm^{+}(X).$
$A \xbe \xdi (X),$ $A \xcv B \xbe \xdm^{+}(X)$ $ \xch $ $X-(A \xcv B) \xce
\xdf (X),$ so by $( \xdm^{++})(1)$
$A=A-(X-(A \xcv B)) \xbe \xdi (X-(X-(A \xcv B)))= \xdi (A \xcv B).$
Likewise $B \xbe \xdi (A \xcv B),$
so $(2*s)$ does not hold for $A \xcv B.$ $ \xcz $
\\[3ex]

\bfa

$\hspace{0.01em}$

% (+++ Orig. No.:  Fact Independence-eM +++)

\label{Fact Independence-eM}

$(eM \xdi )$ and $(eM \xdf )$ are formally independent, though intuitively
equivalent.

\efa

\subparagraph{
Proof
}

$\hspace{0.01em}$

% (+++ Orig.:  Proof +++)

Let $U:=\{x,y,z\},$ $X:=\{x,z\},$ $ \xdy:= \xdp (U)-\{ \xCQ \}$

(1) Let $ \xdf (U):=\{A \xcc U:z \xbe A\},$ $ \xdf (Y)=\{Y\}$ for all $Y
\xcb U.$
$ \xCf (Opt),$ $ \xCf (iM)$ hold, $(eM \xdi )$ holds trivially, so does
$(< \xbo *s),$
but $(eM \xdf )$ fails for $U$ and $X.$

(2) Let $ \xdf (X):=\{\{z\},X\},$ $ \xdf (Y):=\{Y\}$ for all $Y \xcc U,$
$Y \xEd X.$
$ \xCf (Opt),$ $ \xCf (iM),$ $(< \xbo *s)$ hold trivially, $(eM \xdf )$
holds by
Fact \ref{Fact 1-element} (page \pageref{Fact 1-element}). $(eM \xdi )$ fails,
as $\{x\} \xbe \xdi
(X),$ but $\{x\} \xce \xdi (U).$

$ \xcz $
\\[3ex]

\bfa

$\hspace{0.01em}$

% (+++ Orig. No.:  Fact Level-n-n+1 +++)

\label{Fact Level-n-n+1}

A level $n$ system is strictly weaker than a level $n+1$ system.

\efa

\subparagraph{
Proof
}

$\hspace{0.01em}$

% (+++ Orig.:  Proof +++)

Consider $U:=\{1, \Xl,n+1\},$ $ \xdy:= \xdp (U)-\{ \xCQ \}.$ Let $ \xdi
(U):=\{ \xCQ \} \xcv \{\{x\}:x \xbe U\},$
$ \xdi (X):=\{ \xCQ \}$ for $X \xEd U.$
$ \xCf (iM),$ $(eM \xdi ),$ $(eM \xdf )$ hold trivially.
$(n*s)$ holds trivially for $X \xEd U,$ but also for $U.$ $((n+1)*s)$ does
not hold
for $U.$ $ \xcz $
\\[3ex]

\br

$\hspace{0.01em}$

% (+++ Orig. No.:  Remark Infin +++)

\label{Remark Infin}

Note that our schemata allow us to generate infintely many new rules, here
is
an example:

Start with A, add $s_{1,1},$ $s_{1,2}$ two sets small in $A \xcv s_{1,1}$
$(A \xcv s_{1,2}$ respectively).
Consider now $A \xcv s_{1,1} \xcv s_{1,2}$ and $s_{2}$ s.t. $s_{2}$ is
small in $A \xcv s_{1,1} \xcv s_{1,2} \xcv s_{2}.$
Continue with $s_{3,1},$ $s_{3,2}$ small in $A \xcv s_{1,1} \xcv s_{1,2}
\xcv s_{2} \xcv s_{3,1}$ etc.

Without additional properties, this system creates a new rule, which is
not equivalent to any usual rules.

$ \xcz $
\\[3ex]
\subsubsection{
The finite versions
}

\er

\bfa

$\hspace{0.01em}$

% (+++ Orig. No.:  Fact I-n +++)

\label{Fact I-n}

(1) $(I_{n})$ $+$ $(eM \xdi )$ $ \xch $ $( \xdm^{+}_{n}),$

(2) $(I_{n})$ $+$ $(eM \xdi )$ $ \xch $ $(CM_{n}),$

(3) $(I_{n})$ $+$ $(eM \xdi )$ $ \xch $ $(OR_{n}).$

\efa

\subparagraph{
Proof
}

$\hspace{0.01em}$

% (+++ Orig.:  Proof +++)

(1)

Let $X_{1} \xcc  \Xl  \xcc X_{n},$ so
$X_{n}=X_{1} \xcv (X_{2}-X_{1}) \xcv  \Xl  \xcv (X_{n}-X_{n-1}).$ Let
$X_{i} \xbe \xdf (X_{i+1}),$ so $X_{i+1}-X_{i} \xbe \xdi (X_{i+1}) \xcc
\xdi (X_{n})$
by $(eM \xdi )$ for $1 \xck i \xck n-1,$ so by $(I_{n})$ $X_{1} \xbe
\xdm^{+}(X_{n}).$

(2)

Suppose $ \xba \xcn \xbb_{1}, \Xl, \xba \xcn \xbb_{n-1},$ but $ \xba \xcu
\xbb_{1} \xcu  \Xl  \xcu \xbb_{n-2} \xcn \xCN \xbb_{n-1}.$
Then $M( \xba \xcu \xCN \xbb_{1}), \Xl,M( \xba \xcu \xCN \xbb_{n-1}) \xbe
\xdi (M( \xba )),$ and
$M( \xba \xcu \xbb_{1} \xcu  \Xl  \xcu \xbb_{n-2} \xcu \xbb_{n-1}) \xbe
\xdi (M( \xba \xcu \xbb_{1} \xcu  \Xl  \xcu \xbb_{n-2})) \xcc \xdi (M(
\xba ))$ by $(eM \xdi ).$
But $M( \xba )=M( \xba \xcu \xCN \xbb_{1}) \xcv  \Xl  \xcv M( \xba \xcu
\xCN \xbb_{n-1}) \xcv M( \xba \xcu \xbb_{1} \xcu  \Xl  \xcu \xbb_{n-2}
\xcu \xbb_{n-1})$ is
now the union of $n$ small subsets, $contradiction.$

(3)

Let $ \xba_{1} \xcn \xbb, \Xl, \xba_{n-1} \xcn \xbb,$ so $M( \xba_{i}
\xcu \xCN \xbb ) \xbe \xdi (M( \xba_{i}))$ for $1 \xck i \xck n-1,$ so
$M( \xba_{i} \xcu \xCN \xbb ) \xbe \xdi (M( \xba_{1} \xco  \Xl  \xco
\xba_{n-1}))$ for $1 \xck i \xck n-1$ by $(eM \xdi ),$ so
$M(( \xba_{1} \xco  \Xl  \xco \xba_{n-1}) \xcu \xbb )$ $=$ $M( \xba_{1}
\xco  \Xl  \xco \xba_{n-1})- \xcv \{M( \xba_{i} \xcu \xCN \xbb ):1 \xck i
\xck n-1\}$ $ \xce $
$ \xdi (M( \xba_{1} \xco  \Xl  \xco \xba_{n-1}))$ by $(I_{n}),$ so $
\xba_{1} \xco  \Xl  \xco \xba_{n-1} \xcN \xCN \xbb.$

$ \xcz $
\\[3ex]

In the following example, $(OR_{n}),$ $( \xdm^{+}_{n}),$ $(CM_{n})$ hold,
but $( \xdi_{n})$ fails, so
by Fact \ref{Fact I-n} (page \pageref{Fact I-n})  $( \xdi_{n})$ is strictly
stronger
than $(OR_{n}),$ $( \xdm^{+}_{n}),$ $(CM_{n}).$

\be

$\hspace{0.01em}$

% (+++ Orig. No.:  Example Not-I-n +++)

\label{Example Not-I-n}

Let $n \xcg 3.$

Consider $X:=\{1, \Xl,n\},$ $ \xdy:= \xdp (X)-\{ \xCQ \},$
$ \xdi (X):=\{ \xCQ \} \xcv \{\{i\}:1 \xck i \xck n\},$ and for all $Y
\xcb X$ $ \xdi (Y):=\{ \xCQ \}.$

$ \xCf (Opt),$ $ \xCf (iM),$ $(eM \xdi ),$ $(eM \xdf )$ (by Fact \ref{Fact
1-element} (page \pageref{Fact 1-element}) ),
$(1*s),$ $(2*s)$ hold, $(I_{n})$ fails, of course.

(1) $(OR_{n})$ holds:

Suppose $ \xba_{1} \xcn \xbb, \Xl, \xba_{n-1} \xcn \xbb,$ $ \xba_{1}
\xco  \Xl  \xco \xba_{n-1} \xcn \xCN \xbb.$

Case 1: $ \xba_{1} \xco  \Xl  \xco \xba_{n-1} \xcl \xCN \xbb,$ then for
all $i$ $ \xba_{i} \xcl \xCN \xbb,$ so for no $i$ $ \xba_{i} \xcn \xbb $
by $(1*s)$ and thus $(AND_{1}),$ $contradiction.$

Case 2: $ \xba_{1} \xco  \Xl  \xco \xba_{n-1} \xcL \xCN \xbb,$ then $M(
\xba_{1} \xco  \Xl  \xco \xba_{n-1})=X,$ and there is exactly
1 $k \xbe X$ s.t. $k \xcm \xbb.$ Fix this $k.$
By prerequisite, $ \xba_{i} \xcn \xbb.$ If $M( \xba_{i})=X,$ $ \xba_{i}
\xcl \xbb $ cannot be, so there must be
exactly 1 $k' $ s.t. $k' \xcm \xCN \xbb,$ but $card(X) \xcg 3,$
$contradiction.$ So $M( \xba_{i}) \xcb X,$ and $ \xba_{i} \xcl \xbb,$ so
$M( \xba_{i})= \xCQ $ or
$M( \xba_{i})=\{k\}$ for all $i,$ so $M( \xba_{1} \xco  \Xl  \xco
\xba_{n-1}) \xEd X,$ $contradiction.$

(2) $( \xdm^{+}_{n})$ holds:

$( \xdm^{+}_{n})$ is a consequence of $( \xdm^{+}_{ \xbo }),$ (3) so it
suffices to show that the latter
holds. Let $X_{1} \xbe \xdf (X_{2}),$ $X_{2} \xbe \xdf (X_{3}).$ Then
$X_{1}=X_{2}$ or $X_{2}=X_{3},$ so the result is
trivial.

(3) $(CM_{n})$ holds:

Suppose $ \xba \xcn \xbb_{1}, \Xl, \xba \xcn \xbb_{n-1},$ $ \xba \xcu
\xbb_{1} \xcu  \Xl  \xcu \xbb_{n-2} \xcn \xCN \xbb_{n-1}.$

Case 1: For all $i,$ $1 \xck i \xck n-2,$ $ \xba \xcl \xbb_{i},$ then $M(
\xba \xcu \xbb_{1} \xcu  \Xl  \xcu \xbb_{n-2})=M( \xba ),$
so $ \xba \xcn \xbb_{n-1}$ and $ \xba \xcn \xCN \xbb_{n-1},$
$contradiction.$

Case 2: There is $i,$ $1 \xck i \xck n-2,$ $ \xba \xcL \xbb_{i},$ then $M(
\xba )=X,$
$M( \xba \xcu \xbb_{1} \xcu  \Xl  \xcu \xbb_{n-2}) \xcb M( \xba ),$ so $
\xba \xcu \xbb_{1} \xcu  \Xl  \xcu \xbb_{n-2} \xcl \xCN \xbb_{n-1}.$
$Card(M( \xba \xcu \xbb_{1} \xcu  \Xl  \xcu \xbb_{n-2})) \xcg n-(n-2)=2,$
so $card(M( \xCN \xbb_{n-1})) \xcg 2,$ so
$ \xba \xcN \xbb_{n-1},$ $contradiction.$

$ \xcz $
\\[3ex]
\subsubsection{
The $\xbo$ version
}

\ee

\bfa

$\hspace{0.01em}$

% (+++ Orig. No.:  Fact CM-Omega +++)

\label{Fact CM-Omega}

$(CM_{ \xbo })$ $ \xcj $ $( \xdm^{+}_{ \xbo })$ (4)

\efa

\subparagraph{
Proof
}

$\hspace{0.01em}$

% (+++ Orig.:  Proof +++)

`` $ \xch $ ''

Suppose all sets are definable.

Let $A,B \xbe \xdi (X),$
$X=M( \xba ),$ $A=M( \xba \xcu \xCN \xbb ),$ $B=M( \xba \xcu \xCN \xbb '
),$ so $ \xba \xcn \xbb,$ $ \xba \xcn \xbb ',$ so by $(CM_{ \xbo })$
$ \xba \xcu \xbb ' \xcn \xbb,$ so $A-B=M( \xba \xcu \xbb ' \xcu \xCN \xbb
) \xbe \xdi (M( \xba \xcu \xbb ' ))= \xdi (X-$B).

`` $ \xci $ ''

Let $ \xba \xcn \xbb,$ $ \xba \xcn \xbb ',$ so $M( \xba \xcu \xCN \xbb )
\xbe \xdi (M( \xba )),$ $M( \xba \xcu \xCN \xbb ' ) \xbe \xdi (M( \xba
)),$ so
by prerequisite $M( \xba \xcu \xCN \xbb ' )-M( \xba \xcu \xCN \xbb )=M(
\xba \xcu \xbb \xcu \xCN \xbb ' )$ $ \xbe $
$ \xdi (M( \xba )-M( \xba \xcu \xCN \xbb ))= \xdi (M( \xba \xcu \xbb )),$
so $ \xba \xcu \xbb \xcn \xbb '.$

$ \xcz $
\\[3ex]

\bfa

$\hspace{0.01em}$

% (+++ Orig. No.:  Fact I-Omega +++)

\label{Fact I-Omega}

(1) $(I_{ \xbo })$ $+$ $(eM \xdi )$ $ \xch $ $(OR_{ \xbo }),$

(2) $(I_{ \xbo })$ $+$ $(eM \xdi )$ $ \xch $ $( \xdm^{+}_{ \xbo })$ (1),

(3) $(I_{ \xbo })$ $+$ $(eM \xdf )$ $ \xch $ $( \xdm^{+}_{ \xbo })$ (2),

(4) $(I_{ \xbo })$ $+$ $(eM \xdi )$ $ \xch $ $( \xdm^{+}_{ \xbo })$ (3),

(5) $(I_{ \xbo })$ $+$ $(eM \xdf )$ $ \xch $ $( \xdm^{+}_{ \xbo })$ (4)
(and thus, by
Fact \ref{Fact CM-Omega} (page \pageref{Fact CM-Omega}), $(CM_{ \xbo })).$

\efa

\subparagraph{
Proof
}

$\hspace{0.01em}$

% (+++ Orig.:  Proof +++)

(1)

Let $ \xba \xcn \xbb,$ $ \xba ' \xcn \xbb $ $ \xch $ $M( \xba \xcu \xCN
\xbb ) \xbe \xdi (M( \xba )),$ $M( \xba ' \xcu \xCN \xbb ) \xbe \xdi (M(
\xba ' )),$
so by $(eM \xdi )$ $M( \xba \xcu \xCN \xbb ) \xbe \xdi (M( \xba \xco \xba
' )),$ $M( \xba ' \xcu \xCN \xbb ) \xbe \xdi (M( \xba \xco \xba ' )),$
so $M(( \xba \xco \xba ' ) \xcu \xCN \xbb ) \xbe \xdi (M( \xba \xco \xba '
))$ by $(I_{ \xbo }),$ so $ \xba \xco \xba ' \xcn \xbb.$

(2)

Let $A \xcc X \xcc Y,$
$A \xbe \xdi (Y),$ $X-A \xbe \xdi (X) \xcc_{(eM \xdi )} \xdi (Y)$ $ \xch $
$X=(X-A) \xcv A \xbe \xdi (Y)$ by $(I_{ \xbo }).$

(3)

Let $A \xcc X \xcc Y,$
let $A \xbe \xdi (Y),$ $Y-X \xbe \xdi (Y)$ $ \xch $ $A \xcv (Y-X) \xbe
\xdi (Y)$ by $(I_{ \xbo })$ $ \xch $
$X-A=Y-(A \xcv (Y-X)) \xbe \xdf (Y)$ $ \xch $ $X-A \xbe \xdf (X)$ by $(eM
\xdf ).$

(4)

Let $A \xcc X \xcc Y,$ $A \xbe \xdf (X),$ $X \xbe \xdf (Y),$ so
$Y-X \xbe \xdi (Y),$ $X-A \xbe \xdi (X) \xcc_{(eM \xdi )} \xdi (Y)$ $ \xch
$ $Y-A=(Y-X) \xcv (X-A) \xbe \xdi (Y)$ by $( \xdi_{ \xbo })$ $ \xch $
$A \xbe \xdf (Y).$

(5)

Let $A,B \xcc X,$
$A,B \xbe \xdi (X)$ $ \xch_{(I_{ \xbo })}$ $A \xcv B \xbe \xdi (X)$ $ \xch
$ $X-(A \xcv B) \xbe \xdf (X),$ but $X-(A \xcv B) \xcc X-$B,
so $X-(A \xcv B) \xbe \xdf (X-$B) by $(eM \xdf ),$ so $A-B=(X-B)-(X-(A
\xcv B)) \xbe \xdi (X-$B).

$ \xcz $
\\[3ex]

We give three examples of independence of the various versions of
$( \xdm^{+}_{ \xbo }).$

\be

$\hspace{0.01em}$

% (+++ Orig. No.:  Example Versions-M-Omega +++)

\label{Example Versions-M-Omega}

All numbers refer to the versions of $( \xdm^{+}_{ \xbo }).$

For easier reading, we re-write for $A \xcc X \xcc Y$

$( \xdm^{+}_{ \xbo })(1):$ $A \xbe \xdf (X),$ $A \xbe \xdi (Y)$ $ \xch $
$X \xbe \xdi (Y),$

$( \xdm^{+}_{ \xbo })(2):$ $X \xbe \xdf (Y),$ $A \xbe \xdi (Y)$ $ \xch $
$A \xbe \xdi (X).$

We give three examples. Investigating all possibilities exhaustively
seems quite tedious, and might best be done with the help of a computer.
Fact \ref{Fact 1-element} (page \pageref{Fact 1-element})  will be used
repeatedly.

 \xEI

 \xDH

(1), (2), (4) fail, (3) holds:

Let $Y:=\{a,b,c\},$ $ \xdy:= \xdp (Y)-\{ \xCQ \},$ $ \xdf
(Y):=\{\{a,c\},$ $\{b,c\},$ $Y\}$

Let $X:=\{a,b\},$ $ \xdf (X):=\{\{a\},$ $X\},$ $A:=\{a\},$ and $ \xdf
(Z):=\{Z\}$ for all $Z \xEd X,Y.$

$ \xCf (Opt),$ $ \xCf (iM),$ $(eM \xdi ),$ $(eM \xdf )$ hold, $(I_{ \xbo
})$ fails, of course.

(1) fails: $A \xbe \xdf (X),$ $A \xbe \xdi (Y),$ $X \xce \xdi (Y).$

(2) fails: $\{a,c\} \xbe \xdf (Y),$ $\{a\} \xbe \xdi (Y),$ but $\{a\} \xce
\xdi (\{a,c\}).$

(3) holds: If $X_{1} \xbe \xdf (X_{2}),$ $X_{2} \xbe \xdf (X_{3}),$ then
$X_{1}=X_{2}$ or $X_{2}=X_{3},$ so (3) holds
trivially (note that $X \xce \xdf (Y)).$

(4) fails: $\{a\},\{b\} \xbe \xdi (Y),$ $\{a\} \xce \xdi (Y-\{b\})= \xdi
(\{a,c\})=\{ \xCQ \}.$

 \xDH

(2), (3), (4) fail, (1) holds:

Let $Y:=\{a,b,c\},$ $ \xdy:= \xdp (Y)-\{ \xCQ \},$ $ \xdf
(Y):=\{\{a,b\},$ $\{a,c\},$ $Y\}$

Let $X:=\{a,b\},$ $ \xdf (X):=\{\{a\},$ $X\},$ and $ \xdf (Z):=\{Z\}$ for
all $Z \xEd X,Y.$

$ \xCf (Opt),$ $ \xCf (iM),$ $(eM \xdi ),$ $(eM \xdf )$ hold, $(I_{ \xbo
})$ fails, of course.

(1) holds:

Let $X_{1} \xbe \xdf (X_{2}),$ $X_{1} \xbe \xdi (X_{3}),$ we have to show
$X_{2} \xbe \xdi (X_{3}).$
If $X_{1}=X_{2},$ then this is trivial. Consider $X_{1} \xbe \xdf
(X_{2}).$
If $X_{1} \xEd X_{2},$ then $X_{1}$ has to be $\{a\}$ or
$\{a,b\}$ or $\{a,c\}.$ But none of these are in $ \xdi (X_{3})$ for any
$X_{3},$ so the
implication is trivially true.

(2) fails: $\{a,c\} \xbe \xdf (Y),$ $\{c\} \xbe \xdi (Y),$ $\{c\} \xce
\xdi (\{a,c\}).$

(3) fails: $\{a\} \xbe \xdf (X),$ $X \xbe \xdf (Y),$ $\{a\} \xce \xdf
(Y).$

(4) fails: $\{b\},\{c\} \xbe \xdi (Y),$ $\{c\} \xce \xdi (Y-\{b\})= \xdi
(\{a,c\})=\{ \xCQ \}.$

 \xDH

(1), (2), (4) hold, (3) fails:

Let $Y:=\{a,b,c\},$ $ \xdy:= \xdp (Y)-\{ \xCQ \},$ $ \xdf
(Y):=\{\{a,b\},$ $\{a,c\},$ $Y\}$

Let $ \xdf (\{a,b\}):=\{\{a\},\{a,b\}\},$ $ \xdf
(\{a,c\}):=\{\{a\},\{a,c\}\},$
and $ \xdf (Z):=\{Z\}$ for all other $Z.$

$ \xCf (Opt),$ $ \xCf (iM),$ $(eM \xdi ),$ $(eM \xdf )$ hold, $(I_{ \xbo
})$ fails, of course.

(1) holds:

Let $X_{1} \xbe \xdf (X_{2}),$ $X_{1} \xbe \xdi (X_{3}),$ we have to show
$X_{2} \xbe \xdi (X_{3}).$ Consider $X_{1} \xbe \xdi (X_{3}).$
If $X_{1}=X_{2},$ this is trivial. If $ \xCQ \xEd X_{1} \xbe \xdi
(X_{3}),$ then $X_{1}=\{b\}$ or $X_{1}=\{c\},$ but
then by $X_{1} \xbe \xdf (X_{2})$ $X_{2}$ has to be $\{b\},$ or $\{c\},$
so $X_{1}=X_{2}.$

(2) holds:
Let $X_{1} \xcc X_{2} \xcc X_{3},$ let $X_{2} \xbe \xdf (X_{3}),$ $X_{1}
\xbe \xdi (X_{3}),$ we have to show $X_{1} \xbe \xdi (X_{2}).$
If $X_{1}= \xCQ,$ this is trivial, likewise if $X_{2}=X_{3}.$ Otherwise
$X_{1}=\{b\}$ or $X_{1}=\{c\},$ and $X_{3}=Y.$ If $X_{1}=\{b\},$ then
$X_{2}=\{a,b\},$ and the condition
holds, likewise if $X_{1}=\{c\},$ then $X_{2}=\{a,c\},$ and it holds
again.

(3) fails: $\{a\} \xbe \xdf (\{a,c\}),$ $\{a,c\} \xbe \xdf (Y),$ $\{a\}
\xce \xdf (Y).$

(4) holds:

If $A,B \xbe \xdi (X),$ and $A \xEd B,$ $A,B \xEd \xCQ,$ then $X=Y$ and
e.g. $A=\{c\},$ $B=\{b\},$ and
$\{c\} \xbe \xdi (Y-\{b\})= \xdi (\{a,c\}).$

 \xEJ
$ \xcz $
\\[3ex]
\subsubsection{
Rational Monotony
}

\ee

\bfa

$\hspace{0.01em}$

% (+++ Orig. No.:  Fact M-plus-plus +++)

\label{Fact M-plus-plus}

The three versions of $( \xdm^{++})$ are equivalent.

(We assume closure of the domain under set difference.
For the third version of $( \xdm^{++}),$ we use $ \xCf (iM).)$

\efa

\subparagraph{
Proof
}

$\hspace{0.01em}$

% (+++ Orig.:  Proof +++)

For (1) and (2), we have $A,B \xcc X,$ for (3) we have $A \xcc X \xcc Y.$
For $A,B \xcc X,$ $(X-B)-((X-A)-B)=A-B$ holds.

$(1) \xch (2):$ Let $A \xbe \xdf (X),$ $B \xce \xdf (X),$ so $X-A \xbe
\xdi (X),$ so by prerequisite
$(X-A)-B \xbe \xdi (X-$B), so $A-B=(X-B)-((X-A)-B) \xbe \xdf (X-$B).

$(2) \xch (1):$ Let $A \xbe \xdi (X),$ $B \xce \xdf (X),$ so $X-A \xbe
\xdf (X),$ so by prerequisite
$(X-A)-B \xbe \xdf (X-$B), so $A-B=(X-B)-((X-A)-B) \xbe \xdi (X-$B).

$(1) \xch (3):$

Suppose $A \xce \xdm^{+}(Y),$ but $X \xbe \xdm^{+}(Y),$ we show $A \xce
\xdm^{+}(X).$ So $A \xbe \xdi (Y),$ $Y-X \xce \xdf (Y),$
so by (1) $A=A-(Y-X) \xbe \xdi (Y-(Y-X))= \xdi (X).$

$(3) \xch (1):$

Suppose $A-B \xce \xdi (X-$B), $B \xce \xdf (X),$ we show $A \xce \xdi
(X).$ By prerequisite $A-B \xbe \xdm^{+}(X-$B),
$X-B \xbe \xdm^{+}(X),$ so by (3) $A-B \xbe \xdm^{+}(X),$ so by $ \xCf
(iM)$ $A \xbe \xdm^{+}(X),$ so $A \xce \xdi (X).$

$ \xcz $
\\[3ex]

\bfa

$\hspace{0.01em}$

% (+++ Orig. No.:  Fact M-RatM +++)

\label{Fact M-RatM}

We assume that all sets are definable by a formula.

$ \xCf (RatM)$ $ \xcj $ $( \xdm^{++})$

\efa

\subparagraph{
Proof
}

$\hspace{0.01em}$

% (+++ Orig.:  Proof +++)

We show equivalence of $ \xCf (RatM)$ with version (1) of $( \xdm^{++}).$

`` $ \xch $ ''

We have $A,B \xcc X,$ so we can write
$X=M( \xbf ),$ $A=M( \xbf \xcu \xCN \xbq ),$ $B=M( \xbf \xcu \xCN \xbq '
).$ $A \xbe \xdi (X),$ $B \xce \xdf (X),$ so
$ \xbf \xcn \xbq,$ $ \xbf \xcN \xCN \xbq ',$ so by $ \xCf (RatM)$ $ \xbf
\xcu \xbq ' \xcn \xbq,$ so
$A-B=M( \xbf \xcu \xCN \xbq )-M( \xbf \xcu \xCN \xbq ' )=M( \xbf \xcu \xbq
' \xcu \xCN \xbq ) \xbe \xdi (M( \xbf \xcu \xbq ' ))= \xdi (X-$B).

`` $ \xci $ ''

Let $ \xbf \xcn \xbq,$ $ \xbf \xcN \xCN \xbq ',$ so $M( \xbf \xcu \xCN
\xbq ) \xbe \xdi (M( \xbf )),$ $M( \xbf \xcu \xCN \xbq ' ) \xce \xdf (M(
\xbf )),$ so
by $( \xdm^{++})$ (1) $M( \xbf \xcu \xbq ' \xcu \xCN \xbq )=M( \xbf \xcu
\xCN \xbq )-M( \xbf \xcu \xCN \xbq ' ) \xbe \xdi (M( \xbf \xcu \xbq ' )),$
so
$ \xbf \xcu \xbq ' \xcn \xbq.$

$ \xcz $
\\[3ex]
\subsection{
Size and principal filter logic
}
\label{Section Principal}

The connection with logical rules was shown in the table of
Definition \ref{Definition Log-Cond-Ref-Size} (page \pageref{Definition
Log-Cond-Ref-Size}).

(1) to (7) of the following proposition (in different notation, as the
more
systematic connections were found only afterwards) was already published
in
 \cite{GS08c}, we give it here in totality to complete the picture.
\index{Proposition Ref-Class-Mu-neu}

\bp

$\hspace{0.01em}$

% (+++ Orig. No.:  Proposition Ref-Class-Mu-neu +++)

\label{Proposition Ref-Class-Mu-neu}

If $f(X)$ is the smallest $ \xCf A$ s.t. $A \xbe \xdf (X),$ then, given
the property on the
left, the one on the right follows.

Conversely, when we define $ \xdf (X):=\{X':f(X) \xcc X' \xcc X\},$ given
the property on
the right, the one on the left follows. For this direction, we assume
that we can use the full powerset of some base set $U$ - as is the case
for
the model sets of a finite language. This is perhaps not too bold, as
we mainly want to stress here the intuitive connections, without putting
too much weight on definability questions.

We assume $ \xCf (iM)$ to hold.

{\footnotesize

\begin{tabular}{|c|c|c|c|}

\hline

(1.1)
\xEH
$(eM\xdi )$
\xEH
$ \xch $
\xEH
$( \xbm wOR)$
\xEP

\cline{1-1}
\cline{3-3}

(1.2)
\xEH
\xEH
$ \xci $
\xEH
\xEP

\hline

(2.1)
\xEH
$(eM\xdi )+(I_\xbo )$
\xEH
$ \xch $
\xEH
$( \xbm OR)$
\xEP

\cline{1-1}
\cline{3-3}

(2.2)
\xEH
\xEH
$ \xci $
\xEH
\xEP

\hline

(3.1)
\xEH
$(eM\xdi )+(I_\xbo )$
\xEH
$ \xch $
\xEH
$( \xbm PR)$
\xEP

\cline{1-1}
\cline{3-3}

(3.2)
\xEH
\xEH
$ \xci $
\xEH
\xEP

\hline

(4.1)
\xEH
$(I \xcv disj )$
\xEH
$ \xch $
\xEH
$( \xbm disjOR)$
\xEP

\cline{1-1}
\cline{3-3}

(4.2)
\xEH
\xEH
$ \xci $
\xEH
\xEP

\hline

(5.1)
\xEH
$(\xdm^+_\xbo) (4)$
\xEH
$ \xch $
\xEH
$( \xbm CM)$
\xEP

\cline{1-1}
\cline{3-3}

(5.2)
\xEH
\xEH
$ \xci $
\xEH
\xEP

\hline

(6.1)
\xEH
$(\xdm^{++})$
\xEH
$ \xch $
\xEH
$( \xbm RatM)$
\xEP

\cline{1-1}
\cline{3-3}

(6.2)
\xEH
\xEH
$ \xci $
\xEH
\xEP

\hline

(7.1)
\xEH
$(I_\xbo )$
\xEH
$ \xch $
\xEH
$( \xbm AND)$
\xEP

\cline{1-1}
\cline{3-3}

(7.2)
\xEH
\xEH
$ \xci $
\xEH
\xEP

\hline

(8.1)
\xEH
$(eM\xdi )+(I_\xbo )$
\xEH
$ \xch $
\xEH
$( \xbm CUT)$
\xEP

\cline{1-1}
\cline{3-3}

(8.2)
\xEH
\xEH
$ \xcI $
\xEH
\xEP

\hline

(9.1)
\xEH
$(eM\xdi )+(I_\xbo )+(\xdm^+_\xbo) (4)$
\xEH
$ \xch $
\xEH
$( \xbm CUM)$
\xEP

\cline{1-1}
\cline{3-3}

(9.2)
\xEH
\xEH
$ \xcI $
\xEH
\xEP

\hline

(10.1)
\xEH
$(eM\xdi )+(I_\xbo )+(eM\xdf )$
\xEH
$ \xch $
\xEH
$( \xbm \xcc \xcd)$
\xEP

\cline{1-1}
\cline{3-3}

(10.2)
\xEH
\xEH
$ \xcI $
\xEH
\xEP

\hline

\end{tabular}

}

\ep

Note that there is no $( \xbm wCM),$ as the conditions $( \xbm  \Xl.)$
imply that the
filter is principal,
and thus that $(I_{ \xbo })$ holds - we cannot ``see'' $ \xCf (wCM)$ alone
with
principal filters.
\index{Proposition Ref-Class-Mu-neu Proof}

\subparagraph{
Proof
}

$\hspace{0.01em}$

% (+++ Orig.:  Proof +++)

(1.1) $(eM \xdi )$ $ \xch $ $( \xbm wOR):$

$X-f(X)$ is small in $X,$ so it is small in $X \xcv Y$ by $(eM \xdi ),$ so
$A:=X \xcv Y-(X-f(X)) \xbe \xdf (X \xcv Y),$ but $A \xcc f(X) \xcv Y,$ and
$f(X \xcv Y)$ is the smallest element
of $ \xdf (X \xcv Y),$ so $f(X \xcv Y) \xcc A \xcc f(X) \xcv Y.$

(1.2) $( \xbm wOR)$ $ \xch $ $(eM \xdi ):$

Let $X \xcc Y,$ $X':=Y-$X. Let $A \xbe \xdi (X),$ so $X-A \xbe \xdf (X),$
so $f(X) \xcc X-$A, so
$f(X \xcv X' ) \xcc f(X) \xcv X' \xcc (X-A) \xcv X' $ by prerequisite, so
$(X \xcv X' )-((X-A) \xcv X' )=A \xbe \xdi (X \xcv X' ).$

(2.1) $(eM \xdi )+(I_{ \xbo })$ $ \xch $ $( \xbm OR):$

$X-f(X)$ is small in $X,$ $Y-f(Y)$ is small in $Y,$ so both are small in
$X \xcv Y$ by
$(eM \xdi ),$ so $A:=(X-f(X)) \xcv (Y-f(Y))$ is small in $X \xcv Y$ by
$(I_{ \xbo }),$ but
$X \xcv Y-(f(X) \xcv f(Y)) \xcc A,$ so $f(X) \xcv f(Y) \xbe \xdf (X \xcv
Y),$ so, as $f(X \xcv Y)$ is the smallest
element of $ \xdf (X \xcv Y),$ $f(X \xcv Y) \xcc f(X) \xcv f(Y).$

(2.2) $( \xbm OR)$ $ \xch $ $(eM \xdi )+(I_{ \xbo }):$

Let again $X \xcc Y,$ $X':=Y-$X. Let $A \xbe \xdi (X),$ so $X-A \xbe \xdf
(X),$ so $f(X) \xcc X-$A. $f(X' ) \xcc X',$
so $f(X \xcv X' ) \xcc f(X) \xcv f(X' ) \xcc (X-A) \xcv X' $ by
prerequisite, so
$(X \xcv X' )-((X-A) \xcv X' )=A \xbe \xdi (X \xcv X' ).$

$(I_{ \xbo })$ holds by definition.

(3.1) $(eM \xdi )+(I_{ \xbo })$ $ \xch $ $( \xbm PR):$

Let $X \xcc Y.$ $Y-f(Y)$ is the largest element of $ \xdi (Y),$ $X-f(X)
\xbe \xdi (X) \xcc \xdi (Y)$ by
$(eM \xdi ),$ so $(X-f(X)) \xcv (Y-f(Y)) \xbe \xdi (Y)$ by $(I_{ \xbo }),$
so by ``largest'' $X-f(X) \xcc Y-f(Y),$
so $f(Y) \xcs X \xcc f(X).$

(3.2) $( \xbm PR)$ $ \xch $ $(eM \xdi )+(I_{ \xbo })$

Let again $X \xcc Y,$ $X':=Y-$X. Let $A \xbe \xdi (X),$ so $X-A \xbe \xdf
(X),$ so $f(X) \xcc X-$A, so
by prerequisite $f(Y) \xcs X \xcc X-$A, so $f(Y) \xcc X' \xcv (X-$A), so
$(X \xcv X' )-(X' \xcv (X-A))=A \xbe \xdi (Y).$

Again, $(I_{ \xbo })$ holds by definition.

(4.1) $(I \xcv disj)$ $ \xch $ $( \xbm disjOR):$

If $X \xcs Y= \xCQ,$ then (1) $A \xbe \xdi (X),B \xbe \xdi (Y) \xch A
\xcv B \xbe \xdi (X \xcv Y)$ and
(2) $A \xbe \xdf (X),B \xbe \xdf (Y) \xch A \xcv B \xbe \xdf (X \xcv Y)$
are equivalent. (By $X \xcs Y= \xCQ,$
$(X-A) \xcv (Y-B)=(X \xcv Y)-(A \xcv B).)$
So $f(X) \xbe \xdf (X),$ $f(Y) \xbe \xdf (Y)$ $ \xch $ (by prerequisite)
$f(X) \xcv f(Y) \xbe \xdf (X \xcv Y).$ $f(X \xcv Y)$
is the smallest element of $ \xdf (X \xcv Y),$ so $f(X \xcv Y) \xcc f(X)
\xcv f(Y).$

(4.2) $( \xbm disjOR)$ $ \xch $ $(I \xcv disj):$

Let $X \xcc Y,$ $X':=Y-$X. Let $A \xbe \xdi (X),$ $A' \xbe \xdi (X' ),$
so $X-A \xbe \xdf (X),$ $X' -A' \xbe \xdf (X' ),$
so $f(X) \xcc X-$A, $f(X' ) \xcc X' -A',$ so $f(X \xcv X' ) \xcc f(X)
\xcv f(X' ) \xcc (X-A) \xcv (X' -A' )$ by
prerequisite, so $(X \xcv X' )-((X-A) \xcv (X' -A' ))=A \xcv A' \xbe \xdi
(X \xcv X' ).$

(5.1) $( \xdm^{+}_{ \xbo })$ $ \xch $ $( \xbm CM):$

$f(X) \xcc Y \xcc X$ $ \xch $ $X-Y \xbe \xdi (X),$ $X-f(X) \xbe \xdi (X)$
$ \xch $ (by $( \xdm^{+}_{ \xbo }),$ (4))
$A:=(X-f(X))-(X-Y) \xbe \xdi (Y)$ $ \xch $
$Y-A=f(X)-(X-Y) \xbe \xdf (Y)$ $ \xch $ $f(Y) \xcc f(X)-(X-Y) \xcc f(X).$

(5.2) $( \xbm CM)$ $ \xch $ $( \xdm^{+}_{ \xbo })$

Let $X-A \xbe \xdi (X),$ so $A \xbe \xdf (X),$ let $B \xbe \xdi (X),$ so
$f(X) \xcc X-B \xcc X,$ so by prerequisite
$f(X-B) \xcc f(X).$
As $A \xbe \xdf (X),$ $f(X) \xcc A,$ so $f(X-B) \xcc f(X) \xcc A \xcs
(X-B)=A-$B, and $A-B \xbe \xdf (X-$B), so
$(X-A)-B=X-(A \xcv B)=(X-B)-(A-B) \xbe \xdi (X-$B), so
$( \xdm^{+}_{ \xbo }),$ (4) holds.

(6.1) $( \xdm^{++})$ $ \xch $ $( \xbm RatM):$

Let $X \xcc Y,$ $X \xcs f(Y) \xEd \xCQ.$ If $Y-X \xbe \xdf (Y),$ then
$A:=(Y-X) \xcs f(Y) \xbe \xdf (Y),$ but by
$X \xcs f(Y) \xEd \xCQ $ $A \xcb f(Y),$ contradicting ``smallest'' of
$f(Y).$ So $Y-X \xce \xdf (Y),$ and
by $( \xdm^{++})$ $X-f(Y)=(Y-f(Y))-(Y-X) \xbe \xdi (X),$ so $X \xcs f(Y)
\xbe \xdf (X),$ so $f(X) \xcc f(Y) \xcs X.$

(6.2) $( \xbm RatM)$ $ \xch $ $( \xdm^{++})$

Let $A \xbe \xdf (Y),$ $B \xce \xdf (Y).$ $B \xce \xdf (Y)$ $ \xch $ $Y-B
\xce \xdi (Y)$ $ \xch $ $(Y-B) \xcs f(Y) \xEd \xCQ.$
Set $X:=Y-$B, so $X \xcs f(Y) \xEd \xCQ,$ $X \xcc Y,$ so $f(X) \xcc f(Y)
\xcs X$ by prerequisite.
$f(Y) \xcc A$ $ \xch $ $f(X) \xcc f(Y) \xcs X=f(Y)-B \xcc A-$B.

(7.1) $( \xdi_{ \xbo })$ $ \xch $ $( \xbm AND)$

Trivial.

(7.2) $( \xbm AND)$ $ \xch $ $( \xdi_{ \xbo })$

Trivial.

(8.1) Let $f(X) \xcc Y \xcc X.$ $Y-f(Y) \xbe \xdi (Y) \xcc \xdi (X)$ by
$(eM \xdi ).$ $f(X) \xcc Y$ $ \xch $
$X-Y \xcc X-f(X) \xbe \xdi (X),$ so by $ \xCf (iM)$ $X-Y \xbe \xdi (X).$
Thus by $(I_{ \xbo })$
$X-f(Y)=(X-Y) \xcv (Y-f(Y)) \xbe \xdi (X),$ so $f(Y) \xbe \xdf (X),$ so
$f(X) \xcc f(Y)$ by definition.

(8.2) $( \xbm CUT)$ is too special to allow to deduce $(eM \xdi ).$
Consider $U:=\{a,b,c\},$ $X:=\{a,b\},$ $ \xdf (X)=\{X,\{a\}\},$ $ \xdf
(Z)=\{Z\}$ for all other
$X \xEd Z \xcc U.$ Then $(eM \xdi )$ fails, as $\{b\} \xbe \xdi (X),$ but
$\{b\} \xce \xdi (U).$
$ \xCf (iM)$ and $(eM \xdf )$ hold. We have to check $f(A) \xcc B \xcc A
\xch f(A) \xcc f(B).$ The only
case where it might fail is $A=X,$ $B=\{a\},$ but it holds there, too.

(9.1) By Fact \ref{Fact Mu-Base} (page \pageref{Fact Mu-Base}),
published as Fact 14 in  \cite{GS08c}, (6),
we have $( \xbm CM)+( \xbm CUT) \xcj ( \xbm CUM),$ so the result follows
from (5.1) and (8.1).

(9.2) Consider the same example as in (8.2). $f(A) \xcc B \xcc A \xch
f(A)=f(B)$ holds
there, too, by the same argument as above.

(10.1) Let $f(X) \xcc Y,$ $f(Y) \xcc X.$ So $f(X),f(Y) \xcc X \xcs Y,$ and
$X-(X \xcs Y) \xbe \xdi (X),$
$Y-(X \xcs Y) \xbe \xdi (Y)$ by $ \xCf (iM).$ Thus $f(X),f(Y) \xbe \xdf (X
\xcs Y)$ by $(eM \xdf )$ and
$f(X) \xcs f(Y) \xbe \xdf (X \xcs Y)$ by $(I_{ \xbo }).$ So $X \xcs
Y-(f(X) \xcs f(Y)) \xbe \xdi (X \xcs Y),$ so
$X \xcs Y-(f(X) \xcs f(Y)) \xbe \xdi (X), \xdi (Y)$ by $(eM \xdi ),$ so
$(X-(X \xcs Y)) \xcv (X \xcs Y-f(X) \xcs f(Y))$ $=$ $X-f(X) \xcs f(Y) \xbe
\xdi (X)$ by $(I_{ \xbo }),$
so $f(X) \xcs f(Y) \xbe \xdf (X),$ likewise $f(X) \xcs f(Y) \xbe \xdf
(Y),$ so
$f(X) \xcc f(X) \xcs f(Y),$ $f(Y) \xcc f(X) \xcs f(Y),$ and $f(X)=f(Y).$

(10.2) Consider again the same example as in (8.2), we have to show
that $f(A) \xcc B,$ $f(B) \xcc A$ $ \xch $ $f(A)=f(B).$ The only
interesting case is when one
of $A,B$ is $X,$ but not both. Let e.g. $A=X.$ We then have $f(X)=\{a\},$
$f(B)=B \xcc X,$
and $f(X)=\{a\} \xcc B,$ so $B=\{a\},$ and the condition holds.

$ \xcz $
\\[3ex]

The product size defined by principal filters is discussed in
Section \ref{Section Product-Size} (page \pageref{Section Product-Size}).

% ******* BEGIN LATEX SOURCE FILE 4-5-pref.tex *******
%
% Uebers. aus Karltex File: 4-5-pref.m
%
%
\chapter{Preferential structures - Part I}
\label{Chapter Pref}

This chapter, Part $I,$ is dedicated to the basic case without
conditions for the domain.

The following chapter, Part II, will treat the case with supplementary
conditions for the domain, as well as applications and special
cases.

Higher preferential structures will be treated in the next but one
chapter.
\section{Introduction}
\label{Section 2.2.1}

After the present section, we will treat in
Section \ref{Section Without-Domain} (page \pageref{Section Without-Domain}) 
the case
without conditions on the domain,
and in
Section \ref{Section With-Domain} (page \pageref{Section With-Domain})  the case
with the usual conditions on the domain, in particular
closure under finite unions and finite intersections.

But, first some general remarks
\subsection{
Remarks on nonmonotonic logics and preferential semantics}
\label{Section NML-and-Pref}

Nonmonotonic logics were, historically, studied from two different points
of
view: the syntactic side, where rules like (AND), (CUM) (see below,
Definition \ref{Definition Log-Cond-Ref-Size} (page \pageref{Definition
Log-Cond-Ref-Size}) )
were postulated for their naturalness in reasoning,
and from the semantic side, by the introduction of preferential structures
(see Definition \ref{Definition Pref-Str} (page \pageref{Definition Pref-Str})
and Definition \ref{Definition Pref-Log} (page \pageref{Definition Pref-Log})
below). This work was done on the one hand side by
Gabbay  \cite{Gab85},
Makinson  \cite{Mak94}, and others, and for the second approach by
Shoham and others, see  \cite{Sho87b},  \cite{BS85}.
Both approaches were brought together by Kraus, Lehmann, Magidor and
others, see  \cite{KLM90},  \cite{LM92}, in their completeness
results.

A preferential structure $ \xdm $ defines a logic $ \xcn $ by $T \xcn \xbf
$ iff $ \xbf $ holds in all
$ \xdm -$minimal models of $T.$ This is made precise in
Definition \ref{Definition Pref-Str} (page \pageref{Definition Pref-Str})  and
Definition \ref{Definition Pref-Log} (page \pageref{Definition Pref-Log}) 
below.
At the same time, $ \xdm $ defines also a model set function, by assigning
to the
set of models of $T$ the set of its minimal models. As logics can speak
only
about definable model sets (here the model set defined by $T),$ $ \xdm $
defines
a function from the definable sets of models to arbitrary model sets:
$ \xbm_{ \xdm }: \xdD ( \xdl ) \xcp \xdp (M( \xdl )).$ This is the general
framework, within which we will
work most of the time. Different logics and situations (see e.g.
Plausibility
Logic, Section \ref{Section 2.2.9} (page \pageref{Section 2.2.9}), but also
update situations, etc.,
Section \ref{Section 2.3} (page \pageref{Section 2.3}) ), will force us to
generalize, we then
consider
functions $f: \xdy \xcp \xdp (W),$ where $W$ is an arbitrary set, and $
\xdy \xcc \xdp (W).$

This Chapter is about representation proofs in the realm of preferential
and
related structures and concerns mainly the following points:
\begin{enumerate}

\item the importance of closure properties of the domain, in particular
      under finite unions,
\item the conditions affected by lack of definability preservation,
\item the limit version of preferential structures,
\item the problems and solutions for ``hidden dimensions'', i.e. dimensions
      not directly observable.

\end{enumerate}

Concerning (1), the main new result is probably the fact that, in the
absence of
closure under finite unions, Cumulativity fans out to an infinity of
different
conditions. We also separate here clearly the main proof technique from
simplifications possible in the case of closure under finite unions.

Concerning (2), we examine in a systematic way conditions affected by
absence of definability preservation, and use now more carefully crafted
proofs
which can be used in the cases with and those without definability
preservation,
achieving thus a simplification and a conceptually clearer approach.

Concerning (3), we introduce the concept of an algebraic limit, to
separate
logical problems (i.e. due to lack of definability preservation) from
algebraic
ones of the limit variant. Again, this results in a better insight into
problems
and solutions of this variant.

Concerning (4), we describe a problem common to several representation
questions, where we cannot directly observe all dimensions of a result.

Conceptually, one can subsume (2) and (4) under a more general notion of
``blurred observation'', but the solutions appear to be sufficiently
different to
merit a separate treatment - at least at the present stage of
investigation.

Throughout the text, we emphasize those points in proofs where arbitrary
choices are made, where we go beyond that which we know, where we ``loose
ignorance'', as we think that an ideal proof is also in this aspect as
general
as possible, and avoids such losses, or, at least, reveals them clearly.
E.g., ignorance is lost when we complete in arbitrary manner a partial
order
to a complete one, or when we choose arbitrarily some copy to be smaller
than
suitable other ones. The authors thinks that this, perhaps seemingly
pedantic,
attention will reveal itself as fruitful in future work.

The present text is a continuation of the second author's  \cite{Sch04}. Many
results
presented here were already shown in  \cite{Sch04}, but in a less
systematic, more
ad hoc way. The emphasis of this book is on describing the problems and
machinery ``behind the scene'', on a clear separation of general problems
and
possible simplifications. In particular, we take a systematic look at
closure of the domain (mainly under finite unions), conditions affected by
lack of definability preservation, and an analysis of the limit variant,
leading to the notion of an algebraic limit.

Perhaps the single most important new result is the fact that Cumulativity
fans out to an infinity of diffferent conditions in the absence of closure
under finite unions
(see Example \ref{Example Inf-Cum-Alpha} (page \pageref{Example Inf-Cum-Alpha})
).

The systematic investigation of $H(U,u)$ and $H(U),$ see Definition
\ref{Definition HU-All} (page \pageref{Definition HU-All})
is also new.

Many proofs have been systematized, prerequisites have been better
isolated,
and the results are now more general, so they can be re-used more easily.
In particular, we often first look at the case without closure under
finite
union, and obtain in a second step the case with this closure as a
simplification of the former.

The cases of ``hidden dimensions'' are now looked at in a more systematic
way,
and infinite variants are discussed and solved for the first time - to the
authors' knowledge. In particular, such hidden dimensions are present in
update situations, where we only observe the outcome, and can conclude
about
starting and intermediate steps only indirectly. We will see that this
question
leads to an (it seems) non-trivial problem, and how to circumvent it in a
natural way, using the limit approach.

The separation of the limit variant into structural limit, algebraic
limit, and logical limit allows us to better identify problems, and also
see that and how a number of problems here are just those arising also
from lack
of definability preservation.

Finally, we solve some (to the authors' knowledge) open representation
problems: Problems around Aqvist's deontic logic, and Booth's revision
approach - see Section \ref{Section Aqvist} (page \pageref{Section Aqvist})  and
Section \ref{Section Booth} (page \pageref{Section Booth}).

The core of the completeness proofs consists in a general proof strategy,
which
is, essentially, a mathematical reformulation of the things we have to do.
For instance, in general preferential structures, if $x \xbe X- \xbm (X),$
we need to
``minimize'' $x$ by some other element. Yet, we do not know by which one.
This
leads in a natural way to consider copies of $x,$ one for each $x' \xbe
X,$ and to
minimize such $ \xBc x,x'  \xBe $ by $x' $ - somehow. Essentially, this is all
there
is to do
in the most general case. As $x$ may be in several $X- \xbm (X),$ we have
to do the
construction for all such $X$ simultaneously, leading us to consider the
product
and choice functions. This is the basic proof construction - the
mathematical
counterpart of the representation idea. Of course, the same has to be done
for
above $x',$ so we will have $ \xBc x,f \xBe,$ $ \xBc x',f'  \xBe,$ etc. Now,
there is a
problem: do
we make all such $ \xBc x',f'  \xBe $ smaller than the $ \xBc x,f \xBe,$ only
one, some of
them? We
simply do not know, they all give the same result in the basic
construction,
as we are basically interested only in the first coordinate - and will see
in
the outcome only this. Choosing one possibility is completely arbitrary,
and
such a choice is a loss of ignorance. An ideal proof should not do this,
it
should not commit beyond the necessary, it should ``preserve ignorance''.
There is
an easy way out: we just
construct for all possibilities one structure - just as a classical theory
may
have more than one model. This inflation of structures is the honest
solution.
Of course, once we are aware of the problem, we can construct just one
structure, but should remember the arbitrariness of our decision. This
attention is not purely academic or pedantic, as we will see immediately
when
we try to make the construction transitive. If we make all copies of $x' $
smaller
than $ \xBc x,f \xBe,$ then the structure cannot be made transitive by simple
closure: we
did not pay enough attention to what can happen ``in future''. The solution
is
not to consider simple functions for copies, but trees of elements and
functions, giving us complete control.

So far, we did not care about prerequisites, we just looked at the
necessary
construction.

This will be done in a second step, where we will see that sufficiently
strong
prerequisites, especially about closure of the domain of $ \xbm,$ in
particular
whether for $X,$ $X' $ in the domain $X \xcv X' $ also is in the domain,
can simplify the
construction considerably.

Thus, we have a core - the basic construction -, an initial part -
eventual
simplifications due to closure properties -, but are not yet finished: so
far,
we looked only at model sets (or, more generally, just arbitrary sets) and
the
algebraic properties of the $ \xbm -$functions to represent, and
we still have to translate our results to logic. As long as our $ \xbm
-$functions
preserve definability, i.e., $ \xbm (M(T)),$ $M(T)$ the set of models of
some theory $T,$
is again $M(T' )$ for some, usually other, theory $T',$ this is simple.
But, of
course, this is a very strong assumption in the general infinite case.
Usually,
$ \xbm (M(T))$ will only generate a theory, but it will not be all of the
models of
this theory - logic is too poor to see that something is missing. So our
observation of the result is blurred.

Thus, there is still a final part to consider. We have shown in  \cite{Sch04}
that
the general case is impossible to characterize by traditional means, and
we
need other means - we have to work with ``small'' model sets: our
observation
can be a bit off the true situation, but not too much. This, again, can be
seen as a problem of ignorance: we do not know the exact result, but only
that it is not too far off from what we see. Thus, characterization will
be exactly this: it gives a sort of upper limit of how far we can be off,
and anything within those limits is possible, it does not give an exact
result.

This last part will be considered in a more systematic way, too. In
particular,
we will look at conditions which might be affected, and at those which
will not
be affected. At the same time, we take care to make our basic construction
sufficiently general so we do not need that $ \xbm (X)$ is an element of
the domain
of $ \xbm.$

To summarize, we clearly distinguish here three parts: a core part with
the
essential construction, an initial part with possible simplifications when
domain closure conditions are sufficient, and a final part concerning the
sharpness of our observations.

A problem which is conceptually similar to the definability preservation
question is that of ``hidden dimensions''. E.g. in update situations (but
also
in situations of a not necessarily symmetric distance based revision), we
may
see only the result, the last dimension, and initial and intermediate
steps
are hidden from direct observation. More precisely, we may know that the
endpoints of preferred developments all are in some set $X,$ but have no
direct
information where intermediate points are. Again, our tools of observation
are
not sufficiently sharp to see the precise picture. This may generate
non-trivial
problems, especially when we may have infinite descending chains. A
natural
solution is to consider here (partly or totally) the limit approach, where
those things hold, which hold in the limit - i.e. the further we
``approach''
the (nonexisting) minimum, these properties will finally hold.

This is another subject of the present text:

We will introduce the concepts of
the structural, the algebraic, and the logical limit, and will see that
this
allows us to separate problems in this usually quite difficult case. Some
problems are simply due to the fact that a seemingly nice structural limit
does not have nice algebraic properties any more, so it should not be
considered. So, to have a ``good'' limit, the limit should not only capture
the
idea of a structural limit, but its algebraic counterpart should also
capture
the essential algebraic properties of the minimal choice functions.
Other problems are due to the fact that the nice algebraic limit
does not translate to nice logical properties, and we will see that this
is
often due to the same problems we saw in the absence of definability
preservation.

Thus, in a way, we come back in a cycle to the same problems again. This
is
one of the reasons the book form seems adequate for our results and
problems:
they are often interconnected, and a unified presentation seems the best.

It might be useful to emphasize the parallel investigation in the minimal
and
the limit variant:

For minimal preferential structures, we have:

 \xEI
 \xDH logical laws or descriptions like $ \xba \xcn \xba $ - they are the
(imperfect - by
definability preservation problems) reflection of the abstract
description,
 \xDH abstract or algebraic semantics, like $ \xbm (A) \xcc A$ - they are
the abstract
description of the foundation,
 \xDH structural semantics - they are the intuitive foundation.
 \xEJ

Likewise, for the limit situation, we have:

 \xEI
 \xDH structural limits - they are again the foundation,
 \xDH resulting abstract behaviour, which, again, has to be an abstract or
algebraic limit, resulting from the structural limit,
 \xDH a logical limit, which reflects the abstract limit, and may be
plagued by
definability preservation problems etc. when going from the model to the
logics side.
 \xEJ

Note that these clear distinctions have some philosophical importance,
too.
The structures need an intuitive or philosophical justification, why do we
describe preference by transitive relations, why do we admit copies, etc.?
The resulting algebraic choice functions are void of such questions.
\subsection{Basic definitions}
\label{Section 2.2.1.1}

The following two definitions make preferential structures precise. We
first
give the algebraic definition, and then the definition of the consequence
relation generated by an preferential structure. In the algebraic
definition,
the set $U$ is an arbitrary set, in the application to logic, this will be
the
set of classical models of the underlying propositional language.

In both cases, we first present the simpler variant without copies, and
then
the one with copies. (Note that e.g.  \cite{KLM90},  \cite{LM92} use labelling
functions instead, the
version without copies corresponds to injective labelling functions, the
one
with copies to the general case. These are just different ways of
speaking.) We
will discuss the difference between the
version without and the version with copies below, where we show that the
version with copies is strictly more expressive than the version without
copies,
and that transitivity of the relation adds new properties in the case
without
copies. When we summarize our own results below (see
Section \ref{Section 2.2.3} (page \pageref{Section 2.2.3}) ), we will
mention that, in the general case with copies, transitivity can be added
without changing properties.

We give here the ``minimal version'', the much more complicated ``limit
version''
is presented and discussed in
Section \ref{Section Limit} (page \pageref{Section Limit}).
Recall the intuition that the relation
$ \xeb $ expresses ``normality'' or ``importance'' - the $ \xeb -$smaller, the
more normal or
important. The smallest elements are those which count.
\index{Definition Pref-Str}

\bd

$\hspace{0.01em}$

% (+++ Orig. No.:  Definition Pref-Str +++)

\label{Definition Pref-Str}

Fix $U \xEd \xCQ,$ and consider arbitrary $X.$
Note that this $X$ has not necessarily anything to do with $U,$ or $ \xdu
$ below.
Thus, the functions $ \xbm_{ \xdm }$ below are in principle functions from
$V$ to $V$ - where $V$
is the set theoretical universe we work in.

Note that we work here often with copies of elements (or models).
In other areas of logic, most authors work with valuation functions. Both
definitions - copies or valuation functions - are equivalent, a copy
$ \xBc x,i \xBe $ can be seen as a state $ \xBc x,i \xBe $ with valuation $x.$
In the
beginning
of research on preferential structures, the notion of copies was widely
used, whereas e.g. [KLM90] used that of valuation functions. There is
perhaps
a weak justification of the former terminology. In modal logic, even if
two states have the same valid classical formulas, they might still be
distinguishable by their valid modal formulas. But this depends on the
fact
that modality is in the object language. In most work on preferential
stuctures, the consequence relation is outside the object language, so
different states with same valuation are in a stronger sense copies of
each other.

\ed

 \xEh

 \xDH Preferential models or structures.

 \xEh

 \xDH The version without copies:

A pair $ \xdm:= \xBc U, \xeb  \xBe $ with $U$ an arbitrary set, and $ \xeb $ an
arbitrary binary relation
on $U$ is called a preferential model or structure.

 \xDH The version with copies:

A pair $ \xdm:= \xBc  \xdu, \xeb  \xBe $ with $ \xdu $ an arbitrary set of
pairs,
and $ \xeb $ an arbitrary binary
relation on $ \xdu $ is called a preferential model or structure.

If $ \xBc x,i \xBe  \xbe \xdu,$ then $x$ is intended to be an element of $U,$
and
$i$ the index of the
copy.

We sometimes also need copies of the relation $ \xeb,$ we will then
replace $ \xeb $
by one or several arrows $ \xba $ attacking non-minimal elements, e.g. $x
\xeb y$ will
be written $ \xba:x \xcp y,$ $ \xBc x,i \xBe  \xeb  \xBc y,i \xBe $ will be
written $ \xba
: \xBc x,i \xBe  \xcp  \xBc y,i \xBe,$ and
finally we might have $ \xBc  \xba,k \xBe :x \xcp y$ and $ \xBc  \xba,k \xBe :
\xBc x,i \xBe  \xcp
\xBc y,i \xBe,$ etc.

 \xEj

 \xDH Minimal elements, the functions $ \xbm_{ \xdm }$

 \xEh

 \xDH The version without copies:

Let $ \xdm:= \xBc U, \xeb  \xBe,$ and define

$ \xbm_{ \xdm }(X)$ $:=$ $\{x \xbe X:$ $x \xbe U$ $ \xcu $ $ \xCN \xcE x'
\xbe X \xcs U.x' \xeb x\}.$

$ \xbm_{ \xdm }(X)$ is called the set of minimal elements of $X$ (in $
\xdm ).$

Thus, $ \xbm_{ \xdm }(X)$ is the set of elements such that there is no
smaller one
in $X.$

 \xDH The version with copies:

Let $ \xdm:= \xBc  \xdu, \xeb  \xBe $ be as above. Define

$ \xbm_{ \xdm }(X)$ $:=$ $\{x \xbe X:$ $ \xcE  \xBc x,i \xBe  \xbe \xdu. \xCN
\xcE
\xBc x',i'  \xBe  \xbe \xdu (x' \xbe X$ $ \xcu $ $ \xBc x',i'  \xBe ' \xeb  \xBc
x,i \xBe )\}.$

Thus, $ \xbm_{ \xdm }(X)$ is the projection on the first coordinate of the
set of elements
such that there is no smaller one in $X.$

Again, by abuse of language, we say that $ \xbm_{ \xdm }(X)$ is the set of
minimal elements
of $X$ in the structure. If the context is clear, we will also write just
$ \xbm.$

We sometimes say that $ \xBc x,i \xBe $ ``kills'' or ``minimizes'' $ \xBc y,j
\xBe $ if
$ \xBc x,i \xBe  \xeb  \xBc y,j \xBe.$ By abuse of language we also say a set
$X$ kills or
minimizes a set
$Y$ if for all $ \xBc y,j \xBe  \xbe \xdu,$ $y \xbe Y$ there is $ \xBc x,i \xBe 
\xbe \xdu,$
$x \xbe X$ s.t. $ \xBc x,i \xBe  \xeb  \xBc y,j \xBe.$

$ \xdm $ is also called injective or 1-copy, iff there is always at most
one copy
$ \xBc x,i \xBe $ for each $x.$ Note that the existence of copies corresponds to
a
non-injective labelling function - as is often used in nonclassical
logic, e.g. modal logic.

 \xEj

 \xEj

We say that $ \xdm $ is transitive, irreflexive, etc., iff $ \xeb $ is.

Note that $ \xbm (X)$ might well be empty, even if $X$ is not.
\index{Definition Pref-Log}

\bd

$\hspace{0.01em}$

% (+++ Orig. No.:  Definition Pref-Log +++)

\label{Definition Pref-Log}

We define the consequence relation of a preferential structure for a
given propositional language $ \xdl.$

 \xEh

 \xDH

 \xEh

 \xDH If $m$ is a classical model of a language $ \xdl,$ we say by abuse
of language

$ \xBc m,i \xBe  \xcm \xbf $ iff $m \xcm \xbf,$

and if $X$ is a set of such pairs, that

$X \xcm \xbf $ iff for all $ \xBc m,i \xBe  \xbe X$ $m \xcm \xbf.$

 \xDH If $ \xdm $ is a preferential structure, and $X$ is a set of $ \xdl
-$models for a
classical propositional language $ \xdl,$ or a set of pairs $ \xBc m,i \xBe,$
where the $m$ are
such models, we call $ \xdm $ a classical preferential structure or model.

 \xEj

 \xDH

Validity in a preferential structure, or the semantical consequence
relation
defined by such a structure:

Let $ \xdm $ be as above.

We define:

$T \xcm_{ \xdm } \xbf $ iff $ \xbm_{ \xdm }(M(T)) \xcm \xbf,$ i.e. $
\xbm_{ \xdm }(M(T)) \xcc M( \xbf ).$

$ \xdm $ will be called definability preserving iff for all $X \xbe \xdD_{
\xdl }$ $ \xbm_{ \xdm }(X) \xbe \xdD_{ \xdl }.$

 \xEj

As $ \xbm_{ \xdm }$ is defined on $ \xdD_{ \xdl },$ but need by no means
always result in some new
definable set, this is (and reveals itself as a quite strong) additional
property.
\index{Example NeedCopies}

\ed

\be

$\hspace{0.01em}$

% (+++ Orig. No.:  Example NeedCopies +++)

\label{Example NeedCopies}

This simple example illustrates the
importance of copies. Such examples seem to have appeared for the first
time
in print in  \cite{KLM90}, but can probably be
attibuted to folklore.

Consider the propositional language $ \xdl $ of two propositional
variables $p,q$, and
the classical preferential model $ \xdm $ defined by

$m \xcm p \xcu q,$ $m' \xcm p \xcu q,$ $m_{2} \xcm \xCN p \xcu q,$ $m_{3}
\xcm \xCN p \xcu \xCN q,$ with $m_{2} \xeb m$, $m_{3} \xeb m' $, and
let $ \xcm_{ \xdm }$ be its consequence relation. (m and $m' $ are
logically identical.)

Obviously, $Th(m) \xco \{ \xCN p\} \xcm_{ \xdm } \xCN p$, but there is no
complete theory $T' $ s.t.
$Th(m) \xco T' \xcm_{ \xdm } \xCN p$. (If there were one, $T' $ would
correspond to $m$, $m_{2},$ $m_{3},$
or the missing $m_{4} \xcm p \xcu \xCN q$, but we need two models to kill
all copies of $m.)$
On the other hand, if there were just one copy of $m,$ then one other
model,
i.e. a complete theory would suffice. More formally, if we admit at most
one
copy of each model in a structure $ \xdm,$ $m \xcM T,$ and $Th(m) \xco T
\xcm_{ \xdm } \xbf $ for some $ \xbf $ s.t.
$m \xcm \xCN \xbf $ - i.e. $m$ is not minimal in the models of $Th(m) \xco
T$ - then there is a
complete $T' $ with $T' \xcl T$ and $Th(m) \xco T' \xcm_{ \xdm } \xbf $,
i.e. there is $m'' $ with $m'' \xcm T' $ and
$m'' \xeb m.$ $ \xcz $
\\[3ex]

\ee

We define now two additional properties of the relation, smoothness and
rankedness.
\index{Definition Smooth}

\bd

$\hspace{0.01em}$

% (+++ Orig. No.:  Definition Smooth +++)

\label{Definition Smooth}

Let $ \xdy \xcc \xdp (U).$ (In applications to logic, $ \xdy $ will be $
\xdD_{ \xdl }.)$

A preferential structure $ \xdm $ is called $ \xdy -$smooth iff for every
$X \xbe \xdy $ every
element
$x \xbe X$ is either minimal in $X$ or above an element, which is minimal
in $X.$ More
precisely:

 \xEh

 \xDH The version without copies:

If $x \xbe X \xbe \xdy,$ then either $x \xbe \xbm (X)$ or there is $x'
\xbe \xbm (X).x' \xeb x.$

 \xDH The version with copies:

If $x \xbe X \xbe \xdy,$ and $ \xBc x,i \xBe  \xbe \xdu,$ then either there is
no
$ \xBc x',i'  \xBe  \xbe \xdu,$ $x' \xbe X,$
$ \xBc x',i'  \xBe  \xeb  \xBc x,i \xBe $ or there is $ \xBc x',i'  \xBe  \xbe
\xdu,$
$ \xBc x',i'  \xBe  \xeb  \xBc x,i \xBe,$ $x' \xbe X,$ s.t. there is
no $ \xBc x'',i''  \xBe  \xbe \xdu,$ $x'' \xbe X,$
with $ \xBc x'',i''  \xBe  \xeb  \xBc x',i'  \xBe.$

 \xEj

When considering the models of a language $ \xdl,$ $ \xdm $ will be
called smooth iff
it is $ \xdD_{ \xdl }-$smooth; $ \xdD_{ \xdl }$ is the default.

Obviously, the richer the set $ \xdy $ is, the stronger the condition $
\xdy -$smoothness
will be.
\index{Fact Rank-Base}

\ed

\bfa

$\hspace{0.01em}$

% (+++ Orig. No.:  Fact Rank-Base +++)

\label{Fact Rank-Base}

Let $ \xeb $ be an irreflexive, binary relation on $X,$ then the following
two conditions
are equivalent:

(1) There is $ \xbO $ and an irreflexive, total, binary relation $ \xeb '
$ on $ \xbO $ and a
function $f:X \xcp \xbO $ s.t. $x \xeb y$ $ \xcj $ $f(x) \xeb ' f(y)$ for
all $x,y \xbe X.$

(2) Let $x,y,z \xbe X$ and $x \xcT y$ wrt. $ \xeb $ (i.e. neither $x \xeb
y$ nor $y \xeb x),$ then $z \xeb x$ $ \xch $ $z \xeb y$
and $x \xeb z$ $ \xch $ $y \xeb z.$
\index{Fact Rank-Base Proof}

\efa

\subparagraph{
Proof
}

$\hspace{0.01em}$

% (+++ Orig.:  Proof +++)

$(1) \xch (2)$: Let $x \xcT y,$ thus neither $fx \xeb ' fy$ nor $fy \xeb
' fx,$ but then $fx=fy.$ Let now
$z \xeb x,$ so $fz \xeb ' fx=fy,$ so $z \xeb y.$ $x \xeb z$ $ \xch $ $y
\xeb z$ is similar.

$(2) \xch (1)$: For $x \xbe X$ let $[x]:=\{x' \xbe X:x \xcT x' \},$ and $
\xbO:=\{[x]:x \xbe X\}.$ For $[x],[y] \xbe \xbO $
let $[x] \xeb ' [y]: \xcj x \xeb y.$ This is well-defined: Let $x \xcT x'
,$ $y \xcT y' $ and $x \xeb y,$ then
$x \xeb y' $ and $x' \xeb y'.$ Obviously, $ \xeb ' $ is an irreflexive,
total binary relation.
Define $f:X \xcp \xbO $ by $fx:=[x],$ then $x \xeb y \xcj [x] \xeb ' [y]
\xcj fx \xeb ' fy.$ $ \xcz $
\\[3ex]
\index{Definition Rank-Rel}

\bd

$\hspace{0.01em}$

% (+++ Orig. No.:  Definition Rank-Rel +++)

\label{Definition Rank-Rel}

We call an irreflexive, binary relation $ \xeb $ on $X,$ which satisfies
(1)
(equivalently (2)) of Fact \ref{Fact Rank-Base} (page \pageref{Fact Rank-Base})
, ranked.
By abuse of language, we also call a preferential structure $ \xBc X, \xeb  \xBe
$
ranked, iff
$ \xeb $ is.

\ed

The first condition says that if $x \xbe X$ is not a minimal element
of $X,$ then there is $x' \xbe \xbm (X)$ $x' \xeb x.$
In the finite case without copies, smoothness is a trivial consequence of
transitivity and lack of cycles. But note that in the other cases infinite
descending chains might still exist, even if the smoothness condition
holds,
they are just ``short-circuited'': we might have such chains, but below
every
element in the chain is a minimal element. In the authors' opinion,
smoothness
is difficult to justify as a structural property (or, in a more
philosophical
spirit, as a property of the world): why should we always have such
minimal
elements below non-minimal ones? Smoothness has, however, a justification
from
its consequences. Its attractiveness comes from two sides:

First, it generates a very valuable logical property, cumulativity (CUM):
If $ \xdm $ is smooth, and $ \ol{ \ol{T} }$
is the set of $ \xcm_{ \xdm }-consequences,$ then for $T \xcc \ol{T' }
\xcc \ol{ \ol{T} }$ $ \xch $ $ \ol{ \ol{T} }= \ol{ \ol{T' } }.$

Second, for certain approaches, it facilitates completeness proofs, as we
can
look directly at ``ideal'' elements, without having to bother about
intermediate
stages. See in particular the work by Lehmann and his co-authors,  \cite{KLM90},
 \cite{LM92}.

``Smoothness'', or, as it is also called,
``stopperedness'' seems - in the authors' opinion - a misnamer. We think it
should
better be called something like ``weak transitivity'': consider the case
where
$a \xee b \xee c,$ but $c \xeB a,$ with $c \xbe \xbm (X).$ It is then not
necessarily the case that
$a \xee c,$ but there is $c' $ ``sufficiently close to $c$ '', i.e. in $
\xbm (X),$ s.t. $a \xee c'.$
Results and proof techniques underline this idea. First, in the general
case
with copies, and in the smooth case (in the presence of $( \xcv )!),$
transitivity
does not add new properties, it is ``already present'', second, the
construction
of smoothness by sequences $ \xbs $ (see below in
Section \ref{Section 2.2.4} (page \pageref{Section 2.2.4}) )
is very close in spirit to a transitive construction.

The second condition, rankedness, seems easier to justify already as a
property
of the structure. It says that, essentially, the elements are ordered in
layers:
If a and $b$ are not comparable, then they are in the same layer. So, if
$c$ is
above (below) a, it will also be above (below) $b$ - like pancakes or
geological
strata. Apart from the triangle
inequality (and leaving aside cardinality questions), this is then just a
distance from some imaginary, ideal point. Again, this property has
important
consequences on the resulting model choice functions and consequence
relations,
making proof techniques for the non-ranked and the ranked case very
different.

$ \xdy $ can have certain properties, in classical propositional logic for
instance,
if $ \xdy $ is the set of formula defined model sets, then $ \xdy $ is
closed under
complements, finite unions and finite intersection. If $ \xdy $ is the set
of
theory defined model sets, $ \xdy $ is closed under finite unions,
arbitrary
intersections, but not complements any more.

The careful consideration of closure conditions of the domain was
motivated by Lehmann's Plausibility Logic, see  \cite{Leh92a},
and re-motivated by the work of Arieli and Avron, see  \cite{AA00}.
In both cases, the language does not have a
built-in ``or'' - resulting in absence $( \xcv )$ of the domain.

When trying to show completeness of Lehmann's system, the second author
noted
the
importance of the closure of the domain under $( \xcv ),$ see  \cite{Sch96-3}.
The work
of Arieli and Avron incited him to look at this property in a more
systematic
way which lead to the discovery of
Example \ref{Example Inf-Cum-Alpha} (page \pageref{Example Inf-Cum-Alpha}), and
thus of the
enormous strength of closure of the domain under finite unions, and, more
generally, of the importance of domain closure conditions.

In the resulting completeness proofs again, a strategy of ``divide and
conquer''
is useful. This helps us to unify
(or extend) our past completeness proofs for the smooth case
in the following way: We will identify more clearly than in the past a
more
or less simple algebraic property - $ \xCf (HU),$ $ \xCf (HU,u)$ etc. -
which allows us to split
the proofs into two parts. The first part (see
Section \ref{Section Pref-Details} (page \pageref{Section Pref-Details}) ) shows
validity of
the property, and
this demonstration depends on closure properties, the second part shows
how to
construct a representing structure using the algebraic property. This
second
part will be totally independent from closure properties, and is
essentially an
``administrative'' way to use the property for a construction. This split
approach
allows us thus to isolate the demonstration of the used property from the
construction itself, bringing both parts clearer to light, and simplifying
the proofs, by using common parts.

The reader will see that the successively more complicated conditions $
\xCf (HU),$
$ \xCf (HU,u)$ reflect well the successively more complicated situations
of
representation:

$ \xCf (HU):$ smooth (and transitive) structures in the presence of $(
\xcv ),$

$ \xCf (HU,u):$ smooth structures in the absence of $( \xcv ),$

This comparison becomes clearer when we see that in the final, most
complicated
case, we will have to carry around all the history of minimization,
$ \xBc Y_{0}, \Xl,Y_{n} \xBe,$
necessary for transitivity, which could be summarized in the first case
with
finite unions. Thus, from an abstract point of view, it is a very natural
development.

In the rest of this
Section \ref{Section 2.2.1} (page \pageref{Section 2.2.1}), we will only
describe the
problems to solve, without giving
a solution. This will be done in the next sections. Moreover, we will
asume
that we have precise knowledge of $f,$ i.e. what we see as $f(X)$ for $X
\xbe \xdy $ is
really the result, and not some approximation - as we will permit later,
in
Section \ref{Section 2.2.10} (page \pageref{Section 2.2.10}).

So this part is a leisurely description of problems and things to do.
We start with the most general case, arbitrary preferential structures,
turn
to transitive such structures, then to smooth, then to smooth transitive
ones,
and conclude by ranked and ranked smooth structures.

Throughout, we will try to preserve ignorance, i.e. not assume anything we
are
not forced to assume. This will become clearer in a moment. Once we have
understood the problem, we will sometimes just gloss over it by choosing
one
solution, but we should always be conscious that there is a problem.

We will consider here choice functions $f: \xdy \xcp \xdp (W),$ where $
\xdy \xcc \xdp (W),$ and
the problems to represent them by various preferential structures.

We will see the following basic representation problems and the
constructions to
solve them, i.e. to find representing structures for
% #[
% #[
% &' (1) General preferential structures
% &' (1) General preferential structures
% &' (2) General transitive preferential structures
% &' (2) General transitive preferential structures
% &' (3) Smooth preferential structures
% &' (3) Smooth preferential structures
% &' (4) Smooth transitive preferential structures
% &' (4) Smooth transitive preferential structures
% &' (5) Ranked preferential structures
% &' (5) Ranked preferential structures
% &' (6) Smooth ranked preferential structures
% &' (6) Smooth ranked preferential structures
% #]
% #]
\begin{enumerate}

\item General preferential structures
\item General transitive preferential structures
\item Smooth preferential structures
\item Smooth transitive preferential structures
\item Ranked preferential structures
\item Smooth ranked preferential structures

\end{enumerate}

The problems of and solutions to the ranked case are quite different from
the
first four cases. In particular, the situation when $ \xCQ \xEd U \xbe
\xdy,$ but $f(U)= \xCQ $ does
not present major difficulties in cases (1) - (4), but is quite nasty in
the
last case.
\subsubsection{
The situation
}

We work in some universe $W,$ there is a function $f: \xdy \xcp \xdp (W),$
where $ \xdy \xcc \xdp (W),$
$f$ will have certain properties, and perhaps $ \xdy,$ too, and we will
try to
represent $f$ by a preferential structure $ \xdz $ of a certain type, i.e.
we want
$f= \xbm_{ \xdz },$ with $ \xbm_{ \xdz }$ the $ \xbm -$function or choice
function of a preferential
structure $ \xdz.$ Note that the codomain of $f$ is not necessarily a
subset of $ \xdy $ -
so we have to pay attention not to apply $f$ twice.

Before we go into details, we give now an overview of the results.
\index{Proposition Pref-Representation-With-Ref}

The following table summarizes representation by
preferential structures. The positive implications on the right are shown
in
Proposition \ref{Proposition Alg-Log} (page \pageref{Proposition Alg-Log}) 
(going via the $ \xbm -$functions),
those on the left
are shown in the respective representation theorems.

``singletons'' means that the domain must contain all singletons, ``1 copy''
or `` $ \xcg 1$ copy'' means that the structure may contain only 1 copy for
each point,
or several, `` $( \xbm \xCQ )$ '' etc. for the preferential structure mean
that the
$ \xbm -$function of the structure has to satisfy this property.

We call a characterization ``normal'' iff it is a
universally quantified boolean combination (of any fixed, but perhaps
infinite, length) of rules of the usual form. We do not go into
details here.

In the second column from the left `` $ \xch $ '' means, for instance for
the
smooth case, that for any $ \xdy $ closed under finite unions, and any
choice function
$f$ which satisfies the conditions
in the left hand column, there is a (here $ \xdy -$smooth) preferential
structure $ \xdx $ which represents it, i.e. for all $Y \xbe \xdy $ $f(Y)=
\xbm_{ \xdx }(Y),$ where
$ \xbm_{ \xdx }$ is the model choice function of the structure $ \xdx.$
The inverse arrow $ \xci $ means that the model choice function for any
smooth $ \xdx $
defined on such $ \xdy $ will satisfy the conditions on the left.
\label{Proposition Pref-Representation-With-Ref}

{\scriptsize

% \begin{tabular*}{15.9cm}{|c|c|c|c|c|}
\begin{tabular}{|c|c|c|c|c|}

\hline

$\xbm-$ function
\xEH
\xEH
Pref.Structure
\xEH
\xEH
Logic
\xEP

\hline

$(\xbm \xcc)$
\xEH
$\xcj$
\xEH
reactive
\xEH
$\xcj$
\xEH
$(LLE)+(CCL)+$
\xEP

\xEH
Proposition \ref{Proposition Eta-Rho-Repres}
\xEH
\xEH
Proposition \ref{Proposition Higher-Repr}
\xEH
$(SC)$
\xEP

\xEH
page \pageref{Proposition Eta-Rho-Repres}
\xEH
\xEH
page \pageref{Proposition Higher-Repr}
\xEH
\xEP

\hline

$(\xbm \xcc)+(\xbm CUM)$
\xEH
$\xch$ $(\xcs)$
\xEH
reactive +
\xEH
\xEH
\xEP

\xEH
Proposition \ref{Proposition Level-3-Repr}
\xEH
essentially smooth
\xEH
\xEH
\xEP

\xEH
page \pageref{Proposition Level-3-Repr}
\xEH
\xEH
\xEH
\xEP

\hline

$(\xbm \xcc)+(\xbm \xcc \xcd)$
\xEH
$\xch$
\xEH
reactive +
\xEH
$\xcj$
\xEH
$(LLE)+(CCL)+$
\xEP

\xEH
Proposition \ref{Proposition Level-3-Repr}
\xEH
essentially smooth
\xEH
Proposition \ref{Proposition Higher-Repr}
\xEH
$(SC)+(\xcc \xcd)$
\xEP

\xEH
page \pageref{Proposition Level-3-Repr}
\xEH
\xEH
page \pageref{Proposition Higher-Repr}
\xEH
\xEP

\hline

$(\xbm \xcc)+(\xbm CUM)+(\xbm \xcc \xcd)$
\xEH
$\xci$
\xEH
reactive +
\xEH
\xEH
\xEP

\xEH
Fact \ref{Fact X-Sub-X'}
\xEH
essentially smooth
\xEH
\xEH
\xEP

\xEH
page \pageref{Fact X-Sub-X'}
\xEH
\xEH
\xEH
\xEP

\hline

$(\xbm \xcc)+(\xbm PR)$
\xEH
$\xci$
\xEH
general
\xEH
$\xch$ $(\xbm dp)$
\xEH
$(LLE)+(RW)+$
\xEP

\xEH
Fact \ref{Fact Pref-Sound}
\xEH
\xEH
\xEH
$(SC)+(PR)$
\xEP

\xEH
page \pageref{Fact Pref-Sound}
\xEH
\xEH
\xEH
\xEP

\cline{2-2}
\cline{4-4}

\xEH
$\xch$
\xEH
\xEH
$\xci$
\xEH
\xEP

\xEH
Proposition \ref{Proposition Pref-Complete}
\xEH
\xEH
\xEH
\xEP

\xEH
page \pageref{Proposition Pref-Complete}
\xEH
\xEH
\xEH
\xEP

\cline{2-2}
\cline{4-4}

\xEH
\xEH
\xEH
$\xcH$ without $(\xbm dp)$
\xEH
\xEP

\xEH
\xEH
\xEH
Example \ref{Example Pref-Dp}
\xEH
\xEP

\xEH
\xEH
\xEH
page \pageref{Example Pref-Dp}
\xEH
\xEP

\cline{4-5}

\xEH
\xEH
\xEH
$\xcJ$ without $(\xbm dp)$
\xEH
any ``normal''
\xEP

\xEH
\xEH
\xEH
Proposition 5.2.15
\xEH
characterization
\xEP

\xEH
\xEH
\xEH
in \cite{Sch04}
\xEH
of any size
\xEP

\hline

$(\xbm \xcc)+(\xbm PR)$
\xEH
$\xci$
\xEH
transitive
\xEH
$\xch$ $(\xbm dp)$
\xEH
$(LLE)+(RW)+$
\xEP

\xEH
Fact \ref{Fact Pref-Sound}
\xEH
\xEH
\xEH
$(SC)+(PR)$
\xEP

\xEH
page \pageref{Fact Pref-Sound}
\xEH
\xEH
\xEH
\xEP

\cline{2-2}
\cline{4-4}

\xEH
$\xch$
\xEH
\xEH
$\xci$
\xEH
\xEP

\xEH
Proposition \ref{Proposition Pref-Complete-Trans}
\xEH
\xEH
\xEH
\xEP

\xEH
page \pageref{Proposition Pref-Complete-Trans}
\xEH
\xEH
\xEH
\xEP

\cline{2-2}
\cline{4-4}

\xEH
\xEH
\xEH
$\xcH$ without $(\xbm dp)$
\xEH
\xEP

\xEH
\xEH
\xEH
Example \ref{Example Pref-Dp}
\xEH
\xEP

\xEH
\xEH
\xEH
page \pageref{Example Pref-Dp}
\xEH
\xEP

\cline{4-5}

\xEH
\xEH
\xEH
$\xcj$ without $(\xbm dp)$
\xEH
using ``small''
\xEP

\xEH
\xEH
\xEH
Proposition 5.2.5, 5.2.11
\xEH
exception sets
\xEP

\xEH
\xEH
\xEH
in \cite{Sch04}
\xEH
\xEP

\hline

$(\xbm \xcc)+(\xbm PR)+(\xbm CUM)$
\xEH
$\xci$
\xEH
smooth
\xEH
$\xch$ $(\xbm dp)$
\xEH
$(LLE)+(RW)+$
\xEP

\xEH
Fact \ref{Fact Smooth-Sound}
\xEH
\xEH
Proposition \ref{Proposition D-5.4.1}
\xEH
$(SC)+(PR)+$
\xEP

\xEH
page \pageref{Fact Smooth-Sound}
\xEH
\xEH
page \pageref{Proposition D-5.4.1}
\xEH
$(CUM)$
\xEP

\cline{2-2}
\cline{4-4}

\xEH
$\xch$ $(\xcv)$
\xEH
\xEH
$\xci$ $(\xcv)$
\xEH
\xEP

\xEH
Proposition 3.3.4
\xEH
\xEH
Proposition \ref{Proposition D-5.4.1}
\xEH
\xEP

\xEH
in \cite{Sch04},
\xEH
\xEH
page \pageref{Proposition D-5.4.1}
\xEH
\xEP

\xEH
Proposition \ref{Proposition D-4.5.3}
\xEH
\xEH
\xEH
\xEP

\xEH
page \pageref{Proposition D-4.5.3}
\xEH
\xEH
\xEH
\xEP

\cline{2-2}
\cline{4-4}

\xEH
$\xcH$ without $(\xcv)$
\xEH
\xEH
$\xcH$ without $(\xbm dp)$
\xEH
\xEP

\xEH
See \cite{Sch04}
\xEH
\xEH
Example \ref{Example Pref-Dp}
\xEH
\xEP

\xEH
\xEH
\xEH
page \pageref{Example Pref-Dp}
\xEH
\xEP

\hline

$(\xbm \xcc)+(\xbm PR)+(\xbm CUM)$
\xEH
$\xci$
\xEH
smooth+transitive
\xEH
$\xch$ $(\xbm dp)$
\xEH
$(LLE)+(RW)+$
\xEP

\xEH
Fact \ref{Fact Smooth-Sound}
\xEH
\xEH
Proposition \ref{Proposition D-5.4.1}
\xEH
$(SC)+(PR)+$
\xEP

\xEH
page \pageref{Fact Smooth-Sound}
\xEH
\xEH
page \pageref{Proposition D-5.4.1}
\xEH
$(CUM)$
\xEP

\cline{2-2}
\cline{4-4}

\xEH
$\xch$ $(\xcv)$
\xEH
\xEH
$\xci$ $(\xcv)$
\xEH
\xEP

\xEH
Proposition 3.3.8
\xEH
\xEH
Proposition \ref{Proposition D-5.4.1}
\xEH
\xEP

\xEH
in \cite{Sch04},
\xEH
\xEH
page \pageref{Proposition D-5.4.1}
\xEH
\xEP

\xEH
Proposition \ref{Proposition D-4.5.3}
\xEH
\xEH
\xEH
\xEP

\xEH
page \pageref{Proposition D-4.5.3}
\xEH
\xEH
\xEH
\xEP

\cline{2-2}
\cline{4-4}

\xEH
\xEH
\xEH
$\xcH$ without $(\xbm dp)$
\xEH
\xEP

\xEH
\xEH
\xEH
Example \ref{Example Pref-Dp}
\xEH
\xEP

\xEH
\xEH
\xEH
page \pageref{Example Pref-Dp}
\xEH
\xEP

\cline{4-5}

\xEH
\xEH
\xEH
$\xcj$ without $(\xbm dp)$
\xEH
using ``small''
\xEP

\xEH
\xEH
\xEH
Proposition 5.2.9, 5.2.11
\xEH
exception sets
\xEP

\xEH
\xEH
\xEH
in \cite{Sch04}
\xEH
\xEP

\hline

$(\xbm\xcc)+(\xbm=)+(\xbm PR)+$
\xEH
$\xci$
\xEH
ranked, $\xcg 1$ copy
\xEH
\xEH
\xEP

$(\xbm=')+(\xbm\xFO)+(\xbm\xcv)+$
\xEH
Fact \ref{Fact Rank-Hold}
\xEH
\xEH
\xEH
\xEP

$(\xbm\xcv ')+(\xbm\xbe)+(\xbm RatM)$
\xEH
page \pageref{Fact Rank-Hold}
\xEH
\xEH
\xEH
\xEP

\hline

$(\xbm\xcc)+(\xbm=)+(\xbm PR)+$
\xEH
$\xcH$
\xEH
ranked
\xEH
\xEH
\xEP

$(\xbm\xcv)+(\xbm\xbe)$
\xEH
Example \ref{Example Rank-Copies}
\xEH
\xEH
\xEH
\xEP

\xEH
page \pageref{Example Rank-Copies}
\xEH
\xEH
\xEH
\xEP

\hline

$(\xbm\xcc)+(\xbm=)+(\xbm \xCQ)$
\xEH
$\xcj$, $(\xcv)$
\xEH
ranked,
\xEH
\xEH
\xEP

\xEH
Proposition 3.10.11
\xEH
1 copy + $(\xbm \xCQ)$
\xEH
\xEH
\xEP

\xEH
in \cite{Sch04}
\xEH
\xEH
\xEH
\xEP

\hline

$(\xbm\xcc)+(\xbm=)+(\xbm \xCQ)$
\xEH
$\xcj$, $(\xcv)$
\xEH
ranked, smooth,
\xEH
\xEH
\xEP

\xEH
Proposition 3.10.11
\xEH
1 copy + $(\xbm \xCQ)$
\xEH
\xEH
\xEP

\xEH
in \cite{Sch04}
\xEH
\xEH
\xEH
\xEP

\hline

$(\xbm\xcc)+(\xbm=)+(\xbm \xCQ fin)+$
\xEH
$\xcj$, $(\xcv)$, singletons
\xEH
ranked, smooth,
\xEH
\xEH
\xEP

$(\xbm\xbe)$
\xEH
Proposition 3.10.12
\xEH
$\xcg$ 1 copy + $(\xbm \xCQ fin)$
\xEH
\xEH
\xEP

\xEH
in \cite{Sch04}
\xEH
\xEH
\xEH
\xEP

\hline

$(\xbm\xcc)+(\xbm PR)+(\xbm \xFO)+$
\xEH
$\xcj$, $(\xcv)$, singletons
\xEH
ranked
\xEH
$\xcH$ without $(\xbm dp)$
\xEH
$(RatM), (RatM=)$,
\xEP

$(\xbm \xcv)+(\xbm\xbe)$
\xEH
Proposition 3.10.14
\xEH
$\xcg$ 1 copy
\xEH
Example \ref{Example Rank-Dp}
\xEH
$(Log\xcv), (Log\xcv ')$
\xEP

\xEH
in \cite{Sch04}
\xEH
\xEH
page \pageref{Example Rank-Dp}
\xEH
\xEP

\cline{4-5}

\xEH
\xEH
\xEH
$\xcJ$ without $(\xbm dp)$
\xEH
any ``normal''
\xEP

\xEH
\xEH
\xEH
Proposition 5.2.16
\xEH
characterization
\xEP

\xEH
\xEH
\xEH
in \cite{Sch04}
\xEH
of any size
\xEP

\hline

\end{tabular}

}

\section{Preferential structures without domain conditions}
\label{Section Without-Domain}
\subsection{General discussion}

We treat in this Section the general case without conditions on the
domain.
We will see that it is more difficult than when we can impose the
usual conditions (closure under finite intersections and finite unions).
The latter case will be dealt with briefly (as most of it was already
done in  \cite{Sch04}) in
Section \ref{Section With-Domain} (page \pageref{Section With-Domain}).
\subsubsection{General preferential structures}
\label{Section GPS-Without-Domain}

We give now just three simple facts to put the reader in the
right mood for what follows.
\index{Fact Pref-Sound}

\bfa

$\hspace{0.01em}$

% (+++ Orig. No.:  Fact Pref-Sound +++)

\label{Fact Pref-Sound}

$( \xbm \xcc )$ and $( \xbm PR)$ hold in all preferential structures.
\index{Fact Pref-Sound Proof}

\efa

\subparagraph{
Proof
}

$\hspace{0.01em}$

% (+++ Orig.:  Proof +++)

Trivial. The central argument is: if $x,y \xbe X \xcc Y,$ and $x \xeb y$
in $X,$ then also
$x \xeb y$ in $Y.$

$ \xcz $
\\[3ex]
\index{Fact Smooth-Sound}

\bfa

$\hspace{0.01em}$

% (+++ Orig. No.:  Fact Smooth-Sound +++)

\label{Fact Smooth-Sound}

$( \xbm \xcc ),$ $( \xbm PR),$ and $( \xbm CUM)$ hold in all smooth
preferential structures.
\index{Fact Smooth-Sound Proof}

\efa

\subparagraph{
Proof
}

$\hspace{0.01em}$

% (+++ Orig.:  Proof +++)

By Fact \ref{Fact Pref-Sound} (page \pageref{Fact Pref-Sound}), we only have to
show $( \xbm CUM).$
By Fact \ref{Fact Mu-Base} (page \pageref{Fact Mu-Base}), $( \xbm CUT)$ follows
from $( \xbm PR),$
so it remains to show
$( \xbm CM).$ So suppose $ \xbm (X) \xcc Y \xcc X,$ we have to show $ \xbm
(Y) \xcc \xbm (X).$ Let
$x \xbe X- \xbm (X),$ so there is $x' \xbe X,$ $x' \xeb x,$ by smoothness,
there must be $x'' \xbe \xbm (X),$
$x'' \xeb x,$ so $x'' \xbe Y,$ and $x \xce \xbm (Y).$ The proof for the
case with copies is
analogous.
\index{Example Pref-Dp}

\be

$\hspace{0.01em}$

% (+++ Orig. No.:  Example Pref-Dp +++)

\label{Example Pref-Dp}

This example was first given in [Sch92]. It shows
that condition $ \xCf (PR)$ may fail in preferential structures which are
not
definability preserving.

Let $v( \xdl ):=\{p_{i}:i \xbe \xbo \},$ $n,n' \xbe M_{ \xdl }$ be defined
by $n \xcm \{p_{i}:i \xbe \xbo \},$
$n' \xcm \{ \xCN p_{0}\} \xcv \{p_{i}:0<i< \xbo \}.$

Let $ \xdm:= \xBc M_{ \xdl }, \xeb  \xBe $ where only $n \xeb n',$ i.e. just two
models are
comparable. Note that the structure is transitive and smooth.
Thus, by Fact \ref{Fact Smooth-Sound} (page \pageref{Fact Smooth-Sound})  $(
\xbm \xcc ),$ $( \xbm PR),$
$( \xbm CUM)$ hold.

Let $ \xbm:= \xbm_{ \xdm },$ and $ \xcn $ be defined as usual by $ \xbm
.$

Set $T:= \xCQ,$ $T':=\{p_{i}:0<i< \xbo \}.$ We have $M_{T}=M_{ \xdl },$
$f(M_{T})=M_{ \xdl }-\{n' \},$ $M_{T' }=\{n,n' \},$
$f(M_{T' })=\{n\}.$ So by the result of Example \ref{Example Not-Def} (page
\pageref{Example Not-Def}),
$f$ is not
definability preserving, and, furthermore,
$ \ol{ \ol{T} }= \ol{T},$ $ \ol{ \ol{T' } }= \ol{\{p_{i}:i< \xbo \}},$
$so$ $p_{0} \xbe \ol{ \ol{T \xcv T' } },$ $but$ $ \ol{ \ol{ \ol{T} } \xcv
T' }= \ol{ \ol{T} \xcv T' }= \ol{T' },$ $so$ $p_{0} \xce
 \ol{ \ol{ \ol{T} } \xcv T' },$
contradicting $ \xCf (PR),$ which holds in all definability preserving
preferential structures $ \xcz $
\\[3ex]

\ee

We know from Fact \ref{Fact Pref-Sound} (page \pageref{Fact Pref-Sound})
that $f$ has to satisfy $( \xbm \xcc )$ and $( \xbm PR).$
Let then $u \xbe U \xbe \xdy.$

If $u \xbe f(U),$ then, for this $u,$ this $U$ there is nothing to do, we
just have to
take care that at least one copy of $u$ will be minimal in $U.$

If $u \xce f(U),$ then $u$ must be minimized by some $u' \xbe U$ - more
precisely:
It might well be that in all smaller $U' $ is not minimized: $U' \xcB U,$
$U' \xbe \xdy,$
$u \xbe U' \xch u \xbe f(U' ).$ If e.g. $U=\{u,u',u'' \},$ and $U'
=\{u,u' \},$ $U'' =\{u,u'' \}$ with $U',U'' \xbe \xdy,$
then there cannot be $u \xeb u,$ nor $u' \xeb u,$ nor $u'' \xeb u,$ and we
have to make copies of
$u,$ so that only in $U,$ but neither in $U' $ nor in $U'' $ all copies
are minimized.
Thus, what we have to do, is to create $ \xBc u,u \xBe,$ $ \xBc u,u'  \xBe,$ $
\xBc u,u''  \xBe,$ and
to make
$u \xeb  \xBc u,u \xBe,$ $u' \xeb  \xBc u,u'  \xBe,$ $u'' \xeb  \xBc u,u'' 
\xBe $ (or something
similar). Thus, in the presence
of full $U,$ all copies will be minimized, but in all $U' \xcB U$ at least
one copy of
$u$ is not minimized. When we look now at our construction, we note the
following:
(a) $u$ is minimized in $U,$ (b) we took no commitment for other $U' \xcc
U.$ Thus, we
might not know anything about such other $U',$ and leave this question
totally
open - we preserved our ignorance, (c) the construction is independent of
all
other $U' $ - except that any $U' $ with $U \xcc U' $ will also minimize
$u,$ but this was
an inevitable fact. Note that there might also be $U' \xcc U$ with $u \xbe
U' -f(U' ),$
but no minimal one (wrt. $ \xcc ).$

We can see the problem of copies also as preservation of ignorance, and
can
also solve it with many structures - as is done in  \cite{SGMRT00}.

We have to do this construction now also for other $u$ and $U,$ so we will
perhaps
introduce copies for other elements, too, suppose we have copies $ \xBc u',y
\xBe $
and
$ \xBc u',y'  \xBe $ for the above $u'.$ As $ \xbm $ is insensitive to the
particular index, and
it wants only at least one copy of $u' $ to be smaller than $ \xBc u,u'  \xBe,$
we
have a
problem we cannot decide: Shall we make
$ \xBc u',y \xBe  \xeb  \xBc u,u'  \xBe,$ or
$ \xBc u',y'  \xBe  \xeb  \xBc u,u'  \xBe,$ or
both? Deciding for one solution would go beyond our knowledge (though it
would
do no harm, representation would be preserved), and we would not preserve
our
ignorance. The only honest solution to the problem is to admit that we do
not know, and branch into all possible cases, i.e. for any nonempty subset
$X' $
of the copies of $u',$ we make all copies
$ \xBc u',y \xBe  \xbe X' $ smaller than $ \xBc u,u'  \xBe.$ Thus, we
construct many structures instead of one, and say: the real one is one of
them,
but we do not know.

Note that all
these structures will be different, as points which are logically
indiscernible
will be different from an order theoretic point of view. We should also
note the
parallel here to Kripke models for modal logic, where the standard
construction
works with complete consistent theories in the full $ \xcX -$language,
with nested
$ \xcX ' $s etc., where we might see the differences between two points
only when
following arrows to some depths. Here, the situation is similar: $ \xBc u,u' 
\xBe
\xee  \xBc u',y \xBe $
and $ \xBc u,u'  \xBe  \xee  \xBc u',y'  \xBe $ are the same on level 0: it is
in both cases
$u.$ On level 1,
they are the same two, as we see in both cases $u',$ and only in level 2,
they
may begin to look differently: $y$ and $y' $ may choose different
successors.

Once we are aware of the problem - i.e. we do not know
enough for a decision - we can, of course, choose one sufficient for our
representation purpose. But it is important to see the arbitrariness of
the
decision we take. The natural solution will then be to decide for making
ALL
copies of $u' $ smaller than $ \xBc u,u'  \xBe.$

Perhaps $U$ is not the only set s.t. $u \xbe U-f(U),$ but there is also
$U' $ with
$u \xbe U' -f(U' ).$ In this case, we repeat the construction for $U',$
and choose
now for each copy of $u$ (at least) one element in $U,$ and one in $U',$
which
minimizes this copy. Then, in the presence of all elements of $U,$ or all
elements of $U',$ all copies will be minimized. The solution is thus to
consider all copies $ \xBc u,g \xBe,$ where $g \xbe \xbP \{X \xbe \xdy:u \xbe
X-f(X)\}.$

We will define $ \xBc u,g \xBe  \xee  \xBc y,h \xBe $ iff $y \xbe ran(g)$ - this
is the
adequate
condition - and forget about $h,$ but keep in mind that we took here an
arbitrary decision, and should, to preserve ignorance, branch into all
possibilities. We will see its importance immediately, now.
\subsubsection{Transitive preferential structures}
\label{Section 2.2.1.3}
%  Section (3.3):  Transitive preferential structures
%  Section (3.3):  Transitive preferential structures
% %
% ====================================================

The Example \ref{Example Trans-1} (page \pageref{Example Trans-1})  shows that
we cannot just make the
above construction
transitive and preserve representation. This is an illustration of the
fact
that we have to be careful about excessive relations.

The new construction avoids this, as it ``looks ahead'':

Seen from one fixed point, an arbitrary relation is a graph where we
identify
$ \xee $ with $ \xcp,$ i.e. $u \xee x$ will be written $u \xcp x.$ The
picture is perhaps easier to
read when we write this graph as a tree, repeating nodes when necessary.
So, from the starting point $u,$ we can go to $x$ and $x',$ from $x$ to
$y$ and $y',$ from
$x' $ to $w$ and $w' $ and $w'' $ etc. So we can write the tree of all
direct and indirect
successors of $u,$ $t(u),$ and if $x$ is a direct successor of $u,$ then
the tree for $x,$
$t(x),$ will be a subtree of $t(u),$ beginning at the successor $x$ of the
root in
$t(u).$

This gives as now a method to control all direct and indirect successors
of
an element. We write as index above tree, and define $ \xBc u,t(u) \xBe  \xee
\xBc x,t(x) \xBe $ iff
$t(x)$ is the subtree of $t(u)$ which begins at the direct successor $x$
of $u.$ In
the next step, we make the relation transitive, of course, we now have to
see
that this can be done without destroying representation, and we will use
in
our construction special choice functions, which always choose for $u \xbe
Y-f(u)$
$u$ itself - this is allowed, and they will do what we want. The details
are
given in the formal construction below.
\subsubsection{Smooth structures}
\label{Section 2.2.1.4}
%  Section (3.4):  Smooth structures
%  Section (3.4):  Smooth structures
% %
% ===================================

In analogy to Case (1), and with the same argument, we will consider
choice
functions $g \xbe \xbP \{f(X):x \xbe X-f(X)\}.$

(In the final construction, we will construct simultaneously for all $u,U$
s.t.
$u \xbe f(U)$ a $U-$minimal copy, so in the following intuitive
discussion, it will
suffice to find minimal $u,$ $x,$ etc. with the required properties. This
remark is
destined for readers who wonder how this will all fit together. We should
also
note that we will again be in the dilemma which copy to make smaller, and
will
do so for all candidates - violating our principle of preserving
ignorance. Yet,
as before, as long as we are aware of it, it will do no harm.)

To see the new problem arising now, we start with $U,$ and
suppose that $u \xbe f(U).$ Let now $u \xbe X-f(X),$ then we have to find
$x \xbe f(X)$ below $u.$
First, $x$ must not be in $U,$ as we would have destroyed minimality of
$u$ in $U,$
this is analogous to Case (1), so we need $f(X)-U \xEd \xCQ.$ But let now
$u \xbe f(Y),$ $x \xbe Y.$
In Case (1), it
was sufficient to find another copy of $u,$ which is minimal in $Y.$ Now,
we have
to do more: to find an $y \xbe f(Y),$ $y$ below $u,$ so smoothness will
hold. We will
call the following process the ``repairing process for $u,$ $x,$ and $Y$ ''.
Suppose then that $u \xbe f(Y),$ and $x \xbe f(Y)$ for $Y \xbe \xdy.$
Then we have destroyed
minimality of $u$ in $Y,$ but have repaired smoothness immediately again
by finding
the minimal $x.$ The situation is different if $x \xbe Y-f(Y)$ (and there
was no $x' \xeb u$
$x' \xbe f(Y)$ chosen at the same time). Then we have destroyed minimality
of
$u$ in $Y,$ without repairing smoothness, and we have to repair it now by
finding suitable $y \xeb u,$ $y \xbe f(Y).$ Of course, $y$ must not be in
$U,$ as this would
destroy minimality of $u$ in $U.$

Thus, we have to find for all $Y$ with $u \xbe f(Y),$ $x \xbe Y-f(Y)$ some
$y \xbe f(Y),$ $y \xeb u,$ $y \xce U.$
Note that this repair process is individual, i.e. we do not have to find
one
universal $y$ which repairs lost minimality for ALL such $Y$ at the same
time, but
it suffices to do it one by one, individually for every single such $Y.$

But now, the solutions $y$ for such $Y$ may have introduced new problems:
Not only
$x$ is below $u,$ but also $y$ is below $u.$ If there is now $Z \xbe \xdy
$ s.t. $u \xbe f(Z),$ and
$y \xbe Z-f(Z),$ then we have to do the same repairing process for $u,$
$y,$ $Z:$ find
suitable $z \xbe f(Z)$ below $u,$ $z \xce U,$ etc. So we will have an
infinite repairing
process, where each step may introduce new problems, which will be
repaired in
the next step.

To illustrate that the problem is still a bit more complicated, we make a
definition, and see that we have to avoid in above situation not only $U,$
but
$H(U,u),$ to be defined now.

$H(U,u)_{0}$ $:=$ $U,$

$H(U,u)_{ \xba +1}$ $:=$ $H(U,u)_{ \xba }$ $ \xcv $ $ \xcV \{X:$ $u \xbe
X$ $ \xcu $ $ \xbm (X) \xcc H(U,u)_{ \xba }\},$

$H(U,u)_{ \xbl }$ $:=$ $ \xcV \{H(U,u)_{ \xba }: \xba < \xbl \}$ for
$limit( \xbl ),$

$H(U,u)$ $:=$ $ \xcV \{H(U,u)_{ \xba }: \xba < \xbk \}$ for $ \xbk $
sufficiently big

$(card(Z)$ suffices, as the procedure trivializes, when we cannot add any
new
elements).

$ \xCf (HU,u)$ is the property:

$u \xbe \xbm (U),$ $u \xbe Y- \xbm (Y)$ $ \xch $ $ \xbm (Y) \xcC H(U,u)$ -
of course for all $u$ and $U.$

$(U,Y \xbe \xdy ).$

\bfa

$\hspace{0.01em}$

% (+++ Orig. No.:  Fact D-3.4.1 +++)

\label{Fact D-3.4.1}

$ \xCf (HU,u)$ holds in smooth structures.

\efa

The proof is given in Fact \ref{Fact D-4.4.3} (page \pageref{Fact D-4.4.3}) 
(2).

We note now that we have to consider our principle of preserving ignorance
again: We can choose first arbitrary $y \xbe f(Y)$ to repair for $u,$ $x,$
$Y.$ So which
one we choose is - a priori - an arbitrary choice. Yet, this choice might
have repercussions later, as different $y$ and $y' $ chosen to repair for
$u,$ $x,$ $Y$
might force different repairs for $u,$ $y,$ $Z$ or $u,$ $y',$ $Z',$ as
$Z$ might be such
that $u \xbe f(Z),$ $y \xbe Z-f(Z),$ and $Z' $ s.t. $u \xbe f(Z),$ $y'
\xbe Z-f(Z' ),$ and it might be
possible to find suitable $z \xbe f(Z),$ $z \xce H(U,u),$ but no suitable
$z',$ etc. So
the, at first sight, arbitrary choice might reveal an impasse later on. We
will see that we can easily solve the problem in the not necessarily
transitive
case, but we do not see at the time of writing any easy solution in the
transitive case, if the domain is not necessarily closed under finite
unions.
\subsubsection{Transitive smooth structures}
\label{Section 2.2.1.5}
%  Section (3.5):  The transitive smooth case
%  Section (3.5):  The transitive smooth case
% %
% ============================================

The basic, now more complicated, situation to consider is the following:

Let again $u \xbe f(U),$ $u \xbe X-f(X),$ we have to find $x \xbe f(X)$ -
outside $H(U,u)$ as in
Case (3). Thus, we need again $f(X)-H(U,u) \xEd \xCQ.$ Again, we have to
repair
all damage done, i.e. for all $u,$ $x,$ $Y$ as discussed in Case (3), the
infinite
repair process discussed there.

Suppose now that $x \xbe Y-f(Y),$ so we have to find $y \xbe f(Y),$
outside $H(U,u)$ by
transitivity of the relation, as $y \xeb x \xeb u,$ and, in addition
outside $H(X,x),$ as in
Case (3), now for $X$ and $x.$ Thus, we need $f(Y)-(H(U,u) \xcv H(X,x))
\xEd \xCQ.$ Moreover,
we have to do the same for all elements $y$ introduced by the above
repairing
process. Again, we have to do repairing: $y \xeb u,$ and $y \xeb x,$ so
for all $Y' $ s.t.
$u \xbe f(Y' ),$ $y \xbe Y' -f(Y' )$ we have to repair for $u,$ $y,$ $Y'
,$ and if $x \xbe f(Y' ),$
$y \xbe Y' -f(Y' )$ we have to repair for $x,$ $y,$ $Y',$ creating new
smaller elements, etc.

If $y \xbe Z-f(Z),$ we have to find $z \xbe f(Z),$ outside $H(U,u),$
$H(X,x),$ $H(Y,y),$ etc.,
so the further we go down, the longer the condition will be.
Thus, we need $f(Z)-(H(U,u) \xcv H(X,x) \xcv H(Y,y)) \xEd \xCQ.$ And,
again we have to repair,
for $ \xCf u,z,$ $ \xCf x,z,$ and $ \xCf y,z.$

And so on.

Note again the arbitrariness of choice, when there is not
a unique solution, i.e. no unique $x,$ $y,$ $z$ etc. This has to be
considered when
we want to respect preservation of ignorance, but also an early wrong
choice
might lead to an impasse, leading to backtracking to this early wrong
choice.

We will see that the closure of the domain under $( \xcv )$ makes all this
easily
possible, but the authors do not see an easy solution in the absence of
$( \xcv )$ at the time of writing - the problem is an initial potentially
wrong
choice, which we do not see how to avoid other than by trying.

So we will give here only a formal negative result by an example,
see Example \ref{Example D-4.5.1} (page \pageref{Example D-4.5.1}),
and essentially repeat the result given in
 \cite{Sch04} using $( \xcv ),$ see
Proposition \ref{Proposition D-4.5.3} (page \pageref{Proposition D-4.5.3}),
presented in
Section \ref{Section With-Domain} (page \pageref{Section With-Domain}).
\subsubsection{Ranked structures}
\label{Section 2.2.1.6}

We give here some definitions, and show elementary facts about ranked
structures. We also prove a general abstract nonsense fact about extending
relations, to be used here and also later on.

The crucial fact will be
Lemma \ref{Lemma 1-infin} (page \pageref{Lemma 1-infin}), it shows
that we can do with either one or infinitely many copies of each model.
The reason behind it is the following: Suppose we have exactly two copies
of one
model, $m,$ $m',$ where $m$ and $m' $ have the same logical properties.
If, e.g.,
$m \xeb m',$ then, as we consider only minimal elements,
$m' $ will be ``invisible''.
If $m$ and $m' $ are incomparable, then, by rankedness (modularity), they
will
have the same elements above (and below) themselves: they have the same
behavior in the preferential structure.
An immediate consequence is the ``singleton property'' of
Lemma \ref{Lemma 1-infin} (page \pageref{Lemma 1-infin}) :
One element
suffices to destroy minimality, and it suffices to look at pairs (and
singletons).

We first note the following trivial

\bfa

$\hspace{0.01em}$

% (+++ Orig. No.:  Fact D-4.6.1 +++)

\label{Fact D-4.6.1}

In a ranked structure, smoothness and the condition

\index{$(\xbm \xCQ)$}
$( \xbm \xCQ )$ $X \xEd \xCQ $ $ \xch $ $ \xbm (X) \xEd \xCQ $

are (almost) equivalent.

\efa

\subparagraph{
Proof
}

$\hspace{0.01em}$

% (+++ Orig.:  Proof +++)

Suppose $( \xbm \xCQ )$ holds, and let $x \xbe X- \xbm (X),$ $x' \xbe \xbm
(X).$ Then $x' \xeb x$ by rankedness.
Conversely, if the structure is smooth and there is an element $x \xbe X$
in the
structure (recall that structures may have ``gaps'', but this condition is a
minor point, which we shall neglect here - this is the precise meaning of
``almost'' ), then either $x \xbe \xbm (X)$ or there is $x' \xeb x,$ $x'
\xbe \xbm (X),$ so $ \xbm (X) \xEd \xCQ.$ $ \xcz $
\\[3ex]
\index{Fact Rank-Auxil}

\bfa

$\hspace{0.01em}$

% (+++ Orig. No.:  Fact Rank-Auxil +++)

\label{Fact Rank-Auxil}

In the presence of $( \xbm =)$ and $( \xbm \xcc ),$ $f(Y) \xcs (X-f(X))
\xEd \xCQ $ is equivalent to
$f(Y) \xcs X \xEd \xCQ $ and $f(Y) \xcs f(X)= \xCQ.$
\index{Fact Rank-Auxil Proof}

\efa

\subparagraph{
Proof
}

$\hspace{0.01em}$

% (+++ Orig.:  Proof +++)

$f(Y) \xcs (X-f(X))$ $=$ $(f(Y) \xcs X)-(f(Y) \xcs f(X)).$

`` $ \xci $ '': Let $f(Y) \xcs X \xEd \xCQ,$ $f(Y) \xcs f(X)= \xCQ,$ so
$f(Y) \xcs (X-f(X)) \xEd \xCQ.$

`` $ \xch $ '': Suppose $f(Y) \xcs (X-f(X)) \xEd \xCQ,$ so $f(Y) \xcs X
\xEd \xCQ.$ Suppose $f(Y) \xcs f(X) \xEd \xCQ,$ so
by $( \xbm \xcc )$ $f(Y) \xcs X \xcs Y \xEd \xCQ,$ so
by $( \xbm =)$ $f(Y) \xcs X \xcs Y=f(X \xcs Y),$ and $f(X) \xcs X \xcs Y
\xEd \xCQ,$ so by $( \xbm =)$
$f(X) \xcs X \xcs Y=f(X \xcs Y),$ so $f(X) \xcs Y=f(Y) \xcs X$ and $f(Y)
\xcs (X-f(X))= \xCQ.$

$ \xcz $
\\[3ex]
\index{Fact Rank-Trans}

\bfa

$\hspace{0.01em}$

% (+++ Orig. No.:  Fact Rank-Trans +++)

\label{Fact Rank-Trans}

If $ \xeb $ on $X$ is ranked, and free of cycles, then $ \xeb $ is
transitive.
\index{Fact Rank-Trans Proof}

\efa

\subparagraph{
Proof
}

$\hspace{0.01em}$

% (+++ Orig.:  Proof +++)

Let $x \xeb y \xeb z.$ If $x \xcT z,$ then $y \xee z,$ resulting in a
cycle of length 2. If $z \xeb x,$ then
we have a cycle of length 3. So $x \xeb z.$ $ \xcz $
\\[3ex]

The following Fact is essentially Fact 3.10.8 of  \cite{Sch04}.
\index{Fact Rank-Hold}

\bfa

$\hspace{0.01em}$

% (+++ Orig. No.:  Fact Rank-Hold +++)

\label{Fact Rank-Hold}

In all ranked structures, $( \xbm \xcc ),$ $( \xbm =),$ $( \xbm PR),$ $(
\xbm =' ),$ $( \xbm \xFO ),$ $( \xbm \xcv ),$ $( \xbm \xcv ' ),$
$( \xbm \xbe ),$ $( \xbm RatM)$ will hold, if the corresponding closure
conditions are
satisfied.
\index{Fact Rank-Hold Proof}

\efa

\subparagraph{
Proof
}

$\hspace{0.01em}$

% (+++ Orig.:  Proof +++)

$( \xbm \xcc )$ and $( \xbm PR)$ hold in all preferential structures.

$( \xbm =)$ and $( \xbm =' )$ are trivial.

$( \xbm \xcv )$ and $( \xbm \xcv ' ):$ All minimal copies of elements in
$f(Y)$ have the same rank.
If some $y \xbe f(Y)$ has all its minimal copies killed by an element $x
\xbe X,$ by
rankedness, $x$ kills the rest, too.

$( \xbm \xbe ):$ If $f(\{a\})= \xCQ,$ we are done. Take the minimal
copies of a in $\{a\},$ they are
all killed by one element in $X.$

$( \xbm \xFO ):$ Case $f(X)= \xCQ:$ If below every copy of $y \xbe Y$
there is a copy of some $x \xbe X,$
then $f(X \xcv Y)= \xCQ.$ Otherwise $f(X \xcv Y)=f(Y).$ Suppose now $f(X)
\xEd \xCQ,$ $f(Y) \xEd \xCQ,$ then
the minimal ranks decide: if they are equal, $f(X \xcv Y)=f(X) \xcv f(Y),$
etc.

$( \xbm RatM):$ Let $X \xcc Y,$ $y \xbe X \xcs f(Y) \xEd \xCQ,$ $x \xbe
f(X).$ By rankedness, $y \xeb x,$ or
$y \xcT x,$ $y \xeb x$ is impossible, as $y \xbe X,$ so $y \xcT x,$ and $x
\xbe f(Y).$

$ \xcz $
\\[3ex]
\index{Definition 1-infin}

\bd

$\hspace{0.01em}$

% (+++ Orig. No.:  Definition 1-infin +++)

\label{Definition 1-infin}

Let $ \xdz = \xBc  \xdx, \xeb  \xBe $ be a preferential structure. Call $ \xdz $
$1- \xca $ over $Z,$
iff for all $x \xbe Z$ there are exactly one or infinitely many copies of
$x,$ i.e.
for all $x \xbe Z$ $\{u \xbe \xdx:$ $u= \xBc x,i \xBe $ for some $i\}$ has
cardinality 1 or $ \xcg \xbo.$

\ed

The following Lemma is Lemma 3.10.4 of  \cite{Sch04}.
\index{Lemma 1-infin}

\bl

$\hspace{0.01em}$

% (+++ Orig. No.:  Lemma 1-infin +++)

\label{Lemma 1-infin}

Let $ \xdz = \xBc  \xdx, \xeb  \xBe $ be a preferential structure and
$f: \xdy \xcp \xdp (Z)$ with $ \xdy \xcc \xdp (Z)$ be represented by $
\xdz,$ i.e. for $X \xbe \xdy $ $f(X)= \xbm_{ \xdz }(X),$
and $ \xdz $ be ranked and free of cycles. Then there is a structure $
\xdz ' $, $1- \xca $ over
$Z,$ ranked and free of cycles, which also represents $f.$
\index{Lemma 1-infin Proof}

\el

\subparagraph{
Proof
}

$\hspace{0.01em}$

% (+++ Orig.:  Proof +++)

We construct $ \xdz ' = \xBc  \xdx ', \xeb '  \xBe.$

Let $A:=\{x \xbe Z$: there is some $ \xBc x,i \xBe  \xbe \xdx,$ but for all $
\xBc x,i \xBe
\xbe \xdx $ there is
$ \xBc x,j \xBe  \xbe \xdx $ with $ \xBc x,j \xBe  \xeb  \xBc x,i \xBe \},$

let $B:=\{x \xbe Z$: there is some $ \xBc x,i \xBe  \xbe \xdx,$ s.t. for no $
\xBc x,j \xBe
\xbe \xdx $ $ \xBc x,j \xBe  \xeb  \xBc x,i \xBe \},$

let $C:=\{x \xbe Z$: there is no $ \xBc x,i \xBe  \xbe \xdx \}.$

Let $c_{i}:i< \xbk $ be an enumeration of $C.$ We introduce for each such
$c_{i}$ $ \xbo $ many
copies $ \xBc c_{i},n \xBe :n< \xbo $ into $ \xdx ',$ put all $ \xBc c_{i},n
\xBe $ above all
elements in $ \xdx,$ and order
the $ \xBc c_{i},n \xBe $ by $ \xBc c_{i},n \xBe  \xeb '  \xBc c_{i' },n'  \xBe
$ $: \xcj $ $(i=i' $ and
$n>n' )$ or $i>i'.$ Thus, all $ \xBc c_{i},n \xBe $ are
comparable.

If $a \xbe A,$ then there are infinitely many copies of a in $ \xdx,$ as
$ \xdx $ was
cycle-free, we put them all into $ \xdx '.$
If $b \xbe B,$ we choose exactly one such minimal element $ \xBc b,m \xBe $
(i.e.
there
is no $ \xBc b,n \xBe  \xeb  \xBc b,m \xBe )$ into $ \xdx ',$ and omit all other
elements. (For definiteness, assume in all applications $m=0.)$
For all elements from A and $B,$ we take the restriction of the order $
\xeb $ of $ \xdx.$
This is the new structure $ \xdz '.$

Obviously, adding the $ \xBc c_{i},n \xBe $ does not introduce cycles,
irreflexivity
and
rankedness are preserved. Moreover, any substructure of a cycle-free,
irreflexive,
ranked structure also has these properties, so $ \xdz ' $ is $1- \xca $
over $Z,$ ranked and
free of cycles.

We show that $ \xdz $ and $ \xdz ' $ are equivalent. Let then $X \xcc Z,$
we have to prove
$ \xbm (X)= \xbm ' (X)$ $( \xbm:= \xbm_{ \xdz }$, $ \xbm ':= \xbm_{
\xdz ' }).$

Let $z \xbe X- \xbm (X).$ If $z \xbe C$ or $z \xbe A,$ then $z \xce \xbm '
(X).$ If $z \xbe B,$
let $ \xBc z,m \xBe $ be the chosen element. As $z \xce \xbm (X),$ there is $x
\xbe
X$ s.t. some $ \xBc x,j \xBe  \xeb  \xBc z,m \xBe.$
$x$ cannot be in $C.$ If $x \xbe A,$ then also $ \xBc x,j \xBe  \xeb '  \xBc z,m
\xBe $. If
$x \xbe B,$ then there is some
$ \xBc x,k \xBe $ also in $ \xdx '.$ $ \xBc x,j \xBe  \xeb  \xBc x,k \xBe $ is
impossible. If $ \xBc x,k \xBe
\xeb  \xBc x,j \xBe,$ then $ \xBc z,m \xBe  \xee  \xBc x,k \xBe $
by transitivity. If $ \xBc x,k \xBe  \xcT  \xBc x,j \xBe $, then also $ \xBc z,m
\xBe  \xee  \xBc x,k \xBe $ by
rankedness. In any
case, $ \xBc z,m \xBe  \xee '  \xBc x,k \xBe,$ and thus $z \xce \xbm ' (X).$

Let $z \xbe X- \xbm ' (X).$ If $z \xbe C$ or $z \xbe A,$ then $z \xce \xbm
(X).$ Let $z \xbe B,$ and some $ \xBc x,j \xBe  \xeb '  \xBc z,m \xBe.$
$x$ cannot be in $C,$ as they were sorted on top, so $ \xBc x,j \xBe $ exists in
$
\xdx $ too and
$ \xBc x,j \xBe  \xeb  \xBc z,m \xBe.$ But if any other $ \xBc z,i \xBe $ is
also minimal in $ \xdz $
among the $ \xBc z,k \xBe,$
then by rankedness also $ \xBc x,j \xBe  \xeb  \xBc z,i \xBe,$ as $ \xBc z,i
\xBe  \xcT  \xBc z,m \xBe,$ so $z
\xce \xbm (X).$ $ \xcz $
\\[3ex]

\bn

$\hspace{0.01em}$

% (+++ Orig. No.:  Notation D-4.6.2 +++)

\label{Notation D-4.6.2}

We fix the following notation: $A:=\{x \xbe Z:f(x)= \xCQ \}$ and $B:=Z-A$
(here
and in future we sometimes write $f(x)$ for $f(\{x\}),$ likewise $f(x,x'
)=x$ for
$f(\{x,x' \})=\{x\},$ etc., when the meaning is obvious).

\en

\bco

$\hspace{0.01em}$

% (+++ Orig. No.:  Corollary D-4.6.5 +++)

\label{Corollary D-4.6.5}

If $f$ can be represented by a ranked $ \xdz $ free of cycles, then there
is $ \xdz ',$ which
is also ranked and cycle-free, all $b \xbe B$ occur in 1 copy, all $a \xbe
A$ $ \xca $ often.
$ \xcz $
\\[3ex]

\eco

The following Example was presented in Fact 3.10.13 of  \cite{Sch04}.
\index{Example Rank-Copies}

\be

$\hspace{0.01em}$

% (+++ Orig. No.:  Example Rank-Copies +++)

\label{Example Rank-Copies}

This example shows that the conditions $( \xbm \xcc )+( \xbm PR)+( \xbm
=)+( \xbm \xcv )+( \xbm \xbe )$
can be satisfied, and still representation by a ranked structure
is impossible.

Consider $ \xbm (\{a,b\})= \xCQ,$ $ \xbm (\{a\})=\{a\},$ $ \xbm
(\{b\})=\{b\}.$ The conditions
$( \xbm \xcc )+( \xbm PR)+( \xbm =)+( \xbm \xcv )+( \xbm \xbe )$
hold trivially. This is representable, e.g. by $a_{1} \xed b_{1} \xed
a_{2} \xed b_{2} \Xl $ without
transitivity. (Note that rankedness implies transitivity,
$a \xec b \xec c,$ but not for $a=c.)$ But this cannot be represented by a
ranked
structure: As $ \xbm (\{a\}) \xEd \xCQ,$ there must be a copy $a_{i}$ of
minimal rank, likewise for
$b$ and some $b_{i}.$ If they have the same rank, $ \xbm
(\{a,b\})=\{a,b\},$ otherwise it will be
$\{a\}$ or $\{b\}.$

$ \xcz $
\\[3ex]

\ee

In the general situation we have possibly $U \xEd \xCQ,$ but $f(U)= \xCQ
.$

In this case, we only know that below each $u \xbe U,$ there must be
infinitely many
$u' \xbe U$ or infinitely many copies of such $u' \xbe U.$ (It is only in
such cases that
we need copies for representation in ranked structures,
see Lemma \ref{Lemma 1-infin} (page \pageref{Lemma 1-infin}).)
Thus, the amount of information we have is very small. It is not
surprising
that representation problems are now difficult, as we will see below (see
Section \ref{Section 2.2.6} (page \pageref{Section 2.2.6}) ),
and we will not go into more details here.
\subsection{Detailed discussion}
\label{Section Pref-Details}
\subsubsection{General preferential structures}
\label{Section 2.2.2}

The material in this Section is taken from  \cite{Sch04},
Section 3.2.1 there,
the result was already shown in  \cite{Sch92} with the
same methods.
\index{Proposition Pref-Complete}

\bp

$\hspace{0.01em}$

% (+++ Orig. No.:  Proposition Pref-Complete +++)

\label{Proposition Pref-Complete}

Let $ \xbm: \xdy \xcp \xdp (U)$ satisfy $( \xbm \xcc )$ and $( \xbm PR).$
Then there is a preferential
structure $ \xdx $ s.t. $ \xbm = \xbm_{ \xdx }.$ See e.g.  \cite{Sch04}.
\index{Proposition Pref-Complete Proof}

\ep

\subparagraph{
Proof
}

$\hspace{0.01em}$

% (+++ Orig.:  Proof +++)

The preferential structure is defined in
Construction \ref{Construction Pref-Base} (page \pageref{Construction
Pref-Base}), Claim \ref{Claim Pref-Rep-Base} (page \pageref{Claim
Pref-Rep-Base})  shows representation.
The construction is basic for much of
the rest of the material on non-ranked structures.

% {\LARGE karl-search= Start Definition Y-Pi-x }

\index{Definition Y-Pi-x}

\bd

$\hspace{0.01em}$

% (+++ Orig. No.:  Definition Y-Pi-x +++)

\label{Definition Y-Pi-x}

For $x \xbe Z,$ let $ \xdy_{x}:=\{Y \xbe \xdy $: $x \xbe Y- \xbm (Y)\},$
$ \xbP_{x}:= \xbP \xdy_{x}.$

\ed

Note that $ \xCQ \xce \xdy_{x}$, $ \xbP_{x} \xEd \xCQ,$ and that $
\xbP_{x}=\{ \xCQ \}$ iff $ \xdy_{x}= \xCQ.$

% karl-search= End Definition Y-Pi-x
\vspace{7mm}

% *************************************

\vspace{7mm}

% {\LARGE karl-search= Start Claim Mu-f }

\index{Claim Mu-f}

\bc

$\hspace{0.01em}$

% (+++ Orig. No.:  Claim Mu-f +++)

\label{Claim Mu-f}

Let $ \xbm: \xdy \xcp \xdp (Z)$ satisfy $( \xbm \xcc )$ and $( \xbm PR),$
and let $U \xbe \xdy.$
Then $x \xbe \xbm (U)$ $ \xcj $ $x \xbe U$ $ \xcu $ $ \xcE f \xbe
\xbP_{x}.ran(f) \xcs U= \xCQ.$

% karl-search= End Claim Mu-f
\vspace{7mm}

% *************************************

\vspace{7mm}

% {\LARGE karl-search= Start Claim Mu-f Proof }

\index{Claim Mu-f Proof}

\ec

\subparagraph{
Proof
}

$\hspace{0.01em}$

% (+++ Orig.:  Proof +++)

Case 1: $ \xdy_{x}= \xCQ,$ thus $ \xbP_{x}=\{ \xCQ \}.$
`` $ \xch $ '': Take $f:= \xCQ.$
`` $ \xci $ '': $x \xbe U \xbe \xdy,$ $ \xdy_{x}= \xCQ $ $ \xch $ $x \xbe
\xbm (U)$ by definition of $ \xdy_{x}.$

Case 2: $ \xdy_{x} \xEd \xCQ.$
`` $ \xch $ '': Let $x \xbe \xbm (U) \xcc U.$ It suffices to show $Y \xbe
\xdy_{x}$ $ \xch $ $Y-U \xEd \xCQ.$ But if $Y \xcc U$ and
$Y \xbe \xdy_{x}$, then $x \xbe Y- \xbm (Y),$ contradicting $( \xbm PR).$
`` $ \xci $ '': If $x \xbe U- \xbm (U),$ then $U \xbe \xdy_{x}$, so $
\xcA f \xbe \xbP_{x}.ran(f) \xcs U \xEd \xCQ.$
$ \xcz $
\\[3ex]

% karl-search= End Claim Mu-f Proof
\vspace{7mm}

% *************************************

\vspace{7mm}

% {\LARGE karl-search= Start Construction Pref-Base }

\index{Construction Pref-Base}

\bcs

$\hspace{0.01em}$

% (+++ Orig. No.:  Construction Pref-Base +++)

\label{Construction Pref-Base}

Let $ \xdx:=\{ \xBc x,f \xBe :x \xbe Z$ $ \xcu $ $f \xbe \xbP_{x}\},$ and
$ \xBc x',f'  \xBe  \xeb  \xBc x,f \xBe $ $: \xcj $ $x' \xbe ran(f).$
Let $ \xdz:= \xBc  \xdx, \xeb  \xBe.$

% karl-search= End Construction Pref-Base
\vspace{7mm}

% *************************************

\vspace{7mm}

% {\LARGE karl-search= Start Claim Pref-Rep-Base }

\index{Claim Pref-Rep-Base}

\ecs

\bc

$\hspace{0.01em}$

% (+++ Orig. No.:  Claim Pref-Rep-Base +++)

\label{Claim Pref-Rep-Base}

For $U \xbe \xdy,$ $ \xbm (U)= \xbm_{ \xdz }(U).$

% karl-search= End Claim Pref-Rep-Base
\vspace{7mm}

% *************************************

\vspace{7mm}

% {\LARGE karl-search= Start Claim Pref-Rep-Base Proof }

\index{Claim Pref-Rep-Base Proof}

\ec

\subparagraph{
Proof
}

$\hspace{0.01em}$

% (+++ Orig.:  Proof +++)

By Claim \ref{Claim Mu-f} (page \pageref{Claim Mu-f}), it suffices to show that
for all $U \xbe
\xdy $
$x \xbe \xbm_{ \xdz }(U)$ $ \xcj $ $x \xbe U$ and $ \xcE f \xbe
\xbP_{x}.ran(f) \xcs U= \xCQ.$ So let $U \xbe \xdy.$
`` $ \xch $ '': If $x \xbe \xbm_{ \xdz }(U),$ then there is $ \xBc x,f \xBe $
minimal in $ \xdx \xex U$ (recall from
Definition \ref{Definition Alg-Base} (page \pageref{Definition Alg-Base})  that
$ \xdx \xex U:=\{ \xBc x,i \xBe  \xbe \xdx
:x \xbe U\}),$
so $x \xbe U,$ and there
is no $ \xBc x',f'  \xBe  \xeb  \xBc x,f \xBe,$ $x' \xbe U,$ so by $ \xbP_{x' }
\xEd \xCQ $
there is no $x' \xbe ran(f),$ $x' \xbe U,$ but then
$ran(f) \xcs U= \xCQ.$
`` $ \xci $ '': If $x \xbe U,$ and there is $f \xbe \xbP_{x}$, $ran(f)
\xcs U= \xCQ,$ then $ \xBc x,f \xBe $ is minimal in
$ \xdx \xex U.$
$ \xcz $ (Claim \ref{Claim Pref-Rep-Base} (page \pageref{Claim Pref-Rep-Base}) 
and Proposition \ref{Proposition Pref-Complete} (page \pageref{Proposition
Pref-Complete}) )
\\[3ex]

% karl-search= End Claim Pref-Rep-Base Proof
\vspace{7mm}

% *************************************

\vspace{7mm}

\subsubsection{Transitive preferential structures}
\label{Section 2.2.3}
%  Section (4.3):  Transitive preferential structures
%  Section (4.3):  Transitive preferential structures
% %
% ====================================================

The material in this Section is taken from  \cite{Sch04},
Section 3.2.2 there,
the result was already shown in  \cite{Sch92} with different
methods.

We show here:
\index{Proposition Pref-Complete-Trans}

\bp

$\hspace{0.01em}$

% (+++ Orig. No.:  Proposition Pref-Complete-Trans +++)

\label{Proposition Pref-Complete-Trans}

Let $ \xbm: \xdy \xcp \xdp (U)$ satisfy $( \xbm \xcc )$ and $( \xbm PR).$
Then there is a transitive
preferential structure $ \xdx $ s.t. $ \xbm = \xbm_{ \xdx }.$ See e.g.
 \cite{Sch04}.
\index{Proposition Pref-Complete-Trans Proof}

\ep

\subparagraph{
Proof
}

$\hspace{0.01em}$

% (+++ Orig.:  Proof +++)

% {\LARGE karl-search= Start Discussion Pref-Trans }

\index{Discussion Pref-Trans}

\paragraph{
Discussion:
}

$\hspace{0.01em}$

% (+++ Orig.:  Discussion: +++)

\label{Section Discussion:}

The Construction \ref{Construction Pref-Base} (page \pageref{Construction
Pref-Base})
(also used in  \cite{Sch92}) cannot be made transitive as it is,
this
will be shown below in Example \ref{Example Trans-1} (page \pageref{Example
Trans-1}). The second
construction
in  \cite{Sch92} is a
special one, which is transitive, but uses heavily lack of smoothness.
(For
completeness' sake, we give a similar proof
in Proposition \ref{Proposition Equiv-Trans} (page \pageref{Proposition
Equiv-Trans}).) We present here a more
flexibel and more adequate construction, which avoids a certain excess in
the
relation $ \xeb $ of the construction in Proposition \ref{Proposition
Equiv-Trans} (page \pageref{Proposition Equiv-Trans}) :
There, too many elements $ \xBc y,g \xBe $
are smaller than some $ \xBc x,f \xBe,$ as the relation is independent from
$g.$
This excess
prevents transitivity.

We refine now the construction of the relation, to have better control
over
successors.

Recall that a tree of height $ \xck \xbo $ seems the right way to encode
the successors of
an element, as far as transitivity is concerned (which speaks only about
finite
chains). Now, in the basic construction, different copies have different
successors, chosen by different functions (elements of the cartesian
product).
As it suffices to make one copy of the successor smaller than the element
to be
minimized, we do the following: Let $ \xBc x,g \xBe,$ with $g \xbe \xbP \{X:x
\xbe
X-f(X)\}$ be one of the
elements of the standard construction. Let $ \xBc x',g'  \xBe $ be s.t. $x' \xbe
ran(g),$ then we
make again copies $ \xBc x,g,g'  \xBe,$ etc. for each such $x' $ and $g',$ and
make only $ \xBc x',g'  \xBe,$
but not some other $ \xBc x',g''  \xBe $ smaller than $ \xBc x,g,g'  \xBe,$ for
some other
$g'' \xbe \xbP \{X':x' \xbe X' -f(X' )\}.$ Thus, we have a much more
restricted relation, and much
better control over it. More precisely, we make trees, where we mark all
direct
and indirect successors, and each time the choice is made by the
appropriate
choice functions of the cartesian product. An element with its tree is a
successor of another element with its tree, iff the former is an initial
segment of the latter - see the definition in Construction \ref{Construction
Pref-Trees} (page \pageref{Construction Pref-Trees}).

Recall also that transitivity is for free as we can use the element itself
to
minimize it. This is made precise by the use of the trees $tf_{x}$ for a
given element $x$ and choice function $f_{x}.$ But they also serve another
purpose.
The trees $tf_{x}$ are constructed as follows: The root is $x,$ the first
branching is
done according to $f_{x},$ and then we continue with constant choice. Let,
e.g.
$x' \xbe ran(f_{x}),$ we can now always choose $x',$ as it will be a
legal successor of $x' $
itself, being present in all $X' $ s.t. $x' \xbe X' -f(X' ).$ So we have a
tree which
branches once, directly above the root, and is then constant without
branching.
Obviously, this is essentially equivalent to the old construction in the
not
necessarily transitive case. This shows two things: first, the
construction
with trees gives the same $ \xbm $ as the old construction with simple
choice
functions. Second, even if we consider successors of successors, nothing
changes: we are still with the old $x'.$ Consequently, considering the
transitive
closure will not change matters, an element $ \xBc x,tf_{x} \xBe $ will be
minimized
by its
direct successors iff it will be minimized by direct and indirect
successors.
If you like, the trees $tf_{x}$ are the mathematical construction
expressing the
intuition that we know so little about minimization that we have to
consider
suicide a serious possibility - the intuitive reason why transitivity
imposes
no new conditions.

To summarize: Trees seem the right way to encode all the information
needed for
full control over successors for the transitive case. The special trees
$tf_{x}$
show that we have not changed things substantially, i.e. the new $ \xbm
-$functions in
the simple case and for the transitive closure stay the same. We hope that
this
construction will show its usefulness in other contexts, its naturalness
and
generality seem to be a good promise.

We give below the example which shows that the old construction is too
brutal
for transitivity to hold.

Recall that transitivity permits substitution in the following sense:
If (the two copies of) $x$ is killed by $y_{1}$ and $y_{2}$ together, and
$y_{1}$ is killed by
$z_{1}$ and $z_{2}$ together, then $x$ should be killed by $z_{1},$
$z_{2},$ and $y_{2}$ together.

But the old construction substitutes too much: In the old construction,
we considered elements $ \xBc x,f \xBe,$ where $f \xbe \xbP_{x}$, with $ \xBc
y,g \xBe  \xeb
\xBc x,f \xBe $ iff $y \xbe ran(f),$
independent of $g.$ This construction can, in general, not be made
transitive,
as Example \ref{Example Trans-1} (page \pageref{Example Trans-1})  below shows.

The new construction avoids this, as it ``looks ahead'', and not all
elements
$ \xBc y_{1},t_{y_{1}} \xBe $ are smaller than $ \xBc x,t_{x} \xBe,$ where
$y_{1}$ is a child
of $x$ in $t_{x}$ (or $y_{1} \xbe ran(f)).$
The new construction is basically the same as Construction \ref{Construction
Pref-Base} (page \pageref{Construction Pref-Base}),
but avoids to make too many copies smaller than the copy to be killed.

Recall that
we need no new properties of $ \xbm $ to achieve transitivity here, as a
killed
element $x$ might (partially) ``commit suicide'', i.e. for some $i,$ $i' $
$ \xBc x,i \xBe  \xeb  \xBc x,i'  \xBe,$
so we cannot substitute $x$ by any set which does not contain $x:$ In this
simple
situation, if $x \xbe X- \xbm (X),$ we cannot find out whether all copies
of $x$ are killed
by some $y \xEd x,$ $y \xbe X.$ We can assume without loss of generality
that there is an
infinite descending chain of $x-$copies, which are not
killed by other elements. Thus, we cannot replace any $y_{i}$ as above by
any set
which does not contain $y_{i}$, but then substitution becomes trivial, as
any set
substituting $y_{i}$ has to contain $y_{i}$. Thus, we need no new
properties to
achieve transitivity.

% karl-search= End Discussion Pref-Trans
\vspace{7mm}

% *************************************

\vspace{7mm}

% {\LARGE karl-search= Start Example Trans-1 }

\index{Example Trans-1}

\be

$\hspace{0.01em}$

% (+++ Orig. No.:  Example Trans-1 +++)

\label{Example Trans-1}

As we consider only one set in each case, we can index with elements,
instead
of with functions.
So suppose $x,y_{1},y_{2} \xbe X,$ $y_{1},z_{1},z_{2} \xbe Y,$ $x \xce
\xbm (X),$ $y_{1} \xce \xbm (Y),$ and that we
need $y_{1}$ and $y_{2}$ to minimize $x,$ so there are two copies
$ \xBc x,y_{1} \xBe,$ $ \xBc x,y_{2} \xBe,$ likewise
we need $z_{1}$ and $z_{2}$ to minimize $y_{1},$ thus we have $ \xBc x,y_{1}
\xBe
\xee  \xBc y_{1},z_{1} \xBe,$ $ \xBc x,y_{1} \xBe  \xee  \xBc y_{1},z_{2} \xBe
,$
$ \xBc x,y_{2} \xBe  \xee y_{2},$ $ \xBc y_{1},z_{1} \xBe  \xee z_{1},$ $ \xBc
y_{1},z_{2} \xBe  \xee
z_{2}$ (the $z_{i}$ and $y_{2}$ are not killed).
If we take the transitive closure, we have $ \xBc x,y_{1} \xBe  \xee z_{k}$ for
any
$i,k,$ so for any $z_{k}$
$\{z_{k},y_{2}\}$ will minimize all of $x,$ which is not intended. $ \xcz
$
\\[3ex]

% karl-search= End Example Trans-1
\vspace{7mm}

% *************************************

\vspace{7mm}

\ee

The preferential structure is defined in
Construction \ref{Construction Pref-Trees} (page \pageref{Construction
Pref-Trees}), Claim \ref{Claim Tree-Repres-1} (page \pageref{Claim
Tree-Repres-1})  shows representation
for the simple structure, Claim \ref{Claim Tree-Repres-2} (page \pageref{Claim
Tree-Repres-2})
representation
for the transitive closure of the structure.

The main idea is to use the trees $tf_{x}$, whose elements are exactly
the elements
of the range of the choice function $f.$ This makes
Construction \ref{Construction Pref-Base} (page \pageref{Construction
Pref-Base})  and
Construction \ref{Construction Pref-Trees} (page \pageref{Construction
Pref-Trees})
basically equivalent,
and shows that the transitive case is
characterized by the same conditions as the general case. These trees are
defined below in Fact \ref{Fact Pref-Trees} (page \pageref{Fact Pref-Trees}),
(3),
and used in the proofs of Claim \ref{Claim Tree-Repres-1} (page \pageref{Claim
Tree-Repres-1})  and
Claim \ref{Claim Tree-Repres-2} (page \pageref{Claim Tree-Repres-2}).

Again, Construction \ref{Construction Pref-Trees} (page \pageref{Construction
Pref-Trees})  contains the basic idea for
the
treatment of the transitive case. It can certainly be re-used in other
contexts.

% {\LARGE karl-search= Start Construction Pref-Trees }

\index{Construction Pref-Trees}

\bcs

$\hspace{0.01em}$

% (+++ Orig. No.:  Construction Pref-Trees +++)

\label{Construction Pref-Trees}

(1) For $x \xbe Z,$ let $T_{x}$ be the set of trees $t_{x}$ s.t.

(a) all nodes are elements of $Z,$

(b) the root of $t_{x}$ is $x,$

(c) $height(t_{x}) \xck \xbo,$

(d) if $y$ is an element in $t_{x}$, then there is $f \xbe \xbP_{y}:=
\xbP \{Y \xbe \xdy $: $y \xbe Y- \xbm (Y)\}$
s.t. the set of children of $y$ is $ran(f).$

(2) For $x,y \xbe Z,$ $t_{x} \xbe T_{x}$, $t_{y} \xbe T_{y}$, set $t_{x}
\xem t_{y}$ iff $y$ is a (direct) child
of the root $x$ in $t_{x}$, and $t_{y}$ is the subtree of $t_{x}$
beginning at $y.$

(3) Let $ \xdz $ $:=$ $ \xBc $ $\{ \xBc x,t_{x} \xBe :$ $x \xbe Z,$ $t_{x} \xbe
T_{x}\}$,
$ \xBc x,t_{x} \xBe  \xee  \xBc y,t_{y} \xBe $ iff $t_{x} \xem t_{y}$ $ \xBe.$

% karl-search= End Construction Pref-Trees
\vspace{7mm}

% *************************************

\vspace{7mm}

% {\LARGE karl-search= Start Fact Pref-Trees }

\index{Fact Pref-Trees}

\ecs

\bfa

$\hspace{0.01em}$

% (+++ Orig. No.:  Fact Pref-Trees +++)

\label{Fact Pref-Trees}

(1) The construction ends at some $y$ iff $ \xdy_{y}= \xCQ,$ consequently
$T_{x}=\{x\}$ iff $ \xdy_{x}= \xCQ.$ (We identify the tree of height 1
with its root.)

(2) If $ \xdy_{x} \xEd \xCQ,$ $tc_{x}$, the totally ordered tree of
height $ \xbo,$ branching with $card=1,$
and with all elements equal to $x$ is an element of $T_{x}.$ Thus, with
(1), $T_{x} \xEd \xCQ $
for any $x.$

(3) If $f \xbe \xbP_{x}$, $f \xEd \xCQ,$ then the tree $tf_{x}$ with
root $x$ and otherwise
composed of the subtrees $t_{y}$ for $y \xbe ran(f),$ where $t_{y}:=y$ iff
$ \xdy_{y}= \xCQ,$
and $t_{y}:=tc_{y}$ iff $ \xdy_{y} \xEd \xCQ,$ is an element of $T_{x}$.
(Level 0 of $tf_{x}$ has $x$ as
element, the $t_{y}' s$ begin at level 1.)

(4) If $y$ is an element in $t_{x}$ and $t_{y}$ the subtree of $t_{x}$
starting at
$y,$ then $t_{y} \xbe T_{y}$.

(5) $ \xBc x,t_{x} \xBe  \xee  \xBc y,t_{y} \xBe $ implies $y \xbe ran(f)$ for
some $f \xbe
\xbP_{x}.$
$ \xcz $
\\[3ex]

% karl-search= End Fact Pref-Trees
\vspace{7mm}

% *************************************

\vspace{7mm}

% {\LARGE karl-search= Start Claim Tree-Repres-1 }

\index{Claim Tree-Repres-1}

\efa

Claim \ref{Claim Tree-Repres-1} (page \pageref{Claim Tree-Repres-1})  shows
basic representation.

\bc

$\hspace{0.01em}$

% (+++ Orig. No.:  Claim Tree-Repres-1 +++)

\label{Claim Tree-Repres-1}

$ \xcA U \xbe \xdy. \xbm (U)= \xbm_{ \xdz }(U)$

% karl-search= End Claim Tree-Repres-1
\vspace{7mm}

% *************************************

\vspace{7mm}

% {\LARGE karl-search= Start Claim Tree-Repres-1 Proof }

\index{Claim Tree-Repres-1 Proof}

\ec

\subparagraph{
Proof
}

$\hspace{0.01em}$

% (+++ Orig.:  Proof +++)

By Claim \ref{Claim Mu-f} (page \pageref{Claim Mu-f}), it suffices to show that
for all $U \xbe
\xdy $
$x \xbe \xbm_{ \xdz }(U)$ $ \xcj $ $x \xbe U$ $ \xcu $ $ \xcE f \xbe
\xbP_{x}.ran(f) \xcs U= \xCQ.$
Fix $U \xbe \xdy.$
`` $ \xch $ '': $x \xbe \xbm_{ \xdz }(U)$ $ \xch $ ex. $ \xBc x,t_{x} \xBe $
minimal
in $ \xdz \xex U,$ thus $x \xbe U$ and there is no $ \xBc y,t_{y} \xBe  \xbe
\xdz,$
$ \xBc y,t_{y} \xBe  \xeb  \xBc x,t_{x} \xBe,$ $y \xbe U.$ Let $f$ define the
set of children
of the root
$x$ in $t_{x}$. If $ran(f) \xcs U \xEd \xCQ,$ if $y \xbe U$ is a child
of $x$ in $t_{x}$, and if $t_{y}$ is the subtree
of $t_{x}$ starting at $y,$ then $t_{y} \xbe T_{y}$ and $ \xBc y,t_{y} \xBe 
\xeb
\xBc x,t_{x} \xBe,$ contradicting minimality of
$ \xBc x,t_{x} \xBe $ in $ \xdz \xex U.$ So $ran(f) \xcs U= \xCQ.$
`` $ \xci $ '': Let $x \xbe U.$ If $ \xdy_{x}= \xCQ,$ then the tree $x$
has no $ \xem -$successors, and $ \xBc x,x \xBe $ is
$ \xee -$minimal in $ \xdz.$ If $ \xdy_{x} \xEd \xCQ $ and $f \xbe
\xbP_{x}$ s.t. $ran(f) \xcs U= \xCQ,$ then $ \xBc x,tf_{x} \xBe $ is $ \xee
-$minimal
in $ \xdz \xex U.$
$ \xcz $
\\[3ex]

% karl-search= End Claim Tree-Repres-1 Proof
\vspace{7mm}

% *************************************

\vspace{7mm}

% {\LARGE karl-search= Start Claim Tree-Repres-2 }

\index{Claim Tree-Repres-2}

We consider now the transitive closure of $ \xdz.$ (Recall that $
\xeb^{*}$ denotes the
transitive closure of $ \xeb.)$ Claim \ref{Claim Tree-Repres-2} (page
\pageref{Claim Tree-Repres-2})
shows that transitivity does not
destroy what we have achieved. The trees $tf_{x}$ will play a crucial role
in the
demonstration.

\bc

$\hspace{0.01em}$

% (+++ Orig. No.:  Claim Tree-Repres-2 +++)

\label{Claim Tree-Repres-2}

Let $ \xdz ' $ $:=$ $ \xBc $ $\{ \xBc x,t_{x} \xBe :$ $x \xbe Z,$ $t_{x} \xbe
T_{x}\}$,
$ \xBc x,t_{x} \xBe  \xee  \xBc y,t_{y} \xBe $ iff $t_{x} \xem^{*}t_{y}$ $
\xBe.$

Then $ \xbm_{ \xdz }= \xbm_{ \xdz ' }.$

% karl-search= End Claim Tree-Repres-2
\vspace{7mm}

% *************************************

\vspace{7mm}

% {\LARGE karl-search= Start Claim Tree-Repres-2 Proof }

\index{Claim Tree-Repres-2 Proof}

\ec

\subparagraph{
Proof
}

$\hspace{0.01em}$

% (+++ Orig.:  Proof +++)

Suppose there is $U \xbe \xdy,$ $x \xbe U,$ $x \xbe \xbm_{ \xdz }(U),$ $x
\xce \xbm_{ \xdz ' }(U).$
Then there must be an element $ \xBc x,t_{x} \xBe  \xbe \xdz $ with no $ \xBc
x,t_{x} \xBe
\xee  \xBc y,t_{y} \xBe $ for any $y \xbe U.$
Let $f \xbe \xbP_{x}$ determine the set of children of $x$ in $t_{x}$,
then $ran(f) \xcs U= \xCQ,$
consider $tf_{x}.$ As all elements $ \xEd x$ of $tf_{x}$ are already in
$ran(f),$ no element of $tf_{x}$
is in $U.$ Thus there is no $ \xBc z,t_{z} \xBe  \xeb^{*} \xBc x,tf_{x} \xBe $
in $ \xdz $
with $z \xbe U,$ so $ \xBc x,tf_{x} \xBe $ is minimal
in $ \xdz ' \xex U,$ contradiction.
$ \xcz $ (Claim \ref{Claim Tree-Repres-2} (page \pageref{Claim Tree-Repres-2}) 
and Proposition \ref{Proposition Pref-Complete-Trans} (page \pageref{Proposition
Pref-Complete-Trans}) )
\\[3ex]

% karl-search= End Claim Tree-Repres-2 Proof
\vspace{7mm}

% *************************************

\vspace{7mm}

\index{Proposition Equiv-Trans}

We give now the direct proof, which we cannot adapt to the smooth case.
Such
easy results must be part of the folklore, but we give them for
completeness'
sake.

\bp

$\hspace{0.01em}$

% (+++ Orig. No.:  Proposition Equiv-Trans +++)

\label{Proposition Equiv-Trans}

In the general case, every preferential structure is equivalent to a
transitive one - i.e. they have the same $ \xbm -$functions.
\index{Proposition Equiv-Trans Proof}

\ep

\subparagraph{
Proof
}

$\hspace{0.01em}$

% (+++ Orig.:  Proof +++)

If $ \xBc a,i \xBe  \xee  \xBc b,j \xBe,$ we create
an infinite descending chain of new copies $ \xBc b, \xBc j,a,i,n \xBe  \xBe,$
$n \xbe \xbo
,$ where
$ \xBc b, \xBc j,a,i,n \xBe  \xBe  \xee  \xBc b, \xBc j,a,i,n'  \xBe  \xBe $ if
$n' >n,$ and make $ \xBc a,i \xBe  \xee
\xBc b, \xBc j,a,i,n \xBe  \xBe $ for all
$n \xbe \xbo,$ but cancel the pair $ \xBc a,i \xBe  \xee  \xBc b,j \xBe $ from
the relation
(otherwise, we would
not have achieved anything), but $ \xBc b,j \xBe $ stays as element in the set.
Now, the relation is trivially transitive, and all these $ \xBc b, \xBc j,a,i,n
\xBe  \xBe $
just kill themselves, there is no need to minimize them by anything else.
We just continued $ \xBc a,i \xBe  \xee  \xBc b,j \xBe $ in a way it cannot
bother us. For the
$ \xBc b,j \xBe,$ we
do of course the same thing again. So, we have full equivalence, i.e. the
$ \xbm -$functions of both structures are identical (this is trivial to
see). $ \xcz $
\\[3ex]
\subsubsection{Smooth structures}
\label{Section 2.2.4}
%  Section (4.4):  Smooth structures
%  Section (4.4):  Smooth structures
% %
% ===================================
\paragraph{Introduction}
\label{Section 2.2.4.1}
%  Subsubsection (4.4.1):  Introduction
%  Subsubsection (4.4.1):  Introduction
% %
% ===================================
\paragraph{Cumulativity without $(\xcv)$}
\label{Section 2.2.4.2}
%  Subsubsection (4.4.2):  Cumulativity without (%v)
%  Subsubsection (4.4.2):  Cumulativity without (%v)
% %
% ================================================
\index{Comment Cum-Union}

We show here that, without sufficient closure
properties, there is an infinity of versions of cumulativity, which
collapse
to usual cumulativity when the domain is closed under finite unions.
Closure properties thus reveal themselves as a powerful tool to show
independence of properties.

We then show positive results for the smooth and the transitive
smooth case.

We work in some fixed arbitrary set $Z,$ all sets considered will be
subsets of $Z.$

Unless said otherwise, we use without further
mentioning $( \xbm PR)$ and $( \xbm \xcc ).$

Note that $( \xbm PR)$ and $( \xbm \xcc )$ entail $ \xbm (A \xcv B) \xcc
\xbm (A) \xcv \xbm (B)$ whenever
$ \xbm $ is defined for A, $B,$ $A \xcv B.$ $( \xbm (A \xcv B) \xcs A \xcc
\xbm (A),$ $ \xbm (A \xcv B) \xcs B \xcc \xbm (B),$
by $( \xbm PR),$ but $ \xbm (A \xcv B) \xcc A \xcv B$ by $( \xbm \xcc ).)$
\index{Definition Cum-Alpha}

\bd

$\hspace{0.01em}$

% (+++ Orig. No.:  Definition Cum-Alpha +++)

\label{Definition Cum-Alpha}

For any ordinal $ \xba,$ we define

$( \xbm Cum \xba ):$

If for all $ \xbb \xck \xba $ $ \xbm (X_{ \xbb }) \xcc U \xcv \xcV \{X_{
\xbg }: \xbg < \xbb \}$ hold, then so does
$ \xcS \{X_{ \xbg }: \xbg \xck \xba \} \xcs \xbm (U) \xcc \xbm (X_{ \xba
}).$

$( \xbm Cumt \xba ):$

If for all $ \xbb \xck \xba $ $ \xbm (X_{ \xbb }) \xcc U \xcv \xcV \{X_{
\xbg }: \xbg < \xbb \}$ hold, then so does
$X_{ \xba } \xcs \xbm (U) \xcc \xbm (X_{ \xba }).$

( `` $t$ '' stands for transitive, see Fact \ref{Fact Cum-Alpha} (page
\pageref{Fact Cum-Alpha}), (2.2)
below.)

$( \xbm Cum \xca )$ and $( \xbm Cumt \xca )$ will be the class of all $(
\xbm Cum \xba )$ or $( \xbm Cumt \xba )$ -
read their ``conjunction'', i.e. if we say that $( \xbm Cum \xca )$ holds,
we mean that
all $( \xbm Cum \xba )$ hold.
\index{Note Cum-Alpha}

\ed

Note:

The first conditions thus have the form:

$( \xbm Cum0)$ $ \xbm (X_{0}) \xcc U$ $ \xch $ $X_{0} \xcs \xbm (U) \xcc
\xbm (X_{0}),$

$( \xbm Cum1)$ $ \xbm (X_{0}) \xcc U,$ $ \xbm (X_{1}) \xcc U \xcv X_{0}$ $
\xch $ $X_{0} \xcs X_{1} \xcs \xbm (U) \xcc \xbm (X_{1}),$

$( \xbm Cum2)$ $ \xbm (X_{0}) \xcc U,$ $ \xbm (X_{1}) \xcc U \xcv X_{0},$
$ \xbm (X_{2}) \xcc U \xcv X_{0} \xcv X_{1}$ $ \xch $

$ \xDC $ $X_{0} \xcs X_{1} \xcs X_{2} \xcs \xbm (U) \xcc \xbm (X_{2}).$

$( \xbm Cumt \xba )$ differs from $( \xbm Cum \xba )$ only in the
consequence, the intersection contains
only the last $X_{ \xba }$ - in particular, $( \xbm Cum0)$ and $( \xbm
Cumt0)$ coincide.

Recall that condition $( \xbm Cum1)$ is the crucial condition in  \cite{Leh92a},
which
failed, despite $( \xbm CUM),$ but which has to hold in all smooth models.
This
condition $( \xbm Cum1)$ was the starting point of the investigation.

We briefly mention some major results on above conditions, taken from
Fact \ref{Fact Cum-Alpha} (page \pageref{Fact Cum-Alpha})  and
shown there - we use the same numbering:

(1.1) $( \xbm Cum \xba )$ $ \xch $ $( \xbm Cum \xbb )$ for all $ \xbb \xck
\xba $

(1.2) $( \xbm Cumt \xba )$ $ \xch $ $( \xbm Cumt \xbb )$ for all $ \xbb
\xck \xba $

(2.1) All $( \xbm Cum \xba )$ hold in smooth preferential structures

(2.2) All $( \xbm Cumt \xba )$ hold in transitive smooth preferential
structures

(3.1) $( \xbm Cum \xbb )$ $+$ $( \xcv )$ $ \xch $ $( \xbm Cum \xba )$ for
all $ \xbb \xck \xba $

(3.2) $( \xbm Cumt \xbb )$ $+$ $( \xcv )$ $ \xch $ $( \xbm Cumt \xba )$
for all $ \xbb \xck \xba $

(5.2) $( \xbm Cum \xba )$ $ \xch $ $( \xbm CUM)$ for all $ \xba $

(5.3) $( \xbm CUM)$ $+$ $( \xcv )$ $ \xch $ $( \xbm Cum \xba )$ for all $
\xba $
\index{Definition HU-All}

The following inductive definition of $H(U,u)$ and of the property $ \xCf
(HU,u)$
concerns closure under $( \xbm Cum \xca ),$ its main
property is formulated in Fact \ref{Fact HU} (page \pageref{Fact HU}), its main
interest
is its use in the proof of Proposition \ref{Proposition D-4.4.6} (page
\pageref{Proposition D-4.4.6}).

\bd

$\hspace{0.01em}$

% (+++ Orig. No.:  Definition HU-All +++)

\label{Definition HU-All}

$(H(U,u)_{ \xba },$ $H(U)_{ \xba },$ $ \xCf (HU,u),$ $ \xCf (HU).)$

$H(U,u)_{0}$ $:=$ $U,$

$H(U,u)_{ \xba +1}$ $:=$ $H(U,u)_{ \xba }$ $ \xcv $ $ \xcV \{X:$ $u \xbe
X$ $ \xcu $ $ \xbm (X) \xcc H(U,u)_{ \xba }\},$

$H(U,u)_{ \xbl }$ $:=$ $ \xcV \{H(U,u)_{ \xba }: \xba < \xbl \}$ for
$limit( \xbl ),$

$H(U,u)$ $:=$ $ \xcV \{H(U,u)_{ \xba }: \xba < \xbk \}$ for $ \xbk $
sufficiently big $(card(Z)$ suffices, as

the procedure trivializes, when we cannot add any new elements).

$ \xCf (HU,u)$ is the property:

$u \xbe \xbm (U),$ $u \xbe Y- \xbm (Y)$ $ \xch $ $ \xbm (Y) \xcC H(U,u)$ -
of course for all $u$ and $U.$
$(U,Y \xbe \xdy ).$

Thus, $(HU,u)$ entails $ \xbm (U) \xcc H(U,u),$ $u \xbe \xbm (U) \xcs Y$ $
\xch $ $u \xbe \xbm (Y).$

For the case with $( \xcv ),$ we further define, independent of $u,$

$H(U)_{0}$ $:=$ $U,$

$H(U)_{ \xba +1}$ $:=$ $H(U)_{ \xba }$ $ \xcv $ $ \xcV \{X:$ $ \xbm (X)
\xcc H(U)_{ \xba }\},$

$H(U)_{ \xbl }$ $:=$ $ \xcV \{H(U)_{ \xba }: \xba < \xbl \}$ for $limit(
\xbl ),$

$H(U)$ $:=$ $ \xcV \{H(U)_{ \xba }: \xba < \xbk \}$ again for $ \xbk $
sufficiently big

$ \xCf (HU)$ is the property:

$u \xbe \xbm (U),$ $u \xbe Y- \xbm (Y)$ $ \xch $ $ \xbm (Y) \xcC H(U)$ -
of course for all $U.$
$(U,Y \xbe \xdy ).$

Thus, (HU) entails $ \xbm (Y) \xcc H(U)$ $ \xch $ $ \xbm (U) \xcs Y \xcc
\xbm (Y).$

\ed

Obviously, $H(U,u) \xcc H(U),$ so $(HU) \xch (HU,u).$
\index{Example Inf-Cum-Alpha}

\be

$\hspace{0.01em}$

% (+++ Orig. No.:  Example Inf-Cum-Alpha +++)

\label{Example Inf-Cum-Alpha}

This important example shows that the conditions $( \xbm Cum \xba )$ and
$( \xbm Cumt \xba )$ defined in Definition \ref{Definition Cum-Alpha} (page
\pageref{Definition Cum-Alpha})
are all different in the absence of $( \xcv ),$ in its
presence they all collapse (see Fact \ref{Fact Cum-Alpha} (page \pageref{Fact
Cum-Alpha})  below).
More precisely,
the following (class of) examples shows that the $( \xbm Cum \xba )$
increase in
strength. For any finite or infinite ordinal $ \xbk >0$ we construct an
example s.t.

(a) $( \xbm PR)$ and $( \xbm \xcc )$ hold

(b) $( \xbm CUM)$ holds

(c) $( \xcS )$ holds

(d) $( \xbm Cumt \xba )$ holds for $ \xba < \xbk $

(e) $( \xbm Cum \xbk )$ fails.

\ee

\subparagraph{
Proof
}

$\hspace{0.01em}$

% (+++ Orig.:  Proof +++)

We define a suitable base set and a non-transitive binary relation $ \xeb
$
on this set, as well as a suitable set $ \xdx $ of subsets, closed under
arbitrary
intersections, but not under finite unions, and define $ \xbm $ on these
subsets
as usual in preferential structures by $ \xeb.$ Thus, $( \xbm PR)$ and $(
\xbm \xcc )$ will hold.
It will be immediate that $( \xbm Cum \xbk )$ fails, and we will show that
$( \xbm CUM)$ and
$( \xbm Cumt \xba )$ for $ \xba < \xbk $ hold by examining the cases.

For simplicity, we first define a set of generators for $ \xdx,$ and
close under
$( \xcS )$ afterwards. The set $U$ will have a special position, it is the
``useful''
starting point to construct chains corresponding to above definitions of
$( \xbm Cum \xba )$ and $( \xbm Cumt \xba ).$

In the sequel,
$i,j$ will be successor ordinals, $ \xbl $ etc. limit ordinals, $ \xba,$
$ \xbb,$ $ \xbk $ any ordinals,
thus e.g. $ \xbl \xck \xbk $ will imply that $ \xbl $ is a limit ordinal $
\xck \xbk,$ etc.

The base set and the relation $ \xeb:$

$ \xbk >0$ is fixed, but arbitrary. We go up to $ \xbk >0.$

The base set is $\{a,b,c\}$ $ \xcv $ $\{d_{ \xbl }: \xbl \xck \xbk \}$ $
\xcv $ $\{x_{ \xba }: \xba \xck \xbk +1\}$ $ \xcv $ $\{x'_{ \xba }: \xba
\xck \xbk \}.$
$a \xeb b \xeb c,$ $x_{ \xba } \xeb x_{ \xba +1},$ $x_{ \xba } \xeb x'_{
\xba },$ $x'_{0} \xeb x_{ \xbl }$ (for any $ \xbl )$ - $ \xeb $ is NOT
transitive.

The generators:

$U:=\{a,c,x_{0}\} \xcv \{d_{ \xbl }: \xbl \xck \xbk \}$ - i.e. $ \Xl
.\{d_{ \xbl }:lim( \xbl ) \xcu \xbl \xck \xbk \},$

$X_{i}:=\{c,x_{i},x'_{i},x_{i+1}\}$ $(i< \xbk ),$

$X_{ \xbl }:=\{c,d_{ \xbl },x_{ \xbl },x'_{ \xbl },x_{ \xbl +1}\} \xcv
\{x'_{ \xba }: \xba < \xbl \}$ $( \xbl < \xbk ),$

$X'_{ \xbk }:=\{a,b,c,x_{ \xbk },x'_{ \xbk },x_{ \xbk +1}\}$ if $ \xbk $
is a successor,

$X'_{ \xbk }:=\{a,b,c,d_{ \xbk },x_{ \xbk },x'_{ \xbk },x_{ \xbk +1}\}
\xcv \{x'_{ \xba }: \xba < \xbk \}$ if $ \xbk $ is a limit.

Thus, $X'_{ \xbk }=X_{ \xbk } \xcv \{a,b\}$ if $X_{ \xbk }$ were defined.

Note that there is only one $X'_{ \xbk },$ and $X_{ \xba }$ is defined
only for $ \xba < \xbk,$ so we will
not have $X_{ \xba }$ and $X'_{ \xba }$ at the same time.

Thus, the values of the generators under $ \xbm $ are:

$ \xbm (U)=U,$

$ \xbm (X_{i})=\{c,x_{i}\},$

$ \xbm (X_{ \xbl })=\{c,d_{ \xbl }\} \xcv \{x'_{ \xba }: \xba < \xbl \},$

$ \xbm (X'_{i})=\{a,x_{i}\}$ $(i>0,$ $i$ has to be a successor),

$ \xbm (X'_{ \xbl })=\{a,d_{ \xbl }\} \xcv \{x'_{ \xba }: \xba < \xbl \}.$

(We do not assume that the domain is closed under $ \xbm.)$

Intersections:

We consider first pairwise intersections:

(1) $U \xcs X_{0}=\{c,x_{0}\},$

(2) $U \xcs X_{i}=\{c\},$ $i>0,$

(3) $U \xcs X_{ \xbl }=\{c,d_{ \xbl }\},$

(4) $U \xcs X'_{i}=\{a,c\}$ $(i>0),$

(5) $U \xcs X'_{ \xbl }=\{a,c,d_{ \xbl }\},$

(6) $X_{i} \xcs X_{j}:$

(6.1) $j=i+1$ $\{c,x_{i+1}\},$

(6.2) else $\{c\},$

(7) $X_{i} \xcs X_{ \xbl }:$

(7.1) $i< \xbl $ $\{c,x'_{i}\},$

(7.2) $i= \xbl +1$ $\{c,x_{ \xbl +1}\},$

(7.3) $i> \xbl +1$ $\{c\},$

(8) $X_{ \xbl } \xcs X_{ \xbl ' }:$ $\{c\} \xcv \{x'_{ \xba }: \xba \xck
min( \xbl, \xbl ' )\}.$

As $X'_{ \xbk }$ occurs only once, $X_{ \xba } \xcs X'_{ \xbk }$ etc. give
no new results.

Note that $ \xbm $ is constant on all these pairwise intersections.

Iterated intersections:

As $c$ is an element of all sets, sets of the type $\{c,z\}$ do not give
any
new results. The possible subsets of $\{a,c,d_{ \xbl }\}:$ $\{c\},$
$\{a,c\},$ $\{c,d_{ \xbl }\}$ exist
already. Thus, the only source of new sets via iterated intersections is
$X_{ \xbl } \xcs X_{ \xbl ' }=\{c\} \xcv \{x'_{ \xba }: \xba \xck min(
\xbl, \xbl ' )\}.$ But, to intersect them, or with some
old sets, will not generate any new sets either. Consequently, the example
satisfies $( \xcS )$ for $ \xdx $ defined by $U,$ $X_{i}$ $(i< \xbk ),$
$X_{ \xbl }$ $( \xbl < \xbk ),$ $X'_{ \xbk },$ and above
paiwise intersections.

We will now verify the positive properties. This is tedious, but
straightforward, we have to check the different cases.

Validity of $( \xbm CUM):$

Consider the prerequisite $ \xbm (X) \xcc Y \xcc X.$ If $ \xbm (X)=X$ or
if $X- \xbm (X)$ is a singleton,
$X$ cannot give a violation of $( \xbm CUM).$ So we are left with the
following
candidates for $X:$

(1) $X_{i}:=\{c,x_{i},x'_{i},x_{i+1}\},$ $ \xbm (X_{i})=\{c,x_{i}\}$

Interesting candidates for $Y$ will have 3 elements, but they will all
contain
a. (If $ \xbk < \xbo:$ $U=\{a,c,x_{0}\}.)$

(2) $X_{ \xbl }:=\{c,d_{ \xbl },x_{ \xbl },x'_{ \xbl },x_{ \xbl +1}\} \xcv
\{x'_{ \xba }: \xba < \xbl \},$ $ \xbm (X_{ \xbl })=\{c,d_{ \xbl }\} \xcv
\{x'_{ \xba }: \xba < \xbl \}$

The only sets to contain $d_{ \xbl }$ are $X_{ \xbl },$ $U,$ $U \xcs X_{
\xbl }.$ But $a \xbe U,$ and
$U \xcs X_{ \xbl }$ is finite. $(X_{ \xbl }$ and $X'_{ \xbl }$ cannot be
present at the same time.)

(3) $X'_{i}:=\{a,b,c,x_{i},x'_{i},x_{i+1}\},$ $ \xbm (X'_{i})=\{a,x_{i}\}$

a is only in $U,$ $X'_{i},$ $U \xcs X'_{i}=\{a,c\},$ but $x_{i} \xce U,$
as $i>0.$

(4) $X'_{ \xbl }:=\{a,b,c,d_{ \xbl },x_{ \xbl },x'_{ \xbl },x_{ \xbl +1}\}
\xcv \{x'_{ \xba }: \xba < \xbl \},$ $ \xbm (X'_{ \xbl })=\{a,d_{ \xbl }\}
\xcv \{x'_{ \xba }: \xba < \xbl \}$

$d_{ \xbl }$ is only in $X'_{ \xbl }$ and $U,$ but $U$ contains no $x'_{
\xba }.$

Thus, $( \xbm CUM)$ holds trivially.

$( \xbm Cumt \xba )$ hold for $ \xba < \xbk:$

To simplify language, we say that we reach $Y$ from $X$ iff $X \xEd Y$ and
there is a
sequence $X_{ \xbb },$ $ \xbb \xck \xba $ and $ \xbm (X_{ \xbb }) \xcc X
\xcv \xcV \{X_{ \xbg }: \xbg < \xbb \},$ and $X_{ \xba }=Y,$ $X_{0}=X.$
Failure of $( \xbm Cumt \xba )$ would then mean that there are $X$ and
$Y,$ we can reach
$Y$ from $X,$ and $x \xbe ( \xbm (X) \xcs Y)- \xbm (Y).$ Thus, in a
counterexample, $Y= \xbm (Y)$ is
impossible, so none of the intersections can be such $Y.$

To reach $Y$ from $X,$ we have to get started from $X,$ i.e. there must be
$Z$ s.t.
$ \xbm (Z) \xcc X,$ $Z \xcC X$ (so $ \xbm (Z) \xEd Z).$ Inspection of the
different cases shows that
we cannot reach any set $Y$ from any case of the intersections, except
from
(1), (6.1), (7.2).

If $Y$ contains a globally minimal element (i.e. there is no smaller
element in
any set), it can only be reached from any $X$ which already contains this
element. The globally minimal elements are a, $x_{0},$ and the $d_{ \xbl
},$ $ \xbl \xck \xbk.$

By these observations, we see that $X_{ \xbl }$ and $X'_{ \xbk }$ can only
be reached from $U.$
From no $X_{ \xba }$ $U$ can be reached, as the globally minimal a is
missing. But
$U$ cannot be reached from $X'_{ \xbk }$ either, as the globally minimal
$x_{0}$ is missing.

When we look at the relation $ \xeb $ defining $ \xbm,$ we see that we
can reach $Y$ from $X$
only by going upwards, adding bigger elements. Thus, from $X_{ \xba },$ we
cannot reach
any $X_{ \xbb },$ $ \xbb < \xba,$ the same holds for $X'_{ \xbk }$ and
$X_{ \xbb },$ $ \xbb < \xbk.$ Thus, from $X'_{ \xbk },$ we
cannot go anywhere interesting (recall that the intersections are not
candidates
for a $Y$ giving a contradiction).

Consider now $X_{ \xba }.$ We can go up to any $X_{ \xba +n},$ but not to
any $X_{ \xbl },$ $ \xba < \xbl,$ as
$d_{ \xbl }$ is missing, neither to $X'_{ \xbk },$ as a is missing. And we
will be stopped by
the first $ \xbl > \xba,$ as $x_{ \xbl }$ will be missing to go beyond
$X_{ \xbl }.$ Analogous observations
hold for the remaining intersections (1), (6.1), (7.2). But in all these
sets we
can reach, we will not destroy minimality of any element of $X_{ \xba }$
(or of the
intersections).

Consequently, the only candidates for failure will all start with $U.$ As
the only
element of $U$ not globally minimal is $c,$ such failure has to have $c
\xbe Y- \xbm (Y),$ so
$Y$ has to be $X'_{ \xbk }.$ Suppose we omit one of the $X_{ \xba }$ in
the sequence going up to
$X'_{ \xbk }.$ If $ \xbk \xcg \xbl > \xba,$ we cannot reach $X_{ \xbl }$
and beyond, as $x'_{ \xba }$ will be missing.
But we cannot go to $X_{ \xba +n}$ either, as $x_{ \xba +1}$ is missing.
So we will be stopped
at $X_{ \xba }.$ Thus, to see failure, we need the full sequence
$U=X_{0},$ $X'_{ \xbk }=Y_{ \xbk },$
$Y_{ \xba }=X_{ \xba }$ for $0< \xba < \xbk.$

$( \xbm Cum \xbk )$ fails:

The full sequence $U=X_{0},$ $X'_{ \xbk }=Y_{ \xbk },$ $Y_{ \xba }=X_{
\xba }$ for $0< \xba < \xbk $ shows this, as
$c \xbe \xbm (U) \xcs X'_{ \xbk },$ but $c \xce \xbm (X'_{ \xbk }).$

Consequently, the example satisfies $( \xcS ),$ $( \xbm CUM),$ $( \xbm
Cumt \xba )$ for $ \xba < \xbk,$ and
$( \xbm Cum \xbk )$ fails.

$ \xcz $
\\[3ex]
\index{Fact Cum-Alpha}

\bfa

$\hspace{0.01em}$

% (+++ Orig. No.:  Fact Cum-Alpha +++)

\label{Fact Cum-Alpha}

We summarize some properties of $( \xbm Cum \xba )$ and $( \xbm Cumt \xba
)$ - sometimes with some
redundancy. Unless said otherwise, $ \xba,$ $ \xbb $ etc. will be
arbitrary ordinals.

For (1) to (6) $( \xbm PR)$ and $( \xbm \xcc )$ are assumed to hold, for
(7) only
$( \xbm \xcc ).$

(1) Downward:

(1.1) $( \xbm Cum \xba )$ $ \xch $ $( \xbm Cum \xbb )$ for all $ \xbb \xck
\xba $

(1.2) $( \xbm Cumt \xba )$ $ \xch $ $( \xbm Cumt \xbb )$ for all $ \xbb
\xck \xba $

(2) Validity of $( \xbm Cum \xba )$ and $( \xbm Cumt \xba )$:

(2.1) All $( \xbm Cum \xba )$ hold in smooth preferential structures

(2.2) All $( \xbm Cumt \xba )$ hold in transitive smooth preferential
structures

(2.3) $( \xbm Cumt \xba )$ for $0< \xba $ do not necessarily hold in
smooth structures without
transitivity, even in the presence of $( \xcS )$

(3) Upward:

(3.1) $( \xbm Cum \xbb )$ $+$ $( \xcv )$ $ \xch $ $( \xbm Cum \xba )$ for
all $ \xbb \xck \xba $

(3.2) $( \xbm Cumt \xbb )$ $+$ $( \xcv )$ $ \xch $ $( \xbm Cumt \xba )$
for all $ \xbb \xck \xba $

(3.3) $\{( \xbm Cumt \xbb ): \xbb < \xba \}$ $+$ $( \xbm CUM)$ $+$ $( \xcS
)$ $ \xcP $ $( \xbm Cum \xba )$ for $ \xba >0.$

(4) Connection $( \xbm Cum \xba )/( \xbm Cumt \xba )$:

(4.1) $( \xbm Cumt \xba )$ $ \xch $ $( \xbm Cum \xba )$

(4.2) $( \xbm Cum \xba )$ $+$ $( \xcS )$ $ \xcP $ $( \xbm Cumt \xba )$

(4.3) $( \xbm Cum \xba )$ $+$ $( \xcv )$ $ \xch $ $( \xbm Cumt \xba )$

(5) $( \xbm CUM)$ and $( \xbm Cumi)$:

(5.1) $( \xbm CUM)$ $+$ $( \xcv )$ entail:

(5.1.1) $ \xbm (A) \xcc B$ $ \xch $ $ \xbm (A \xcv B)= \xbm (B)$

(5.1.2) $ \xbm (X) \xcc U,$ $U \xcc Y$ $ \xch $ $ \xbm (Y \xcv X)= \xbm
(Y)$

(5.1.3) $ \xbm (X) \xcc U,$ $U \xcc Y$ $ \xch $ $ \xbm (Y) \xcs X \xcc
\xbm (U)$

(5.2) $( \xbm Cum \xba )$ $ \xch $ $( \xbm CUM)$ for all $ \xba $

(5.3) $( \xbm CUM)$ $+$ $( \xcv )$ $ \xch $ $( \xbm Cum \xba )$ for all $
\xba $

(5.4) $( \xbm CUM)$ $+$ $( \xcs )$ $ \xch $ $( \xbm Cum0)$

(6) $( \xbm CUM)$ and $( \xbm Cumt \xba )$:

(6.1) $( \xbm Cumt \xba )$ $ \xch $ $( \xbm CUM)$ for all $ \xba $

(6.2) $( \xbm CUM)$ $+$ $( \xcv )$ $ \xch $ $( \xbm Cumt \xba )$ for all $
\xba $

(6.3) $( \xbm CUM)$ $ \xcP $ $( \xbm Cumt \xba )$ for all $ \xba >0$

(7) $( \xbm Cum0)$ $ \xch $ $( \xbm PR)$
\index{Fact Cum-Alpha Proof}

\efa

\subparagraph{
Proof
}

$\hspace{0.01em}$

% (+++ Orig.:  Proof +++)

We prove these facts in a different order: (1), (2), (5.1), (5.2), (4.1),
(6.1),
(6.2), (5.3), (3.1), (3.2), (4.2), (4.3), (5.4), (3.3), (6.3), (7).

(1.1)

For $ \xbb < \xbg \xck \xba $ set $X_{ \xbg }:=X_{ \xbb }.$ Let the
prerequisites of $( \xbm Cum \xbb )$ hold. Then for
$ \xbg $ with $ \xbb < \xbg \xck \xba $ $ \xbm (X_{ \xbg }) \xcc X_{ \xbb
}$ by $( \xbm \xcc ),$ so the prerequisites
of $( \xbm Cum \xba )$ hold, too, so by $( \xbm Cum \xba )$ $ \xcS \{X_{
\xbd }: \xbd \xck \xbb \} \xcs \xbm (U)$ $=$
$ \xcS \{X_{ \xbd }: \xbd \xck \xba \} \xcs \xbm (U)$ $ \xcc $ $ \xbm (X_{
\xba })$ $=$ $ \xbm (X_{ \xbb }).$

(1.2)

Analogous.

(2.1)

Proof by induction.

$( \xbm Cum0)$ Let $ \xbm (X_{0}) \xcc U,$ suppose there is $x \xbe \xbm
(U) \xcs (X_{0}- \xbm (X_{0})).$ By smoothness,
there is $y \xeb x,$ $y \xbe \xbm (X_{0}) \xcc U,$ $contradiction$ (The
same arguments works for copies: all
copies of $x$ must be minimized by some $y \xbe \xbm (X_{0}),$ but at
least one copy of $x$
has to be minimal in $U.)$

Suppose $( \xbm Cum \xbb )$ hold for all $ \xbb < \xba.$ We show $( \xbm
Cum \xba ).$ Let the prerequisites
of $( \xbm Cum \xba )$ hold, then those for $( \xbm Cum \xbb ),$ $ \xbb <
\xba $ hold, too. Suppose there is
$x \xbe \xbm (U) \xcs \xcS \{X_{ \xbg }: \xbg \xck \xba \}- \xbm (X_{ \xba
}).$ So by $( \xbm Cum \xbb )$ for $ \xbb < \xba $ $x \xbe \xbm (X_{ \xbb
})$
for all $ \xbb < \xba,$
moreover $x \xbe \xbm (U).$ By smoothness, there is $y \xbe \xbm (X_{ \xba
}) \xcc U \xcv \xcV \{X_{ \xbb ' }: \xbb ' < \xba \},$ $y \xeb x,$
but this is a contradiction. The same argument works again for copies.

(2.2)

We use the following Fact:
Let, in a smooth transitive structure, $ \xbm (X_{ \xbb })$ $ \xcc $ $U
\xcv \xcV \{X_{ \xbg }: \xbg < \xbb \}$ for all $ \xbb \xck \xba,$
and let $x \xbe \xbm (U).$ Then there is no $y \xeb x,$ $y \xbe U \xcv
\xcV \{X_{ \xbg }: \xbg \xck \xba \}.$

Proof of the Fact by induction:
Suppose such $y \xbe U \xcv X_{0}$ exists. $y \xbe U$ is impossible. Let
$y \xbe X_{0},$ by
$ \xbm (X_{0}) \xcc U,$ $y \xbe X_{0}- \xbm (X_{0}),$ so
there is $z \xbe \xbm (X_{0}),$ $z \xeb y,$ so $z \xeb x$ by transitivity,
but $ \xbm (X_{0}) \xcc U.$
Let the result hold for all $ \xbb < \xba,$ but fail for $ \xba,$
so $ \xCN \xcE y \xeb x.y \xbe U \xcv \xcV \{X_{ \xbg }: \xbg < \xba \},$
but $ \xcE y \xeb x.y \xbe U \xcv \xcV \{X_{ \xbg }: \xbg \xck \xba \},$
so $y \xbe X_{ \xba }.$
If $y \xbe \xbm (X_{ \xba }),$ then $y \xbe U \xcv \xcV \{X_{ \xbg }: \xbg
< \xba \},$ but this is impossible, so
$y \xbe X_{ \xba }- \xbm (X_{ \xba }),$ let by smoothness $z \xeb y,$ $z
\xbe \xbm (X_{ \xba }),$ so by transitivity $z \xeb x,$ $contradiction.$
The result is easily modified for the case with copies.

Let the prerequisites of $( \xbm Cumt \xba )$ hold, then those of the Fact
will hold,
too. Let now $x \xbe \xbm (U) \xcs (X_{ \xba }- \xbm (X_{ \xba })),$ by
smoothness, there must be $y \xeb x,$
$y \xbe \xbm (X_{ \xba }) \xcc U \xcv \xcV \{X_{ \xbg }: \xbg < \xba \},$
contradicting the Fact.

(2.3)

Let $ \xba >0,$ and consider the following structure over $\{a,b,c\}:$
$U:=\{a,c\},$
$X_{0}:=\{b,c\},$ $X_{ \xba }:= \Xl:=X_{1}:=\{a,b\},$ and their
intersections, $\{a\},$ $\{b\},$ $\{c\},$ $ \xCQ $ with
the order $c \xeb b \xeb a$ (without transitivity). This is preferential,
so $( \xbm PR)$ and
$( \xbm \xcc )$ hold.
The structure is smooth for $U,$ all $X_{ \xbb },$ and their
intersections.
We have $ \xbm (X_{0}) \xcc U,$ $ \xbm (X_{ \xbb }) \xcc U \xcv X_{0}$ for
all $ \xbb \xck \xba,$ so $ \xbm (X_{ \xbb }) \xcc U \xcv \xcV \{X_{ \xbg
}: \xbg < \xbb \}$
for all $ \xbb \xck \xba $ but $X_{ \xba } \xcs \xbm (U)=\{a\} \xcC \{b\}=
\xbm (X_{ \xba })$ for $ \xba >0.$

(5.1)

(5.1.1) $ \xbm (A) \xcc B$ $ \xch $ $ \xbm (A \xcv B) \xcc \xbm (A) \xcv
\xbm (B) \xcc B$ $ \xch_{( \xbm CUM)}$ $ \xbm (B)= \xbm (A \xcv B).$

(5.1.2) $ \xbm (X) \xcc U \xcc Y$ $ \xch $ (by (5.1.1)) $ \xbm (Y \xcv X)=
\xbm (Y).$

(5.1.3) $ \xbm (Y) \xcs X$ $=$ (by (5.1.2)) $ \xbm (Y \xcv X) \xcs X$ $
\xcc $ $ \xbm (Y \xcv X) \xcs (X \xcv U)$ $ \xcc $ (by $( \xbm PR))$
$ \xbm (X \xcv U)$ $=$ (by (5.1.1)) $ \xbm (U).$

(5.2)

Using (1.1), it suffices to show $( \xbm Cum0)$ $ \xch $ $( \xbm CUM).$
Let $ \xbm (X) \xcc U \xcc X.$ By $( \xbm Cum0)$ $X \xcs \xbm (U) \xcc
\xbm (X),$ so by $ \xbm (U) \xcc U \xcc X$ $ \xch $ $ \xbm (U) \xcc \xbm
(X).$
$U \xcc X$ $ \xch $ $ \xbm (X) \xcs U \xcc \xbm (U)$ by $( \xbm PR),$ but
also $ \xbm (X) \xcc U,$ so $ \xbm (X) \xcc \xbm (U).$

(4.1)

Trivial.

(6.1)

Follows from (4.1) and (5.2).

(6.2)

Let the prerequisites of $( \xbm Cumt \xba )$ hold.

We first show by induction $ \xbm (X_{ \xba } \xcv U) \xcc \xbm (U).$

Proof:

$ \xba =0:$ $ \xbm (X_{0}) \xcc U$ $ \xch $ $ \xbm (X_{0} \xcv U)= \xbm
(U)$ by (5.1.1).
Let for all $ \xbb < \xba $ $ \xbm (X_{ \xbb } \xcv U) \xcc \xbm (U) \xcc
U.$ By prerequisite,
$ \xbm (X_{ \xba }) \xcc U \xcv \xcV \{X_{ \xbb }: \xbb < \xba \},$ thus $
\xbm (X_{ \xba } \xcv U)$ $ \xcc $ $ \xbm (X_{ \xba }) \xcv \xbm (U)$ $
\xcc $ $ \xcV \{U \xcv X_{ \xbb }: \xbb < \xba \},$

moreover for all $ \xbb < \xba $ $ \xbm (X_{ \xbb } \xcv U) \xcc U \xcc
X_{ \xba } \xcv U,$
so $ \xbm (X_{ \xba } \xcv U) \xcs (U \xcv X_{ \xbb })$ $ \xcc $ $ \xbm
(U)$ by (5.1.3), thus $ \xbm (X_{ \xba } \xcv U) \xcc \xbm (U).$

Consequently, under the above prerequisites, we have $ \xbm (X_{ \xba }
\xcv U)$ $ \xcc $ $ \xbm (U)$ $ \xcc $
$U$ $ \xcc $ $U \xcv X_{ \xba },$ so by $( \xbm CUM)$ $ \xbm (U)= \xbm
(X_{ \xba } \xcv U),$ and, finally,
$ \xbm (U) \xcs X_{ \xba }= \xbm (X_{ \xba } \xcv U) \xcs X_{ \xba } \xcc
\xbm (X_{ \xba })$ by $( \xbm PR).$

Note that finite unions take us over the limit step, essentially, as all
steps collapse, and $ \xbm (X_{ \xba } \xcv U)$ will always be $ \xbm
(U),$ so there are no real
changes.

(5.3)

Follows from (6.2) and (4.1).

(3.1)

Follows from (5.2) and (5.3).

(3.2)

Follows from (6.1) and (6.2).

(4.2)

Follows from (2.3) and (2.1).

(4.3)

Follows from (5.2) and (6.2).

(5.4)

$ \xbm (X) \xcc U$ $ \xch $ $ \xbm (X) \xcc U \xcs X \xcc X$ $ \xch $ $
\xbm (X \xcs U)= \xbm (X)$ $ \xch $
$X \xcs \xbm (U)=(X \xcs U) \xcs \xbm (U) \xcc \xbm (X \xcs U)= \xbm (X)$

(3.3)

See Example \ref{Example Inf-Cum-Alpha} (page \pageref{Example Inf-Cum-Alpha}).

(6.3)

See Example \ref{Example Inf-Cum-Alpha} (page \pageref{Example Inf-Cum-Alpha}).

(7)

Trivial. Let $X \xcc Y,$ so by $( \xbm \xcc )$ $ \xbm (X) \xcc X \xcc Y,$
so by $( \xbm Cum0)$ $X \xcs \xbm (Y) \xcc \xbm (X).$

$ \xcz $
\\[3ex]
\index{Fact Cum-Alpha-HU}

\bfa

$\hspace{0.01em}$

% (+++ Orig. No.:  Fact Cum-Alpha-HU +++)

\label{Fact Cum-Alpha-HU}

Assume $( \xbm \xcc ).$

We have for $( \xbm Cum \xca )$ and $ \xCf (HU,u)$:

(1) $x \xbe \xbm (Y),$ $ \xbm (Y) \xcc H(U,x)$ $ \xch $ $Y \xcc H(U,x)$

(2) $( \xbm Cum \xca )$ $ \xch $ $ \xCf (HU,u)$

(3) $ \xCf (HU,u)$ $ \xch $ $( \xbm Cum \xca )$
\index{Fact Cum-Alpha-HU Proof}

\efa

\subparagraph{
Proof
}

$\hspace{0.01em}$

% (+++ Orig.:  Proof +++)

(1)

Trivial by definition of $H(U,x).$

(2)

Let $x \xbe \xbm (U),$ $x \xbe Y,$ $ \xbm (Y) \xcc H(U,x)$ (and thus $Y
\xcc H(U,x)$ by definition).
Thus, we have a sequence $X_{0}:=U,$ $ \xbm (X_{ \xbb }) \xcc U \xcv \xcV
\{X_{ \xbg }: \xbg < \xbb \},$ $x \xbe X_{ \xbb },$
and $Y=X_{ \xba }$ for some $ \xba $
(after $X_{0},$ enumerate arbitrarily $H(U,x)_{1},$ then $H(U,x)_{2},$
etc., do nothing at
limits). So $x \xbe \xcS \{X_{ \xbg }: \xbg \xck \xba \} \xcs \xbm (U)
\xcc \xbm (X_{ \xba })= \xbm (Y)$ by $( \xbm Cum \xca ).$

Remark: The same argument shows that we can replace `` $x \xbe X$ ''
equivalently by
`` $x \xbe \xbm (X)$ '' in the definition of $H(U,x)_{ \xba +1},$ as was
done in Definition 3.7.5
in  \cite{Sch04}.

(3)

Suppose $( \xbm Cum \xba )$ fails, we show that then so does $ \xCf
(HU,u)$ for $u=x.$
As $( \xbm Cum \xba )$ fails, for
all $ \xbb \xck \xba $ $ \xbm (X_{ \xbb }) \xcc U \xcv \xcV \{X_{ \xbg }:
\xbg < \xbb \},$ but there is $x \xbe \xcS \{X_{ \xbg }: \xbg \xck \xba \}
\xcs \xbm (U),$
$x \xce \xbm (X_{ \xba }).$ Thus for all $ \xbb \xck \xba $ $ \xbm (X_{
\xbb }) \xcc X_{ \xbb } \xcc H(U,x),$ moreover $x \xbe \xbm (U),$
$x \xbe X_{ \xba }- \xbm (X_{ \xba }),$ but $ \xbm (X_{ \xba }) \xcc
H(U,x),$ so $ \xCf (HU,u)$ fails for $u=x.$

$ \xcz $
\\[3ex]
\index{Fact HU}

\bfa

$\hspace{0.01em}$

% (+++ Orig. No.:  Fact HU +++)

\label{Fact HU}

We continue to show results for $H(U)$ and $H(U,u).$

Let A, $X,$ $U,$ $U',$ $Y$ and all $A_{i}$ be in $ \xdy.$

(0) $H(U)$ and $H(U,u)$

(0.1) $H(U,u) \xcc H(U),$

(0.2) $(HU) \xch (HU,u),$

(0.3) $( \xcv )$ $+$ $( \xbm PR)$ entail $H(U) \xcc H(U,u)$ for $u \xbe
\xbm (U),$

(0.4) $( \xcv )$ $+$ $( \xbm PR)$ entail $(HU,u) \xch (HU),$

(1) $( \xbm \xcc )$ and $ \xCf (HU)$ entail:

(1.1) $( \xbm PR),$

(1.2) $( \xbm CUM),$

(3) $( \xbm \xcc )$ and $( \xbm PR)$ entail:

(3.1) $A= \xcV \{A_{i}:i \xbe I\}$ $ \xch $ $ \xbm (A) \xcc \xcV \{ \xbm
(A_{i}):i \xbe I\},$

(3.2) $U \xcc H(U),$ and $U \xcc U' \xch H(U) \xcc H(U' ),$

(3.3) $ \xbm (U \xcv Y)-H(U) \xcc \xbm (Y)$ - if $ \xbm (U \xcv Y)$ is
defined, in particular, if $( \xcv )$
holds.

(4) $( \xcv ),$ $( \xbm \xcc ),$ $( \xbm PR),$ $( \xbm CUM)$ entail:

(4.1) $H(U)=H_{1}(U),$

(4.2) $U \xcc A,$ $ \xbm (A) \xcc H(U)$ $ \xch $ $ \xbm (A) \xcc U,$

(4.3) $ \xbm (Y) \xcc H(U)$ $ \xch $ $Y \xcc H(U)$ and $ \xbm (U \xcv Y)=
\xbm (U),$

(4.4) $x \xbe \xbm (U),$ $x \xbe Y- \xbm (Y)$ $ \xch $ $Y \xcC H(U)$ (and
thus $ \xCf (HU)),$

(4.5) $Y \xcC H(U)$ $ \xch $ $ \xbm (U \xcv Y) \xcC H(U).$

(5) $( \xcv ),$ $( \xbm \xcc ),$ $ \xCf (HU)$ entail

(5.1) $H(U)=H_{1}(U),$

(5.2) $U \xcc A,$ $ \xbm (A) \xcc H(U)$ $ \xch $ $ \xbm (A) \xcc U,$

(5.3) $ \xbm (Y) \xcc H(U)$ $ \xch $ $Y \xcc H(U)$ and $ \xbm (U \xcv Y)=
\xbm (U),$

(5.4) $x \xbe \xbm (U),$ $x \xbe Y- \xbm (Y)$ $ \xch $ $Y \xcC H(U),$

(5.5) $Y \xcC H(U)$ $ \xch $ $ \xbm (U \xcv Y) \xcC H(U).$
\index{Fact HU Proof}

\efa

\bfa

$\hspace{0.01em}$

% (+++ Orig. No.:  Fact HU Proof +++)

\label{Fact HU Proof}

(0.1) and (0.2) trivial by definition.

(0.3) Proof by induction.
$H(U)_{0}=H(U,u)_{0}$ is trivial.
Suppose $H(U)_{ \xbb }=H(U,u)_{ \xbb }$ has been shown for $ \xbb < \xba
.$
The limit step is trivial, so suppose $ \xba = \xbb +1.$
Let $X$ be such that $ \xbm (X) \xcc H(U)_{ \xbb }=H(U,u)_{ \xbb },$ so $X
\xcc H(U)_{ \xba }.$ Consider $X \xcv U,$ so
$u \xbe X \xcv U,$ $ \xbm (X \xcv U)$ is defined and by $( \xbm PR)$ and
$( \xbm \xcc )$
$ \xbm (X \xcv U) \xcc \xbm (X) \xcv \xbm (U) \xcc H(U)_{ \xbb }=H(U,u)_{
\xbb },$ so $X \xcv U \xcc H(U,u)_{ \xba }.$

(0.4) Immediate by (0.3).

(1.1) By $ \xCf (HU),$ if $ \xbm (Y) \xcc H(U),$ then $ \xbm (U) \xcs Y
\xcc \xbm (Y).$ But, if $Y \xcc U,$ then
$ \xbm (Y) \xcc H(U)$ by $( \xbm \xcc ).$

(1.2) Let $ \xbm (U) \xcc X \xcc U.$ Then by (1.1) $ \xbm (U)= \xbm (U)
\xcs X \xcc \xbm (X).$ By
$ \xbm (U) \xcc X$ and $( \xbm \xcc )$
$ \xbm (U) \xcc U \xcc H(X),$ so by (HU) and $X \xcc U$ and $( \xbm \xcc
),$
$ \xbm (X)= \xbm (X) \xcs U \xcc \xbm (U)$ by $( \xbm \xcc ).$

(3.1) $ \xbm (A) \xcs A_{j} \xcc \xbm (A_{j}) \xcc \xcV \xbm (A_{i}),$ so
by $ \xbm (A) \xcc A= \xcV A_{i}$ $ \xbm (A) \xcc \xcV \xbm (A_{i}).$

(3.2) trivial.

(3.3) $ \xbm (U \xcv Y)-H(U)$ $ \xcc_{(3.2)}$ $ \xbm (U \xcv Y)-U$ $ \xcc
$ (by $( \xbm \xcc ))$
$ \xbm (U \xcv Y) \xcs Y$ $ \xcc_{( \xbm PR)}$ $ \xbm (Y).$

(4.1) We show that, if $X \xcc H_{2}(U),$ then $X \xcc H_{1}(U),$ more
precisely, if $ \xbm (X) \xcc H_{1}(U),$
then already $X \xcc H_{1}(U),$ so the construction stops already at
$H_{1}(U).$
Suppose then $ \xbm (X) \xcc \xcV \{Y: \xbm (Y) \xcc U\},$ and let $A:=X
\xcv U.$ We show that $ \xbm (A) \xcc U,$ so
$X \xcc A \xcc H_{1}(U).$ Let $a \xbe \xbm (A).$ By $( \xbm PR),$ $( \xbm
\xcc ),$
$ \xbm (A) \xcc \xbm (X) \xcv \xbm (U).$ If $a \xbe \xbm (U) \xcc U,$ we
are done. If $a \xbe \xbm (X),$ there is $Y$ s.t. $ \xbm (Y) \xcc U$ and
$a \xbe Y,$ so $a \xbe \xbm (A) \xcs Y.$
By Fact \ref{Fact Cum-Alpha} (page \pageref{Fact Cum-Alpha}), (5.1.3), we have
for $Y$ s.t.
$ \xbm (Y) \xcc U$ and $U \xcc A$ $ \xbm (A) \xcs Y \xcc \xbm (U).$ Thus
$a \xbe \xbm (U),$ and we are done again.

(4.2) Let $U \xcc A,$ $ \xbm (A) \xcc H(U)=H_{1}(U)$ by (4.1). So $ \xbm
(A)$ $=$ $ \xcV \{ \xbm (A) \xcs Y: \xbm (Y) \xcc U\}$ $ \xcc $
$ \xbm (U)$ $ \xcc $ $U,$ again by Fact \ref{Fact Cum-Alpha} (page \pageref{Fact
Cum-Alpha}),
(5.1.3).

(4.3) Let $ \xbm (Y) \xcc H(U),$ then by $ \xbm (U) \xcc H(U)$ and $( \xbm
PR),$ $( \xbm \xcc ),$
$ \xbm (U \xcv Y) \xcc \xbm (U) \xcv \xbm (Y) \xcc H(U),$ so by (4.2) $
\xbm (U \xcv Y) \xcc U$ and $U \xcv Y \xcc H(U).$
Moreover, $ \xbm (U \xcv Y) \xcc U \xcc U \xcv Y$ $ \xch_{( \xbm CUM)}$ $
\xbm (U \xcv Y)= \xbm (U).$

(4.4) If not, $Y \xcc H(U),$ so $ \xbm (Y) \xcc H(U),$ so $ \xbm (U \xcv
Y)= \xbm (U)$ by (4.3),
but $x \xbe Y- \xbm (Y)$ $ \xch_{( \xbm PR)}$ $x \xce \xbm (U \xcv Y)=
\xbm (U),$ $contradiction.$

(4.5) $ \xbm (U \xcv Y) \xcc H(U)$ $ \xch_{(4.3)}$ $U \xcv Y \xcc H(U).$

(5) Trivial by (1) and (4).

$ \xcz $
\\[3ex]
\paragraph{The representation result}
\label{Section 2.2.4.3}

\efa

We turn now to the representation result and its proof.

We adapt Proposition 3.7.15 in  \cite{Sch04} and its proof. All we
need is $ \xCf (HU,u)$ and
$( \xbm \xcc ).$ We modify the proof of Remark 3.7.13 (1) in  \cite{Sch04}
(now Remark \ref{Remark D-4.4.4} (page \pageref{Remark D-4.4.4}) )
so we will not need $( \xcs )$ any more.
We will give the full proof, although its essential elements have already
been
published, for three reasons: First, the new version will need less
prerequisites than the old proof does (closure under finite intersections
is
not needed any more, and replaced by $ \xCf (HU,u)).$ Second, we will more
clearly
separate the requirements to do the construction from the construction
itself,
thus splitting the proof neatly into two parts.

We show how to work with $( \xbm \xcc )$ and $ \xCf (HU,u)$ only. Thus,
once we have shown $( \xbm \xcc )$
and $ \xCf (HU,u),$ we have finished the substantial side, and enter the
administrative
part, which will not use any prerequisites about domain closure any more.
At the same time, this gives a uniform proof of the difficult part for the
case
with and without $( \xcv ),$ in the former case we can even work with the
stronger
$H(U).$ The easy direction of the former parts needs a proof of the
stronger
$H(U),$ but this is easy.

Note that, in the presence of $( \xbm \xcc ),$
$(HU,u) \xch ( \xbm Cum \xca )$ and $( \xbm Cum0) \xch ( \xbm PR),$ by
Fact \ref{Fact Cum-Alpha-HU} (page \pageref{Fact Cum-Alpha-HU}), (3) and
Fact \ref{Fact Cum-Alpha} (page \pageref{Fact Cum-Alpha}), (7), so
$ \xCf (HU,u)$ entails $( \xbm PR),$
so we can use it in our context, where $ \xCf (HU,u)$ will be the central
property.

\bfa

$\hspace{0.01em}$

% (+++ Orig. No.:  Fact D-4.4.3 +++)

\label{Fact D-4.4.3}

$ \xCf (HU,u)$ holds in all smooth models.

\efa

\subparagraph{
Proof
}

$\hspace{0.01em}$

% (+++ Orig.:  Proof +++)

(1) Trivial by definition.

(2) Suppose not. So let $x \xbe \xbm (U),$ $x \xbe Y- \xbm (Y),$ $ \xbm
(Y) \xcc H(U,x).$
By smoothness, there is $x_{1} \xbe \xbm (Y),$ $x \xee x_{1},$
and let $ \xbk_{1}$ be the least $ \xbk $ s.t. $x_{1} \xbe H(U,x)_{
\xbk_{1}}.$ $ \xbk_{1}$ is not a
limit, and $x_{1} \xbe U'_{x_{1}}- \xbm (U'_{x_{1}})$ with $x \xbe
U'_{x_{1}}$ by definition of $H(U,x)$
for some $U'_{x_{1}},$ so as $x_{1} \xce \xbm (U'_{x_{1}}),$
there must be (by smoothness) some other
$x_{2} \xbe \xbm (U'_{x_{1}}) \xcc H(U,x)_{ \xbk_{1}-1}$ with $x \xee
x_{2}.$ Continue with $x_{2},$ we thus construct
a descending chain of ordinals, which cannot be infinite, so there must be
$x_{n} \xbe \xbm (U'_{x_{n}}) \xcc U,$ $x \xee x_{n},$ contradicting
minimality of $x$ in $U.$
(More precisely, this works for all copies of $x.)$
$ \xcz $
\\[3ex]

We first show two basic facts and then turn to the main result,
Proposition \ref{Proposition D-4.4.6} (page \pageref{Proposition D-4.4.6}).

\bd

$\hspace{0.01em}$

% (+++ Orig. No.:  Definition D-4.4.4 +++)

\label{Definition D-4.4.4}

For $x \xbe Z,$ let $ \xdw_{x}:=\{ \xbm (Y)$: $Y \xbe \xdy $ $ \xcu $ $x
\xbe Y- \xbm (Y)\},$ $ \xbG_{x}:= \xbP \xdw_{x}$, and $K:=\{x \xbe Z$: $
\xcE X \xbe \xdy.x \xbe \xbm (X)\}.$

\ed

\br

$\hspace{0.01em}$

% (+++ Orig. No.:  Remark D-4.4.4 +++)

\label{Remark D-4.4.4}

(1) $x \xbe K$ $ \xch $ $ \xbG_{x} \xEd \xCQ,$

(2) $g \xbe \xbG_{x}$ $ \xch $ $ran(g) \xcc K.$

\er

\subparagraph{
Proof
}

$\hspace{0.01em}$

% (+++ Orig.:  Proof +++)

(1) We give two proofs, the first uses $( \xbm Cum0),$ the second the
stronger
(by Fact \ref{Fact Cum-Alpha-HU} (page \pageref{Fact Cum-Alpha-HU})  (3))
$ \xCf (HU,u).$

(a) We have to show that $Y \xbe \xdy,$ $x \xbe Y- \xbm (Y)$ $ \xch $ $
\xbm (Y) \xEd \xCQ.$ Suppose then $x \xbe \xbm (X),$
this exists, as $x \xbe K,$ and $ \xbm (Y)= \xCQ,$ so $ \xbm (Y) \xcc X,$
$x \xbe Y,$ so by $( \xbm Cum0)$ $x \xbe \xbm (Y).$

(b) Consider $H(X,x),$ suppose $ \xbm (Y)= \xCQ,$ $x \xbe Y,$ so $Y \xcc
H(X,x),$ so
$x \xbe \xbm (Y)$ by $(HU,u).$

(2) By definition, $ \xbm (Y) \xcc K$ for all $Y \xbe \xdy.$
$ \xcz $
\\[3ex]

\bc

$\hspace{0.01em}$

% (+++ Orig. No.:  Claim D-4.4.5 +++)

\label{Claim D-4.4.5}

Let $U \xbe \xdy,$ $x \xbe K.$ Then

(1) $x \xbe \xbm (U)$ $ \xcj $ $x \xbe U$ $ \xcu $ $ \xcE f \xbe
\xbG_{x}.ran(f) \xcs U= \xCQ,$

(2) $x \xbe \xbm (U)$ $ \xcj $ $x \xbe U$ $ \xcu $ $ \xcE f \xbe
\xbG_{x}.ran(f) \xcs H(U,x)= \xCQ.$

\ec

\subparagraph{
Proof
}

$\hspace{0.01em}$

% (+++ Orig.:  Proof +++)

(1)

Case 1: $ \xdw_{x}= \xCQ,$ thus $ \xbG_{x}=\{ \xCQ \}.$

`` $ \xch $ '': Take $f:= \xCQ.$

`` $ \xci $ '': $x \xbe U \xbe \xdy,$ $ \xdw_{x}= \xCQ $ $ \xch $ $x \xbe
\xbm (U)$ by definition of $ \xdw_{x}.$

Case 2: $ \xdw_{x} \xEd \xCQ.$

`` $ \xch $ '': Let $x \xbe \xbm (U) \xcc U.$ Consider $H(U,x).$ If $ \xbm
(Y) \xbe \xdw_{x},$ then $x \xbe Y- \xbm (Y),$ so
by $(HU,u)$ for $H(U,x)$ $ \xbm (Y)-H(U,x) \xEd \xCQ,$ but $ \xbm (U)
\xcc U \xcc H(U,x).$

`` $ \xci $ '': If $x \xbe U- \xbm (U),$ so $ \xbm (U) \xbe \xdw_{x}$,
moreover $ \xbG_{x} \xEd \xCQ $ by
Remark \ref{Remark D-4.4.4} (page \pageref{Remark D-4.4.4}), (1)
and thus $ \xbm (U) \xEd \xCQ,$ so $ \xcA f \xbe \xbG_{x}.ran(f) \xcs U
\xEd \xCQ.$

(2): The Case 1 is as for (1).

Case 2: `` $ \xch $ '' was shown already in Case 1.

`` $ \xci $ '': Let $x \xbe U- \xbm (U),$ then by $x \xbe K$ $ \xbm (U)
\xEd \xCQ $ (see proof of
Remark \ref{Remark D-4.4.4} (page \pageref{Remark D-4.4.4}) ), moreover $ \xbm
(U) \xcc U \xcc H(U,x),$
so
$ \xcA f \xbe \xbG_{x}.ran(f) \xcs H(U,x) \xEd \xCQ.$

$ \xcz $ (Claim \ref{Claim D-4.4.5} (page \pageref{Claim D-4.4.5}) )
\\[3ex]

The following Proposition \ref{Proposition D-4.4.6} (page \pageref{Proposition
D-4.4.6})
is the main positive result of
Section \ref{Section 2.2.4} (page \pageref{Section 2.2.4})  and
shows how to characterize smooth structures in the absence of closure
under
finite unions. The strategy of the proof follows closely the proof of
Proposition 3.3.4 in  \cite{Sch04}.

\bp

$\hspace{0.01em}$

% (+++ Orig. No.:  Proposition D-4.4.6 +++)

\label{Proposition D-4.4.6}

\label{Proposition Smooth-Complete}

Let $ \xbm: \xdy \xcp \xdp (Z).$
Then there is a $ \xdy -$smooth preferential structure $ \xdz,$ s.t. for
all $X \xbe \xdy $
$ \xbm (X)= \xbm_{ \xdz }(X)$ iff $ \xbm $ satisfies $( \xbm \xcc )$ and $
\xCf (HU,u)$ above.

In particular, we need no prerequisites about domain closure.

\ep

\subparagraph{
Proof
}

$\hspace{0.01em}$

% (+++ Orig.:  Proof +++)

`` $ \xch $ '' $ \xCf (HU,u)$ was shown in
Fact \ref{Fact D-4.4.3} (page \pageref{Fact D-4.4.3}).

Outline of `` $ \xci $ '': We first define a structure $ \xdz $ which
represents $ \xbm,$ but is
not necessarily $ \xdy -$smooth, refine it to $ \xdz ' $ and show that $
\xdz ' $ represents $ \xbm $
too, and that $ \xdz ' $ is $ \xdy -$smooth.

In the structure $ \xdz ',$ all pairs destroying smoothness in $ \xdz $
are successively
repaired, by adding minimal elements: If $ \xBc y,j \xBe $ is not minimal, and
has
no minimal
$ \xBc x,i \xBe $ below it, we just add one such $ \xBc x,i \xBe.$ As the
repair process
might itself
generate such ``bad'' pairs, the process may have to be repeated infinitely
often.
Of course, one has to take care that the representation property is
preserved.

\bcs

$\hspace{0.01em}$

% (+++ Orig. No.:  Construction D-4.4.1 +++)

\label{Construction D-4.4.1}

(Construction of $ \xdz )$

Let $ \xdx $ $:=$ $\{ \xBc x,g \xBe $: $x \xbe K,$ $g \xbe \xbG_{x}\},$ $ \xBc
x',g'  \xBe
\xeb  \xBc x,g \xBe $ $: \xcj $ $x' \xbe ran(g),$ $ \xdz:= \xBc  \xdx, \xeb 
\xBe.$

\ecs

\bc

$\hspace{0.01em}$

% (+++ Orig. No.:  Claim D-4.4.7 +++)

\label{Claim D-4.4.7}

$ \xcA U \xbe \xdy. \xbm (U)= \xbm_{ \xdz }(U)$

\ec

\subparagraph{
Proof
}

$\hspace{0.01em}$

% (+++ Orig.:  Proof +++)

Case 1: $x \xce K.$ Then $x \xce \xbm (U)$ and $x \xce \xbm_{ \xdz }(U).$

Case 2: $x \xbe K.$

By Claim \ref{Claim D-4.4.5} (page \pageref{Claim D-4.4.5}), (1)
it suffices to show that for all $U \xbe \xdy $
$x \xbe \xbm_{ \xdz }(U)$ $ \xcj $ $x \xbe U$ $ \xcu $ $ \xcE f \xbe
\xbG_{x}.ran(f) \xcs U= \xCQ.$
Fix $U \xbe \xdy.$

`` $ \xch $ '': $x \xbe \xbm_{ \xdz }(U)$ $ \xch $ ex. $ \xBc x,f \xBe $ minimal
in
$ \xdx \xex U,$ thus $x \xbe U$ and there is no
$ \xBc x',f'  \xBe  \xeb  \xBc x,f \xBe,$ $x' \xbe U,$ $x' \xbe K.$ But if $x'
\xbe K,$ then
by
Remark \ref{Remark D-4.4.4} (page \pageref{Remark D-4.4.4}),
(1), $ \xbG_{x' } \xEd \xCQ,$
so we find suitable $f'.$ Thus, $ \xcA x' \xbe ran(f).x' \xce U$ or $x'
\xce K.$ But $ran(f) \xcc K,$ so
$ran(f) \xcs U= \xCQ.$

`` $ \xci $ '': If $x \xbe U,$ $f \xbe \xbG_{x}$ s.t. $ran(f) \xcs U= \xCQ
,$ then $ \xBc x,f \xBe $ is minimal in $ \xdx \xex U.$
$ \xcz $ (Claim \ref{Claim D-4.4.7} (page \pageref{Claim D-4.4.7}) )
\\[3ex]

We will use in the construction of the refined structure $ \xdz ' $ the
following
definition:

\bd

$\hspace{0.01em}$

% (+++ Orig. No.:  Definition D-4.4.5 +++)

\label{Definition D-4.4.5}

$ \xbs $ is called $x-$admissible sequence iff

1. $ \xbs $ is a sequence of length $ \xck \xbo,$ $ \xbs =\{ \xbs_{i}:i
\xbe \xbo \},$

2. $ \xbs_{o} \xbe \xbP \{ \xbm (Y)$: $Y \xbe \xdy $ $ \xcu $ $x \xbe Y-
\xbm (Y)\},$

3. $ \xbs_{i+1} \xbe \xbP \{ \xbm (X)$: $X \xbe \xdy $ $ \xcu $ $x \xbe
\xbm (X)$ $ \xcu $ $ran( \xbs_{i}) \xcs X \xEd \xCQ \}.$

By 2., $ \xbs_{0}$ minimizes $x,$ and by 3., if $x \xbe \xbm (X),$ and
$ran( \xbs_{i}) \xcs X \xEd \xCQ,$ i.e. we
have destroyed minimality of $x$ in $X,$ $x$ will be above some $y$
minimal in $X$ to
preserve smoothness.

Let $ \xbS_{x}$ be the set of $x-$admissible sequences, for $ \xbs \xbe
\xbS_{x}$ let $ \wt{ \xbs }:= \xcV \{ran( \xbs_{i}):i \xbe \xbo \}.$

\ed

\bcs

$\hspace{0.01em}$

% (+++ Orig. No.:  Construction D-4.4.2 +++)

\label{Construction D-4.4.2}

(Construction of $ \xdz ' )$

Note that by Remark \ref{Remark D-4.4.4} (page \pageref{Remark D-4.4.4}),
(1), $ \xbS_{x} \xEd \xCQ,$ if $x \xbe K$
(this does $ \xbs_{0},$ $ \xbs_{i+1}$ is trivial as by prerequisite $ \xbm
(X) \xEd \xCQ ).$

Let $ \xdx ' $ $:=$ $\{ \xBc x, \xbs  \xBe $: $x \xbe K$ $ \xcu $ $ \xbs \xbe
\xbS_{x}\}$ and $ \xBc x', \xbs '  \xBe  \xeb '  \xBc x, \xbs  \xBe $ $: \xcj $
$x' \xbe \wt{
\xbs }$.
Finally, let $ \xdz ':= \xBc  \xdx ', \xeb '  \xBe,$ and $ \xbm ':= \xbm_{ \xdz
' }.$

\ecs

It is now easy to show that $ \xdz ' $ represents $ \xbm,$ and that $
\xdz ' $ is smooth.
For $x \xbe \xbm (U),$ we construct a special $x-$admissible sequence $
\xbs^{x,U}$ using the
properties of $H(U,x)$ as described at the beginning of this section.

\bc

$\hspace{0.01em}$

% (+++ Orig. No.:  Claim D-4.4.8 +++)

\label{Claim D-4.4.8}

For all $U \xbe \xdy $ $ \xbm (U)= \xbm_{ \xdz }(U)= \xbm ' (U).$

\ec

\subparagraph{
Proof
}

$\hspace{0.01em}$

% (+++ Orig.:  Proof +++)

If $x \xce K,$ then $x \xce \xbm_{ \xdz }(U),$ and $x \xce \xbm ' (U)$ for
any $U.$ So assume $x \xbe K.$ If $x \xbe U$ and
$x \xce \xbm_{ \xdz }(U),$ then for all $ \xBc x,f \xBe  \xbe \xdx,$ there is
$ \xBc x',f'  \xBe  \xbe \xdx $ with $ \xBc x',f'  \xBe  \xeb  \xBc x,f \xBe $
and
$x' \xbe U.$ Let now $ \xBc x, \xbs  \xBe  \xbe \xdx ',$ then $ \xBc x, \xbs_{0}
\xBe  \xbe
\xdx,$ and let
$ \xBc x',f'  \xBe  \xeb  \xBc x, \xbs_{0} \xBe $ in $ \xdz $ with
$x' \xbe U.$ As $x' \xbe K,$ $ \xbS_{x' } \xEd \xCQ,$ let $ \xbs ' \xbe
\xbS_{x' }.$
Then $ \xBc x', \xbs '  \xBe  \xeb '  \xBc x, \xbs  \xBe $ in $ \xdz '.$ Thus $x
\xce \xbm '
(U).$
Thus, for all $U \xbe \xdy,$ $ \xbm ' (U) \xcc \xbm_{ \xdz }(U)= \xbm
(U).$

It remains to show $x \xbe \xbm (U) \xch x \xbe \xbm ' (U).$

Assume $x \xbe \xbm (U)$ (so $x \xbe K),$ $U \xbe \xdy,$ we will
construct minimal $ \xbs,$ i.e. show that
there is $ \xbs^{x,U} \xbe \xbS_{x}$ s.t. $ \wt{ \xbs^{x,U}} \xcs U= \xCQ
.$ We construct this $ \xbs^{x,U}$ inductively, with the
stronger property that $ran( \xbs^{x,U}_{i}) \xcs H(U,x)= \xCQ $ for all
$i \xbe \xbo.$

$ \ul{ \xbs^{x,U}_{0}:}$

$x \xbe \xbm (U),$ $x \xbe Y- \xbm (Y)$ $ \xch $ $ \xbm (Y)-H(U,x) \xEd
\xCQ $ by $ \xCf (HU,u)$ for $H(U,x).$
Let $ \xbs^{x,U}_{0}$ $ \xbe $ $ \xbP \{ \xbm (Y)-H(U,x):$ $Y \xbe \xdy,$
$x \xbe Y- \xbm (Y)\},$ so $ran( \xbs^{x,U}_{0}) \xcs H(U,x)= \xCQ.$

$ \ul{ \xbs^{x,U}_{i} \xch \xbs^{x,U}_{i+1}:}$

By the induction hypothesis, $ran( \xbs^{x,U}_{i}) \xcs H(U,x)= \xCQ.$
Let $X \xbe \xdy $ be s.t. $x \xbe \xbm (X),$
$ran( \xbs^{x,U}_{i}) \xcs X \xEd \xCQ.$ Thus $X \xcC H(U,x),$ so $ \xbm
(X)-H(U,x) \xEd \xCQ $ by
Fact \ref{Fact Cum-Alpha-HU} (page \pageref{Fact Cum-Alpha-HU}), (1).
Let $ \xbs^{x,U}_{i+1}$ $ \xbe $ $ \xbP \{ \xbm (X)-H(U,x):$ $X \xbe \xdy
,$ $x \xbe \xbm (X),$ $ran( \xbs^{x,U}_{i}) \xcs X \xEd \xCQ \},$ so $ran(
\xbs^{x,U}_{i+1}) \xcs H(U,x)= \xCQ.$
The construction satisfies the $x-$admissibility condition.
$ \xcz $
\\[3ex]

It remains to show:

\bc

$\hspace{0.01em}$

% (+++ Orig. No.:  Claim D-4.4.9 +++)

\label{Claim D-4.4.9}

$ \xdz ' $ is $ \xdy -$smooth.

\ec

\subparagraph{
Proof
}

$\hspace{0.01em}$

% (+++ Orig.:  Proof +++)

Let $X \xbe \xdy,$ $ \xBc x, \xbs  \xBe  \xbe \xdx ' \xex X.$

Case 1, $x \xbe X- \xbm (X):$ Then $ran( \xbs_{0}) \xcs \xbm (X) \xEd \xCQ
,$ let $x' \xbe ran( \xbs_{0}) \xcs \xbm (X).$ Moreover,
$ \xbm (X) \xcc K.$ Then for all $ \xBc x', \xbs '  \xBe  \xbe \xdx ' $
$ \xBc x', \xbs '  \xBe  \xeb  \xBc x, \xbs  \xBe.$ But $ \xBc x', \xbs^{x',X}
\xBe $ as
constructed in the proof of
Claim \ref{Claim D-4.4.8} (page \pageref{Claim D-4.4.8})
is minimal in $ \xdx ' \xex X.$

Case 2, $x \xbe \xbm (X)= \xbm_{ \xdz }(X)= \xbm ' (X):$ If $ \xBc x, \xbs  \xBe
$
is minimal in $ \xdx ' \xex X,$ we are done.
So suppose there is $ \xBc x', \xbs '  \xBe  \xeb  \xBc x, \xbs  \xBe,$ $x'
\xbe X.$ Thus
$x' \xbe \wt{ \xbs }.$ Let
$x' \xbe ran( \xbs_{i}).$ So $x \xbe \xbm (X)$ and $ran( \xbs_{i}) \xcs X
\xEd \xCQ.$ But
$ \xbs_{i+1} \xbe \xbP \{ \xbm (X' )$: $X' \xbe \xdy $ $ \xcu $ $x \xbe
\xbm (X' )$ $ \xcu $ $ran( \xbs_{i}) \xcs X' \xEd \xCQ \},$ so $X$ is one
of the $X',$
moreover $ \xbm (X) \xcc K,$ so there is $x'' \xbe \xbm (X) \xcs ran(
\xbs_{i+1}) \xcs K,$ so for all
$ \xBc x'', \xbs ''  \xBe  \xbe \xdx ' $
$ \xBc x'', \xbs ''  \xBe  \xeb  \xBc x, \xbs  \xBe.$ But again $ \xBc x'',
\xbs^{x'',X} \xBe $ as
constructed in the proof
of Claim \ref{Claim D-4.4.8} (page \pageref{Claim D-4.4.8})
is minimal in $ \xdx ' \xex X.$

$ \xcz $ (Claim \ref{Claim D-4.4.9} (page \pageref{Claim D-4.4.9})  and
Proposition \ref{Proposition D-4.4.6} (page \pageref{Proposition D-4.4.6}) )
\\[3ex]

We conclude this section by showing that we cannot improve substantially.

\bp

$\hspace{0.01em}$

% (+++ Orig. No.:  Proposition D-4.4.10 +++)

\label{Proposition D-4.4.10}

There is no fixed size characterization of $ \xbm -$functions which are
representable
by smooth structures, if the domain is not closed under finite unions.

\ep

\subparagraph{
Proof
}

$\hspace{0.01em}$

% (+++ Orig.:  Proof +++)

Suppose we have a fixed size characterization, which allows to distinguish
$ \xbm -$functions on domains which are not necessarily closed under
finite unions,
and which can be represented by smooth structures, from those which cannot
be
represented in this way. Let the characterization have $ \xba $ parameters
for sets,
and consider
Example \ref{Example Inf-Cum-Alpha} (page \pageref{Example Inf-Cum-Alpha})
with $ \xbk = \xbb +1,$ $ \xbb > \xba $
(as a cardinal). This structure
cannot be represented, as $( \xbm Cum \xbk )$ fails - see
Fact \ref{Fact Cum-Alpha} (page \pageref{Fact Cum-Alpha}), (2.1).
As we have only $ \xba $ parameters, at
least one of the $X_{ \xbg }$ is not mentioned, say $X_{ \xbd }.$ Without
loss of generality,
we may assume that
$ \xbd = \xbd ' +1.$ We change now the structure, and erase one pair of
the relation,
$x_{ \xbd } \xeb x_{ \xbd +1}.$ Thus, $ \xbm (X_{ \xbd })=\{c,x_{ \xbd
},x_{ \xbd +1}\}.$ But now we cannot go any more from
$X_{ \xbd ' }$ to $X_{ \xbd ' +1}=X_{ \xbd },$ as $ \xbm (X_{ \xbd }) \xcC
X_{ \xbd ' }.$ Consequently, the only chain showing
that $( \xbm Cum \xca )$ fails is interrupted - and we have added no new
possibilities,
as inspection of cases shows. $(x_{ \xbd +1}$ is now globally minimal, and
increasing
$ \xbm (X)$ cannot introduce new chains, only interrupt chains.) Thus, $(
\xbm Cum \xca )$ holds
in the modified example, and it is thus representable by a smooth
structure, as
above proposition shows. As we did not touch any of the parameters, the
truth
value of the characterization is unchanged, which was negative. So the
``characterization'' cannot be correct.
$ \xcz $
\\[3ex]
\subsubsection{The transitive smooth case}
\label{Section 2.2.5}
\label{Section CS-Smooth-Trans}

Unfortunately, $( \xbm Cumt \xca )$ is a necessary but not sufficient
condition for
smooth transitive structures, as can be seen in the following example.

\be

$\hspace{0.01em}$

% (+++ Orig. No.:  Example D-4.5.1 +++)

\label{Example D-4.5.1}

We assume no closure whatever.

$U:=\{u_{1},u_{2},u_{3},u_{4}\},$ $ \xbm (U):=\{u_{3},u_{4}\}$

$Y_{1}:=\{u_{4},v_{1},v_{2},v_{3},v_{4}\},$ $ \xbm
(Y_{1}):=\{v_{3},v_{4}\}$

$Y_{2,1}:=\{u_{2},v_{2},v_{4}\},$ $ \xbm (Y_{2,1}):=\{u_{2},v_{2}\}$

$Y_{2,2}:=\{u_{1},v_{1},v_{3}\},$ $ \xbm (Y_{2,2}):=\{u_{1},v_{1}\}$

For no $ \xCf A,B$ $ \xbm (A) \xcc B$ $(A \xEd B),$ so the prerequisite of
$( \xbm Cumt \xba )$ is false,
and $( \xbm Cumt \xba )$ holds, but there is no smooth transitive
representation possible:
Consider $Y_{1}.$
If $u_{4} \xee v_{3},$ then $Y_{2,2}$ makes this impossible, if $u_{4}
\xee v_{4},$ then $Y_{2,1}$ makes this
impossible.

$ \xcz $
\\[3ex]

\ee

\br

$\hspace{0.01em}$

% (+++ Orig. No.:  Remark D-4.5.1 +++)

\label{Remark D-4.5.1}

(1)
The situation does not change when we have copies, the same argument will
still
work: There is a $U-$minimal copy $ \xBc u_{4},i \xBe,$ by smoothness and
$Y_{1},$
there must be a
$Y_{1}-$minimal copy, e.g. $ \xBc v_{3},j \xBe  \xeb  \xBc u_{4},i \xBe.$ By
smoothness and
$Y_{2,2},$ there must be a
$Y_{2,2}-$minimal $ \xBc u_{1},k \xBe $ or $ \xBc v_{1},l \xBe $ below $ \xBc
v_{3},j \xBe.$ But
$v_{1}$ is in $Y_{1},$ contradicting
minimality of $ \xBc v_{3},j \xBe,$ $u_{1}$ is in $U,$ contadicting minimality
of
$ \xBc u_{4},i \xBe $ by
transitivity. If we choose $ \xBc v_{4},j \xBe $ minimal below
$ \xBc u_{4},i \xBe,$ we will work with
$Y_{2,1}$ instead of $Y_{2,2}.$

(2)
We can also close under arbitrary intersections, and the example will
still
work: We have to consider $U \xcs Y_{1},$ $U \xcs Y_{2,1},$ $U \xcs
Y_{2,2},$ $Y_{2,1} \xcs Y_{2,2},$ $Y_{1} \xcs Y_{2,1},$ $Y_{1} \xcs
Y_{2,2},$
there are no further intersections to consider. We may assume $ \xbm
(A)=A$ for all
these intersections (working with copies). But then $ \xbm (A) \xcc B$
implies $ \xbm (A)=A$
for all sets, and all $( \xbm Cumt \xba )$ hold again trivially.

(3) If we had finite unions, we could form $A:=U \xcv Y_{1} \xcv Y_{2,1}
\xcv Y_{2,2},$ then $ \xbm (A)$ would
have to be a subset of $\{u_{3}\}$ by $( \xbm PR),$ so by $( \xbm CUM)$
$u_{4} \xce \xbm (U),$ a contradiction.
Finite unions allow us to ``look ahead'', without $( \xcv ),$ we see
desaster only at
the end - and have to backtrack, i.e. try in our example $Y_{2,1},$ once
we have seen
impossibility via $Y_{2,2},$ and discover impossibility again at the end.
$ \xcz $
\\[3ex]
\subsubsection{General ranked structures}
\label{Section 2.2.6}
%  Section (4.6):  General ranked structures
%  Section (4.6):  General ranked structures
% %
% ===========================================
\paragraph{Introduction}
\label{Section 2.2.6.1}
%  Subsection (4.6.1):  Introduction
%  Subsection (4.6.1):  Introduction
% %
% ===================================

\er

Fix $f: \xdy \xcp \xdp (Z).$
\paragraph{The general case}
\label{Section 2.2.6.2}
%  Subsection (4.6.2):  The general case
%  Subsection (4.6.2):  The general case
% %
% =======================================

We summarize in the following
Lemma \ref{Lemma D-4.6.6} (page \pageref{Lemma D-4.6.6})
our results for the general
ranked case, many of them trivial.

\bl

$\hspace{0.01em}$

% (+++ Orig. No.:  Lemma D-4.6.6 +++)

\label{Lemma D-4.6.6}

We assume here for simplicity that all elements occur in the structure.

(1) If $ \xbm (X)= \xCQ,$ then each element $x \xbe X$ either has
infinitely many copies, or
below each copy of each $x,$ there is an infinite descending chain of
other
elements.

(2) If there is no $X$ s.t. $x \xbe \xbm (X),$ then we can make infinitely
many copies
of $x.$

(3) There is no simple way to detect whether there is for all $x$ some $X$
s.t.
$x \xbe \xbm (X).$ More precisely: there is no normal finite
characterization of ranked
structures, in which each $x$ in the domain occurs in at least one $ \xbm
(X).$

Suppose in the sequel that for each $x$ there is some $X$ s.t. $x \xbe
\xbm (X).$ (This is
the hard case.)

(4) If the language is finite, then $X \xEd \xCQ $ implies $ \xbm (X) \xEd
\xCQ.$

Suppose now the language to be infinite.

(5) If we admit all theories, then $ \xbm (M(T))=M(T)$ for all complete
theories.

(6) It is possible to have $ \xbm (M( \xbf ))= \xCQ $ for all formulas $
\xbf,$ even though all
models occur in exactly one copy.

(7) If the domain is sufficiently rich, then we cannot have $ \xbm (X)=
\xCQ $ for ``many''
$X.$

(8) We see that a small domain (see Case (6)) can have many $X$ with $
\xbm (X)= \xCQ,$
but if the domain is too dense (see Case (7)), then we cannot have many
$ \xbm (X)= \xCQ.$
(We do not know any criterion to distinguish poor from rich domains.)

(9) If we have all pairs in the domain, we can easily construct the
ranking.

\el

\subparagraph{
Proof
}

$\hspace{0.01em}$

% (+++ Orig.:  Proof +++)

(1), (2), (4), (5), (9) are trivial, there is nothing to show for (8).

(3)
Suppose there is a normal characterization $ \xbF $ of such structures,
where each
element $x$ occurs at least once in a set $X$ s.t. $x \xbe \xbm (X).$ Such
a
characterization will be a finite boolean combination of set expressions $
\xbF,$
universally quantified, in the spirit of $ \xCf (AND),$ $ \xCf (RM)$ etc.

We consider a realistic counterexample - an infinite propositional
language
and the sets definable by formulas. We do not necessarily assume
definability
preservation, and work with full equality of results.

Take an infinite propositional language $p_{i}:i< \xbo.$ Choose an
arbitrary model $m,$
say $m \xcm p_{i}:i< \xbo.$

Now, determine the height of any model
$m' $ as follows: $ht(m' ):=the$ first $p_{i}$ s.t. $m(p_{i}) \xEd m'
(p_{i}),$ in our example then
the first $p_{i}$ s.t. $m' \xcm \xCN p_{i}.$ Thus, only $m$ has infinite
height, essentially,
the more different $m' $ is from $m$ (in an alphabetical order), the lower
it is.

Make now $ \xbo $ many copies of $m,$ in infinite descending order, which
you put on
top of the rest.

$ \xbF $ has to fail for some instantiation, as $ \xdx $ does not have the
desired
property. Write this instantiation of $ \xbF $ wlog. as a disjunction of
conjunctions: $ \xcO ( \xcU \xbf_{i,j}).$

Each (consistent, or non-empty) component $ \xbf_{i,j}$ has finite height,
more precisely: the minimum of all heigts of its models (which is a finite
height).
Thus, $ \xcn ( \xbf_{i,j})$ will be just the minimally high models of $
\xbf_{i,j}$ in this order.

Modify now $ \xdx $ s.t. $m$ has only 1 copy, and is just $(+1$ suffices)
above the
minimum
of all the finitely many $ \xbf_{i,j}.$ Then none of the $ \xcn (
\xbf_{i,j})$ is affected, and $m$
has now finite height, say $h,$ and is a minimum in any $M( \xbf ' )$
where $ \xbf ' =$ the
conjunction of the first $h$ values of $m.$

(Remark:
Obviously, there are two easy generalizations for this ranking: First, we
can
go beyond $ \xbo $ (but also stay below $ \xbo ),$ second, instead of
taking just one $m$
as a scale, and which has maximal height, we can take a set $M$ of models:
$ht(m' )$
is then the first $p_{i}$ where $m' (p_{i})$ is different from $ \xCf all$
$m \xbe M.$ Note that in
this
case, in general, not all levels need to be filled. If e.g., $m_{0},$
$m_{1} \xbe M,$ and
$m_{0}(p_{0})=false,$ $m_{1}(p_{0})=true,$ then level 0 will be empty.)

(6)
Let the $p_{i}$ again define an infinite language. Denote by $p^{+}_{i}$
the set of all $+p_{j}$,
where $j>i.$ Let $T$ be the usual tree of models (each model is a branch)
for the
$p_{i},$ with an artificial root $*.$ Let the first model $(=$ branch ) be
$*^{+},$ i.e. the
leftest branch in the obvious way of drawing it. Next, we choose $ \xCN
p^{+}_{0},$ i.e.
we go right, and then all the way left. Next, we consider the 4 sequences
of
$+/-p_{0},$ $+/-p_{1},$ two of them were done already, both ending in
$p^{+}_{1}$, and choose the
remaining two, both ending in $ \xCN p^{+}_{1},$ i.e. the positive
prolongations of
$p_{0}, \xCN p_{1}$ and $ \xCN p_{0}, \xCN p_{1}.$ Thus, at each level, we
take all possible prolongations,
the positive ones were done already, and we count those, which begin
negatively,
and then continue positively. Each formula has in this counting
arbitrarily
big models.

This is not yet a full enumeration of all models, e.g. the branch with all
models negative will never be enumerated. But it suffices for our
purposes.

Reverse the order so far constructed, and put the models not enumerated on
top.
Then all models are considered, and each formula has arbitrarily small
models,
thus $ \xbm ( \xbf )= \xCQ $ for all $ \xbf.$

(7)
Let the domain contain all singletons, and let the structure be without
copies.
The latter can be seen by considering singletons. Suppose now there is a
set $X$
in the domain s.t. $ \xbm (X)= \xCQ.$ Thus, each $x \xbe X$ must have
infinitely many $x' \xbe X$
$x' \xeb x.$ Suppose $ \xdp (X)$ is a subset of the domain. Then there
must be infinite
$Y \xbe \xdp (X)$ s.t. $ \xbm (Y) \xEd \xCQ:$ Suppose not. Let $ \xeb $
be the ranking order. Choose
arbitrary $x \xbe X.$ Consider $X':=\{x' \xbe X:x \xeb x' \},$ then $x
\xbe \xbm (X' ),$ and not all such
$X' $ can be finite - assuming $X$ is big enough, e.g. uncountable.

$ \xcz $
\\[3ex]

We conclude by giving an example of a definability preserving non-compact
preferential logic - in answer to a question by D.Makinson (personal
communication):

\be

$\hspace{0.01em}$

% (+++ Orig. No.:  Example D-4.6.1 +++)

\label{Example D-4.6.1}

Take an infinite language, $p_{i},$ $i< \xbo.$
Fix one model, $m,$ which makes $p_{0}$ true (and, say, for
definiteness, all the others true, too), and $m' $ which is just
like $m,$ but it makes $p_{0}$ false.
Well-order all the other $p_{0}-models,$ and all the other
$ \xCN p_{0}-models$ separately.

Construct now the following ranked structure:

On top, put $m,$ directly below it $m'.$ Further down put the bloc of
the other $ \xCN p_{0}-models,$ and at the bottom the bloc of the other
$p_{0}-models.$

As the structure is well-ordered, it is definability preserving
(singletons
are definable).

Let $T$ be the theory defined by $m,m',$ then $T \xcn \xCN p_{0}.$

Let $ \xbf $ be s.t. $M(T) \xcc M( \xbf ),$ then $M( \xbf )$ contains a
$p_{0}-model$
other than $m,$ so $ \xbf \xcn p_{0}.$

$ \xcz $
\\[3ex]
\subsubsection{Smooth ranked structures}
\label{Section 2.2.7}
%  Section (4.7):  Smooth ranked structures
%  Section (4.7):  Smooth ranked structures
% %
% ==========================================

\ee

We assume that all elements occur in the structure, so smoothness and
$ \xbm (X) \xEd \xCQ $ for $X \xEd \xCQ $ coincide.

The following abstract definition is motivated by:

$ \xCd (u),$ the set of $u' \xbe W$ which have same rank as $u,$

$ \xeb (u),$ the set of $u' \xbe W$ which have lower rank than $u,$

$ \xee (u),$ the set of $u' \xbe W$ which have higher rank than $u,$

all other $u' \xbe W$ will by default have unknown rank in comparison.

We can diagnose e.g. $u' \xbe \xCd (u)$ if $u,u' \xbe \xbm (X)$ for some
$X,$ and $u' \xbe \xee (u)$ if
$u \xbe \xbm (X)$ and $u' \xbe X- \xbm (X)$ for some $X.$

If we sometimes do not know more, we will have to consider also
$ \xec (u)$ and $ \xed (u)$ - this will be needed in
Section \ref{Section 2.3.2} (page \pageref{Section 2.3.2}), where we
will have only incomplete information, due to hidden dimensions.

All other $u' \xbe W$ will by default have unknown rank in comparison.

\bd

$\hspace{0.01em}$

% (+++ Orig. No.:  Definition D-4.7.1 +++)

\label{Definition D-4.7.1}

 \xEh

 \xDH

Define for each $u \xbe W$ three subsets of $W$ $ \xCd (u),$ $ \xeb (u),$
and $ \xee (u).$ Let $ \xdo $ be
the set of all these subsets, i.e. $ \xdo:=\{ \xCd (u),$ $ \xeb (u),$ $
\xee (u):$ $u \xbe W\}$

 \xDH

We say that $ \xdo $ is generated by a choice function $f$

iff

(1) $ \xcA U \xbe \xdy \xcA x,x' \xbe f(U)$ $x' \xbe \xCd (x),$

(2) $ \xcA U \xbe \xdy \xcA x \xbe f(U) \xcA x' \xbe U-f(U)$ $x' \xbe \xee
(x)$

 \xDH

$ \xdo $ is said to be representable by a ranking iff there is a function
$f:W \xcp  \xBc O, \xen  \xBe $ into a total order $ \xBc O, \xen  \xBe $ s.t.

(1) $u' \xbe \xCd (u)$ $ \xch $ $f(u' )=f(u)$

(2) $u' \xbe \xeb (u)$ $ \xch $ $f(u' ) \xen f(u)$

(3) $u' \xbe \xee (u)$ $ \xch $ $f(u' ) \xeq f(u)$

 \xDH

Let $ \xdc ( \xdo )$ be the closure of $ \xdo $ under the following
operations:

 \xEI

 \xDH $u \xbe \xCd (u),$

 \xDH if $u' \xbe \xCd (u),$ then $ \xCd (u)= \xCd (u' ),$ $ \xeb (u)=
\xeb (u' ),$ $ \xee (u)= \xee (u' ),$

 \xDH $u' \xbe \xeb (u)$ iff $u \xbe \xee (u' ),$

 \xDH $u \xbe \xeb (u' ),$ $u' \xbe \xeb (u'' )$ $ \xch $ $u \xbe \xeb
(u'' ),$

or, equivalently,

 \xDH $u \xbe \xeb (u' )$ $ \xch $ $ \xeb (u' ) \xcc \xeb (u).$

 \xEJ

 \xEj

\ed

Note that we will generally loose much ignorance in applying the next two
Facts.

\bfa

$\hspace{0.01em}$

% (+++ Orig. No.:  Fact D-4.7.1 +++)

\label{Fact D-4.7.1}

A partial (strict) order on $W$ can be extended to a total (strict) order.

\efa

\subparagraph{
Proof
}

$\hspace{0.01em}$

% (+++ Orig.:  Proof +++)

Take an arbitrary enumeration of all pairs a, $b$ of $W:$
$ \xBc a,b \xBe _{i}:$ $i \xbe \xbk.$
Suppose all $ \xBc a,b \xBe _{j}$ for $j<i$ have been ordered, and we have no
information
if $a \xeb b$ or $a \xCd b$ or $a \xee b.$ Choose arbitrarily $a \xeb b.$
A contradiction would be
a (finite) cycle involving $ \xeb.$ But then we would have known already
that
$b \xec a.$
$ \xcz $
\\[3ex]

We use now a generalized abstract nonsense result, taken from  \cite{LMS01},
which
must also be part of the folklore:

\bfa

$\hspace{0.01em}$

% (+++ Orig. No.:  Fact D-4.7.2 +++)

\label{Fact D-4.7.2}

\label{Fact Abs-Rel-Ext}

Given a set $X$ and a binary relation $R$ on $X,$ there exists a total
preorder (i.e.
a total, reflexive, transitive relation) $S$ on $X$ that extends $R$ such
that

$ \xcA x,y \xbe X(xSy,ySx \xch xR^{*}y)$

where $R^{*}$ is the reflexive and transitive closure of $R.$

\efa

\subparagraph{
Proof
}

$\hspace{0.01em}$

% (+++ Orig.:  Proof +++)

Define $x \xDd y$ iff $xR^{*}y$ and $yR^{*}x.$
The relation $ \xDd $ is an equivalence relation.
Let $[x]$ be the equivalence class of $x$ under $ \xDd.$ Define $[x] \xec
[y]$ iff $xR^{*}y.$
The definition of $ \xec $ does not depend on the representatives $x$ and
$y$ chosen.
The relation $ \xec $ on equivalence classes is a partial order.

Let $ \xck $ be any total order on these equivalence classes that extends
$ \xec,$ by above
Fact 1.08.

Define $ \xCf xSy$ iff $[x] \xck [y].$
The relation $S$ is total (since $ \xck $ is total) and transitive
(since $ \xck $ is transitive) and is therefore a total preorder.
It extends $R$ by the definition of $ \xec $ and the fact that $ \xck $
extends $ \xec.$
Suppose now $ \xCf xSy$ and $ \xCf ySx.$ We have $[x] \xck [y]$ and $[y]
\xck [x]$
and therefore $[x]=[y]$ by antisymmetry. Therefore $x \xDd y$ and
$xR^{*}y.$
$ \xcz $
\\[3ex]

\bfa

$\hspace{0.01em}$

% (+++ Orig. No.:  Fact D-4.7.3 +++)

\label{Fact D-4.7.3}

$ \xdo $ can be represented by a ranking iff in $ \xdc ( \xdo )$ the sets
$ \xCd (u),$ $ \xeb (u),$ $ \xee (u)$
are pairwise disjoint.

\efa

\subparagraph{
Proof
}

$\hspace{0.01em}$

% (+++ Orig.:  Proof +++)

(Outline)
By the construction of $ \xdc ( \xdo )$ and disjointness, there are no
cycles involving
$ \xeb.$ Extend the relation by
Fact \ref{Fact D-4.7.2} (page \pageref{Fact D-4.7.2}).
Let the $ \xCd (u)$ be the equivalence classes. Define $ \xCd (u) \xen
\xCd (u' )$ iff $u \xbe \xeb (u' ).$
$ \xcz $
\\[3ex]

\bp

$\hspace{0.01em}$

% (+++ Orig. No.:  Proposition D-4.7.4 +++)

\label{Proposition D-4.7.4}

Let $f: \xdy \xcp \xdp (W).$ $f$ is representable by a smooth ranked
structure iff in $ \xdc ( \xdo )$
the sets $ \xCd (u),$ $ \xeb (u),$ $ \xee (u)$ are pairwise disjoint,
where $ \xdo $ is the system
generated by $f,$ as in
Definition \ref{Definition D-4.7.1} (page \pageref{Definition D-4.7.1}).

\ep

\subparagraph{
Proof
}

$\hspace{0.01em}$

% (+++ Orig.:  Proof +++)

If the sets are not pairwise disjoint, we have a cycle. If not,
use Fact \ref{Fact D-4.7.3} (page \pageref{Fact D-4.7.3}).
$ \xcz $
\\[3ex]
\chapter{Preferential structures - Part II}
\label{Chapter Pref-II}
\section{Simplifications by domain conditions, logical properties}
\label{Section With-Domain}
\subsection{Introduction}
\label{Section 2.2.8.1}

We examine here simplifications made possible by stronger closure
conditions of
the domain $ \xdy,$ in particular $( \xcv ).$

For general preferential structures, there is nothing to show - there were
no
prerequisites about closure of the domain.

The smooth case is more interesting. The work for the not necessarily
transitive case was done already, and, as we did not know how to do
better,
we give now directly the result for the smooth transitive case,
using in an essential way $( \xcv ).$
\subsection{Smooth structures}
\label{Section 2.2.8.2}

For completeness' sake and for the reader's convenience, we
will just repeat here our result from  \cite{Sch04}, with the
slight improvement that we do not need $( \xcs )$ any more, and the
codomain
need not be $ \xdy $ any more.
The central condition is, of course, $( \xcv ),$ which we use now as we
prepare
the classical propositional case, where we have $ \xco.$
\index{CS-Smooth-Trans}

\paragraph{
Discussion of the smooth and transitive case
}

$\hspace{0.01em}$

% (+++ Orig.:  Discussion of the smooth and transitive case +++)

\label{Section Discussion of the smooth and transitive case}

In a certain way, it is not surprising that transitivity does not impose
stronger conditions in the smooth case either. Smoothness is itself a weak
kind of
transitivity: If an element is not minimal, then there is a minimal
element
below it, i.e., $x \xee y$ with $y$ not minimal is possible, there is $z'
\xeb y,$ but
then there is $z$ minimal with $x \xee z.$ This is ``almost'' $x \xee z',$
transitivity.

To obtain representation,
we will combine here the ideas of the smooth, but not necessarily
transitive
case with those of the general transitive case - as the reader will have
suspected. Thus, we will index again with trees, and work with (suitably
adapted) admissible sequences for the construction of the trees. In the
construction of the admissible sequences, we were careful to repair all
damage
done in previous steps. We have to add now reparation of all damage done
by
using transitivity, i.e., the transitivity of the relation might destroy
minimality, and we have to construct minimal elements below all elements
for
which we thus destroyed minimality. Both cases are combined by considering
immediately all $Y$ s.t. $x \xbe Y-H(U).$ Of course, the properties
described in
Fact \ref{Fact HU} (page \pageref{Fact HU})  play again a central role.

The (somewhat complicated) construction will be commented on in more
detail
below.

Note that even beyond Fact \ref{Fact HU} (page \pageref{Fact HU}),
closure of the domain under finite
unions is used in the construction of the trees. This - or something like
it -
is necessary, as we have to respect the hulls of all elements treated so
far
(the predecessors), and not only of the first element, because of
transitivity.
For the same reason, we need more bookkeeping, to annotate all the hulls
(or
the union of the respective $U' $s) of all predecessors to be respected.
One can
perhaps do with a weaker operation than union - i.e. just look at the
hulls of
all U's separately, to obtain a transitive construction where unions are
lacking, see the case of plausibility logic below - but we have not
investigated
this problem.

To summarize: we combine the ideas from the transitive general case and
the
simple smooth case, using the crucial
Fact \ref{Fact HU} (page \pageref{Fact HU})  to show that the construction
goes through. The construction leaves still some freedom, and
modifications
are possible as indicated below in the course of the proof. The
construction is
perhaps the most complicated in the entire book, as it combines several
ideas,
some of which are already somewhat involved. If necessary, the proof can
certainly still be elaborated, and its main points (use of a suitable
$H(U)$ to
avoid $U,$ successive repair of damage done in the construction, trees as
indexing) may probably be used in other contexts, too.

\paragraph{
The construction:
}

$\hspace{0.01em}$

% (+++ Orig.:  The construction: +++)

\label{Section The construction:}

Recall that $ \xdy $ will be closed under finite unions
in this section, and let again $ \xbm: \xdy \xcp \xdp (Z).$

Proposition \ref{Proposition D-4.5.3} (page \pageref{Proposition D-4.5.3})  is
the representation result for the smooth transitive case.

\bp

$\hspace{0.01em}$

% (+++ Orig. No.:  Proposition D-4.5.3 +++)

\label{Proposition D-4.5.3}

Let $ \xdy $ be closed under finite unions, and $ \xbm: \xdy \xcp \xdp
(Z).$
Then there is a $ \xdy -$smooth transitive preferential structure $ \xdz
,$ s.t. for all
$X \xbe \xdy $ $ \xbm (X)= \xbm_{ \xdz }(X)$ iff $ \xbm $ satisfies $(
\xbm \xcc ),$ $( \xbm PR),$ $( \xbm CUM).$

\ep

\subparagraph{
Proof
}

$\hspace{0.01em}$

% (+++ Orig.:  Proof +++)

\paragraph{
The idea:
}

$\hspace{0.01em}$

% (+++ Orig.:  The idea: +++)

\label{Section The idea:}

We have to adapt Construction \ref{Construction D-4.4.2} (page
\pageref{Construction D-4.4.2})
(using $x-$admissible sequences) to the
transitive situation, and to our construction with trees. If
$ \xBc  \xCQ,x \xBe $ is the root,
$ \xbs_{0} \xbe \xbP \{ \xbm (Y):x \xbe Y- \xbm (Y)\}$ determines some
children of the root.
To preserve smoothness, we have to compensate and add
other children by the $ \xbs_{i+1}:$ $ \xbs_{i+1} \xbe \xbP \{ \xbm (X):x
\xbe \xbm (X),$ $ran( \xbs_{i}) \xcs X \xEd \xCQ \}.$
On the other hand, we have to pursue the same construction for the
children so
constructed. Moreover, these indirect children have to be added to those
children
of the root, which have to be compensated (as the first children are
compensated
by $ \xbs_{1})$ to preserve smoothness. Thus, we build the tree in a
simultaneous vertical
and horizontal induction.

This construction can be simplified, by considering immediately all $Y
\xbe \xdy $ s.t.
$x \xbe Y \xcC H(U)$ - independent of whether $x \xce \xbm (Y)$ (as done
in $ \xbs_{0}),$ or whether
$x \xbe \xbm (Y),$ and some child $y$ constructed before is in $Y$ (as
done in the $ \xbs_{i+1}),$ or
whether $x \xbe \xbm (Y),$ and some indirect child $y$ of $x$ is in $Y$
(to take care of
transitivity, as indicated above). We make this simplified construction.

There are two ways to proceed. First, we can take as $ \xej^{*}$ in the
trees
the transitive closure of $ \xej.$ Second, we can deviate from the idea
that
children are chosen by selection functions $f,$ and take nonempty subsets
of
elements instead, making more elements children than in the first case. We
take
the first alternative, as it is more in the spirit of the construction.

We will suppose for simplicity that $Z=K$ - the general case in easy to
obtain
by a technique similar to that in
Section \ref{Section 2.2.2} (page \pageref{Section 2.2.2}), but
complicates the picture.

\index{$t_x$}
For each $x \xbe Z,$ we construct trees $t_{x}$, which will be used to
index
different copies of $x,$ and control the relation $ \xeb.$

These trees $t_{x}$ will have the following form:

(a) the root of $t$ is $ \xBc  \xCQ,x \xBe $ or
$ \xBc U,x \xBe $ with $U \xbe \xdy $ and $x \xbe \xbm (U),$

(b) all other nodes are pairs $ \xBc Y,y \xBe,$ $Y \xbe \xdy,$ $y \xbe \xbm
(Y),$

(c) $ht(t) \xck \xbo,$

(d) if $ \xBc Y,y \xBe $ is an element in $t_{x},$ then there is some $ \xdy (y)
\xcc \{W \xbe \xdy:y \xbe W\},$ and
$f \xbe \xbP \{ \xbm (W):W \xbe \xdy (y)\}$ s.t. the set of children of
$ \xBc Y,y \xBe $ is $\{ \xBc Y \xcv W,f(W) \xBe :$ $W \xbe \xdy (y)\}.$

The first coordinate is used for bookkeeping when constructing children,
in
particular for condition (d).

The relation $ \xeb $ will essentially be determined by the subtree
relation.

We first construct the trees $t_{x}$ for those sets $U$ where $x \xbe \xbm
(U),$ and then take
care of the others. In the construction for the minimal elements,
at each level $n>0,$ we may have several ways to choose a selection
function $f_{n}$,
and each such choice leads to the construction of a different tree - we
construct all these trees. (We could also construct only one tree, but
then
the choice would have to be made coherently for different $x,U.$ It is
simpler to
construct more trees than necessary.)

We control the relation by indexing with trees, just as it was done in the
not
necessarily smooth case before.

\bd

$\hspace{0.01em}$

% (+++ Orig. No.:  Definition 3.3.3: +++)

\label{Definition 3.3.3:}

\index{$t/c$}
If $t$ is a tree with root $ \xBc a,b \xBe,$ then t/c will be the same tree,
only with the root $ \xBc c,b \xBe.$

\ed

\bcs

$\hspace{0.01em}$

% (+++ Orig. No.:  Construction 3.3.3: +++)

\label{Construction 3.3.3:}

\index{$T_x$}
(A) The set $T_{x}$ of trees $t$ for fixed $x$:

\index{$T\xbm_x$}
(1) Construction of the set $T \xbm_{x}$ of trees for those sets $U \xbe
\xdy,$ where $x \xbe \xbm (U):$

Let $U \xbe \xdy,$ $x \xbe \xbm (U).$ The trees $t_{U,x} \xbe T \xbm_{x}$
are constructed inductively,
observing simultaneously:

If $ \xBc U_{n+1},x_{n+1} \xBe $ is a child of
$ \xBc U_{n},x_{n} \xBe,$ then

(a) $x_{n+1} \xbe \xbm (U_{n+1})-H(U_{n}),$
and
(b) $U_{n} \xcc U_{n+1}$.

Set $U_{0}:=U,$ $x_{0}:=x.$

Level 0: $ \xBc U_{0},x_{0} \xBe.$

Level $n \xcp n+1$:
Let $ \xBc U_{n},x_{n} \xBe $ be in level $n.$
Suppose $Y_{n+1} \xbe \xdy,$ $x_{n} \xbe Y_{n+1},$ and $Y_{n+1} \xcC
H(U_{n}).$ Note that $ \xbm (U_{n} \xcv Y_{n+1})-H(U_{n}) \xEd \xCQ $ by
Fact \ref{Fact HU} (page \pageref{Fact HU}), (5.5), and
$ \xbm (U_{n} \xcv Y_{n+1})-H(U_{n}) \xcc \xbm (Y_{n+1})$
by Fact \ref{Fact HU} (page \pageref{Fact HU}), (3.3).
Choose $f_{n+1} \xbe \xbP \{ \xbm (U_{n} \xcv Y_{n+1})-H(U_{n}):$ $Y_{n+1}
\xbe \xdy,$ $x_{n} \xbe Y_{n+1} \xcC H(U_{n})\}$ (for the construction
of this tree, at this element), and let the set of
children of $ \xBc U_{n},x_{n} \xBe $ be
$\{ \xBc U_{n} \xcv Y_{n+1},f_{n+1}(Y_{n+1}) \xBe :$ $Y_{n+1} \xbe \xdy,$ $x_{n}
\xbe Y_{n+1} \xcC H(U_{n})\}.$
(If there is no such $Y_{n+1}$, $ \xBc U_{n},x_{n} \xBe $ has no children.)
Obviously, (a) and (b) hold.

\index{$U,x$-tree}
We call such trees $U,x-$trees.

(2) Construction of the set $T'_{x}$ of trees for the nonminimal elements.
Let $x \xbe Z.$ Construct the tree $t_{x}$ as follows (here, one tree per
$x$ suffices for
all $U)$:

Level 0: $ \xBc  \xCQ,x \xBe $

Level 1:
Choose arbitrary $f \xbe \xbP \{ \xbm (U):x \xbe U \xbe \xdy \}.$ Note
that $U \xEd \xCQ \xcp \xbm (U) \xEd \xCQ $ by $Z=K$ (by
Remark \ref{Remark D-4.4.4} (page \pageref{Remark D-4.4.4}), (1)).
Let $\{ \xBc U,f(U) \xBe :x \xbe U \xbe \xdy \}$ be the set of children of $<
\xCQ
,x>.$
This assures that the element will be nonminimal.

Level $>1$:
Let $ \xBc U,f(U) \xBe $ be an element of level 1, as $f(U) \xbe \xbm (U),$
there is
a $t_{U,f(U)} \xbe T \xbm_{f(U)}.$
Graft one of these trees $t_{U,f(U)} \xbe T \xbm_{f(U)}$ at
$ \xBc U,f(U) \xBe $ on the level 1.
This assures that a minimal element will be below it to guarantee
smoothness.

Finally, let $T_{x}:=T \xbm_{x} \xcv T'_{x}.$

(B) The relation $ \xej $ between trees:
For $x,y \xbe Z,$ $t \xbe T_{x}$, $t' \xbe T_{y}$, set $t \xem t' $ iff
for some $Y$
$ \xBc Y,y \xBe $ is a child of
the root $ \xBc X,x \xBe $ in $t,$ and $t' $ is the subtree of $t$ beginning
at this $ \xBc Y,y \xBe.$

(C) The structure $ \xdz $:

Let $ \xdz $ $:=$ $ \xCc $ $\{ \xBc x,t_{x} \xBe :$ $x \xbe Z,$ $t_{x} \xbe
T_{x}\}$
,
$ \xBc x,t_{x} \xBe  \xee  \xBc y,t_{y} \xBe $ iff $t_{x} \xem^{*}t_{y}$ $
\xCe.$

The rest of the proof are simple observations.

\ecs

\bfa

$\hspace{0.01em}$

% (+++ Orig. No.:  Fact 3.3.9 +++)

\label{Fact 3.3.9}

(1) If $t_{U,x}$ is an $U,x-$tree, $ \xBc U_{n},x_{n} \xBe $ an element of
$t_{U,x}$
,
$ \xBc U_{m},x_{m} \xBe $ a direct or
indirect child of $ \xBc U_{n},x_{n} \xBe,$ then $x_{m} \xce H(U_{n}).$

(2) Let $ \xBc Y_{n},y_{n} \xBe $ be an element in $t_{U,x} \xbe T \xbm_{x}$,
$t' $
the subtree
starting at $ \xBc Y_{n},y_{n} \xBe,$
then $t' $ is a $Y_{n},y_{n}-tree.$

(3) $ \xeb $ is free from cycles.

(4) If $t_{U,x}$ is an $U,x-$tree, then
$ \xBc x,t_{U,x} \xBe $ is $ \xeb -$minimal in $ \xdz \xex U.$

(5) No $ \xBc x,t_{x} \xBe,$ $t_{x} \xbe T'_{x}$ is minimal in any $ \xdz \xex
U,$
$U \xbe \xdy.$

(6) Smoothness is respected for the elements of the form
$ \xBc x,t_{U,x} \xBe.$

(7) Smoothness is respected for the elements of the form
$ \xBc x,t_{x} \xBe $ with $t_{x} \xbe T'_{x}.$

(8) $ \xbm = \xbm_{ \xdz }.$

\efa

\subparagraph{
Proof
}

$\hspace{0.01em}$

% (+++ Orig.:  Proof +++)

(1) trivial by (a) and (b).

(2) trivial by (a).

(3) Note that no $ \xBc x,t_{x} \xBe $ $t_{x} \xbe T'_{x}$ can be smaller than
any
other element (smaller
elements require $U \xEd \xCQ $ at the root). So no cycle involves
any such $ \xBc x,t_{x} \xBe.$
Consider now $ \xBc x,t_{U,x} \xBe,$ $t_{U,x} \xbe T \xbm_{x}$. For any
$ \xBc y,t_{V,y} \xBe  \xeb  \xBc x,t_{U,x} \xBe,$ $y \xce H(U)$ by (1),
but $x \xbe \xbm (U) \xcc H(U),$ so $x \xEd y.$

(4) This is trivial by (1).

(5) Let $x \xbe U \xbe \xdy,$ then $f$ as used in the construction of
level 1 of $t_{x}$ chooses
$y \xbe \xbm (U) \xEd \xCQ,$ and some $ \xBc y,t_{U,y} \xBe $ is in $ \xdz \xex
U$
and below
$ \xBc x,t_{x} \xBe.$

(6) Let $x \xbe A \xbe \xdy,$ we have to show that either
$ \xBc x,t_{U,x} \xBe $ is minimal in $ \xdz \xex A,$ or that
there is $ \xBc y,t_{y} \xBe  \xeb  \xBc x,t_{U,x} \xBe $ minimal in $ \xdz \xex
A.$
Case 1, $A \xcc H(U)$: Then $ \xBc x,t_{U,x} \xBe $ is minimal in $ \xdz \xex
A,$
again by (1).
Case 2, $A \xcC H(U)$: Then A is one of the $Y_{1}$ considered for level
1. So there is
$ \xBc U \xcv A,f_{1}(A) \xBe $ in level 1 with $f_{1}(A) \xbe \xbm (A) \xcc A$
by
Fact \ref{Fact HU} (page \pageref{Fact HU}), (3.3). But note that by (1)
all elements below $ \xBc U \xcv A,f_{1}(A) \xBe $ avoid $H(U \xcv A).$ Let $t$
be
the subtree of $t_{U,x}$
beginning at $ \xBc U \xcv A,f_{1}(A) \xBe,$ then by (2) $t$ is one of the $U
\xcv
A,f_{1}(A)-trees,$ and
$ \xBc f_{1}(A),t \xBe $ is minimal in $ \xdz \xex U \xcv A$ by (4), so in $
\xdz
\xex A,$ and
$ \xBc f_{1}(A),t \xBe  \xeb  \xBc x,t_{U,x} \xBe.$

(7) Let $x \xbe A \xbe \xdy,$ $ \xBc x,t_{x} \xBe,$ $t_{x} \xbe T_{x}',$ and
consider the subtree $t$ beginning at
$ \xBc A,f(A) \xBe,$
then $t$ is one of the $A,f(A)-$trees, and
$ \xBc f(A),t \xBe $ is minimal in $ \xdz \xex A$ by (4).

(8) Let $x \xbe \xbm (U).$ Then any $ \xBc x,t_{U,x} \xBe $ is $ \xeb -$minimal
in $
\xdz \xex U$ by (4), so $x \xbe \xbm_{ \xdz }(U).$
Conversely, let $x \xbe U- \xbm (U).$ By (5), no
$ \xBc x,t_{x} \xBe $ is minimal in $U.$ Consider now some
$ \xBc x,t_{V,x} \xBe  \xbe \xdz,$ so $x \xbe \xbm (V).$ As $x \xbe U- \xbm
(U),$
$U \xcC H(V)$
by Fact \ref{Fact HU} (page \pageref{Fact HU}), (5.4). Thus $U$ was
considered in the construction of level 1 of $t_{V,x}.$ Let $t$ be the
subtree of $t_{V,x}$
beginning at $ \xBc V \xcv U,f_{1}(U) \xBe,$ by $ \xbm (V \xcv U)-H(V) \xcc
\xbm
(U)$
(Fact \ref{Fact HU} (page \pageref{Fact HU}), (3.3)), $f_{1}(U) \xbe \xbm (U)
\xcc U,$
and $ \xBc f_{1}(U),t \xBe  \xeb  \xBc x,t_{V,x} \xBe.$

$ \xcz $ (Fact \ref{Fact 3.3.9} (page \pageref{Fact 3.3.9})  and Proposition
\ref{Proposition D-4.5.3} (page \pageref{Proposition D-4.5.3}) )
\\[3ex]
\subsection{Ranked structures}
\label{Section 2.2.8.3}

We summarize for completess' sake results from  \cite{Sch04}:

First two results for the case without copies
(Proposition \ref{Proposition D-5.3.4} (page \pageref{Proposition D-5.3.4})  and
Proposition \ref{Proposition D-5.3.5} (page \pageref{Proposition D-5.3.5}) ).

\bp

$\hspace{0.01em}$

% (+++ Orig. No.:  Proposition D-5.3.4 +++)

\label{Proposition D-5.3.4}

Let $ \xdy \xcc \xdp (U)$ be closed under finite unions.
Then $( \xbm \xcc ),$ $( \xbm \xCQ ),$ $( \xbm =)$ characterize ranked
structures for which for all
$X \xbe \xdy $ $X \xEd \xCQ $ $ \xch $ $ \xbm_{<}(X) \xEd \xCQ $ hold,
i.e. $( \xbm \xcc ),$ $( \xbm \xCQ ),$ $( \xbm =)$ hold in such
structures for $ \xbm_{<},$ and if they hold for some $ \xbm,$ we can
find a ranked relation
$<$ on $U$ s.t. $ \xbm = \xbm_{<}.$ Moreover, the structure can be choosen
$ \xdy -$smooth.

\ep

For the following representation result, we assume only $( \xbm \xCQ
fin),$ but the
domain has to contain singletons.

\bp

$\hspace{0.01em}$

% (+++ Orig. No.:  Proposition D-5.3.5 +++)

\label{Proposition D-5.3.5}

\index{singletons}
Let $ \xdy \xcc \xdp (U)$ be closed under finite unions, and contain
singletons.
Then $( \xbm \xcc ),$ $( \xbm \xCQ fin),$ $( \xbm =),$ $( \xbm \xbe )$
characterize ranked structures for which
for all finite $X \xbe \xdy $ $X \xEd \xCQ $ $ \xch $ $ \xbm_{<}(X) \xEd
\xCQ $ hold, i.e. $( \xbm \xcc ),$ $( \xbm \xCQ fin),$ $( \xbm =),$ $(
\xbm \xbe )$
hold in such structures for $ \xbm_{<},$ and if they hold for some $ \xbm
,$ we can find
a ranked relation $<$ on $U$ s.t. $ \xbm = \xbm_{<}.$

\ep

Note that the prerequisites of
Proposition \ref{Proposition D-5.3.5} (page \pageref{Proposition D-5.3.5})
hold in particular in the case
of ranked structures without copies, where all elements of $U$ are present
in the
structure - we need infinite descending chains to have $ \xbm (X)= \xCQ $
for $X \xEd \xCQ.$

We turn now to the general case, where every element may occur in several
copies.

\bfa

$\hspace{0.01em}$

% (+++ Orig. No.:  Fact D-5.3.6 +++)

\label{Fact D-5.3.6}

(1) $( \xbm \xcc )+( \xbm PR)+( \xbm =)+( \xbm \xcv )+( \xbm \xbe )$ do
not imply representation by a ranked
structure.

(2) The infinitary version of $( \xbm \xFO ):$

\index{$(\xbm \xFO \xca)$}
$( \xbm \xFO \xca )$ $ \xbm ( \xcV \{A_{i}:i \xbe I\})$ $=$ $ \xcV \{ \xbm
(A_{i}):i \xbe I' \}$ for some $I' \xcc I.$

will not always hold in ranked structures.

\efa

We assume again the existence of singletons for the following
representation
result.

\bp

$\hspace{0.01em}$

% (+++ Orig. No.:  Proposition D-5.3.7 +++)

\label{Proposition D-5.3.7}

Let $ \xdy $ be closed under finite unions and contain singletons. Then
$( \xbm \xcc )+( \xbm PR)+( \xbm \xFO )+( \xbm \xcv )+( \xbm \xbe )$
characterize ranked structures.
\subsection{The logical properties with definability preservation}
\label{Section 2.2.8.4}
%  Section (5.4):  The logical properties with dp
%  Section (5.4):  The logical properties with dp
% %
% ================================================

\ep

We repeat for completeness' sake:

\bp

$\hspace{0.01em}$

% (+++ Orig. No.:  Proposition D-5.4.1 +++)

\label{Proposition D-5.4.1}

Let $ \xcn $ be a logic for $ \xdl.$ Set $T^{ \xdm }:=Th( \xbm_{ \xdm
}(M(T))),$
where $ \xdm $ is a preferential structure.

(1) Then there is a (transitive) definability preserving
classical preferential model $ \xdm $ s.t. $ \ol{ \ol{T} }=T^{ \xdm }$ iff

\index{$(LLE)$}
\index{$(CCL)$}
\index{$(SC)$}
\index{$(PR)$}
\index{$(CUM)$}
(LLE) $ \ol{T}= \ol{T' }$ $ \xch $ $ \ol{ \ol{T} }= \ol{ \ol{T' } },$

(CCL) $ \ol{ \ol{T} }$ is classically closed,

(SC) $T \xcc \ol{ \ol{T} },$

(PR) $ \ol{ \ol{T \xcv T' } } \xcc \ol{ \ol{ \ol{T} } \xcv T' }$

for all $T,T' \xcc \xdl.$

(2) The structure can be chosen smooth, iff, in addition

(CUM) $T \xcc \ol{T' } \xcc \ol{ \ol{T} }$ $ \xch $ $ \ol{ \ol{T} }= \ol{
\ol{T' } }$

holds.

\ep

The proof is an immediate consequence of
Proposition \ref{Proposition D-5.4.2} (page \pageref{Proposition D-5.4.2})  and
Proposition \ref{Proposition D-4.5.3} (page \pageref{Proposition D-4.5.3}). $
\xcz $
\\[3ex]

\bp

$\hspace{0.01em}$

% (+++ Orig. No.:  Proposition D-5.4.2 +++)

\label{Proposition D-5.4.2}

Consider for a logic $ \xcn $ on $ \xdl $ the properties

(LLE) $ \ol{T}= \ol{T' }$ $ \xch $ $ \ol{ \ol{T} }= \ol{ \ol{T' } },$

(CCL) $ \ol{ \ol{T} }$ is classically closed,

(SC) $T \xcc \ol{ \ol{T} },$

(PR) $ \ol{ \ol{T \xcv T' } } \xcc \ol{ \ol{ \ol{T} } \xcv T' },$

(CUM) $T \xcc \ol{T' } \xcc \ol{ \ol{T} }$ $ \xch $ $ \ol{ \ol{T} }= \ol{
\ol{T' } }$

for all $T,T' \xcc \xdl,$

and for a function $ \xbm: \xdD_{ \xdl } \xcp \xdp (M_{ \xdl })$ the
properties

\index{$(\xbm dp)$}
\index{$(\xbm \xcc)$}
\index{$(\xbm PR)$}
\index{$(\xbm CUM)$}
$( \xbm dp)$ $ \xbm $ is definability preserving, i.e. $ \xbm (M(T))=M(T'
)$ for some $T' $

$( \xbm \xcc )$ $ \xbm (X) \xcc X,$

$( \xbm PR)$ $X \xcc Y$ $ \xch $ $ \xbm (Y) \xcs X \xcc \xbm (X),$

$( \xbm CUM)$ $ \xbm (X) \xcc Y \xcc X$ $ \xch $ $ \xbm (X)= \xbm (Y)$

for all $X,Y \xbe \xdD_{ \xdl }.$

It then holds:

(a) If $ \xbm $ satisfies $( \xbm dp),$ $( \xbm \xcc ),$ $( \xbm PR),$
then $ \xcn $ defined by $ \ol{ \ol{T} }:=T^{ \xbm }:=$
$Th( \xbm (M(T)))$ satisfies (LLE), (CCL), (SC), (PR).
If $ \xbm $ satisfies in addition $( \xbm CUM),$ then (CUM) will hold,
too.

(b) If $ \xcn $ satisfies (LLE), (CCL), (SC), (PR),
then there is $ \xbm: \xdD_{ \xdl } \xcp \xdp (M_{ \xdl })$ s.t. $ \ol{
\ol{T} }=T^{ \xbm }$
for all $T \xcc \xdl $ and $ \xbm $ satisfies $( \xbm dp),$ $( \xbm \xcc
),$ $( \xbm PR).$
If, in addition, (CUM) holds, then $( \xbm CUM)$ will hold, too.

The proof follows from Proposition \ref{Proposition Alg-Log} (page
\pageref{Proposition Alg-Log}). $ \xcz $
\\[3ex]

\ep

The logical properties of definability preserving ranked structures are
straightforward now, and left to the reader.
\section{$\xda$-ranked structures}
\label{Section A-ranked-Rep}

We do now the completeness proofs for the preferential part of
hierarchical conditionals.
All motivation etc. will be found in
Section \ref{Section Hierarchical-Conditionals} (page \pageref{Section
Hierarchical-Conditionals}).

First the basic semantical definition:
\label{Section Hier-Def-2}
\index{Section Hier-Def-2}

\bd

$\hspace{0.01em}$

% (+++ Orig. No.:  Definition A-ranked +++)

\label{Definition A-ranked}

Let $ \xdA $ be a fixed set, and $ \xda $ a finite, totally ordered (by
$<)$ disjoint cover
by non-empty subsets of $ \xdA.$

For $x \xbe \xdA,$ let $rg(x)$ the unique $A \xbe \xda $ such that $x
\xbe A,$ so $rg(x)<rg(y)$ is
defined in the natural way.

A preferential structure $ \xBc  \xdx, \xeb  \xBe $ $( \xdx $ a set of pairs
$ \xBc x,i \xBe )$ is called $ \xda -$ranked
iff for all $x,x' $ $rg(x)<rg(x' )$ implies $ \xBc x,i \xBe  \xeb  \xBc x',i' 
\xBe $ for all
$ \xBc x,i \xBe, \xBc x',i'  \xBe  \xbe \xdx.$
\subsection{
Representation results for $\xda$-ranked structures
}
\subsubsection{
Discussion
}
\label{Section Hier-ARepr-Disc-1}
\index{Section Hier-ARepr-Disc-1}

\ed

The not necessarily smooth and the smooth case will be treated
differently.

Strangely, the smooth case is simpler, as an added new layer in the proof
settles it. Yet, this is not surprising when looking closer, as minimal
elements
never have
higher rank, and we know from $( \xbm CUM)$ that minimizing by minimal
elements
suffices. All we have to add that any element in the minimal layer
minimizes
any element higher up.

In the simple, not necessarily smooth, case, we have to go deeper into the
original proof to obtain the result.

The following idea, inspired by the treatment of the smooth case, will not
work: Instead of minimizing by arbitrary elements, minimize only by
elements of
minimal rank, as the following example shows. If it worked, we might add
just
another layer to the original proof without $( \xbm \xda ),$
(see Definition \ref{Definition Mu-A} (page \pageref{Definition Mu-A}) ), as in
the smooth case.

\be

$\hspace{0.01em}$

% (+++ Orig. No.:  Example 6.1 +++)

\label{Example 6.1}

Consider the base set $\{a,b,c\},$ $ \xbm (\{a,b,c\})=\{b\},$ $ \xbm
(\{a,b\})=\{a,b\},$
$ \xbm (\{a,c\})= \xCQ,$ $ \xbm (\{b,c\})=\{b\},$ $ \xda $ defined by
$\{a,b\}<\{c\}.$

Obviously, $( \xbm \xda )$ is satisfied. $ \xbm $ can be represented by
the (not transitive!)
relation $a \xeb c \xeb a,$ $b \xeb c,$ which is $ \xda -$ranked.

But trying to minimize a in $\{a,b,c\}$ in the minimal layer will lead to
$b \xeb a,$
and thus $a \xce \xbm (\{a,b\}),$ which is wrong.

$ \xcz $
\\[3ex]

\ee

The proofs of the general and transitive general case are (minor)
adaptations of the proofs in
Section \ref{Section Pref-Details} (page \pageref{Section Pref-Details}).
For the smooth case, we
only have to add a supplementary layer
in the end (Fact \ref{Fact Smooth-to-A} (page \pageref{Fact Smooth-to-A}) ),
which
will make the construction $ \xda -$ranked.

In the following, we will assume the partition $ \xda $ to be given. We
could also
construct it from the properties of $ \xbm,$ but this would need stronger
closure
properties of the domain. The construction of $ \xda $ is more difficult
than the
construction of the ranking in fully ranked structures, as $x \xbe \xbm
(X),$ $y \xbe X- \xbm (X)$
will guarantee only $rg(x) \xck rg(y),$ and not $rg(x)<rg(y),$ as is the
case in the
latter situation. This corresponds to the separate treatment of the $ \xba
$ and
other formulas in the logical version, discussed
in Section \ref{Section Hier-ARepr-Logic} (page \pageref{Section
Hier-ARepr-Logic}).
\subsubsection{
$\xda$-ranked general and transitive structures
}
\paragraph{
Introduction
}
\label{Section Hier-ARepr-General-Intro}
\index{Section Hier-ARepr-General-Intro}

We will show here the following representation result:

Let $ \xda $ be given.

An operation $ \xbm: \xdy \xcp \xdp (Z)$ is representable by an $ \xda
-$ranked preferential
structure iff $ \xbm $ satisfies $( \xbm \xcc ),$ $( \xbm PR),$ $( \xbm
\xda )$ (Proposition \ref{Proposition 6.3} (page \pageref{Proposition 6.3}) ),
and,
moreover, the structure can be chosen transitive (Proposition \ref{Proposition
6.5} (page \pageref{Proposition 6.5}) ).

Note that we carefully avoid any unnecessary assumptions about the domain
$ \xdy \xcc \xdp (Z)$ of the function $ \xbm.$

\bd

$\hspace{0.01em}$

% (+++ Orig. No.:  Definition Mu-A +++)

\label{Definition Mu-A}

We define a new condition:

Let $ \xda $ be given as defined in Definition \ref{Definition A-ranked} (page
\pageref{Definition A-ranked}).

$( \xbm \xda )$ If $X \xbe \xdy,$ $A,A' \xbe \xda,$ $A<A',$ $X \xcs A
\xEd \xCQ,$ $X \xcs A' \xEd \xCQ $ then $ \xbm (X) \xcs A' = \xCQ.$

\ed

This new condition will be central for the modified representation.
\paragraph{The basic, not necessarily transitive, case}
\label{Section Hier-ARepr-General-Intrans-Core}
\index{Section Hier-ARepr-General-Intrans-Core}

\bco

$\hspace{0.01em}$

% (+++ Orig. No.:  Corollary 6.2 +++)

\label{Corollary 6.2}

Let $ \xbm: \xdy \xcp \xdp (Z)$ satisfy $( \xbm \xcc ),$ $( \xbm PR),$ $(
\xbm \xda ),$ and let $U \xbe \xdy.$

If $x \xbe U$ and $ \xcE x' \xbe U.rg(x' )<rg(x),$ then $ \xcA f \xbe
\xbP_{x}.ran(f) \xcs U \xEd \xCQ.$

\eco

\subparagraph{
Proof
}

$\hspace{0.01em}$

% (+++ Orig.:  Proof +++)

By $( \xbm \xda )$ $x \xce \xbm (U),$ thus by Claim \ref{Claim Mu-f} (page
\pageref{Claim Mu-f})
$ \xcA f \xbe \xbP_{x}.ran(f) \xcs U \xEd \xCQ.$ $ \xcz $
\\[3ex]

\bp

$\hspace{0.01em}$

% (+++ Orig. No.:  Proposition 6.3 +++)

\label{Proposition 6.3}

Let $ \xda $ be given.

An operation $ \xbm: \xdy \xcp \xdp (Z)$ is representable by an $ \xda
-$ranked preferential
structure iff $ \xbm $ satisfies $( \xbm \xcc ),$ $( \xbm PR),$ $( \xbm
\xda ).$

\ep

\subparagraph{
Proof
}

$\hspace{0.01em}$

% (+++ Orig.:  Proof +++)

One direction is trivial. The central argument is: If $a \xeb b$ in $X,$
and $X \xcc Y,$
then $a \xeb b$ in $Y,$ too.

We turn to the other direction. The preferential structure is defined in
Construction \ref{Construction 6.1} (page \pageref{Construction 6.1}), Claim
\ref{Claim 6.4} (page \pageref{Claim 6.4})  shows
representation.

\bcs

$\hspace{0.01em}$

% (+++ Orig. No.:  Construction 6.1 +++)

\label{Construction 6.1}

Let $ \xdx:=\{ \xBc x,f \xBe :x \xbe Z$ $ \xcu $ $f \xbe \xbP_{x}\},$ and
$ \xBc x',f'  \xBe  \xeb  \xBc x,f \xBe $ $: \xcr $ $x' \xbe ran(f)$ or $rg(x'
)<rg(x).$

Note that, as $ \xda $ is given, we also know $rg(x).$

Let $ \xdz:= \xBc  \xdx, \xeb  \xBe.$

\ecs

Obviously, $ \xdz $ is $ \xda -$ranked.

\bc

$\hspace{0.01em}$

% (+++ Orig. No.:  Claim 6.4 +++)

\label{Claim 6.4}

For $U \xbe \xdy,$ $ \xbm (U)= \xbm_{ \xdz }(U).$

\ec

\subparagraph{
Proof
}

$\hspace{0.01em}$

% (+++ Orig.:  Proof +++)

By Claim \ref{Claim Mu-f} (page \pageref{Claim Mu-f}), it suffices to show that
for all $U \xbe
\xdy $
$x \xbe \xbm_{ \xdz }(U)$ $ \xcr $ $x \xbe U$ and $ \xcE f \xbe
\xbP_{x}.ran(f) \xcs U= \xCQ.$ So let $U \xbe \xdy.$

`` $ \xcp $ '': If $x \xbe \xbm_{ \xdz }(U),$ then there is
$ \xBc x,f \xBe $ minimal in $ \xdx \xex U$ - where
$ \xdx \xex U:=\{ \xBc x,i \xBe  \xbe \xdx:x \xbe U\}),$ so $x \xbe U,$ and
there
is no $ \xBc x',f'  \xBe  \xeb  \xBc x,f \xBe,$ $x' \xbe U,$ so by $ \xbP_{x' }
\xEd \xCQ $
there is no $x' \xbe ran(f),$ $x' \xbe U,$ but then
$ran(f) \xcs U= \xCQ.$

`` $ \xcq $ '': If $x \xbe U,$ and there is $f \xbe \xbP_{x}$, $ran(f)
\xcs U= \xCQ,$
then by Corollary \ref{Corollary 6.2} (page \pageref{Corollary 6.2}), there is
no $x' \xbe U,$ $rg(x'
)<rg(x),$
so $ \xBc x,f \xBe $ is
minimal in $ \xdx \xex U.$

$ \xcz $ (Claim \ref{Claim 6.4} (page \pageref{Claim 6.4})  and Proposition
\ref{Proposition 6.3} (page \pageref{Proposition 6.3}) )
\\[3ex]
\paragraph{The transitive case}
\label{Section Hier-ARepr-General-Trans}
\index{Section Hier-ARepr-General-Trans}

\bp

$\hspace{0.01em}$

% (+++ Orig. No.:  Proposition 6.5 +++)

\label{Proposition 6.5}

Let $ \xda $ be given.

An operation $ \xbm: \xdy \xcp \xdp (Z)$ is representable by an $ \xda
-$ranked transitive
preferential structure iff $ \xbm $ satisfies $( \xbm \xcc ),$ $( \xbm
PR),$ $( \xbm \xda ).$

\ep

\bcs

$\hspace{0.01em}$

% (+++ Orig. No.:  Construction 6.2 +++)

\label{Construction 6.2}

\index{$T_x$}
(1) For $x \xbe Z,$ let $T_{x}$ be the set of trees $t_{x}$ s.t.

(a) all nodes are elements of $Z,$

(b) the root of $t_{x}$ is $x,$

(c) $height(t_{x}) \xck \xbo,$

(d) if $y$ is an element in $t_{x}$, then there is $f \xbe \xbP_{y}:=
\xbP \{Y \xbe \xdy $: $y \xbe Y- \xbm (Y)\}$
s.t. the set of children of $y$ is $ran(f) \xcv \{y' \xbe Z:rg(y'
)<rg(y)\}.$

(2) For $x,y \xbe Z,$ $t_{x} \xbe T_{x}$, $t_{y} \xbe T_{y}$, set $t_{x}
\xem t_{y}$ iff $y$ is a (direct) child
of the root $x$ in $t_{x}$, and $t_{y}$ is the subtree of $t_{x}$
beginning at $y.$

(3) Let $ \xdz $ $:=$ $ \xBc $ $\{ \xBc x,t_{x} \xBe :$ $x \xbe Z,$ $t_{x} \xbe
T_{x}\}$,
$ \xBc x,t_{x} \xBe  \xee  \xBc y,t_{y} \xBe $ iff $t_{x} \xem t_{y}$ $ \xBe.$

\ecs

\bfa

$\hspace{0.01em}$

% (+++ Orig. No.:  Fact 6.6 +++)

\label{Fact 6.6}

(1) The construction ends at some $y$ iff $ \xdy_{y}= \xCQ $ and there is
no $y' $ s.t.
$rg(y' )<rg(y),$ consequently $T_{x}=\{x\}$ iff $ \xdy_{x}= \xCQ $ and
there are no $x' $ with
lesser rang. (We identify the tree of height 1 with its root.)

(2) We define a special tree $tc_{x}$ for all $x:$ For all nodes $y$ in
$tc_{x},$ the
successors are as follows:
if $ \xdy_{y} \xEd \xCQ,$ then $z$ is an successor iff $z=y$ or
$rg(z)<rg(y);$
if $ \xdy_{y}= \xCQ,$ then $z$ is an successor iff $rg(z)<rg(y).$
(In the first case, we make $f \xbe \xdy_{y}$ always choose $y$ itself.)
$tc_{x}$ is an element of $T_{x}.$ Thus, with (1), $T_{x} \xEd \xCQ $ for
any $x.$
Note: $tc_{x}=x$ iff $ \xdy_{x}= \xCQ $ and $x$ has minimal rang.

(3) If $f \xbe \xbP_{x}$, then the tree $tf_{x}$ with root $x$ and
otherwise
composed of the subtrees $tc_{y}$ for $y \xbe ran(f) \xcv \{y':rg(y'
)<rg(y)\}$
is an element of $T_{x}$.
(Level 0 of $tf_{x}$ has $x$ as element, the $t_{y}' s$ begin at level 1.)

(4) If $y$ is an element in $t_{x}$ and $t_{y}$ the subtree of $t_{x}$
starting at
$y,$ then $t_{y} \xbe T_{y}$.

(5) $ \xBc x,t_{x} \xBe  \xee  \xBc y,t_{y} \xBe $ implies $y \xbe ran(f) \xcv
\{x':rg(x'
)<rg(x)\}$ for some $f \xbe \xbP_{x}.$

$ \xcz $
\\[3ex]

\efa

Claim \ref{Claim 6.7} (page \pageref{Claim 6.7})  shows basic representation.

\bc

$\hspace{0.01em}$

% (+++ Orig. No.:  Claim 6.7 +++)

\label{Claim 6.7}

$ \xcA U \xbe \xdy. \xbm (U)= \xbm_{ \xdz }(U)$

\ec

\subparagraph{
Proof
}

$\hspace{0.01em}$

% (+++ Orig.:  Proof +++)

By Claim \ref{Claim Mu-f} (page \pageref{Claim Mu-f}), it suffices to show that
for all $U \xbe
\xdy $
$x \xbe \xbm_{ \xdz }(U)$ $ \xcr $ $x \xbe U$ $ \xcu $ $ \xcE f \xbe
\xbP_{x}.ran(f) \xcs U= \xCQ.$

Fix $U \xbe \xdy.$

`` $ \xcp $ '': $x \xbe \xbm_{ \xdz }(U)$ $ \xcp $ ex. $ \xBc x,t_{x} \xBe $
minimal
in $ \xdz \xex U,$ thus $x \xbe U$ and there is no
$ \xBc y,t_{y} \xBe  \xbe \xdz,$
$ \xBc y,t_{y} \xBe  \xeb  \xBc x,t_{x} \xBe,$ $y \xbe U.$ Let $f$ define the
first part of
the set of children of
the root $x$ in $t_{x}$. If $ran(f) \xcs U \xEd \xCQ,$ if $y \xbe U$ is
a child of $x$ in $t_{x}$, and if
$t_{y}$ is the subtree of $t_{x}$ starting at $y,$ then $t_{y} \xbe T_{y}$
and
$ \xBc y,t_{y} \xBe  \xeb  \xBc x,t_{x} \xBe,$
contradicting minimality of $ \xBc x,t_{x} \xBe $ in $ \xdz \xex U.$ So $ran(f)
\xcs
U= \xCQ.$

`` $ \xcq $ '': Let $x \xbe U,$ and $ \xcE f \xbe \xbP_{x}.ran(f) \xcs U=
\xCQ.$ By Corollary \ref{Corollary 6.2} (page \pageref{Corollary 6.2}),
there is no $x' \xbe U,$
$rg(x' )<rg(x).$ If $ \xdy_{x}= \xCQ,$ then the tree $tc_{x}$ has no $
\xem -$successors in $U,$ and
$ \xBc x,tc_{x} \xBe $ is $ \xee -$minimal in $ \xdz \xex U.$ If $ \xdy_{x} \xEd
\xCQ $ and $f \xbe \xbP_{x}$ s.t. $ran(f) \xcs U= \xCQ,$ then
$ \xBc x,tf_{x} \xBe $ is again $ \xee -$minimal in $ \xdz \xex U.$

$ \xcz $
\\[3ex]

We consider now the transitive closure of $ \xdz.$ (Recall that $
\xeb^{*}$ denotes the
transitive closure of $ \xeb.)$ Claim \ref{Claim 6.8} (page \pageref{Claim 6.8})
 shows that
transitivity does not
destroy what we have achieved. The trees $tf_{x}$ play a crucial role in
the
demonstration.

\bc

$\hspace{0.01em}$

% (+++ Orig. No.:  Claim 6.8 +++)

\label{Claim 6.8}

Let $ \xdz ' $ $:=$ $ \xBc $ $\{ \xBc x,t_{x} \xBe :$ $x \xbe Z,$ $t_{x} \xbe
T_{x}\}$,
$ \xBc x,t_{x} \xBe  \xee  \xBc y,t_{y} \xBe $ iff $t_{x} \xem^{*}t_{y}$ $
\xBe.$
Then $ \xbm_{ \xdz }= \xbm_{ \xdz ' }.$

\ec

\subparagraph{
Proof
}

$\hspace{0.01em}$

% (+++ Orig.:  Proof +++)

Suppose there is $U \xbe \xdy,$ $x \xbe U,$ $x \xbe \xbm_{ \xdz }(U),$ $x
\xce \xbm_{ \xdz ' }(U).$
Then there must be an element $ \xBc x,t_{x} \xBe  \xbe \xdz $ with no
$ \xBc x,t_{x} \xBe  \xee  \xBc y,t_{y} \xBe $ for any $y \xbe U.$
Let $f \xbe \xbP_{x}$ determine the first part of the set of children of
$x$ in $t_{x}$,
then $ran(f) \xcs U= \xCQ,$ consider $tf_{x}.$
All elements $w \xEd x$ of $tf_{x}$ are already in $ran(f),$ or
$rg(w)<rg(x)$ holds. (Note
that the elements chosen by rang in $tf_{x}$ continue by themselves or by
another
element of even smaller rang, but the rang order is transitive.) But all
$w$ s.t.
$rg(w)<rg(x)$ were already successors at level 1 of $x$ in $tf_{x}.$
By Corollary \ref{Corollary 6.2} (page \pageref{Corollary 6.2}),
there is no $w \xbe U,$ $rg(w)<rg(x).$ Thus, no element $ \xEd x$ of
$tf_{x}$ is in $U.$
Thus there is no $ \xBc z,t_{z} \xBe  \xeb^{*} \xBc x,tf_{x} \xBe $ in $ \xdz $
with $z \xbe
U,$ so
$ \xBc x,tf_{x} \xBe $ is minimal
in $ \xdz ' \xex U,$ contradiction.

$ \xcz $ (Claim \ref{Claim 6.8} (page \pageref{Claim 6.8})  and Proposition
\ref{Proposition 6.5} (page \pageref{Proposition 6.5}) )
\\[3ex]
\subsubsection{
$\xda$-ranked smooth structures
}
\label{Section Hier-ARepr-Smooth-Intro}
\index{Section Hier-ARepr-Smooth-Intro}

All smooth cases have a simple solution. We use one of our existing proofs
for
the not necessarily $ \xda -$ranked case, and add one litte result:

\bfa

$\hspace{0.01em}$

% (+++ Orig. No.:  Fact Smooth-to-A +++)

\label{Fact Smooth-to-A}

Let $( \xbm \xda )$ hold, and let $ \xdz = \xBc  \xdx, \xeb  \xBe $ be a smooth
preferential structure
representing $ \xbm,$ i.e. $ \xbm = \xbm_{ \xdz }.$

Suppose that

$ \xBc x,i \xBe  \xeb  \xBc y,j \xBe $ implies $rg(x) \xck rg(y).$

Define $ \xdz ':= \xBc  \xdx, \xer  \xBe $ where $ \xBc x,i \xBe  \xer  \xBc y,j
\xBe $ iff
$ \xBc x,i \xBe  \xeb  \xBc y,j \xBe $ or $rg(x)<rg(y).$

Then $ \xdz ' $ is $ \xda -$ranked.

$ \xdz ' $ is smooth, too, and $ \xbm_{ \xdz }= \xbm_{ \xdz ' }=: \xbm '
.$

In addition, if $ \xeb $ is free from cycles, so is $ \xer,$ if $ \xeb $
is transitive, so is $ \xer.$

\efa

\subparagraph{
Proof
}

$\hspace{0.01em}$

% (+++ Orig.:  Proof +++)

$ \xda -$rankedness is trivial.

Suppose $ \xBc x,i \xBe $ is $ \xeb -$minimal, but not $ \xer -$minimal. Then
there
must be
$ \xBc y,j \xBe  \xer  \xBc x,i \xBe,$ $ \xBc y,j \xBe  \xeB  \xBc x,i \xBe,$
$y \xbe X,$ so
$rg(y)<rg(x).$ By $( \xbm \xda ),$ all $x \xbe \xbm (X)$ have minimal $
\xda -$rang among the elements of
$X,$ so this is impossible. Thus, $ \xbm -$minimal elements stay $ \xbm '
-$minimal, so
smoothness will also be preserved - remember that we increased the
relation.

By prerequisite, there cannot be any cycle involving only $ \xeb,$ but
the rang
order is free from cycles, too, and $ \xeb $ respects the rang order, so $
\xer $ is free
from cycles.

Let $ \xeb $ be transitive, so is the rang order. But if
$ \xBc x,i \xBe  \xeb  \xBc y,j \xBe $ and
$rg(y)<rg(z)$ for some $ \xBc z,k \xBe,$ then by prerequisite $rg(x) \xck
rg(y),$
so $rg(x)<rg(z),$
so $ \xBc x,i \xBe  \xer  \xBc z,k \xBe $ by definition. Likewise for
$rg(x)<rg(y)$ and
$ \xBc y,j \xBe  \xeb  \xBc z,k \xBe.$

$ \xcz $
\\[3ex]

All that remains to show then is that our constructions of smooth and of
smooth and transitive structures satisfy the condition

$ \xBc x,i \xBe  \xeb  \xBc y,j \xBe $ implies $rg(x) \xck rg(y).$
\index{Proposition A-Smooth-Complete}

\bp

$\hspace{0.01em}$

% (+++ Orig. No.:  Proposition A-Smooth-Complete +++)

\label{Proposition A-Smooth-Complete}

Let $ \xda $ be given.

Let - for simplicity - $ \xdy $ be closed under finite unions, and $ \xbm
: \xdy \xcp \xdp (Z).$
Then there is a $ \xdy -$smooth $ \xda -$ranked preferential structure $
\xdz,$ s.t. for all
$X \xbe \xdy $ $ \xbm (X)= \xbm_{ \xdz }(X)$ iff $ \xbm $ satisfies $(
\xbm \xcc ),$ $( \xbm PR),$ $( \xbm CUM),$ $( \xbm \xda ).$
\index{Proposition A-Smooth-Complete Proof}

\ep

\subparagraph{
Proof
}

$\hspace{0.01em}$

% (+++ Orig.:  Proof +++)

Consider the construction in the proof of Proposition \ref{Proposition
Smooth-Complete} (page \pageref{Proposition Smooth-Complete}).
We have to show that it respects the rang order with respect
to $ \xda,$ i.e. that $ \xBc x', \xbs '  \xBe  \xeb '  \xBc x, \xbs  \xBe $
implies $rg(x' )
\xck rg(x).$ This is easy: By
definition, $x' \xbe \xcV \{ran( \xbs_{i}):i \xbe \xbo \}.$
If $x' \xbe ran( \xbs_{0}),$ then for some $Y$ $x' \xbe \xbm (Y),$ $x \xbe
Y- \xbm (Y),$ so $rg(x' ) \xck rg(x)$
by $( \xbm \xda ).$
If $x' \xbe ran( \xbs_{i}),$ $i>0,$ then for some $X$ $x',x \xbe \xbm
(X),$ so $rg(x)=rg(x' )$ by $( \xbm \xda ).$

$ \xcz $ (Proposition \ref{Proposition A-Smooth-Complete} (page
\pageref{Proposition A-Smooth-Complete}) )
\\[3ex]
\index{Proposition A-Smooth-Trans-Complete}

\bp

$\hspace{0.01em}$

% (+++ Orig. No.:  Proposition A-Smooth-Trans-Complete +++)

\label{Proposition A-Smooth-Trans-Complete}

Let $ \xda $ be given.

Let - for simplicity - $ \xdy $ be closed under finite unions, and $ \xbm
: \xdy \xcp \xdp (Z).$
Then there is a $ \xdy -$smooth $ \xda -$ranked transitive preferential
structure $ \xdz,$ s.t.
for all
$X \xbe \xdy $ $ \xbm (X)= \xbm_{ \xdz }(X)$ iff $ \xbm $ satisfies $(
\xbm \xcc ),$ $( \xbm PR),$ $( \xbm CUM),$ $( \xbm \xda ).$
\index{Proposition A-Smooth-Trans-Complete Proof}

\ep

\subparagraph{
Proof
}

$\hspace{0.01em}$

% (+++ Orig.:  Proof +++)

Consider the construction in the proof of
Proposition \ref{Proposition D-4.5.3} (page \pageref{Proposition D-4.5.3}).

Thus, we only have to show that in $ \xdz $ defined by

$ \xdz $ $:=$ $ \xBc $ $\{ \xBc x,t_{x} \xBe :$ $x \xbe Z,$ $t_{x} \xbe
T_{x}\}$,
$ \xBc x,t_{x} \xBe  \xee  \xBc y,t_{y} \xBe $ iff $t_{x} \xem^{*}t_{y}$ $
\xBe,$
$t_{x} \xem t_{y}$ implies
$rg(y) \xck rg(x).$

But by construction of the trees,
$x_{n} \xbe Y_{n+1},$ and $x_{n+1} \xbe \xbm (U_{n} \xcv Y_{n+1}),$ so
$rg(x_{n+1}) \xck rg(x_{n}).$

$ \xcz $ (Proposition \ref{Proposition A-Smooth-Trans-Complete} (page
\pageref{Proposition A-Smooth-Trans-Complete}) )
\\[3ex]
\subsubsection{
The logical properties with definability preservation
}
\label{Section Hier-ARepr-Logic}
\label{Section Hier-ARepr-Logic-Without-Proof}
\index{Section Hier-ARepr-Logic-Without-Proof}

First, a small fact about the $ \xda.$

\bfa

$\hspace{0.01em}$

% (+++ Orig. No.:  Fact 6.21 +++)

\label{Fact 6.21}

Let $ \xda $ be as above (and thus finite).
Then each $A_{i}$ is equivalent to a formula $ \xba_{i}.$

\efa

\subparagraph{
Proof
}

$\hspace{0.01em}$

% (+++ Orig.:  Proof +++)

We use the standard topology and its compactness.
By definition, each $M(A_{i})$ is closed, by finiteness all unions of such
$M(A_{i})$ are
closed, too, so $ \xdC (M(A_{i}))$ is closed. By compactness, each open
cover $X_{j}:j \xbe J$
of the clopen $M(A_{i})$ contains a finite subcover, so also $ \xcV
\{M(A_{j}):j \xEd i\}$ has
a finite open cover. But the $M( \xbf ),$ $ \xbf $ a formula form a basis
of the closed
sets, so we are done. $ \xcz $
\\[3ex]

\bp

$\hspace{0.01em}$

% (+++ Orig. No.:  Proposition 6.22 +++)

\label{Proposition 6.22}

Let $ \xcn $ be a logic for $ \xdl.$ Set $T^{ \xdm }:=Th( \xbm_{ \xdm
}(M(T))),$
where $ \xdm $ is a preferential structure.

(1) Then there is a (transitive) definability preserving
classical preferential model $ \xdm $ s.t. $ \ol{ \ol{T} }=T^{ \xdm }$ iff

(LLE), (CCL), (SC), (PR) hold for all $T,T' \xcc \xdl.$

(2) The structure can be chosen smooth, iff, in addition

(CUM) holds.

(3) The structure can be chosen $ \xda -$ranked, iff, in addition

$( \xda -$min) $T \xcL \xCN \xba_{i}$ and $T \xcL \xCN \xba_{j},$ $i<j$
implies $ \ol{ \ol{T} } \xcl \xCN \xba_{j}$

holds.

\ep

The proof is an immediate consequence of Proposition \ref{Proposition 6.23}
(page \pageref{Proposition 6.23})
and the respective above results.
This proposition (or its analogue) was mostly already
shown in  \cite{Sch92} and
 \cite{Sch96-1} and is repeated here for completeness' sake, but
with a new and partly stronger proof.

\bp

$\hspace{0.01em}$

% (+++ Orig. No.:  Proposition 6.23 +++)

\label{Proposition 6.23}

Consider for a logic $ \xcn $ on $ \xdl $ the properties

(LLE), (CCL), (SC), (PR), (CUM), $( \xda -$min) hold for all $T,T' \xcc
\xdl.$

and for a function $ \xbm: \xdD_{ \xdl } \xcp \xdp (M_{ \xdl })$ the
properties

$( \xbm dp)$ $ \xbm $ is definability preserving, i.e. $ \xbm (M(T))=M(T'
)$ for some $T' $

$( \xbm \xcc ),$ $( \xbm PR),$ $( \xbm CUM),$ $( \xbm \xda )$

for all $X,Y \xbe \xdD_{ \xdl }.$

It then holds:

(a) If $ \xbm $ satisfies $( \xbm dp),$ $( \xbm \xcc ),$ $( \xbm PR),$
then $ \xcn $ defined by $ \ol{ \ol{T} }:=T^{ \xbm }:=$
$Th( \xbm (M(T)))$ satisfies (LLE), (CCL), (SC), (PR).
If $ \xbm $ satisfies in addition $( \xbm CUM),$ then (CUM) will hold,
too.
If $ \xbm $ satisfies in addition $( \xbm \xda ),$ then $( \xda -$min)
will hold, too.

(b) If $ \xcn $ satisfies (LLE), (CCL), (SC), (PR),
then there is $ \xbm: \xdD_{ \xdl } \xcp \xdp (M_{ \xdl })$ s.t. $ \ol{
\ol{T} }=T^{ \xbm }$
for all $T \xcc \xdl $ and $ \xbm $ satisfies $( \xbm dp),$ $( \xbm \xcc
),$ $( \xbm PR).$
If, in addition, (CUM) holds, then $( \xbm CUM)$ will hold, too.
If, in addition, $( \xda -$min) holds, then $( \xbm \xda )$ will hold,
too.

\ep

\subparagraph{
Proof
}

$\hspace{0.01em}$

% (+++ Orig.:  Proof +++)

All properties except $( \xda -$min) and $( \xbm \xda )$ are shown in
Proposition \ref{Proposition Alg-Log} (page \pageref{Proposition Alg-Log}).
But the remaining two are trivial. $ \xcz $
\\[3ex]
\section{Two sequent calculi}
\label{Section 2.2.9}
%  CHAPTER (6):  COMPARISON TO EXISTING PROOFS AND SPECIAL CASES
%  CHAPTER (6):  COMPARISON TO EXISTING PROOFS AND SPECIAL CASES
% %
% ===============================================================
\subsection{Introduction}
\label{Section 2.2.9.1}
%  Section (6.1):  Introduction
%  Section (6.1):  Introduction
% %
% ==============================

This section serves mainly as a posteriori motivation for our examination
of weak closure conditions of the domain. The second author realized first
when
looking at Lehmann's plausibility logic, that absence of $( \xcv )$ might
be a
problem for representation - see  \cite{Sch96-3} or  \cite{Sch04}.

Beyond motivation, the reader will see here two ``real life'' examples where
closure under $( \xcv )$ is not given, and thus problems arise. So this is
also
a warning against a too naive treatment of representation problems,
neglecting
domain closure issues.
\subsection{Plausibility Logic}
\label{Section 2.2.9.2}
%  Subsection (6.2.1):  Plausibility Logic
%  Subsection (6.2.1):  Plausibility Logic
% %
% =========================================

\paragraph{
Discussion of plausibility logic
}

$\hspace{0.01em}$

% (+++ Orig.:  Discussion of plausibility logic +++)

\label{Section Discussion of plausibility logic}

Plausibility logic was introduced by $D.$ Lehmann  \cite{Leh92a},
 \cite{Leh92b}
as a sequent
calculus in a propositional language without connectives. Thus, a
plausibility
logic language $ \xdl $ is just a set, whose elements correspond to
propositional
variables, and a sequent has the form $X \xcn Y,$ where $X,$ $Y$ are $
\ul{finite}$ subsets of
$ \xdl,$ thus, in the intuitive reading, $ \xcU X \xcn \xcO Y.$ (We use $
\xcn $ instead of the $ \xcl $ used
in  \cite{Leh92a},  \cite{Leh92b} and continue to reserve $
\xcl $ for classical logic.)

\paragraph{
The details:
}

$\hspace{0.01em}$

% (+++ Orig.:  The details: +++)

\label{Section The details:}

\bn

$\hspace{0.01em}$

% (+++ Orig. No.:  Notation D-6.2.1 +++)

\label{Notation D-6.2.1}

We abuse notation, and write $X \xcn a$ for $X \xcn \{a\},$ $X,a \xcn Y$
for $X \xcv \{a\} \xcn Y,$ $ab \xcn Y$ for
$\{a,b\} \xcn Y,$ etc. When discussing plausibility logic, $X,Y,$ etc.
will denote finite
subsets of $ \xdl,$ $a,b,$ etc. elements of $ \xdl.$

We first define the logical properties we will examine.

\en

\bd

$\hspace{0.01em}$

% (+++ Orig. No.:  Definition D-6.2.1 +++)

\label{Definition D-6.2.1}

$X$ and $Y$ will be finite subsets of $ \xdl,$ a, etc. elements of $ \xdl
.$
The base axiom and rules of plausibility logic are
(we use the prefix ``Pl'' to differentiate them from the usual ones):

\index{$(PlI)$}
\index{$(PlRM)$}
\index{$(PlCLM)$}
\index{$(PlCC)$}
\index{$(PlUCC)$}
\index{$PL$}
(PlI) (Inclusion): $X \xcn a$ for all $a \xbe X,$

(PlRM) (Right Monotony): $X \xcn Y$ $ \xch $ $X \xcn a,Y,$

(PlCLM) (Cautious Left Monotony): $X \xcn a,$ $X \xcn Y$ $ \xch $ $X,a
\xcn Y,$

(PlCC) (Cautious Cut): $X,a_{1} \Xl a_{n} \xcn Y,$ and for all $1 \xck i
\xck n$ $X \xcn a_{i},Y$ $ \xch $ $X \xcn Y,$

and as a special case of (PlCC):

(PlUCC) (Unit Cautious Cut): $X,a \xcn Y$, $X \xcn a,Y$ $ \xch $ $X \xcn
Y.$

and we denote by PL, for plausibility logic, the full system, i.e.
$(PlI)+(PlRM)+(PlCLM)+(PlCC).$ $ \xcz $
\\[3ex]

\ed

We now adapt the definition of a preferential model to plausibility logic.
This is the central definition on the semantic side.

\bd

$\hspace{0.01em}$

% (+++ Orig. No.:  Definition D-6.2.2 +++)

\label{Definition D-6.2.2}

Fix a plausibility logic language $ \xdl.$ A model for $ \xdl $ is then
just an arbitrary
subset of $ \xdl.$

If $ \xdm:= \xBc M, \xeb  \xBe $ is a preferential model s.t. $M$ is a set of
(indexed) $ \xdl -$models,
then for a finite set $X \xcc \xdl $ (to be imagined on the left hand side
of $ \xcn $!), we
define

(a) $m \xcm X$ iff $X \xcc m$

(b) $M(X)$ $:=$ $\{m$: $ \xBc m,i \xBe  \xbe M$ for some $i$ and $m \xcm X\}$

(c) $ \xbm (X)$ $:=$ $\{m \xbe M(X)$:
$ \xcE  \xBc m,i \xBe  \xbe M. \xCN \xcE  \xBc m',i'  \xBe  \xbe M$ $(m' \xbe
M(X)$
$ \xcu $ $ \xBc m',i'  \xBe  \xeb  \xBc m,i \xBe )\}$

(d) $X \xcm_{ \xdm }Y$ iff $ \xcA m \xbe \xbm (X).m \xcs Y \xEd \xCQ.$ $
\xcz $
\\[3ex]

\ed

(a) reflects the intuitive reading of $X$ as $ \xcU X,$ and (d) that of
$Y$ as $ \xcO Y$ in
$X \xcn Y.$ Note that $X$ is a set of ``formulas'', and $ \xbm (X)= \xbm_{
\xdm }(M(X)).$

We note as trivial consequences of the definition.

\bfa

$\hspace{0.01em}$

% (+++ Orig. No.:  Fact D-6.2.1 +++)

\label{Fact D-6.2.1}

(a) $a \xcm_{ \xdm }b$ iff for all $m \xbe \xbm (a).b \xbe m$

(b) $X \xcm_{ \xdm }Y$ iff $ \xbm (X) \xcc \xcV \{M(b):b \xbe Y\}$

(c) $m \xbe \xbm (X)$ $ \xcu $ $X \xcc X' $ $ \xcu $ $m \xbe M(X' )$ $
\xch $ $m \xbe \xbm (X' )$. $ \xcz $
\\[3ex]

\efa

We note without proof: $(PlI)+(PlRM)+(PlCC)$ is complete (and sound) for
preferential models - see  \cite{Sch96-3} or  \cite{Sch04} for
a proof.

We note the following fact for smooth preferential models:

\bfa

$\hspace{0.01em}$

% (+++ Orig. No.:  Fact D-6.2.2 +++)

\label{Fact D-6.2.2}

Let $ \xCf U,X,Y$ be any sets, $ \xdm $ be smooth for at least $\{Y,X\}$
and
let $ \xbm (Y) \xcc U \xcv X,$ $ \xbm (X) \xcc U,$ then $X \xcs Y \xcs
\xbm (U) \xcc \xbm (Y).$ (This is, of course,
a special case of $( \xbm Cum1).$

\efa

\be

$\hspace{0.01em}$

% (+++ Orig. No.:  Example D-6.2.1 +++)

\label{Example D-6.2.1}

Let $ \xdl:=\{a,b,c,d,e,f\},$ and
$ \xdx $ $:=$ $\{a \xcn b$, $b \xcn a$, $a \xcn c$, $a \xcn fd$, $dc
\xcn ba$, $dc \xcn e$, $fcba \xcn e\}.$
We show that $ \xdx $ does not have a smooth representation.

\ee

\bfa

$\hspace{0.01em}$

% (+++ Orig. No.:  Fact D-6.2.3 +++)

\label{Fact D-6.2.3}

$ \xdx $ does not entail $a \xcn e.$

\efa

See  \cite{Sch96-3} or  \cite{Sch04} for a proof.

Suppose now that there is a smooth preferential model
$ \xdm = \xBc M, \xeb  \xBe $ for plausibility
logic which represents $ \xcn,$ i.e. for all $ \xCf X,Y$ finite subsets
of $ \xdl $ $X \xcn Y$ iff
$X \xcm_{ \xdm }Y.$ (See
Definition \ref{Definition D-6.2.2} (page \pageref{Definition D-6.2.2})  and
Fact \ref{Fact D-6.2.1} (page \pageref{Fact D-6.2.1}).)

$a \xcn a,$ $a \xcn b,$ $a \xcn c$ implies for $m \xbe \xbm (a)$ $a,b,c
\xbe m.$ Moreover, as $a \xcn df,$ then also
$d \xbe m$ or $f \xbe m.$ As $a \xcN e,$ there must be $m \xbe \xbm (a)$
s.t. $e \xce m.$ Suppose now $m \xbe \xbm (a)$
with $f \xbe m.$ So $a,b,c,f \xbe m,$ thus by $m \xbe \xbm (a)$ and
Fact \ref{Fact D-6.2.1} (page \pageref{Fact D-6.2.1}),
$m \xbe \xbm (a,b,c,f).$ But
$fcba \xcn e,$ so $e \xbe m.$ We thus have shown that $m \xbe \xbm (a)$
and $f \xbe m$ implies $e \xbe m.$
Consequently, there must be $m \xbe \xbm (a)$ s.t. $d \xbe m,$ $e \xce m.$
Thus, in particular, as $cd \xcn e,$ there is $m \xbe \xbm (a),$ $a,b,c,d
\xbe m,$ $m \xce \xbm (cd).$
But by $cd \xcn ab,$ and $b \xcn a,$ $ \xbm (cd) \xcc M(a) \xcv M(b)$ and
$ \xbm (b) \xcc M(a)$ by
Fact \ref{Fact D-6.2.1} (page \pageref{Fact D-6.2.1}).
Let now $T:=M(cd),$ $R:=M(a),$ $S:=M(b),$ and $ \xbm_{ \xdm }$ be the
choice function of the
minimal elements in the structure $ \xdm,$ we then have by $ \xbm (S)=
\xbm_{ \xdm }(M(S))$:

1. $ \xbm_{ \xdm }(T) \xcc R \xcv S,$

2. $ \xbm_{ \xdm }(S) \xcc R,$

3. there is $m \xbe S \xcs T \xcs \xbm_{ \xdm }(R),$ but $m \xce \xbm_{
\xdm }(T),$

but this contradicts above Fact \ref{Fact D-6.2.2} (page \pageref{Fact D-6.2.2})
.

$ \xcz $ (Counterexample $D-6.2.1)$
\\[3ex]
\subsection{A comment on the work by Arieli and Avron}
\label{Section 2.2.9.3}
%  Subsection (6.2.2):  Arieli/Avron
%  Subsection (6.2.2):  Arieli/Avron
% %
% ===================================

We turn to a similar case, published in  \cite{AA00}.
Definitions are due to  \cite{AA00}, for motivation the reader is
referred there.

\bd

$\hspace{0.01em}$

% (+++ Orig. No.:  Definition D-6.2.3 +++)

\label{Definition D-6.2.3}

(1) A Scott consequence relation, abbreviated scr, is a binary relation $
\xcl $
between sets of formulae, that satisfies the following conditions:

(s-R) if $ \xbG \xcs \xbD \xEd \xCQ,$ the $ \xbG \xcl \xbD,$

(M) if $ \xbG \xcl \xbD $ and $ \xbG \xcc \xbG ',$ $ \xbD \xcc \xbD ',$
then $ \xbG ' \xcl \xbD ',$

(C) if $ \xbG \xcl \xbq, \xbD $ and $ \xbG ', \xbq \xcl \xbD ',$ then $
\xbG, \xbG ' \xcl \xbD, \xbD '.$

(2) A Scott cautious consequence relation, abbreviated sccr, is a binary
relation $ \xcn $ between nonempty sets of formulae, that satisfies the
following
conditions:

(s-R) if $ \xbG \xcs \xbD \xEd \xCQ,$ the $ \xbG \xcn \xbD,$

(CM) if $ \xbG \xcn \xbD $ and $ \xbG \xcn \xbq,$ then $ \xbG, \xbq \xcn
\xbD,$

(CC) if $ \xbG \xcn \xbq $ and $ \xbG, \xbq \xcn \xbD,$ then $ \xbG \xcn
\xbD.$

\ed

\be

$\hspace{0.01em}$

% (+++ Orig. No.:  Example D-6.2.2 +++)

\label{Example D-6.2.2}

We have two consequence relations, $ \xcl $ and $ \xcn.$

The rules to consider are

$LCC^{n}$ $ \frac{ \xbG \xcn \xbq_{1}, \xbD  \Xl  \xbG \xcn \xbq_{n}, \xbD
\xbG, \xbq_{1}, \Xl, \xbq_{n} \xcn \xbD }{ \xbG \xcn \xbD }$

$RW^{n}$ $ \frac{ \xbG \xcn \xbq_{i}, \xbD i=1 \Xl n \xbG, \xbq_{1}, \Xl
, \xbq_{n} \xcl \xbf }{ \xbG \xcn \xbf, \xbD }$

Cum $ \xbG, \xbD \xEd \xCQ,$ $ \xbG \xcl \xbD $ $ \xch $ $ \xbG \xcn
\xbD $

RM $ \xbG \xcn \xbD $ $ \xch $ $ \xbG \xcn \xbq, \xbD $

CM $ \frac{ \xbG \xcn \xbq \xbG \xcn \xbD }{ \xbG, \xbq \xcn \xbD }$

$s-R$ $ \xbG \xcs \xbD \xEd \xCQ $ $ \xch $ $ \xbG \xcn \xbD $

$M$ $ \xbG \xcl \xbD,$ $ \xbG \xcc \xbG ',$ $ \xbD \xcc \xbD ' $ $ \xch
$ $ \xbG ' \xcl \xbD ' $

$C$ $ \frac{ \xbG_{1} \xcl \xbq, \xbD_{1} \xbG_{2}, \xbq \xcl \xbD_{2}}{
\xbG_{1}, \xbG_{2} \xcl \xbD_{1}, \xbD_{2}}$

Let $ \xdl $ be any set.
Define now $ \xbG \xcl \xbD $ iff $ \xbG \xcs \xbD \xEd \xCQ.$
Then $s-R$ and $M$ for $ \xcl $ are trivial. For $C:$ If $ \xbG_{1} \xcs
\xbD_{1} \xEd \xCQ $ or $ \xbG_{1} \xcs \xbD_{1} \xEd \xCQ,$ the
result is trivial. If not, $ \xbq \xbe \xbG_{1}$ and $ \xbq \xbe
\xbD_{2},$ which implies the result.
So $ \xcl $ is a scr.

Consider now the rules for a sccr which is $ \xcl -$plausible for this $
\xcl.$
Cum is equivalent to $s-$R, which is essentially (PlI) of Plausibility
Logic.
Consider $RW^{n}.$ If $ \xbf $ is one of the $ \xbq_{i},$ then the
consequence $ \xbG \xcn \xbf, \xbD $ is a case
of one of the other hypotheses. If not, $ \xbf \xbe \xbG,$ so $ \xbG \xcn
\xbf $ by $s-$R, so $ \xbG \xcn \xbf, \xbD $
by RM (if $ \xbD $ is finite). So, for this $ \xcl,$ $RW^{n}$ is a
consequence of $s-R$ $+$ RM.

We are left with $LCC^{n},$ RM, CM, $s-$R, it was shown in  \cite{Sch04}
and  \cite{Sch96-3}
that this does not suffice to guarantee smooth representability, by
failure of
$( \xbm Cum1).$
\section{Blurred observation - absence of definability preservation}
\label{Section 2.2.10}
%  CHAPTER (7):  BLURRED OBSERVATION - DEFINABILITY PRESERVATION
%  CHAPTER (7):  BLURRED OBSERVATION - DEFINABILITY PRESERVATION
% %
% ===============================================================
\subsection{Introduction}
\label{Section 2.2.10.1}
%  Section (7.1):  Introduction
%  Section (7.1):  Introduction
% %
% ==============================

\ee

Lack of definability preservation results in uncertainty. We do not know
exactly the result, but only that it is not too far away from what we
(seem to)
observe.

Thus, we pretend to know more than we know, and, according to our general
policy of not neglecting ignorance, we should branch here into a multitude
of solutions.

We take here a seemingly different way, but, as a matter of fact, we just
describe the boundaries of what is permitted. So, everything which lies in
those
boundaries, is a possible solution, and every such solution should be
considered
as equal, and, again, we should not pretend to know more than we actually
do.
\subsubsection{General remarks, affected conditions}

$\hspace{0.01em}$
%  Subsection (7.1.1)  General remarks, affected conditions
%  Subsection (7.1.1)  General remarks, affected conditions
% %
% =========================================================

We assume now - unless explicitly stated otherwise - $ \xdy \xcc \xdp (Z)$
to be closed
under arbitrary intersections (this is
used for the definition of $ \wt{.})$ and finite unions, and $ \xCQ,Z
\xbe \xdy.$ This holds,
of course, for $ \xdy = \xdD_{ \xdl },$ $ \xdl $ any propositional
language.

The aim of Section \ref{Section 2.2.10} (page \pageref{Section 2.2.10})
is to present
the results of  \cite{Sch04} connected to
problems of definability preservation in a uniform way, stressing the
crucial
condition $ \wt{X} \xcs \wt{Y}= \wt{X \xcs Y}.$ This presentation shall
help and guide future research
concerning similar problems.

For motivation, we first consider the problem with definability
preservation for
the rules

(PR) $ \ol{ \ol{T \xcv T' } }$ $ \xcc $ $ \ol{ \ol{ \ol{T} } \xcv T' },$
and

$( \xcn =)$ $T \xcl T',$ $Con( \ol{ \ol{T' } },T)$ $ \xch $ $ \ol{ \ol{T}
}$ $=$ $ \ol{ \ol{ \ol{T' } } \xcv T}$ holds.

which are consequences of

$( \xbm PR)$ $X \xcc Y$ $ \xch $ $ \xbm (Y) \xcs X \xcc \xbm (X)$ or

$( \xbm =)$ $X \xcc Y,$ $ \xbm (Y) \xcs X \xEd \xCQ $ $ \xch $ $ \xbm (Y)
\xcs X= \xbm (X)$ respectively

and definability preservation.

Example \ref{Example Pref-Dp} (page \pageref{Example Pref-Dp})
showed that in the general case without
definability preservation, (PR) fails, and
the following
Example \ref{Example D-7.1.1} (page \pageref{Example D-7.1.1})
shows that
in the ranked case, $( \xcn =)$ may fail. So failure is not just a
consequence of
the very liberal definition of general preferential structures.

\be

$\hspace{0.01em}$

% (+++ Orig. No.:  Example D-7.1.1 +++)

\label{Example D-7.1.1}

Take $\{p_{i}:i \xbe \xbo \}$ and put $m:=m_{ \xcU p_{i}},$ the model
which makes all $p_{i}$ true, in the top
layer, all the other in the bottom layer. Let $m' \xEd m,$ $T':= \xCQ,$
$T:=Th(m,m' ).$ Then
Then $ \ol{ \ol{T' } }=T',$ so $Con( \ol{ \ol{T' } },T),$ $ \ol{ \ol{T}
}=Th(m' ),$ $ \ol{ \ol{ \ol{T' } } \xcv T}=T.$

$ \xcz $
\\[3ex]

\ee

We remind the reader of
Definition \ref{Definition Def-Clos} (page \pageref{Definition Def-Clos})  and
Fact \ref{Fact Def-Clos} (page \pageref{Fact Def-Clos}), partly taken
from  \cite{Sch04}.

We turn to the central condition.
\subsubsection{The central condition}

$\hspace{0.01em}$
%                                            ~  ~ ~~~~
%                                            ~  ~ ~~~~
%  Subsection (7.1.2) The central condition (X%sY=X%sY)
% %
% =====================================================

We analyze the problem of (PR), seen in
Example \ref{Example D-7.1.1} (page \pageref{Example D-7.1.1})  (1) above,
working in the intended application.

(PR) is equivalent to $M( \ol{ \ol{T} } \xcv T' )$ $ \xcc $ $M( \ol{ \ol{T
\xcv T' } }).$ To show (PR) from $( \xbm PR),$ we argue
as follows, the crucial point is marked by ``?'':

$M( \ol{ \ol{T \xcv T' } })$ $=$ $M(Th( \xbm (M_{T \xcv T' })))$ $=$ $
\wt{ \xbm (M_{T \xcv T' })}$ $ \xcd $ $ \xbm (M_{T \xcv T' })$ $=$ $ \xbm
(M_{T} \xcs M_{T' })$
$ \xcd $ (by $( \xbm PR))$ $ \xbm (M_{T}) \xcs M_{T' }$? $ \wt{ \xbm
(M_{T})} \xcs M_{T' }$ $=$ $M(Th( \xbm (M_{T}))) \xcs M_{T' }$ $=$ $M(
\ol{ \ol{T} }) \xcs M_{T' }$ $=$
$M( \ol{ \ol{T} } \xcv T' ).$ If $ \xbm $ is definability preserving, then
$ \xbm (M_{T})$ $=$ $ \wt{ \xbm (M_{T})},$ so ``?'' above
is equality, and everything is fine. In general, however, we have only
$ \xbm (M_{T})$ $ \xcc $ $ \wt{ \xbm (M_{T})},$ and the argument
collapses.

But it is not necessary to
impose $ \xbm (M_{T})$ $=$ $ \wt{ \xbm (M_{T})},$ as we still have room to
move: $ \wt{ \xbm (M_{T \xcv T' })}$ $ \xcd $ $ \xbm (M_{T \xcv T' }).$
(We do not consider here $ \xbm (M_{T} \xcs M_{T' })$ $ \xcd $ $ \xbm
(M_{T}) \xcs M_{T' }$ as room to move,
as we are now interested
only in questions related to definability preservation.) If we had
$ \wt{ \xbm (M_{T})} \xcs M_{T' }$ $ \xcc $ $ \wt{ \xbm (M_{T}) \xcs M_{T'
}}$ $,$
we could use $ \xbm (M_{T}) \xcs M_{T' }$ $ \xcc $ $ \xbm (M_{T} \xcs
M_{T' })$ $=$
$ \xbm (M_{T \xcv T' })$ and monotony of $ \wt{.}$ to obtain $ \wt{ \xbm
(M_{T})} \xcs M_{T' }$ $ \xcc $ $ \wt{ \xbm (M_{T}) \xcs M_{T' }}$ $ \xcc
$
$ \wt{ \xbm (M_{T} \xcs M_{T' })}$ $=$ $ \wt{ \xbm (M_{T \xcv T' })}.$
If, for instance, $T' =\{ \xbq \},$ we have $ \wt{ \xbm (M_{T})} \xcs
M_{T' }$ $=$ $ \wt{ \xbm (M_{T}) \xcs M_{T' }}$
by Fact \ref{Fact Def-Clos} (page \pageref{Fact Def-Clos})  $(Cl \xcs +).$
Thus, definability preservation

is not the only solution to the problem.

We have seen in Fact \ref{Fact Def-Clos} (page \pageref{Fact Def-Clos})
that $ \wt{X \xcv Y}$ $=$ $ \wt{X} \xcv \wt{Y},$ moreover $ \xCf X-Y$ $=$
$X \xcs \xdC Y$ $( \xdC Y$ the set
complement of $Y),$ so, when considering boolean expressions of model sets
(as we do in usual properties describing logics), the central question is
whether

$( \xCq \xcs )$ $ \wt{X \xcs Y}$ $=$ $ \wt{X} \xcs \wt{Y}$

holds.

We take a closer look at this question.

$ \wt{X \xcs Y}$ $ \xcc $ $ \wt{X} \xcs \wt{Y}$
holds by Fact \ref{Fact Def-Clos} (page \pageref{Fact Def-Clos})  (6).
Using $(Cl \xcv )$ and monotony
of $ \wt{.}$, we have $ \wt{X} \xcs \wt{Y}$ $=$ $ \wt{((X \xcs Y) \xcv
(X-Y))}$ $ \xcs $ $ \wt{((X \xcs Y) \xcv (Y-X))}$ $=$
$( \wt{(X \xcs Y)} \xcv \wt{(X-Y)})$ $ \xcs $ $( \wt{(X \xcs Y)} \xcv
\wt{(Y-X)})$ $=$
$ \wt{X \xcs Y}$ $ \xcv $ $( \wt{X-Y} \xcs \wt{Y-X}),$
thus $ \wt{X} \xcs \wt{Y}$ $ \xcc $ $ \wt{X \xcs Y}$ iff

$( \xCq \xcs ' )$ $ \wt{Y-X}$ $ \xcs $ $ \wt{X-Y}$ $ \xcc $ $ \wt{X \xcs
Y}$

holds.

Intuitively speaking, the condition holds
iff we cannot approximate any element both from $X-Y$ and $X-$Y, which
cannot be
approximated from $X \xcs Y,$ too.

Note that in above Example \ref{Example D-7.1.1} (page \pageref{Example
D-7.1.1})  (1)
$X:= \xbm (M_{T})=M_{ \xdl }-\{n' \},$ $Y:=M_{T' }=\{n,n' \},$
$ \wt{X-Y}=M_{ \xdl },$ $ \wt{Y-X}=\{n' \},$ $ \wt{X \xcs Y}=\{n\},$
and $ \wt{X} \xcs \wt{Y}=\{n,n' \}.$

We consider now particular cases:

(1) If $X \xcs Y= \xCQ,$ then by $ \xCQ \xbe \xdy,$ $( \xCq \xcs )$
holds iff $ \wt{X} \xcs \wt{Y}= \xCQ.$

(2) If $X \xbe \xdy $ and $Y \xbe \xdy,$ then $ \wt{X-Y} \xcc X$
and $ \wt{Y-X} \xcc Y,$
so $ \wt{X-Y} \xcs \wt{Y-X} \xcc X \xcs Y \xcc \wt{X \xcs Y}$
and $( \xCq \xcs )$ trivially holds.

(3) $X \xbe \xdy $ and $ \xdC X \xbe \xdy $ together also suffice - in
these cases
$ \wt{Y-X}$ $ \xcs $ $ \wt{X-Y}$ $=$ $ \xCQ:$ $ \wt{Y-X}= \wt{Y \xcs \xdC
X} \xcc \xdC X,$
and $ \wt{X-Y} \xcc X,$
so $ \wt{Y-X} \xcs \wt{X-Y} \xcc X \xcs \xdC X= \xCQ \xcc \wt{X \xcs Y}.$
(The same holds, of course, for $Y.)$
(In the intended application, such $X$ will be $M( \xbf )$ for some
formula $ \xbf.$ But, a
warning, $ \xbm (M( \xbf ))$ need not again be the $M( \xbq )$ for some $
\xbq.)$

We turn to the properties of various structures and apply our results.
\subsubsection{Application to various structures}

$\hspace{0.01em}$
%  Subsection (7.1.3) Application to various structures
%  Subsection (7.1.3) Application to various structures
% %
% =====================================================

We now take a look at other frequently used logical conditions. First, in
the
context on nonmonotonic logics, the following rules will always hold in
smooth
preferential structures, even if we consider full theories, and not
necessarily
definability preserving structures:

\bfa

$\hspace{0.01em}$

% (+++ Orig. No.:  Fact D-7.1.2 +++)

\label{Fact D-7.1.2}

Also for full theories, and not necessarily definability preserving
structures
hold:

(1) $ \xCf (LLE),$ $ \xCf (RW),$ $ \xCf (AND),$ $ \xCf (REF),$ by
definition and $( \xbm \xcc ),$

(2) $ \xCf (OR),$

(3) $ \xCf (CM)$ in smooth structures,

(4) the infinitary version of $ \xCf (CUM)$ in smooth structures.
In definability preserving structures, but also when considering only
formulas
hold:

(5) $ \xCf (PR),$

(6) $( \xcn =)$ in ranked structures.

\efa

\subparagraph{
Proof
}

$\hspace{0.01em}$

% (+++ Orig.:  Proof +++)

We use the corresponding algebraic properties. The result then follows
from
Proposition \ref{Proposition Alg-Log} (page \pageref{Proposition Alg-Log}). $
\xcz $
\\[3ex]

We turn to theory revision. The following definition and example, taken
from
 \cite{Sch04}
shows, that the usual AGM axioms for theory revision fail in distance
based
structures in the general case, unless we require definability
preservation.
See Chapter \ref{Chapter TR} (page \pageref{Chapter TR})
for discussion and motivation.

\bd

$\hspace{0.01em}$

% (+++ Orig. No.:  Definition D-7.1.2 +++)

\label{Definition D-7.1.2}

We summarize the AGM postulates $(K*7)$ and $(K*8)$ in $(*4):$

$(*4)$ If $T*T' $ is consistent with $T'',$ then $T*(T' \xcv T'' )$ $=$ $
\ol{(T*T' ) \xcv T'' }.$

\ed

\be

$\hspace{0.01em}$

% (+++ Orig. No.:  Example D-7.1.2 +++)

\label{Example D-7.1.2}

Consider an infinite propositional language $ \xdl.$

Let $X$ be an infinite set of models, $m,$ $m_{1},$ $m_{2}$ be models for
$ \xdl.$
Arrange the models of $ \xdl $ in the real plane s.t. all $x \xbe X$ have
the same
distance $<2$ (in the real plane) from $m,$ $m_{2}$ has distance 2 from
$m,$ and $m_{1}$ has
distance 3 from $m.$

Let $T,$ $T_{1},$ $T_{2}$ be complete (consistent) theories, $T' $ a
theory with infinitely
many models, $M(T)=\{m\},$ $M(T_{1})=\{m_{1}\},$ $M(T_{2})=\{m_{2}\}.$
$M(T' )=X \xcv \{m_{1},m_{2}\},$ $M(T'' )=\{m_{1},m_{2}\}.$

Assume $Th(X)=T',$ so $X$ will not be definable by a theory.

Then $M(T) \xfA M(T' )=X,$ but $T*T' =Th(X)=T'.$ So $T*T' $ is consistent
with $T'',$ and
$ \ol{(T*T' ) \xcv T'' }=T''.$ But $T' \xcv T'' =T'',$ and $T*(T' \xcv
T'' )=T_{2} \xEd T'',$ contradicting $(*4).$

$ \xcz $
\\[3ex]

\ee

We show now that the version with formulas only holds here, too, just as
does
above (PR), when we consider formulas only - this is needed below for $T''
$ only.
This was already shown in  \cite{Sch04}, we give now a proof based
on our new
principles.

\bfa

$\hspace{0.01em}$

% (+++ Orig. No.:  Fact D-7.1.3 +++)

\label{Fact D-7.1.3}

$(*4)$ holds when considering only formulas.

\efa

\subparagraph{
Proof
}

$\hspace{0.01em}$

% (+++ Orig.:  Proof +++)

When we fix the left hand side, the structure is ranked, so $Con(T*T',T''
)$
implies $(M_{T} \xfA M_{T' }) \xcs M_{T'' } \xEd \xCQ $ by $T'' =\{ \xbq
\}$ and thus $M_{T} \xfA M_{T' \xcv T'' }$ $=$ $M_{T} \xfA (M_{T' } \xcs
M_{T'' })$ $=$
$(M_{T} \xfA M_{T' }) \xcs M_{T'' }.$
So $M(T*(T' \xcv T'' ))$ $=$ $ \wt{M_{T} \xfA M_{T' \xcv T'' }}$ $=$ $
\wt{(M_{T} \xfA M_{T' }) \xcs M_{T'' }}$ $=$ (by $T'' =\{ \xbq \},$ see
above)
$ \wt{(M_{T} \xfA M_{T' })} \xcs \wt{M_{T'' }}$ $=$ $ \wt{(M_{T} \xfA
M_{T' })} \xcs M_{T'' }$ $=$ $M((T*T' ) \xcv T'' ),$
and $T*(T' \xcv T'' )= \ol{(T*T' ) \xcv T'' }.$
$ \xcz $
\\[3ex]
\subsection{
General and smooth structures without definability preservation
}
\label{Section 2.2.10.2}
\subsubsection{Introduction}

$\hspace{0.01em}$
%  Subsection (7.2.1):  General and smooth structures, Introduction
%  Subsection (7.2.1):  General and smooth structures, Introduction
% %
% ==================================================================

Note that in Sections 3.2 and 3.3 of  \cite{Sch04},
as well as in Proposition 4.2.2 of  \cite{Sch04} we
have characterized $ \xbm: \xdy \xcp \xdy $ or $ \xfA: \xdy \xCK \xdy
\xcp \xdy,$ but a closer inspection of the
proofs shows that the destination can as well be assumed $ \xdp (Z),$
consequently
we can simply re-use above algebraic representation results also for the
not
definability preserving case. (Note that the easy direction of all
these results work for destination $ \xdp (Z),$ too.) In particular, also
the
proof for the not definability preserving case of revision in  \cite{Sch04} can
be simplified - but we will not go into details here.

$( \xcv )$ and $( \xcS )$ are again assumed to hold now - we need $( \xcS
)$ for $ \wt{.}$.

The central functions and conditions to consider are summarized in the
following definition.

\bd

$\hspace{0.01em}$

% (+++ Orig. No.:  Definition D-7.2.1 +++)

\label{Definition D-7.2.1}

Let $ \xbm: \xdy \xcp \xdy,$ we define $ \xbm_{i}: \xdy \xcp \xdp (Z):$

$ \xbm_{0}(U)$ $:=$ $\{x \xbe U:$ $ \xCN \xcE Y \xbe \xdy (Y \xcc U$ and
$x \xbe Y- \xbm (Y))\},$

$ \xbm_{1}(U)$ $:=$ $\{x \xbe U:$ $ \xCN \xcE Y \xbe \xdy ( \xbm (Y) \xcc
U$ and $x \xbe Y- \xbm (Y))\},$

$ \xbm_{2}(U)$ $:=$ $\{x \xbe U:$ $ \xCN \xcE Y \xbe \xdy ( \xbm (U \xcv
Y) \xcc U$ and $x \xbe Y- \xbm (Y))\}$

(note that we use $( \xcv )$ here),

$ \xbm_{3}(U)$ $:=$ $\{x \xbe U:$ $ \xcA y \xbe U.x \xbe \xbm (\{x,y\})\}$

(we use here $( \xcv )$ and that singletons are in $ \xdy ).$

``Small'' is now in the sense of
Definition \ref{Definition Def-Clos} (page \pageref{Definition Def-Clos}).

$( \xbm PR0)$ $ \xbm (U)- \xbm_{0}(U)$ is small,

$( \xbm PR1)$ $ \xbm (U)- \xbm_{1}(U)$ is small,

$( \xbm PR2)$ $ \xbm (U)- \xbm_{2}(U)$ is small,

$( \xbm PR3)$ $ \xbm (U)- \xbm_{3}(U)$ is small.

$( \xbm PR0)$ with its function will be the one to consider for general
preferential
structures, $( \xbm PR2)$ the one for smooth structures.

\ed

\paragraph{
A non-trivial problem
}

$\hspace{0.01em}$

% (+++ Orig.:  A non-trivial problem +++)

\label{Section A non-trivial problem}

Unfortunately, we cannot use
$( \xbm PR0)$ in the smooth case, too, as
Example \ref{Example D-7.2.2} (page \pageref{Example D-7.2.2})  below will show.
This sheds some
doubt on the possibility to find an easy common approach to all cases of
not
definability preserving preferential, and perhaps other, structures. The
next
best guess, $( \xbm PR1)$ will not work either, as the same example shows
- or by
Fact \ref{Fact D-7.2.1} (page \pageref{Fact D-7.2.1})  (10),
if $ \xbm $ satisfies $( \xbm Cum),$ then $ \xbm_{0}(U)= \xbm_{1}(U).$ $(
\xbm PR3)$ and
$ \xbm_{3}$ are used for ranked structures.

We will now see that this first impression of a difficult situation is
indeed
well founded.

First, note that in our context, $ \xbm $ will not necessarily respect $(
\xbm PR).$ Thus,
if e.g. $x \xbe Y- \xbm (Y),$ and $ \xbm (Y) \xcc U,$ we cannot
necessarily conclude that
$x \xce \xbm (U \xcv Y)$ - the fact that $x$ is minimized in $U \xcv Y$
might be hidden by the
bigger $ \xbm (U \xcv Y).$

Consequently, we may have to work with small sets $( \xCf Y$ in the case
of $ \xbm_{2}$ above)
to see the problematic elements - recall that the smaller the set $ \xbm
(X)$ is,
the less it can ``hide'' missing elements - but will need bigger sets $(U
\xcv Y$ in
above example) to recognize the contradiction.

Second, ``problematic'' elements are those involved in a contradiction, i.e.
contradicting the representation conditions. Now, a negation of a
conjunction
is a disjunction of negations, so, generally, we will have to look at
various
possibilities of violated conditions. But the general situation is much
worse,
still.

\be

$\hspace{0.01em}$

% (+++ Orig. No.:  Example D-7.2.1 +++)

\label{Example D-7.2.1}

Look at the ranked case, and assume no closure properties of the domain.
Recall that we might be unable to see $ \xbm (X),$ but see only $ \wt{
\xbm (X)}.$
Suppose we have $ \wt{ \xbm (X_{1})} \xcs (X_{2}- \wt{ \xbm (X_{2})}) \xEd
\xCQ,$ $ \wt{ \xbm (X_{2})} \xcs (X_{3}- \wt{ \xbm (X_{3})}) \xEd \xCQ,$
$ \wt{ \xbm (X_{n-1})} \xcs (X_{n}- \wt{ \xbm (X_{n})}) \xEd \xCQ,$ $
\wt{ \xbm (X_{n})} \xcs (X_{1}- \wt{ \xbm (X_{1})}) \xEd \xCQ,$
which seems to be a contradiction. (It only is a real contradiction if it
still
holds without the closures.)
But, we do not know where the contradiction is situated. It might well be
that
for all but one $i$ really $ \xbm (X_{i}) \xcs (X_{i+1}- \xbm (X_{i+1}))
\xEd \xCQ,$ and not only that
for the closure $ \wt{ \xbm (X_{i})}$ of $ \xbm (X_{i})$ $ \wt{ \xbm
(X_{i})} \xcs (X_{i+1}- \wt{ \xbm (X_{i+1})}) \xEd \xCQ,$ but we might be
unable to find this out. So we have to branch into all possibilities, i.e.
for one, or several $i$ $ \wt{ \xbm (X_{i})} \xcs (X_{i+1}- \wt{ \xbm
(X_{i+1})}) \xEd \xCQ,$ but $ \xbm (X_{i}) \xcs (X_{i+1}- \xbm
(X_{i+1}))= \xCQ.$

$ \xcz $
\\[3ex]

\ee

The situation might even be worse, when those $ \wt{ \xbm (X_{i})} \xcs
(X_{i+1}- \wt{ \xbm (X_{i+1})}) \xEd \xCQ $ are
involved in several cycles, etc. Consequently, it seems very difficult to
describe all possible violations in one concise condition, and thus we
will
examine here only some specific cases, and do not pretend that they are
the
only ones, that other cases are similar, or that our solutions (which
depend
on closure conditions) are the best ones.

\paragraph{
Outline of our solutions in some particular cases
}

$\hspace{0.01em}$

% (+++ Orig.:  Outline of our solutions in some particular cases +++)

\label{Section Outline of our solutions in some particular cases}

The strategy of representation without definability preservation will in
all
cases be very simple: Under sufficient conditions,
among them smallness $( \xbm PRi)$ as described above, the corresponding
function
$ \xbm_{i}$ has all the properties to guarantee representation by a
corresponding
structures, and we can just take our representation theorems for the dp
case,
to show this. Using smallness again, we can show that we have obtained a
sufficient approximation - see
Proposition \ref{Proposition D-7.2.3} (page \pageref{Proposition D-7.2.3}),
Proposition \ref{Proposition D-7.2.4} (page \pageref{Proposition D-7.2.4}),
Proposition \ref{Proposition D-7.3.3} (page \pageref{Proposition D-7.3.3}).

We first show some properties for the $ \xbm_{i},$ $i=0,1,2.$ A
corresponding
result for $ \xbm_{3}$ is given in
Fact \ref{Fact D-7.3.1} (page \pageref{Fact D-7.3.1})  below.
(The conditions and results
are sufficiently different for $ \xbm_{3}$ to make a separation more
natural.)

Property (9) of the following
Fact \ref{Fact D-7.2.1} (page \pageref{Fact D-7.2.1})
fails for $ \xbm_{0}$ and $ \xbm_{1},$ as
Example \ref{Example D-7.2.2} (page \pageref{Example D-7.2.2})
below will show. We will therefore work in the smooth case
with $ \xbm_{2}.$
\subsubsection{Results}

$\hspace{0.01em}$
%  Subsection (7.2.2):  General and smooth structures, The results
%  Subsection (7.2.2):  General and smooth structures, The results
% %
% =================================================================

\bfa

$\hspace{0.01em}$

% (+++ Orig. No.:  Fact D-7.2.1 +++)

\label{Fact D-7.2.1}

(This is partly Fact 5.2.6 in  \cite{Sch04}.)

Recall that $ \xdy $ is closed under $( \xcv ),$ and $ \xbm: \xdy \xcp
\xdy.$
Let $ \xCf A,B,U,U',X,Y$ be elements of $ \xdy $ and the $ \xbm_{i}$ be
defined from $ \xbm $ as in
Definition \ref{Definition D-7.2.1} (page \pageref{Definition D-7.2.1}).
$i$ will here be $0,1,$ or 2, but not 3.

(1) Let $ \xbm $ satisfy $( \xbm \xcc ),$ then $ \xbm_{1}(X) \xcc
\xbm_{0}(X)$ and $ \xbm_{2}(X) \xcc \xbm_{0}(X),$

(2) Let $ \xbm $ satisfy $( \xbm \xcc )$ and $( \xbm Cum),$ then $ \xbm (U
\xcv U' ) \xcc U$ $ \xcj $ $ \xbm (U \xcv U' )= \xbm (U),$

(3) Let $ \xbm $ satisfy $( \xbm \xcc ),$ then $ \xbm_{i}(U) \xcc \xbm
(U),$ and $ \xbm_{i}(U) \xcc U,$

(4) Let $ \xbm $ satisfy $( \xbm \xcc )$ and one of the $( \xbm PRi),$
then
$ \xbm (A \xcv B) \xcc \xbm (A) \xcv \xbm (B),$

(5) Let $ \xbm $ satisfy $( \xbm \xcc )$ and one of the $( \xbm PRi),$
then
$ \xbm_{2}(X) \xcc \xbm_{1}(X),$

(6) Let $ \xbm $ satisfy $( \xbm \xcc ),$ $( \xbm PRi),$ then
$ \xbm_{i}(U) \xcc U' $ $ \xcj $ $ \xbm (U) \xcc U',$

(7) Let $ \xbm $ satisfy $( \xbm \xcc )$ and one of the $( \xbm PRi),$
then
$X \xcc Y,$ $ \xbm (X \xcv U) \xcc X$ $ \xch $ $ \xbm (Y \xcv U) \xcc Y,$

(8) Let $ \xbm $ satisfy $( \xbm \xcc )$ and one of the $( \xbm PRi),$
then
$X \xcc Y$ $ \xch $ $X \xcs \xbm_{i}(Y) \xcc \xbm_{i}(X)$ - so $( \xbm
PR)$ holds for $ \xbm_{i},$
(more precisely, only for $ \xbm_{2}$ we need the prerequisites, in the
other
cases the definition suffices)

(9) Let $ \xbm $ satisfy $( \xbm \xcc ),$ $( \xbm PR2),$ $( \xbm Cum),$
then
$ \xbm_{2}(X) \xcc Y \xcc X$ $ \xch $ $ \xbm_{2}(X)= \xbm_{2}(Y)$ - so $(
\xbm Cum)$ holds for $ \xbm_{2}.$

(10) $( \xbm \xcc )$ and $( \xbm Cum)$ for $ \xbm $ entail $ \xbm_{0}(U)=
\xbm_{1}(U).$

\efa

\subparagraph{
Proof
}

$\hspace{0.01em}$

% (+++ Orig.:  Proof +++)

(1) $ \xbm_{1}(X) \xcc \xbm_{0}(X)$ follows from $( \xbm \xcc )$ for $
\xbm.$ For $ \xbm_{2}:$ By $Y \xcc U,$ $U \xcv Y=U,$ so
$ \xbm (U) \xcc U$ by $( \xbm \xcc ).$

(2) $ \xbm (U \xcv U' ) \xcc U \xcc U \xcv U' $ $ \xch_{( \xbm CUM)}$ $
\xbm (U \xcv U' )= \xbm (U).$

(3) $ \xbm_{i}(U) \xcc U$ by definition. To show $ \xbm_{i}(U) \xcc \xbm
(U),$ take in all three cases
$Y:=U,$ and use for $i=1,2$ $( \xbm \xcc ).$

(4) By definition of $ \xbm_{0},$ we have $ \xbm_{0}(A \xcv B) \xcc A \xcv
B,$ $ \xbm_{0}(A \xcv B) \xcs (A- \xbm (A))= \xCQ,$
$ \xbm_{0}(A \xcv B) \xcs (B- \xbm (B))= \xCQ,$ so $ \xbm_{0}(A \xcv B)
\xcs A \xcc \xbm (A),$ $ \xbm_{0}(A \xcv B) \xcs B \xcc \xbm (B),$ and
$ \xbm_{0}(A \xcv B) \xcc \xbm (A) \xcv \xbm (B).$ By $ \xbm: \xdy \xcp
\xdy $ and $( \xcv ),$ $ \xbm (A) \xcv \xbm (B) \xbe \xdy.$
Moreover, by (3) $ \xbm_{0}(A \xcv B) \xcc \xbm (A \xcv B),$ so $
\xbm_{0}(A \xcv B)$ $ \xcc $
$( \xbm (A) \xcv \xbm (B)) \xcs \xbm (A \xcv B),$ so by (1) $ \xbm_{i}(A
\xcv B) \xcc ( \xbm (A) \xcv \xbm (B)) \xcs \xbm (A \xcv B)$
for $i=0, \Xl,2.$
If $ \xbm (A \xcv B) \xcC \xbm (A) \xcv \xbm (B),$ then $( \xbm (A) \xcv
\xbm (B)) \xcs \xbm (A \xcv B)$ $ \xcb $ $ \xbm (A \xcv B),$
contradicting $( \xbm PRi).$

(5) Let $Y \xbe \xdy,$ $ \xbm (Y) \xcc U,$ $x \xbe Y- \xbm (Y),$ then (by
(4)) $ \xbm (U \xcv Y) \xcc \xbm (U) \xcv \xbm (Y) \xcc U.$

(6) `` $ \xci $ '' by (3). `` $ \xch $ '': By $( \xbm PRi),$ $ \xbm (U)-
\xbm_{i}(U)$ is small, so there
is no $X \xbe \xdy $ s.t.
$ \xbm_{i}(U) \xcc X \xcb \xbm (U).$ If there were $U' \xbe \xdy $ s.t. $
\xbm_{i}(U) \xcc U',$ but $ \xbm (U) \xcC U',$ then
for $X:=U' \xcs \xbm (U) \xbe \xdy,$ $ \xbm_{i}(U) \xcc X \xcb \xbm (U),$
$contradiction.$

(7) $ \xbm (Y \xcv U)= \xbm (Y \xcv X \xcv U) \xcc_{(4)} \xbm (Y) \xcv
\xbm (X \xcv U) \xcc Y \xcv X=Y.$

(8) For $i=0,1:$ Let $x \xbe X- \xbm_{0}(X),$ then there is A s.t. $A \xcc
X,$ $x \xbe A- \xbm (A),$ so
$A \xcc Y.$ The case $i=1$ is similar. We need here only the definitions.
For $i=2:$ Let $x \xbe X- \xbm_{2}(X),$ A s.t. $x \xbe A- \xbm (A),$ $
\xbm (X \xcv A) \xcc X,$ then by
(7) $ \xbm (Y \xcv A) \xcc Y.$

(9) `` $ \xcc $ '': Let $x \xbe \xbm_{2}(X),$ so $x \xbe Y,$ and $x \xbe
\xbm_{2}(Y)$ by (8).
`` $ \xcd $ '': Let $x \xbe \xbm_{2}(Y),$ so $x \xbe X.$ Suppose $x \xce
\xbm_{2}(X),$ so there is $U \xbe \xdy $ s.t.
$x \xbe U- \xbm (U)$ and $ \xbm (X \xcv U) \xcc X.$ Note that by $ \xbm (X
\xcv U) \xcc X$ and (2), $ \xbm (X \xcv U)= \xbm (X).$
Now, $ \xbm_{2}(X) \xcc Y,$ so by (6) $ \xbm (X) \xcc Y,$ thus $ \xbm (X
\xcv U)= \xbm (X) \xcc Y \xcc Y \xcv U \xcc X \xcv U,$
so $ \xbm (Y \xcv U)= \xbm (X \xcv U)= \xbm (X) \xcc Y$ by $( \xbm Cum),$
so $x \xce \xbm_{2}(Y),$ $contradiction.$

(10) $ \xbm_{1}(U) \xcc \xbm_{0}(U)$ by (1). Let $Y$ s.t. $ \xbm (Y) \xcc
U,$ $x \xbe Y- \xbm (Y),$ $x \xbe U.$ Consider
$Y \xcs U,$ $x \xbe Y \xcs U,$ $ \xbm (Y) \xcc Y \xcs U \xcc Y,$ so $ \xbm
(Y)= \xbm (Y \xcs U)$ by $( \xbm Cum),$ and
$x \xce \xbm (Y \xcs U).$ Thus, $ \xbm_{0}(U) \xcc \xbm_{1}(U).$

$ \xcz $
\\[3ex]

\bfa

$\hspace{0.01em}$

% (+++ Orig. No.:  Fact D-7.2.2 +++)

\label{Fact D-7.2.2}

In the presence of $( \xbm \xcc ),$ $( \xbm Cum)$ for $ \xbm,$ we have:
$( \xbm PR0) \xcj ( \xbm PR1),$ and $( \xbm PR2) \xch ( \xbm PR1).$

If $( \xbm PR)$ also holds for $ \xbm,$ then so will $( \xbm PR1) \xch (
\xbm PR2).$

(Recall that $( \xcv )$ and $( \xcs )$ are assumed to hold.)

\efa

\subparagraph{
Proof
}

$\hspace{0.01em}$

% (+++ Orig.:  Proof +++)

$( \xbm PR0) \xcj ( \xbm PR1):$ By Fact \ref{Fact D-7.2.1} (page \pageref{Fact
D-7.2.1}), (10),
$ \xbm_{0}(U)= \xbm_{1}(U)$ if $( \xbm Cum)$ holds for $ \xbm.$

$( \xbm PR2) \xch ( \xbm PR1):$
Suppose $( \xbm PR2)$ holds. By $( \xbm PR2)$ and (5), $ \xbm_{2}(U) \xcc
\xbm_{1}(U),$ so
$ \xbm (U)- \xbm_{1}(U) \xcc \xbm (U)- \xbm_{2}(U).$ By $( \xbm PR2),$ $
\xbm (U)- \xbm_{2}(U)$ is small, then so is
$ \xbm (U)- \xbm_{1}(U),$ so $( \xbm PR1)$ holds.

$( \xbm PR1) \xch ( \xbm PR2):$
Suppose $( \xbm PR1)$ holds, and $( \xbm PR2)$ fails. By failure of $(
\xbm PR2),$ there is
$X \xbe \xdy $ s.t. $ \xbm_{2}(U) \xcc X \xcb \xbm (U).$ Let $x \xbe \xbm
(U)-$X, as $x \xce \xbm_{2}(U),$ there is $Y$ s.t.
$ \xbm (U \xcv Y) \xcc U,$ $x \xbe Y- \xbm (Y).$ Let $Z:=U \xcv Y \xcv X.$
By $( \xbm PR),$ $x \xce \xbm (U \xcv Y),$ and
$x \xce \xbm (U \xcv Y \xcv X).$
Moreover, $ \xbm (U \xcv X \xcv Y) \xcc \xbm (U \xcv Y) \xcv \xbm (X)$
by Fact \ref{Fact D-7.2.1} (page \pageref{Fact D-7.2.1})  (4),
$ \xbm (U \xcv Y) \xcc U,$ $ \xbm (X) \xcc X \xcc \xbm (U) \xcc U$ by
prerequisite, so
$ \xbm (U \xcv X \xcv Y) \xcc U \xcc U \xcv Y \xcc U \xcv X \xcv Y,$ so $
\xbm (U \xcv X \xcv Y)= \xbm (U \xcv Y) \xcc U.$
Thus, $x \xce \xbm_{1}(U),$ and $ \xbm_{1}(U) \xcc X,$ too, a
contradiction.

$ \xcz $
\\[3ex]

Here is an example which shows that
Fact \ref{Fact D-7.2.1} (page \pageref{Fact D-7.2.1}), (9) may fail for
$ \xbm_{0}$ and $ \xbm_{1}.$

\be

$\hspace{0.01em}$

% (+++ Orig. No.:  Example D-7.2.2 +++)

\label{Example D-7.2.2}

Consider $ \xdl $ with $v( \xdl ):=\{p_{i}:i \xbe \xbo \}.$ Let $m \xcM
p_{0},$ let $m' \xbe M(p_{0})$ arbitrary.
Make for each $n \xbe M(p_{0})-\{m' \}$ one copy of $m,$ likewise
of $m',$ set $ \xBc m,n \xBe  \xeb  \xBc m',n \xBe $ for all $n,$ and
$n \xeb  \xBc m,n \xBe,$ $n \xeb  \xBc m',n \xBe $ for all $n.$ The
resulting structure $ \xdz $ is smooth and transitive. Let $ \xdy:=
\xdD_{ \xdl },$ define
$ \xbm (X):= \wt{ \xbm_{ \xdz }(X)}$ for $X \xbe \xdy.$

Let $m' \xbe X- \xbm_{ \xdz }(X).$ Then $m \xbe X,$ or $M(p_{0}) \xcc X.$
In the latter case, as all $m'' $ s.t.
$m'' \xEd m',$ $m'' \xcm p_{0}$ are minimal, $M(p_{0})-\{m' \} \xcc
\xbm_{ \xdz }(X),$ so $m' \xbe \wt{ \xbm_{ \xdz }(X)}= \xbm (X).$ Thus,
as $ \xbm_{ \xdz }(X) \xcc \xbm (X),$ if $m' \xbe X- \xbm (X),$ then $m
\xbe X.$

Define now $X:=M(p_{0}) \xcv \{m\},$ $Y:=M(p_{0}).$

We first show that $ \xbm_{0}$ does not satisfy $( \xbm Cum).$
$ \xbm_{0}(X):=\{x \xbe X: \xCN \xcE A \xbe \xdy (A \xcc X:x \xbe A- \xbm
(A))\}.$ $m \xce \xbm_{0}(X),$ as $m \xce \xbm (X)= \wt{ \xbm_{ \xdz
}(X)}.$
Moreover, $m' \xce \xbm_{0}(X),$ as $\{m,m' \} \xbe \xdy,$ $\{m,m' \}
\xcc X,$ and $ \xbm (\{m,m' \})= \xbm_{ \xdz }(\{m,m' \})=\{m\}.$
So $ \xbm_{0}(X) \xcc Y \xcc X.$ Consider now $ \xbm_{0}(Y).$ As $m \xce
Y,$ for any $A \xbe \xdy,$ $A \xcc Y,$ if
$m' \xbe A,$ then $m' \xbe \xbm (A),$ too, by above argument, so $m' \xbe
\xbm_{0}(Y),$ and $ \xbm_{0}$ does
not satisfy $( \xbm Cum).$

We turn to $ \xbm_{1}.$

By Fact \ref{Fact D-7.2.1} (page \pageref{Fact D-7.2.1})  (1),
$ \xbm_{1}(X) \xcc \xbm_{0}(X),$ so $m,m' \xce \xbm_{1}(X),$ and
again $ \xbm_{1}(X) \xcc Y \xcc X.$
Consider again $ \xbm_{1}(Y).$ As $m \xce Y,$ for any $A \xbe \xdy,$ $
\xbm (A) \xcc Y,$ if $m' \xbe A,$
then $m' \xbe \xbm (A),$ too: if $M(p_{0})-\{m' \} \xcc A,$ then $m' \xbe
\wt{ \xbm_{ \xdz }(A)},$ if $M(p_{0})-\{m' \} \xcC A,$
but $m' \xbe A,$ then either $m' \xbe \xbm_{ \xdz }(A),$ or $m \xbe \xbm_{
\xdz }(A) \xcc \xbm (A),$ but $m \xce Y.$ Thus,
$( \xbm Cum)$ fails for $ \xbm_{1},$ too.

It remains to show that $ \xbm $ satisfies $( \xbm \xcc ),$ $( \xbm Cum),$
$( \xbm PR0),$ $( \xbm PR1).$
Note that by Fact \ref{Fact Cum-Alpha-HU} (page \pageref{Fact Cum-Alpha-HU}) 
(3)
and Proposition \ref{Proposition D-4.4.6} (page \pageref{Proposition D-4.4.6})
$ \xbm_{ \xdz }$ satisfies $( \xbm Cum),$ as $ \xdz $ is smooth.
$( \xbm \xcc )$ is trivial. We show $( \xbm PRi)$ for $i=0,1.$ As $ \xbm_{
\xdz }(A) \xcc \xbm (A),$ by $( \xbm PR)$
and $( \xbm Cum)$ for $ \xbm_{ \xdz },$ $ \xbm_{ \xdz }(X) \xcc
\xbm_{0}(X)$ and $ \xbm_{ \xdz }(X) \xcc \xbm_{1}(X):$
To see this, we note
$ \xbm_{ \xdz }(X) \xcc \xbm_{0}(X):$ Let $x \xbe X- \xbm_{0}(X),$ then
there is $Y$ s.t. $x \xbe Y- \xbm (Y).Y \xcc X,$
but $ \xbm_{ \xdz }(Y) \xcc \xbm (Y),$ so by $Y \xcc X$ and $( \xbm PR)$
for $ \xbm_{ \xdz }$ $x \xce \xbm_{ \xdz }(X).$
$ \xbm_{ \xdz }(X) \xcc \xbm_{1}(X):$ Let $x \xbe X- \xbm_{1}(X),$ then
there is $Y$ s.t. $x \xbe Y- \xbm (Y),$ $ \xbm (Y) \xcc X,$
so $x \xbe Y- \xbm_{ \xdz }(Y)$ and $ \xbm_{ \xdz }(Y) \xcc X.$ $ \xbm_{
\xdz }(X \xcv Y) \xcc \xbm_{ \xdz }(X) \xcv \xbm_{ \xdz }(Y) \xcc X \xcc X
\xcv Y,$ so
$ \xbm_{ \xdz }(X \xcv Y)= \xbm_{ \xdz }(X)$ by $( \xbm Cum)$ for $ \xbm_{
\xdz }.$ $x \xbe Y- \xbm_{ \xdz }(Y)$ $ \xch $ $x \xce \xbm_{ \xdz }(X
\xcv Y)$ by
$( \xbm PR)$ for $ \xbm_{ \xdz },$ so $x \xce \xbm_{ \xdz }(X).$

But by Fact \ref{Fact D-7.2.1} (page \pageref{Fact D-7.2.1}), (3)
$ \xbm_{i}(X) \xcc \xbm (X).$ As by definition, $ \xbm (X)- \xbm_{ \xdz
}(X)$ is small, $( \xbm PRi)$ hold for
$i=0,1.$ It remains to show $( \xbm Cum)$ for $ \xbm.$ Let $ \xbm (X)
\xcc Y \xcc X,$ then
$ \xbm_{ \xdz }(X) \xcc \xbm (X) \xcc Y \xcc X,$ so by $( \xbm Cum)$ for $
\xbm_{ \xdz }$
$ \xbm_{ \xdz }(X)= \xbm_{ \xdz }(Y),$ so by definition of $ \xbm,$ $
\xbm (X)= \xbm (Y).$

(Note that by Fact \ref{Fact D-7.2.1} (page \pageref{Fact D-7.2.1})  (10),
$ \xbm_{0}= \xbm_{1}$ follows from $( \xbm Cum)$
for $ \xbm,$ so we could have demonstrated part of the properties also
differently.)

$ \xcz $
\\[3ex]

\ee

By Fact \ref{Fact D-7.2.1} (page \pageref{Fact D-7.2.1})  (3) and (8) and
Proposition \ref{Proposition Pref-Complete-Trans} (page \pageref{Proposition
Pref-Complete-Trans}), $ \xbm_{0}$
has a representation by a (transitive) preferential structure, if $ \xbm:
\xdy \xcp \xdy $
satisfies $( \xbm \xcc )$ and $( \xbm PR0),$ and $ \xbm_{0}$ is defined as
in
Definition \ref{Definition D-7.2.1} (page \pageref{Definition D-7.2.1}).

We thus have (taken from  \cite{Sch04}, Proposition 5.2.5 there):

\bp

$\hspace{0.01em}$

% (+++ Orig. No.:  Proposition D-7.2.3 +++)

\label{Proposition D-7.2.3}

Let $Z$ be an arbitrary set, $ \xdy \xcc \xdp (Z),$ $ \xbm: \xdy \xcp
\xdy,$ $ \xdy $ closed under arbitrary
intersections and finite unions, and $ \xCQ,Z \xbe \xdy,$ and let $
\wt{.}$ be defined wrt. $ \xdy.$

(a) If $ \xbm $ satisfies $( \xbm \xcc ),$ $( \xbm PR0),$ then there is a
transitive
preferential structure $ \xdz $ over $Z$ s.t. for all $U \xbe \xdy $ $
\xbm (U)= \wt{ \xbm_{ \xdz }(U)}.$

(b) If $ \xdz $ is a preferential structure over $Z$ and
$ \xbm: \xdy \xcp \xdy $ s.t. for all $U \xbe \xdy $ $ \xbm (U)= \wt{
\xbm_{ \xdz }(U)},$ then $ \xbm $ satisfies $( \xbm \xcc ),$ $( \xbm
PR0).$

\ep

\subparagraph{
Proof
}

$\hspace{0.01em}$

% (+++ Orig.:  Proof +++)

(a) Let $ \xbm $ satisfy $( \xbm \xcc ),$ $( \xbm PR0).$ $ \xbm_{0}$ as
defined in
Definition \ref{Definition D-7.2.1} (page \pageref{Definition D-7.2.1})
satisfies properties $( \xbm \xcc ),$ $( \xbm PR)$ by
Fact \ref{Fact D-7.2.1} (page \pageref{Fact D-7.2.1}), (3) and (8).
Thus, by Proposition \ref{Proposition Pref-Complete-Trans} (page
\pageref{Proposition Pref-Complete-Trans}),
there is a transitive
structure $ \xdz $
over $Z$ s.t. $ \xbm_{0}= \xbm_{ \xdz }$, but by $( \xbm PR0)$ $ \xbm
(U)= \wt{ \xbm_{0}(U)}= \wt{ \xbm_{ \xdz }(U)}$ for $U \xbe \xdy.$

(b) $( \xbm \xcc ):$ $ \xbm_{ \xdz }(U) \xcc U,$ so by $U \xbe \xdy $ $
\xbm (U)= \wt{ \xbm_{ \xdz }(U)} \xcc U.$

$( \xbm PR0):$ If $( \xbm PR0)$ is false, there is $U \xbe \xdy $ s.t. for
$U':= \xcV \{Y' - \xbm (Y' ):Y' \xbe \xdy,Y' \xcc U\}$ $ \wt{ \xbm
(U)-U' } \xcb \xbm (U).$ By $ \xbm_{ \xdz }(Y' ) \xcc \xbm (Y' ),$
$Y' - \xbm (Y' ) \xcc Y' - \xbm_{ \xdz }(Y' ).$ No copy of any $x \xbe Y'
- \xbm_{ \xdz }(Y' )$ with $Y' \xcc U,$ $Y' \xbe \xdy $ can be minimal in
$ \xdz \xex U.$ Thus, by $ \xbm_{ \xdz }(U) \xcc \xbm (U),$ $ \xbm_{ \xdz
}(U) \xcc \xbm (U)-U',$ so $ \wt{ \xbm_{ \xdz }(U)} \xcc \wt{ \xbm (U)-U'
} \xcb \xbm (U),$ $contradiction.$

$ \xcz $
\\[3ex]

We turn to the smooth case.

If $ \xbm: \xdy \xcp \xdy $ satisfies $( \xbm \xcc ),$ $( \xbm PR2),$ $(
\xbm CUM)$ and $ \xbm_{2}$ is
defined from $ \xbm $ as in
Definition \ref{Definition D-7.2.1} (page \pageref{Definition D-7.2.1}),
then $ \xbm_{2}$ satisfies $( \xbm \xcc ),$ $( \xbm PR),$ $( \xbm Cum)$ by
Fact \ref{Fact D-7.2.1} (page \pageref{Fact D-7.2.1})  (3),
(8), and (9), and can thus be represented by a (transitive) smooth
structure, by
Proposition \ref{Proposition D-4.5.3} (page \pageref{Proposition D-4.5.3}), and
we finally have (taken from  \cite{Sch04}, Proposition 5.2.9 there):

\bp

$\hspace{0.01em}$

% (+++ Orig. No.:  Proposition D-7.2.4 +++)

\label{Proposition D-7.2.4}

Let $Z$ be an arbitrary set, $ \xdy \xcc \xdp (Z),$ $ \xbm: \xdy \xcp
\xdy,$ $ \xdy $ closed under arbitrary
intersections and finite unions, and $ \xCQ,Z \xbe \xdy,$ and let $
\wt{.}$ be defined wrt. $ \xdy.$

(a) If $ \xbm $ satisfies $( \xbm \xcc ),$ $( \xbm PR2),$ $( \xbm CUM),$
then there is a transitive smooth
preferential structure $ \xdz $ over $Z$ s.t. for all $U \xbe \xdy $ $
\xbm (U)= \wt{ \xbm_{ \xdz }(U)}.$

(b) If $ \xdz $ is a smooth preferential structure over $Z$ and
$ \xbm: \xdy \xcp \xdy $ s.t. for all $U \xbe \xdy $ $ \xbm (U)= \wt{
\xbm_{ \xdz }(U)},$ then $ \xbm $ satisfies
$( \xbm \xcc ),$ $( \xbm PR2),$ $( \xbm CUM).$

\ep

\subparagraph{
Proof
}

$\hspace{0.01em}$

% (+++ Orig.:  Proof +++)

(a) If $ \xbm $ satisfies $( \xbm \xcc ),$ $( \xbm PR2),$ $( \xbm CUM),$
then $ \xbm_{2}$
defined from $ \xbm $ as in Definition \ref{Definition D-7.2.1} (page
\pageref{Definition D-7.2.1})
satisfies $( \xbm \xcc ),$ $( \xbm PR),$
$( \xbm CUM)$ by Fact \ref{Fact D-7.2.1} (page \pageref{Fact D-7.2.1})  (3), (8)
and (9).
Thus, by Proposition \ref{Proposition D-4.5.3} (page \pageref{Proposition
D-4.5.3}),
there is a smooth transitive preferential
structure $ \xdz $ over $Z$ s.t. $ \xbm_{2}= \xbm_{ \xdz }$, but by $(
\xbm PR2)$ $ \xbm (U)= \wt{ \xbm_{2}(U)}= \wt{ \xbm_{ \xdz }(U)}.$

(b)
$( \xbm \xcc ):$ $ \xbm_{ \xdz }(U) \xcc U$ $ \xch $ $ \xbm (U)= \wt{
\xbm_{ \xdz }(U)} \xcc U$ by $U \xbe \xdy.$

$( \xbm PR2):$ If $( \xbm PR2)$ fails, then there is $U \xbe \xdy $ s.t.
for
$U':= \xcV \{Y' - \xbm (Y' ):$ $Y' \xbe \xdy,$ $ \xbm (U \xcv Y' ) \xcc
U\}$ $ \wt{ \xbm (U)-U' } \xcb \xbm (U).$

By $ \xbm_{ \xdz }(Y' ) \xcc \xbm (Y' ),$ $Y' - \xbm (Y' ) \xcc Y' -
\xbm_{ \xdz }(Y' ).$
But no copy of any $x \xbe Y' - \xbm_{ \xdz }(Y' )$ with $ \xbm_{ \xdz }(U
\xcv Y' ) \xcc \xbm (U \xcv Y' ) \xcc U$ can be minimal in $ \xdz \xex U:$
As $x \xbe Y' - \xbm_{ \xdz }(Y' ),$ if $ \xBc x,i \xBe $ is any copy of $x,$
then
there
is $ \xBc y,j \xBe  \xeb  \xBc x,i \xBe,$ $y \xbe Y'.$
Consider now $U \xcv Y'.$ As $ \xBc x,i \xBe $ is not minimal in $ \xdz \xex U
\xcv
Y',$ by smoothness of $ \xdz $
there must be $ \xBc z,k \xBe  \xeb  \xBc x,i \xBe,$
$ \xBc z,k \xBe $ minimal in $ \xdz \xex U \xcv Y'.$ But all minimal elements of
$
\xdz \xex U \xcv Y' $ must be in $ \xdz \xex U,$
so there must be $ \xBc z,k \xBe  \xeb  \xBc x,i \xBe,$ $z \xbe U,$ thus
$ \xBc x,i \xBe $ is not minimal in $ \xdz \xex U.$
Thus by $ \xbm_{ \xdz }(U) \xcc \xbm (U),$ $ \xbm_{ \xdz }(U) \xcc \xbm
(U)-U',$ so $ \wt{ \xbm_{ \xdz }(U)} \xcc \wt{ \xbm (U)-U' } \xcb \xbm
(U),$ $contradiction.$

$( \xbm CUM):$ Let $ \xbm (X) \xcc Y \xcc X.$ Now $ \xbm_{ \xdz }(X) \xcc
\wt{ \xbm_{ \xdz }(X)}= \xbm (X),$ so by smoothness of $ \xdz $
$ \xbm_{ \xdz }(Y)= \xbm_{ \xdz }(X),$ thus $ \xbm (X)= \wt{ \xbm_{ \xdz
}(X)}= \wt{ \xbm_{ \xdz }(Y)}= \xbm (Y).$
$ \xcz $
\\[3ex]
\subsection{Ranked structures}
\label{Section 2.2.10.3}
%  Section (7.3):  Ranked structures
%  Section (7.3):  Ranked structures
% %
% ===================================

We recall from  \cite{Sch04}
and Section \ref{Section 2.2.8.3} (page \pageref{Section 2.2.8.3})
above the basic properties of ranked structures.

We give now an easy version of representation results for ranked
structures
without definability preservation.

\bn

$\hspace{0.01em}$

% (+++ Orig. No.:  Notation D-7.3.1 +++)

\label{Notation D-7.3.1}

We abbreviate $ \xbm (\{x,y\})$ by $ \xbm (x,y)$ etc.

\en

\bfa

$\hspace{0.01em}$

% (+++ Orig. No.:  Fact D-7.3.1 +++)

\label{Fact D-7.3.1}

Let the domain contain singletons and be closed under $( \xcv ).$

Let for $ \xbm: \xdy \xcp \xdy $ hold:

$( \xbm =)$ for finite sets, $( \xbm \xbe ),$ $( \xbm PR3),$ $( \xbm \xCQ
fin).$

Then the following properties hold for $ \xbm_{3}$ as defined in
Definition \ref{Definition D-7.2.1} (page \pageref{Definition D-7.2.1}) :

(1) $ \xbm_{3}(X) \xcc \xbm (X),$

(2) for finite $X,$ $ \xbm (X)= \xbm_{3}(X),$

(3) $( \xbm \xcc ),$

(4) $( \xbm PR),$

(5) $( \xbm \xCQ fin),$

(6) $( \xbm =),$

(7) $( \xbm \xbe ),$

(8) $ \xbm (X)= \wt{ \xbm_{3}(X)}.$

\efa

\subparagraph{
Proof
}

$\hspace{0.01em}$

% (+++ Orig.:  Proof +++)

(1) Suppose not, so $x \xbe \xbm_{3}(X),$ $x \xbe X- \xbm (X),$ so by $(
\xbm \xbe )$ for $ \xbm,$ there is
$y \xbe X,$ $x \xce \xbm (x,y),$ $contradiction.$

(2) By $( \xbm PR3)$ for $ \xbm $ and (1), for finite $U$ $ \xbm (U)=
\xbm_{3}(U).$

(3) $( \xbm \xcc )$ is trivial for $ \xbm_{3}.$

(4) Let $X \xcc Y,$ $x \xbe \xbm_{3}(Y) \xcs X,$ suppose $x \xbe X-
\xbm_{3}(X),$ so there is $y \xbe X \xcc Y,$
$x \xce \xbm (x,y),$ so $x \xce \xbm_{3}(Y).$

(5) $( \xbm \xCQ fin)$ for $ \xbm_{3}$ follows from $( \xbm \xCQ fin)$ for
$ \xbm $ and (2).

(6) Let $X \xcc Y,$ $y \xbe \xbm_{3}(Y) \xcs X,$ $x \xbe \xbm_{3}(X),$ we
have to show $x \xbe \xbm_{3}(Y).$
By (4), $y \xbe \xbm_{3}(X).$
Suppose $x \xce \xbm_{3}(Y).$ So there is $z \xbe Y.x \xce \xbm (x,z).$ As
$y \xbe \xbm_{3}(Y),$ $y \xbe \xbm (y,z).$
As $x \xbe \xbm_{3}(X),$ $x \xbe \xbm (x,y),$ as $y \xbe \xbm_{3}(X),$ $y
\xbe \xbm (x,y).$
Consider $\{x,y,z\}.$ Suppose $y \xce \xbm (x,y,z),$ then by $( \xbm \xbe
)$ for $ \xbm,$ $y \xce \xbm (x,y)$ or
$y \xce \xbm (y,z),$ $contradiction.$
Thus $y \xbe \xbm (x,y,z) \xcs \xbm (x,y).$ As $x \xbe \xbm (x,y),$ and $(
\xbm =)$ for $ \xbm $ and finite sets,
$x \xbe \xbm (x,y,z).$ Recall that $x \xce \xbm (x,z).$ But for finite
sets $ \xbm = \xbm_{3},$ and by
(4) $( \xbm PR)$ holds for $ \xbm_{3},$ so it holds for $ \xbm $ and
finite sets. $contradiction$

(7) Let $x \xbe X- \xbm_{3}(X),$ so there is $y \xbe X.x \xce \xbm (x,y)=
\xbm_{3}(x,y).$

(8) As $ \xbm (X) \xbe \xdy,$ and $ \xbm_{3}(X) \xcc \xbm (X),$ $ \wt{
\xbm_{3}(X)} \xcc \xbm (X),$ so by $( \xbm PR3)$ $ \wt{ \xbm_{3}(X)}= \xbm
(X).$

$ \xcz $
\\[3ex]

\bfa

$\hspace{0.01em}$

% (+++ Orig. No.:  Fact D-7.3.2 +++)

\label{Fact D-7.3.2}

If $ \xdz $ is ranked, and we define $ \xbm (X):= \wt{ \xbm_{ \xdz }(X)},$
and $ \xdz $ has no
copies, then the following hold:

(1) $ \xbm_{ \xdz }(X)=\{x \xbe X: \xcA y \xbe X.x \xbe \xbm (x,y)\},$ so
$ \xbm_{ \xdz }(X)= \xbm_{3}(X)$ for $X \xbe \xdy,$

(2) $ \xbm (X)= \xbm_{ \xdz }(X)$ for finite $X,$

(3) $( \xbm =)$ for finite sets for $ \xbm,$

(4) $( \xbm \xbe )$ for $ \xbm,$

(5) $( \xbm \xCQ fin)$ for $ \xbm,$

(6) $( \xbm PR3)$ for $ \xbm.$

\efa

\subparagraph{
Proof
}

$\hspace{0.01em}$

% (+++ Orig.:  Proof +++)

(1) holds for ranked structures.

(2) and (6) are trivial. (3) and (5) hold for $ \xbm_{ \xdz },$ so by (2)
for $ \xbm.$

(4) If $x \xce \xbm (X),$ then $x \xce \xbm_{ \xdz }(X),$ $( \xbm \xbe )$
holds for $ \xbm_{ \xdz },$ so there is $y \xbe X$ s.t.
$x \xce \xbm_{ \xdz }(x,y)= \xbm (x,y)$ by (2).

$ \xcz $
\\[3ex]

We summarize:

\bp

$\hspace{0.01em}$

% (+++ Orig. No.:  Proposition D-7.3.3 +++)

\label{Proposition D-7.3.3}

Let $Z$ be an arbitrary set, $ \xdy \xcc \xdp (Z),$ $ \xbm: \xdy \xcp
\xdy,$ $ \xdy $ closed under arbitrary
intersections and finite unions, contain singletons, and $ \xCQ,Z \xbe
\xdy,$ and let $ \wt{.}$
be defined wrt. $ \xdy.$

(a) If $ \xbm $ satisfies $( \xbm =)$ for finite sets, $( \xbm \xbe ),$ $(
\xbm PR3),$ $( \xbm \xCQ fin),$
then there is a ranked preferential structure $ \xdz $ without copies
over $Z$ s.t. for all $U \xbe \xdy $ $ \xbm (U)= \wt{ \xbm_{ \xdz }(U)}.$

(b) If $ \xdz $ is a ranked preferential structure over $Z$ without copies
and
$ \xbm: \xdy \xcp \xdy $ s.t. for all $U \xbe \xdy $ $ \xbm (U)= \wt{
\xbm_{ \xdz }(U)},$ then $ \xbm $ satisfies
$( \xbm =)$ for finite sets, $( \xbm \xbe ),$ $( \xbm PR3),$ $( \xbm \xCQ
fin).$

\ep

\subparagraph{
Proof
}

$\hspace{0.01em}$

% (+++ Orig.:  Proof +++)

(a) Let $ \xbm $ satisfy $( \xbm =)$ for finite sets, $( \xbm \xbe ),$ $(
\xbm PR3),$ $( \xbm \xCQ fin),$ then
$ \xbm_{3}$ as defined in Definition \ref{Definition D-7.2.1} (page
\pageref{Definition D-7.2.1})
satisfies properties $( \xbm \xcc ),$ $( \xbm \xCQ fin),$ $( \xbm =),$ $(
\xbm \xbe )$ by
Fact \ref{Fact D-7.3.1} (page \pageref{Fact D-7.3.1}).
Thus, by Proposition \ref{Proposition D-5.3.5} (page \pageref{Proposition
D-5.3.5}),
there is a transitive structure $ \xdz $
over $Z$ s.t. $ \xbm_{3}= \xbm_{ \xdz }$, but by Fact \ref{Fact D-7.3.1} (page
\pageref{Fact D-7.3.1})  (8)
$ \xbm (U)= \wt{ \xbm_{3}(U)}= \wt{ \xbm_{ \xdz }(U)}$ for $U \xbe \xdy.$

(b) This was shown in Fact \ref{Fact D-7.3.2} (page \pageref{Fact D-7.3.2}).

$ \xcz $
\\[3ex]
\section{The limit variant}
\label{Section Limit}
%  CHAPTER (8):  THE LIMIT VARIANT
%  CHAPTER (8):  THE LIMIT VARIANT
% %
% =================================
\subsection{Introduction}
\label{Section 2.2.11.1}
%  Section (8.1):  THE LIMIT VARIANT, INTRODUCTION
%  Section (8.1):  THE LIMIT VARIANT, INTRODUCTION
% %
% =================================================

Distance based semantics give perhaps the clearest motivation for the
limit
variant. For instance,
the Stalnaker/Lewis semantics for counterfactual conditionals defines
$ \xbf > \xbq $ to hold in a (classical) model $m$ iff in those models of
$ \xbf,$ which are
closest to $m,$ $ \xbq $ holds. For this to make sense, we need, of
course, a distance
$d$ on the model set. We call this approach the minimal variant.
Usually, one makes a limit assumption: The set of $ \xbf -$models
closest to $m$ is not empty if $ \xbf $ is consistent - i.e. the $ \xbf
-$models are not
arranged around $m$ in a way that they come closer and closer, without a
minimal
distance. This is, of course, a very strong assumption, and which is
probably
difficult to justify philosophically. It seems to have its only
justification
in the fact that it avoids degenerate cases, where, in above example, for
consistent $ \xbf $ $m \xcm \xbf >FALSE$ holds. As such, this assumption
is unsatisfactory.

Our aim here is to analyze the limit version more closely, in particular,
to see
criteria whether the much more complex limit version can be reduced to the
simpler minimal variant. In the limit version, roughly, $ \xbq $ is a
consequence
of $ \xbf,$ if $ \xbq $ holds ``in the limit'' in all $ \xbf -$models.
That is, iff, ``going sufficiently far down'', $ \xbq $ will become and
stay true.

The problem is not simple, as there are two sides which come into play,
and
sometimes we need both to cooperate to achieve a satisfactory translation.

The first component is what we call the ``algebraic limit'', i.e. we
stipulate
that the limit version should have properties which correspond to the
algebraic
properties of the minimal variant. An exact correspondence cannot always
be
achieved, and we give a translation which seems reasonable.

But once the translation is done, even if it is exact, there might still
be
problems linked to translation to logic.

 \xEh

 \xDH The structural limit:
It is a natural and much more convincing solution to the problem described
above to modify the
basic definition, and work without the rather artificial assumption that
the closest world exists. We adopt
what we call a ``limit approach'', and define $m \xcm \xbf > \xbq $ iff
there is a distance $d' $
such that for all $m' \xcm \xbf $ and $d(m,m' ) \xck d' $ $m' \xcm \xbq.$
Thus, from a certain point
onward, $ \xbq $ becomes and stays true. We will call this definition the
structural
limit, as it is based directly on the structure (the distance on the model
set).

 \xDH The algebraic limit:
The model sets to consider are spheres around $m,$ $S:=\{m' \xbe M( \xbf
):d(m,m' ) \xck d' \}$ for
some $d',$ s.t. $S \xEd \xCQ.$ The system of such $S$ is nested, i.e.
totally ordered by
inclusion; and if $m \xcm \xbf,$ it has a smallest element $\{m\},$ etc.
When we forget the
underlying structure, and consider just the properties of these systems of
spheres around different $m,$ and for different $ \xbf,$ we obtain what
we call the
algebraic limit.

 \xDH The logical limit:
The logical limit speaks about the logical properties which hold ``in the
limit'',
i.e. finally in all such sphere systems.

 \xEj

The interest to investigate this algebraic limit is twofold: first, we
shall
see (for other kinds of structures) that there are reasonable and not so
reasonable algebraic limits. Second, this distinction permits us to
separate
algebraic from logical problems, which have to do with definability of
model
sets, in short definability problems. We will see that
we find common definability problems and also common solutions in the
usual
minimal, and the limit variant.

In particular, the decomposition into three layers on both sides (minimal
and
limit version) can reveal that a
(seemingly) natural notion of structural limit results in algebraic
properties which have not much to do any more with the minimal variant.
So, to
speak about a limit variant, we will demand that this variant is not only
a natural structural limit, but results in a natural abstract limit, too.
Conversely, if the algebraic limit preserves the properties of the minimal
variant, there is hope that it preserves the logical properties, too - not
more
than hope, however, due to definability problems.

We give now the basic definitions for preferential and ranked preferential
structures.

\bd

$\hspace{0.01em}$

% (+++ Orig. No.:  Definition D-8.1.1 +++)

\label{Definition D-8.1.1}

(1) General preferential structures

(1.1) The version without copies:

Let $ \xdm:= \xBc U, \xeb  \xBe.$ Define

$Y \xcc X \xcc U$ is a minimizing initial segment, or MISE, of $X$ iff:

(a) $ \xcA x \xbe X \xcE x \xbe Y.y \xec x$ - where $y \xec x$ stands for
$x \xeb y$ or $x=y$
(i.e. $Y$ is minimizing)
and

(b) $ \xcA y \xbe Y, \xcA x \xbe X(x \xeb y$ $ \xch $ $x \xbe Y)$
(i.e. $Y$ is downward closed or an initial part).

(1.2) The version with copies:

Let $ \xdm:= \xBc  \xdu, \xeb  \xBe $ be as above. Define for $Y \xcc X \xcc
\xdu
$

$Y$ is a minimizing initial segment, or MISE of $X$ iff:

(a) $ \xcA  \xBc x,i \xBe  \xbe X \xcE  \xBc y,j \xBe  \xbe Y. \xBc y,j \xBe 
\xec  \xBc x,i \xBe $

and

(b) $ \xcA  \xBc y,j \xBe  \xbe Y, \xcA  \xBc x,i \xBe  \xbe X$
$( \xBc x,i \xBe  \xeb  \xBc y,j \xBe $ $ \xch $ $ \xBc x,i \xBe  \xbe Y).$

(1.3) For $X \xcc \xdu,$ let $ \xbL (X)$ be the set of MISE of $X.$

(1.4) We say that a set $ \xdx $ of MISE is cofinal in another set of MISE
$ \xdx ' $ (for the same base set $X)$ iff for all $Y' \xbe \xdx ',$
there is $Y \xbe \xdx,$ $Y \xcc Y'.$

(1.5) A MISE $X$ is called definable iff
$\{x: \xcE i. \xBc x,i \xBe  \xbe X\} \xbe \xdD_{ \xdl }.$

(1.6) $T \xcm_{ \xdm } \xbf $ iff there is $Y \xbe \xbL ( \xdu \xex M(T))$
s.t. $Y \xcm \xbf.$

$( \xdu \xex M(T):=\{ \xBc x,i \xBe  \xbe \xdu:x \xbe M(T)\}$ - if there are no
copies, we simplify in
the obvious way.)

(2) Ranked preferential structures

In the case of ranked structures, we may assume without loss of generality
that
the MISE sets have a particularly simple form:

For $X \xcc U$ $A \xcc X$ is MISE iff $X \xEd \xCQ $ and $ \xcA a \xbe A
\xcA x \xbe X(x \xeb a \xco x \xcT a$ $ \xch $ $x \xbe A).$
(A is downward and horizontally closed.)

(3) Theory Revision

Recall that we have a distance $d$ on the model set, and are interested
in $y \xbe Y$ which are close to $X.$

Thus, given $ \xCf X,Y,$ we define analogously:

$B \xcc Y$ is MISE iff

(1) $B \xEd \xCQ $

(2) there is $d' $ s.t. $B:=\{y \xbe Y: \xcE x \xbe X.d(x,y) \xck d' \}$
(we could also have chosen $d(x,y)<d',$ this is not important).

And we define $ \xbf \xbe T*T' $ iff there is $B \xbe \xbL (M(T),M(T' ))$
$B \xcm \xbf.$
\subsection{The algebraic limit}
\label{Section 2.2.11.2}
%  Section (8.2):  THE LIMIT VARIANT, Algebraic Limit
%  Section (8.2):  THE LIMIT VARIANT, Algebraic Limit
% %
% ====================================================

\ed

There are basic problems with the algebraic limit in general preferential
structures.

\be

$\hspace{0.01em}$

% (+++ Orig. No.:  Example D-8.2.1 +++)

\label{Example D-8.2.1}

Let $a \xeb b,$ $a \xeb c,$ $b \xeb d,$ $c \xeb d$ (but $ \xeb $ not
transitive!), then $\{a,b\}$ and $\{a,c\}$ are
such $S$ and $S',$ but there is no $S'' \xcc S \xcs S' $ which is an
initial segment. If, for
instance, in a and $b$ $ \xbq $ holds, in a and $c$ $ \xbq ',$ then ``in
the limit'' $ \xbq $ and $ \xbq ' $
will hold, but not $ \xbq \xcu \xbq '.$ This does not seem right. We
should not be obliged
to give up $ \xbq $ to obtain $ \xbq '.$ $ \xcz $
\\[3ex]

\ee

When we look at the system of such $S$ generated by a preferential
structure and
its algebraic properties, we will therefore require it to be closed under
finite intersections, or at least, that if $S,$ $S' $ are such segments,
then there
must be $S'' \xcc S \xcs S' $ which is also such a segment.

We make this official. Let $ \xbL (X)$ be the set of initial segments of
$X,$ then
we require:

$( \xbL \xcs )$ If $A,B \xbe \xbL (X)$ then there is $C \xcc A \xcs B,$ $C
\xbe \xbL (X).$

More precisely, a limit should be a structural limit in a reasonable sense
-
whatever the underlying structure is -, and the resulting algebraic limit
should respect $( \xbL \xcs ).$

We should not demand too much, either. It would be wrong to demand closure
under
arbitrary intersections, as this would mean that there is an initial
segment
which makes all consequences true - trivializing the very idea of a limit.

But we can make our requirements more precise, and bind the limit variant
closely to the minimal variant, by looking at the algebraic version of
both.

Before we look at deeper problems, we show some basic facts about the
algebraic
limit.

\bfa

$\hspace{0.01em}$

% (+++ Orig. No.:  Fact D-8.2.1 +++)

\label{Fact D-8.2.1}

(Taken from  \cite{Sch04}, Fact 3.4.3, Proposition 3.10.16 there.)

Let the relation $ \xeb $ be transitive. The following hold in the limit
variant of
general preferential structures:

(1) If $A \xbe \xbL (Y),$ and $A \xcc X \xcc Y,$ then $A \xbe \xbL (X).$

(2) If $A \xbe \xbL (Y),$ and $A \xcc X \xcc Y,$ and $B \xbe \xbL (X),$
then $A \xcs B \xbe \xbL (Y).$

(3) If $A \xbe \xbL (Y),$ $B \xbe \xbL (X),$ then there is $Z \xcc A \xcv
B$ $Z \xbe \xbL (Y \xcv X).$

The following hold in the limit variant of ranked structures without
copies,
where the domain is closed under finite unions and contains all finite
sets.

(4) $A,B \xbe \xbL (X)$ $ \xch $ $A \xcc B$ or $B \xcc A,$

(5) $A \xbe \xbL (X),$ $Y \xcc X,$ $Y \xcs A \xEd \xCQ $ $ \xch $ $Y \xcs
A \xbe \xbL (Y),$

(6) $ \xbL ' \xcc \xbL (X),$ $ \xcS \xbL ' \xEd \xCQ $ $ \xch $ $ \xcS
\xbL ' \xbe \xbL (X).$

(7) $X \xcc Y,$ $A \xbe \xbL (X)$ $ \xch $ $ \xcE B \xbe \xbL (Y).B \xcs
X=A$

\efa

\subparagraph{
Proof
}

$\hspace{0.01em}$

% (+++ Orig.:  Proof +++)

(1) trivial.

(2)

(2.1) $A \xcs B$ is closed in $Y:$ Let $ \xBc x,i \xBe  \xbe A \xcs B,$
$ \xBc y,j \xBe  \xeb  \xBc x,i \xBe,$ then $ \xBc y,j \xBe  \xbe A.$ If
$ \xBc y,j \xBe  \xce X,$ then $ \xBc y,j \xBe  \xce A,$ $contradiction.$
So $ \xBc y,j \xBe  \xbe X,$ but then $ \xBc y,j \xBe  \xbe B.$

(2.2) $A \xcs B$ minimizes $Y:$ Let $ \xBc a,i \xBe  \xbe Y.$

(a) If $ \xBc a,i \xBe  \xbe A-B \xcc X,$ then there is
$ \xBc y,j \xBe  \xeb  \xBc a,i \xBe,$ $ \xBc y,j \xBe  \xbe B.$ Xy closure of
A,
$ \xBc y,j \xBe  \xbe A.$

(b) If $ \xBc a,i \xBe  \xce A,$ then there is
$ \xBc a',i'  \xBe  \xbe A \xcc X,$ $ \xBc a',i'  \xBe  \xeb  \xBc a,i \xBe,$
continue by (a).

(3)

Let $Z$ $:=$ $\{ \xBc x,i \xBe  \xbe A$:
$ \xCN \xcE  \xBc b,j \xBe  \xec  \xBc x,i \xBe. \xBc b,j \xBe  \xbe X-B\}$ $
\xcv $
$\{ \xBc y,j \xBe  \xbe B$:
$ \xCN \xcE  \xBc a,i \xBe  \xec  \xBc y,j \xBe. \xBc a,i \xBe  \xbe Y-A\},$
where $ \xec $ stands for $ \xeb $ or $=.$

(3.1) $Z$ minimizes $Y \xcv X:$ We consider $Y,$ $X$ is symmetrical.

(a) We first show: If $ \xBc a,k \xBe  \xbe A-$Z, then there is
$ \xBc y,i \xBe  \xbe Z. \xBc a,k \xBe  \xee  \xBc y,i \xBe.$
Broof: If $ \xBc a,k \xBe  \xbe A-$Z,
then there is $ \xBc b,j \xBe  \xec  \xBc a,k \xBe,$ $ \xBc b,j \xBe  \xbe
X-$B. Then there
is $ \xBc y,i \xBe  \xeb  \xBc b,j \xBe,$ $ \xBc y,i \xBe  \xbe B.$
Xut $ \xBc y,i \xBe  \xbe Z,$ too: If not, there would be
$ \xBc a',k'  \xBe  \xec  \xBc y,i \xBe,$ $ \xBc a',k'  \xBe  \xbe Y- \xCf A,$
but
$ \xBc a',k'  \xBe  \xeb  \xBc a,k \xBe,$ contradicting closure of $ \xCf A.$

(b) If $ \xBc a'',k''  \xBe  \xbe Y- \xCf A,$ there is $ \xBc a,k \xBe  \xbe A,$
$ \xBc a,k \xBe  \xeb  \xBc a'',k''  \xBe.$ If $ \xBc a,k \xBe  \xce Z,$
continue
with (a).

(3.2) $Z$ is closed in $Y \xcv X:$ Let then $ \xBc z,i \xBe  \xbe Z,$
$ \xBc u,k \xBe  \xeb  \xBc z,i \xBe,$ $ \xBc u,k \xBe  \xbe Y \xcv X.$
Suppose $ \xBc z,i \xBe  \xbe A$ - the case
$ \xBc z,i \xBe  \xbe B$ is symmetrical.

(a) $ \xBc u,k \xBe  \xbe Y-A$ cannot be, by closure of $ \xCf A.$

(b) $ \xBc u,k \xBe  \xbe X-B$ cannot be, as $ \xBc z,i \xBe  \xbe Z,$ and by
definition of
$Z.$

(c) If $ \xBc u,k \xBe  \xbe A- \xCf Z,$ then there is
$ \xBc v,l \xBe  \xec  \xBc u,k \xBe,$ $ \xBc v,l \xBe  \xbe X-$B, so $ \xBc
v,l \xBe  \xeb  \xBc z,i \xBe,$
contradicting (b).

(d) If $ \xBc u,k \xBe  \xbe B-$Z, then there is
$ \xBc v,l \xBe  \xec  \xBc u,k \xBe,$ $ \xBc v,l \xBe  \xbe Y- \xCf A,$
contradicting (a).

(4) Suppose not, so there are $a \xbe A- \xCf B,$ $b \xbe B- \xCf A.$ But
if $a \xcT b,$ $a \xbe B$ and $b \xbe A,$
similarly if $a \xeb b$ or $b \xeb a.$

(5) As $A \xbe \xbL (X)$ and $Y \xcc X,$ $Y \xcs A$ is downward and
horizontally closed. As $Y \xcs A \xEd \xCQ,$
$Y \xcs A$ minimizes $Y.$

(6) $ \xcS \xbL ' $ is downward and horizontally closed, as all $A \xbe
\xbL ' $ are. As $ \xcS \xbL ' \xEd \xCQ,$
$ \xcS \xbL ' $ minimizes $X.$

(7) Set $B:=\{b \xbe Y: \xcE a \xbe A.a \xcT b$ or $b \xck a\}$

$ \xcz $
\\[3ex]

We have as immediate logical consequence:

\bfa

$\hspace{0.01em}$

% (+++ Orig. No.:  Fact D-8.2.2 +++)

\label{Fact D-8.2.2}

(Fact 3.4.4 of  \cite{Sch04}.)

If $ \xeb $ is transitive, then in the limit variant hold:

(1) $ \xCf (AND),$

(2) $ \xCf (OR).$

\efa

\subparagraph{
Proof
}

$\hspace{0.01em}$

% (+++ Orig.:  Proof +++)

Let $ \xdz $ be the structure.

(1) Immediate by Fact \ref{Fact D-8.2.1} (page \pageref{Fact D-8.2.1}), (2) -
set $A=B.$

(2) Immediate by Fact \ref{Fact D-8.2.1} (page \pageref{Fact D-8.2.1}), (3).
$ \xcz $
\\[3ex]
\subsection{The logical limit}
\label{Section 2.2.11.3}
%  Section (8.3):  THE LIMIT VARIANT, Logical Limit
%  Section (8.3):  THE LIMIT VARIANT, Logical Limit
% %
% ==================================================
\subsubsection{Translation between the minimal and the limit variant}

$\hspace{0.01em}$
%  Subsection (8.3.1):  LIMIT: Translation between minimal and limit variant
%  Subsection (8.3.1):  LIMIT: Translation between minimal and limit variant
% %
% ===========================================================================

A good example for problems linked to the translation from the algebraic
limit
to the logical limit is the property $( \xbm =)$ of ranked structures:

$( \xbm =)$ $X \xcc Y,$ $ \xbm (Y) \xcs X \xEd \xCQ $ $ \xch $ $ \xbm (Y)
\xcs X= \xbm (X)$

or its logical form

$( \xcn =)$ $T \xcl T',$ $Con( \ol{ \ol{T' } },T)$ $ \xch $ $ \ol{ \ol{T}
}= \ol{ \ol{ \ol{T' } } \xcv T}.$

$ \xbm (Y)$ or its analogue $ \ol{ \ol{T' } }$ (set $X:=M(T),$ $Y:=M(T'
))$ speak about the limit, the
``ideal'', and this, of course, is not what we have in the limit version.
This
limit version was intoduced precisely to avoid speaking about the ideal.

So, first, we have to translate $ \xbm (Y) \xcs X \xEd \xCQ $ to something
else, and the
natural candidate seems to be

$ \xcA B \xbe \xbL (Y).B \xcs X \xEd \xCQ.$

In logical terms, we have replaced the set of consequences of $Y$ by some
$Th(B)$
where $T' \xcc Th(B) \xcc \ol{ \ol{T' } }.$ The conclusion can now be
translated in a similar way to
$ \xcA B \xbe \xbL (Y). \xcE A \xbe \xbL (X).A \xcc B \xcs X$ and $ \xcA A
\xbe \xbL (X). \xcE B \xbe \xbL (Y).B \xcs X \xcc A.$ The total
translation reads now:

$( \xbL =)$ Let $X \xcc Y.$ Then

$ \xCI \xcA B \xbe \xbL (Y).B \xcs X \xEd \xCQ \xCJ $ $ \xch $
$ \xCI \xcA B \xbe \xbL (Y). \xcE A \xbe \xbL (X).A \xcc B \xcs X$ and $
\xcA A \xbe \xbL (X). \xcE B \xbe \xbL (Y).B \xcs X \xcc A \xCJ.$

By Fact \ref{Fact D-8.2.1} (page \pageref{Fact D-8.2.1})  (5) and (7),
we see that this holds in ranked
structures. Thus, the limit reading seems to provide a correct algebraic
limit.

Yet, Example \ref{Example D-8.3.1} (page \pageref{Example D-8.3.1})
below shows the following:

Let $m' \xEd m$ be arbitrary.
For $T':=Th(\{m,m' \}),$ $T:= \xCQ,$ we have $T' \xcl T,$ $ \ol{ \ol{T'
} }=Th(\{m' \}),$ $ \ol{ \ol{T} }=Th(\{m\}),$ $Con( \ol{ \ol{T} },T' ),$
but $Th(\{m\})= \ol{ \ol{ \ol{T} } \xcv T' } \xEd \ol{ \ol{T' } }.$

Thus:

(1) The prerequisite holds, though usually for $A \xbe \xbL (T),$ $A \xcs
M(T' )= \xCQ.$

(2) (PR) fails, which is independent of the prerequisite $Con( \ol{ \ol{T}
},T' ),$ so the
problem is not just due to the prerequsite.

(3) Both inclusions of $( \xcn =)$ fail.

We will see below in
Corollary \ref{Corollary D-8.3.4} (page \pageref{Corollary D-8.3.4})
a sufficient condition to make $( \xcn =)$
hold in ranked structures. It has to do with definability or formulas,
more
precisely, the crucial property is to have sufficiently often
$ \wt{A} \xcs \wt{M(T' )}= \wt{A \xcs M(T' )}$
for $A \xbe \xbL (T)$ - see
Section \ref{Section 2.2.10.1} (page \pageref{Section 2.2.10.1})  for reference.

\be

$\hspace{0.01em}$

% (+++ Orig. No.:  Example D-8.3.1 +++)

\label{Example D-8.3.1}

(Taken from  \cite{Sch04}, Example 3.10.1 (1) there.)

Take an infinite propositional language $p_{i}:i \xbe \xbo.$ We have $
\xbo_{1}$ models (assume
for simplicity CH).

Take the model $m$ which makes all $p_{i}$ true, and put it on top. Next,
going down, take
all models which make $p_{0}$ false, and then all models which make
$p_{0}$ true, but $p_{1}$
false, etc. in a ranked construction.
So, successively more $p_{i}$ will become (and stay) true. Consequently,
$ \xCQ \xcm_{ \xbL }p_{i}$ for all $i.$ But the structure has no minimum,
and the ``logical''
limit $m$ is not in the set wise limit.
Let $T:= \xCQ $ and $m' \xEd m,$ $T':=Th(\{m,m' \}),$ then $ \ol{ \ol{T}
}=Th(\{m\}),$ $ \ol{ \ol{T' } }=Th(\{m' \}),$ and
$ \ol{ \ol{ \ol{T' } } \xcv T}= \ol{ \ol{T' } }=Th(\{m' \})$ and $ \ol{
\ol{ \ol{T} } \xcv T' }= \ol{ \ol{T} }=Th(\{m\}).$ $ \xcz $
\\[3ex]

\ee

This example shows that our translation is not perfect, but it is half the
way. Note that the minimal variant faces the same problems (definability
and
others), so the problems are probably at least not totally due to our
perhaps insufficient translation.

We turn to other rules.

$( \xbL \xcs )$ If $A,B \xbe \xbL (X)$ then there is $C \xcc A \xcs B,$ $C
\xbe \xbL (X)$

seems a minimal requirement for an appropriate limit. It holds in
transitive structures by
Fact \ref{Fact D-8.2.1} (page \pageref{Fact D-8.2.1})  (2).

The central logical condition for minimal smooth structures is

(CUM) $T \xcc T' \xcc \ol{ \ol{T} }$ $ \xch $ $ \ol{ \ol{T} }= \ol{ \ol{T'
} }$

It would again be wrong - using the limit - to translate this only partly
by:
If $T \xcc T' \xcc \ol{ \ol{T} },$ then for all $A \xbe \xbL (M(T))$ there
is $B \xbe \xbL (M(T' ))$ s.t. $A \xcc B$ - and
vice versa.
Now, smoothness is in itself a wrong condition for limit structures, as it
speaks about minimal elements, which we will not necessarily have. This
cannot
guide us. But when we consider a more modest version of cumulativity, we
see
what to do.

(CUMfin) If $T \xcn \xbf,$ then $ \ol{ \ol{T} }= \ol{ \ol{T \xcv \{ \xbf
\}} }.$

This translates into algebraic limit conditions as follows - where
$Y=M(T),$ and
$X=M(T \xcv \{ \xbf \}):$

$( \xbL CUMfin)$ Let $X \xcc Y.$ If there is $B \xbe \xbL (Y)$ s.t. $B
\xcc X,$ then:

$ \xCI \xcA A \xbe \xbL (X) \xcE B' \xbe \xbL (Y).B' \xcc A$ and $ \xcA B'
\xbe \xbL (Y) \xcE A \xbe \xbL (X).A \xcc B' \xCJ.$

Note, that in this version, we do not have the ``ideal'' limit on the left
of the
implication, but one fixed approximation $B \xbe \xbL (Y).$
We can now prove that $( \xbL CUMfin)$ holds in transitive structures:
The first part holds by
Fact \ref{Fact D-8.2.1} (page \pageref{Fact D-8.2.1})  (2),
the second, as $B \xcs B' \xbe \xbL (Y)$ by
Fact \ref{Fact D-8.2.1} (page \pageref{Fact D-8.2.1})  (1).
This is true without additional properties of the structure,
which might at first sight seem surprising. But note that the initial
segments
play a similar role as the set of minimal elements: an initial segment has
to
minimize the other elements, just as the set of minimal elements in the
smooth
case does.

The central algebraic property of minimal preferential structures is

$( \xbm PR)$ $X \xcc Y$ $ \xch $ $ \xbm (Y) \xcs X \xcc \xbm (X)$

This translates naturally and directly to

$( \xbL PR)$ $X \xcc Y$ $ \xch $ $ \xcA A \xbe \xbL (X) \xcE B \xbe \xbL
(Y).B \xcs X \xcc A$

$( \xbL PR)$ holds in transitive structures:
$Y-X \xbe \xbL (Y- \xCf X),$ so the result holds by
Fact \ref{Fact D-8.2.1} (page \pageref{Fact D-8.2.1})  (3).

The central algebraic condition of ranked minimal structures is

$( \xbm =)$ $X \xcc Y,$ $ \xbm (Y) \xcs X \xEd \xCQ $ $ \xch $ $ \xbm (Y)
\xcs X= \xbm (X)$

We saw above how to translate this condition to $( \xbL =),$ we also saw
that $( \xbL =)$
holds
in ranked structures.

We will see in
Corollary \ref{Corollary D-8.3.4} (page \pageref{Corollary D-8.3.4})
that the following logical version holds
in ranked structures:

$T \xcN \xCN \xbg $ implies $ \ol{ \ol{T} }= \ol{ \ol{T \xcv \{ \xbg \}}
}$

We generalize above translation results to a recipe:

Translate

 \xEh

 \xDH $ \xbm (X) \xcc \xbm (Y)$ to $ \xcA B \xbe \xbL (Y) \xcE A \xbe \xbL
(X).A \xcc B,$
and thus

 \xDH $ \xbm (Y) \xcs X \xcc \xbm (X)$ to $ \xcA A \xbe \xbL (X) \xcE B
\xbe \xbL (Y).B \xcs X \xcc A,$

 \xDH $ \xbm (X) \xcc Y$ to $ \xcE A \xbe \xbL (X).A \xcc Y,$
and thus

 \xDH $ \xbm (Y) \xcs X \xEd \xCQ $ to $ \xcA B \xbe \xbL (Y).B \xcs X
\xEd \xCQ $

 \xDH $X \xcc \xbm (Y)$ to $ \xcA B \xbe \xbL (Y).X \xcc B,$

and quantify expressions separately, thus we repeat:

 \xDH $( \xbm CUM)$ $ \xbm (Y) \xcc X \xcc Y$ $ \xch $ $ \xbm (X)= \xbm
(Y)$ translates to

 \xDH $( \xbL CUMfin)$ Let $X \xcc Y.$ If there is $B \xbe \xbL (Y)$ s.t.
$B \xcc X,$ then:

$ \xCI \xcA A \xbe \xbL (X) \xcE B' \xbe \xbL (Y).B' \xcc A$ and $ \xcA B'
\xbe \xbL (Y) \xcE A \xbe \xbL (X).A \xcc B' \xCJ.$

 \xDH $( \xbm =)$ $X \xcc Y,$ $ \xbm (Y) \xcs X \xEd \xCQ $ $ \xch $ $
\xbm (Y) \xcs X= \xbm (X)$ translates to

 \xDH $( \xbL =)$ Let $X \xcc Y.$ If $ \xcA B \xbe \xbL (Y).B \xcs X \xEd
\xCQ,$ then

$ \xCI \xcA A \xbe \xbL (X) \xcE B' \xbe \xbL (Y).B' \xcs X \xcc A,$ and $
\xcA B' \xbe \xbL (Y) \xcE A \xbe \xbL (X).A \xcc B' \xcs X \xCJ.$

 \xEj

We collect now for easier reference the definitions and some algebraic
properties which we saw above to hold:

\bd

$\hspace{0.01em}$

% (+++ Orig. No.:  Definition D-8.3.1 +++)

\label{Definition D-8.3.1}

$( \xbL \xcs )$ If $A,B \xbe \xbL (X)$ then there is $C \xcc A \xcs B,$ $C
\xbe \xbL (X),$

$( \xbL PR)$ $X \xcc Y$ $ \xch $ $ \xcA A \xbe \xbL (X) \xcE B \xbe \xbL
(Y).B \xcs X \xcc A,$

$( \xbL CUMfin)$ Let $X \xcc Y.$ If there is $B \xbe \xbL (Y)$ s.t. $B
\xcc X,$ then:

$ \xCI \xcA A \xbe \xbL (X) \xcE B' \xbe \xbL (Y).B' \xcc A$ and $ \xcA B'
\xbe \xbL (Y) \xcE A \xbe \xbL (X).A \xcc B' \xCJ.$

 \xDH $( \xbL =)$ Let $X \xcc Y.$ If $ \xcA B \xbe \xbL (Y).B \xcs X \xEd
\xCQ,$ then

$ \xCI \xcA A \xbe \xbL (X) \xcE B' \xbe \xbL (Y).B' \xcs X \xcc A,$ and $
\xcA B' \xbe \xbL (Y) \xcE A \xbe \xbL (X).A \xcc B' \xcs X \xCJ.$

\ed

\bfa

$\hspace{0.01em}$

% (+++ Orig. No.:  Fact D-8.3.1 +++)

\label{Fact D-8.3.1}

In transitive structures hold:

(1) $( \xbL \xcs )$

(2) $( \xbL PR)$

(3) $( \xbL CUMfin)$

In ranked structures holds:

(4) $( \xbL =)$

\efa

\subparagraph{
Proof
}

$\hspace{0.01em}$

% (+++ Orig.:  Proof +++)

(1) By Fact \ref{Fact D-8.2.1} (page \pageref{Fact D-8.2.1})  (2).

(2) $Y-X \xbe \xbL (Y- \xCf X),$ so the result holds by
Fact \ref{Fact D-8.2.1} (page \pageref{Fact D-8.2.1})  (3).

(3) By Fact \ref{Fact D-8.2.1} (page \pageref{Fact D-8.2.1})  (1) and (2).

(4) By Fact \ref{Fact D-8.2.1} (page \pageref{Fact D-8.2.1})  (5) and (7).

$ \xcz $
\\[3ex]

To summarize the discussion:

Just as in the minimal case, the algebraic laws may hold, but not the
logical
ones, due in both cases to definability problems. Thus, we cannot expect a
clean
proof of correspondence. But we can argue that we did a correct
translation, which shows its limitation, too. The part with $ \xbm (X)$
and $ \xbm (Y)$ on
both sides of $ \xcc $ is obvious, we will have a perfect correspondence.
The part
with $X \xcc \xbm (Y)$ is obvious, too. The problem is in the part with $
\xbm (X) \xcc Y.$ As we
cannot use the limit, but only its approximation, we are limited here to
one
(or finitely many) consequences of $T,$ if $X=M(T),$ so we obtain only $T
\xcn \xbf,$ if
$Y \xcc M( \xbf ),$ and if there is $A \xbe \xbL (X).A \xcc Y.$

We consider a limit only appropriate, if it is an algebraic limit which
preserves algebraic properties of the minimal version in above
translation.

The advantage of such limits is that they allow - with suitable caveats -
to
show that they preserve the logical properties of the minimal variant, and
thus are equivalent to the minimal case (with, of course, perhaps a
different
relation). Thus, they allow a straightforward trivialization.
\subsubsection{Logical properties of the limit variant}

$\hspace{0.01em}$
%  Subsection (8.3.2):  LIMIT: Logical properties
%  Subsection (8.3.2):  LIMIT: Logical properties
% %
% ================================================

We begin with some simple logical facts about the limit version.

We abbreviate $ \xbL (T):= \xbL (M(T))$ etc., assume transitivity.

\bfa

$\hspace{0.01em}$

% (+++ Orig. No.:  Fact D-8.3.2 +++)

\label{Fact D-8.3.2}

(1) $A \xbe \xbL (T)$ $ \xch $ $M( \ol{ \ol{T} }) \xcc \wt{A}$

(2) $M( \ol{ \ol{T} })$ $=$ $ \xcS \{ \wt{A}:A \xbe \xbL (T)\}$

(2a) $M( \ol{ \ol{T' } }) \xcm \xbs $ $ \xch $ $ \xcE B \xbe \xbL (T' ).
\wt{B} \xcm \xbs $

(3) $M( \ol{ \ol{T' } }) \xcs M(T) \xcm \xbs $ $ \xch $ $ \xcE B \xbe \xbL
(T' ). \wt{B} \xcs M(T) \xcm \xbs.$

\efa

\subparagraph{
Proof
}

$\hspace{0.01em}$

% (+++ Orig.:  Proof +++)

(1)

Note that $A \xcm \xbf $ $ \xch $ $T \xcn \xbf $ by definition, see
Definition \ref{Definition D-8.1.1} (page \pageref{Definition D-8.1.1}).

Let $M( \ol{ \ol{T} }) \xcC \wt{A},$ so there is $ \xbf,$ $ \wt{A} \xcm
\xbf,$ so $A \xcm \xbf,$ but $M( \ol{ \ol{T} }) \xcM \xbf,$ so $T \xcN
\xbf,$ $contradiction.$

(2) `` $ \xcc $ '' by (1). `` $ \xcd $ '': Let $x \xbe \xcS \{ \wt{A}:A \xbe
\xbL (T)\}$ $ \xch $ $ \xcA A \xbe \xbL (T).x \xcm Th(A)$ $ \xch $ $x \xcm
\ol{ \ol{T} }.$

(2a) $M( \ol{ \ol{T' } }) \xcm \xbs $ $ \xch $ $T' \xcn \xbs $ $ \xch $ $
\xcE B \xbe \xbL (T' ).B \xcm \xbs.$ But $B \xcm \xbs $ $ \xch $ $ \wt{B}
\xcm \xbs.$

(3) $M( \ol{ \ol{T' } }) \xcs M(T) \xcm \xbs $ $ \xch $ $ \ol{ \ol{T' } }
\xcv T \xcl \xbs $ $ \xch $ $ \xcE \xbt_{1} \Xl  \xbt_{n} \xbe \ol{ \ol{T'
} }$ s.t. $T \xcv \{ \xbt_{1}, \Xl, \xbt_{n}\} \xcl
 \xbs,$
so $ \xcE B \xbe \xbL (T' ).Th(B) \xcv T \xcl \xbs.$ So $M(Th(B)) \xcs
M(T) \xcm \xbs $ $ \xch $ $ \wt{B} \xcs M(T) \xcm \xbs.$

$ \xcz $
\\[3ex]

We saw in Example \ref{Example D-8.3.1} (page \pageref{Example D-8.3.1})
and its discussion the problems which
might arise in the limit version, even if the algebraic behaviour is
correct.

This analysis leads us to consider the following facts:

\bfa

$\hspace{0.01em}$

% (+++ Orig. No.:  Fact D-8.3.3 +++)

\label{Fact D-8.3.3}

(1) Let $ \xcA B \xbe \xbL (T' ) \xcE A \xbe \xbL (T).A \xcc B \xcs M(T),$
then $ \ol{ \ol{ \ol{T' } } \xcv T} \xcc \ol{ \ol{T} }.$

Let, in addition, $\{B \xbe \xbL (T' ): \wt{B} \xcs \wt{M(T)}= \wt{B \xcs
M(T)}\}$ be cofinal in $ \xbL (T' ).$ Then

(2) $Con( \ol{ \ol{T' } },T)$ implies $ \xcA A \xbe \xbL (T' ).A \xcs M(T)
\xEd \xCQ.$

(3) $ \xcA A \xbe \xbL (T) \xcE B \xbe \xbL (T' ).B \xcs M(T) \xcc A$
implies $ \ol{ \ol{T} } \xcc \ol{ \ol{ \ol{T' } } \xcv T}.$

\efa

Note that $M(T)= \wt{M(T)},$ so we could also have written $ \wt{B} \xcs
M(T)= \wt{B \xcs M(T)},$ but above
way of writing stresses more the essential condition $ \wt{X} \xcs \wt{Y}=
\wt{X \xcs Y}.$

\subparagraph{
Proof
}

$\hspace{0.01em}$

% (+++ Orig.:  Proof +++)

(1)
Let $ \ol{ \ol{T' } } \xcv T \xcl \xbs,$ so $ \xcE B \xbe \xbL (T' ).
\wt{B} \xcs M(T) \xcm \xbs $ by
Fact \ref{Fact D-8.3.2} (page \pageref{Fact D-8.3.2}), (3) above (using
compactness). Thus $ \xcE A \xbe \xbL (T).A \xcc B \xcs M(T) \xcm \xbs $
by prerequisite, so $ \xbs \xbe \ol{ \ol{T} }.$

(2)
Let $Con( \ol{ \ol{T' } },T),$ so $M( \ol{ \ol{T' } }) \xcs M(T) \xEd \xCQ
.$ $M( \ol{ \ol{T' } })= \xcS \{ \wt{A}:A \xbe \xbL (T' )\}$ by
Fact \ref{Fact D-8.3.2} (page \pageref{Fact D-8.3.2})  (2), so
$ \xcA A \xbe \xbL (T' ). \wt{A} \xcs M(T) \xEd \xCQ.$ As cofinally often
$ \wt{A} \xcs M(T)= \wt{A \xcs M(T)},$
$ \xcA A \xbe \xbL (T' ). \wt{A \xcs M(T)} \xEd \xCQ,$ so $ \xcA A \xbe
\xbL (T' ).A \xcs M(T) \xEd \xCQ $ by $ \wt{ \xCQ }= \xCQ.$

(3)
Let $ \xbs \xbe \ol{ \ol{T} },$ so $T \xcn \xbs,$ so $ \xcE A \xbe \xbL
(T).A \xcm \xbs,$ so $ \xcE B \xbe \xbL (T' ).B \xcs M(T) \xcc A$ by
prerequisite, so $ \xcE B \xbe \xbL (T' ). \xCI B \xcs M(T) \xcc A$ and $
\wt{B} \xcs \wt{M(T)}= \wt{B \xcs M(T)} \xCJ.$ So for
such $B$ $ \wt{B} \xcs \wt{M(T)}= \wt{B \xcs M(T)} \xcc \wt{A} \xcm \xbs
.$ By
Fact \ref{Fact D-8.3.2} (page \pageref{Fact D-8.3.2})  (1)
$M( \ol{ \ol{T' } }) \xcc \wt{B},$ so $M( \ol{ \ol{T' } }) \xcs M(T) \xcm
\xbs,$ so $ \ol{ \ol{T' } } \xcv T \xcl \xbs.$

$ \xcz $
\\[3ex]

We obtain now as easy corollaries of a more general situation the
following
properties shown in  \cite{Sch04} by direct proofs. Thus, we have
the trivialization
results shown there.

\bco

$\hspace{0.01em}$

% (+++ Orig. No.:  Corollary D-8.3.4 +++)

\label{Corollary D-8.3.4}

Let the structure be transitive.

(1) Let $\{B \xbe \xbL (T' ): \wt{B} \xcs \wt{M(T)}= \wt{B \xcs M(T)}\}$
be cofinal in $ \xbL (T' ),$ then

(PR) $T \xcl T' $ $ \xch $ $ \ol{ \ol{T} }$ $ \xcc $ $ \ol{ \ol{ \ol{T' }
} \xcv T}$ holds.

(2) $ \ol{ \ol{ \xbf \xcu \xbf ' } } \xcc \ol{ \ol{ \ol{ \xbf } } \xcv \{
\xbf ' \}}$ holds.

If the structure is ranked, then also:

(3) Let $\{B \xbe \xbL (T' ): \wt{B} \xcs \wt{M(T)}= \wt{B \xcs M(T)}\}$
be cofinal in $ \xbL (T' ),$ then

$( \xcn =)$ $T \xcl T',$ $Con( \ol{ \ol{T' } },T)$ $ \xch $ $ \ol{ \ol{T}
}$ $=$ $ \ol{ \ol{ \ol{T' } } \xcv T}$ holds.

(4) $T \xcN \xCN \xbg $ $ \xch $ $ \ol{ \ol{T} }= \ol{ \ol{T \xcv \{ \xbg
\}} }$ holds.

\eco

\subparagraph{
Proof
}

$\hspace{0.01em}$

% (+++ Orig.:  Proof +++)

(1) $ \xcA A \xbe \xbL (M(T)) \xcE B \xbe \xbL (M(T' )).B \xcs M(T) \xcc
A$ by
Fact \ref{Fact D-8.3.1} (page \pageref{Fact D-8.3.1})  (2).
So the result follows from
Fact \ref{Fact D-8.3.3} (page \pageref{Fact D-8.3.3})  (3).

(2) Set $T':=\{ \xbf \},$ $T:=\{ \xbf, \xbf ' \}.$ Then for $B \xbe \xbL
(T' )$ $ \wt{B} \xcs M(T)= \wt{B} \xcs M( \xbf ' )= \wt{B \xcs M( \xbf '
)}$
by Fact \ref{Fact Def-Clos} (page \pageref{Fact Def-Clos})  $(Cl \xcs +),$
so the result follows by (1).

(3) Let $Con( \ol{ \ol{T' } },T),$ then by
Fact \ref{Fact D-8.3.3} (page \pageref{Fact D-8.3.3})  (2)
$ \xcA A \xbe \xbL (T' ).A \xcs M(T) \xEd \xCQ,$ so by
Fact \ref{Fact D-8.3.1} (page \pageref{Fact D-8.3.1})  (4)
$ \xcA B \xbe \xbL (T' ) \xcE A \xbe \xbL (T).A \xcc B \xcs M(T),$ so $
\ol{ \ol{ \ol{T' } } \xcv T} \xcc \ol{ \ol{T} }$
by Fact \ref{Fact D-8.3.3} (page \pageref{Fact D-8.3.3})  (1).

The other direction follows from (1).

(4) Set $T:=T' \xcv \{ \xbg \}.$ Then for $B \xbe \xbL (T' )$ $ \wt{B}
\xcs M(T)= \wt{B} \xcs M( \xbg )= \wt{B \xcs M( \xbg )}$ again
by Fact \ref{Fact Def-Clos} (page \pageref{Fact Def-Clos})  $(Cl \xcs +),$
so the result follows from (3).

$ \xcz $
\\[3ex]

We summarize for easier reference here our main positive logical results
on the
limit variant of general preferential structures where each model occurs
in one
copy only, Proposition 3.4.7 and Proposition 3.10.19 from  \cite{Sch04}:

\bp

$\hspace{0.01em}$

% (+++ Orig. No.:  Proposition D-8.3.5 +++)

\label{Proposition D-8.3.5}

Let the relation be transitive. Then

(1) Every instance of the the limit version, where the definable closed
minimizing sets are cofinal in the closed minimizing sets, is equivalent
to
an instance of the minimal version.

(2) If we consider only formulas on the left of $ \xcn,$ the resulting
logic of the
limit version can also be generated by the minimal version of a (perhaps
different) preferential structure. Moreover, the structure can be chosen
smooth.

\ep

\bp

$\hspace{0.01em}$

% (+++ Orig. No.:  Proposition D-8.3.6 +++)

\label{Proposition D-8.3.6}

When considering just formulas, in the ranked case without copies, $ \xbL
$ is
equivalent to $ \xbm $ - so $ \xbL $ is trivialized in this case. More
precisely:

Let a logic $ \xbf \xcn \xbq $ be given by the limit variant without
copies, i.e. by
Definition \ref{Definition D-8.1.1} (page \pageref{Definition D-8.1.1}).
Then there is a ranked structure, which gives
exactly the same logic, but interpreted in the minimal variant.

(As Example 3.10.2 in  \cite{Sch04} has shown, this is NOT
necessarily true if we
consider full theories $T$ and $T \xcn \xbq.)$

\ep

This shows that there is an important difference between considering full
theories and considering just formulas (on the left of $ \xcn ).$ If we
consider full
theories, we can ``grab'' single models, and thus determine the full order.
As
long as we restrict ourselves to formulas, we are much more shortsighted,
and see only a blurred picture. In particular, we can make sequences of
models
to converge to some model, but put this model elsewhere. Suitable such
manipulations will pass unobserved by formulas. The example also shows
that
there are structures whose limit version for theories is unequal to any
minimal structure.

(The negative results for the general not definability preserving minimal
case
apply also to the general limit case - see Section 5.2.3 in
 \cite{Sch04} for details.)
\chapter{Higher preferential structures}
\section{
Introduction
}
\label{Section Reac-GenPref-Intro}
\index{Section Reac-GenPref-Intro}
\index{Definition Generalized preferential structure}

\bd

$\hspace{0.01em}$

% (+++ Orig. No.:  Definition Generalized preferential structure +++)

\label{Definition Generalized preferential structure}

An IBR is called a generalized preferential structure iff the origins of
all
arrows are points. We will usually write $x,y$ etc. for points, $ \xba,$
$ \xbb $ etc. for
arrows.
\index{Definition Level-n-Arrow}

\ed

\bd

$\hspace{0.01em}$

% (+++ Orig. No.:  Definition Level-n-Arrow +++)

\label{Definition Level-n-Arrow}

Consider a generalized preferential structure $ \xdx.$

(1) Level $n$ arrow:

Definition by upward induction.

If $ \xba:x \xcp y,$ $x,y$ are points, then $ \xba $ is a level 1 arrow.

If $ \xba:x \xcp \xbb,$ $x$ is a point, $ \xbb $ a level $n$ arrow, then
$ \xba $ is a level $n+1$ arrow.
$(o( \xba )$ is the origin, $d( \xba )$ is the destination of $ \xba.)$

$ \xbl ( \xba )$ will denote the level of $ \xba.$

(2) Level $n$ structure:

$ \xdx $ is a level $n$ structure iff all arrows in $ \xdx $ are at most
level $n$ arrows.

We consider here only structures of some arbitrary but finite level $n.$

(3) We define for an arrow $ \xba $ by induction $O( \xba )$ and $D( \xba
).$

If $ \xbl ( \xba )=1,$ then $O( \xba ):=\{o( \xba )\},$ $D( \xba ):=\{d(
\xba )\}.$

If $ \xba:x \xcp \xbb,$ then $D( \xba ):=D( \xbb ),$ and $O( \xba
):=\{x\} \xcv O( \xbb ).$

Thus, for example, if $ \xba:x \xcp y,$ $ \xbb:z \xcp \xba,$ then $O(
\xbb ):=\{x,z\},$ $D( \xbb )=\{y\}.$
\index{Comment Gen-Pref}

\ed

\bcom

$\hspace{0.01em}$

% (+++ Orig. No.:  Comment Gen-Pref +++)

\label{Comment Gen-Pref}

A counterargument to $ \xba $ is NOT an argument for $ \xCN \xba $ (this
is asking for too
much), but just showing one case where $ \xCN \xba $ holds. In
preferential structures,
an argument for $ \xba $ is a set of level 1 arrows, eliminating $ \xCN
\xba -$models. A
counterargument is one level 2 arrow, attacking one such level 1 arrow.

Of course, when we have copies, we may need many successful attacks, on
all
copies, to achieve the goal. As we may have copies of level 1 arrows, we
may need many level 2 arrows to destroy them all.
\index{Example Inf-Level}

\ecom

We will not consider here diagrams with arbitrarily high levels. One
reason is
that diagrams like the following will have an unclear meaning:

\be

$\hspace{0.01em}$

% (+++ Orig. No.:  Example Inf-Level +++)

\label{Example Inf-Level}

$ \xBc  \xba,1 \xBe :x \xcp y,$

$ \xBc  \xba,n+1 \xBe :x \xcp  \xBc  \xba,n \xBe $ $(n \xbe \xbo ).$

Is $y \xbe \xbm (X)?$
\index{Definition Valid-Arrow}

\ee

\bd

$\hspace{0.01em}$

% (+++ Orig. No.:  Definition Valid-Arrow +++)

\label{Definition Valid-Arrow}

Let $ \xdx $ be a generalized preferential structure of (finite) level
$n.$

We define (by downward induction):

(1) Valid $X-to-Y$ arrow:

Let $X,Y \xcc \xdP ( \xdx ).$

$ \xba \xbe \xdA ( \xdx )$ is a valid $X-to-Y$ arrow iff

(1.1) $O( \xba ) \xcc X,$ $D( \xba ) \xcc Y,$

(1.2) $ \xcA \xbb:x' \xcp \xba.(x' \xbe X$ $ \xch $ $ \xcE \xbg:x''
\xcp \xbb.( \xbg $ is a valid $X-to-Y$ arrow)).

We will also say that $ \xba $ is a valid arrow in $X,$ or just valid in
$X,$ iff $ \xba $ is a
valid $X-to-X$ arrow.

(2) Valid $X \xch Y$ arrow:

Let $X \xcc Y \xcc \xdP ( \xdx ).$

$ \xba \xbe \xdA ( \xdx )$ is a valid $X \xch Y$ arrow iff

(2.1) $o( \xba ) \xbe X,$ $O( \xba ) \xcc Y,$ $D( \xba ) \xcc Y,$

(2.2) $ \xcA \xbb:x' \xcp \xba.(x' \xbe Y$ $ \xch $ $ \xcE \xbg:x''
\xcp \xbb.( \xbg $ is a valid $X \xch Y$ arrow)).

(Note that in particular $o( \xbg ) \xbe X,$ and that $o( \xbb )$ need not
be in $X,$ but
can be in the bigger $Y.)$
\index{Fact Higher-Validity}

\ed

\bfa

$\hspace{0.01em}$

% (+++ Orig. No.:  Fact Higher-Validity +++)

\label{Fact Higher-Validity}

(1) If $ \xba $ is a valid $X \xch Y$ arrow, then $ \xba $ is a valid
$Y-to-Y$ arrow.

(2) If $X \xcc X' \xcc Y' \xcc Y \xcc \xdP ( \xdx )$ and $ \xba \xbe \xdA
( \xdx )$ is a valid $X \xch Y$ arrow, and
$O( \xba ) \xcc Y',$ $D( \xba ) \xcc Y',$ then $ \xba $ is a valid $X'
\xch Y' $ arrow.
\index{Fact Higher-Validity Proof}

\efa

\subparagraph{
Proof
}

$\hspace{0.01em}$

% (+++ Orig.:  Proof +++)

Let $ \xba $ be a valid $X \xch Y$ arrow. We show (1) and (2) together by
downward induction
(both are trivial).

By prerequisite $o( \xba ) \xbe X \xcc X',$ $O( \xba ) \xcc Y' \xcc Y,$
$D( \xba ) \xcc Y' \xcc Y.$

Case 1: $ \xbl ( \xba )=n.$ So $ \xba $ is a valid $X' \xch Y' $ arrow,
and a valid $Y-to-Y$ arrow.

Case 2: $ \xbl ( \xba )=n-1.$ So there is no $ \xbb:x' \xcp \xba,$ $y
\xbe Y,$ so $ \xba $ is a valid
$Y-to-Y$ arrow. By $Y' \xcc Y$ $ \xba $ is a valid $X' \xch Y' $ arrow.

Case 3: Let the result be shown down to $m,$ $n>m>1,$ let $ \xbl ( \xba
)=m-1.$
So $ \xcA \xbb:x' \xcp \xba (x' \xbe Y$ $ \xch $ $ \xcE \xbg:x'' \xcp
\xbb (x'' \xbe X$ and $ \xbg $ is a valid $X \xch Y$ arrow)).
By induction hypothesis $ \xbg $ is a valid $Y-to-Y$ arrow, and a valid
$X' \xch Y' $ arrow. So $ \xba $ is a valid $Y-to-Y$ arrow, and by $Y'
\xcc Y,$ $ \xba $ is a valid
$X' \xch Y' $ arrow.

$ \xcz $
\\[3ex]
\index{Definition Higher-Mu}

\bd

$\hspace{0.01em}$

% (+++ Orig. No.:  Definition Higher-Mu +++)

\label{Definition Higher-Mu}

Let $ \xdx $ be a generalized preferential structure of level $n,$ $X \xcc
\xdP ( \xdx ).$

$ \xbm (X):=\{x \xbe X:$ $ \xcE  \xBc x,i \xBe. \xCN \xcE $ valid $X-to-X$
arrow $
\xba:x' \xcp  \xBc x,i \xBe \}.$
\index{Comment Smooth-Gen}

\ed

\bcom

$\hspace{0.01em}$

% (+++ Orig. No.:  Comment Smooth-Gen +++)

\label{Comment Smooth-Gen}

The purpose of smoothness is to guarantee cumulativity. Smoothness
achieves Cumulativity by mirroring all information present in $X$ also in
$ \xbm (X).$
Closer inspection shows that smoothness does more than necessary. This is
visible when there are copies (or, equivalently, non-injective labelling
functions). Suppose we have two copies of $x \xbe X,$
$ \xBc x,i \xBe $ and $ \xBc x,i'  \xBe,$ and there
is $y \xbe X,$ $ \xba: \xBc y,j \xBe  \xcp  \xBc x,i \xBe,$ but there is no
$ \xba ': \xBc y',j'  \xBe  \xcp  \xBc x,i'  \xBe,$ $y' \xbe X.$ Then
$ \xba: \xBc y,j \xBe  \xcp  \xBc x,i \xBe $
is irrelevant, as $x \xbe \xbm (X)$ anyhow. So mirroring
$ \xba: \xBc y,j \xBe  \xcp  \xBc x,i \xBe $ in $ \xbm (X)$ is not
necessary, i.e. it is not necessary to have some
$ \xba ': \xBc y',j'  \xBe  \xcp  \xBc x,i \xBe,$ $y' \xbe \xbm (X).$

On the other hand, Example \ref{Example Need-Smooth} (page \pageref{Example
Need-Smooth})  shows that,
if we want smooth structures to correspond to the property $( \xbm CUM),$
we
need at least some valid arrows from $ \xbm (X)$ also for higher level
arrows.
This ``some'' is made precise (essentially) in Definition \ref{Definition
X-Sub-X'} (page \pageref{Definition X-Sub-X'}).

From a more philosophical point of view,
when we see the (inverted) arrows of preferential structures as attacks on
non-minimal elements, then we should see smooth structures as always
having
attacks also from valid (minimal) elements. So, in general structures,
also
attacks from non-valid elements are valid, in smooth structures we always
also have attacks from valid elements.

The analogon to usual smooth structures, on level 2, is then that any
successfully attacked level 1 arrow is also attacked from a minimal point.
\index{Definition X-Sub-X'}

\ecom

\bd

$\hspace{0.01em}$

% (+++ Orig. No.:  Definition X-Sub-X' +++)

\label{Definition X-Sub-X'}

Let $ \xdx $ be a generalized preferential structure.

$X \xes X' $ iff

(1) $X \xcc X' \xcc \xdP ( \xdx ),$

(2) $ \xcA x \xbe X' -X$ $ \xcA  \xBc x,i \xBe $
$ \xcE \xba:x' \xcp  \xBc x,i \xBe ( \xba $ is a valid $X \xch X' $ arrow),

(3) $ \xcA x \xbe X$ $ \xcE  \xBc x,i \xBe $

$ \xDC $ $( \xcA \xba:x' \xcp  \xBc x,i \xBe (x' \xbe X' $ $ \xch $ $ \xcE \xbb
:x'' \xcp \xba.( \xbb $ is a valid $X \xch X' $ arrow))).

Note that (3) is not simply the negation of (2):

Consider a level 1 structure. Thus all level 1 arrows are valid, but the
source
of the arrows must not be neglected.

(2) reads now: $ \xcA x \xbe X' -X$
$ \xcA  \xBc x,i \xBe $ $ \xcE \xba:x' \xcp  \xBc x,i \xBe.x' \xbe X$

(3) reads: $ \xcA x \xbe X$ $ \xcE  \xBc x,i \xBe $
$ \xCN \xcE \xba:x' \xcp  \xBc x,i \xBe.x' \xbe X' $

This is intended: intuitively, $X= \xbm (X' ),$ and minimal elements must
not be
attacked at all, but non-minimals must be attacked from $X$ - which is a
modified
version of smoothness.
\index{Remark X-Sub-X'}

\ed

\br

$\hspace{0.01em}$

% (+++ Orig. No.:  Remark X-Sub-X' +++)

\label{Remark X-Sub-X'}

We note the special case of Definition \ref{Definition X-Sub-X'} (page
\pageref{Definition X-Sub-X'})  for level
3 structures,
as it will be used later. We also write it immediately for the intended
case
$ \xbm (X) \xes X,$ and explicitly with copies.

$x \xbe \xbm (X)$ iff

(1) $ \xcE  \xBc x,i \xBe  \xcA  \xBc  \xba,k \xBe : \xBc y,j \xBe  \xcp  \xBc
x,i \xBe $

$ \xDC (y \xbe X$ $ \xch $
$ \xcE  \xBc  \xbb ',l'  \xBe : \xBc z',m'  \xBe  \xcp  \xBc  \xba,k \xBe.$

$ \xDC \xDC (z' \xbe \xbm (X)$ $ \xcu $
$ \xCN \xcE  \xBc  \xbg ',n'  \xBe : \xBc u',p'  \xBe  \xcp  \xBc  \xbb ',l' 
\xBe.u' \xbe X))$

See Diagram \ref{Diagram Essential-Smooth-3-1-2} (page \pageref{Diagram
Essential-Smooth-3-1-2}).

$x \xbe X- \xbm (X)$ iff

(2) $ \xcA  \xBc x,i \xBe  \xcE  \xBc  \xba ',k'  \xBe : \xBc y',j'  \xBe  \xcp 
\xBc x,i \xBe $

$ \xDC (y' \xbe \xbm (X)$ $ \xcu $

$ \xDC \xDC (a)$ $ \xCN \xcE  \xBc  \xbb ',l'  \xBe : \xBc z',m'  \xBe  \xcp
\xBc \xba ',k'
\xBe.z' \xbe X$

$ \xDC \xDC or$

$ \xDC \xDC (b)$ $ \xcA  \xBc  \xbb ',l'  \xBe : \xBc z',m'  \xBe  \xcp \xBc
\xba ',k' \xBe
$

$ \xDC \xDC \xDC (z' \xbe X$ $ \xch $
$ \xcE  \xBc  \xbg ',n'  \xBe : \xBc u',p'  \xBe  \xcp  \xBc  \xbb ',l'  \xBe
.u' \xbe \xbm (X))$ )

See Diagram \ref{Diagram Essential-Smooth-3-2} (page \pageref{Diagram
Essential-Smooth-3-2}).

\vspace{15mm}

\begin{diagram}

\label{Diagram Essential-Smooth-3-1-2}
\index{Diagram Essential-Smooth-3-1-2}

\centering
\setlength{\unitlength}{1mm}
{\renewcommand{\dashlinestretch}{30}
\begin{picture}(100,130)(0,0)
\put(50,80){\circle{80}}
\path(10,80)(90,80)

\put(100,100){{\xssc $X$}}
\put(100,60){{\xssc $\xbm (X)$}}

\path(50,100)(50,60)
\path(48.5,63)(50,60)(51.5,63)
\put(50,101){\circle*{0.3}}
\put(50,59){\circle*{0.3}}
\put(50,57){{\xssc $ \xBc x,i \xBe $}}
\put(50,102){{\xssc $ \xBc y,j \xBe $}}
\put(52,90){{\xssc $ \xBc \xba,k \xBe $}}

\path(20,70)(48,70)
\path(46,71)(48,70)(46,69)
\put(19,70){\circle*{0.3}}
\put(13,67){{\xssc $ \xBc z',m' \xBe $}}
\put(26,71){{\xssc $ \xBc \xbb',l' \xBe $}}

\path(30,30)(30,68)
\path(29,66)(30,68)(31,66)
\put(30,29){\circle*{0.3}}
\put(31,53){{\xssc $ \xBc \xbg ',n' \xBe $}}
\put(30,26.5){{\xssc $ \xBc u',p' \xBe $}}

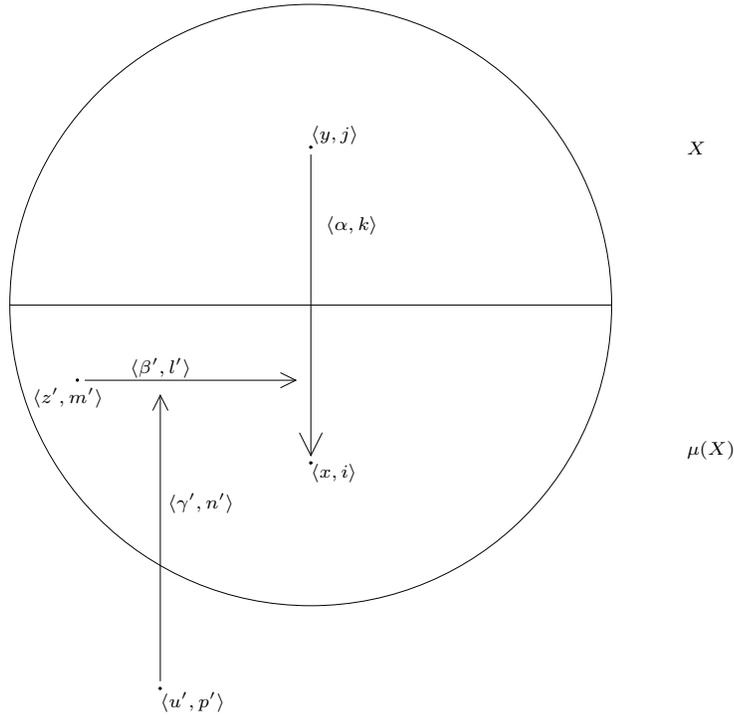
\put(30,10) {{\rm\bf Case 3-1-2}}

\end{picture}
}
\end{diagram}

\vspace{4mm}

\vspace{10mm}

\begin{diagram}

\label{Diagram Essential-Smooth-3-2}
\index{Diagram Essential-Smooth-3-2}

\centering
\setlength{\unitlength}{1mm}
{\renewcommand{\dashlinestretch}{30}
\begin{picture}(100,110)(0,0)
\put(50,60){\circle{80}}
\path(10,60)(90,60)

\put(100,80){{\xssc $X$}}
\put(100,40){{\xssc $\xbm (X)$}}

\path(50,80)(50,40)
\path(48.5,77)(50,80)(51.5,77)
\put(50,81){\circle*{0.3}}
\put(50,39){\circle*{0.3}}
\put(50,37){{\xssc $ \xBc y',j' \xBe $}}
\put(50,82){{\xssc $ \xBc x,i \xBe $}}
\put(52,70){{\xssc $ \xBc \xba',k' \xBe $}}

\path(20,70)(48,70)
\path(46,71)(48,70)(46,69)
\put(19,70){\circle*{0.3}}
\put(13,67){{\xssc $ \xBc z',m' \xBe $}}
\put(26,71){{\xssc $ \xBc \xbb',l' \xBe $}}

\path(30,30)(30,68)
\path(29,66)(30,68)(31,66)
\put(30,29){\circle*{0.3}}
\put(31,53){{\xssc $ \xBc \xbg ',n' \xBe $}}
\put(31,26.5){{\xssc $ \xBc u',p' \xBe $}}

\put(30,0) {{\rm\bf Case 3-2}}

\end{picture}
}
\end{diagram}

\vspace{4mm}

\index{Fact X-Sub-X'}

\er

\bfa

$\hspace{0.01em}$

% (+++ Orig. No.:  Fact X-Sub-X' +++)

\label{Fact X-Sub-X'}

(1) If $X \xes X',$ then $X= \xbm (X' ),$

(2) $X \xes X',$ $X \xcc X'' \xcc X' $ $ \xch $ $X \xes X''.$ (This
corresponds to $( \xbm CUM).)$

(3) $X \xes X',$ $X \xcc Y',$ $Y \xes Y',$ $Y \xcc X' $ $ \xch $ $X=Y.$
(This corresponds to $( \xbm \xcc \xcd ).)$
\index{Fact X-Sub-X' Proof}

\efa

\subparagraph{
Proof
}

$\hspace{0.01em}$

% (+++ Orig.:  Proof +++)

\subparagraph{
Proof
}

$\hspace{0.01em}$

% (+++ Orig.:  Proof +++)

(1) Trivial by Fact \ref{Fact Higher-Validity} (page \pageref{Fact
Higher-Validity})  (1).

(2)

We have to show

(a) $ \xcA x \xbe X'' -X$ $ \xcA  \xBc x,i \xBe $
$ \xcE \xba:x' \xcp  \xBc x,i \xBe ( \xba $ is a valid $X \xch X'' $ arrow), and

(b) $ \xcA x \xbe X$ $ \xcE  \xBc x,i \xBe $
$( \xcA \xba:x' \xcp  \xBc x,i \xBe (x' \xbe X'' $ $ \xch $ $ \xcE \xbb:x'' \xcp
\xba.( \xbb $ is a valid $X \xch X'' $ arrow))).

Both follow from the corresponding condition for $X \xch X',$ the
restriction of the
universal quantifier, and Fact \ref{Fact Higher-Validity} (page \pageref{Fact
Higher-Validity})  (2).

(3)

Let $x \xbe X-$Y.

(a) By $x \xbe X \xes X',$ $ \xcE  \xBc x,i \xBe $ s.t.
$( \xcA \xba:x' \xcp  \xBc x,i \xBe (x' \xbe X' $ $ \xch $ $ \xcE \xbb:x'' \xcp
\xba.( \xbb $ is a valid $X \xch X' $ arrow))).

(b) By $x \xce Y \xes \xcE \xba_{1}:x' \xcp  \xBc x,i \xBe $ $ \xba_{1}$ is a
valid
$Y \xch Y' $ arrow, in
particular $x' \xbe Y \xcc X'.$ Moreover, $ \xbl ( \xba_{1})=1.$

So by (a) $ \xcE \xbb_{2}:x'' \xcp \xba_{1}.( \xbb_{2}$ is a valid $X \xch
X' $ arrow), in particular $x'' \xbe X \xcc Y',$
moreover $ \xbl ( \xbb_{2})=2.$

It follows by induction from the definition of valid $A \xch B$ arrows
that

$ \xcA n \xcE \xba_{2m+1},$ $ \xbl ( \xba_{2m+1})=2m+1,$ $ \xba_{2m+1}$ a
valid $Y \xch Y' $ arrow and

$ \xcA n \xcE \xbb_{2m+2},$ $ \xbl ( \xbb_{2m+2})=2m+2,$ $ \xbb_{2m+2}$ a
valid $X \xch X' $ arrow,

which is impossible, as $ \xdx $ is a structure of finite level.

$ \xcz $
\\[3ex]
\index{Definition Totally-Smooth}

\bd

$\hspace{0.01em}$

% (+++ Orig. No.:  Definition Totally-Smooth +++)

\label{Definition Totally-Smooth}

Let $ \xdx $ be a generalized preferential structure, $X \xcc \xdP ( \xdx
).$

$ \xdx $ is called totally smooth for $X$ iff

(1) $ \xcA \xba:x \xcp y \xbe \xdA ( \xdx )(O( \xba ) \xcv D( \xba ) \xcc
X$ $ \xch $ $ \xcE \xba ':x' \xcp y.x' \xbe \xbm (X))$

(2) if $ \xba $ is valid, then there must also exist such $ \xba ' $ which
is valid.

(y a point or an arrow).

If $ \xdy \xcc \xdP ( \xdx ),$ then $ \xdx $ is called $ \xdy -$totally
smooth iff for all $X \xbe \xdy $
$ \xdx $ is totally smooth for $X.$
\index{Example Totally-Smooth}

\ed

\be

$\hspace{0.01em}$

% (+++ Orig. No.:  Example Totally-Smooth +++)

\label{Example Totally-Smooth}

$X:=\{ \xba:a \xcp b,$ $ \xba ':b \xcp c,$ $ \xba '':a \xcp c,$ $ \xbb
:b \xcp \xba ' \}$ is not totally smooth,

$X:=\{ \xba:a \xcp b,$ $ \xba ':b \xcp c,$ $ \xba '':a \xcp c,$ $ \xbb
:b \xcp \xba ',$ $ \xbb ':a \xcp \xba ' \}$ is totally smooth.
\index{Example Need-Smooth}

\ee

\be

$\hspace{0.01em}$

% (+++ Orig. No.:  Example Need-Smooth +++)

\label{Example Need-Smooth}

Consider $ \xba ':a \xcp b,$ $ \xba '':b \xcp c,$ $ \xba:a \xcp c,$ $
\xbb:a \xcp \xba.$

Then $ \xbm (\{a,b,c\})=\{a\},$ $ \xbm (\{a,c\})=\{a,c\}.$
Thus, $( \xbm CUM)$ does not hold in this structure.
Note that there is no valid arrow from $ \xbm (\{a,b,c\})$ to $c.$
\index{Definition Essentially-Smooth}

\ee

\bd

$\hspace{0.01em}$

% (+++ Orig. No.:  Definition Essentially-Smooth +++)

\label{Definition Essentially-Smooth}

Let $ \xdx $ be a generalized preferential structure, $X \xcc \xdP ( \xdx
).$

$ \xdx $ is called essentially smooth for $X$ iff $ \xbm (X) \xes X.$

If $ \xdy \xcc \xdP ( \xdx ),$ then $ \xdx $ is called $ \xdy
-$essentially smooth iff for all $X \xbe \xdy $
$ \xbm (X) \xes X.$
\index{Example Total-vs-Essential}

\ed

\be

$\hspace{0.01em}$

% (+++ Orig. No.:  Example Total-vs-Essential +++)

\label{Example Total-vs-Essential}

It is easy to see that we can distinguish total and essential smoothness
in richer structures, as the following Example shows:

We add an accessibility relation $R,$ and consider only those models which
are accessible.

Let e.g. $a \xcp b \xcp  \xBc c,0 \xBe,$ $ \xBc c,1 \xBe,$ without
transitivity. Thus, only
$c$ has two
copies. This structure is essentially smooth, but of course not totally
so.

Let now mRa, mRb, $mR \xBc c,0 \xBe,$ $mR \xBc c,1 \xBe,$ $m' Ra,$ $m' Rb,$
$m' R \xBc c,0 \xBe.$

Thus, seen from $m,$ $ \xbm (\{a,b,c\})=\{a,c\},$ but seen from $m',$ $
\xbm (\{a,b,c\})=\{a\},$
but $ \xbm (\{a,c\})=\{a,c\},$ contradicting (CUM).

$ \xcz $
\\[3ex]
\section{
The general case
}
\label{Section Reac-GenPref-NonSmooth}
\index{Section Reac-GenPref-NonSmooth}

\ee

The idea to solve the representation problem illustrated by Example \ref{Example
Need-Pr} (page \pageref{Example Need-Pr})
is to use the points $c$ and $d$ as bases for counterarguments against $
\xba:b \xcp a$ - as
is possible in IBRS. We do this now. We will obtain a representation
for logics weaker than $ \xCf P$ by generalized preferential structures.

We will now prove a representation theorem, but will make it more general
than
for preferential structures only. For this purpose, we will introduce some
definitions first.

\bd

$\hspace{0.01em}$

% (+++ Orig. No.:  Definition Eta-Rho-Structure +++)

\label{Definition Eta-Rho-Structure}

Let $ \xbh, \xbr: \xdy \xcp \xdp (U).$

(1) If $ \xdx $ is a simple structure:

$ \xdx $ is called an attacking structure relative to $ \xbh $
representing $ \xbr $ iff

$ \xbr (X)$ $=$ $\{x \xbe \xbh (X):$ there is no valid $X-to- \xbh (X)$
arrow $ \xba:x' \xcp x\}$

for all $X \xbe \xdy.$

(2) If $ \xdx $ is a structure with copies:

$ \xdx $ is called an attacking structure relative to $ \xbh $
representing $ \xbr $ iff

$ \xbr (X)$ $=$ $\{x \xbe \xbh (X):$ there is $ \xBc x,i \xBe $ and no valid
$X-to-
\xbh (X)$ arrow
$ \xba: \xBc x',i'  \xBe  \xcp  \xBc x,i \xBe \}$

for all $X \xbe \xdy.$

Obviously, in those cases $ \xbr (X) \xcc \xbh (X)$ for all $X \xbe \xdy
.$

Thus, $ \xdx $ is a preferential structure iff $ \xbh $ is the identity.

See Diagram \ref{Diagram Eta-Rho-1} (page \pageref{Diagram Eta-Rho-1})

\vspace{10mm}

\begin{diagram}

\label{Diagram Eta-Rho-1}
\index{Diagram Eta-Rho-1}

\centering
\setlength{\unitlength}{1mm}
{\renewcommand{\dashlinestretch}{30}
\begin{picture}(150,100)(0,0)

\put(30,50){\circle{50}}
\put(100,50){\circle{40}}
\put(100,50){\circle{70}}

\put(30,50){\circle*{1}}
\put(70,50){\circle*{1}}

\path(31,50)(69,50)
\path(66.6,50.9)(69,50)(66.6,49.1)

\put(30,60){\xssc{$X$}}
\put(100,80){\xssc{$\xbh(X)$}}
\put(100,60){\xssc{$\xbr(X)$}}

\put(30,10) {{\rm\bf Attacking structure}}

\end{picture}
}

\end{diagram}

\vspace{4mm}

\ed

(Note that it does not seem very useful to generalize the notion of
smoothness
from preferential structures to general attacking structures, as, in the
general
case, the minimizing set $X$ and the result $ \xbr (X)$ may be disjoint.)

The following result is the first positive representation result of this
paper, and shows that we can obtain (almost) anything with level 2
structures.

\bp

$\hspace{0.01em}$

% (+++ Orig. No.:  Proposition Eta-Rho-Repres +++)

\label{Proposition Eta-Rho-Repres}

Let $ \xbh, \xbr: \xdy \xcp \xdp (U).$ Then there is an attacking level
2 structure relative
to $ \xbh $ representing $ \xbr $ iff

(1) $ \xbr (X) \xcc \xbh (X)$ for all $X \xbe \xdy,$

(2) $ \xbr ( \xCQ )= \xbh ( \xCQ )$ if $ \xCQ \xbe \xdy.$

\ep

(2) is, of course, void for preferential structures.

\subparagraph{
Proof
}

$\hspace{0.01em}$

% (+++ Orig.:  Proof +++)

(A) The construction

We make a two stage construction.

(A.1) Stage 1.

In stage one, consider (almost as usual)

$ \xdu:= \xBc  \xdx,\{ \xba_{i}:i \xbe I\} \xBe $ where

$ \xdx $ $:=$ $\{ \xBc x,f \xBe :$ $x \xbe U,$ $f \xbe \xbP \{X \xbe \xdy:x \xbe
\xbh (X)- \xbr (X)\}\},$

$ \xba:x' \xcp  \xBc x,f \xBe $ $: \xcj $ $x' \xbe ran(f).$ Attention: $x' \xbe
X,$
not $x' \xbe \xbr (X)!$

(A.2) Stage 2.

Let $ \xdx ' $ be the set of all $ \xBc x,f,X \xBe $ s.t. $ \xBc x,f \xBe  \xbe
\xdx $ and

(a) either $X$ is some dummy value, say $*$

or

(b) all of the following $(1)-(4)$ hold:

(1) $X \xbe \xdy,$

(2) $x \xbe \xbr (X),$

(3) there is $X' \xcc X,$ $x \xbe \xbh (X' )- \xbr (X' ),$ $X' \xbe \xdy
,$
(thus $ran(f) \xcs X \xEd \xCQ $ by definition),

(4) $ \xcA X'' \xbe \xdy.(X \xcc X'',$ $x \xbe \xbh (X'' )- \xbr (X'' )$
$ \xch $ $(ran(f) \xcs X'' )-X \xEd \xCQ ).$

(Thus, $f$ chooses in (4) for $X'' $ also outside $X.$ If there is no such
$X'',$ (4) is
void, and only $(1)-(3)$ need to hold, i.e. we may take any $f$ with
$ \xBc x,f \xBe  \xbe \xdx.)$

See Diagram \ref{Diagram Condition-rho-eta} (page \pageref{Diagram
Condition-rho-eta}).

\vspace{10mm}

\begin{diagram}

\label{Diagram Condition-rho-eta}
\index{Diagram Condition-rho-eta}

\centering
\setlength{\unitlength}{0.00083333in}
{\renewcommand{\dashlinestretch}{30}
\begin{picture}(2727,2755)(0,-500)
\put(1304.562,-118.818){\arc{3584.794}{4.1827}{5.2252}}
\put(1356.377,-664.180){\arc{3978.201}{4.0248}{5.3839}}
\put(1370.659,430.666){\arc{1327.952}{3.9224}{5.4136}}
\put(1344,1343){\ellipse{2672}{2672}}
\put(1334,1293){\ellipse{1178}{1178}}
\put(1289,1338){\ellipse{1802}{1802}}
\put(2200,1950){\circle*{30}}
\put(2094,1800){{\xssc $f(X'')$}}
\put(2214,1318){{\xssc $\rho(X)$}}
\put(2074,2608){{\xssc $X''$}}
\put(1709,2228){{\xssc $X$}}
\put(2624,698) {{\xssc $\rho(X'')$}}
\put(1334,828) {{\xssc $\rho(X')$}}
\put(1404,1928){{\xssc $X'$}}
\put(1514,1433){\circle*{30}}
\put(1544,1403){{\xssc $x$}}
\put(100,2800){{\xssc For simplicity, $\xbh(X)=X$ here}}

\put(150,-400) {{\rm\bf The complicated case}}

\end{picture}
}
\end{diagram}

\vspace{4mm}

Note: If $(1)-(3)$ are satisfied for $x$ and $X,$ then we will find $f$
s.t. $ \xBc x,f \xBe  \xbe \xdx,$
and $ \xBc x,f,X \xBe $ satisfies $(1)-(4):$ As $X \xcB X'' $ for $X'' $ as in
(4),
we find $f$ which
chooses for such $X'' $ outside of $X.$

So for any $ \xBc x,f \xBe  \xbe \xdx,$ there is $ \xBc x,f,* \xBe,$ and maybe
also some
$ \xBc x,f,X \xBe $ in $ \xdx '.$

Let again for any $x',$ $ \xBc x,f,X \xBe  \xbe \xdx ' $

$ \xba:x' \xcp  \xBc x,f,X \xBe $ $: \xcj $ $x' \xbe ran(f)$

(A.3) Adding arrows.

Consider $x' $ and $ \xBc x,f,X \xBe.$

If $X=*,$ or $x' \xce X,$ we do nothing, i.e. leave a simple arrow
$ \xba:x' \xcp  \xBc x,f,X \xBe $ $ \xcj $ $x' \xbe ran(f).$

If $X \xbe \xdy,$ and $x' \xbe X,$ and $x' \xbe ran(f),$ we make $X$ many
copies of the attacking
arrow and have then: $ \xBc  \xba,x''  \xBe :x' \xcp  \xBc x,f,X \xBe $ for all
$x'' \xbe X.$

In addition, we add attacks on the $ \xBc  \xba,x''  \xBe :$ $ \xBc  \xbb,x'' 
\xBe :x''
\xcp  \xBc  \xba,x''  \xBe $ for all $x'' \xbe X.$

The full structure $ \xdz $ is thus:

$ \xdx ' $ is the set of elements.

If $x' \xbe ran(f),$ and $X=*$ or $x' \xce X$ then $ \xba:x' \xcp
\xBc x,f,X \xBe $

if $x' \xbe ran(f),$ and $X \xEd *$ and $x' \xbe X$ then

(a) $ \xBc  \xba,x''  \xBe :x' \xcp  \xBc x,f,X \xBe $ for all $x'' \xbe X,$

(b) $ \xBc  \xbb,x''  \xBe :x'' \xcp  \xBc  \xba,x''  \xBe $ for all $x'' \xbe
X.$

See Diagram \ref{Diagram Structure-rho-eta} (page \pageref{Diagram
Structure-rho-eta}).

\vspace{10mm}

\begin{diagram}

\label{Diagram Structure-rho-eta}
\index{Diagram Structure-rho-eta}

\centering
\setlength{\unitlength}{1mm}
{\renewcommand{\dashlinestretch}{30}
\begin{picture}(150,100)(0,0)

\put(30,50){\circle{50}}
\put(100,50){\circle{40}}
\put(100,50){\circle{70}}

\put(30,50){\circle*{1}}
\put(90,50){\circle*{1}}
\put(92,50){{\xssc $ \xBc  x, f, X \xBe $}}
\put(32,50) {{\xssc $x'$}}

\put(60,20){\arc{84.5}{-2.34}{-0.8}}
\put(70,62){{\xssc $ \xBc  a, x_0'' \xBe $}}
\put(60,80){\arc{84.5}{0.8}{2.34}}
\put(70,38){{\xssc $ \xBc  a, x_1'' \xBe $}}

\path(88.9,51.7)(90,50)(88,50.5)
\path(88,49.5)(90,50)(88.9,48.3)

\put(35,65){\circle*{1}}
\put(35,35){\circle*{1}}

\put(36,65) {{\xssc $x_0''$}}
\put(36,35){{\xssc $x''_1$}}

\path(35,65)(40,58)
\path(39.7,59.2)(40,58)(39.0,58.7)
\put(25,62){{\xssc $ \xBc  \beta, x_1'' \xBe $}}

\path(35,35)(40,42)
\path(39,41.3)(40,42)(39.7,40.8)
\put(25,38) {{\xssc $ \xBc  \beta, x_o'' \xBe $}}

\put(15,50){\xssc{$X$}}
\put(100,80){\xssc{$\xbh(X)$}}
\put(100,60){\xssc{$\xbr(X)$}}

\put(30,10) {{\rm\bf Attacking structure}}

\end{picture}
}

\end{diagram}

\vspace{4mm}

(B) Representation

We have to show that this structure represents $ \xbr $ relative to $ \xbh
.$

Let $y \xbe \xbh (Y),$ $Y \xbe \xdy.$

Case 1. $y \xbe \xbr (Y).$

We have to show that there is $ \xBc y,g,Y''  \xBe $ s.t. there is no valid
$ \xba:y' \xcp  \xBc y,g,Y''  \xBe,$ $y' \xbe Y.$ In Case 1.1 below, $Y'' $
will
be $*,$ in Case 1.2,
$Y'' $ will be $Y,$ $g$ will be chosen suitably.

Case 1.1. There is no $Y' \xcc Y,$ $y \xbe \xbh (Y' )- \xbr (Y' ),$ $Y'
\xbe \xdy.$

So for all $Y' $ with $y \xbe \xbh (Y' )- \xbr (Y' )$ $Y' -Y \xEd \xCQ.$
Let $g \xbe \xbP \{Y' -Y:y \xbe \xbh (Y' )- \xbr (Y' )\}.$
Then $ran(g) \xcs Y= \xCQ,$ and $ \xBc y,g,* \xBe $ is not attacked from $Y.$
$( \xBc y,g \xBe $ was already not attacked in $ \xdx.)$

Case 1.2. There is $Y' \xcc Y,$ $y \xbe \xbh (Y' )- \xbr (Y' ),$ $Y' \xbe
\xdy.$

Let now $ \xBc y,g,Y \xBe  \xbe \xdx ',$ s.t. $g(Y'' ) \xce Y$ if $Y \xcc Y'',$
$y
\xbe \xbh (Y'' )- \xbr (Y'' ),$ $Y'' \xbe \xdy.$
As noted above, such $g$ and thus $ \xBc y,g,Y \xBe $ exist. Fix $ \xBc y,g,Y
\xBe.$

Consider any $y' \xbe ran(g).$ If $y' \xce Y,$ $y' $ does not attack
$ \xBc y,g,Y \xBe $ in $Y.$ Suppose $y' \xbe Y.$ We had made $Y$ many copies
$ \xBc  \xba,y''  \xBe,$ $y'' \xbe Y$ with $ \xBc  \xba,y''  \xBe :y' \xcp 
\xBc y,g,Y \xBe $ and had
added the level 2
arrows $ \xBc  \xbb,y''  \xBe :y'' \xcp  \xBc  \xba,y''  \xBe $ for $y'' \xbe
Y.$
So all copies $ \xBc  \xba,y''  \xBe $ are destroyed in $Y.$ This was done for
all
$y' \xbe Y,$ $y' \xbe ran(g),$ so $ \xBc y,g,Y \xBe $ is now not (validly)
attacked
in $Y$
any more.

Case 2. $y \xbe \xbh (Y)- \xbr (Y).$

Let $ \xBc y,g,Y'  \xBe $ (where $Y' $ can be $*)$ be any copy of $y,$ we have
to
show that
there is $z \xbe Y,$ $ \xba:z \xcp  \xBc y,g,Y'  \xBe,$ or some
$ \xBc  \xba,z'  \xBe :z \xcp  \xBc y,g,Y'  \xBe,$ $z' \xbe Y',$ which is not
destroyed by
some
level 2 arrow $ \xBc  \xbb,z'  \xBe :z' \xcp  \xBc  \xba,z'  \xBe,$ $z' \xbe
Y.$

As $y \xbe \xbh (Y)- \xbr (Y),$ $ran(g) \xcs Y \xEd \xCQ,$ so there is $z
\xbe ran(g) \xcs Y.$ Fix such $z.$
(We will modify the choice of $z$ only in Case 2.2.2 below.)

Case 2.1. $Y' =*.$

As $z \xbe ran(g),$ $ \xba:z \xcp  \xBc y,g,* \xBe.$ (There were no level 2
arrows
introduced for
this copy.)

Case 2.2. $Y' \xEd *.$

So $ \xBc y,g,Y'  \xBe $ satisfies the conditions $(1)-(4)$ of (b) at the
beginning
of the
proof.

If $z \xce Y',$ we are done, as $ \xba:z \xcp  \xBc y,g,Y'  \xBe,$ and there
were
no level 2
arrows introduced in this case.
If $z \xbe Y',$ we had made $Y' $ many copies $ \xBc  \xba,z'  \xBe,$
$ \xBc  \xba,z'  \xBe :z \xcp  \xBc y,g,Y'  \xBe,$ one for
each $z' \xbe Y'.$ Each $ \xBc  \xba,z'  \xBe $ was destroyed by
$ \xBc  \xbb,z'  \xBe :z' \xcp  \xBc  \xba,z'  \xBe,$ $z' \xbe Y'.$

Case 2.2.1. $Y' \xcC Y.$

Let $z'' \xbe Y' -$Y, then
$ \xBc  \xba,z''  \xBe :z \xcp  \xBc y,g,Y'  \xBe $ is destroyed only by
$ \xBc  \xbb,z''  \xBe :z'' \xcp  \xBc  \xba,z''  \xBe $ in $Y',$ but not in
$Y,$ as $z''
\xce Y,$ so
$ \xBc y,g,Y'  \xBe $ is
attacked by $ \xBc  \xba,z''  \xBe :z \xcp  \xBc y,g,Y'  \xBe,$ valid in $Y.$

Case 2.2.2. $Y' \xcB Y$ $(Y=Y' $ is impossible, as $y \xbe \xbr (Y' ),$ $y
\xce \xbr (Y)).$

Then there was by definition (condition (b) (4)) some $z' \xbe (ran(g)
\xcs Y)-Y' $ and
$ \xba:z' \xcp  \xBc y,g,Y'  \xBe $ is valid, as $z' \xce Y'.$ (In this case,
there
are no copies of $ \xba $ and no level 2 arrows.)

$ \xcz $
\\[3ex]

\bco

$\hspace{0.01em}$

% (+++ Orig. No.:  Corollary Eta-Rho +++)

\label{Corollary Eta-Rho}

(1) We cannot distinguish general structures of level 2 from those of
higher
levels by their $ \xbr -$functions relative to $ \xbh.$

(2) Let $U$ be the universe, $ \xdy \xcc \xdp (U),$ $ \xbm: \xdy \xcp
\xdp (U).$
Then any $ \xbm $ satisfying $( \xbm \xcc )$ can be represented by a level
2 preferential
structure. (Choose $ \xbh =identity.)$

Again, we cannot distinguish general structures of level 2 from those of
higher
levels by their $ \xbm -$functions.

$ \xcz $
\\[3ex]

\eco

A remark on the function $ \xbh:$

We can also obtain the function $ \xbh $ via arrows. Of course, then we
need
positive arrows (not only negative arrows against negative arrows, as we
first
need to have something positive).

If $ \xbh $ is the identity, we can make a positive arrow from each point
to itself.
Otherwise, we can connect every point to every point by a positive arrow,
and then choose those we really want in $ \xbh $ by a choice function
obtained from
arrows just as we obtained $ \xbr $ from arrows.
\section{
Discussion of the totally smooth case
}
\label{Section Reac-GenPref-TotalSmooth}
\index{Section Reac-GenPref-TotalSmooth}

\bfa

$\hspace{0.01em}$

% (+++ Orig. No.:  Fact Val-Arrow +++)

\label{Fact Val-Arrow}

Let $X,Y \xbe \xdy,$ $ \xdx $ a level $n$ structure.
Let $ \xBc  \xba,k \xBe : \xBc x,i \xBe  \xcp  \xBc y,j \xBe,$ where
$ \xBc y,j \xBe $ may itself be (a copy of) an arrow.

(1) Let $n>1,$ $X \xcc Y,$ $ \xBc  \xba,k \xBe  \xbe X$ a level $n-1$ arrow in $
\xdx \xex X.$
If $ \xBc  \xba,k \xBe $ is valid in
$ \xdx \xex Y,$ then it is valid in $ \xdx \xex X.$

(2) Let $ \xdx $ be totally smooth, $ \xbm (X) \xcc Y,$ $ \xbm (Y) \xcc
X,$
$ \xBc  \xba,k \xBe  \xbe \xdx \xex X \xcs Y,$ then $ \xBc  \xba,k \xBe $
is valid in $ \xdx \xex X$ iff it is valid in $ \xdx \xex Y.$

Note that we will also sometimes write $X$ for $ \xdx \xex X,$ when the
context is clear.

\efa

\subparagraph{
Proof
}

$\hspace{0.01em}$

% (+++ Orig.:  Proof +++)

(1)
If $ \xBc  \xba,k \xBe $ is not valid in $ \xdx \xex X,$ then there must be a
level
$n$ arrow
$ \xBc  \xbb,r \xBe : \xBc z,s \xBe  \xcp  \xBc  \xba,k \xBe $ in $ \xdx \xex X
\xcc \xdx \xex Y.$
$ \xBc  \xbb,r \xBe $ must be valid in $ \xdx \xex X$ and $ \xdx \xex Y,$
as there are no level $n+1$ arrows.
So $ \xBc  \xba,k \xBe $ is not valid in $ \xdx \xex Y,$ $contradiction.$

(2)
By downward induction.
Case $n:$ $ \xBc  \xba,k \xBe  \xbe \xdx \xex X \xcs Y,$ so it is valid in both
as
there are no level $n+1$
arrows.
Case $m \xcp m-1:$ Let $ \xBc  \xba,k \xBe  \xbe \xdx \xex X \xcs Y$ be a level
$m-1$ arrow valid in $ \xdx \xex X,$ but not
in $ \xdx \xex Y.$ So there must be a level $m$ arrow
$ \xBc  \xbb,r \xBe : \xBc z,s \xBe  \xcp  \xBc  \xba,k \xBe $ valid in $ \xdx
\xex Y.$
By total smoothness, we may assume $z \xbe \xbm (Y) \xcc X,$ so
$ \xBc  \xbb,r \xBe  \xbe \xdx \xex X$ is valid by
induction hypothesis. So
$ \xBc  \xba,k \xBe $ is not valid in $ \xdx \xex X,$ $contradiction.$

$ \xcz $
\\[3ex]

\bco

$\hspace{0.01em}$

% (+++ Orig. No.:  Corollary Total-Mu +++)

\label{Corollary Total-Mu}

Let $X,Y \xbe \xdy,$ $ \xdx $ a totally smooth level $n$ structure, $
\xbm (X) \xcc Y,$ $ \xbm (Y) \xcc X.$
Then $ \xbm (X)= \xbm (Y).$

\eco

\subparagraph{
Proof
}

$\hspace{0.01em}$

% (+++ Orig.:  Proof +++)

Let $x \xbe \xbm (X)- \xbm (Y).$ Then by $ \xbm (X) \xcc Y,$ $x \xbe Y,$
so there must be for all
$ \xBc x,i \xBe  \xbe \xdx $
an arrow $ \xBc  \xba,k \xBe : \xBc y,j \xBe  \xcp  \xBc x,i \xBe $ valid in $
\xdx \xex Y,$ wlog. $y
\xbe \xbm (Y) \xcc X$ by total
smoothness. So by Fact \ref{Fact Val-Arrow} (page \pageref{Fact Val-Arrow}),
(2),
$ \xBc  \xba,k \xBe $ is valid in $ \xdx \xex X.$ This
holds for all $ \xBc x,i \xBe,$ so $x \xce \xbm (X),$ $contradiction.$
$ \xcz $
\\[3ex]

\bfa

$\hspace{0.01em}$

% (+++ Orig. No.:  Fact Level-bigger-2 +++)

\label{Fact Level-bigger-2}

There are situations satisfying $( \xbm \xcc )+( \xbm CUM)+( \xcs )$ which
cannot be represented
by level 2 totally smooth preferential structures.

\efa

The proof is given in the following example.

\be

$\hspace{0.01em}$

% (+++ Orig. No.:  Example Level-bigger-2 +++)

\label{Example Level-bigger-2}

Let $Y:=\{x,y,y' \},$ $X:=\{x,y\},$ $X':=\{x,y' \}.$ Let $ \xdy:= \xdp
(Y).$
Let $ \xbm (Y):=\{y,y' \},$ $ \xbm (X):= \xbm (X' ):=\{x\},$ and $ \xbm
(Z):=Z$ for all other sets.

Obviously, this satisfies $( \xcs ),$ $( \xbm \xcc ),$ and $( \xbm CUM).$

Suppose $ \xdx $ is a totally smooth level 2 structure representing $ \xbm
.$

So $ \xbm (X)= \xbm (X' ) \xcc Y- \xbm (Y),$ $ \xbm (Y) \xcc X \xcv X'.$
Let $ \xBc x,i \xBe $ be minimal in $ \xdx \xex X.$
As $ \xBc x,i \xBe $ cannot be minimal in $ \xdx \xex Y,$ there must be
$ \xba: \xBc z,j \xBe  \xcp  \xBc x,i \xBe,$ valid in
$ \xdx \xex Y.$

Case 1: $z \xbe X'.$

So $ \xba \xbe \xdx \xex X'.$
If $ \xba $ is valid in $ \xdx \xex X',$ there must be
$ \xba ': \xBc x',i'  \xBe  \xcp  \xBc x,i \xBe,$ $x' \xbe \xbm (X' ),$ valid
in
$ \xdx \xex X',$ and thus in $ \xdx \xex X,$ by $ \xbm (X)= \xbm (X' )$
and
Fact \ref{Fact Val-Arrow} (page \pageref{Fact Val-Arrow})  (2). This is
impossible, so there must be
$ \xbb: \xBc x',i'  \xBe  \xcp \xba,$ $x' \xbe \xbm (X' ),$ valid in $ \xdx \xex
X'.$
As $ \xbb $ is in $ \xdx \xex Y$ and $ \xdx $ a level $ \xck 2$ structure,
$ \xbb $ is valid in $ \xdx \xex Y,$ so
$ \xba $ is not valid in $ \xdx \xex Y,$ $contradiction.$

Case 2: $z \xbe X.$

$ \xba $ cannot be valid in $ \xdx \xex X,$ so there must be
$ \xbb: \xBc x',i'  \xBe  \xcp \xba,$ $x' \xbe \xbm (X),$
valid in $ \xdx \xex X.$ Again, as $ \xbb $ is in $ \xdx \xex Y$ and $
\xdx $ a level $ \xck 2$ structure, $ \xbb $ is
valid in $ \xdx \xex Y,$ so $ \xba $ is not valid in $ \xdx \xex Y,$
$contradiction.$

$ \xcz $
\\[3ex]

\ee

It is unknown to the authors whether an analogon is true for essential
smoothness, i.e. whether there are examples of such $ \xbm $ function
which need
at least level 3 essentially smooth structures for representation.
Proposition \ref{Proposition Level-3-Repr} (page \pageref{Proposition
Level-3-Repr})  below shows
that such structures suffice, but we do not
know whether level 3 is necessary.

\bfa

$\hspace{0.01em}$

% (+++ Orig. No.:  Fact Level-3-Solution +++)

\label{Fact Level-3-Solution}

Above Example \ref{Example Level-bigger-2} (page \pageref{Example
Level-bigger-2})
can be solved by a totally smooth level 3 structure:

Let $ \xba_{1}:x \xcp y,$ $ \xba_{2}:x \xcp y',$ $ \xba_{3}:y \xcp x,$
$ \xbb_{1}:y \xcp \xba_{2},$ $ \xbb_{2}:y' \xcp \xba_{1},$ $ \xbb_{3}:y
\xcp \xba_{3},$ $ \xbb_{4}:x \xcp \xba_{3},$
$ \xbg_{1}:y' \xcp \xbb_{3},$ $ \xbg_{2}:y' \xcp \xbb_{4}.$

See Diagram \ref{Diagram Smooth-Level-3} (page \pageref{Diagram Smooth-Level-3})
.

\vspace{10mm}

\begin{diagram}

\label{Diagram Smooth-Level-3}
\index{Diagram Smooth-Level-3}

\centering
\setlength{\unitlength}{0.00083333in}
{\renewcommand{\dashlinestretch}{30}
\begin{picture}(3584,3681)(0,-500)
\put(1980.858,1684.979){\arc{2896.137}{3.0863}{5.1276}}
\path(503.037,1724.492)(535.000,1605.000)(563.029,1725.476)
\put(613.871,2084.678){\arc{1045.415}{1.8998}{5.0662}}
\path(671.312,2573.903)(795.000,2575.000)(685.347,2632.239)
\put(1755.902,1748.510){\arc{3480.502}{1.0969}{3.3803}}
\path(69.787,2036.400)(65.000,2160.000)(11.059,2048.688)
\put(1915.859,2762.471){\arc{1397.908}{3.4506}{5.9440}}
\path(1268.561,3097.293)(1250.000,2975.000)(1323.927,3074.171)
\put(1892.516,1773.709){\arc{3365.589}{4.8901}{7.4299}}
\path(2313.635,3433.788)(2190.000,3430.000)(2300.872,3375.161)
\path(605,1595)(2535,2930)
\path(2453.376,2837.062)(2535.000,2930.000)(2419.243,2886.407)
\path(600,1540)(2545,245)
\put(550,1560){\circle*{30}}
\path(2428.488,286.533)(2545.000,245.000)(2461.741,336.476)
\path(2560,2870)(1820,765)
\path(1831.496,888.158)(1820.000,765.000)(1888.100,868.259)
\path(2530,270)(1745,2345)
\put(2570,220){\circle*{30}}
\put(2600,2950){\circle*{30}}
\path(1815.520,2243.378)(1745.000,2345.000)(1759.401,2222.148)
\put(3315,2745){{\xssc $\gamma_1$}}
\put(2570,2845){{\xssc $y$}}
\put(2575,60)  {{\xssc $y'$}}
\put(1960,950) {{\xssc $\beta_1$}}
\put(1915,1985){{\xssc $\beta_2$}}
\put(2030,2445){{\xssc $\alpha_1$}}
\put(520,1435) {{\xssc $x$}}
\put(1960,3510){{\xssc $\beta_3$}}
\put(750,185)  {{\xssc $\gamma_2$}}
\put(170,2530) {{\xssc $\beta_4$}}
\put(630,2110) {{\xssc $\alpha_3$}}
\put(1920,500) {{\xssc $\alpha_2$}}

\put(100,-400) {{\rm\bf Solution by smooth level 3 structure}}

\end{picture}
}
\end{diagram}

\vspace{4mm}

\efa

The subdiagram generated by $X$ contains $ \xba_{1},$ $ \xba_{3},$ $
\xbb_{3},$ $ \xbb_{4}.$
$ \xba_{1},$ $ \xbb_{3},$ $ \xbb_{4}$ are valid, so $ \xbm (X)=\{x\}.$

The subdiagram generated by $X' $ contains $ \xba_{2}.$
$ \xba_{2}$ is valid, so $ \xbm (X' )=\{x\}.$

In the full diagram, $ \xba_{3},$ $ \xbb_{1},$ $ \xbb_{2},$ $ \xbg_{1},$ $
\xbg_{2}$ are valid, so $ \xbm (Y)=\{y,y' \}.$

$ \xcz $
\\[3ex]

\br

$\hspace{0.01em}$

% (+++ Orig. No.:  Remark Need-Mucd +++)

\label{Remark Need-Mucd}

Example \ref{Example Mu-Cum-Cd} (page \pageref{Example Mu-Cum-Cd})  together
with Corollary \ref{Corollary Total-Mu} (page \pageref{Corollary Total-Mu}) 
show
that $( \xbm \xcc )$ and $( \xbm CUM)$ without $( \xcs )$ do not guarantee
representability by a
level $n$ totally smooth structure.
\section{
The essentially smooth case
}
\label{Section Reac-GenPref-EssSmooth}
\index{Section Reac-GenPref-EssSmooth}

\er

\bd

$\hspace{0.01em}$

% (+++ Orig. No.:  Definition ODPi +++)

\label{Definition ODPi}

Let $ \xbm: \xdy \xcp \xdp (U)$ and $ \xdx $ be given, let
$ \xba: \xBc y,j \xBe  \xcp  \xBc x,i \xBe  \xbe \xdx.$

Define

$ \xdO ( \xba )$ $:=$ $\{Y \xbe \xdy:x \xbe Y- \xbm (Y),y \xbe \xbm
(Y)\},$

$ \xdD ( \xba )$ $:=$ $\{X \xbe \xdy:x \xbe \xbm (X),y \xbe X\},$

$ \xbP ( \xdO, \xba )$ $:=$ $ \xbP \{ \xbm (Y):Y \xbe \xdO ( \xba )\},$

$ \xbP ( \xdD, \xba )$ $:=$ $ \xbP \{ \xbm (X):X \xbe \xdD ( \xba )\}.$

\ed

\bl

$\hspace{0.01em}$

% (+++ Orig. No.:  Lemma Level-3-Constr +++)

\label{Lemma Level-3-Constr}

Let $U$ be the universe, $ \xbm: \xdy \xcp \xdp (U).$
Let $ \xbm $ satisfy $( \xbm \xcc )+( \xbm \xcc \xcd ).$

Let $ \xdx $ be a level 1 preferential structure,
$ \xba: \xBc y,j \xBe  \xcp  \xBc x,i \xBe,$ $ \xdO ( \xba ) \xEd \xCQ,$
$ \xdD ( \xba ) \xEd \xCQ.$

We can modify $ \xdx $ to a level 3 structure $ \xdx ' $ by introducing
level 2 and level 3
arrows s.t. no copy of $ \xba $ is valid in any $X \xbe \xdD ( \xba ),$
and in every $Y \xbe \xdO ( \xba )$ at
least one copy of $ \xba $ is valid. (More precisely, we should write $
\xdx ' \xex X$ etc.)

Thus, in $ \xdx ',$

(1) $ \xBc x,i \xBe $ will not be minimal in any $Y \xbe \xdO ( \xba ),$

(2) if $ \xba $ is the only arrow minimizing
$ \xBc x,i \xBe $ in $X \xbe \xdD ( \xba ),$ $ \xBc x,i \xBe $ will now be
minimal in $X.$

The construction is made independently for all such arrows $ \xba \xbe
\xdx.$

(This is probably the main technical result of the paper.)

\el

\subparagraph{
Proof
}

$\hspace{0.01em}$

% (+++ Orig.:  Proof +++)

(1) The construction

Make $ \xbP ( \xdD, \xba )$ many copies of $ \xba:$
$\{ \xBc  \xba,f \xBe :f \xbe \xbP ( \xdD, \xba )\},$ all
$ \xBc  \xba,f \xBe : \xBc y,j \xBe  \xcp  \xBc x,i \xBe.$ Note that
$ \xBc  \xba,f \xBe  \xbe X$ for all $X \xbe \xdD ( \xba )$ and
$ \xBc  \xba,f \xBe  \xbe Y$ for all
$Y \xbe \xdO ( \xba ).$

Add to the structure
$ \xBc  \xbb,f,X_{r},g \xBe : \xBc f(X_{r}),i_{r} \xBe  \xcp  \xBc  \xba,f \xBe
,$ for any $X_{r} \xbe
\xdD ( \xba ),$ and
$g \xbe \xbP ( \xdO, \xba )$ (and some or all $i_{r}$ - this does not
matter).

For all $Y_{s} \xbe \xdO ( \xba ):$

if $ \xbm (Y_{s}) \xcC X_{r}$ and $f(X_{r}) \xbe Y_{s},$ then add to the
structure
$ \xBc  \xbg,f,X_{r},g,Y_{s} \xBe : \xBc g(Y_{s}),j_{s} \xBe  \xcp \xBc
\xbb,f,X_{r},g \xBe
$ (again for all or some $j_{s}),$

if $ \xbm (Y_{s}) \xcc X_{r}$ or $f(X_{r}) \xce Y_{s},$
$ \xBc  \xbg,f,X_{r},g,Y_{s} \xBe $ is not added.

See Diagram \ref{Diagram Essential-Smooth-Repr} (page \pageref{Diagram
Essential-Smooth-Repr}).

\vspace{30mm}

\begin{diagram}

\label{Diagram Essential-Smooth-Repr}
\index{Diagram Essential-Smooth-Repr}

\centering
\setlength{\unitlength}{1mm}
{\renewcommand{\dashlinestretch}{30}
\begin{picture}(110,150)(0,0)
\put(50,80){\ellipse{120}{80}}
\put(50,100){\ellipse{120}{80}}
\path(20,134.3)(20,65.5)
\path(80,114.3)(80,45.5)
\put(50,142){\xssc{$X \xbe \xdD(\xba)$}}
\put(50,36){\xssc{$Y \xbe \xdO(\xba)$}}
\put(10,138){\xssc{$\xbm(X)$}}
\put(90,44){\xssc{$\xbm(Y)$}}

\path(10,82)(90,82)
\put(9.2,82){\circle*{0.3}}
\put(90.8,82){\circle*{0.3}}
% \path(12,83)(10,82)(12,81)
\path(12.8,81)(10,82)(12.8,83)
\put(5,78){{\xssc $ \xBc x,i \xBe $}}
\put(85,78){{\xssc $ \xBc y,j \xBe $}}
\put(55,78){{\xssc $ \xBc \xba,f \xBe $}}

\path(8,100)(35,82)
\path(32.5,82.6)(35,82)(33.5,84.2)
\put(7.2,100.5){\circle*{0.3}}
\put(7.7,101){{\xssc $ \xBc f(X_r),i_r \xBe $}}
\put(4,90){{\xssc $ \xBc \xbb,f,X_r,g \xBe $}}

\path(25,89)(95,95)
\path(27.9,88.2)(25,89)(27.7,90.26)
\put(95.8,95.2){\circle*{0.3}}
\put(90,98){{\xssc $ \xBc g(Y_s),j_s \xBe $}}
\put(55,95){{\xssc $ \xBc \xbg,f,X_r,g,Y_s \xBe $}}

\put(10,15) {{\rm\bf The construction}}

\end{picture}
}

\end{diagram}

\vspace{4mm}

(2)

Let $X_{r} \xbe \xdD ( \xba ).$ We have to show that no
$ \xBc  \xba,f \xBe $ is valid in $X_{r}.$ Fix $f.$

$ \xBc  \xba,f \xBe $ is in $X_{r},$ so we have to show that for at least one $g
\xbe \xbP ( \xdO, \xba )$
$ \xBc  \xbb,f,X_{r},g \xBe $ is valid in $X_{r},$ i.e. that for this $g,$ no
$ \xBc  \xbg,f,X_{r},g,Y_{s} \xBe : \xBc g(Y_{s}),j_{s} \xBe  \xcp  \xBc 
\xbb,f,X_{r},g \xBe,$
$Y_{s} \xbe \xdO ( \xba )$
attacks $ \xBc  \xbb,f,X_{r},g \xBe $ in $X_{r}.$

We define $g.$ Take $Y_{s} \xbe \xdO ( \xba ).$

Case 1: $ \xbm (Y_{s}) \xcc X_{r}$ or $f(X_{r}) \xce Y_{s}:$ choose
arbitrary $g(Y_{s}) \xbe \xbm (Y_{s}).$

Case 2: $ \xbm (Y_{s}) \xcC X_{r}$ and $f(X_{r}) \xbe Y_{s}:$ Choose
$g(Y_{s}) \xbe \xbm (Y_{s})-X_{r}.$

In Case 1, $ \xBc  \xbg,f,X_{r},g,Y_{s} \xBe $ does not exist, so it cannot
attack
$ \xBc  \xbb,f,X_{r},g \xBe.$

In Case 2,
$ \xBc  \xbg,f,X_{r},g,Y_{s} \xBe : \xBc g(Y_{s}),j_{s} \xBe  \xcp \xBc
\xbb,f,X_{r},g \xBe
$ is not in $X_{r},$ as $g(Y_{s}) \xce X_{r}.$

Thus, no $ \xBc  \xbg,f,X_{r},g,Y_{s} \xBe : \xBc g(Y_{s}),j_{s} \xBe  \xcp \xBc
\xbb
,f,X_{r},g \xBe,$ $Y_{s} \xbe \xdO ( \xba )$
attacks $ \xBc  \xbb,f,X_{r},g \xBe $
in $X_{r}.$

So $ \xcA  \xBc  \xba,f \xBe : \xBc y,j \xBe  \xcp  \xBc x,i \xBe $

$ \xDC \xDC y \xbe X_{r}$ $ \xch $
$ \xcE  \xBc  \xbb,f,X_{r},g \xBe : \xBc f(X_{r}),i_{r} \xBe  \xcp \xBc \xba,f
\xBe $

$ \xDC \xDC \xDC \xDC (f(X_{r}) \xbe \xbm (X_{r})$ $ \xcu $
$ \xCN \xcE  \xBc  \xbg,f,X_{r},g,Y_{s} \xBe : \xBc g(Y_{s}),j_{s} \xBe  \xcp
\xBc \xbb
,f,X_{r},g \xBe.g(Y_{s}) \xbe X_{r}).$

But $ \xBc  \xbb,f,X_{r},g \xBe $ was constructed only for
$ \xBc  \xba,f \xBe,$ so was $ \xBc  \xbg,f,X_{r},g,Y_{s} \xBe,$ and
there was no other $ \xBc  \xbg,i \xBe $ attacking
$ \xBc  \xbb,f,X_{r},g \xBe,$ so we are done.

(3)

Let $Y_{s} \xbe \xdO ( \xba ).$ We have to show that at least one
$ \xBc  \xba,f \xBe $ is valid in $Y_{s}.$

We define $f \xbe \xbP ( \xdD, \xba ).$ Take $X_{r}.$

If $ \xbm (X_{r}) \xcC Y_{s},$ choose $f(X_{r}) \xbe \xbm (X_{r})-Y_{s}.$
If $ \xbm (X_{r}) \xcc Y_{s},$ choose arbitrary $f(X_{r}) \xbe \xbm
(X_{r}).$

All attacks on $ \xBc x,f \xBe $ have the form
$ \xBc  \xbb,f,X_{r},g \xBe : \xBc f(X_{r}),i_{r} \xBe  \xcp \xBc \xba,f \xBe,$
$X_{r} \xbe
\xdD ( \xba ),$
$g \xbe \xbP ( \xdO, \xba ).$ We have to show that they are either not in
$Y_{s},$ or that they
are themselves attacked in $Y_{s}.$

Case 1: $ \xbm (X_{r}) \xcC Y_{s}.$ Then $f(X_{r}) \xce Y_{s},$ so
$ \xBc  \xbb,f,X_{r},g \xBe : \xBc f(X_{r}),i_{r} \xBe  \xcp  \xBc  \xba,f \xBe
$ is not
in $Y_{s}$ (for no $g).$

Case 2: $ \xbm (X_{r}) \xcc Y_{s}.$ Then $ \xbm (Y_{s}) \xcC X_{r}$ by $(
\xbm \xcc \xcd )$ and $f(X_{r}) \xbe Y_{s},$ so
$ \xBc  \xbb,f,X_{r},g \xBe : \xBc f(X_{r}),i_{r} \xBe  \xcp  \xBc  \xba,f \xBe
$ is in $Y_{s}$ (for
all $g).$ Take any $g \xbe \xbP ( \xdO, \xba ).$ As
$ \xbm (Y_{s}) \xcC X_{r}$ and $f(X_{r}) \xbe Y_{s},$
$ \xBc  \xbg,f,X_{r},g,Y_{s} \xBe : \xBc g(Y_{s}),j_{s} \xBe  \xcp \xBc
\xbb,f,X_{r},g \xBe
$ is defined, and
$g(Y_{s}) \xbe \xbm (Y_{s}),$ so it is in $Y_{s}$ (for all $g).$ Thus,
$ \xBc  \xbb,f,X_{r},g \xBe $ is attacked in $Y_{s}.$

Thus, for this $f,$ all $ \xBc  \xbb,f,X_{r},g \xBe $ are either not in $Y_{s},$
or
attacked in $Y_{s},$
thus for this $f,$ $ \xBc  \xba,f \xBe $ is valid in $Y_{s}.$

So for this $ \xBc x,i \xBe $

$ \xcE  \xBc  \xba,f \xBe : \xBc y,j \xBe  \xcp  \xBc x,i \xBe.y \xbe \xbm
(Y_{s})$ $ \xcu $

$ \xDC \xDC $ (a)
$ \xCN \xcE  \xBc  \xbb,f,X_{r},g \xBe : \xBc f(X_{r}),i \xBe  \xcp \xBc \xba,f
\xBe
.f(X_{r}) \xbe Y_{s}$

$ \xDC \xDC \xDC $ or

$ \xDC \xDC $ (b)
$ \xcA  \xBc  \xbb,f,X_{r},g \xBe : \xBc f(X_{r}),i \xBe  \xcp  \xBc  \xba,f
\xBe $

$ \xDC \xDC \xDC \xDC (f(X_{r}) \xbe Y_{s}$ $ \xch $

$ \xDC \xDC \xDC \xDC \xcE  \xBc  \xbg,f,X_{r},g,Y_{s} \xBe : \xBc
g(Y_{s}),j_{s} \xBe  \xcp
\xBc \xbb,f,X_{r},g \xBe.g(Y_{s}) \xbe \xbm (Y_{s})).$

As we made copies of $ \xba $ only, introduced only $ \xbb ' $s attacking
the $ \xba -$copies, and
$ \xbg ' $s attacking the $ \xbb ' $s, the construction is independent for
different $ \xba ' $s.

$ \xcz $
\\[3ex]

\bp

$\hspace{0.01em}$

% (+++ Orig. No.:  Proposition Level-3-Repr +++)

\label{Proposition Level-3-Repr}

Let $U$ be the universe, $ \xbm: \xdy \xcp \xdp (U).$

Then any $ \xbm $ satisfying $( \xbm \xcc ),$ $( \xcs ),$ $( \xbm CUM)$
(or, alternatively, $( \xbm \xcc )$ and
$( \xbm \xcc \xcd ))$ can be represented by a level 3 essentially smooth
structure.

\ep

\subparagraph{
Proof
}

$\hspace{0.01em}$

% (+++ Orig.:  Proof +++)

In stage one, consider as usual
$ \xdu:= \xBc  \xdx,\{ \xba_{i}:i \xbe I\} \xBe $ where
$ \xdx $ $:=$ $\{ \xBc x,f \xBe :$ $x \xbe U,$ $f \xbe \xbP \{ \xbm (X):$ $X
\xbe
\xdy,$ $x \xbe X- \xbm (X)\}\},$ and set
$ \xba: \xBc x',f'  \xBe  \xcp  \xBc x,f \xBe $ $: \xcj $ $x' \xbe ran(f).$

For stage two:

Any level 1 arrow
$ \xba: \xBc y,j \xBe  \xcp  \xBc x,i \xBe $ was introduced in stage one by some
$Y \xbe \xdy
$
s.t. $y \xbe \xbm (Y),$ $x \xbe Y- \xbm (Y).$ Do the construction
of Lemma \ref{Lemma Level-3-Constr} (page \pageref{Lemma Level-3-Constr})  for
all
level 1 arrows of $ \xdx $ in parallel or successively.

We have to show that the resulting structure represents $ \xbm $ and is
essentially smooth. (Level 3 is obvious.)

(1) Representation

Suppose $x \xbe Y- \xbm (Y).$ Then there was in stage 1 for all
$ \xBc x,i \xBe $ some $ \xba: \xBc y,j \xBe  \xcp  \xBc x,i \xBe,$
$y \xbe \xbm (Y).$ We examine the $y.$

If there is no $X$ s.t. $x \xbe \xbm (X),$ $y \xbe X,$ then there were no
$ \xbb ' $s and $ \xbg ' $s
introduced for this
$ \xba: \xBc y,j \xBe  \xcp  \xBc x,i \xBe,$ so $ \xba $ is valid.

If there is $X$ s.t. $x \xbe \xbm (X),$ $y \xbe X,$ consider
$ \xba: \xBc y,j \xBe  \xcp  \xBc x,i \xBe.$ So $X \xbe \xdD ( \xba ),$
$Y \xbe \xdO ( \xba ),$ so we did the construction
of Lemma \ref{Lemma Level-3-Constr} (page \pageref{Lemma Level-3-Constr}), and
by its
result, $ \xBc x,i \xBe $ is not minimal in $Y.$

Thus, in both cases, $ \xBc x,i \xBe $ is successfully attacked in $Y,$ and no
$ \xBc x,i \xBe $
is a minimal element in $Y.$

Suppose $x \xbe \xbm (X)$ (we change notation to conform
to Lemma \ref{Lemma Level-3-Constr} (page \pageref{Lemma Level-3-Constr}) ).
Fix $ \xBc x,i \xBe.$

If there is no $ \xba: \xBc y,j \xBe  \xcp  \xBc x,i \xBe,$ $y \xbe X,$ then
$ \xBc x,i \xBe $ is minimal in $X,$ and we are
done.

If there is $ \xba $ or
$ \xBc  \xba,k \xBe : \xBc y,j \xBe  \xcp  \xBc x,i \xBe,$ $y \xbe X,$ then $
\xba $ originated from
stage one
through some $Y$ s.t. $x \xbe Y- \xbm (Y),$ and $y \xbe \xbm (Y).$ (Note
that stage 2 of the
construction did not
introduce any new level 1 arrows - only copies of existing level 1
arrows.) So
$X \xbe \xdD ( \xba ),$ $Y \xbe \xdO ( \xba ),$ so we did the construction
of Lemma \ref{Lemma Level-3-Constr} (page \pageref{Lemma Level-3-Constr}), and
by its result,
$ \xBc x,i \xBe $ is minimal in $X,$ and we are done again.

In both cases, all
$ \xBc x,i \xBe $ are minimal elements in $X.$

(2) Essential smoothness. We have to show the conditions of
Definition \ref{Definition X-Sub-X'} (page \pageref{Definition X-Sub-X'}). We
will, however, work with the
reformulation given in Remark \ref{Remark X-Sub-X'} (page \pageref{Remark
X-Sub-X'}).

Case (1), $x \xbe \xbm (X).$

Case (1.1), there is $ \xBc x,i \xBe $ with no
$ \xBc  \xba,f \xBe : \xBc y,j \xBe  \xcp  \xBc x,i \xBe,$ $y \xbe X.$
There is nothing to show.

Case (1.2), for all $ \xBc x,i \xBe $ there is
$ \xBc  \xba,f \xBe : \xBc y,j \xBe  \xcp  \xBc x,i \xBe,$ $y \xbe X.$

$ \xba $ was introduced in stage 1 by some $Y$ s.t.
$x \xbe Y- \xbm (Y),$ $y \xbe X \xcs \xbm (Y),$ so $X \xbe \xdD ( \xba ),$
$Y \xbe \xdO ( \xba ).$
In the proof of Lemma \ref{Lemma Level-3-Constr} (page \pageref{Lemma
Level-3-Constr}), at the end of (2),
it was shown that

$ \xDC \xDC \xcE  \xBc  \xbb,f,X_{r},g \xBe : \xBc f(X_{r}),i_{r} \xBe  \xcp
\xBc \xba,f
\xBe $

$ \xDC \xDC \xDC \xDC (f(X_{r}) \xbe \xbm (X_{r})$ $ \xcu $

$ \xDC \xDC \xDC \xDC \xCN \xcE  \xBc  \xbg,f,X_{r},g,Y_{s} \xBe : \xBc
g(Y_{s}),j_{s} \xBe \xcp \xBc \xbb,f,X_{r},g \xBe.g(Y_{s}) \xbe X_{r}).$

By $f(X_{r}) \xbe \xbm (X_{r}),$ condition (1) of Remark \ref{Remark X-Sub-X'}
(page \pageref{Remark X-Sub-X'})  is true.

Case (2), $x \xce \xbm (Y).$ Fix
$ \xBc x,i \xBe.$ (We change notation back to $Y.)$

In stage 1, we constructed
$ \xba: \xBc y,j \xBe  \xcp  \xBc x,i \xBe,$ $y \xbe \xbm (Y),$ so $Y \xbe \xdO
( \xba ).$

If $ \xdD ( \xba )= \xCQ,$ then there is no attack on $ \xba,$ and the
condition (2) of
Remark \ref{Remark X-Sub-X'} (page \pageref{Remark X-Sub-X'})  is trivially
true.

If $ \xdD ( \xba ) \xEd \xCQ,$ we did the construction
of Lemma \ref{Lemma Level-3-Constr} (page \pageref{Lemma Level-3-Constr}), so

$ \xcE  \xBc  \xba,f \xBe : \xBc y,j \xBe  \xcp  \xBc x,i \xBe.y \xbe \xbm
(Y_{s})$ $ \xcu $

$ \xDC \xDC $ (a) $ \xCN \xcE  \xBc  \xbb,f,X_{r},g \xBe : \xBc f(X_{r}),i \xBe
\xcp \xBc \xba,f \xBe.f(X_{r}) \xbe Y_{s}$

$ \xDC \xDC \xDC $ or

$ \xDC \xDC $ (b) $ \xcA  \xBc  \xbb,f,X_{r},g \xBe : \xBc f(X_{r}),i \xBe \xcp
\xBc \xba,f \xBe $

$ \xDC \xDC \xDC \xDC (f(X_{r}) \xbe Y_{s}$ $ \xch $

$ \xDC \xDC \xDC \xDC \xcE  \xBc  \xbg,f,X_{r},g,Y_{s} \xBe : \xBc
g(Y_{s}),j_{s}
\xBe \xcp \xBc \xbb,f,X_{r},g \xBe.g(Y_{s}) \xbe \xbm (Y_{s}).$

As the only attacks on $ \xBc  \xba,f \xBe $ had the form
$ \xBc  \xbb,f,X_{r},g \xBe,$ and $g(Y_{s}) \xbe \xbm (Y_{s}),$
condition (2) of Remark \ref{Remark X-Sub-X'} (page \pageref{Remark X-Sub-X'}) 
is satisfied.

$ \xcz $
\\[3ex]

As said after Example \ref{Example Level-bigger-2} (page \pageref{Example
Level-bigger-2}),
we do not know if level 3 is necessary for
representation. We also do not know whether the same can be achieved with
level 3, or higher, totally smooth structures.
\section{
Translation to logic
}
\label{Section Reac-GenPref-Logic}
\index{Section Reac-GenPref-Logic}

We turn to the translation to logics.

\bp

$\hspace{0.01em}$

% (+++ Orig. No.:  Proposition Higher-Repr +++)

\label{Proposition Higher-Repr}

Let $ \xcn $ be a logic for $ \xdl.$ Set $T^{ \xdm }:=Th( \xbm_{ \xdm
}(M(T))),$ where $ \xdm $ is a
generalized preferential structure, and $ \xbm_{ \xdm }$ its choice
function. Then

(1) there is a level 2 preferential structure $ \xdm $ s.t. $ \ol{ \ol{T}
}=T^{ \xdm }$ iff
(LLE), (CCL), (SC) hold for all $T,T' \xcc \xdl.$

(2) there is a level 3 essentially smooth preferential structure
$ \xdm $ s.t. $ \ol{ \ol{T} }=T^{ \xdm }$ iff
(LLE), (CCL), (SC), $( \xcc \xcd )$ hold for all $T,T' \xcc \xdl.$

\ep

\subparagraph{
Proof
}

$\hspace{0.01em}$

% (+++ Orig.:  Proof +++)

The proof is an immediate consequence of Corollary \ref{Corollary Eta-Rho} (page
\pageref{Corollary Eta-Rho})
(2),
Fact \ref{Fact X-Sub-X'} (page \pageref{Fact X-Sub-X'}),
Proposition \ref{Proposition Level-3-Repr} (page \pageref{Proposition
Level-3-Repr}),
and Proposition \ref{Proposition Alg-Log} (page \pageref{Proposition Alg-Log}) 
(10) and (11).

(More precisely, for (2): Let $ \xdm $ be an essentially smooth
structure, then by
Definition \ref{Definition Essentially-Smooth} (page \pageref{Definition
Essentially-Smooth})  for all $X$ $ \xbm (X) \xes
X.$ Consider $( \xbm CUM).$
So by Fact \ref{Fact X-Sub-X'} (page \pageref{Fact X-Sub-X'})  (2)
$ \xbm (X' ) \xcc X'' \xcc X' $ $ \xch $ $ \xbm (X' ) \xes X'',$
so by Fact \ref{Fact X-Sub-X'} (page \pageref{Fact X-Sub-X'})  (1)
$ \xbm (X' )= \xbm (X'' ).$
$( \xbm \xcc \xcd )$ is analogous, using
Fact \ref{Fact X-Sub-X'} (page \pageref{Fact X-Sub-X'})  (3).

$ \xcz $
\\[3ex]

We leave aside the generalization of preferential structures to attacking
structures relative to $ \xbh,$ as this can cause problems, without
giving real
insight: It might well be that $ \xbr (X) \xcC \xbh (X),$ but, still, $
\xbr (X)$ and $ \xbh (X)$ might
define the same theory - see
e.g. Example \ref{Example Pref-Dp} (page \pageref{Example Pref-Dp}).

% ******* BEGIN LATEX SOURCE FILE 6-9-dur.tex *******
%
% Uebers. aus Karltex File: 6-9-dur.m
%
%
\chapter{Deontic logic and hierarchical conditionals}
\section{Semantics of deontic logic}
\label{Section Deon-Semantik}
\subsection{
Introductory remarks
}

We see some relation of ``better'' as central for obligations.
Obligations determine what is ``better'' and what
is ``worse'', conversely, given such a relation of
``better'', we can define obligations.
The problems lie, in our opinion, in the fact that an adequate treatment
of such a relation is somewhat complicated, and leads to many
ramifications.

On the other hand, obligations have sufficiently many things
in common with facts so we can in a useful way say that an
obligation is satisfied in a situation, and one can also define a notion
of derivation for obligations.

Our approach is almost exclusively semantical.
\subsubsection{
Context
}

The problem with formalisation using logic is that the natural
movements in the application area being formalised do not exactly
correspond to natural movements in the logic being used as a tool of
formalisation. Put differently, we may be able to express statement A
of the application area by a formula $ \xbf $ of the logic, but the
subtleties of the way A is manipulated in the application area cannot
be matched in the formal logic used.
This gives rise to paradoxes. To resolve the paradoxes one needs to
improve the logic.
So the progress in the formalisation program depends on the state of
development of logic itself.
Recent serious advances in logical tools by the authors of this paper
enable us to offer some better formalisations possibilities for the
notion of obligation. This is what we offer in this paper.

Historically, articles on Deontic Logic include collections of
problems, see e.g.  \cite{MDW94}, semantical
approaches, see e.g.  \cite{Han69},
and others, like  \cite{CJ02}.

Our basic idea is to see obligations as tightly connected to some
relation of ``better''. An immediate consequence is that negation, which
inverses such a relation, behaves differently in the case of obligations
and of classical logic. $('' And'' $ and ``or'' seem to show analogue
behaviour in
both logics.)
The relation of ``better'' has to be treated with some caution, however,
and we introduce and investigate local and global properties about
``better'' of
obligations. Most of these properties coincide in sufficiently nice
situations, in others, they are different.

We do not come to a final conclusion which properties obligations should
or should not have, perhaps this will be answered in future, perhaps there
is no
universal answer. We provide a list of ideas which seem reasonable to us.

Throughout, we work in a finite (propositional) setting.
\subsubsection{
Central idea
}

We see a tight connection between obligations and a relation of
``morally'' better between situations. Obligations are there to
guide us for ``better'' actions, and, conversely, given some
relation of ``better'', we can define obligations.

The problems lie in the fact that a simple approach via quality is not
satisfactory. We examine a number of somewhat more subtle ideas, some
of them also using a notion of distance.
\subsubsection{
A common property of facts and obligations
}

We are fully aware that an obligation has a conceptually different status
than a fact. The latter is true or false, an obligation has been set by
some
instance as a standard for behaviour, or whatever.

Still, we will say that an obligation holds in a situation, or that a
situation
satisfies an obligation. If the letter is in the mail box, the obligation
to
post it is satisfied, if it is in the trash bin, the obligation is not
satisfied. In some set of worlds, the obligation is satisfied, in the
complement
of this set, it is not satisfied. Thus, obligations behave in this respect
like facts, and we put for this aspect all distinctions between facts and
obligations aside.

Thus, we will treat obligations most of the time as subsets of the set of
all
models, but also sometimes as formulas. As we work mostly in a finite
propositional setting, both are interchangeable.

We are $ \xCf not$ concerned here with actions to fulfill obligations,
developments
or so, just simply situations which satisfy or not obligations.

This is perhaps also the right place to point out that one has to
distinguish
between facts that $ \xCf hold$ in ``good'' situations (they will be closed
under
arbitrary right weakening), and obligations which describe what $ \xCf
should$
be, they will not be closed under arbitrary right weakening. This article
is only about the latter.
\subsubsection{
Derivations of obligations
}

Again, we are aware that ``deriving'' obligations is different from
``deriving''
facts. Derivation of facts is supposed to conclude from truth to truth,
deriving obligations will be about concluding what can also be considered
an obligation, given some set of ``basic'' obligations. The parallel is
sufficiently strong to justify the double use of the word ``derive''.

Very roughly, we will say that conjunctions (or intersections) and
disjunctions (unions) of obligations lead in a reasonable way to derived
obligations, but negations do not. We take the Ross paradox (see below)
very seriously, it was, as a matter of fact, our starting point to avoid
it in a reasonable notion of derivation.

We mention two simple postulates derived obligations should probably
satisfy.

 \xEh

 \xDH Every original obligation should also be a derived obligation,
corresponding to $ \xba, \xbb \xcn \xba.$

 \xDH A derived obligation should not be trivial, i.e. neither empty nor
$U,$
the universe we work in.

 \xEj

The last property is not very important from an algebraic point of
view, and easily satisfiable, so we will not give it too much attention.
\subsubsection{
Orderings and obligations
}

There is, in our opinion, a deep connection between obligations and
orderings (and, in addition, distances), which works both ways.

First, given a set of obligations, we can say that one situation is
``better'' than a seond situation, if the first satisfies ``more'' obligations
than the second does. ``More'' can be measured by the set of obligations
satisfied, and then by the subset/superset relation, or by the number of
obligations.
Both are variants of the Hamming distance. ``Distance'' between two
situations can be measured by the set or number of obligations in which
they differ (i.e. one situation satisfies them, the other not).
In both cases, we will call the variants the set or the counting variant.

This is also the deeper reason why we have to be careful with negation.
Negation inverses such orderings, if $ \xbf $ is better than $ \xCN \xbf
,$ then $ \xCN \xbf $ is
worse than $ \xCN \xCN \xbf = \xbf.$ But in some reasonable sense $ \xcu
$ and $ \xco $ preserve the
ordering, thus they are compatibel with obligations.

Conversely, given a relation (of quality), we might for instance require
that obligations are closed under improvement. More subtle requirements
might work with distances.
The relations of quality and distance can be given abstractly (as the
notion of size used for ``soft'' obligations), or as above by a starting set
of
obligations. We will also define important
auxiliary concepts on such abstract relations.
\subsubsection{
Derivation revisited
}

A set of ``basic'' obligations generates an ordering and a distance between
situations, ordering and distance can be used to define properties
obligations should have. It is thus natural to define obligations derived
from the basic set as those sets of situations which satisfy the
desirable properties of obligations defined via the order and distance
generated
by the basic set of obligations.
Our derivation is thus a two step procedure: first generate the order and
distance, which define suitable sets of situations.

We will call properties which are defined without using distances global
properties (like closure under improving quality), properties involving
distance (like being a neighbourhood) local properties.
\subsubsection{
Relativization
}

An important problem is relativization. Suppose $ \xdo $ is a set of
obligations for
all possible situations, e.g. $O$ is the obligation to post the letter,
and
$O' $ is the obligation to water the flowers. Ideally, we will do both.
Suppose
we consider now a subset, where we cannot do both (e.g. for lack of time).
What are our obligations in this subset? Are they just the restrictions to
the subset? Conversely, if $O$ is an obligation for a subset of all
situations, is then some $O' $ with $O \xcc O' $ an obligation for the set
of all
situations?

In more complicated requirements, it might be reasonable e.g. to choose
the ideal situations still in the big set, even if they are not in the
subset to be considered, but use an approximation inside the subset.
Thus, relativizations present a non-trivial problem with many possible
solutions, and it seems doubtful whether a universally adequate solution
can be
found.
\subsubsection{
Numerous possibilities
}

Seeing the possibilities presented so far (set or counting order,
set or counting distance, various relativizations), we can already guess
that there are numerous possible reasonable approaches to what an
obligation
is or should be. Consequently, it seems quite impossible to pursue all
these combinations in detail. Thus, we concentrate mostly on one
combination, and leave it to the reader to fill in details for the others,
if (s)he is so interested.

We will also treat the defeasible case here. Perhaps somewhat
surprisingly,
this is straightforward, and largely due to the fact that there is one
natural
definition of ``big'' sets for the product set, given that ``big'' sets are
defined for the components. So there are no new possibilities to deal with
here.

The multitude of possibilities is somewhat disappointing. It may,
of course, be due to an incapacity of the present authors to find $ \xCf
the$
right notion. But it may also be a genuine problem of ambiguous
intuitions,
and thus generate conflicting views and misunderstandings on the one side,
and loopholes for not quite honest argumentation in practical juridical
reasoning on the other hand.
\subsubsection{
Notation
}

$ \xdp (X)$ will be the powerset of $X,$ $A \xcc B$ will mean that A is a
subset of $B,$ or
$A=B,$ $A \xcb B$ that A is a proper subset of $B.$
\subsubsection{
Overview
}

We will work in a finite propositional setting, so there is a trivial
and 1-1 correspondence between formulas and model sets. Thus, we can
just work with model sets - which implies, of course, that obligations
will be robust under logical reformulation. So we will formulate most
results only for sets.

 \xEI

 \xDH In Section \ref{Section Definitions} (page \pageref{Section Definitions})
,
we give the basic definitions, together with
some simple results about those definitions.

 \xDH Section \ref{Section Phil} (page \pageref{Section Phil})  will
present a more philosophical discussion, with more
examples, and we will argue that our definitions are relevant for our
purpose.

 \xDH As said already, there seems to be a multitude of possible and
reasonable definitions of what an obligation can or should be, so we
limit our formal investigation to a few cases, this is given in
Section \ref{Section Exam} (page \pageref{Section Exam}).

 \xDH In Section \ref{Section What} (page \pageref{Section What}), we give
a tentative definition of an obligation.

 \xEJ

(The concept of
neighbourhood semantics is not new, and was already
introduced by D.Scott,  \cite{Sco70}, and R.Montague,  \cite{Mon70}. Further
investigations showed that it was also used by O.Pacheco,  \cite{Pac07},
precisely to avoid unwanted weakening for obligations. We came the other
way
and started with the concept of independent strengthening, see
Definition \ref{Definition Delta-O} (page \pageref{Definition Delta-O}), and
introduced the abstract
concept of neighbourhood semantics only at the end. This is one of the
reasons we also have different descriptions which turned out to be
equivalent:
we came from elsewhere.)
\subsection{
Basic definitions
}
\label{Section Definitions}

We give here all definitions needed for our treatment of obligations.
The reader may continue immediately
to Section \ref{Section Phil} (page \pageref{Section Phil}), and come
back to the present section whenever necessary.

Intuitively, $U$ is supposed to be the set of
models of some propositional language $ \xdl,$ but we will stay purely
algebraic
whenever possible. $U' \xcc U$ is some subset.

For ease of writing, we will here and later sometimes work with
propositional
variables, and also identify models with the formula describing them,
still
this is just shorthand for arbitrary sets and elements. pq will stand for
$p \xcu q,$ etc.

If a set $ \xdo $ of obligations is given, these will be just arbitrary
subsets
of the universe $U.$ We will also say that $ \xdo $ is over $U.$

Before we deepen the discussion of more conceptual aspects, we give some
basic
definitions (for which we claim no originality). We will need them quite
early,
and think it is better to put them together here, not to be read
immediately,
but so the reader can leaf back and find them easily.

We work here with a notion of size (for the defeasible case), a notion $d$
of
distance, and a quality relation $ \xec.$ The latter two can be given
abstractly, but
may also be defined from a set of (basic) obligations.

We use these notions to describe properties obligations should, in our
opinion,
have. A careful analysis will show later interdependencies between the
different
properties.
\subsubsection{
Product Size
}
\label{Section Product-Size}

For each $U' \xcc U$ we suppose an abstract notion of size to be given.
We may assume this notion to be a filter or an ideal. Coherence properties
between the filters/ideals of different $U',$ $U'' $ will be left open,
the reader may assume them to be the conditions of the system $P$ of
preferential
logic,

see Section \ref{Section Basic-definitions-and-results} (page \pageref{Section
Basic-definitions-and-results}).

Given such notions of size on $U' $ and $U'',$ we will also need a notion
of size
on $U' \xCK U''.$ We take the evident solution:

\bd

$\hspace{0.01em}$

% (+++ Orig. No.:  Definition Size-Product +++)

\label{Definition Size-Product}

Let a notion of ``big subset'' be defined by a principal filter for all $X
\xcc U$ and
all $X' \xcc U'.$ Thus, for all $X \xcc U$ there exists a fixed principal
filter
$ \xdf (X) \xcc \xdp (X),$ and likewise for all $X' \xcc U'.$ (This is
the situation in the case
of preferential structures, where $ \xdf (X)$ is generated by $ \xbm (X),$
the set of
minimal elements of $X.)$

Define now $ \xdf (X \xCK X' )$ as generated by $\{A \xCK A':$ $A \xbe
\xdf (X),$ $A' \xbe \xdf (X' )\},$ i.e.
if A is the smallest element of $ \xdf (X),$ $A' $ the smallest element of
$ \xdf (X' ),$
then $ \xdf (X \xCK X' ):=\{B \xcc X \xCK X':$ $A \xCK A' \xcc B\}.$

\ed

\bfa

$\hspace{0.01em}$

% (+++ Orig. No.:  Fact Product-Pref +++)

\label{Fact Product-Pref}

If $ \xdf (X)$ and $ \xdf (X' )$ are generated by preferential structures
$ \xec_{X},$ $ \xec_{X' },$
then $ \xdf (X \xCK X' )$ is
generated by the product structure defined by

$ \xBc x,x'  \xBe  \xec_{X \xCK X' } \xBc y,y'  \xBe $ $: \xcj $ $x \xec_{X}y$
and $x'
\xec_{X' }y'.$

\efa

\subparagraph{
Proof
}

$\hspace{0.01em}$

% (+++ Orig.:  Proof +++)

We will omit the indices of the orderings when this causes no confusion.

Let $A \xbe \xdf (X),$ $A' \xbe \xdf (X' ),$ i.e. A minimizes $X,$ $A' $
minimizes $X'.$ Let
$ \xBc x,x'  \xBe  \xbe X \xCK X',$ then there are $a \xbe A,$ $a' \xbe A' $
with
$a \xec x,$ $a' \xec x',$ so
$ \xBc a,a'  \xBe  \xec  \xBc x,x'  \xBe.$

Conversely, suppose $U \xcc X \xCK X',$ $U$ minimizes $X \xCK X'.$ but
there is no $A \xCK A' \xcc U$
s.t. $A \xbe \xdf (X),$ $A' \xbe \xdf (X' ).$ Assume $A= \xbm (X),$ $A' =
\xbm (X' ),$ so there is
$ \xBc a,a'  \xBe  \xbe \xbm (X) \xCK \xbm (X' ),$ $ \xBc a,a'  \xBe  \xce U.$
But only
$ \xBc a,a'  \xBe  \xec  \xBc a,a'  \xBe,$ and $U$ does not
minimize $X \xCK X',$ $contradiction.$

$ \xcz $
\\[3ex]

Note that a natural modification of our definition:

There is $A \xbe \xdf (X)$ s.t.
for all $a \xbe A$ there is a (maybe varying) $A'_{a} \xbe \xdf (X' ),$
and
$U:=\{ \xBc a,a'  \xBe :$ $a \xbe A,$ $a' \xbe A'_{a}\}$ as generating sets

will result in the same
definition, as our filters are principal, and thus stable under arbitrary
intersections.
\subsubsection{
Distance
}

We consider a set of sequences $ \xbS,$ for $x \xbe \xbS $ $x:I \xcp S,$
$I$ a finite index set, $S$
some set.
Often, $S$ will be $\{0,1\},$ $x(i)=1$ will mean that $x \xbe i,$ when $I
\xcc \xdp (U)$ and $x \xbe U.$
For abbreviation, we will call this (unsystematically, often context will
tell) the $ \xbe -$case.
Often, $I$ will be written $ \xdo,$ intuitively, $O \xbe \xdo $ is then
an obligation, and
$x(O)=1$ means $x \xbe O,$ or $x$ ``satisfies'' the obligation $O.$

\bd

$\hspace{0.01em}$

% (+++ Orig. No.:  Definition O-x +++)

\label{Definition O-x}

In the $ \xbe -$case, set $ \xdo (x):=\{O \xbe \xdo:x \xbe O\}.$

\ed

\bd

$\hspace{0.01em}$

% (+++ Orig. No.:  Definition Hamming-Distance +++)

\label{Definition Hamming-Distance}

Given $x,y \xbe \xbS,$ the Hamming distance comes in two flavours:

$d_{s}(x,y):=\{i \xbe I:x(i) \xEd y(i)\},$ the set variant,

$d_{c}(x,y):=card(d_{s}(x,y)),$ the counting variant.

We define $d_{s}(x,y) \xck d_{s}(x',y' )$ iff $d_{s}(x,y) \xcc d_{s}(x'
,y' ),$

thus, $s-$distances are not always comparabel.

\ed

\bfa

$\hspace{0.01em}$

% (+++ Orig. No.:  Fact Distance +++)

\label{Fact Distance}

(1) In the $ \xbe -$case, we have $d_{s}(x,y)= \xdo (x) \xeY \xdo (y),$
where $ \xeY $ is the symmetric
set difference.

(2) $d_{c}$ has the normal addition, set union takes the role of addition
for $d_{s},$
$ \xCQ $ takes the role of 0 for $d_{s},$
both are distances in the following sense:

(2.1) $d(x,y)=0$ if $x=y,$ but not conversely,

(2.2) $d(x,y)=d(y,x),$

(2.3) the triangle inequality holds, for the set variant in the form
$d_{s}(x,z) \xcc d_{s}(x,y) \xcv d_{s}(y,z).$

(If $d(x,y)=0$ $ \xcH $ $x=y$ poses a problem, one can always consider
equivalence
classes.)

\efa

\subparagraph{
Proof
}

$\hspace{0.01em}$

% (+++ Orig.:  Proof +++)

(2.1) Suppose $U=\{x,y\},$ $ \xdo =\{U\},$ then $ \xdo (X)= \xdo (Y),$ but
$x \xEd y.$

(2.3) If $i \xce d_{s}(x,y) \xcv d_{s}(y,z),$ then $x(i)=y(i)=z(i),$ so
$x(i)=z(i)$ and
$i \xce d_{s}(x,z).$

The others are trivial.

$ \xcz $
\\[3ex]

\bd

$\hspace{0.01em}$

% (+++ Orig. No.:  Definition Between +++)

\label{Definition Between}

(1) We can define for any distance $d$ with some minimal requirements a
notion of
``between''.

If the codomain of $d$ has an ordering $ \xck,$ but no addition, we
define:

$ \xBc x,y,z \xBe _{d}$ $: \xcj $ $d(x,y) \xck d(x,z)$ and $d(y,z) \xck d(x,z).$

If the codomain has a commutative addition, we define

$ \xBc x,y,z \xBe _{d}$ $: \xcj $ $d(x,z)=d(x,y)+d(y,z)$ - in $d_{s}$ $+$ will
be
replaced by $ \xcv,$ i.e.

$ \xBc x,y,z \xBe _{s}$ $: \xcj $ $d(x,z)=d(x,y) \xcv d(y,z).$

For above two Hamming distances, we will write $ \xBc x,y,z \xBe _{s}$ and
$ \xBc x,y,z \xBe _{c}.$

(2) We further define:

$[x,z]_{d}:=\{y \xbe X: \xBc x,y,x \xBe _{d}\}$ - where $X$ is the set we work
in.

We will write $[x,z]_{s}$ and $[x,z]_{c}$ when appropriate.

(3) For $x \xbe U,$ $X \xcc U$ set $x \xFO_{d}X$ $:=$ $\{x' \xbe X: \xCN
\xcE x'' \xEd x' \xbe X.d(x,x' ) \xcg d(x,x'' )\}.$

Note that, if $X \xEd \xCQ,$ then $x \xFO X \xEd \xCQ.$

We omit the index when this does not cause confusion. Again, when
adequate,
we write $ \xFO_{s}$ and $ \xFO_{c}.$

\ed

For problems with characterizing ``between'' see
 \cite{Sch04}.

\bfa

$\hspace{0.01em}$

% (+++ Orig. No.:  Fact Between +++)

\label{Fact Between}

(0) $ \xBc x,y,z \xBe _{d}$ $ \xcj $ $ \xBc z,y,x \xBe _{d}.$

Consider the situation of a set of sequences $ \xbS.$

Let $A:=A_{ \xbs, \xbs '' }:=\{ \xbs ': \xcA i \xbe I( \xbs (i)= \xbs ''
(i) \xcp \xbs ' (i)= \xbs (i)= \xbs '' (i))\}.$
Then

(1) If $ \xbs ' \xbe A,$ then $d_{s}( \xbs, \xbs '' )=d_{s}( \xbs, \xbs
' ) \xcv d_{s}( \xbs ', \xbs '' ),$ so $ \xBc  \xbs, \xbs ', \xbs ''  \xBe
_{s}.$

(2) If $ \xbs ' \xbe A$ and $S$ consists of 2 elements (as in classical
2-valued
logic), then $d_{s}( \xbs, \xbs ' )$ and $d_{s}( \xbs ', \xbs '' )$ are
disjoint.

(3) $[ \xbs, \xbs '' ]_{s}=A.$

(4) If, in addition, $S$ consists of 2 elements, then $[ \xbs, \xbs ''
]_{c}=A.$

\efa

\subparagraph{
Proof
}

$\hspace{0.01em}$

% (+++ Orig.:  Proof +++)

(0) Trivial.

(1) `` $ \xcc $ '' follows from Fact \ref{Fact Distance} (page \pageref{Fact
Distance}), (2.3).

Conversely, if e.g. $i \xbe d_{s}( \xbs, \xbs ' ),$ then
by prerequisite $i \xbe d_{s}( \xbs, \xbs '' ).$

(2) Let $i \xbe d_{s}( \xbs, \xbs ' ) \xcs d_{s}( \xbs ', \xbs '' ),$
then $ \xbs (i) \xEd \xbs ' (i)$ and $ \xbs ' (i) \xEd \xbs '' (i),$
but then by $card(S)=2$ $ \xbs (i)= \xbs '' (i),$ but $ \xbs ' \xbe A,$
$contradiction.$

We turn to (3) and (4):

If $ \xbs ' \xce A,$ then there is $i' $ s.t. $ \xbs (i' )= \xbs '' (i' )
\xEd \xbs ' (i' ).$ On the other hand,
for all $i$ s.t. $ \xbs (i) \xEd \xbs '' (i)$ $i \xbe d_{s}( \xbs, \xbs '
) \xcv d_{s}( \xbs ', \xbs '' ).$ Thus:

(3) By (1) $ \xbs ' \xbe A$ $ \xch $ $ \xBc  \xbs, \xbs ', \xbs ''  \xBe _{s}.$
Suppose $ \xbs ' \xce A,$ so there is $i' $ s.t.
$i' \xbe d_{s}( \xbs, \xbs ' )-d_{s}( \xbs, \xbs '' ),$ so
$ \xBc  \xbs, \xbs ', \xbs ''  \xBe _{s}$ cannot be.

(4) By (1) and (2) $ \xbs ' \xbe A$ $ \xch $
$ \xBc  \xbs, \xbs ', \xbs ''  \xBe _{c}.$ Conversely, if $ \xbs ' \xce A,$ then
$card(d_{s}( \xbs, \xbs ' ))+card(d_{s}( \xbs ', \xbs '' )) \xcg
card(d_{s}( \xbs, \xbs '' ))+2.$

$ \xcz $
\\[3ex]

\bd

$\hspace{0.01em}$

% (+++ Orig. No.:  Definition L-And +++)

\label{Definition L-And}

Given a finite propositional laguage $ \xdl $ defined by the set $v( \xdl
)$ of
propositional
variables, let $ \xdl_{ \xcu }$ be the set of all consistent conjunctions
of
elements from $v( \xdl )$ or their negations. Thus, $p \xcu \xCN q \xbe
\xdl_{ \xcu }$ if $p,q \xbe v( \xdl ),$ but
$p \xco q,$ $ \xCN (p \xcu q) \xce \xdl_{ \xcu }.$ Finally, let $ \xdl_{
\xco \xcu }$ be the set of all (finite)
disjunctions of formulas from $ \xdl_{ \xcu }.$ (As we will later not
consider all
formulas from $ \xdl_{ \xcu },$ this will be a real restriction.)

Given a set of models $M$ for a finite language $ \xdl,$ define
$ \xbf_{M}:= \xcU \{p \xbe v( \xdl ): \xcA m \xbe M.m(p)=v\} \xcu \xcU \{
\xCN p:p \xbe v( \xdl ), \xcA m \xbe M.m(p)=f\} \xbe \xdl_{ \xcu }.$
(If there are no such $p,$ set $ \xbf_{M}:=TRUE.)$

This is the strongest $ \xbf \xbe \xdl_{ \xcu }$ which holds in $M.$

\ed

\bfa

$\hspace{0.01em}$

% (+++ Orig. No.:  Fact Hamming-Neighbourhood +++)

\label{Fact Hamming-Neighbourhood}

If $x,y$ are models, then $[x,y]=M( \xbf_{\{x,y\}}).$ $ \xcz $
\\[3ex]

\efa

\subparagraph{
Proof
}

$\hspace{0.01em}$

% (+++ Orig.:  Proof +++)

$m \xbe [x,y]$ $ \xcj $ $ \xcA p(x \xcm p,y \xcm p \xch m \xcm p$ and $x
\xcM p,y \xcM p \xch m \xcM p),$
$m \xcm \xbf_{\{x,y\}}$ $ \xcj $ $m \xcm \xcU \{p:x(p)=y(p)=v\} \xcu \xcU
\{ \xCN p:x(p)=y(p)=f\}.$
\subsubsection{
Quality and closure
}

\bd

$\hspace{0.01em}$

% (+++ Orig. No.:  Definition Closed +++)

\label{Definition Closed}

Given any relation $ \xec $ (of quality), we say that $X \xcc U$ is
(downward) closed (with
respect to $ \xec )$ iff $ \xcA x \xbe X \xcA y \xbe U(y \xec x$ $ \xch $
$y \xbe X).$

\ed

(Warning, we follow the preferential tradition, ``smaller'' will
mean ``better''.)

We think that being closed is a desirable property for obligations: what
is
at least as good as one element in the obligation should be ``in'', too.

\bfa

$\hspace{0.01em}$

% (+++ Orig. No.:  Fact Subset-Closure +++)

\label{Fact Subset-Closure}

Let $ \xec $ be given.

(1) Let $D \xcc U' \xcc U'',$ $D$ closed in $U'',$ then $D$ is also
closed in $U'.$

(2) Let $D \xcc U' \xcc U'',$ $D$ closed in $U',$ $U' $ closed in $U''
,$ then $D$ is
closed in $U''.$

(3) Let $D_{i} \xcc U' $ be closed for all $i \xbe I,$ then so are $ \xcV
\{D_{i}:i \xbe I\}$ and $ \xcS \{D_{i}:i \xbe I\}.$

\efa

\subparagraph{
Proof
}

$\hspace{0.01em}$

% (+++ Orig.:  Proof +++)

(1) Trivial.

(2) Let $x \xbe D \xcc U',$ $x' \xec x,$ $x' \xbe U'',$ then $x' \xbe U'
$ by closure of $U'',$ so $x' \xbe D$ by
closure of $U'.$

(3) Trivial.

$ \xcz $
\\[3ex]

We may have an abstract relation $ \xec $ of quality on the domain, but we
may also
define it from the structure of the sequences, as we will do now.

\bd

$\hspace{0.01em}$

% (+++ Orig. No.:  Definition Quality +++)

\label{Definition Quality}

Consider the case of sequences.

Given a relation $ \xec $ (of quality) on the codomain, we extend this to
sequences
in $ \xbS:$

$x \xCq y$ $: \xcj $ $ \xcA i \xbe I(x(i) \xCq y(i))$

$x \xec y$ $: \xcj $ $ \xcA i \xbe I(x(i) \xec y(i))$

$x \xeb y$ $: \xcj $ $ \xcA i \xbe I(x(i) \xec y(i))$ and $ \xcE i \xbe
I(x(i) \xeb y(i))$

In the $ \xbe -$case, we will consider $x \xbe i$ better than $x \xce i.$
As we have only two
values, true/false, it is easy to count the positive and negative cases
(in more complicated situations, we might be able to multiply), so we have
an analogue of the two Hamming distances, which we might call the Hamming
quality relations.

Let $ \xdo $ be given now.

(Recall that we follow the preferential tradition, ``smaller'' will
mean ``better''.)

$x \xCq_{s}y$ $: \xcj $ $ \xdo (x)= \xdo (y),$

$x \xec_{s}y$ $: \xcj $ $ \xdo (y) \xcc \xdo (x),$

$x \xeb_{s}y$ $: \xcj $ $ \xdo (y) \xcb \xdo (x),$

$x \xCq_{c}y$ $: \xcj $ $card( \xdo (x))=card( \xdo (y)),$

$x \xec_{c}y$ $: \xcj $ $card( \xdo (y)) \xck card( \xdo (x)),$

$x \xeb_{c}y$ $: \xcj $ $card( \xdo (y))<card( \xdo (x)).$

\ed

The requirement of closure causes a problem for
the counting approach: Given e.g. two obligations $O,$ $O',$ then any two
elements
in just one obligation have the same quality, so if one is in, the other
should
be, too. But this prevents now any of the original obligations to have the
desirable property of closure. In the counting case, we will obtain a
ranked
structure, where elements satisfy
0, 1, 2, etc. obligations, and we are unable to differentiate inside those
layers. Moreover, the set variant seems to be closer to logic, where we do
not
count the propositional variables which hold in a model, but consider them
individually. For these reasons, we will not pursue the counting approach
as
systematically as the set approach. One should, however, keep in mind that
the
counting variant gives a ranking relation of quality, as all qualities are
comparable, and the set variant does not. A ranking seems to be
appreciated
sometimes in the literature, though we are not really sure why.

Of particular interest is the combination of $d_{s}$ and $ \xec_{s}$
$(d_{c}$ and $ \xec_{c})$
respectively - where by $ \xec_{s}$ we also mean $ \xeb_{s}$ and $
\xCq_{s},$ etc. We turn to
this now.

\bfa

$\hspace{0.01em}$

% (+++ Orig. No.:  Fact Quality-Distance +++)

\label{Fact Quality-Distance}

We work in the $ \xbe -$case.

(1) $x \xec_{s}y$ $ \xch $ $d_{s}(x,y)= \xdo (x)- \xdo (y)$

Let $a \xeb_{s}b \xeb_{s}c.$ Then

(2) $d_{s}(a,b)$ and $d_{s}(b,c)$ are not comparable,

(3) $d_{s}(a,c)=d_{s}(a,b) \xcv d_{s}(b,c),$ and thus $b \xbe [a,c]_{s}.$

This does not hold in the counting variant, as Example \ref{Example Count} (page
\pageref{Example Count})  shows.

(4) Let $x \xeb_{s}y$ and $x' \xeb_{s}y$ with $x,x' $ $
\xeb_{s}-$incomparabel. Then $d_{s}(x,y)$ and $d_{s}(x',y)$
are incomparable.

(This does not hold in the counting variant, as then all distances are
comparable.)

(5) If $x \xeb_{s}z,$ then for all $y \xbe [x,z]_{s}$ $x \xec_{s}y
\xec_{s}z.$

\efa

\subparagraph{
Proof
}

$\hspace{0.01em}$

% (+++ Orig.:  Proof +++)

(1) Trivial.

(2) We have $ \xdo (c) \xcb \xdo (b) \xcb \xdo (a),$ so the results
follows from (1).

(3) By definition of $d_{s}$ and (1).

(4) $x$ and $x' $ are $ \xec_{s}-incomparable,$ so there are $O \xbe \xdo
(x)- \xdo (x' ),$
$O' \xbe \xdo (x' )- \xdo (x).$

As $x,x' \xeb_{s}y,$ $O,O' \xce \xdo (y),$ so $O \xbe d_{s}(x,y)-d_{s}(x'
,y),$
$O' \xbe d_{s}(x',y)-d_{s}(x,y).$

(5) $x \xeb_{s}z$ $ \xch $ $ \xdo (z) \xcb \xdo (x),$ $d_{s}(x,z)= \xdo
(x)- \xdo (z).$ By prerequisite
$d_{s}(x,z)=d_{s}(x,y) \xcv d_{s}(y,z).$
Suppose $x \xeC_{s}y.$ Then there is $i \xbe \xdo (y)- \xdo (x) \xcc
d_{s}(x,y),$ so
$i \xce \xdo (x)- \xdo (z)=d_{s}(x,z),$ $contradiction.$

Suppose $y \xeC_{s}z.$ Then there is $i \xbe \xdo (z)- \xdo (y) \xcc
d_{s}(y,z),$ so
$i \xce \xdo (x)- \xdo (z)=d_{s}(x,z),$ $contradiction.$

$ \xcz $
\\[3ex]

\be

$\hspace{0.01em}$

% (+++ Orig. No.:  Example Count +++)

\label{Example Count}

In this and similar examples, we will use the model notation. Some
propositional variables $p,$ $q,$ etc. are given, and models are described
by
$p \xCN qr,$ etc. Moreover, the propositional variables are the
obligations, so
in this example we have the obligations $p,$ $q,$ $r.$

Consider $x:= \xCN p \xCN qr,$ $y:=pq \xCN r,$ $z:= \xCN p \xCN q \xCN r.$
Then $y \xeb_{c}x \xeb_{c}z,$ $d_{c}(x,y)=3,$
$d_{c}(x,z)=1,$ $d_{c}(z,y)=2,$ so $x \xce [y,z]_{c}.$ $ \xcz $
\\[3ex]

\ee

\bd

$\hspace{0.01em}$

% (+++ Orig. No.:  Definition Quality-Extension +++)

\label{Definition Quality-Extension}

Given a quality relation $ \xeb $ between elements, and a distance $d,$ we
extend the
quality relation to sets and define:

(1) $x \xeb Y$ $: \xcj $ $ \xcA y \xbe (x \xFO Y).x \xeb y.$ (The closest
elements - i.e. there are no
closer ones - of $Y,$ seen from $x,$ are less good than $x.)$

analogously $X \xeb y$ $: \xcj $ $ \xcA x \xbe (y \xFO X).x \xeb y$

(2) $X \xeb_{l}Y$ $: \xcj $ $ \xcA x \xbe X.x \xeb Y$ and $ \xcA y \xbe
Y.X \xeb y$
(X is locally better than $Y).$

When necessary, we will write $ \xeb_{l,s}$ or $ \xeb_{l,c}$ to
distinguish the
set from the counting variant.

For the next definition, we use the notion of size: $ \xeA \xbf $ iff for
almost all $ \xbf $
holds i.e. the set of exceptions is small.

(3) $X \xDc_{l}Y$ $: \xcj $ $ \xeA x \xbe X.x \xeb Y$ and $ \xeA y \xbe
Y.X \xeb y.$

We will likewise write $ \xDc_{l,s}$ etc.

This definition is supposed to capture quality difference under minimal
change,
the ``ceteris paribus'' idea: $X \xeb_{l} \xdC X$ should hold for an
obligation $X.$
Minimal change is coded by $ \xFO,$ and ``ceteris paribus'' by minimal
change.

\ed

\bfa

$\hspace{0.01em}$

% (+++ Orig. No.:  Fact General-Obligation +++)

\label{Fact General-Obligation}

If $X \xeb_{l} \xdC X,$ and $x \xbe U$ an optimal point (there is no
better one), then $x \xbe X.$

\efa

\subparagraph{
Proof
}

$\hspace{0.01em}$

% (+++ Orig.:  Proof +++)

If not, then take $x' \xbe X$ closest to $x,$ this must be better than
$x,$ contradiction.
$ \xcz $
\\[3ex]

\bfa

$\hspace{0.01em}$

% (+++ Orig. No.:  Fact Local-Global +++)

\label{Fact Local-Global}

Take the set version.

If $X \xeb_{l,s} \xdC X,$ then $X$ is downward $ \xeb_{s}-closed.$

\efa

\subparagraph{
Proof
}

$\hspace{0.01em}$

% (+++ Orig.:  Proof +++)

Suppose $X \xeb_{l,s} \xdC X,$ but $X$ is not downward closed.

Case 1: There are $x \xbe X,$ $y \xce X,$ $y \xCq_{s}x.$ Then $y \xbe x
\xFO_{s} \xdC X,$ but $x \xeB y,$ $contradiction.$

Case 2:
There are $x \xbe X,$ $y \xce X,$
$y \xeb_{s}x.$ By $X \xeb_{l,s} \xdC X,$ the elements in $X$ closest to
$y$ must be better than $y.$
Thus, there is $x' \xeb_{s}y,$ $x' \xbe X,$ with minimal distance from
$y.$ But then
$x' \xeb_{s}y \xeb_{s}x,$ so $d_{s}(x',y)$ and $d_{s}(y,x)$ are
incomparable by
Fact \ref{Fact Quality-Distance} (page \pageref{Fact Quality-Distance}), so $x$
is among those
with minimal distance from $y,$ so $X \xeb_{l,s} \xdC X$ does not hold. $
\xcz $
\\[3ex]

\be

$\hspace{0.01em}$

% (+++ Orig. No.:  Example Dependent-2 +++)

\label{Example Dependent-2}

We work with the set variant.

This example shows that $ \xec_{s}-$closed does not imply $X \xeb_{l,s}
\xdC X,$ even
if $X$ contains the best elements.

Let $ \xdo:=\{p,q,r,s\},$ $U':=\{x:=p \xCN q \xCN r \xCN s,$ $y:= \xCN
pq \xCN r \xCN s,$ $x':=pqrs\},$ $X:=\{x,x' \}.$
$x' $ is the best element of $U',$ so $X$ contains the best elements,
and $X$ is downward closed in $U',$
as $x$ and $y$ are not comparable. $d_{s}(x,y)=\{p,q\},$ $d_{s}(x'
,y)=\{p,r,s\},$ so the
distances from $y$ are not comparable, so $x$ is among the closest
elements in $X,$
seen from $y,$ but $x \xeB_{s}y.$

The lack of comparability is essential here, as the following Fact shows.

$ \xcz $
\\[3ex]

\ee

We have, however, for the counting variant:

\bfa

$\hspace{0.01em}$

% (+++ Orig. No.:  Fact Count-Closed +++)

\label{Fact Count-Closed}

Consider the counting variant. Then

If $X$ is downward closed, then $X \xeb_{l,c} \xdC X.$

\efa

\subparagraph{
Proof
}

$\hspace{0.01em}$

% (+++ Orig.:  Proof +++)

Take any $x \xbe X,$ $y \xce X.$ We have $y \xec_{c}x$ or $x \xeb_{c}y,$
as any two elements are
$ \xec_{c}-$comparabel. $y \xec_{c}x$ contradicts closure, so $x
\xeb_{c}y,$ and $X \xeb_{l,c} \xdC X$ holds
trivially. $ \xcz $
\\[3ex]
\subsubsection{
Neighbourhood
}

\bd

$\hspace{0.01em}$

% (+++ Orig. No.:  Definition Neighbourhood +++)

\label{Definition Neighbourhood}

Given a distance $d,$ we define:

(1) Let $X \xcc Y \xcc U',$ then $Y$ is a neighbourhood of $X$ in $U' $
iff

$ \xcA y \xbe Y \xcA x \xbe X(x$ is closest to $y$ among all $x' $ with
$x' \xbe X$ $ \xch $ $[x,y] \xcs U' \xcc Y).$

(Closest means that there are no closer ones.)

When we also have a quality relation $ \xeb,$ we define:

(2) Let $X \xcc Y \xcc U',$ then $Y$ is an improving neighbourhood of $X$
in $U' $ iff

$ \xcA y \xbe Y \xcA x((x$ is closest to $y$ among all $x' $ with $x' \xbe
X$ and $x' \xec y)$ $ \xch $ $[x,y] \xcs U' \xcc Y).$

When necessary, we will have to say for (3) and (4) which variant, i.e.
set or counting, we mean.

\ed

\bfa

$\hspace{0.01em}$

% (+++ Orig. No.:  Fact Neighbourhood +++)

\label{Fact Neighbourhood}

(1) If $X \xcc X' \xcc \xbS,$ and $d(x,y)=0$ $ \xch $ $x=y,$ then $X$ and
$X' $ are Hamming neighbourhoods
of $X$ in $X'.$

(2) If $X \xcc Y_{j} \xcc X' \xcc \xbS $ for $j \xbe J,$ and all $Y_{j}$
are Hamming Neighbourhoods
of $X$ in $X',$ then so are $ \xcV \{Y_{j}:j \xbe J\}$ and $ \xcS
\{Y_{j}:j \xbe J\}.$

\efa

\subparagraph{
Proof
}

$\hspace{0.01em}$

% (+++ Orig.:  Proof +++)

(1) is trivial (we need here that $d(x,y)=0$ $ \xch $ $x=y).$

(2) Trivial.

$ \xcz $
\\[3ex]
\subsubsection{
Unions of intersections and other definitions
}

\bd

$\hspace{0.01em}$

% (+++ Orig. No.:  Definition ui +++)

\label{Definition ui}

Let $ \xdo $ over $U$ be given.

$X \xcc U' $ is $ \xCf (ui)$ (for union of intersections) iff there is a
family $ \xdo_{i} \xcc \xdo,$
$i \xbe I$ s.t. $X=( \xcV \{ \xcS \xdo_{i}:i \xbe I\}) \xcs U'.$

\ed

Unfortunately, this definition is not very useful for simple
relativization.

\bd

$\hspace{0.01em}$

% (+++ Orig. No.:  Definition Validity +++)

\label{Definition Validity}

Let $ \xdo $ be over $U.$ Let $ \xdo ' \xcc \xdo.$ Define for $m \xbe U$
and $ \xbd: \xdo ' \xcp 2=\{0,1\}$

$m \xcm \xbd $ $: \xcj $ $ \xcA O \xbe \xdo ' (m \xbe O \xcj \xbd (O)=1)$

\ed

\bd

$\hspace{0.01em}$

% (+++ Orig. No.:  Definition Independence +++)

\label{Definition Independence}

Let $ \xdo $ be over $U.$

$ \xdo $ is independent iff $ \xcA \xbd: \xdo \xcp 2. \xcE m \xbe U.m
\xcm \xbd.$

\ed

Obviously, independence does not inherit downward to subsets of $U.$

\bd

$\hspace{0.01em}$

% (+++ Orig. No.:  Definition Delta-O +++)

\label{Definition Delta-O}

This definition is only intended for the set variant.

Let $ \xdo $ be over $U.$

$ \xdd ( \xdo ):=\{X \xcc U':$ $ \xcA \xdo ' \xcc \xdo $ $ \xcA \xbd:
\xdo ' \xcp 2$

$ \xDC $ $(( \xcE m,m' \xbe U,$ $m,m' \xcm \xbd,$ $m \xbe X,m' \xce X)$ $
\xch $ $( \xcE m'' \xbe X.m'' \xcm \xbd \xcu m'' \xeb_{s}m' ))\}$

This property expresses that we can satisfy obligations independently: If
we
respect $O,$ we can, in addition,
respect $O',$ and if we are hopeless kleptomaniacs, we may still not be a
murderer. If $X \xbe \xdd ( \xdo ),$ we can go from $U-X$
into $X$ by improving on all $O \xbe \xdo,$ which we have not fixed by $
\xbd,$ if $ \xbd $ is
not too rigid.
\subsection{
Philosophical discussion of obligations
}
\label{Section Phil}

\ed

We take now a closer look at obligations, in particular at the
ramifications
of the treatment of the relation ``better''. Some aspects of obligations
will also
need a notion of distance, we call them local properties of obligations.
\subsubsection{
A fundamental difference between facts and obligations: asymmetry and
negation
}

There is an important difference between facts and obligations. A
situation
which satisfies an obligation is in some sense ``good'', a situation which
does not, is in some sense ``bad''. This is not true of facts. Being
``round'' is a priori not better than ``having corners'' or vice versa. But
given
the obligation to post the letter, the letter in the mail box is ``good'',
the
letter in he trash bin is ``bad''. Consequently, negation has to play
different role for obligations and for facts.

This is a fundamental property, which can also be
found in orders, planning (we move towards the goal or not), reasoning
with
utility (is $ \xbf $ or $ \xCN \xbf $ more useful?), and probably others,
like perhaps the
Black Raven paradox.

We also think that the Ross paradox (see below) is a true paradox, and
should be
avoided. A closer look shows that this paradox involves arbitrary
weakening,
in particular by the ``negation'' of an obligation. This was a starting
point of our analysis.

``Good'' and ``bad'' cannot mean that any situation satisfying obligation $O$
is better than any situation not satisfying $O,$ as the following example
shows.

\be

$\hspace{0.01em}$

% (+++ Orig. No.:  Example 3-Obligations +++)

\label{Example 3-Obligations}

If we have three independent
and equally strong obligations, $O,$ $O',$ $O'',$ then a situation
satisfying
$O$ but neither $O' $ nor $O'' $ will not be better than one satisfying
$O' $ and $O'',$
but not $O.$

\ee

We have to introduce some kind of ``ceteris paribus''. All other
things being equal, a situation satisfying $O$ is better than a situation
not satisfying $O,$
see Section \ref{Section Ceteris} (page \pageref{Section Ceteris}).

\be

$\hspace{0.01em}$

% (+++ Orig. No.:  Example Ross-Paradox +++)

\label{Example Ross-Paradox}

The original version of the Ross paradox reads: If we have the obligation
to post the letter, then we have the obligation to post or burn the
letter.
Implicit here is the background knowledge that burning the letter implies
not
to post it, and is even worse than not posting it.

We prefer a modified version, which works with two independent
obligations:
We have the obligation to post the letter, and we have the obligation to
water
the plants. We conclude by unrestricted weakening that we have the
obligation to
post the letter or $ \xCf not$ to water the plants. This is obvious
nonsense.

It is not the ``or'' itself which is the problem. For instance, in case of
an
accident, to call an ambulance or to help the victims by giving first aid
is a
perfectly reasonable obligation. It is the negation of the obligation to
water the plants which is the problem. More generally,
it must not be that the system of suitable sets is closed under
arbitrary supersets,
otherwise we have closure under arbitrary right weakening, and thus the
Ross
paradox. Notions like ``big subset'' or ``small exception sets'' from the
semantics of nonmonotonic logics are closed under
supersets, so they are not suitable.
\subsubsection{
``And'' and ``or'' for obligations
}

\ee

``Not'' behaves differently for facts and for obligations. If $O$ and $O' $
are
obligations, can $O \xcu O' $ be considered an obligation? We think, yes.
``Ceteris paribus'', satisfying $O$ and $O' $ together is better than not to
do
so. If is the obligation to post the letter, $O' $ to water the plants,
then
doing both is good, and better than doing none, or only one. Is $O \xco O'
$
an obligation? Again, we think, yes. Satisfying one (or even both, a
non-exclusive or) is better than doing nothing. We might not have enough
time
to do both, so we do our best, and water the plants or post the letter.
Thus, if $ \xba $ and $ \xbb $ are obligations, then so will be $ \xba
\xcu \xbb $ and $ \xba \xco \xbb,$ but not
anything involving $ \xCN \xba $ or $ \xCN \xbb.$ (In a non-trivial
manner, leaving aside
tautologies and contradictions which have to be considered separately.)
To summarize: ``and'' and ``or'' preserve the asymmetry, ``not'' does not,
therefore
we can combine obligations using ``and'' and ``or'', but not ``not''. Thus,
a reasonable notion of derivation of obligations will work with $ \xcu $
and $ \xco,$ but
not with $ \xCN.$

We should not close under inverse $ \xcu,$ i.e. if $ \xbf \xcu \xbf ' $
is an obligation, we
should not conclude that $ \xbf $ and $ \xbf ' $ separately are
obligations, as the following
example shows.

\be

$\hspace{0.01em}$

% (+++ Orig. No.:  Example Not-Inv-And +++)

\label{Example Not-Inv-And}

Let $p$ stand for: post letter, $w:$ water plants, $s:$ strangle
grandmother.

Consider now $ \xbf \xcu \xbf ',$ where $ \xbf =p \xco ( \xCN p \xcu \xCN
w),$ $ \xbf ' =p \xco ( \xCN p \xcu w \xcu s).$
$ \xbf \xcu \xbf ' $ is equivalent to $p$ - though it is perhaps a bizarre
way to express
the obligation to post the letter. $ \xbf $ leaves us the possibility not
to water the
plants, and $ \xbf ' $ to strangle the grandmother, and neither seem good
obligations.
$ \xcz $
\\[3ex]

\ee

\br

$\hspace{0.01em}$

% (+++ Orig. No.:  Remark Not-Inv-And +++)

\label{Remark Not-Inv-And}

This is particularly important in the case of soft obligations, as we see
now,
when we try to apply the rules of preferential reasoning to obligations.

One of the rules of preferential reasoning is the $ \xCf (OR)$ rule:

$ \xbf \xcn \xbq,$ $ \xbf ' \xcn \xbq $ $ \xch $ $ \xbf \xco \xbf ' \xcn
\xbq.$

Suppose we have $ \xbf \xcn \xbq ' \xcu \xbq '',$ and $ \xbf ' \xcn \xbq
'.$ We might be tempted to split
$ \xbq ' \xcu \xbq '' $ - as $ \xbq ' $ is a ``legal'' obligation, and
argue: $ \xbf \xcn \xbq ' \xcu \xbq '',$ so
$ \xbf \xcn \xbq ',$ moreover $ \xbf ' \xcn \xbq ',$ so $ \xbf \xco \xbf
' \xcn \xbq '.$ The following example shows
that this is not always justified.

\er

\be

$\hspace{0.01em}$

% (+++ Orig. No.:  Example Drugs +++)

\label{Example Drugs}

Consider the following obligations for a physician:

Let $ \xbf ' $ imply that the patient has no heart disease, and if $ \xbf
' $ holds,
we should give drug $ \xCf A$ or (not drug $ \xCf A,$ but drug $ \xCf B),$
abbreviated $A \xco ( \xCN A \xcu B).$
$( \xCf B$ is considered dangerous for people with heart problems.)

Let $ \xbf $ imply that the patient has heart problems. Here, the
obligation is
$(A \xco ( \xCN A \xcu B)) \xcu (A \xco ( \xCN A \xcu \xCN B)),$
equivalent to $ \xCf A.$

The false conclusion would then be $ \xbf ' \xcn A \xco ( \xCN A \xcu B),$
and $ \xbf \xcn A \xco ( \xCN A \xcu B),$
so $ \xbf \xco \xbf ' \xcn A \xco ( \xCN A \xcu B),$ so in both situation
we should either give $ \xCf A$ or
$B,$ but $B$ is dangerous in ``one half'' of the situations.

$ \xcz $
\\[3ex]

\ee

We captured this idea about ``and'' and ``or'' in
Definition \ref{Definition ui} (page \pageref{Definition ui}).
\subsubsection{
Ceteris paribus - a local poperty
}
\label{Section Ceteris}

Basically, the set of points ``in'' an obligation has to be better than the
set
of ``exterior'' points. As
above Example \ref{Example 3-Obligations} (page \pageref{Example 3-Obligations})
with three obligations shows, demanding
that any element inside is better than any element outside, is too strong.
We use instead the ``ceteris paribus'' idea.

``All other things being equal'' seems to play a crucial role in
understanding obligations. Before we try to analyse it, we look for other
concepts which have something to do with it.

The Stalnaker/Lewis semantics for counterfactual conditionals also works
with
some kind of ``ceteris paribus''.
``If it were to rain, $I$ would use an umbrella'' means something like:
``If it were to rain, and there were not a very strong wind'' (there is no
such
wind now), ``if $I$ had an umbrella'' (I have one now), etc., i.e. if things
were mostly as they are now, with the exception that now it does not rain,
and in the situation $I$ speak about it rains, then $I$ will use an
umbrella.

But also theory revision in the AGM sense contains - at least as objective
-
this idea: Change things as little as possible to incorporate some new
information in a consistent way.

When looking at the ``ceteris paribus'' in obligations, a natural
interpretation
is to read it as ``all other obligations being unchanged'' (i.e. satisfied
or not as before). This is then just a Hamming distance considering the
obligations (but not other information).

Then, in particular, if $ \xdo $ is a family of obligations, and
if $x$ and $x' $ are in the same subset $ \xdo ' \xcc \xdo $ of
obligations, then an
obligation derived from $ \xdo $ should not separate them. More precisely,
if
$x \xbe O \xbe \xdo \xcj x' \xbe O \xbe \xdo,$ and $D$ is a derived
obligation, then $x \xbe D \xcj x' \xbe D.$

\be

$\hspace{0.01em}$

% (+++ Orig. No.:  Example Spaghetti +++)

\label{Example Spaghetti}

If the only obligation is not to kill, then it should not be
derivable not to kill and to eat spaghetti.

\ee

Often, this is impossible, as obligations are not independent. In this
case,
but also in other situations, we can push ``ceteris paribus'' into an
abstract distance $d$ (as in the Stalnaker/Lewis semantics), which we
postulate
as given, and say that satisfying an obligation makes things better when
going from ``outside'' the obligation to the $d-$closest situation ``inside''.
Conversely, whatever the analysis of ``ceteris paribus'', and given a
quality
order on the situations, we can now define an obligation as a formula
which
(perhaps among other criteria)
``ceteris paribus'' improves the situation when we go from ``outside'' the
formula ``inside''.

A simpler way to capture ``ceteris paribus'' is to connect it directly to
obligations,
see Definition \ref{Definition Delta-O} (page \pageref{Definition Delta-O}).
This is probably too much tied
to independence (see below), and thus too rigid.
\subsubsection{
Hamming neighbourhoods
}

A combination concept is a Hamming neighbourhood:

$X$ is called a Hamming neighbourhood of the best cases iff for any $x
\xbe X$
and $y$ a best case with minimal distance from $x,$ all elements between
$x$ and $y$ are in $X.$

For this, we need a notion of distance (also to define ``between'' ).
This was made precise
in Definition \ref{Definition Hamming-Distance} (page \pageref{Definition
Hamming-Distance})  and
Definition \ref{Definition Neighbourhood} (page \pageref{Definition
Neighbourhood}).
\subsubsection{
Global and mixed global/local properties of obligations
}

We look now at some global properties (or mixtures of global and local)
which
seem desirable for obligations:

 \xEh

 \xDH Downward closure

Consider the following example:

\be

$\hspace{0.01em}$

% (+++ Orig. No.:  Example Not-Global +++)

\label{Example Not-Global}

Let $U':=\{x,x',y,y' \}$ with $x':=pqrs,$ $y':=pqr \xCN s,$ $x:= \xCN
p \xCN qr \xCN s,$ $y:= \xCN p \xCN q \xCN r \xCN s.$

Consider $X:=\{x,x' \}.$

The counting version:

Then $x' $ has quality 4 (the best), $y' $ has quality 3, $x$ has 1, $y$
has 0.

$d_{c}(x',y' )=1,$ $d_{c}(x,y)=1,$ $d_{c}(x,y' )=2.$

\ee

Then above ``ceteris paribus'' criterion is satisfied, as $y' $ and $x$ do
not
``see'' each other, so
$X \xeb_{l,c} \xdC X.$

But $X$ is not downward closed, below $x \xbe X$ is a better element $y'
\xce X.$

This seems an argument against $X$ being an obligation.

The set version:

We still have $x' \xeb_{s}y' \xeb_{s}x \xeb_{s}y.$ As shown in
Fact \ref{Fact Quality-Distance} (page \pageref{Fact Quality-Distance}),
$d_{s}(x,y)$ (and also
$d_{s}(x',y' ))$ and $d_{s}(x,y' )$ are not comparable, so our argument
collapses.

As a matter of fact, we have the result that the ``ceteris paribus''
criterion entails downward closure in the set variant, see
Fact \ref{Fact Local-Global} (page \pageref{Fact Local-Global}).

$ \xcz $
\\[3ex]

Note that a sufficiently rich domain (put elements between $y' $ and $x)$
will make
this local condition (for $ \xeb )$ a global one, so we have here a domain
problem. Domain problems are discussed e.g. in
 \cite{Sch04} and  \cite{GS08a}.

 \xDH Best states

It seems also reasonable to postulate that obligations contain all best
states. In particular, obligations have then to be consistent - under the
condition that best states exist. We are aware that this point can be
debated,
there is, of course, an easy technical way out: we take, when necessary,
unions of obligations to cover the set of ideal cases.
So obligations will be certain ``neighbourhoods'' of the ``best'' situations.

We think, that some such notion of neighbourhood is a good
candidate for a semantics:

 \xEI

 \xDH
A system of neighbourhoods is not necessarily closed under supersets.

 \xDH
Obligations express something like an approximation to the ideal case
where all obligations (if possible, or, as many as possible) are
satisfied,
so we try to be close to the ideal. If we satisfy an obligation, we are
(relatively) close, and stay so as long as the obligation is satisfied.

 \xDH
The notion of neighbourhood expresses the idea of being close, and
containing
everything which is sufficiently close.
Behind ``containing everything which is sufficiently close'' is the idea
of being in some sense convex. Thus, ``convex'' or ``between''
is another basic notion
to be investigated. See here also the discussion of ``between''
in  \cite{Sch04}.

 \xEJ

 \xEj
\subsubsection{
Soft obligations
}

``Soft'' obligations are obligations which have exceptions. Normally,
one is obliged to do $O,$ but there are cases where one is not obliged.
This is like soft rules, as ``Birds fly'' (but penguins do not), where
exceptions are not explicitly mentioned.

The semantic notions of size are very useful here, too. We will content
ourselves that soft obligations satisfy the postulates of usual
obligations
everywhere except on a small set of cases. For instance, a soft obligation
$O$
should be downward closed ``almost'' everywhere, i.e. for a small subset
of pairs $ \xBc a,b \xBe $ in $U \xCK U$ we accept that $a \xeb b,$ $b \xbe O,$
$a
\xce O.$ We transplanted
a suitable and cautious notion of size from the components to the
product in Definition \ref{Definition Size-Product} (page \pageref{Definition
Size-Product}).

When we look at the requirement to contain the best cases, we might have
to
soften this, too. We will admit that a small set of the ideal cases might
be
excluded. Small can be relative to all cases, or only to all ideal cases.

Soft obligations generate an ordering which takes care of exceptions, like
the normality ordering of birds will take care of penguins: within the
set of pengins, non-flying animals are the normal ones.
Based on this ordering, we define ``derived soft obligations'', they may
have (a small set of) exceptions with respect to this ordering.
\subsubsection{
Overview of different types of obligations
}

 \xEh

 \xDH Hard obligations. They hold without exceptions, as in the Ten
Commandments. You should not kill.

 \xEh

 \xDH
In the simplest case, they apply everywhere and can be combined
arbitrarily,
i.e. for any $ \xdo ' \xcc \xdo $ there is a model where all $O \xbe \xdo
' $ hold, and no $O' \xbe \xdo - \xdo '.$

 \xDH
In a more complicated case, not all combinations are possible. This is the
same as considering just an arbitrary subset of $U$ with the same set $
\xdo $ of
obligations. This case is very similar to the case of conditional
obligations
(which might not be defined outside a subset of $U),$ and we treat them
together.

A good example is the Considerate Assassin:

\be

$\hspace{0.01em}$

% (+++ Orig. No.:  Example Considerate-Assassin +++)

\label{Example Considerate-Assassin}

Normally, one should not offer a cigarette to someone, out of respect for
his health. But the considerate assassin might do so nonetheless, on the
cynical reasoning that the victim's health is going to suffer anyway:

(1) One should not kill, $ \xCN k.$

(2) One should not offer cigarettes, $ \xCN o.$

(3) The assassin should offer his victim a cigarette before killing him,
if $k,$ then $o.$

Here, globally, $ \xCN k$ and $ \xCN o$ is best, but among $k-$worlds, $o$
is better than $ \xCN o.$
The model ranking is $ \xCN k \xcu \xCN o \xeb \xCN k \xcu o \xeb k \xcu o
\xeb k \xcu \xCN o.$

\ee

 \xEj

Recall that an obligation for the whole set need not be an obligation for
a
subset any more, as it need not contain all best states. In this case, we
may have to take a union with other obligations.

 \xDH Soft obligations.

Many obligations have exceptions. Consider the following example:

\be

$\hspace{0.01em}$

% (+++ Orig. No.:  Example Library +++)

\label{Example Library}

You are in
a library. Of course, you should not pour water on a book. But if the book
has caught fire, you should pour water on it to prevent worse damage. In
stenographic style these obligations read: ``Do not pour water on books''.
``If a book is on fire, do pour water on it.'' It is like ``birds fly'', but
``penguins do not fly'', ``soft'' or nonmonotonic obligations, which have
exceptions, which are not formulated in the original obligation, but added
as exceptions.

\ee

We could have formulated the library obligation also without exceptions:
``When you are in a library, and the book is not on fire, do not pour water
on
it.''
``When you are in a library, and the book is on fire, pour water on it.''
This formulation avoids exceptions. Conditional obligations behave like
restricted quantifiers: they apply in a subset of all possible cases.

We treat now the considerate assassin case as an obligation
(not to offer) with exceptions.
Consider the full set $U,$ and
consider the obligation $ \xCN o.$ This is not downward
closed, as $k \xcu o$ is better than $k \xcu \xCN o.$
Downward closure will only hold for ``most'' cases, but not for all.

 \xDH

Contrary-to-duty obligations.

Contrary-to-duty obligations are about different degrees of fulfillment.
If you should ideally not have any fence, but are not willing or able to
fulfill this obligation (e.g. you have a dog which might stray), then you
should
at least paint it white to make it less conspicuous.
This is also a conditional obligation.
Conditional, as it specifies what has to be done if there is a fence.
The new aspect in contrary-to-duty obligations is the different degree of
fulfillment.

We will not treat contrary-to-duty obligations here, as they do not seem
to
have any import on our basic ideas and solutions.

 \xDH
A still more complicated case is when the language of obligations is not
uniform, i.e. there are subsets $V \xcc U$ where obligations are defined,
which are not defined in $ \xCf U- \xCf V.$

We will not pursue this case here.

 \xEj
\subsubsection{
Summary of the philosophical remarks
}

 \xEh

 \xDH It seems justifiable to say that an obligation is satisfied or holds
in a certain situation.

 \xDH Obligations are fundamentally asymmetrical, thus negation has to be
treated with care. ``Or'' and ``and'' behave as for facts.

 \xDH Satisfying obligations improves the situation with respect to some
given
grading - ceteris paribus.

 \xDH ``Ceteris paribus'' can be defined by minimal change with respect to
other
obligations, or by an abstract distance.

 \xDH Conversely, given a grading and some distance, we can define an
obligation
locally as describing an improvement with respect to this grading when
going
from ``outside'' to the closest point ``inside'' the obligation.

 \xDH Obligations should also have global properties: they should be
downward (i.e. under increasing quality) closed, and cover the set
of ideal cases.

 \xDH The properties of ``soft'' obligations, i.e. with exceptions, have to
be
modified appropriately. Soft obligations generate an ordering, which in
turn may generate other obligations, where exceptions to the ordering
are permitted.

 \xDH Quality and distance can be defined from an existing set of
obligations
in the set or the counting variant. Their behaviour is quite different.

 \xDH We distinguished various cases of obligations, soft and hard, with
and without all possibilities, etc.

 \xEj

Finally, we should emphasize that the notions of distance, quality, and
size
are in principle independent, even if they may be based on a common
substructure.
\subsection{
Examination of the various cases
}
\label{Section Exam}

We will concentrate here on the set version of hard obligations.
\subsubsection{
Hard obligations for the set approach
}
\label{Section Hard-Obligations}
\paragraph{
Introduction
}

We work here in the set version, the $ \xbe -$case, and examine mostly
the set version only.

We will assume a set $ \xdo $ of obligations to be given.
We define the relation $ \xeb:= \xeb_{ \xdo }$ as described in
Definition \ref{Definition Hamming-Distance} (page \pageref{Definition
Hamming-Distance}), and the distance $d$ is the
Hamming distance
based on $ \xdo.$
\paragraph{
The not necessarily independent case
}

\be

$\hspace{0.01em}$

% (+++ Orig. No.:  Example Dependent-1 +++)

\label{Example Dependent-1}

Work in the set variant. We show that $X$ $ \xec_{s}-closed$ does not
necessarily imply that $X$ contains all $ \xec_{s}-best$ elements.

Let $ \xdo:=\{p,q\},$ $U':=\{p \xCN q, \xCN pq\},$ then all elements of
$U' $ have best quality
in $U',$ $X:=\{p \xCN q\}$ is closed, but does not contain all best
elements. $ \xcz $
\\[3ex]

\ee

\be

$\hspace{0.01em}$

% (+++ Orig. No.:  Example Dependent-3 +++)

\label{Example Dependent-3}

Work in the set variant. We show that $X$ $ \xec_{s}-closed$ does not
necessarily imply that $X$ is a neighbourhood of the best elements, even
if $X$ contains them.

Consider $x:=pq \xCN rstu,$ $x':= \xCN pqrs \xCN t \xCN u,$ $x'':=p \xCN
qr \xCN s \xCN t \xCN u,$
$y:=p \xCN q \xCN r \xCN s \xCN t \xCN u,$ $z:=pq \xCN r \xCN s \xCN t
\xCN u.$ $U:=\{x,x',x'',y,z\},$ the $ \xeb_{s}-$best elements
are $x,x',x'',$ they are contained in $X:=\{x,x',x'',z\}.$
$d_{s}(z,x)=\{s,t,u\},$
$d_{s}(z,x' )=\{p,r,s\},$ $d_{s}(z,x'' )=\{q,r\},$ so $x'' $ is one of the
best
elements closest to $z.$ $d(z,y)=\{q\},$ $d(y,x'' )=\{r\},$ so $[z,x''
]=\{z,y,x'' \},$ $y \xce X,$ but $X$
is downward closed. $ \xcz $
\\[3ex]

\ee

\bfa

$\hspace{0.01em}$

% (+++ Orig. No.:  Fact Global-Dependent-S +++)

\label{Fact Global-Dependent-S}

Work in the set variant.

Let $X \xEd \xCQ,$ $X$ $ \xec_{s}-closed.$ Then

(1) $X$ does not necessarily contain all best elements.

Assume now that $X$ contains, in addition, all best elements. Then

(2) $X \xeb_{l,s} \xdC X$ does not necessarily hold.

(3) $X$ is (ui).

(4) $X \xbe \xdd ( \xdo )$ does not necessarily hold.

(5) $X$ is not necessarily a neighbourhood of the best elements.

(6) $X$ is an improving neighbourhood of the best elements.

\efa

\subparagraph{
Proof
}

$\hspace{0.01em}$

% (+++ Orig.:  Proof +++)

(1) See Example \ref{Example Dependent-1} (page \pageref{Example Dependent-1})

(2) See Example \ref{Example Dependent-2} (page \pageref{Example Dependent-2})

(3)
If there is $m \xbe X,$ $m \xce O$ for all $O \xbe \xdo,$ then by closure
$X=U,$ take $ \xdo_{i}:= \xCQ.$

For $m \xbe X$ let $ \xdo_{m}:=\{O \xbe \xdo:m \xbe O\}.$ Let $X':= \xcV
\{ \xcS \xdo_{m}:m \xbe X\}.$

$X \xcc X':$ trivial, as $m \xbe X \xcp m \xbe \xcS \xdo_{m} \xcc X'.$

$X' \xcc X:$ Let $m' \xbe \xcS \xdo_{m}$ for some $m \xbe X.$ It suffices
to show that $m' \xec_{s}m.$
$m' \xbe \xcS \xdo_{m}= \xcS \{O \xbe \xdo:m \xbe O\},$ so for all $O
\xbe \xdo $ $(m \xbe O \xcp m' \xbe O).$

(4) Consider Example \ref{Example Dependent-2} (page \pageref{Example
Dependent-2}),
let $dom( \xbd )=\{r,s\},$ $ \xbd (r)= \xbd (s)=0.$ Then $x,y \xcm \xbd,$
but $x' \xcM \xbd $ and $x \xbe X,$ $y \xce X,$
but there is no $z \xbe X,$ $z \xcm \xbd $ and $z \xeb y,$ so $X \xce \xdd
( \xdo ).$

(5) See Example \ref{Example Dependent-3} (page \pageref{Example Dependent-3}).

(6) By Fact \ref{Fact Quality-Distance} (page \pageref{Fact Quality-Distance}),
(5).

$ \xcz $
\\[3ex]

\bfa

$\hspace{0.01em}$

% (+++ Orig. No.:  Fact Global-Dependent-Converse-S +++)

\label{Fact Global-Dependent-Converse-S}

Work in the set variant

(1.1) $X \xeb_{l,s} \xdC X$ implies that $X$ is $ \xec_{s}-closed.$

(1.2) $X \xeb_{l,s} \xdC X$ $ \xch $ $X$ contains all best elements

(2.1) $X$ is (ui) $ \xch $ $X$ is $ \xec_{s}-closed.$

(2.2) $X$ is (ui) does not necessarily imply that $X$ contains all $
\xec_{s}-best$
elements.

(3.1) $X \xbe \xdd ( \xdo )$ $ \xch $ $X$ is $ \xec_{s}-closed$

(3.2) $X \xbe \xdd ( \xdo )$ implies that $X$ contains all $
\xec_{s}-best$ elements.

(4.1) $X$ is an improving neighbourhood of the $ \xec_{s}-best$ elements $
\xch $ $X$ is
$ \xec_{s}-closed.$

(4.2) $X$ is an improving neighbourhood of the best elements $ \xch $ $X$
contains all best elements.

\efa

\subparagraph{
Proof
}

$\hspace{0.01em}$

% (+++ Orig.:  Proof +++)

(1.1) By Fact \ref{Fact Local-Global} (page \pageref{Fact Local-Global}).

(1.2) By Fact \ref{Fact General-Obligation} (page \pageref{Fact
General-Obligation}).

(2.1)
Let $O \xbe \xdo,$ then $O$ is downward closed (no $y \xce O$ can be
better than $x \xbe O).$
The rest follows from Fact \ref{Fact Subset-Closure} (page \pageref{Fact
Subset-Closure})  (3).

(2.2) Consider Example \ref{Example Dependent-1} (page \pageref{Example
Dependent-1}), $p$ is (ui) (formed in
$U!),$ but $p \xcs X$ does
not contain $ \xCN pq.$

(3.1)
Let $X \xbe \xdd ( \xdo ),$ but let $X$ not be closed.
Thus, there are $m \xbe X,$ $m' \xec_{s}m,$ $m' \xce X.$

Case 1: Suppose $m' \xCq m.$ Let $ \xbd_{m}: \xdo \xcp 2,$ $
\xbd_{m}(O)=1$ iff $m \xbe O.$ Then $m,m' \xcm \xbd_{m},$
and there cannot be any $m'' \xcm \xbd_{m},$ $m'' \xeb_{s}m',$ so $X \xce
\xdd ( \xdo ).$

Case 2: $m' \xeb_{s}m.$ Let $ \xdo ':=\{O \xbe \xdo:m \xbe O \xcj m'
\xbe O\},$ $dom( \xbd )= \xdo ',$ $ \xbd (O):=1$ iff $m \xbe O$ for
$O \xbe \xdo '.$ Then $m,m' \xcm \xbd.$ If there is $O \xbe \xdo $ s.t.
$m' \xce O,$ then by $m' \xec_{s}m$ $m \xce O,$
so $O \xbe \xdo '.$ Thus for all $O \xce dom( \xbd ).m' \xbe O.$ But then
there is no $m'' \xcm \xbd,$ $m'' \xeb_{s}m',$
as $m' $ is already optimal among the $n$ with $n \xcm \xbd.$

(3.2) Suppose $X \xbe \xdd ( \xdo ),$ $x' \xbe U-X$ is a best element,
take $ \xbd:= \xCQ,$ $x \xbe X.$
Then there must be $x'' \xeb x',$ $x'' \xbe X,$ but this is impossible as
$x' $ was best.

(4.1) By Fact \ref{Fact Quality-Distance} (page \pageref{Fact Quality-Distance})
, (4) all minimal elements
have incomparabel
distance. But if $z \xec y,$ $y \xbe X,$ then either $z$ is
minimal or it is above a minimal element, with minimal distance from $y,$
so $z \xbe X$
by Fact \ref{Fact Quality-Distance} (page \pageref{Fact Quality-Distance})  (3).

(4.2) Trivial.

$ \xcz $
\\[3ex]
\paragraph{
The independent case
}

Assume now the system to be independent, i.e. all combinations of $ \xdo $
are present.

Note that there is now only one minimal element, and the notions of
Hamming neighbourhood of the best elements
and improving Hamming neighbourhood of the best elements coincide.

\bfa

$\hspace{0.01em}$

% (+++ Orig. No.:  Fact Global-Independent-S +++)

\label{Fact Global-Independent-S}

Work in the set variant.

Let $X \xEd \xCQ,$ $X$ $ \xec_{s}-closed.$ Then

(1) $X$ contains the best element.

(2) $X \xeb_{l,s} \xdC X$

(3) $X$ is (ui).

(4) $X \xbe \xdd ( \xdo )$

(5) $X$ is a (improving) Hamming neighbourhood of the best elements.

\efa

\subparagraph{
Proof
}

$\hspace{0.01em}$

% (+++ Orig.:  Proof +++)

(1) Trivial.

(2) Fix $x \xbe X,$ let $y$ be closest to $x,$ $y \xce X.$ Suppose $x \xeB
y,$ then there must be
$O \xbe \xdo $ s.t. $y \xbe O,$ $x \xce O.$ Choose $y' $ s.t. $y' $ is
like $y,$ only $y' \xce O.$ If $y' \xbe X,$ then
by closure $y \xbe X,$ so $y' \xce X.$ But $y' $ is closer to $x$ than $y$
is, $contradiction.$

Fix $y \xbe U-$X. Let $x$ be closest to $y,$ $x \xbe X.$ Suppose $x \xeB
y,$ then there is $O \xbe \xdo $
s.t. $y \xbe O,$ $x \xce O.$ Choose $x' $ s.t. $x' $ is like $x,$ only $x'
\xbe O.$ By closure of $X,$
$x' \xbe X,$ but $x' $ is closer to $y$ than $x$ is, $contradiction.$

(3) By Fact \ref{Fact Global-Dependent-S} (page \pageref{Fact
Global-Dependent-S})  (3)

(4)
Let $X$ be closed, and $ \xdo ' \xcc \xdo,$ $ \xbd: \xdo ' \xcp 2,$
$m,m' \xcm \xbd,$ $m \xbe X,$ $m' \xce X.$
Let $m'' $ be s.t. $m'' \xcm \xbd,$ and for all $O \xbe \xdo -dom( \xbd
)$ $m'' \xbe O.$ This exists by
independence.
Then $m'' \xec_{s}m',$ but also $m'' \xec_{s}m,$ so $m'' \xbe X.$ Suppose
$m'' \xCq m',$ then
$m' \xec_{s}m'',$ so $m' \xbe X,$ contradiction, so $m'' \xeb_{s}m'.$

(5) Trivial by (1), the remark preceding this Fact, and
Fact \ref{Fact Global-Dependent-S} (page \pageref{Fact Global-Dependent-S}) 
(6).

\bfa

$\hspace{0.01em}$

% (+++ Orig. No.:  Fact Global-Independent-Converse-S +++)

\label{Fact Global-Independent-Converse-S}

Work in the set variant.

(1) $X \xeb_{l,s} \xdC x$ $ \xch $ $X$ is $ \xec_{s}-closed,$

(2) $X$ is (ui) $ \xch $ $X$ is $ \xec_{s}-closed,$

(3) $X \xbe \xdd ( \xdo )$ $ \xch $ $X$ is $ \xec_{s}-closed,$

(4) $X$ is a (improving) neighbourhood of the best elements $ \xch $ $X$
is $ \xec_{s}-closed.$

\efa

\subparagraph{
Proof
}

$\hspace{0.01em}$

% (+++ Orig.:  Proof +++)

(1) Suppose there are $x \xbe X,$ $y \xbe U-$X, $y \xeb x.$ Choose them
with minimal distance.
If $card(d_{s}(x,y))>1,$ then there is
$z,$ $y \xeb_{s}z \xeb_{s}x,$ $z \xbe X$ or $z \xbe U-$X, contradicting
minimality. So $card(d_{s}(x,y))=1.$
So $y$ is among the closest elements of $U-X$ seen from $x,$ but then by
prerequisite
$x \xeb y,$ $contradiction.$

(2) By Fact \ref{Fact Global-Dependent-Converse-S} (page \pageref{Fact
Global-Dependent-Converse-S})  (2.1).

(3) By Fact \ref{Fact Global-Dependent-Converse-S} (page \pageref{Fact
Global-Dependent-Converse-S})  (3.1).

(4) There is just one best element $z,$ so if $x \xbe X,$ then $[x,z]$
contains all
$y$ $y \xeb x$ by Fact \ref{Fact Quality-Distance} (page \pageref{Fact
Quality-Distance})  (3).

$ \xcz $
\\[3ex]

The $ \xdd ( \xdo )$ condition seems to be adequate only for the
independent situation,
so we stop considering it now.

\bfa

$\hspace{0.01em}$

% (+++ Orig. No.:  Fact Int-Union +++)

\label{Fact Int-Union}

Let $X_{i} \xcc U,$ $i \xbe I$ a family of sets, we note the following
about closure under unions and intersections:

(1) If the $X_{i}$ are downward closed, then so are their unions and
intersections.

(2) If the $X_{i}$ are (ui), then so are their unions and intersections.

\efa

\subparagraph{
Proof
}

$\hspace{0.01em}$

% (+++ Orig.:  Proof +++)

Trivial. $ \xcz $
\\[3ex]

We do not know whether $ \xeb_{l,s}$ is preserved under unions and
intersections,
it does not seem an easy problem.

\bfa

$\hspace{0.01em}$

% (+++ Orig. No.:  Fact Relativization +++)

\label{Fact Relativization}

(1) Being downward closed is preserved while going to subsets.

(2) Containing the best elements is not preserved (and thus neither the
neighbourhood property).

(3) The $ \xdd ( \xdo )$ property is not preserved.

(4) $ \xec_{l,s}$ is not preserved.

\efa

\subparagraph{
Proof
}

$\hspace{0.01em}$

% (+++ Orig.:  Proof +++)

(4) Consider Example \ref{Example Not-Global} (page \pageref{Example
Not-Global}), and
eliminate $y$ from $U',$ then the closest to $x$ not in $X$ is $y',$
which is better.

$ \xcz $
\\[3ex]
\subsubsection{
Remarks on the counting case
}

\br

$\hspace{0.01em}$

% (+++ Orig. No.:  Remark +++)

\label{Remark}

In the counting variant all qualities are comparabel. So if $X$
is closed, it will contain all minimal elements.

\er

\be

$\hspace{0.01em}$

% (+++ Orig. No.:  Example H-N-Local +++)

\label{Example H-N-Local}

We measure distance by counting.

Consider $a:= \xCN p \xCN q \xCN r \xCN s,$ $b:= \xCN p \xCN q \xCN rs,$
$c:= \xCN p \xCN qr \xCN s,$ $d:=pqr \xCN s,$
let $U:=\{a,b,c,d\},$ $X:=\{a,c,d\}.$ $d$ is the best element,
$[a,d]=\{a,d,c\},$
so $X$ is an improving Hamming neighbourhood, but $b \xeb a,$ so $X
\xeB_{l,c} \xdC X.$

$ \xcz $
\\[3ex]

\ee

\bfa

$\hspace{0.01em}$

% (+++ Orig. No.:  Fact Local-H-N +++)

\label{Fact Local-H-N}

We measure distances by counting.

$X \xeb_{l,c} \xdC X$ does not necessarily imply that $X$ is an improving
Hamming
neighbourhood of the best elements.

\efa

\subparagraph{
Proof
}

$\hspace{0.01em}$

% (+++ Orig.:  Proof +++)

Consider Example \ref{Example Not-Global} (page \pageref{Example Not-Global}).
There $X \xeb_{l,c} \xdC X.$
$x' $ is the best element, and
$y' \xbe [x',x],$ but $y' \xce X.$ $ \xcz $
\\[3ex]
\subsection{
What is an obligation?
}
\label{Section What}

The reader will probably not expect a final definition. All we can do is
to give a tentative definition, which, in all probability, will not be
satisfactory in all cases.

\bd

$\hspace{0.01em}$

% (+++ Orig. No.:  Definition Obligation +++)

\label{Definition Obligation}

We decide for the set relation and distance.

(1) Hard obligation

A hard obligation has the following properties:

(1.1) It contains all ideal cases in the set considered.

\ed

(1.2) It is closed under increasing quality,
Definition \ref{Definition Closed} (page \pageref{Definition Closed})

(1.3) It is an improving neighbourhood of the ideal cases (this also
implies (1.1)),
Definition \ref{Definition Neighbourhood} (page \pageref{Definition
Neighbourhood})

We are less committed to:

(1.4) It is ceteris paribus improving,
Definition \ref{Definition Quality-Extension} (page \pageref{Definition
Quality-Extension})

An obligation $O$ is a derived obligation of a system $ \xdo $ of
obligations iff
it is a hard obligation based on the set variant of the order and distance
generated by $ \xdo.$

(2) Soft obligations

A set is a soft obligation iff it satisfies the soft versions of
above postulates. The notion of size has to be given, and is transferred
to products as described in Definition \ref{Definition Size-Product} (page
\pageref{Definition Size-Product}).
More precisely, strict universal quantifiers are transformed into their
soft variant ``almost all'', and the other operators are left as they are.
Of course, one might also want to use a mixture of soft and hard
conditions, e.g. we might want to have all ideal cases, but renounce
on closure for a small set of pairs
$ \xBc x,x'  \xBe.$

An obligation $O$ is derived from $ \xdo $ iff it is a soft obligation
based
on the set variant of the order and distance generated by the
translation of $ \xdo $ into their hard versions. (I.e. exceptions will
be made explicit.)

\bfa

$\hspace{0.01em}$

% (+++ Orig. No.:  Fact Garbage-In +++)

\label{Fact Garbage-In}

Let $O \xbe \xdo,$ then $ \xdo \xcn O$ in the independent set case.

\efa

\subparagraph{
Proof
}

$\hspace{0.01em}$

% (+++ Orig.:  Proof +++)

We check $(1.1)-(1.3)$ of Definition \ref{Definition Obligation} (page
\pageref{Definition Obligation}).

(1.1) holds by independence.

(1.2) If $x \xbe O,$ $x' \xce O,$ then $x' \xeC_{s}x.$

(1.3) By Fact \ref{Fact Global-Dependent-S} (page \pageref{Fact
Global-Dependent-S})  (6).

Note that (1.4) will also hold by
Fact \ref{Fact Global-Independent-S} (page \pageref{Fact Global-Independent-S}) 
(2).

$ \xcz $
\\[3ex]

\bco

$\hspace{0.01em}$

% (+++ Orig. No.:  Corollary Classical-Consequence +++)

\label{Corollary Classical-Consequence}

Every derived obligation is a classical consequence of the original set
of obligations in the independent set case.

\eco

\subparagraph{
Proof
}

$\hspace{0.01em}$

% (+++ Orig.:  Proof +++)

This follows
from Fact \ref{Fact Global-Independent-S} (page \pageref{Fact
Global-Independent-S})  (3)
and Fact \ref{Fact Garbage-In} (page \pageref{Fact Garbage-In}).

\be

$\hspace{0.01em}$

% (+++ Orig. No.:  Example Burnt-Letter +++)

\label{Example Burnt-Letter}

The Ross paradox is not a derived obligation.

\ee

\subparagraph{
Proof
}

$\hspace{0.01em}$

% (+++ Orig.:  Proof +++)

Suppose we have the alphabet $p,q$ and the obligations $\{p,q\},$ let
$R:=p \xco \xCN q.$ This is not not closed, as $ \xCN p \xcu q \xeb \xCN p
\xcu \xCN q \xbe R.$
$ \xcz $
\\[3ex]
\subsection{
Conclusion
}

Obligations differ from facts in the behaviour of negation, but not of
conjunction and disjunction. The Ross paradox originates, in our opinion,
from the differences in negation. Central to the treatment of obligations
seems to be a relation of ``better'', which can generate obligations, but
also be generated by obligations. The connection between obligations and
this relation of ``better'' seems to be somewhat complicated and leads to
a number of ramifications. A tentative definition of a derivation of
obligations is given.
\section{A comment on work by Aqvist}
\label{Section 2.4.2}
\label{Section Aqvist}
\subsection{
Introduction
}

The article  \cite{Aqv00} discusses three systems, which are
presented now in
outlines. (When necessary, we will give details later in this section.)

(1) The systems $ \xCf Hm,$ where $m \xbe \xbo.$
The (Kripke style) semantics has for each $i \xbe \xbo $ a subset
$opt_{i}$ of the model
set, s.t. the $opt_{i}$ form a disjoint cover of the model set, all
$opt_{i}$ for all
$i \xck m$ are not-empty, and all other $opt_{i}$ are empty. The $opt_{i}$
interpret new
frame constants $Q_{i}.$ The $opt_{i}$ describe intuitively levels of
perfection,
where $opt_{1}$ is the best, and $opt_{m}$ the worst level.

(2) The dyadic deontic logics $ \xCf Gm,$ $m \xbe \xbo.$
The semantics is again given by a cover $opt_{i}$ as above, and in
addition, a
function best, which assigns to each formula $ \xbf $ the ``best'' models of
$ \xbf,$
i.e. those which are in the best $opt_{i}$ set.
The language has the $Q_{i}$ operators, and a new binary operator $O( \xbf
/ \xbq )$
(and its dual $P(./.)),$ which expresses that in the best $ \xbq -$models
$ \xbf $ holds.
Note that there is no explicit ``best'' operator in the language.

(3) The dyadic deontic logic $G.$
The semantics does not contain the $opt_{i}$ any more, but still the
``best''
operator as in case (2), which now corresponds to a ranking of the models
(sufficient axioms are given). The language does not have the $Q_{i}$ any
more,
but contains the $O$ (and $P)$ operator, which is interpreted in the
natural way:
$O( \xbf / \xbq )$ holds iff the best $ \xbq -$models are $ \xbf -$models.
Note that again there is no explicit ``best'' operator in the language.

Thus, it corresponds to putting the $ \xcn -$relation of ranked models in
the object
language.

In particular, there is no finiteness restriction any more, as in cases
(1) and
(2) - and here lies the main difference.

Aqvist gives a Theorem (Theorem 6) which shows (among other things) that

If $G$ is a $ \xCf G-$sentence, provable in $G,$ then it is also provable
in $ \xCf Gm.$

The converse is left open, and Aqvist thought to have found proof using
another
result of his, but there was a loophole, as he had found out.

We close this hole here.
\subsection{
There are (at least) two solutions
}

(1) We take a detour via a language which contains an operator $ \xbb $ to
be
interpreted by ``best''.

(2) We work with the original language.

As a matter of fact, Aqvist's paper contains already an almost complete
solution of type (1), just the translation part is lacking. We give a
slightly
different proof (which will be self contained), which works with the
original
(i.e. possibly infinite) model set, but reduces the quantity of levels to
a
finite number. Thus, the basic idea is the same, our technique is more
specific
to the problem at hand, thus less versatile, but also less ``brutal''.

Yet, we can also work with the original language, even if the operator $O(
\xbf / \xbq )$
will not allow to describe ``best'' as the operator $ \xbb $ can, we can
approximate
``best'' sufficiently well to suit our purposes. (Note also that $O$ allows,
using all formulas, to approximate ``best'' from above: $best( \xbf )$ $=$
$ \xcS \{M( \xbq ):O( \xbq / \xbf )\}.)$

We may describe the difference between solution (1) and (2) as follows:
Solution (1) will preserve the exact ``best'' value of some - but not
necessarily
all, and this cannot be really improved - sets, solution (2) will not even
allow
this, but still, we stay sufficiently close to the original ``best'' value,
so
the formula at hand, and its subformulas, will preserve their truth
values.

In both cases, the basic idea will be the same:
We have a $G-$model for $ \xbf,$ and now construct a $ \xCf Gm-$model
(i.e. with finitely many
levels) for $ \xbf.$ $ \xbf $ is a finite entity, containing only
finitely many
propositional variables, say $p_{1}, \Xl,p_{n}.$ Then we look at all set
of the type
$ \xCL p_{1} \xcu  \Xl  \xcu \xCL p_{n},$ where $ \xCL $ is nothing or the
negation symbol. This is basically
how fine our structure will have to be. If we work with $ \xbb,$ we have
to get
better, as $ \xbb (p_{i})$ will, in general, be something new, so will $
\xbb (p_{i} \xcu \xCN \xbb (p_{i})),$
etc., thus we have to take the nesting of $ \xbb ' $s in $ \xbf $ into
account, too. (This
will be done below.) If we work directly with the $O( \xba / \xbb )$
operator, then
we need not go so far down as $O( \xba / \xbb )$ will always evaluate to
(universal)
true or false. If e.g. $ \xbf $ contains $O( \xba / \xbb )$ and $O( \xba '
/ \xbb ),$ then we will try
to define the ``best'' elements of $M( \xbb )$ as $M( \xbb \xcu \xba \xcu
\xba ' ),$ etc. We have to
be careful to make a ranking, i.e. $O(.,true)$ will give the lowest layer,
and if $ \xbq \xcl \xbq ',$ $O( \xbf / \xbq ),$ $O( \xCN \xbf / \xbq ),$
then the rank of $ \xbq ' $ will be strictly
smaller than the rank of $ \xbq,$ etc. This is all straightforward, but a
bit
tedious.

As said, we will take the first approach, which seems a bit ``cleaner''.

We repeat now the definitions of  \cite{Aqv00} only as far as
necessary to
understand the subsequent pages. In particular, we will not introduce the
axiom systems, as we will work on the semantic side only.

All systems are propositional.

\bd

$\hspace{0.01em}$

% (+++ Orig. No.:  Definition D-6.3.A.1 +++)

\label{Definition D-6.3.A.1}

The systems $ \xCf Hm,$ $m \xbe \xbo.$

The language:

A set $ \xCf Prop$ of propositional variables, True, the usual Boolean
operators,
$N$ (universal necessity), a set $\{Q_{i}:i \xbe \xbo \}$ of systematic
frame constants (zero
place connectives like True) (and the duals False, $M).$

The semantics:

$M= \xBc W,V,\{opt_{i}:i \xbe \xbo \},m \xBe,$ where $W$ is a set of possible
worlds, $V$ a valuation as
usual, each $opt_{i}$ is a subset of $W.$

Validity:

Let $x \xbe W.$

$M,x \xcm p$ iff $x \xbe V(p)$ for $x \xbe W,$ $p \xbe Prop,$

$M,x \xcm True,$

$M,x \xcm N \xbf $ iff for all $y \xbe W$ $M,y \xcm \xbf,$

$M,x \xcm Q_{i}$ iff $x \xbe opt_{i}.$

Conditions on the $opt_{i}:$

(a) $opt_{i} \xcs opt_{j}= \xCQ $ if $i \xEd j,$

(b) $opt_{1} \xcv  \Xl  \xcv opt_{m}=W,$

(c) $opt_{i} \xEd \xCQ $ for $i \xck m,$

(d) $opt_{i}= \xCQ $ for $i>m.$

\ed

\bd

$\hspace{0.01em}$

% (+++ Orig. No.:  Definition D-6.3.A.2 +++)

\label{Definition D-6.3.A.2}

The systems $ \xCf Gm,$ $m \xbe \xbo.$

The language:

It is just like that for the systems $ \xCf Hm,$ with, in addition, a new
binary
connective $O( \xbf / \xbq ),$ with the meaning that the ``best'' $ \xbq
-$worlds satisfy $ \xbf $
(and its dual).

The semantics:

It is also just like the one for $ \xCf Hm,$ with, in addition, a function
$B$ (for
``best'' ), assigning to each formula $ \xbf $ a set of worlds, the best $
\xbf -$worlds, thus
$M= \xBc W,V,\{opt_{i}:i \xbe \xbo \},B,m \xBe.$

Validity:

$M,x \xcm O( \xbf, \xbq )$ iff $B( \xbq ) \xcc M( \xbf ).$

Conditions:

We have to connect $B$ to the $opt_{i},$ this is done by the condition $(
\xbg 0)$ in the
obvious way:

$( \xbg 0)$ $x \xbe B( \xbf )$ iff $M,x \xcm \xbf $ and for each $y \xbe
W,$ if $M,y \xcm \xbf,$ then $x$ is at least
as good as $y$ - i.e. there is no strictly lower $opt_{i}-level$ where $
\xbf $ holds.

\ed

\bd

$\hspace{0.01em}$

% (+++ Orig. No.:  Definition D-6.3.A.3 +++)

\label{Definition D-6.3.A.3}

The system $G.$

The language:

It is just like that for the systems $ \xCf Gm,$ but without the $Q_{i}.$

The semantics:

It is now just $M= \xBc W,V,B \xBe.$

Validity:

As for $ \xCf Gm.$

Conditions:

We have no $opt_{i}$ now, which are replaced by suitable conditions on
$B,$ which
make $B$ the choice function of a ranked structure:

$( \xbs_{0})$ $M( \xbf )=M( \xbf ' )$ $ \xcp $ $B( \xbf )=B( \xbf ' ),$

$( \xbs_{1})$ $B( \xbf ) \xcc M( \xbf ),$

$( \xbs_{2})$ $B( \xbf ) \xcs M( \xbq ) \xcc B( \xbf \xcu \xbq ),$

$( \xbs_{3})$ $M( \xbf ) \xEd \xCQ $ $ \xcp $ $B( \xbf ) \xEd \xCQ,$

$( \xbs_{4})$ $B( \xbf ) \xcs M( \xbq ) \xEd \xCQ $ $ \xcp $ $B( \xbf \xcu
\xbq ) \xcc B( \xbf ) \xcs M( \xbq ).$
\subsection{
Outline
}

\ed

We first show that the language (and logic) $G$ may necessitate infinitely
many
levels (whereas the languages $ \xCf Gm$ $(m \xbe \xbo )$ admit only
finitely many ones). Thus,
when trying to construct a $ \xCf Gm-$model for a $G-$formula $ \xbf,$ we
have to
construct from a possibly ``infinite'' function $B$ via finitely many
$opt_{i}$
levels a new function $B' $ with only finitely many levels, which is
sufficient
for the formula $ \xbf $ under consideration. Crucial for the argument is
that $ \xbf $ is
a finite formula, and we need only a limited degree of discernation for
any
finite formula. Most of the argument is standard reasoning about ranked
structures, as it is common for sufficiently strong preferential
structures.

We will first reformulate the problem slightly using directly an operator
$ \xbb,$
which will be interpreted by the semantic function $B,$ and which results
in a
slightly richer language, as the operators $O$ and $P$ can be expressed
using
$B,$ but not necessarily conversely. Thus, we show slightly more than what
is sufficient to solve Aqvist's problem.

We make this official:

\bd

$\hspace{0.01em}$

% (+++ Orig. No.:  Definition D-6.3.A.4 +++)

\label{Definition D-6.3.A.4}

The language $G' $ is like the language $G,$ without $O(./.)$ and
$P(./.),$ but with
a new unary operator $ \xbb,$ which is interpreted by the semantic
function $B.$

\ed

\br

$\hspace{0.01em}$

% (+++ Orig. No.:  Remark D-6.3.A.1 +++)

\label{Remark D-6.3.A.1}

(1) $O( \xbf / \xbq )$ can thus be expressed by $N( \xbb ( \xbq ) \xcp
\xbf )$ - $N$ universal necessity. In
particular, and this is the important aspect of $ \xbb,$ $ \xbb ( \xbf )$
is now (usually) a
non-trivial set of models, whereas $O$ and $P$ result always in $ \xCQ $
or the set of
all models.

(2) Note that by $( \xbs 3),$ if $M( \xbf ) \xEd \xCQ,$ then $M( \xbb (
\xbf )) \xEd \xCQ,$ but it may very well
be that $M( \xbb ( \xbf ))=M( \xbf )$ - all $ \xbf -$models may be equally
good. $(M( \xbf )$ is the set
of $ \xbf -$models.)

$ \xcz $
\\[3ex]

\er

We use now an infinite theory to force $B$ to have infinitely many levels.

\be

$\hspace{0.01em}$

% (+++ Orig. No.:  Example D-6.3.A.1 +++)

\label{Example D-6.3.A.1}

We construct $ \xbo $ many levels for the models of the language
$\{p_{i}:i \xbe \xbo \},$ going
downward:

Top level: $M(p_{0})$

second level: $M( \xCN p_{0} \xcu p_{1})$

third level: $M( \xCN p_{0} \xcu \xCN p_{1} \xcu p_{2}),$ etc.

We can express via $ \xbb,$ and even via $O(./.),$ that these levels are
all
distinct, e.g. by $ \xbb (p_{0} \xco ( \xCN p_{0} \xcu p_{1})) \xcm \xCN
p_{0},$ or by $O( \xCN p_{0}/p_{0} \xco ( \xCN p_{0} \xcu p_{1})),$ etc.
So we will necessarily have $ \xbo $ many non-empty $opt_{i}-levels,$ for
any $B$ which
satisfies the condition connecting the $opt_{i}$ and $B,$ condition $(
\xbg_{0}).$

$ \xcz $
\\[3ex]

\ee

We work now in the fixed $G-$model $ \xbG $ (with the $ \xbb $ operator.)
Let in the sequel
all $X,$ $Y,$ $X_{i}$ be non-empty (to avoid trivial cases) model sets. By
prerequisite,
$B$ satisfies $( \xbs 0)-( \xbs 4).$

\bd

$\hspace{0.01em}$

% (+++ Orig. No.:  Definition D-6.3.A.5 +++)

\label{Definition D-6.3.A.5}

(1) $X \xeb Y$ iff $B(X \xcv Y) \xcs Y= \xCQ,$

(2) $X \xCq Y$ iff $B(X \xcv Y) \xcs X \xEd \xCQ $ and $B(X \xcv Y) \xcs Y
\xEd \xCQ,$ i.e. iff $X \xeB Y$ and $Y \xeB X.$

\ed

\br

$\hspace{0.01em}$

% (+++ Orig. No.:  Remark D-6.3.A.2 +++)

\label{Remark D-6.3.A.2}

(1) $ \xeb $ and $ \xCq $ behave nicely:

(1.1) $ \xeb $ and $ \xCq $ are transitive, $ \xCq $ is reflexive and
symmetric, and for no $X$
$X \xeb X.$

(1.2) $ \xeb $ and $ \xCq $ cooperate: $X \xeb Y \xCq Z$ $ \xcp $ $X \xeb
Z$ and $X \xee Y \xCq Z$ $ \xcp $ $X \xee Z$

Thus, $ \xeb $ and $ \xCq $ give a ranking.

(2.1) $B(X \xcv Y)=B(X)$ if $X \xeb Y$

(2.2) $B(X \xcv Y)=B(X) \xcv B(Y)$ if $X \xCq Y$

(3) $B(X_{1} \xcv  \Xl  \xcv X_{n})$ $=$ $ \xcV \{B(X_{i}):$ $ \xCN \xcE
X_{j}(X_{j} \xeb X_{i})\}$ $1 \xck i,j \xck n,$ - i.e.
$B(X_{1} \xcv  \Xl  \xcv X_{n})$ is the union of the $B(X_{i})$ with
minimal $ \xeb -$rank.

$ \xcz $
\\[3ex]

\er

We come to the main idea and its formalization.

Let a fixed $ \xbf $ be given. Let $ \xbf $ contain the propositional
variables
$p_{1}, \Xl,p_{m},$ and the operator $ \xbb $ to a depth of nesting $n.$

\bd

$\hspace{0.01em}$

% (+++ Orig. No.:  Definition D-6.3.A.6 +++)

\label{Definition D-6.3.A.6}

(1) Define by induction:

Elementary sets of degree 0:
all intersections of type $M( \xCL p_{1}) \xcs M( \xCL p_{2}) \xcs  \Xl
\xcs M( \xCL p_{m}),$ where $ \xCL p_{j}$ is
either $p_{j}$ or $ \xCN p_{j}.$

Elementary sets of degree $i+1:$
If $s$ is an elementary set of degree $i,$ then $B(s)$ and $s-B(s)$ are
elementary
sets of degree $i+1.$

(2) Unions of degree $i$ are either $ \xCQ $ or (arbitrary) unions of
elementary
sets of degree $i.$

This definition goes up to $i=n,$ though we could continue for all $i \xbe
\xbo,$ but
this is not necessary.

(3) $ \xde $ is the set of all elementary sets of degree $n.$

\ed

\br

$\hspace{0.01em}$

% (+++ Orig. No.:  Remark D-6.3.A.3 +++)

\label{Remark D-6.3.A.3}

(1) for any degree $i,$ the elementary sets of degree $i$ form a disjoint
cover of
$M_{ \xdl }$ - the set of all classical models of the language.

(2) The elementary sets of degree $i+1$ form a refinement of the
elementary sets
of degree $i.$

(3) The set of unions of degree $i$ is closed under union, intersection,
set
difference, i.e. if $ \xCf X,Y$ are unions of degree $i,$ so are $X \xcv
Y,$ $ \xCf X- \xCf Y,$ etc.

(4) If $X$ is a union of degree $i,$ $B(X)$ and $X-B(X)$ are unions of
degree $i+1.$

(5) A union of degree $i$ is also a union of degree $j$ if $j>i.$

$ \xcz $
\\[3ex]

\er

We construct now the new $ \xCf Gm-$structure. We first define the
$opt_{i}$ levels, and
from these levels, the new function $B'.$ Thus, $( \xbg_{0})$ will hold
automatically.
As we know from above example, we may forcibly loose some discerning
power, so
a word where it is lost and where it goes may be adequate. The
construction will
not look inside $B(s)$ and $s-B(s)$ for $s \xbe \xde.$ For $B(s),$ this
is not necessary, as
we know that all elements inside $B(s)$ are on the same level
$(B(s)=B(B(s))),$
inside $s-B(s),$ the original function $B$ may well be able to discern
still
(even infinitely many) different levels, but our fomula $ \xbf $ does not
permit us
to look down at these details - $s-B(s)$ is treated as an atom, one chunk
without
any details inside.

\bd

$\hspace{0.01em}$

% (+++ Orig. No.:  Definition D-6.3.A.7 +++)

\label{Definition D-6.3.A.7}

We define the rank of $X \xbe \xde,$ and the $opt_{i}$ and $B'.$ This is
the central definition.
Let $X \xEd \xCQ,$ $X \xbe \xde.$

(1) Set $rank(X)=0,$ iff there is no $Y \xeb X,$ $Y \xbe \xde.$

Set $rank(X)=i+1,$ iff $X \xbe \xde -\{Z:rank(z) \xck i\}$ and there is no
$Y \xeb X,$
$Y \xbe \xde -\{Z:rank(z) \xck i\}.$

So, $rank(X)$ is the $ \xeb -$level of $X.$

(2) Set $opt_{i}:= \xcV \{X \xbe \xde:rank(X)=i\},$ and $opt_{i}= \xCQ $
iff there is no $X$ s.t. $rank(X)=i.$

(3) $B' (A)$ $:=$ $A \xcs opt_{i},$ where $i$ is the smallest $j$ s.t. $A
\xcs opt_{j} \xEd \xCQ,$ for $A \xEd \xCQ.$

\ed

\br

$\hspace{0.01em}$

% (+++ Orig. No.:  Remark D-6.3.A.4 +++)

\label{Remark D-6.3.A.4}

(1) The $opt_{i}$ form again a disjoint cover of $M_{ \xdl },$ all
$opt_{i}$ are $ \xEd \xCQ $ up to some $k,$
and $ \xCQ $ beyond.

(2) $( \xbg_{0})$ will now hold by definition.

(3) $B' $ is not necessarily $B,$ but sufficiently close.

\er

\bl

$\hspace{0.01em}$

% (+++ Orig. No.:  Lemma D-6.3.A.5 +++)

\label{Lemma D-6.3.A.5}

(Main result)

For any union $X$ of degree $k<n,$ $B(X)=B' (X).$

\el

\subparagraph{
Proof
}

$\hspace{0.01em}$

% (+++ Orig.:  Proof +++)

Write $X$ as a union of degree $n-1,$ and let $ \xdx:=\{X' \xcc X:X' $ is
of degree $n-1\}.$
Note that the construction of $ \xeb / \xCq /opt$ splits all $X' \xbe \xdx
$ into $B(X' )$ and
$X' -B(X' ),$ both of degree $n,$ and that both parts are always present,
there will
not be any isolated $X' -B(X' )$ without its counterpart $B(X' ).$

Let $X' \xbe \xdx.$ Then $(X' -B(X' )) \xcs B' (X)= \xCQ,$ as the
opt-level of $B(X' )$ is better than
the opt-level of $X' -B(X' ).$ Obviously, also $(X' -B(X' )) \xcs B(X)=
\xCQ.$ Thus, $B(X)$ and
$B' (X)$ are the union of certain $B(X' ),$ $X' \xbe \xdx.$
Suppose $B(X' ) \xcc B(X)$ for some $X' \xbe \xdx.$ Then for no $X'' \xbe
\xdx $
$B(B(X' ) \xcv B(X'' )) \xcs B(X' )= \xCQ,$ so $B(X' )$ has minimal
opt-level in $X,$ and
$B(X' ) \xcc B' (X).$ Conversely, let $B(X' ) \xcC B(X).$ Then there is
$X'' \xbe \xdx $ s.t.
$B(B(X' ) \xcv B(X'' )) \xcs B(X' )= \xCQ,$ so $B(X' )$ has not minimal
opt-level in $X,$ and
$B(X' ) \xcs B' (X)= \xCQ.$ $ \xcz $
\\[3ex]

\bco

$\hspace{0.01em}$

% (+++ Orig. No.:  Corollary D-6.3.A.6 +++)

\label{Corollary D-6.3.A.6}

If $ \xbf ' $ is built up from the ingredients of $ \xbf,$ then $[ \xbf '
]=[ \xbf ' ]' $ - where
$[ \xbf ' ]$ $([ \xbf ' ]' )$ is the set of models where $ \xbf ' $ holds
in the original (new)
structure.

\eco

\subparagraph{
Proof
}

$\hspace{0.01em}$

% (+++ Orig.:  Proof +++)

Let $ \xbF $ be the set of formulas which are built up from the
ingredients of $ \xbf,$
i.e. using (some of) the propositional variables of $ \xbf,$ and up to
the nesting
depth of $ \xbb $ of $ \xbf.$

Case 0: Let $p \xbe \xbF $ be a propositional variable. Then $[p]=[p]',$
as we did not
change the classical model. Moreover, $[p]$ is a union of degree 0.

Case 1: Let $ \xbf ', \xbf '' \xbe \xbF,$ and let $[ \xbf ' ]=[ \xbf '
]',$ $[ \xbf '' ]=[ \xbf '' ]' $ be unions of degree
$k' $ and $k'' $ respectively, and let $k:=max(k',k'' ) \xck n.$ Then
both are unions of degree
$k,$ and so are $[ \xbf ' \xcu \xbf '' ],$ $[ \xCN \xbf ' ]$ etc., and $[
\xbf ' \xcu \xbf '' ]=[ \xbf ' \xcu \xbf '' ]',$ as $ \xcu $ is
interpreted by intersection, etc. Let now $k<n.$ We have to show that
$[ \xbb ( \xbf ' )]=[ \xbb ( \xbf ' )]'.$ But $[ \xbb ( \xbf ' )]=B([
\xbf ' ])=B' ([ \xbf ' ])=B' ([ \xbf ' ]' )=[ \xbb ( \xbf ' )]' $ by
above Lemma and induction hypothesis.

$ \xcz $
\\[3ex]

\be

$\hspace{0.01em}$

% (+++ Orig. No.:  Example D-6.3.A.2 +++)

\label{Example D-6.3.A.2}

We interpret $ \xbf = \xbb (p).$
The elementary sets of degree 0 are $M(p),$ $M( \xCN p),$ those of degree
1 are
$M( \xbb (p)),$ $M(p \xcu \xCN \xbb (p)),$ $M( \xbb ( \xCN p)),$ $M( \xCN
p \xcu \xCN \xbb ( \xCN p)),$ suppose the construction
of $ \xeb,$ $ \xCq $ and the opt-levels results in (omitting the $ \xCf
M' $s for clarity)

$p \xcu \xCN \xbb (p),$ $ \xCN p \xcu \xCN \xbb ( \xCN p)$ on the worst
non-empty level

$ \xbb (p)$

$ \xbb ( \xCN p)$ on the best level.

When we calculate now $B' (M(p))$ via opt, we decompose $p$ in its
components $ \xbb (p)$
and $p- \xbb (p),$ and see that $ \xbb (p)$ is on a better level than $p-
\xbb (p),$ so
$ \xbb ' (M(p))= \xbb (M(p)),$ as it should be. $ \xcz $
\\[3ex]

\ee

We finally have:
\subsection{
$Gm \xcl A$ implies $G \xcl A$ (Outline)
}

We show that $ \xcA m(Gm \xcl A)$ implies $G \xcl A$ $( \xCf A$ a $ \xCf
G-$formula).

As L.Aqvist has given equivalent semantics for both systems, we can argue
semantically. We turn the problem round: If $G \xcL A,$ then there is $m$
s.t. $Gm \xcL A.$
Or, in other words, if there is a $ \xCf G-$model and a point in it, where
$ \xCf A$ does
not
hold, then we find an analogue $ \xCf Gm-$model, or, still differently, if
$ \xbf $ is any
$ \xCf G-$formula, and there is a $ \xCf G-$model $ \xbG $ and a point $x$
in $ \xbG,$ where $ \xbf $ holds,
then we find $m$ and $ \xCf Gm-$model $ \xbD,$ and a point $y$ in $ \xbD
$ where $ \xbf $ holds.

By prerequisite, $ \xbf $ contains some propositional variables, perhaps
the
(absolute) quantifiers $M$ and $N,$ usual connectives, and the binary
operators
$O$ and $P.$ Note that the function ``best'' intervenes only in the
interpretation
of $O$ and $P.$ Moreover, the axioms $ \xbs_{i}$ express that best defines
a ranking, e.g.
in the sense of ranked models in preferential reasoning. In addition, $
\xbs_{3}$ is
a limit condition (which essentially excludes unbounded descending
chains).

Let $ \xbf $ contain $n$ symbols. We thus use ``best'' for at most $n$
different $ \xbq,$ where
$O( \xbf ' / \xbq )$ (or $P( \xbf ' / \xbq ))$ is a subformula of $ \xbf
.$

We introduce now more structure into the $ \xCf G-$model. We make $m:=n+1$
layers in
$ \xCf G' $s
universe, where the first $n$ layers are those mentioned in $ \xbf.$ More
precisely,
we put the best $ \xbq -$models for each $ \xbq $ mentioned as above in
its layer -
respecting relations between layers when needed (this is possible, as the
$ \xbs_{i}$
are sufficiently strong), and put all other $ \xbq -$models somewhere
above. The fact
that we have one supplementary layer (which, of course, we put on top)
guarantees that we can do so. The $opt_{i}$ will be the layers.

We then have a $ \xCf Gm-$model (if we take a little care, so nothing gets
empty
prematurely), and $ \xbf $ will hold in our new structure.

$ \xcz $
\\[3ex]
\section{
Hierarchical conditionals
}
\label{Section Hierarchical-Conditionals}
\subsection{
Introduction
}
\label{Section Hier-In}
\index{Section Hier-In}
\label{Section Hier-In-O-Hist-1}
\index{Section Hier-In-O-Hist-1}
\subsubsection{
Description of the problem
}

We often see a hierarchy of situations, e.g.:

 \xEh

 \xDH it is better to prevent an accident than to help the victims,

 \xDH it is better to prove a difficult theorem than to prove an easy
lemma,

 \xDH it is best not to steal, but if we have stolen, we should return
the stolen object to its legal owner, etc.

 \xEj

On the other hand, it is sometimes impossible to achieve the best
objective.

We might have seen the accident happen from far away, so we were unable to
interfere in time to prevent it, but we can still run to the scene and
help
the victims.

We might have seen friends last night and had a drink too many, so today's
headaches will not allow us to do serious work, but we can still prove a
little lemma.

We might have needed a hammer to smash the windows of a car involved in an
accident, so we stole it from a building site, but will return it
afterwards.

We see in all cases:

- a hierarchy of situations

- not all situations are possible or accessible for an agent.

In addition, we often have implicitly a ``normality'' relation:

Normally, we should help the victims, but there might be situations where
not: This would expose ourselves to a very big danger, or this would
involve
neglecting another, even more important task (we are supervisor in a
nuclear
power plant  \Xl.), etc.

Thus, in all ``normal'' situations where an accident seems imminent, we
should
try to prevent it. If this is impossible, in all ``normal'' situations, we
should
help the victims, etc.

We combine these three ideas

 \xEh
 \xDH normality,
 \xDH hierarchy,
 \xDH accessibility
 \xEj

in the present paper.

Note that it might be well possible to give each situation a numerical
value and decide by this value what is right to do - but humans do not
seem
to think this way, and we want to formalize human common sense reasoning.

Before we begin the formal part, we elaborate above situations with more
examples.

 \xEI

 \xDH

We might have the overall intention to advance computer science.

So we apply for the job of head of department of computer science at
Stanford,
and promise every tenured scientist his own laptop.

Unfortunately, we do not get the job, but become head of computer
science department at the local community college. The college does not
have research as priority, but we can still do our best to achieve our
overall intention, by, say buying good books for the library, or buy
computers
for those still active in research, etc.

So, it is reasonable to say that, even if we failed in the best possible
situation - it was not accessible to us - we still succeeded in another
situation, so we achieved the overall goal.

 \xDH

The converse is also possible, where better solutions become possible,
as is illustrated by the following example.

The daughter and her husband say to have the overall intention to start a
family
life with a house of their own, and children.

Suppose the mother now asks her daughter: You have been married now for
two
years, how come you are not pregnant?

Daughter - we cannot afford a baby now, we had to take a huge mortgage to
buy
our house and we both have to work.

Mother - $I$ shall pay off your mortgage. Get on with it!

In this case, what was formerly inaccessible, is now accessible, and
if the daughter was serious about her intentions - the mother can begin
to look for baby carriages.

Note that we do not distinguish here how the situations change, whether
by our own doing, or by someone else's doing, or by some events not
controlled
by anyone.

 \xDH

Consider the following hierarchy of obligations making fences as
unobtrusive as possible, involving contrary to duty obligations.

 \xEh

 \xDH You should have no fence (main duty).

 \xDH If this is impossible (e.g. you have a dog which might invade
neighbours'
property), it should be less than 3 feet high (contrary to duty, but
second
best choice).

 \xDH If this is impossible too (e.g. your dog might jump over it), it
should be
white (even more contrary to duty, but still better than nothing).

 \xDH If all is impossible, you should get the neighbours' consent (etc.).

 \xEj

 \xEJ
\subsubsection{
Outline of the solution
}

The last example can be modelled as follows $( \xbm (x)$ is the minimal
models of $x):$

Layer 1: $ \xbm (True):$ all best models have no fence.

Layer 2: $ \xbm (fence):$ all best models with a
fence are less than 3 ft. high.

Layer 3: $ \xbm (fence$ and more than 3 ft. high): all best models with a
tall
fence have a white fence.

Layer 4: $ \xbm (fence$ and non-white and $ \xcg 3$ ft): in all best
models with a
non-white fence taller than 3 feet, you have permission

Layer 5: all the rest

This will be modelled by a corresponding $ \xda -$structure.

In summary:

 \xEh

 \xDH We have a hierarchy of situations, where one group (e.g. preventing
accidents) is strictly better than another group (e.g. helping victims).

 \xDH Within each group, preferences are not so clear (first help person
A,
or person $B,$ first call ambulance, etc.?).

 \xDH We have a subset of situations which are attainable, this can be
modelled
by an accessibility relation which tells us which situations are possible
or
can be reached.

 \xEj

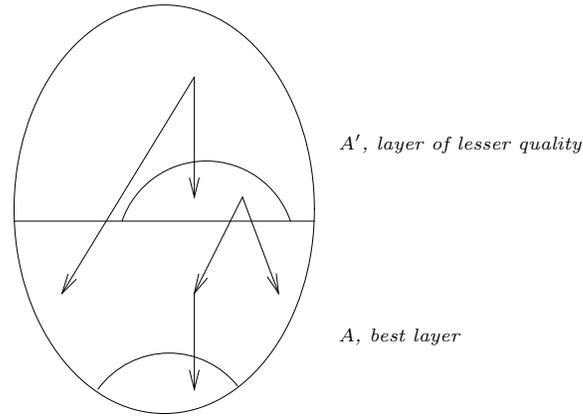
\vspace{20mm}

\begin{diagram}

\label{Diagram A-Ranked}
\index{Diagram A-Ranked}

\centering
\setlength{\unitlength}{0.00083333in}
{\renewcommand{\dashlinestretch}{30}
\begin{picture}(2390,3581)(0,0)
\put(1212.000,2028.000){\arc{1110.000}{3.4719}{5.9529}}
\put(979.651,835.201){\arc{1095.700}{3.7717}{5.6071}}
\put(949,2283){\ellipse{1874}{2550}}
\path(12,2208)(1887,2208)
\path(1137,3108)(1137,2358)
\path(1107.000,2478.000)(1137.000,2358.000)(1167.000,2478.000)
\path(1137,3108)(312,1758)
\path(348.976,1876.037)(312.000,1758.000)(400.172,1844.750)
\path(1437,2358)(1662,1758)
\path(1591.775,1859.826)(1662.000,1758.000)(1647.955,1880.893)
\path(1437,2358)(1137,1758)
\path(1163.833,1878.748)(1137.000,1758.000)(1217.498,1851.915)
\path(1137,1758)(1137,1158)
\path(1107.000,1278.000)(1137.000,1158.000)(1167.000,1278.000)
\put(2037,2658){{\xssc $A'$, layer of lesser quality}}
\put(2037,1458){{\xssc $A$, best layer}}
\put(-700,800){{\xssc
Each layer behaves inside like any preferential structure.}}
\put(-700,600){{\xssc
Amongst each other, layers behave like ranked structures.}}

\put(100,200) {{\rm\bf $\xda-$ ranked structure}}

\end{picture}
}
\end{diagram}

\vspace{4mm}

We combine all three ideas, consider what we call $ \xda -$ranked
structures,
structures which are organized in levels $A_{1},$ $A_{2},$ $A_{3},$ etc.,
where all
elements of $A_{1}$ are better than any element of $A_{2}$ - this is
basically
rankedness -, and where inside each $A_{i}$ we have an arbitrary relation
of
preference. Thus, an $ \xda -$ranked structure is between a simple
preferential
structure and a fully ranked structure.

See Diagram \ref{Diagram A-Ranked} (page \pageref{Diagram A-Ranked}).

Remark: It is not at all necessary that the rankedness relation between
the different layers and the relation inside the layers express the
same concept. For instance, rankedness may express deontic preference,
whereas the inside relation expresses normality or some usualness.

In addition, we have an accessibility relation $R,$ which tells us which
situations are reachable.

It is perhaps easiest to motivate the precise choice of modelling
by layered (or contrary to duty) obligations.

For any point $t,$ let $R(t):=\{s:tRs\},$ the set of $R-$reachable points
from $t.$ Given a preferential structure $ \xdx:= \xBc X, \xeb  \xBe,$ we can
relativize $ \xdx $
by considering only those points in $X,$ which are reachable from $t.$

Let $X' \xcc X,$ and $ \xbm (X' )$ the minimal points of $X,$ we will now
consider
$ \xbm (X' ) \xcs R(t)$ - attention, not: $ \xbm (X' \xcs R(t))!$ This
choice is motivated
by the following: norms are universal, and do not depend on one's
situation $t.$

If $ \xdx $ describes a simple obligation, then we are obliged to $Y$ iff
$ \xbm (X' ) \xcs R(t) \xEd \xCQ,$ and $ \xbm (X' ) \xcs R(t) \xcc Y.$
The first clause excludes
obligations to the unattainable. We can write this as follows, supposing
that $X' $ is the set of models of $ \xbf ',$ and $Y$ is the set of
models of $ \xbq:$

$m \xcm \xbf ' > \xbq.$

Thus, we put the usual consequence relation $ \xcn $ into the object
language
as $>,$ and relativize to the attainable (from $m).$

If an $ \xda -$ranked structure has two or more layers, then we are, if
possible,
obliged to fulfill the lower obligation, e.g. prevent an accident,
but if this is impossible, we are obliged to fulfill the upper
obligation, e.g. help the victims, etc.

See Diagram \ref{Diagram Pischinger} (page \pageref{Diagram Pischinger}).

\vspace{10mm}

\begin{diagram}

\label{Diagram Pischinger}
\index{Diagram Pischinger}

\centering
\setlength{\unitlength}{1mm}
{\renewcommand{\dashlinestretch}{30}
\begin{picture}(150,150)(0,0)

\path(30,10)(30,110)(110,110)(110,10)(30,10)
\path(30,30)(110,30)
\path(30,50)(110,50)
\path(30,70)(110,70)
\path(30,90)(110,90)

\put(70,10){\arc{20}{-3.14}{0}}
\put(70,30){\arc{20}{-3.14}{0}}
\put(70,50){\arc{20}{-3.14}{0}}
\put(70,70){\arc{20}{-3.14}{0}}
\put(70,90){\arc{20}{-3.14}{0}}

\put(60,80){\circle{30}}

\path(5,80)(25,80)
\path(22.3,81)(25,80)(22.3,79)
\put(5,80){\circle*{1}}
\put(25,80){\circle*{1}}

\put(15,75){\xssc{$R$}}
\put(5,75){\xssc{$t$}}
\put(25,75){\xssc{$s$}}

\put(10,130){\xssc{The overall structure is visible from $t$}}
\put(10,123){\xssc{Only the inside of the circle is visible from $s$}}
\put(10,116){\xssc{Half-circles are the sets of minimal elements of layers}}

\put(20,3) {{\rm\bf $\xda-$ ranked structure and accessibility}}

\end{picture}
}

\end{diagram}

\vspace{4mm}

Let now, for simplicity, $ \xdB $ be a subset of the union of all layers
A, and
let $ \xdB $ be the set of models of $ \xbb.$ This can be done, as the
individual
subset can be found by considering $A \xcs \xdB,$ and call the whole
structure
$ \xBc  \xda, \xdB  \xBe.$

Then we say that $m$ satisfies $ \xBc  \xda, \xdB  \xBe $ iff in the lowest
layer A
where
$ \xbm (A) \xcs R(m) \xEd \xCQ $ $ \xbm (A) \xcs R(m) \xcc \xdB.$

When we want a terminology closer to usual conditionals, we may
write e.g. $(A_{1}>B_{1};A_{2}>B_{2}; \Xl.)$ expressing that the best is
$A_{1},$ and
then $B_{1}$ should hold, the second best is $A_{2},$ then $B_{2}$ should
hold, etc.
(The $B_{i}$ are just $A_{i} \xcs \xdB.)$
See Diagram \ref{Diagram C-Validity} (page \pageref{Diagram C-Validity}).
\label{Section Hier-In-O-Hist-2}
\index{Section Hier-In-O-Hist-2}
\subsection{
Formal modelling and summary of results
}

% {\LARGE karl-search= Start Hier-In-Summary }

\label{Section Hier-In-Summary}
\index{Section Hier-In-Summary}

We started with an investigation of ``best fulfillment'' of abstract
requirements,
and contrary to duty obligations. - See
also  \cite{Gab08} and  \cite{Gab08a}.

It soon became evident that semi-ranked preferential structures give a
natural
semantics to contrary to duty obligations, just as simple preferential
structures give a natural semantics to simple obligations - the latter
goes
back to Hansson  \cite{Han69}.

A semi-ranked - or $ \xda -$ranked preferential
structure, as we will call them later, as they are based on a system of
sets $ \xda $ - has a finite number of layers, which amongst them are
totally ordered
by a ranking, but the internal ordering is just any (binary) relation.
It thus has stronger properties than a simple preferential structure, but
not as strong ones as a (totally) ranked structure.

The idea is to put the (cases of the) strongest obligation at
the bottom, and the weaker ones more towards the top. Then, fulfillment of
a strong obligation makes the whole obligation automatically satisfied,
and
the weaker ones are forgotten.

Beyond giving a natural semantics to contrary to duty obligations,
semi-ranked
structures seem very useful for other questions of knowledge
representation.
For instance, any blackbird might seem a more normal bird than any
penguin,
but we might not be so sure within each set of birds.

Thus, this generalization of preferential semantics seems very natural and
welcome.

The second point of this paper is to make some, but not necessarily all,
situations accessible to each point of departure. Thus, if we imagine
agent $ \xCf a$
to be at point $p,$ some fulfillments of the obligation, which are
reachable
to agent $a' $ from point $p' $ might just be impossible to reach for him.
Thus, we introduce a second relation, of accessibility in the intuitive
sense, denoting situations which can be reached. If this relation is
transitive, then we have restrictions on the set of reachable situations:
if $p$ is accessible from $p',$ and $p$ can access situation $s,$ then so
can
$p',$ but not necessarily the other way round.

On the formal side, we characterize:

(1) $ \xda -$ranked structures,

(2) satisfaction of an $ \xda -$ranked conditional once an accessibility
relation
between the points $p,$ $p',$ etc. is given.

For the convience of the reader, we now state the main formal results of
this paper - together with the more unusual definitions.

On (1):

Let $ \xdA $ be a fixed set, and $ \xda $ a finite, totally ordered (by
$<)$ disjoint cover
by non-empty subsets of $ \xdA.$

For $x \xbe \xdA,$ let $rg(x)$ the unique $A \xbe \xda $ such that $x
\xbe A,$ so $rg(x)<rg(y)$ is
defined in the natural way.

A preferential structure $ \xBc  \xdx, \xeb  \xBe $ $( \xdx $ a set of pairs
$ \xBc x,i \xBe )$ is called $ \xda -$ranked
iff for all $x,x' $ $rg(x)<rg(x' )$ implies $ \xBc x,i \xBe  \xeb  \xBc x',i' 
\xBe $ for all
$ \xBc x,i \xBe, \xBc x',i'  \xBe  \xbe \xdx.$
See Definition \ref{Definition Pref-Str} (page \pageref{Definition Pref-Str}) 
for the definition of preferential
structures, and Diagram \ref{Diagram A-Ranked} (page \pageref{Diagram A-Ranked})
 for an illustration.

We then have:

Let $ \xcn $ be a logic for $ \xdl.$ Set $T^{ \xdm }:=Th( \xbm_{ \xdm
}(M(T))),$ and $ \ol{ \ol{T} }:=\{ \xbf:T \xcn \xbf \}.$
where $ \xdm $ is a preferential structure.

(1) Then there is a (transitive) definability preserving
classical preferential model $ \xdm $ s.t. $ \ol{ \ol{T} }=T^{ \xdm }$ iff

(LLE), (CCL), (SC), (PR) hold for all $T,T' \xcc \xdl.$

(2) The structure can be chosen smooth, iff, in addition

(CUM) holds.

(3) The structure can be chosen $ \xda -$ranked, iff, in addition

$( \xda -$min) $T \xcL \xCN \xba_{i}$ and $T \xcL \xCN \xba_{j},$ $i<j$
implies $ \ol{ \ol{T} } \xcl \xCN \xba_{j}$

holds.

See Definition \ref{Definition Pref-Log} (page \pageref{Definition Pref-Log}) 
for the logic defined by a
preferential
structure, Definition \ref{Definition Log-Cond-Ref-Size} (page
\pageref{Definition Log-Cond-Ref-Size})
for the logical conditions,
Definition \ref{Definition Smooth} (page \pageref{Definition Smooth})  for
smoothness.

On (2)

Given a transitive accessibility relation $R,$ $R(m):=\{x:mRx\}.$

Given $ \xda $ as above, let $ \xdB \xcc \xdA $ be the set of ``good''
points in $ \xdA,$ and
set $ \xdc:= \xBc  \xda, \xdB  \xBe.$

We define:

(1) $ \xbm ( \xda ):= \xcV \{ \xbm (A_{i}):i \xbe I\}$

(warning: this is NOT $ \xbm ( \xdA ))$

(2) $ \xda_{m}:=R(m) \xcs \xdA,$

(3) $ \xbm ( \xda_{m}):= \xcV \{ \xbm (A_{i}) \xcs R(m):i \xbe I\}$

(3a) $ \xbn ( \xda_{m})$ $:=$ $ \xbm ( \xbm ( \xda_{m}))$

(thus $ \xbn ( \xda_{m})$ $=$ $\{a \xbe \xdA:$ $ \xcE A \xbe \xda (a \xbe
\xbm (A),$ $a \xbe R(m),$ and

$ \xDC \xDC \xDC \xCN \xcE a' ( \xcE A' \xbe \xda (a' \xbe \xbm (A' ),$
$a' \xbe R(m),$ $a' \xeb a\}.$

(4) $m \xcm \xdc $ $: \xcr $ $ \xbn ( \xda_{m})) \xcc \xdB.$

See Diagram \ref{Diagram C-Validity} (page \pageref{Diagram C-Validity})

Then the following hold:

Let $m,m' \xbe M,$ $A,A' \xbe \xda,$ $ \xdA $ be the set of models of $
\xba.$

(1) $m \xcm \xcX \xCN \xba,$ mRm' $ \xch $ $m' \xcm \xcX \xCN \xba $

(2) $ \xCf mRm',$ $ \xbn ( \xda_{m}) \xcs A \xEd \xCQ,$ $ \xbn (
\xda_{m' }) \xcs A' \xEd \xCQ,$ $ \xch $ $A \xck A' $ (in the ranking)

(3) $ \xCf mRm',$ $ \xbn ( \xda_{m}) \xcs A \xEd \xCQ,$ $ \xbn (
\xda_{m' }) \xcs A' \xEd \xCQ,$ $m \xcm \xdc,$ $m' \xcM \xdc, \xch $
$A<A' $

Conversely, these conditions suffice to construct an accessibility
relation between $M$ and $ \xdA $ satisfying them, so they are sound and
complete.

% karl-search= End Hier-In-DovIn-Motiv
\vspace{7mm}

% *************************************

\vspace{7mm}

% karl-search= End Hier-In-DovIn
\vspace{7mm}

% *************************************

\vspace{7mm}

% karl-search= End Hier-In-Summary
\vspace{7mm}

% *************************************

\vspace{7mm}

\subsection{
Overview
}

We next point out some connections with other domains of artificial
intelligence and computer science.

We then put our work in perspective with a summary of logical and
semantical conditions for nonmonotonic and related logics, and present
basic defintions for preferential structures.

Next, we will give special definitions for our framework.

We then start the main formal part, and prove representation results
for $ \xda -$ranked structures, first for the general case, then for
the smooth case. The general case needs more work, as we have to
do a (minor) modification of the not $ \xda -$ranked case. The smooth case
is
easy, we simply have to append a small construction. Both proofs are
given in full detail, in order to make the text self-contained.

Finally, we characterize changes due to restricted accessibility.
\label{Section Hier-Def-1}
\index{Section Hier-Def-1}

\bd

$\hspace{0.01em}$

% (+++ Orig. No.:  Definition Hier-2.1 +++)

\label{Definition Hier-2.1}

We have the usual framework of preferential structures, i.e.
either a set with a possibly non-injective labelling function, or,
equivalently,
a set of possible worlds with copies. The relation of the
preferential structure will be fixed, and will not depend on the point $m$
from
where we look at it.

Next, we have a set $ \xdA,$ and a finite, disjoint cover $A_{i}:i<n$ of
$ \xda,$
with a relation ``of quality'' $<,$ $ \xda $ will denote the $A_{i}$ (and
thus $ \xdA ),$ and $<,$
i.e. $ \xda = \xBc \{A_{i}:i \xbe I\},< \xBe.$

By Fact \ref{Fact 6.21} (page \pageref{Fact 6.21}), we may assume that all
$A_{i}$ are described
by a formula.

Finally, we have $ \xdB \xcc \xdA,$ the subset of ``good'' elements of $
\xdA $ - which
we also assume to be described by a formula.

In addition, we have a binary relation of accessibility, $R,$ which we
assume
transitive - modal operators will be defined relative to $R.$ $R$
determines which
part of the preferential structure is visible.

Let $R(s):=\{t:sRt\}.$

\ed

\bd

$\hspace{0.01em}$

% (+++ Orig. No.:  Definition 2.2 +++)

\label{Definition 2.2}

We repeat here from the introduction, and assume $A_{i}=M( \xba_{i}),$
$B=M( \xbb ),$ and
$ \xbm $ expresses the minimality of the preferential structure.

$t \xcm \xba_{i}> \xbb: \xcj $ $ \xbm (A_{i}) \xcs R(t) \xcc B,$

we will also abuse notation and just write

$t \xcm A_{i}>B$ in this case.

We then define:

$t \xcm \xdc $ iff at the smallest $i$ s.t. $ \xbm (A_{i}) \xcs R(t) \xEd
\xCQ,$ $ \xbm (A_{i}) \xcs R(t) \xcc \xdB $ holds.

\ed

This motivates Definition \ref{Definition A-ranked} (page \pageref{Definition
A-ranked}).
\label{Section Hier-Def-3}
\index{Section Hier-Def-3}

Note that automatically for $X \xcc \xdA,$ $ \xbm (X) \xcc A_{j}$ when
$j$ is the smallest $i$ s.t.
$X \xcs A_{i} \xEd \xCQ.$

The idea is now to make the $A_{i}$ the layers, and ``trigger'' the first
layer
$A_{j}$ s.t. $ \xbm (A_{j}) \xcs R(x) \xEd \xCQ,$ and check whether $
\xbm (A_{j}) \xcs R(x) \xcc B_{j}.$ A suitable
ranked structure will automatically find this $A_{j}.$

More definitions and results for such $ \xda $ and $ \xdc $ will be found
in Section \ref{Section Hier-CondRepr} (page \pageref{Section Hier-CondRepr}).
\subsection{
Connections with other concepts
}
\subsubsection{
Hierarchical conditionals and programs
}

% {\LARGE karl-search= Start Hier-Cond }

\label{Section Hier-Cond}
\index{Section Hier-Cond}

Our situation is now very similar to a sequence of computer program
instructions:

if $A_{1}$ then do $B_{1};$

else if $A_{2}$ then do $B_{2};$

else if $A_{3}$ then do $B_{3};$

where we can see the $B_{i}$ as subroutines.

We can deepen this analogy in two directions:

(1) connect it to Update

(2) put an imperative touch to it.

In both cases, we differentiate between different degrees of fulfillment
of $ \xdc:$ the lower the level is which is fulfilled, the better.

(1) We can consider all threads of reachability which lead to a model $m$
where
$m \xcm \xdc.$ Then we take as best threads those which lead to the best
fulfillment of
$ \xdc.$ So degree of fulfillment gives the order by which we should do
the update.
(This is then not update in the sense that we choose the most normal
developments, but rather we actively decide for the most desirable ones.)
We will not pursue this line any further here, but leave it for future
research.

(2): We introduce an imperative operator, say!.! means that one should
fulfill $ \xdc $ as best as possible by suitable choices.
We will elaborate this now.

First, we can easily compare the degree of satisfaction of $ \xdc $ of two
models:

\bd

$\hspace{0.01em}$

% (+++ Orig. No.:  Definition 3.1 +++)

\label{Definition 3.1}

Let $m,m' \xcm \xdc,$ and define $m<m' $ $: \xch $ $ \xbm ( \xbm (
\xda_{m}) \xcv \xbm ( \xda_{m' })) \xcs \xbm ( \xda_{m' })= \xCQ.$ $(
\xbm $ is,
as usual, relative to some fixed $ \xck_{t}.)$

\ed

For two sets of models, $X,$ $X',$ the situation does not seem so easy.
So suppose
that $X,X' \xcm \xdc.$ First, we
have to decide how to compare this, we do by the maximum: $X<X' $ iff the
worst
satisfaction of all $x \xbe X$ is better than the worst satisfaction in
$X'.$
More precisely, we look at all $ \xbg ( \xdc )$ for all $x \xbe X,$ take
the maximum (which
exists, as $ \xda $ is finite), and then compare the maxima for $X$ and
for $X'.$

Suppose now that there are points where we can make decisions $('' free$
$will'' ),$
let $m$ be such a point. We introduce a new relation $D,$ and let mDm' iff
we can
decide to go from $m$ to $m'.$ The relation $D$ expresses this
possibility - it is
our definition of ``free will''.

\bd

$\hspace{0.01em}$

% (+++ Orig. No.:  Definition 3.2 +++)

\label{Definition 3.2}

Consider now some formula $ \xbf,$ and define

$m \xcm! \xbf $ $: \xch $ $D(m) \xcs M( \xbf )<D(m) \xcs M( \xCN \xbf )$

(as defined in Definition \ref{Definition 3.1} (page \pageref{Definition 3.1})
).

% karl-search= End Hier-Cond
\vspace{7mm}

% *************************************

\vspace{7mm}

\subsubsection{
Connection with Theory Revision
}

% {\LARGE karl-search= Start Hier-TR }

\label{Section Hier-TR}
\index{Section Hier-TR}

\ed

In particular, the situation of contrary to duty obligations
(see Section \ref{Section Hier-In} (page \pageref{Section Hier-In}) ) shows an
intuitive
similarity to revision. You have the duty not to have a fence. If this is
impossible (read: inconsistent), then it should be white. So the duty is
revised.

But there is also a formal analogy: As is well known, AGM revision (with
fixed left hand side $K)$ corresponds to a ranked order of models, where
models of $K$ have lowest rank (or: distance 0 from $K-$models). The
structures we consider $( \xda -$rankings) are partially ranked, i.e.
there is only a partial ranked preference, inside the layers, nothing is
said about the ordering. This partial ranking is natural, as we have only
a limited number of cases to consider.

But we use the revision order (based on $K,$ so it really is a $ \xck_{K}$
relation)
differently: We do not revise $K,$ but use only the order to choose the
first
layer which has non-empty intersection with the set of possible cases.
Still, the spirit (and formal apparatus) of revision is there, just used
somewhat differently. The $K-$relation expresses here deontic quality, and
if
the best situation is impossible, we choose the second best, etc.

Theory revision with variable $K$ is expressed by a distance between
models
(see  \cite{LMS01}),
where $K* \xbf $ is defined by the set of $ \xbf $ models which have
minimal distance from
the set of $K$ models.

We can now generalize our idea of layered structure to a partial distance
as
follows: For instance, $d(K,A)$ is defined, $d(K,B)$ too, and we know that
all A models with minimal distance to $K$ have smaller distance than the
$B$ models
with minimal distance to $K.$ But we do NOT know a precise distance for
other A
models, we can sometimes compare, but not always. We may also know that
all
A models are closer to $K$ than any $B$ model is, but for a and $a',$
both A models,
we might not know if one or the other is closer to $K,$ or is they have
the
same distance.

% karl-search= End Hier-TR
\vspace{7mm}

% *************************************

\vspace{7mm}

The representation results for $ \xda -$ranked structures were shown
already
in Section \ref{Section A-ranked-Rep} (page \pageref{Section A-ranked-Rep}),
so we can turn immediately to the following:
\subsection{
Formal results and representation for hierarchical conditionals
}
\label{Section Hier-CondRepr}
\index{Section Hier-CondRepr}

We look here at the following problem:

Given

(1.1) a finite, ordered partition $ \xda $ of $ \xdA,$ $ \xda = \xBc
\{A_{i}:i \xbe I\},< \xBe $

(1.2) a normality relation $ \xeb,$ which is an $ \xda -$ranking,
defining a
choice function $ \xbm $ on subsets of $ \xdA,$ (so, obviously,
$A<A' $ iff $ \xbm (A \xcv A' ) \xcs A' = \xCQ ),$

(1.3) a subset $ \xdB \xcc \xdA,$ and we set $ \xdc:= \xBc  \xda, \xdB  \xBe $
(thus, the $B_{i}$ are just $A_{i} \xcs \xdB,$
this way of writing saves a little notation),

(2.1) a set of models $M,$

(2.2) an accessibility relation $R$ on $M,$ with some finite upper bound
on
$R-$chains,

(2.3) an unknown extension of $R$ to pairs $(m,a),$ $m \xbe M,$ $a \xbe
\xdA,$

(3.1) a notion of validity $m \xcm \xdc,$ for $m \xbe M,$ defined by
$m \xcm \xdc $ iff
$\{a \xbe \xdA:$ $ \xcE A \xbe \xda (a \xbe \xbm (A),$ $a \xbe R(m),$ and

$ \xDC \xDC \xDC \xCN \xcE a' ( \xcE A' \xbe \xda (a' \xbe \xbm (A' ),$
$a' \xbe R(m),$ $a' \xeb a\}$ $ \xcc $ $ \xdB,$

(3.2) a subset $M' $ of $M$

give a criterion which decides whether it is possible to construct the
extension of $R$ to pairs $(m,a)$ s.t. $ \xcA m \xbe M.(m \xbe M' $ $ \xcj
$ $m \xcm \xdc ).$

We first show some elementary facts on the situation, and give the
criterion
in Proposition \ref{Proposition 7.4} (page \pageref{Proposition 7.4}), together
with the proof that it does
what is wanted.

\bfa

$\hspace{0.01em}$

% (+++ Orig. No.:  Fact 7.1 +++)

\label{Fact 7.1}

Reachability for a transitive relation is characterized by

$y \xbe R(x)$ $ \xcp $ $R(y) \xcc R(x)$

\efa

\subparagraph{
Proof
}

$\hspace{0.01em}$

% (+++ Orig.:  Proof +++)

Define directly xRz iff $z \xbe R(x).$ This does it. $ \xcz $
\\[3ex]

Let now $S$ be a set with an accessibility relation $R',$ generated by
transitive
closure from the intransitive subrelation $R.$ All modal notation will be
relative
to this $R.$

Let $ \xdA =M( \xba ),$ $A_{i}=M( \xba_{i}),$ the latter is justified by
Fact \ref{Fact 6.21} (page \pageref{Fact 6.21}).

\bd

$\hspace{0.01em}$

% (+++ Orig. No.:  Definition 7.1 +++)

\label{Definition 7.1}

(1) $ \xbm ( \xda ):= \xcV \{ \xbm (A_{i}):i \xbe I\}$

(warning: this is NOT $ \xbm ( \xdA ))$

(2) $ \xda_{m}:=R(m) \xcs \xdA,$

(3) $ \xbm ( \xda_{m}):= \xcV \{ \xbm (A_{i}) \xcs R(m):i \xbe I\}$

(3a) $ \xbn ( \xda_{m})$ $:=$ $ \xbm ( \xbm ( \xda_{m}))$

(thus $ \xbn ( \xda_{m})$ $=$ $\{a \xbe \xdA:$ $ \xcE A \xbe \xda (a \xbe
\xbm (A),$ $a \xbe R(m),$ and

$ \xDC \xDC \xDC \xCN \xcE a' ( \xcE A' \xbe \xda (a' \xbe \xbm (A' ),$
$a' \xbe R(m),$ $a' \xeb a\}.$

(4) $m \xcm \xdc $ $: \xcr $ $ \xbn ( \xda_{m})) \xcc \xdB.$

See Diagram \ref{Diagram C-Validity} (page \pageref{Diagram C-Validity})

\vspace{20mm}

\begin{diagram}

\label{Diagram C-Validity}
\index{Diagram C-Validity}

\centering
\setlength{\unitlength}{0.00083333in}
{\renewcommand{\dashlinestretch}{30}
\begin{picture}(3841,3665)(0,200)
\put(2407.000,2628.000){\arc{1110.000}{3.4719}{5.9529}}
\put(2407.000,1578.000){\arc{1110.000}{3.4719}{5.9529}}
\put(2398.149,661.355){\arc{1503.312}{3.9463}{5.4292}}
\put(-24.763,2227.012){\arc{4826.591}{5.7407}{6.7841}}
\put(2369,2283){\ellipse{1874}{2550}}
\path(1507,2808)(3232,2808)
\path(1507,1758)(3232,1758)
\path(907,3263)(2562,1963)
\path(2449.100,2013.534)(2562.000,1963.000)(2486.163,2060.718)
\path(157,3248)(2072,1213)
\path(1967.915,1279.831)(2072.000,1213.000)(2011.611,1320.950)
\path(272,3258)(847,3258)
\put(210,3258){\circle*{20}}
\put(880,3258){\circle*{20}}
\path(727.000,3228.000)(847.000,3258.000)(727.000,3288.000)
\put(3457,3108){{\xssc $A''$}}
\put(3457,2208){{\xssc $A'$}}
\put(3457,1308){{\xssc $A$}}
\put(4000,1093){{\xssc $m \xcm \xdc$}}
\put(4000,8033){{\xssc $m' \xcM \xdc$}}
\put(2407,2883){{\xssc $\mu(A'')$}}
\put(2407,1833){{\xssc $\mu(A')$}}
\put(2407,1093){{\xssc $\mu(A)$}}
\put(132,3318){{\xssc $m$}}
\put(882,3313){{\xssc $m'$}}
\put(1972,3518){{\xssc $B$}}
\put(10,800){{\xssc Here, the ``best'' element $m$ sees is in
$B$, so $\xdc$ holds in $m$.}}
\put(10,600){{\xssc The ``best'' element $m'$ sees is not in
$B$, so $\xdc$ does not hold in $m'$.}}

\put(100,250) {{\rm\bf Validity of $\xdc$ from $m$ and $m'$}}

\end{picture}
}
\end{diagram}

\vspace{4mm}

\ed

We have the following Fact for $m \xcm \xdc:$

\bfa

$\hspace{0.01em}$

% (+++ Orig. No.:  Fact 7.2 +++)

\label{Fact 7.2}

Let $m,m' \xbe M,$ $A,A' \xbe \xda.$

(1) $m \xcm \xcX \xCN \xba,$ mRm' $ \xch $ $m' \xcm \xcX \xCN \xba $

(2) $ \xCf mRm',$ $ \xbn ( \xda_{m}) \xcs A \xEd \xCQ,$ $ \xbn (
\xda_{m' }) \xcs A' \xEd \xCQ,$ $ \xch $ $A \xck A' $

(3) $ \xCf mRm',$ $ \xbn ( \xda_{m}) \xcs A \xEd \xCQ,$ $ \xbn (
\xda_{m' }) \xcs A' \xEd \xCQ,$ $m \xcm \xdc,$ $m' \xcM \xdc, \xch $
$A<A' $

\efa

\subparagraph{
Proof
}

$\hspace{0.01em}$

% (+++ Orig.:  Proof +++)

Trivial. $ \xcz $
\\[3ex]

\bfa

$\hspace{0.01em}$

% (+++ Orig. No.:  Fact 7.3 +++)

\label{Fact 7.3}

We can conclude from above properties that there are no arbitrarily long
$R-$chains of models $m,$ changing from $m \xcm \xdc $ to $m \xcM \xdc $
and back.

\efa

\subparagraph{
Proof
}

$\hspace{0.01em}$

% (+++ Orig.:  Proof +++)

Trivial: By Fact \ref{Fact 7.2} (page \pageref{Fact 7.2}), (3), any change from
$ \xcm \xdc $
to $ \xcM \xdc $ results in a
strict increase in rank. $ \xcz $
\\[3ex]

We solve now the representation task described at the beginning of
Section \ref{Section Hier-CondRepr} (page \pageref{Section Hier-CondRepr}), all
we need
are the properties shown in Fact \ref{Fact 7.2} (page \pageref{Fact 7.2}).

(Note that constructing $R$ between the different $m,$ $m' $ is trivial:
we could just
choose the empty relation.)

\bp

$\hspace{0.01em}$

% (+++ Orig. No.:  Proposition 7.4 +++)

\label{Proposition 7.4}

If the properties of Fact \ref{Fact 7.2} (page \pageref{Fact 7.2})  hold, we can
extend $R$ to
solve the
representation problem described at the beginning of
this Section \ref{Section Hier-CondRepr} (page \pageref{Section Hier-CondRepr})
.

\ep

\subparagraph{
Proof
}

$\hspace{0.01em}$

% (+++ Orig.:  Proof +++)

By induction on $R.$ This is possible, as the depth of $R$ on $M$ was
assumed to
be finite.

\bcs

$\hspace{0.01em}$

% (+++ Orig. No.:  Construction 7.1 +++)

\label{Construction 7.1}

We choose now elements as possible, which ones are chosen exactly does not
matter.

$X_{i}:=\{b_{i},c_{i}\}$ iff $ \xbm (A_{i}) \xcs \xdB \xEd \xCQ $ and $
\xbm (A_{i})- \xdB \xEd \xCQ,$ $b_{i} \xbe \xbm (A_{i}) \xcs \xdB,$
$c_{i} \xbe \xbm (A_{i})- \xdB.$

$X_{i}:=\{c_{i}\}$ iff $ \xbm (A_{i}) \xcs \xdB = \xCQ $ and $ \xbm
(A_{i})- \xdB \xEd \xCQ,$ $c_{i} \xbe \xbm (A_{i})- \xdB $

$X_{i}:=\{b_{i}\}$ iff $ \xbm (A_{i}) \xcs \xdB \xEd \xCQ $ and $ \xbm
(A_{i})- \xdB = \xCQ,$ $b_{i} \xbe \xbm (A_{i}) \xcs \xdB,$

$X_{i}:= \xCQ $ iff $ \xbm (A_{i})= \xCQ.$

Case 1:

Let $m$ be $R-$minimal and $m \xcm \xdc.$ Let $i_{0}$ be the first $i$
s.t. $b_{i} \xbe X_{i},$ make
$ \xbg (m):=i_{0},$ and make $R(m):=\{b_{i_{0}}\} \xcv \xcV
\{X_{i}:i>i_{0}\}.$ This makes $ \xdc $ hold.
(This leaves us as many possibilities open as possible - remember we have
to
decrease the set of reachable elements now.)

Case 2:

Let $m$ be $R-$minimal and $m \xcM \xdc.$ Let $i_{0}$ be the first $i$
s.t. $c_{i} \xbe X_{i},$ make
$ \xbg (m):=i_{0},$ and make $R(m):= \xcV \{X_{i}:i \xcg i_{0}\}.$ This
makes $ \xdc $ false.

Let all $R-$predecessors of $m$ be determined, and $i:=max\{ \xbg (m' ):m'
Rm\}.$

Case 3: $m \xcm \xdc.$ Let $j$ be the smallest $i' \xcg i$ with $ \xbm
(A_{i' }) \xcs \xdB \xEd \xCQ.$
Let $R(m):=\{b_{j}\} \xcv \xcV \{X_{k}:k>j\},$ and $ \xbg (m):=j.$

Case 4: $m \xcM \xdc.$

Case 4.1: For all $m' Rm$ with $i= \xbg (m' )$ $m' \xcM \xdc.$

Take one such $m' $ and set $R(m):=R(m' ),$ $ \xbg (m):=i.$

Case 4.2: There is $m' Rm$ with $i= \xbb (m' )$ $m' \xcm \xdc.$

Let $j$ be the smallest $i' >i$ with $ \xbm (A_{i' })- \xdB \xEd \xCQ.$
Let $R(m):= \xcV \{X_{k}:k \xcg j\}.$
(Remark: To go from $ \xcm $ to $ \xcM,$ we have to go higher in the
hierarchy.)

Obviously, validity is done as it should be.
It remains to show that the sets of reachable elements decrease with $R.$

\ecs

\bfa

$\hspace{0.01em}$

% (+++ Orig. No.:  Fact 7.5 +++)

\label{Fact 7.5}

In above construction, if mRm', then $R(m' ) \xcc R(m).$

\efa

\subparagraph{
Proof
}

$\hspace{0.01em}$

% (+++ Orig.:  Proof +++)

By induction, considering $R.$ $ \xcz $ (Fact \ref{Fact 7.5} (page \pageref{Fact
7.5})  and
Proposition \ref{Proposition 7.4} (page \pageref{Proposition 7.4}) )
\\[3ex]

We consider an example for illustration.

\be

$\hspace{0.01em}$

% (+++ Orig. No.:  Example 7.1 +++)

\label{Example 7.1}

Let $a_{1}Ra_{2}RcRc_{1},$ $b_{1}Rb_{2}Rb_{3}RcRd_{1}Rd_{2}.$

Let $ \xdc =(A_{1}>B_{1}, \Xl,A_{n}>B_{n})$ with the $C_{i}$ consistency
with $ \xbm (A_{i}).$

Let $ \xbm (A_{2}) \xcs B_{2}= \xCQ,$ $ \xbm (A_{3}) \xcc B_{3},$ and for
the other $i$ hold neither of these two.

Let $a_{1},a_{2},b_{2},c_{1},d_{2} \xcm \xdc,$ the others $ \xcM \xdc.$

Let $ \xbm (A_{1})=\{a_{1,1},a_{1,2}\},$ with $a_{1,1} \xbe B_{1},$
$a_{1,2} \xce B_{1},$

$ \xbm (A_{2})=\{a_{2,1}\},$ $ \xbm (A_{3})=\{a_{3,1}\}$ (there is no
reason to differentiate),

and the others like $ \xbm (A_{1}).$ Let $ \xbm A:= \xcV \{ \xbm (A_{i}):i
\xck n\}.$

We have to start at $a_{1}$ and $b_{1},$ and make $R(x)$ progressively
smaller.

Let $R(a_{1}):= \xbm A-\{a_{1,2}\},$ so $a_{1} \xcm \xdc.$ Let
$R(a_{2})=R(a_{1}),$ so again $a_{2} \xcm \xdc.$

Let $R(b_{1}):= \xbm A-\{a_{1,1}\},$ so $b_{1} \xcM \xdc.$ We now have to
take $a_{1,2}$ away, but
$a_{2,1}$ too to be able to change. So let
$R(b_{2}):=R(b_{1})-\{a_{1,2},a_{2,1}\},$ so we
begin at $ \xbm (A_{3}),$ which is a (positive) singleton.
Then let $R(b_{3}):=R(b_{2})-\{a_{3,1}\}.$

We can choose $R(c):=R(b_{3}),$ as $R(b_{3}) \xcc R(a_{2}).$

Let $R(c_{1}):=R(c)-\{a_{4,2}\}$ to make $ \xdc $ hold again. Let
$R(d_{1}):=R(c),$ and $R(d_{2}):=R(c_{1}).$

$ \xcz $
\\[3ex]
\chapter{Theory update and theory revision}
\label{Section 2.3}
\section{Update}
\subsection{Introduction}
\label{Section 2.3.1}

\ee

We will treat here problems due to lack of information, i.e.
we can ``see'' some dimensions, but not all.
\subsection{Hidden dimensions}
\label{Section 2.3.2}
\subsubsection{Introduction}
\label{Section 2.3.2.1}

We look here at situations where only one dimension is visible in the
results,
and the other ones stay hidden. This is e.g. the case when we can observe
only
the outcome of developments, but not the developments themselves.

It was the authors' intention to treat here the general infinite case, and
then
show that the problems treated in  \cite{BLS99} and  \cite{LMS01}
(the not necessarily symmetric case there) are special cases thereof.
Unfortunately, we failed in the attempt to solve the general
infinite case, it seems that one needs new and quite different methods to
solve
it, so we will just describe modestly what we see that can be done, what
the
problems seem to be, and conclude with a very short remark on the
situation
described in  \cite{BLS99}.
\subsubsection{Situation, definitions, and basic results}
\label{Section 2.3.2.2}

In several situations, we can observe directly only one dimension of a
problem.
In a classical ranked structure, we ``see'' everything about an optimal
model.
It is there, and fully described. But look at a ranking of developments,
where
we can observe only the outcome. The earlier dimensions remain hidden, and
when
we see the outcome of the ``best'' developments, we do not directly see the
threads which led there. A similar case is theory revision based on not
necessarily symmetric distances, where we cannot ``look back'' from the
result to
see the closest elements of the former theory (see  \cite{LMS01} for
details).

The non-definable case and the case of hidden dimensions are different
aspects
of a common problem: In the case of non-definability, any not too small
subset
might generate what we see, in the case of hidden dimensions, any thread
ending in an optimal outcome might be optimal.

\paragraph{
The situation, more formally
}

$\hspace{0.01em}$

% (+++ Orig.:  The situation, more formally +++)

\label{Section The situation, more formally}

The universe is a finite or even infinite product $ \xbP U_{i},$ $i \xbe
I.$ We will see that
the finite case is already sufficiently nasty, so we will consider only
the finite situation. If $X \xcc \xbP U_{i},$ then possible results will
be projections
on some fixed coordinate, say $j,$ of the best $ \xbs \xbe X,$ where best
is determined by
some ranking $ \xeb $ on $ \xbP U_{i},$ $ \xbp_{j}( \xbm (X)).$

As input, we will usually not have arbitrary $X \xcc \xbP U_{i},$ but
again some product
$X:= \xbP U'_{i},$ with $U'_{i} \xcc U_{i}.$ Here is the main problem: we
cannot use as input
arbitrary sets, but only products. We will see that this complication will
hide almost all information in sufficiently nasty situations.

We will make now some reasonable assumptions:

First, without loss of generality, we will always take the last dimension
as
outcome. Obviously, this does not change the general picture.

Second, the difficult situations are those where (some of) the $U_{i}$ are
infinite.
We will take the infinite propositional case with theories as input as
motivation, and assume that for each $U_{i},$ we can choose any finite
$U'_{i}$ as input,
and the possible $U'_{i}$ are closed under intersections and finite
unions.

Of course, $ \xbP U'_{i} \xcv \xbP U''_{i}$ need not be a product $ \xbP
V_{i}$ - here lies the main
problem, the domain is not closed under finite unions.

We will see that the case presents serious difficulties, and which the
authors
do not know how to solve. The basic problem is that we do not have enough
information to construct an order.

\bn

$\hspace{0.01em}$

% (+++ Orig. No.:  Notation D-6.3.1 +++)

\label{Notation D-6.3.1}

(.)! will be the projection on the fixed, (last) coordinate, so for $X
\xcc \xbP U_{i},$
$i \xck n,$ $X!:=\{ \xbs (n): \xbs \xbe X\},$ and analogously, $ \xbs!:=
\xbs (n)$ for $ \xbs \xbe \xbP U_{i}.$

To avoid excessive parentheses, we will also write $ \xbm X!$ for $( \xbm
(X))!,$ where
$ \xbm $ is the minimal choice operator, etc.

For any set of sequences $X,$ [X] will be the smallest product which
contains $X.$

Likewise, for $ \xbs,$ $ \xbs ',$ $[ \xbs, \xbs ' ]$ will be the
product (which will always be
defined as we have singletons and finite unions in the components),
$ \xbP \{ \xbs (i), \xbs ' (i)\}$ $(i \xbe I).$ Analogously, $( \xbs,
\xbs ' ):=[ \xbs, \xbs ' ]-\{ \xbs, \xbs ' \}.$ (The
interval notation is intended.)

If $ \xbs \xbe X,$ a $ \xbs -$cover of $X$ will be a set $\{X_{k}:k \xbe
K\}$ s.t. $ \xbs \xce X_{k}$ for all $k,$ and
$ \xcV \{X_{k}:k \xbe K\} \xcv \{ \xbs \}=X.$

A (finite) sequence $ \xbs = \xbs_{1}, \Xl, \xbs_{n}$ will also be
written $( \xbs_{1}, \Xl, \xbs_{n}).$

The Hamming distance is very important in our situation. We can define
the Hamming distance between two (finite) sequences as the number of
arguments
where they disagree, or as the set of those arguments. In our context, it
does
not matter which definition we choose. `` $H-$closest to $ \xbs $ '' means
thus
``closest to $ \xbs,$ measured by (one of) the Hamming $distance(s)$ ''.

\en

\paragraph{
Direct and indirect observation
}

$\hspace{0.01em}$

% (+++ Orig.:  Direct and indirect observation +++)

\label{Section Direct and indirect observation}

If $ \xbs!= \xbs '!,$ we cannot compare them directly, as we always see
the same
projection. If $ \xbs! \xEd \xbs '!,$ we can compare them, if e.g.
$X:=\{ \xbs, \xbs ' \}$ is in
the domain of $ \xbm $ (which usually is not the case, but only if they
differ only
in the last coordinate), $ \xbm (X)!$ can be $\{ \xbs!\},$ $\{ \xbs '
!\},$ $\{ \xbs!, \xbs '!\},$ so $ \xbs $ is
better, $ \xbs ' $ is better, or they are equal.

Thus, to compare $ \xbs $ and $ \xbs ' $ if $ \xbs!= \xbs '!,$ we have
to take a detour, via some
$ \xbt $ with $ \xbt! \xEd \xbs!.$ This is illustrated by the following
example:

\be

$\hspace{0.01em}$

% (+++ Orig. No.:  Example D-6.3.1 +++)

\label{Example D-6.3.1}

Let a set A be given, $B:=\{b,b' \},$ and $ \xbm (A \xCK B)!=\{b\}.$
Of course, without any further information, we have no
idea if for all $a \xbe A$ $ \xBc a,b \xBe $ is optimal, if only for some, etc.

Consider now the following situation: Let $a,a' \xbe A$ be given, $A'
:=\{a,a' \},$ and
$C:=\{b,c,c' \}$ - where $c,c' $ need not be in $B.$ Let $A' \xCK C,$
$\{a\} \xCK C,$ $\{a' \} \xCK C$ be
defined and observable, and

(1) $ \xbm (A' \xCK C)!=\{b,c\},$

(2) $ \xbm (\{a\} \xCK C)!=\{b,c\},$

(3) $ \xbm (\{a' \} \xCK C)!=\{b,c' \}.$

Then by (3) $(a',b) \xCd (a',c' ) \xeb (a',c),$ by (2) $(a,b) \xCd
(a,c) \xeb (a,c' ),$ by (1)
$(a',c' )$ cannot be optimal, so neither is $(a',b),$ thus $(a,b)$ must
be better
than $(a',b),$ as one of $(a,b)$ and $(a',b)$ must be optimal.

Thus, in an indirect way, using elements $c$ and $c' $ not in $B,$ we may
be able to
find out which pairs in $A \xCK B$ are really optimal. Obviously,
arbitrarily long
chains might be involved, and such chains may also lead to contradictions.

$ \xcz $
\\[3ex]

\ee

So it seems difficult to find a general way to find the best pairs - and
much
will depend on the domain, how much we can compare, etc. It might also
well
be that we cannot find out - so we will branch into all possibilities, or
choose
arbitrarily one - loosing ignorance, of course.

\bd

$\hspace{0.01em}$

% (+++ Orig. No.:  Definition of the relation +++)

\label{Definition of the relation}

Let $\{X_{k}:k \xbe K\}$ be a $ \xbs -$cover of $X.$

Then we can define $ \xeb $ and $ \xec $ in the following cases:

(1) $ \xcA k( \xbm X_{k}! \xcC \xbm X!)$ $ \xcp $ $ \xbs \xeb \xbs ' $ for
all $ \xbs ' \xbe X,$ $ \xbs ' \xEd \xbs,$

(2) $ \xbm X! \xcC \xcV \xbm X_{k}!$ $ \xcp $ $ \xbs \xeb \xbs ' $ for all
$ \xbs ' \xbe X$ s.t. $ \xbs '! \xce \xbm (X)!,$ and
$ \xbs \xec \xbs ' $ for all $ \xbs ' \xbe X$ s.t. $ \xbs '! \xbe \xbm
X!.$

Explanation of the latter: If $ \xbs '! \xbe X!,$ there is some $ \xbs ''
$ s.t. $ \xbs '' $ is one of
the best, and $ \xbs ''!= \xbs ' $ - but we are not sure which one. So we
really need
$ \xec $ here.

Of course, we also know that for each $x \xbe \xbm X!$ there is some $
\xbs $ s.t. $ \xbs!=x$ and
$ \xbs \xec \xbs ' $ for all $ \xbs ' \xbe X,$ etc., but we do not know
which $ \xbs.$

\ed

We describe now the two main problems.

We want to construct a representing relation. In particular,

(a) if $ \xbt \xbe X,$ $ \xbt! \xce \xbm X!,$ we will need some $ \xbs
\xbe X,$ $ \xbs \xeb \xbt,$

and

(b) if $ \xbt! \xbe \xbm X!,$ and $x' \xEd \xbt!,$ $x' \xbe \xbm X!,$ we
will need some $ \xbs \xbe X,$ $ \xbs \xec \xbt.$

Consider now some such candidate, and the smallest set containing them,
$X:=[ \xbt, \xbs ].$

Problem (1): $ \xbs \xeb \xbt $ might well be the case, but there is
$ \xbt ' \xbe ( \xbt, \xbs ),$ $ \xbt '!= \xbt!,$ and $ \xbt ' \xeb
\xbs $ - so we will not see this, as $ \xbm X!=\{ \xbt '!\}.$

Problem (2): $ \xbs \xeb \xbt $ might well be the case, but there is
$ \xbs ' \xbe ( \xbt, \xbs ),$ $ \xbs '!= \xbs!,$ $ \xbs ' \xeb \xbt $
- so we will not see this, as already
$ \xbm [ \xbs ', \xbt ]!= \xbs!.$

Problem (2) can be overcome for our purposes, by choosing a $H-$closest $
\xbs $
s.t. $ \xbs \xeb \xbt,$ and if (1) is not a problem, we can see now that
$ \xbs \xeb \xbt:$ for all
suitable $ \xbr $ we have $ \xbm [ \xbr, \xbt ]!= \xbt!,$ and only $
\xbm [ \xbt, \xbs ]!= \xbs!.$

Problem (1) can be overcome if we know that some $ \xbs $ minimizes all $
\xbt,$ e.g.
if $ \xbs! \xbe \xbm X!,$ $ \xbt! \xce \xbm X!.$ We choose then for some
$ \xbt $ some such $ \xbs,$ and
consider $[ \xbt, \xbs ].$ There might now already be some $ \xbs ' \xbe
( \xbt, \xbs ),$ $ \xbs '!= \xbs!,$
which minimizes $ \xbt $ and all elements of $[ \xbs ', \xbt ],$ so we
have to take the
$H-$closest one to avoid Problem (2), work with this one, and we will see
that
$ \xbs ' \xeb \xbt.$

Both cases illustrate the importance of choosing $H-$closest elements,
which is,
of course, possible by finiteness of the index set.

Problem (1) can, however, be unsolvable if we take the limit approach, as
we
may be forced to do (see discussion in Section \ref{Section 2.3.2} (page
\pageref{Section 2.3.2}) ),
as we will not necessarily
have any more such $ \xbs $ which minimize all $ \xbt $ - see Example
\ref{Example 6.3.4} (page \pageref{Example 6.3.4}), Case (4).

But in the other, i.e. non-problematic, cases, this is the approach to
take:

Choose for $ \xbt $ some $ \xbs $
s.t. $ \xbs \xec \xbr $ for all $ \xbr $ - if this exists - then consider
$[ \xbt, \xbs ],$ choose in
$[ \xbt, \xbs ]$ the $H-$closest (measured from $ \xbt )$ $ \xbr $ s.t. $
\xbr $ is smaller or equal to all
$ \xbr ' \xbe [ \xbr, \xbt ].$ This has to exist, as $ \xbs $ is a
candidate, and we work with a
finite product, so the $H-$closest such has to exist.

We will then use a Cover Property:

\paragraph{
Property 6.3.1
}

$\hspace{0.01em}$

% (+++ Orig.:  Property 6.3.1 +++)

\label{Section Property 6.3.1}

Let $ \xbt \xbe X,$ $ \xbt! \xce \xbm X!,$ $ \xbS_{x}:=\{ \xbs \xbe X:
\xbs!=x\}$ for some fixed $x \xbe \xbm X!.$
Let $ \xbS_{x} \xcc \xcV \xbS_{i}$ for $i \xbe I,$ and $X_{i} \xcd
\xbS_{i} \xcv \{ \xbs \},$ then:
There is $i_{0} \xbe I$ s.t. $ \xbt! \xce \xbm X_{i}!.$

(This will also hold in the limit case.)

This allows to find witnesses, if possible:

Let $ \xbt \xbe X,$ $ \xbt! \xce \xbm X!,$ then $ \xbt $ is minimized by
one $ \xbs $ with $ \xbs! \xbe \xbm X!.$ Take such
$ \xbs $ with $H-$closest to $ \xbt,$ this exists by Cover Property. Then
we see that
$ \xbt! \xce \xbm [ \xbs, \xbt ]!,$ and that $ \xbs $ does it, so we
define $ \xbs \xeb \xbt.$

We will also have something similar for 2 elements $ \xbs, \xbs ' \xbe
X,$ $ \xbs! \xEd \xbs '!,$
$ \xbs!, \xbs! \xbe \xbm X!:$ We will seek $ \xbs,$ $ \xbs ' $ which
have minimal $H-$distance among
all $ \xbt,$ $ \xbt ' $ with the same endpoints s.t. $ \xbs \xec \xbr,$
$ \xbs ' \xec \xbr $ for all $ \xbr \xbe [ \xbs, \xbs ' ].$
This must exist in the minimal variant, but, warning, we are not sure that
they are really the minimal ones (by the order) - see Example \ref{Example
6.3.4} (page \pageref{Example 6.3.4}), (1).

\paragraph{
Small sets and easy properties
}

$\hspace{0.01em}$

% (+++ Orig.:  Small sets and easy properties +++)

\label{Section Small sets and easy properties}

As a first guess, one might think that small sets always suffice to
define the relation. This is not true. First, if $ \xbs!= \xbs '!,$ then
$[ \xbs, \xbs ' ]$ will
give us no information.

The following (also negative) example is more instructive:

\be

$\hspace{0.01em}$

% (+++ Orig. No.:  Example D-6.3.2 +++)

\label{Example D-6.3.2}

Consider $X:=\{0,1,2\} \xCK \{0,1\},$ with the ordering:
$(0,0) \xeb (2.1) \xeb (1,0) \xeb (0,1) \xCq (1,1) \xCq (2,0).$ We want to
find out that $(0,0) \xeb (1,1)$
using the results $ \xbm Y!.$

Consider $X':=\{0,1\} \xCK \{0,1\},$ the smallest set containing both.
Any $(0,0)-$cover
of $X' $ has to contain some $X'' $ with $(1,0) \xbe X'',$ but then $
\xbm X''!=\{0\},$ so we cannot
see that $(0,0) \xeb (1,1).$

But consider the $(0,0)-$cover $\{X',X'' \}$ of $X$ with $X':=\{1,2\}
\xCK \{0,1\},$ $X'':=\{(0,1)\},$
then $ \xbm X'!= \xbm X''!=\{1\},$ but $ \xbm X!=\{0\},$ so we see that
$(0,0) \xeb (1,1).$

$ \xcz $
\\[3ex]

\ee

(But we can obtain the same result through a chain of small sets:
$(0,0) \xeb (2,1):$ look at the $(0,0)-$cover $\{\{2\} \xCK
\{0,1\},\{(0,1)\}\}$ of $\{0,2\} \xCK \{0,1\},$
$(2,1) \xeb (2,0)$ is obvious,
$(2,0) \xec (1,1):$ look at the $(2,0)-$cover $\{\{1\} \xCK
\{0,1\},\{(2,1)\}\}$ of $\{1,2\} \xCK \{0,1\},$ we
see that $ \xbm (\{1\} \xCK \{0,1\})!= \xbm \{(2,1)\}!=\{1\},$ and $ \xbm
(\{1,2\} \xCK \{0,1\})!=\{0,1\}.)$

\br

$\hspace{0.01em}$

% (+++ Orig. No.:  Remark Immediate +++)

\label{Remark Immediate}

The following are immediate:

(1) $ \xcS \{ \xbP \{X^{j}_{i}:i \xbe I\}:j \xbe J\}$ $=$ $ \xbP \{ \xcS
\{X^{j}_{i}:j \xbe J\}:i \xbe I\}$

(2) $ \xbm X! \xcc X!$

(3) $X= \xcV X_{k}$ $ \xcp $ $ \xbm X! \xcc \xcV ( \xbm X_{k}!)$

(4) $X \xcc X' $ $ \xcp $ $X!- \xbm X! \xcc X'!- \xbm X'!$ is in general
wrong, though
$X \xcc X' $ $ \xcp $ $X- \xbm X \xcc X' - \xbm X' $ is correct: There
might be
a new minimal thread in $X',$ which happens to have the same endpoint
as discarded threads in $X.$
\subsubsection{A first serious problem}
\label{Section 2.3.2.3}

\er

We describe now the first serious problem, and indicate how to work around
it.
So, we think this can still be avoided in a reasonable way, it is the next
one,
where we meet essentially the same difficulty, and where we see no
solution
or detour.

We will complete the ranked relation (or, better, what we see of it) as
usual,
i.e. close under transitivity,
it has to be free of cycles, etc. Usually, the relation will not be fully
determined, so, forgetting about ignorance, we complete it using our
abstract nonsense lemma Fact \ref{Fact D-4.7.2} (page \pageref{Fact D-4.7.2}).
But there is a
problem, as the
following sketch of an example shows. We have to take care that we do not
have
infinite descending chains of threads $ \xbs $ with the same endpoint, as
the
endpoint will then not appear anymore, by definition of minimality. More
precisely:

\be

$\hspace{0.01em}$

% (+++ Orig. No.:  Example D-6.3.3 +++)

\label{Example D-6.3.3}

Let $ \xbs!= \xbt!$ and $ \xbs '!= \xbt '!,$ and $ \xbs ' \xeb \xbs,$
$ \xbt \xeb \xbt ',$ then choosing $ \xbs \xec \xbt $ will
result in $ \xbs ' \xeb \xbt '.$ Likewise, $ \xbs '' \xec \xbt '' $ with
$ \xbs ''!= \xbt ''!$ may lead to $ \xbr ' \xeb \xbs ' $
(with $ \xbs '!= \xbr '!),$ etc, and to an infinite descending chain.
This problem
seems difficult to solve, and the authors do not know if there is
always a solution, probably not, as we need just $ \xbo $ many pairs to
create
an infinite descending chain, but the universe can be arbitrarily big.

\ee

\paragraph{
Two ways to work around this problem:
}

$\hspace{0.01em}$

% (+++ Orig.:  Two ways to work around this problem: +++)

\label{Section Two ways to work around this problem:}

(1)
We may neglect such infinite descending chains, and do as if they were not
there. Thus, we will take a limit approach (see Section \ref{Section Limit}
(page \pageref{Section Limit}) ) locally,
within the concerned $ \xbs!.$ This might be justified intuitively by the
fact
that we are really interested only in one dimension, and do not really
care about how we got there, so we neglect the other dimensions somewhat.

(2)
We may consider immediately the global limit approach. In this case, it
will
probably be the best strategy to consider formula defined model sets, in
order
to be able to go easily to the logical version - see
Section \ref{Section Limit} (page \pageref{Section Limit})

For (1):

We continue to write $ \xbm X!$ for the observed result, though this will
not
be a minimum any more.

We then have:

$ \xcA \xbs \xbe X \xcA x \xbe \xbm X! \xcE \xbt ( \xbt!=x \xcu \xbt \xec
\xbs ),$

but we will not have any more:

$ \xcA x \xbe \xbm X! \xcA \xbs \xbe X \xcE \xbt \xbe X( \xbt!=x \xcu
\xcA \xbs \xbe X. \xbt \xec \xbs )$

Note that in the limit case:

(1) $X \xEd \xCQ $ $ \xcp $ $ \xbm X! \xEd \xCQ $

(2) If $x,x' \xbe \xbm X!,$ and $x$ is a limit (i.e. there is an infitite
descending chain
of $ \xbs ' $s with $ \xbs =x),$ then so is $x'.$

(3) By finiteness of the Hamming distance, e.g. there will be for all $
\xbs \xbe X$
and all $x \xbe \xbm X!$ cofinally often some $H-$closest (to $ \xbs )$ $
\xbt \xbe X$ s.t. $ \xbt!=x$ and
$ \xbt \xec \xbs.$ (Still, this will not help us much, unfortunately.)
\subsubsection{A second serious problem}
\label{Section 2.3.2.4}

This problem can be seen in a (class of) unpleasant $example(s),$ it has
to do
with the fact that interesting elements may be hidden by other elements.

\be

$\hspace{0.01em}$

% (+++ Orig. No.:  Example 6.3.4 +++)

\label{Example 6.3.4}

We discuss various variants of this example, and first present the common
parts.

The domain $U$ will be $2 \xCK \xbo \xCK 2.$ $X \xcc U$ etc. will be legal
products. We define
suitable orders, and will see that we can hardly get any information about
the order.

$ \xbs_{i}:=(0,i,0),$ $ \xbt_{i}:=(1,i,1).$

The key property is:

If $ \xbs_{i}, \xbt_{j} \xbe X \xcc U$ for some $i,j \xbe \xbo,$ then $
\xbs_{k} \xbe X$ $ \xcr $ $ \xbt_{k} \xbe X.$
(Proof: Let e.g. $ \xbs_{k}=(0,k,0) \xbe X,$ then, as $(1,j,1) \xbe X,$
$(1,k,1) \xbe X.)$

In all variants, we make a top layer consisting of all $(1,i,0)$ and
$(0,i,1),$
$i \xbe \xbo.$ They will play a minor role. All other elements will be
strictly
below this top layer.

We turn to the variants.

(1)
Let $ \xbs_{i} \xCq \xbt_{i}$ for all $i,$ and $ \xbs_{i} \xCq \xbt_{i}
\xeb \xbs_{i+1} \xCq \xbt_{i+1},$ close under transitivity.
Thus, there is a minimal pair, $ \xbs_{0},$ $ \xbt_{0},$ and we have an
``ascending ladder'' of
the other $ \xbs_{i},$ $ \xbt_{i}.$
But we cannot see which pair is at the bottom.

Proof: Let $X \xcc U.$ If $X$ contains only elements of the type $(a,b,0)$
or $(a,b,1),$
the result is trivial: $ \xbm X!=\{0\},$ etc. So suppose not. If $X$
contains no $ \xbs_{i}$ and
no $ \xbt_{i},$ the result is trivial again, and gives no information
about our
question. The same is true, if $X$ contains only $ \xbs_{i}' $s or only $
\xbt_{i}' s,$ but not
both. If $X$ contains some $ \xbs_{i},$ and some $ \xbt_{j},$ let $k$ be
the smallest such index
$i$ or $j,$ then, by above remark, $ \xbs_{k}$ and $ \xbt_{k}$ will be in
$X,$ so $ \xbm X!=\{0,1\}$ in any
case, and this will not tell us anything. (More formally, consider in
parallel
the modified order, where the pair $ \xbs_{1},$ $ \xbt_{1}$ is the
smallest one, instead of
$ \xbs_{0},$ $ \xbt_{0},$ and order the others as before, then this
different order will give
the same result in all cases.)

(2)
Order the elements as in (1), only make a descending ladder, instead of an
ascending one. Thus, $ \xbm U!= \xCQ,$ but $ \xbm X!$ in variant (1) and
in variant (2) will
not be different for any finite $X.$

Proof: As above, the interesting case in when there is some $ \xbs_{i}$
and some $ \xbt_{j}$ in
$X,$ we now take the biggest such index $k,$ then $ \xbs_{k}$ and $
\xbt_{k}$ will both be in $X,$
and thus $ \xbm X!=\{0,1\}.$ Consequently, only infinite $X$ allow us to
distinguish
variant (1) from variant (2). But they give no information how the
relation
is precisely defined, only that there are infinite descending chains.

Thus, finite $X$ do not allow us to distinguish between absence and
presence
of infinite descending chains.

(3)
Order the $ \xbs_{i}$ in an ascending chain, let $ \xbs_{i} \xeb
\xbt_{i},$ and no other comparisons
between the $ \xbs_{i}$ and the $ \xbt_{j}.$ Thus, this is not a ranked
relation. Again, the
only interesting case is where some $ \xbs_{i} \xbe X,$ and some $
\xbt_{j} \xbe X.$ But by above
reasoning, we have always $ \xbm X!=\{0\}.$

Modifying the order such that $ \xbt_{i} \xeb \xbs_{i+1}$ (making one
joint ascending chain)
is a ranked order, with the same results $ \xbm X!.$

Thus, we cannot distinguish ranked from non-ranked relations.

(4)
Let now $ \xbs_{i} \xeb \xbt_{i},$ and $ \xbt_{i+1} \xeb \xbs_{i},$ close
under transitivity. As we have an
infinite descending chain, we take the limit approach.

Let $X$ be finite, and some $ \xbs_{i}=(0,i,0) \xbe X,$ then $ \xbm
X!=\{0\}.$ Thus, we can never
see that some specific $ \xbt_{j}$ minimizes some $ \xbs_{k}.$

Proof: All elements of the top layer are minimized by $ \xbs_{i}.$ If
there is no
$ \xbt_{j} \xbe X,$ we are done. Let $ \xbt_{k}=(1,k,1)$ be the smallest $
\xbt_{j} \xbe X,$ then by the above
$ \xbs_{k}=(0,k,0) \xbe X,$ too, so $1 \xce \xbm X!.$

But we can see the converse: Consider $X:=[ \xbs_{i}, \xbt_{i}],$ then $
\xbm X!=\{0\},$ as
we just saw. The only other $ \xbs \xbe X$ with $ \xbs!=0$ is $(1,i,0),$
but we see that
$ \xbt_{i} \xeb \xbs,$ as $ \xbm [ \xbt_{i}, \xbs ]=\{1\}$ and $[
\xbt_{i}, \xbs ]$ contains only 2 elements, so $ \xbs_{i}$
must minimize $ \xbt_{i}.$

Consequently, we have the information $ \xbs_{i} \xeb \xbt_{i},$ but not
any information about
any $ \xbt_{i} \xeb \xbs_{j}.$ Yet, taking the limit approach applied to
$U,$ we see that there
must be below each $ \xbs_{j}$ some $ \xbt_{i},$ otherwise we would see
only $\{0\}$ as result.
But we do not know how to choose the $ \xbt_{i}$ below the $ \xbs_{j}.$

$ \xcz $
\\[3ex]
\subsubsection{Resume}
\label{Section 2.3.2.5}

\ee

The problem is that we might not have enough information to construct the
order. Taking any completion is not sure to work, as it might indirectly
contradict existing limits or non-limits, which give only very scant
information (there is or is not an infinite descending chain), but do not
tell us anything about the order.

For such situations, cardinalities of the sets involved might play a role.

Unfortunately, the authors have no idea how to attack these problems. It
seems
quite certain that the existing techniques in the domain will not help us.

So this might be a reason (or pretext?) to look at simplifications.

It might be intuitively justifiable to impose a continuity condition:

(Cont) If $ \xbt $ is between $ \xbs $ and $ \xbs ' $ in the Hamming
distance, then so is its
ranking, i.e. $ \xbr ( \xbs ) \xck \xbr ( \xbs ' )$ $ \xcp $ $ \xbr ( \xbs
) \xck \xbr ( \xbt ) \xck \xbr ( \xbs ' ),$ if $ \xbr $
is the ranking function.

This new condition seems to solve above problems, and can perhaps be
justified
intuitively. But we have to be aware of the temptation to hide our
ignorance
and inability to solve problems behind flimsy intuitions. Yet, has someone
ever come up with an intuitive justification of definability preservation?
(The
desastrous consequences of its absence were discussed in  \cite{Sch04}.)

On the other hand, this continuity or interpolation property might be
useful
in other circumstances, too, where we can suppose that small changes have
small
effects, see the second author's  \cite{Sch95-2} or Chapter 4 in
 \cite{Sch97-2} for a broader and deeper discussion.

This additional property should probably be further investigated.

A word of warning: this condition is NOT compatibel with distance
based reasoning: The distance from 0 to 2 is the same as the distance from
1 to 3, but the set $[(0,2),(1,3)]$ contains $(1,2),$ with a strictly
smaller
distance. A first idea is then not to work on states, but on differences
between states, considering the sum. But this does not work either. The
sequences $(0,1),$ $(1,0)$ contain between them the smaller one $(0,0),$
and the
bigger one $(1,1).$ Thus, the condition can be countertintuitive, and the
authors
do not know if there are sufficiently many natural scenarios where it
holds.

It seems, however, that we do not need the full condition, but the
following
would suffice: If $ \xbs $ is such that there is $ \xbt $ with $ \xbs!
\xEd \xbt!,$ and $ \xbt \xeb \xbs,$ then
we find $ \xbt ' $ s.t. $ \xbt!= \xbt '!$ and $ \xbt ' $ is smaller than
all sequences in $[ \xbs, \xbt ]$ (and
perhaps a similar condition for $ \xec ).$ Then, we would see smaller
elements in
Example \ref{Example 6.3.4} (page \pageref{Example 6.3.4}). This is a limit
condition, and is similar
to the condition
that there are cofinally many definable initial segments in the limit
approach,
see the discussion there - or in  \cite{Sch04}.
\vspace{.3cm}
\subsubsection{A comment on former work}
\label{Section 2.3.2.6}

We make here some very short remarks on our joint article with S.Berger
and
D.Lehmann,  \cite{BLS99}.

The perhaps central definition of the article is Definition 3.10 in 
\cite{BLS99}.

First, we see that the relation $R$ is generated only by sets where one
includes
the other. Thus, if the domain does not contain such sets, the relation
will
be void, and trivial examples show that completeness may collaps in such
cases.

More subtly, case 3. (a) is surprising, as $ \xcA i([B^{i}_{1} \Xl
B^{i}_{n-1}.C] \xcC [A_{1} \Xl.A_{n-1}.C])$ is
expected, and not the (in result) much weaker condition given.

The proof shows that the very definition of a patch allows to put the
interesting sequence in all elements of the patch (attention: elements of
the patch are not mutually disjoint, they are only disjoint from the
starting
set), so we can work with the intersection. The proof uses, however, that
singletons are in the domain, and that the elements of arbitrary patches
are
in the domain. In particular, this will generally not be true in the
infinite
case, as patches work with set differences.

Thus, this proof is another good illustration that we may pay in domain
conditions what we win in representation conditions.
\section{Theory revision}
\label{Section 2.5}
\label{Chapter TR}

We begin with a very succinct introduction into AGM theory revision,
and the subsequent results by the second author as published
in  \cite{Sch04}. It is not supposed to be a self-contained
introduction,
but to help the reader recollect the situation.

Section \ref{Section Booth} (page \pageref{Section Booth})  describes very
briefly parts of the work
by Booth and co-authors, and then solves a representation problem
in the infinite case left open by Booth et al.
\subsection{Introduction to theory revision}
\label{Section 2.5.1}

Recall from the introduction that theory revision was invented in order
to ``fuse'' together two separately consistent, but together inconsistent
theories or formulas to a consistent result.
The by far best known approach is that by Alchourron, Gardenfors, and
Makinson, and know as the AGM approach, see  \cite{AGM85}.
They formulated ``rationality postulates'' for
various variants of theory revision, which we give now in a very
succinct form.
Lehmann, Magidor, Schlechta, see  \cite{LMS01},
gave a distance semantics for theory revision, this is further
elaborated in  \cite{Sch04}, and presented here in very brief
outline, too.
\index{Definition AGM}

\bd

$\hspace{0.01em}$

% (+++ Orig. No.:  Definition AGM +++)

\label{Definition AGM}

We present in parallel the logical
and the semantic (or purely algebraic) side. For the latter, we work in
some
fixed universe $U,$ and the intuition is $U=M_{ \xdl },$ $X=M(K),$ etc.,
so, e.g. $A \xbe K$
becomes $X \xcc B,$ etc.

(For reasons of readability, we omit most caveats about definability.)

$K_{ \xcT }$ will denote the inconsistent theory.

We consider two functions, - and $*,$ taking a deductively closed theory
and a
formula as arguments, and returning a (deductively closed) theory on the
logics
side. The algebraic counterparts work on definable model sets. It is
obvious
that $ \xCf (K-1),$ $(K*1),$ $ \xCf (K-6),$ $(K*6)$ have vacuously true
counterparts on the
semantical side. Note that $K$ $ \xCf (X)$ will never change, everything
is relative
to fixed $K$ $ \xCf (X).$ $K* \xbf $ is the result of revising $K$ with $
\xbf.$ $K- \xbf $ is the result of
subtracting enough from $K$ to be able to add $ \xCN \xbf $ in a
reasonable way, called
contraction.

Moreover,
let $ \xck_{K}$ be a relation on the formulas relative to a deductively
closed theory $K$
on the formulas of $ \xdl,$ and $ \xck_{X}$ a relation on $ \xdp (U)$ or
a suitable subset of $ \xdp (U)$
relative to fixed $X.$ When the context is clear, we simply write $ \xck
.$
$ \xck_{K}$ $( \xck_{X})$ is called a relation of epistemic entrenchment
for $K$ $ \xCf (X).$

The following table presents the ``rationality postulates'' for contraction
(-),
revision $(*)$ and epistemic entrenchment. In AGM tradition, $K$ will be a
deductively closed theory, $ \xbf, \xbq $ formulas. Accordingly, $X$ will
be the set of
models of a theory, $A,B$ the model sets of formulas.

In the further development, formulas $ \xbf $ etc. may sometimes also be
full
theories. As the transcription to this case is evident, we will not go
into
details.

\renewcommand{\arraystretch}{1.2}

{\small
% \begin{tabular*}{15cm}{|c@{\extracolsep\fill}|c|c|c|}
\begin{tabular}{|c|c|c|c|}

\hline

\multicolumn{4}{|c|} {Contraction, $K-\xbf $} \xEP

\hline

$(K-1)$ \xEH $K-\xbf $ is deductively closed \xEH \xEH \xEP

\hline

$(K-2)$ \xEH $K-\xbf $ $ \xcc $ $K$ \xEH $(X \xDN 2)$ \xEH $X \xcc X \xDN A$
\xEP

\hline

$(K-3)$ \xEH $\xbf  \xce K$ $ \xch $ $K-\xbf =K$ \xEH $(X \xDN 3)$ \xEH $X \xcC
A$ $
\xch $ $X \xDN A=X$ \xEP

\hline

$(K-4)$ \xEH $ \xcL \xbf $ $ \xch $ $\xbf  \xce K-\xbf $ \xEH $(X \xDN 4)$ \xEH
$A \xEd
U$ $ \xch $ $X \xDN A \xcC A$ \xEP

\hline

$(K-5)$ \xEH $K \xcc \ol{(K-\xbf ) \xcv \{\xbf \}}$ \xEH $(X \xDN 5)$ \xEH $(X
\xDN
A) \xcs A$ $ \xcc $ $X$ \xEP

\hline

$(K-6)$ \xEH $ \xcl \xbf  \xcr \xbq $ $ \xch $ $K-\xbf =K-\xbq $ \xEH \xEH \xEP

\hline

$(K-7)$ \xEH $(K-\xbf ) \xcs (K-\xbq )  \xcc  $ \xEH
$(X \xDN 7)$ \xEH $X \xDN (A \xcs B)  \xcc  $ \xEP

\xEH $K-(\xbf  \xcu \xbq ) $ \xEH
\xEH $(X \xDN A) \xcv (X \xDN B)$ \xEP

\hline

$(K-8)$ \xEH $\xbf  \xce K-(\xbf  \xcu \xbq )  \xch  $ \xEH
$(X \xDN 8)$ \xEH $X \xDN (A \xcs B) \xcC A  \xch  $ \xEP

\xEH $K-(\xbf  \xcu \xbq ) \xcc K-\xbf $ \xEH
\xEH $X \xDN A \xcc X \xDN (A \xcs B)$ \xEP

\hline
\hline

\multicolumn{4}{|c|} {Revision, $K*\xbf $} \xEP

\hline

$(K*1)$ \xEH $K*\xbf $ is deductively closed \xEH - \xEH \xEP

\hline

$(K*2)$ \xEH $\xbf  \xbe K*\xbf $ \xEH $(X \xfA 2)$ \xEH $X \xfA A \xcc A$
\xEP

\hline

$(K*3)$ \xEH $K*\xbf $ $ \xcc $ $ \ol{K \xcv \{\xbf \}}$ \xEH $(X \xfA 3)$ \xEH
$X
\xcs A \xcc X \xfA A$ \xEP

\hline

$(K*4)$ \xEH $ \xCN \xbf  \xce K  \xch $ \xEH
$(X \xfA 4)$ \xEH $X \xcs A \xEd \xCQ   \xch $ \xEP

\xEH $\ol{K \xcv \{\xbf \}}  \xcc  K*\xbf $ \xEH
\xEH $X \xfA A \xcc X \xcs A$ \xEP

\hline

$(K*5)$ \xEH $K*\xbf =K_{ \xcT }$ $ \xch $ $ \xcl \xCN \xbf $ \xEH $(X \xfA 5)$
\xEH $X \xfA A= \xCQ $ $ \xch $ $A= \xCQ $ \xEP

\hline

$(K*6)$ \xEH $ \xcl \xbf  \xcr \xbq $ $ \xch $ $K*\xbf =K*\xbq $ \xEH - \xEH
\xEP

\hline

$(K*7)$ \xEH $K*(\xbf  \xcu \xbq )  \xcc $ \xEH
$(X \xfA 7)$ \xEH $(X \xfA A) \xcs B  \xcc  $ \xEP

\xEH $\ol{(K*\xbf ) \xcv \{\xbq \}}$ \xEH
\xEH $X \xfA (A \xcs B)$ \xEP

\hline

$(K*8)$ \xEH $ \xCN \xbq  \xce K*\xbf  \xch $ \xEH
$(X \xfA 8)$ \xEH $(X \xfA A) \xcs B \xEd \xCQ \xch $ \xEP

\xEH $\ol{(K*\xbf ) \xcv \{\xbq \}} \xcc K*(\xbf  \xcu \xbq )$ \xEH
\xEH $ X \xfA (A \xcs B) \xcc (X \xfA A) \xcs B$ \xEP

\hline
\hline

\multicolumn{4}{|c|} {Epistemic entrenchment} \xEP

\hline

$(EE1)$ \xEH $ \xck_{K}$ is transitive \xEH
$(EE1)$ \xEH $ \xck_{X}$ is transitive \xEP

\hline

$(EE2)$ \xEH $\xbf  \xcl \xbq   \xch  \xbf  \xck_{K}\xbq $ \xEH
$(EE2)$ \xEH $A \xcc B  \xch  A \xck_{X}B$ \xEP

\hline

$(EE3)$ \xEH $ \xcA  \xbf,\xbq  $ \xEH
$(EE3)$ \xEH $ \xcA A,B $ \xEP
\xEH $ (\xbf  \xck_{K}\xbf  \xcu \xbq $ or $\xbq  \xck_{K}\xbf  \xcu \xbq )$
\xEH
\xEH $ (A \xck_{X}A \xcs B$ or $B \xck_{X}A \xcs B)$ \xEP

\hline

$(EE4)$ \xEH $K \xEd K_{ \xcT }  \xch  $ \xEH
$(EE4)$ \xEH $X \xEd \xCQ   \xch  $ \xEP
\xEH $(\xbf  \xce K$ iff $ \xcA  \xbq.\xbf  \xck_{K}\xbq )$ \xEH
\xEH $(X \xcC A$ iff $ \xcA B.A \xck_{X}B)$ \xEP

\hline

$(EE5)$ \xEH $ \xcA \xbq.\xbq  \xck_{K}\xbf   \xch   \xcl \xbf $ \xEH
$(EE5)$ \xEH $ \xcA B.B \xck_{X}A  \xch  A=U$ \xEP

\hline

\end{tabular}
}
\\

\index{Remark TR-Rank}

\ed

\br

$\hspace{0.01em}$

% (+++ Orig. No.:  Remark TR-Rank +++)

\label{Remark TR-Rank}

(1) Note that $(X \xfA 7)$ and $(X \xfA 8)$ express a central condition
for ranked
structures, see Section 3.10: If we note $X \xfA.$ by $f_{X}(.),$ we then
have:
$f_{X}(A) \xcs B \xEd \xCQ $ $ \xch $ $f_{X}(A \xcs B)=f_{X}(A) \xcs B.$

(2) It is trivial to see that AGM revision cannot be defined by an
individual
distance (see Definition 2.3.5 below):
Suppose $X \xfA Y$ $:=$ $\{y \xbe Y:$ $ \xcE x_{y} \xbe X( \xcA y' \xbe
Y.d(x_{y},y) \xck d(x_{y},y' ))\}.$
Consider $a,b,c.$ $\{a,b\} \xfA \{b,c\}=\{b\}$ by $(X \xfA 3)$ and $(X
\xfA 4),$ so $d(a,b)<d(a,c).$
But on the other hand $\{a,c\} \xfA \{b,c\}=\{c\},$ so $d(a,b)>d(a,c),$
$contradiction.$
\index{Proposition AGM-Equiv}

\er

\bp

$\hspace{0.01em}$

% (+++ Orig. No.:  Proposition AGM-Equiv +++)

\label{Proposition AGM-Equiv}

Contraction, revision, and epistemic entrenchment are interdefinable by
the
following equations, i.e., if the defining side has the respective
properties,
so will the defined side.

\renewcommand{\arraystretch}{1.5}

{\scriptsize
% \begin{tabular*}{15.0cm}{|c|c|}
\begin{tabular}{|c|c|}

\hline

$K*\xbf:= \ol{(K- \xCN \xbf )} \xcv {\xbf }$  \xEH  $X \xfA A:= (X \xDN  \xdC A)
\xcs A$ \xEP

\hline

$K-\xbf:= K \xcs (K* \xCN \xbf )$  \xEH  $X \xDN A:= X \xcv (X \xfA  \xdC A)$
\xEP

\hline

$K-\xbf:=\{\xbq  \xbe K:$ $(\xbf <_{K}\xbf  \xco \xbq $ or $ \xcl \xbf )\}$ \xEH

$
X \xDN A:=
\left\{
\begin{array}{rcl}
X & iff & A=U, \\
 \xcS \{B: X \xcc B \xcc U, A<_{X}A \xcv B\} & & otherwise \\
\end{array}
\right. $
\xEP

\hline

$
\xbf  \xck_{K}\xbq: \xcr
\left\{
\begin{array}{l}
\xcl \xbf  \xcu \xbq  \\
or \\
\xbf  \xce K-(\xbf  \xcu \xbq ) \\
\end{array}
\right. $
\xEH

$
A \xck_{X}B: \xcr
\left\{
\begin{array}{l}
A,B=U  \\
or \\
X \xDN (A \xcs B) \xcC A \\
\end{array}
\right. $
\xEP

\hline

\end{tabular}
}

\index{Intuit-Entrench}
\paragraph{A remark on intuition}

\ep

The idea of epistemic entrenchment is that $ \xbf $ is more entrenched
than $ \xbq $
(relative to $K)$ iff $M( \xCN \xbq )$ is closer to $M(K)$ than $M( \xCN
\xbf )$ is to $M(K).$ In
shorthand, the more we can twiggle $K$ without reaching $ \xCN \xbf,$ the
more $ \xbf $ is
entrenched. Truth is maximally entrenched - no twiggling whatever will
reach
falsity. The more $ \xbf $ is entrenched,
the more we are certain about it. Seen this way, the properties of
epistemic
entrenchment relations are very natural (and trivial): As only the closest
points of $M( \xCN \xbf )$ count (seen from $M(K)),$ $ \xbf $ or $ \xbq $
will be as entrenched as
$ \xbf \xcu \xbq,$ and there is a logically strongest $ \xbf ' $ which is
as entrenched as $ \xbf $ -
this is just the sphere around $M(K)$ with radius $d(M(K),M( \xCN \xbf
)).$
\index{Definition Distance}

\bd

$\hspace{0.01em}$

% (+++ Orig. No.:  Definition Distance +++)

\label{Definition Distance}

$d:U \xCK U \xcp Z$ is called a pseudo-distance on $U$ iff (d1) holds:

(d1) $Z$ is totally ordered by a relation $<.$

If, in addition, $Z$ has a $<-$smallest element 0, and (d2) holds, we say
that $d$
respects identity:

(d2) $d(a,b)=0$ iff $a=b.$

If, in addition, (d3) holds, then $d$ is called symmetric:

(d3) $d(a,b)=d(b,a).$

(For any $a,b \xbe U.)$

Note that we can force the triangle inequality to hold trivially (if we
can
choose the values in the real numbers): It suffices to choose the values
in
the set $\{0\} \xcv [0.5,1],$ i.e. in the interval from 0.5 to 1, or as 0.
\index{Definition Dist-Indiv-Coll}

\ed

\bd

$\hspace{0.01em}$

% (+++ Orig. No.:  Definition Dist-Indiv-Coll +++)

\label{Definition Dist-Indiv-Coll}

We define the collective and the individual variant of choosing the
closest
elements in the second operand by two operators,
$ \xfA, \xfB: \xdp (U) \xCK \xdp (U) \xcp \xdp (U):$

Let $d$ be a distance or pseudo-distance.

$X \xfA Y$ $:=$ $\{y \xbe Y:$ $ \xcE x_{y} \xbe X. \xcA x' \xbe X, \xcA y'
\xbe Y(d(x_{y},y) \xck d(x',y' )\}$

(the collective variant, used in theory revision)

and

$X \xfB Y$ $:=$ $\{y \xbe Y:$ $ \xcE x_{y} \xbe X. \xcA y' \xbe
Y(d(x_{y},y) \xck d(x_{y},y' )\}$

(the individual variant, used for counterfactual conditionals and theory
update).

Thus, $A \xfA_{d}B$ is the subset of $B$ consisting of all $b \xbe B$ that
are closest to A.
Note that, if $ \xCf A$ or $B$ is infinite, $A \xfA_{d}B$ may be empty,
even if $ \xCf A$ and $B$ are not
empty. A condition assuring nonemptiness will be imposed when necessary.
\index{Definition Dist-Repr}

\ed

\bd

$\hspace{0.01em}$

% (+++ Orig. No.:  Definition Dist-Repr +++)

\label{Definition Dist-Repr}

An operation $ \xfA: \xdp (U) \xCK \xdp (U) \xcp \xdp (U)$ is
representable iff there is a
pseudo-distance $d:U \xCK U \xcp Z$ such that

$A \xfA B$ $=$ $A \xfA_{d}B$ $:=$ $\{b \xbe B:$ $ \xcE a_{b} \xbe A \xcA
a' \xbe A \xcA b' \xbe B(d(a_{b},b) \xck d(a',b' ))\}.$
\index{Definition TR*d}

\ed

The following is the central definition, it describes the way a revision
$*_{d}$ is
attached to a pseudo-distance $d$ on the set of models.

\bd

$\hspace{0.01em}$

% (+++ Orig. No.:  Definition TR*d +++)

\label{Definition TR*d}

$T*_{d}T' $ $:=$ $Th(M(T) \xfA_{d}M(T' )).$

$*$ is called representable iff there is a pseudo-distance $d$ on the set
of models
s.t. $T*T' =Th(M(T) \xfA_{d}M(T' )).$
\index{Fact AGM-In-Dist}

\ed

\bfa

$\hspace{0.01em}$

% (+++ Orig. No.:  Fact AGM-In-Dist +++)

\label{Fact AGM-In-Dist}

A distance based revision satisfies the AGM postulates provided:

(1) it respects identity, i.e. $d(a,a)<d(a,b)$ for all $a \xEd b,$

(2) it satisfies a limit condition: minima exist,

(3) it is definability preserving.

(It is trivial to see that the first two are necessary,
and Example \ref{Example TR-Dp} (page \pageref{Example TR-Dp})  (2)
below shows the necessity of (3). In particular, (2) and (3) will hold for
finite languages.)
\index{Fact AGM-In-Dist Proof}

\efa

\subparagraph{
Proof
}

$\hspace{0.01em}$

% (+++ Orig.:  Proof +++)

We use $ \xfA $ to abbreviate $ \xfA_{d}.$ As a matter of fact, we show
slightly more, as
we admit also full theories on the right of $*.$

$(K*1),$ $(K*2),$ $(K*6)$ hold by definition, $(K*3)$ and $(K*4)$ as $d$
respects
identity, $(K*5)$ by existence of minima.

It remains to show $(K*7)$ and $(K*8),$ we do them together, and show:
If $T*T' $ is consistent with $T'',$ then $T*(T' \xcv T'' )$ $=$ $
\ol{(T*T' ) \xcv T'' }.$

Note that $M(S \xcv S' )=M(S) \xcs M(S' ),$ and that $M(S*S' )=M(S) \xfA
M(S' ).$ (The latter
is only true if $ \xfA $ is definability preserving.)
By prerequisite, $M(T*T' ) \xcs M(T'' ) \xEd \xCQ,$ so $(M(T) \xfA M(T'
)) \xcs M(T'' ) \xEd \xCQ.$
Let $A:=M(T),$ $B:=M(T' ),$ $C:=M(T'' ).$ `` $ \xcc $ '': Let $b \xbe A
\xfA (B \xcs C).$
By prerequisite, there is $b' \xbe (A \xfA B) \xcs C.$ Thus $d(A,b' ) \xcg
d(A,B \xcs C)=d(A,b).$
As $b \xbe B,$ $b \xbe A \xfA B,$ but $b \xbe C,$ too. `` $ \xcd $ '': Let
$b' \xbe (A \xfA B) \xcs C.$ Thus $d(A,b' )=$
$d(A,B) \xck d(A,B \xcs C),$ so by $b' \xbe B \xcs C$ $b' \xbe A \xfA (B
\xcs C).$
We conclude $M(T) \xfA (M(T' ) \xcs M(T'' ))$ $=$ $(M(T) \xfA M(T' )) \xcs
M(T'' ),$ thus that
$T*(T' \xcv T'' )= \ol{(T*T' ) \xcv T'' }.$

$ \xcz $
\\[3ex]
\index{Definition TR-Umgeb}

\bd

$\hspace{0.01em}$

% (+++ Orig. No.:  Definition TR-Umgeb +++)

\label{Definition TR-Umgeb}

For $X,Y \xEd \xCQ,$ set $U_{Y}(X):=\{z:d(X,z) \xck d(X,Y)\}.$
\index{Fact TR-Umgeb}

\ed

\bfa

$\hspace{0.01em}$

% (+++ Orig. No.:  Fact Tr-Umgeb +++)

\label{Fact Tr-Umgeb}

Let $X,Y,Z \xEd \xCQ.$ Then

(1) $U_{Y}(X) \xcs Z \xEd \xCQ $ iff $(X \xfA (Y \xcv Z)) \xcs Z \xEd \xCQ
,$

(2) $U_{Y}(X) \xcs Z \xEd \xCQ $ iff $ \xdC Z \xck_{X} \xdC Y$ - where $
\xck_{X}$ is epistemic entrenchement relative
to $X.$
\index{Fact TR-Umgeb Proof}

\efa

\subparagraph{
Proof
}

$\hspace{0.01em}$

% (+++ Orig.:  Proof +++)

(1) Trivial.

(2) $ \xdC Z \xck_{X} \xdC Y$ iff $X \xDN ( \xdC Z \xcs \xdC Y) \xcC \xdC
Z.$ $X \xDN ( \xdC Z \xcs \xdC Y)$ $=$ $X \xcv (X \xfA \xdC ( \xdC Z \xcs
\xdC Y))$ $=$
$X \xcv (X \xfA (Z \xcv Y)).$ So $X \xDN ( \xdC Z \xcs \xdC Y) \xcC \xdC
Z$ $ \xcj $ $(X \xcv (X \xfA (Z \xcv Y))) \xcs Z \xEd \xCQ $ $ \xcj $
$X \xcs Z \xEd \xCQ $ or $(X \xfA (Z \xcv Y)) \xcs Z \xEd \xCQ $ $ \xcj $
$d(X,Z) \xck d(X,Y).$

$ \xcz $
\\[3ex]
\index{Definition TR-Dist}

\bd

$\hspace{0.01em}$

% (+++ Orig. No.:  Definition TR-Dist +++)

\label{Definition TR-Dist}

Let $U \xEd \xCQ,$ $ \xdy \xcc \xdp (U)$ satisfy $( \xcs ),$ $( \xcv ),$
$ \xCQ \xce \xdy.$

Let $A,B,X_{i} \xbe \xdy,$ $ \xfA: \xdy \xCK \xdy \xcp \xdp (U).$

Let $*$ be a revision function defined for
arbitrary consistent theories on both sides. (This is thus a slight
extension of
the AGM framework, as AGM work with formulas only on the right of $*.)$

{\scriptsize

% \begin{tabular*}{15.0cm}{|c|c|c|}
\begin{tabular}{|c|c|c|}

\hline

\xEH
\xEH
$(*Equiv)$
\xEP
\xEH
\xEH
$ \xcm T \xcr S,$ $ \xcm T' \xcr S',$ $\xch$ $T*T' =S*S',$
\xEP

\hline

\xEH
\xEH
$(*CCL)$
\xEP
\xEH
\xEH
$T*T' $ is a consistent, deductively closed theory,
\xEP

\hline

\xEH
$( \xfA Succ)$
\xEH
$(*Succ)$
\xEP
\xEH
$A \xfA B \xcc B$
\xEH
$T' \xcc T*T',$
\xEP

\hline

\xEH
$( \xfA Con)$
\xEH
$(*Con)$
\xEP
\xEH
$A \xcs B \xEd \xCQ $ $ \xch $ $A \xfA B=A \xcs B$
\xEH
$Con(T \xcv T') $ $\xch$ $T*T' = \ol{T \xcv T' },$
\xEP

\hline

Intuitively,
\xEH
$( \xfA Loop)$
\xEH
$(*Loop)$
\xEP
Using symmetry
\xEH
\xEH
\xEP
$d(X_{0},X_{1}) \xck d(X_{1},X_{2}),$
\xEH
$(X_{1} \xfA (X_{0} \xcv X_{2})) \xcs X_{0} \xEd \xCQ,$
\xEH
$Con(T_{0},T_{1}*(T_{0} \xco T_{2})),$
\xEP
$d(X_{1},X_{2}) \xck d(X_{2},X_{3}),$
\xEH
$(X_{2} \xfA (X_{1} \xcv X_{3})) \xcs X_{1} \xEd \xCQ,$
\xEH
$Con(T_{1},T_{2}*(T_{1} \xco T_{3})),$
\xEP
$d(X_{2},X_{3}) \xck d(X_{3},X_{4})$
\xEH
$(X_{3} \xfA (X_{2} \xcv X_{4})) \xcs X_{2} \xEd \xCQ,$
\xEH
$Con(T_{2},T_{3}*(T_{2} \xco T_{4}))$
\xEP
\Xl
\xEH
\Xl
\xEH
\Xl
\xEP
$d(X_{k-1},X_{k}) \xck d(X_{0},X_{k})$
\xEH
$(X_{k} \xfA (X_{k-1} \xcv X_{0})) \xcs X_{k-1} \xEd \xCQ $
\xEH
$Con(T_{k-1},T_{k}*(T_{k-1} \xco T_{0}))$
\xEP
$\xch$
\xEH
$\xch$
\xEH
$\xch$
\xEP
$d(X_{0},X_{1}) \xck d(X_{0},X_{k}),$
\xEH
$(X_{0} \xfA (X_{k} \xcv X_{1})) \xcs X_{1} \xEd \xCQ$
\xEH
$Con(T_{1},T_{0}*(T_{k} \xco T_{1}))$
\xEP

i.e. transitivity, or absence of
\xEH
\xEH
\xEP

loops involving $<$
\xEH
\xEH
\xEP

\hline

\end{tabular}

}

\index{Proposition TR-Alg-Log}

\ed

\bp

$\hspace{0.01em}$

% (+++ Orig. No.:  Proposition TR-Alg-Log +++)

\label{Proposition TR-Alg-Log}

The following connections between the logical and the algebraic side might
be the most interesting ones. We will consider in all cases also the
variant
with full theories.

Given $*$ which respects logical equivalence, let $M(T) \xfA M(T'
):=M(T*T' ),$
conversely, given $ \xfA,$ let $T*T':=Th(M(T) \xfA M(T' )).$ We then
have:

{\small

% \begin{tabular*}{13.0cm}{|c|c|c|c|}
\begin{tabular}{|c|c|c|c|}

\hline

(1.1)
\xEH
$(K*7)$
\xEH
$\xch$
\xEH
$(X \xfA 7)$
\xEP

\cline{1-1}
\cline{3-3}

(1.2)
\xEH
\xEH
$\xci$ $(\xbm dp)$
\xEH
\xEP

\cline{1-1}
\cline{3-3}

(1.3)
\xEH
\xEH
$\xci$ B is the model set for some $\xbf$
\xEH
\xEP

\cline{1-1}
\cline{3-3}

(1.4)
\xEH
\xEH
$\xcI$ in general
\xEH
\xEP

\hline

(2.1)
\xEH
$(*Loop)$
\xEH
$\xch$
\xEH
$(\xfA Loop)$
\xEP

\cline{1-1}
\cline{3-3}

(2.2)
\xEH
\xEH
$\xci$ $(\xbm dp)$
\xEH
\xEP

\cline{1-1}
\cline{3-3}

(2.3)
\xEH
\xEH
$\xci$ all $X_i$ are the model sets for some $\xbf_i$
\xEH
\xEP

\cline{1-1}
\cline{3-3}

(2.4)
\xEH
\xEH
$\xcI$ in general
\xEH
\xEP

\hline

\end{tabular}

}

\index{Proposition TR-Alg-Log Proof}

\ep

\subparagraph{
Proof
}

$\hspace{0.01em}$

% (+++ Orig.:  Proof +++)

(1)

We consider the equivalence of $T*(T' \xcv T'' ) \xcc \ol{(T*T' ) \xcv T''
}$ and
$(M(T) \xfA M(T' )) \xcs M(T'' ) \xcc M(T) \xfA (M(T' ) \xcs M(T'' )).$

(1.1)

$(M(T) \xfA M(T' )) \xcs M(T'' )$ $=$ $M(T*T' ) \xcs M(T'' )$ $=$ $M((T*T'
) \xcv T'' )$ $ \xcc_{(K*7)}$
$M(T*(T' \xcv T'' ))$ $=$ $M(T) \xfA M(T' \xcv T'' )$ $=$ $M(T) \xfA (M(T'
) \xcs M(T'' )).$

(1.2)

$T*(T' \xcv T'' )$ $=$ $Th(M(T) \xfA M(T' \xcv T'' ))$ $=$ $Th(M(T) \xfA
(M(T' ) \xcs M(T'' )))$ $ \xcc_{(X \xfA 7)}$
$Th((M(T) \xfA M(T' )) \xcs M(T'' ))$ $=_{( \xbm dp)}$
$ \ol{Th(M(T) \xfA M(T' )) \xcv T'' }$ $=$ $ \ol{Th(M(T*T' ) \xcv T'' }$
$=$ $ \ol{(T*T' ) \xcv T'' }.$

(1.3)

Let $T'' $ be equivalent to $ \xbf ''.$ We can then replace the use of $(
\xbm dp)$
in the proof of (1.2) by Fact \ref{Fact Log-Form} (page \pageref{Fact Log-Form})
 (3).

(1.4)

By Example \ref{Example TR-Dp} (page \pageref{Example TR-Dp})  (2), $(K*7)$ may
fail, though $(X \xfA
7)$ holds.

(2.1) and (2.2):

$Con(T_{0},T_{1}*(T_{0} \xco T_{2}))$ $ \xcj $ $M(T_{0}) \xcs
M(T_{1}*(T_{0} \xco T_{2})) \xEd \xCQ.$

$M(T_{1}*(T_{0} \xco T_{2}))$ $=$ $M(Th(M(T_{1}) \xfA M(T_{0} \xco
T_{2})))$ $=$ $M(Th(M(T_{1}) \xfA (M(T_{0}) \xcv M(T_{2}))))$ $=_{( \xbm
dp)}$
$M(T_{1}) \xfA (M(T_{0}) \xcv (T_{2})),$ so
$Con(T_{0},T_{1}*(T_{0} \xco T_{2}))$ $ \xcj $ $M(T_{0}) \xcs (M(T_{1})
\xfA (M(T_{0}) \xcv (T_{2}))) \xEd \xCQ.$

Thus, all conditions translate one-to-one, and we use $( \xfA Loop)$ and
$(*Loop)$
to go back and forth.

(2.3):

Let $A:=M(Th(M(T_{1}) \xfA (M(T_{0}) \xcv M(T_{2})))),$ $A':=M(T_{1})
\xfA (M(T_{0}) \xcv (T_{2})),$ then we do
not need $A=A',$ it suffices to have $M(T_{0}) \xcs A \xEd \xCQ $ $ \xcj
$ $M(T_{0}) \xcs A' \xEd \xCQ.$ $A= \wt{A' },$ so
we can use Fact \ref{Fact Log-Form} (page \pageref{Fact Log-Form})  (4), if
$T_{0}$ is equivalent to
some $ \xbf_{0}.$

This has to hold for all $T_{i},$ so all $T_{i}$ have to be equivalent to
some $ \xbf_{i}.$

(2.4):

By Proposition \ref{Proposition TR-Alg-Repr} (page \pageref{Proposition
TR-Alg-Repr}), all distance defined $ \xfA $
satisfy
$( \xfA Loop).$ By Example \ref{Example TR-Dp} (page \pageref{Example TR-Dp}) 
(1), $(*Loop)$ may fail.

$ \xcz $
\\[3ex]
\index{Proposition TR-Representation-With-Ref}

The following table summarizes representation of theory revision
functions by structures with a distance.

By ``pseudo-distance'' we mean here a pseudo-distance which respects
identity, and is symmetrical.

$( \xfA \xCQ )$ means that if $X,Y \xEd \xCQ,$ then $X \xfA_{d}Y \xEd
\xCQ.$
\label{Proposition TR-Representation-With-Ref}

{\scriptsize

% \begin{tabular*}{15.9cm}{|c|c|c|c|c|}
\begin{tabular}{|c|c|c|c|c|}

\hline

$\xfA-$ function
\xEH
\xEH
Distance Structure
\xEH
\xEH
$*-$ function
\xEP

\hline

$(\xfA Succ)+(\xfA Con)+$
\xEH
$\xcj$, $(\xcv)+(\xcs)$
\xEH
pseudo-distance
\xEH
$\xcj$ $(\xbm dp)+(\xfA\xCQ)$
\xEH
$(*Equiv)+(*CCL)+(*Succ)+$
\xEP

$(\xfA Loop)$
\xEH
Proposition \ref{Proposition TR-Alg-Repr}
\xEH
\xEH
Proposition \ref{Proposition TR-Log-Repr}
\xEH
$(*Con)+(*Loop)$
\xEP

\xEH
page \pageref{Proposition TR-Alg-Repr}
\xEH
\xEH
page \pageref{Proposition TR-Log-Repr}
\xEH
\xEP

\cline{1-2}
\cline{4-4}

any finite
\xEH
$\xcJ$
\xEH
\xEH
$\xcH$ without $(\xbm dp)$
\xEH
\xEP

characterization
\xEH
Proposition \ref{Proposition Hamster}
\xEH
\xEH
Example \ref{Example TR-Dp}
\xEH
\xEP

\xEH
page \pageref{Proposition Hamster}
\xEH
\xEH
page \pageref{Example TR-Dp}
\xEH
\xEP

\hline

\end{tabular}

}

\index{Example TR-Dp}

The following Example \ref{Example TR-Dp} (page \pageref{Example TR-Dp})  shows
that, in general, a
revision operation
defined on models via a pseudo-distance by $T*T':=Th(M(T) \xfA_{d}M(T'
))$ might not
satisfy $(*Loop)$ or $(K*7),$ unless we require $ \xfA_{d}$ to preserve
definability.

\be

$\hspace{0.01em}$

% (+++ Orig. No.:  Example TR-Dp +++)

\label{Example TR-Dp}

Consider an infinite propositional language $ \xdl.$

Let $X$ be an infinite set of models, $m,$ $m_{1},$ $m_{2}$ be models for
$ \xdl.$
Arrange the models of $ \xdl $ in the real plane s.t. all $x \xbe X$ have
the same
distance $<2$ (in the real plane) from $m,$ $m_{2}$ has distance 2 from
$m,$ and $m_{1}$ has
distance 3 from $m.$

Let $T,$ $T_{1},$ $T_{2}$ be complete (consistent) theories, $T' $ a
theory with infinitely
many models, $M(T)=\{m\},$ $M(T_{1})=\{m_{1}\},$ $M(T_{2})=\{m_{2}\}.$ The
two variants diverge now
slightly:

(1) $M(T' )=X \xcv \{m_{1}\}.$ $T,T',T_{2}$ will be pairwise
inconsistent.

(2) $M(T' )=X \xcv \{m_{1},m_{2}\},$ $M(T'' )=\{m_{1},m_{2}\}.$

Assume in both cases $Th(X)=T',$ so $X$ will not be definable by a
theory.

Now for the results:

Then $M(T) \xfA M(T' )=X,$ but $T*T' =Th(X)=T'.$

(1) We easily verify
$Con(T,T_{2}*(T \xco T)),$ $Con(T_{2},T*(T_{2} \xco T_{1})),$
$Con(T,T_{1}*(T \xco T)),$ $Con(T_{1},T*(T_{1} \xco T' )),$
$Con(T,T' *(T \xco T)),$ and conclude by Loop (i.e. $(*Loop))$
$Con(T_{2},T*(T' \xco T_{2})),$ which
is wrong.

(2) So $T*T' $ is consistent with $T'',$ and
$ \ol{(T*T' ) \xcv T'' }=T''.$ But $T' \xcv T'' =T'',$ and $T*(T' \xcv
T'' )=T_{2} \xEd T'',$ contradicting $(K*7).$

$ \xcz $
\\[3ex]
\index{Proposition TR-Alg-Repr}

\ee

\bp

$\hspace{0.01em}$

% (+++ Orig. No.:  Proposition TR-Alg-Repr +++)

\label{Proposition TR-Alg-Repr}

Let $U \xEd \xCQ,$ $ \xdy \xcc \xdp (U)$ be closed under finite $ \xcs $
and finite $ \xcv,$ $ \xCQ \xce \xdy.$

(a) $ \xfA $ is representable by a symmetric pseudo-distance $d:U \xCK U
\xcp Z$ iff $ \xfA $
satisfies $( \xfA Succ)$ and $( \xfA Loop)$ in Definition \ref{Definition
TR-Dist} (page \pageref{Definition TR-Dist}).

(b) $ \xfA $ is representable by an identity respecting symmetric
pseudo-distance
$d:U \xCK U \xcp Z$ iff $ \xfA $ satisfies $( \xfA Succ),$ $( \xfA Con),$
and $( \xfA Loop)$
in Definition \ref{Definition TR-Dist} (page \pageref{Definition TR-Dist}).

See  \cite{LMS01} or  \cite{Sch04}.
\index{Proposition TR-Log-Repr}

\ep

\bp

$\hspace{0.01em}$

% (+++ Orig. No.:  Proposition TR-Log-Repr +++)

\label{Proposition TR-Log-Repr}

Let $ \xdl $ be a propositional language.

(a) A revision operation $*$ is representable by a symmetric consistency
and
definability preserving pseudo-distance iff $*$ satisfies $(*Equiv),$
$(*CCL),$ $(*Succ),$ $(*Loop).$

(b) A revision operation $*$ is representable by a symmetric consistency
and
definability preserving, identity respecting pseudo-distance iff $*$
satisfies $(*Equiv),$ $(*CCL),$ $(*Succ),$ $(*Con),$ $(*Loop).$

See  \cite{LMS01} or  \cite{Sch04}.
\index{Example WeakTR}

\ep

\be

$\hspace{0.01em}$

% (+++ Orig. No.:  Example WeakTR +++)

\label{Example WeakTR}

This example shows the expressive weakness of revision based on
distance: not all distance relations can be reconstructed from
the revision operator. Thus, a revision operator does not allow
to ``observe'' all distances relations, so transitivity of $ \xck $ cannot
necessarily be captured in a short condition, requiring
arbitrarily long conditions, see Proposition \ref{Proposition Hamster} (page
\pageref{Proposition Hamster}).

Note that even when the pseudo-distance is a real distance, the
resulting revision operator $ \xfA_{d}$ does not always permit to
reconstruct the
relations of the distances: revision is a coarse instrument to investigate
distances.

Distances with common start (or end, by symmetry) can always be
compared by looking at the result of revision:

$a \xfA_{d}\{b,b' \}=b$ iff $d(a,b)<d(a,b' ),$

$a \xfA_{d}\{b,b' \}=b' $ iff $d(a,b)>d(a,b' ),$

$a \xfA_{d}\{b,b' \}=\{b,b' \}$ iff $d(a,b)=d(a,b' ).$

This is not the case with arbitrary distances $d(x,y)$ and $d(a,b),$
as this example will show.

\ee

We work in the real plane, with the standard distance, the angles have 120
degrees. $a' $ is closer to $y$ than $x$ is to $y,$ a is closer to $b$
than $x$ is to $y,$
but $a' $ is farther away from $b' $ than $x$ is from $y.$ Similarly for
$b,b'.$
But we cannot distinguish the situation $\{a,b,x,y\}$ and the
situation $\{a',b',x,y\}$ through $ \xfA_{d}.$ (See Diagram \ref{Diagram WeakTR}
(page \pageref{Diagram WeakTR}) ):

Seen from a, the distances are in that order: $y,b,x.$

Seen from $a',$ the distances are in that order: $y,b',x.$

Seen from $b,$ the distances are in that order: $y,a,x.$

Seen from $b',$ the distances are in that order: $y,a',x.$

Seen from $y,$ the distances are in that order: $a/b,x.$

Seen from $y,$ the distances are in that order: $a' /b',x.$

Seen from $x,$ the distances are in that order: $y,a/b.$

Seen from $x,$ the distances are in that order: $y,a' /b'.$

Thus, any $c \xfA_{d}C$ will be the same in both situations (with a
interchanged with
$a',$ $b$ with $b' ).$ The same holds for any $X \xfA_{d}C$ where $X$ has
two elements.

Thus, any $C \xfA_{d}D$ will be the same in both situations, when we
interchange a with
$a',$ and $b$ with $b'.$ So we cannot determine by $ \xfA_{d}$ whether
$d(x,y)>d(a,b)$ or not.
$ \xcz $
\\[3ex]

\vspace{10mm}

\begin{diagram}

\label{Diagram WeakTR}
\index{Diagram WeakTR}

\centering
\setlength{\unitlength}{1mm}
{\renewcommand{\dashlinestretch}{30}
\begin{picture}(130,90)(0,0)

\put(5,50){\line(1,0){30}}
\put(35,50){\line(1,2){6}}
\put(35,50){\line(1,-2){6}}

\put(5,50){\circle*{1.5}}
\put(35,50){\circle*{1.5}}

\put(41,62){\circle*{1.5}}
\put(41,38){\circle*{1.5}}

\put(5,47){$x$}
\put(32,47){$y$}
\put(43,61){$a$}
\put(43,37){$b$}

\put(65,50){\line(1,0){35}}
\put(100,50){\line(1,2){12}}
\put(100,50){\line(1,-2){12}}

\put(65,50){\circle*{1.5}}
\put(100,50){\circle*{1.5}}

\put(112,74){\circle*{1.5}}
\put(112,26){\circle*{1.5}}

\put(65,47){$x$}
\put(97,47){$y$}
\put(114,73){$a'$}
\put(114,25){$b'$}

\put(29,10){Indiscernible by revision}

\end{picture}

}

\end{diagram}

\vspace{4mm}

\index{Proposition Hamster}

\bp

$\hspace{0.01em}$

% (+++ Orig. No.:  Proposition Hamster +++)

\label{Proposition Hamster}

There is no finite characterization of distance based $ \xfA -$operators.

(Attention: this is, of course, false when we fix the left hand side:
the AGM axioms give a finite characterization. So this also shows the
strength of being able to change the left hand side.)

See  \cite{Sch04}.
\subsection{
Booth revision
}
\label{Section Booth}
\subsubsection{
Introduction
}

\ep

This material is due to Booth and co-authors.
\paragraph{
The problem we solve
}

Booth and his co-authors have shown in very interesting papers, see
 \cite{BN06} and  \cite{BCMG06}, that
many new approaches to theory revision (with fixed $K)$ can be represented
by two
relations, $<$ and $ \xej,$ where $<$ is the usual ranked relation, and $
\xej $ is a
sub-relation of $<.$ They have, however, left open a characterization of
the
infinite case, which we treat here.

The, for us, main definition they give is (in slight modification, we use
the
strict subrelations):

\bd

$\hspace{0.01em}$

% (+++ Orig. No.:  Definition Booth-K-phi +++)

\label{Definition Booth-K-phi}

Given $K,$ and $<$ and $ \xej,$ we define

$K \xDN \xbf:=Th(\{w:w \xej w' $ for some $w' \xbe min(M( \xCN \xbf
),<)\}),$

i.e. $K \xDN \xbf $ is given by all those worlds, which are below the
closest $ \xbf -$worlds,
as seen from $K.$

\ed

We want to characterize $K \xDN \xbf,$ for fixed $K.$ Booth et al. have
done the finite
case by working with complete consistent formulas, i.e. single models. We
want to do the infinite case without using complete consistent theories,
i.e.
in the usual style of completeness results in the area.

Our approach is basically semantic, though we use
sometimes the language of logic, on the one hand to show how to
approximate
with formulas a single model, and on the other hand when we use classical
compactness. This is, however, just a matter of speaking, and we could
translate it into model sets, too, but we do not think that we would win
much
by doing so. Moreoever, we will treat only the formula case, as this seems
to be the most interesting (otherwise the problem of approximation by
formulas
would not exist), and restrict ourselves to the definability preserving
case.
The more general case is left open, for a young researcher who
wants to sharpen his tools by solving it. Another open problem is to treat
the same question for variable $K,$ for distance based revision.
\paragraph{
The framework
}

For the reader's convenience, and to put our work a bit more into
perspective, we repeat now some of the definitions and results given
by Booth and his co-authors.

Consequently, all material in this section is due to Booth and his
co-authors.

$ \xck $ will be a total pre-order, anchored on $M(K),$ the models of $K,$
i.e.
$M(K)=min(W, \xck ),$ the set of $ \xck -$minimal worlds.

We have a second binary relation $ \xec $ on $W,$ which is a reflexive
subrelation of $ \xck.$

\bd

$\hspace{0.01em}$

% (+++ Orig. No.:  Definition Booth-Context +++)

\label{Definition Booth-Context}

(1) $( \xck, \xec )$ is a $K-$context iff $ \xck $ is a total pre-order
on $W,$ anchored on $M(K),$
and $ \xec $ is a reflexive sub-relation of $ \xck.$

(2) $K \xDN \xbf:=Th(\{w:w \xec w' $ for some $w' \xbe min(M( \xCN \xbf
), \xck )\})$ is called a basic removal
operator.

\ed

\bt

$\hspace{0.01em}$

% (+++ Orig. No.:  Theorem Booth-Basic-Removal +++)

\label{Theorem Booth-Basic-Removal}

Basic removal is characterzed by:

$ \xCf (B1)$ $K \xDN \xbf =Cn(K \xDN \xbf )$ - $ \xCf Cn$ classical
consequence,

$ \xCf (B2)$ $ \xbf \xce K \xDN \xbf,$

$ \xCf (B3)$ If $ \xcm \xbf \xcr \xbf ',$ then $K \xDN \xbf =K \xDN \xbf
',$

$ \xCf (B4)$ $K \xDN \xcT =K,$

$ \xCf (B5)$ $K \xDN \xbf \xcc Cn(K \xcv \{ \xCN \xbf \}),$

$ \xCf (B6)$ if $ \xbs \xbe K \xDN ( \xbs \xcu \xbf ),$ then $ \xbs \xbe K
\xDN ( \xbs \xcu \xbf \xcu \xbq ),$

$ \xCf (B7)$ if $ \xbs \xbe K \xDN ( \xbs \xcu \xbf ),$ then $K \xDN \xbf
\xcc K \xDN ( \xbs \xcu \xbf ),$

$ \xCf (B8)$ $(K \xDN \xbs ) \xcs (K \xDN \xbf ) \xcc K \xDN ( \xbs \xcu
\xbf ),$

$ \xCf (B9)$ if $ \xbf \xce K \xDN ( \xbs \xcu \xbf ),$ then $K \xDN (
\xbs \xcu \xbf ) \xcc K \xDN \xbf.$

\et

$ \xCf (B1)- \xCf (B3)$ belong to the basic AGM contraction postulates, $
\xCf (B4)- \xCf (B5)$ are
weakened versions of another basic AGM postulate:

$ \xCf (Vacuity)$ If $ \xbf \xce K,$ then $K \xDN \xbf =K$

which does not necessarily hold for basic removal operators.

The same holds for the remaining two basic AGM contraction postulates:

$ \xCf (Inclusion)$ $K \xDN \xbf \xcc K$

$ \xCf (Recovery)$ $K \xcc Cn((K \xDN \xbf ) \xcv \{ \xbf \}).$

The main definition towards the completeness result of Booth et al. is:

\bd

$\hspace{0.01em}$

% (+++ Orig. No.:  Definition Booth-C(K,-) +++)

\label{Definition Booth-C}

Given $K$ and $ \xDN,$ the structure $C(K, \xDN )$ is defined by:

$( \xck )$ $w \xck w' $ iff $ \xCN \xba \xce K \xDN ( \xCN \xba \xcu \xCN
\xba ' )$ and

$( \xec )$ $w \xec w' $ iff $ \xCN \xba \xce K \xDN \xCN \xba ',$

where $ \xba $ is a formula which holds exactly in $w,$ analogously for
$w' $ and $ \xba '.$

\ed

Booth et al. then give a long list of Theorems showing equivalence between
various postulates, and conditions on the orderings $ \xck $ and $ \xec.$
This, of
course, shows the power of their approach.

We give three examples:

\bcd

$\hspace{0.01em}$

% (+++ Orig. No.:  Conditions Booth-div +++)

\label{Conditions Booth-div}

(c) If (for each $i=1,2)$ $w_{i} \xck w' $ for all $w',$ then $w_{1} \xec
w_{2}.$

(d) If $w_{1} \xck w_{2}$ for all $w_{2},$ then $w_{1} \xec w_{2}$ for all
$w_{2}.$

(e) If $w_{1} \xec w_{2},$ then $w_{1}=w_{2}$ or $w_{1} \xck w' $ for all
$w'.$

\ecd

\bt

$\hspace{0.01em}$

% (+++ Orig. No.:  Theorem Booth-Conditions-Postulates +++)

\label{Theorem Booth-Conditions-Postulates}

Let $ \xDN $ be a basic removal operator as defined above.

(1) $ \xDN $ satisfies one half of $ \xCf (Vacuity):$ If $ \xbf \xce K,$
then $K \xcc K \xDN \xbf,$

(2.1) If $( \xck, \xec )$ satisfies (c), then $ \xDN $ satisfies $ \xCf
(Vacuity).$

(2.2) If $ \xDN $ satisfies $ \xCf (Vacuity),$ then $C(K, \xDN )$ satsfies
(c).

(3.1) If $( \xck, \xec )$ satisfies (d), then $ \xDN $ satisfies $ \xCf
(Inclusion).$

(3.2) If $ \xDN $ satisfies $ \xCf (Inclusion),$ then $C(K, \xDN )$
satsfies (d).

(4.1) If $( \xck, \xec )$ satisfies (e), then $ \xDN $ satisfies $ \xCf
(Recovery).$

(4.2) If $ \xDN $ satisfies $ \xCf (Recovery),$ then $C(K, \xDN )$
satsfies (e).

(5) The following are equivalent:

(5.1) $ \xDN $ is a full AGM contraction operator,

(5.2) $ \xDN $ satisfies $ \xCf (B1)- \xCf (B9),$ $ \xCf (Inclusion),$ and
$ \xCf (Recovery)$

(5.3) $ \xDN $ is generated by some $( \xck, \xec )$ satisfying (d) and
(e).
\subsubsection{
Construction and proof
}
\label{Section Booth-Proof@}

\et

We change perspective a little, and work directly with a ranked relation,
so we
forget about the (fixed) $K$ of revision, and have an equivalent, ranked
structure. We are then interested in an operator $ \xbn,$ which returns a
model set
$ \xbn ( \xbf ):= \xbn (M( \xbf )),$ where $ \xbn ( \xbf ) \xcs M( \xbf )$
is given by a ranked relation $<,$ and
$ \xbn ( \xbf )-M( \xbf )$ $:=$ $\{x \xce M( \xbf ): \xcE y \xbe \xbn (
\xbf ) \xcs M( \xbf )(x \xej y)\},$ and $ \xej $ is an arbitrary
subrelation of $<.$ The essential problem is to find such $y,$ as we have
only
formulas to find it. (If we had full theories, we could just look at all
$Th(\{y\})$ whether $x \xbe \xbn (Th(\{y\})).)$ There is still some more
work to do, as we
have to connect the two relations, and simply taking a ready
representation
result will not do, as we shall see.

We first introduce some notation, then a set of conditions, and formulate
the
representation result. Soundness will be trivial. For completeness, we
construct first the ranked relation $<,$ show that it does what it should
do,
and then the subrelation $ \xej.$

\bn

$\hspace{0.01em}$

% (+++ Orig. No.:  Notation Booth-Nota-1 +++)

\label{Notation Booth-Nota-1}

We set

$ \xbm^{+}(X):= \xbn (X) \xcs X$

$ \xbm^{-}(X):= \xbn (X)-X$

where $X:=M( \xbf )$ for some $ \xbf.$

\en

\bcd

$\hspace{0.01em}$

% (+++ Orig. No.:  Conditions Booth-Cond-1 +++)

\label{Conditions Booth-Cond-1}

$( \xbm^{-}1)$ $Y \xcs \xbm^{-}(X) \xEd \xCQ $ $ \xcp $ $ \xbm^{+}(Y) \xcs
X= \xCQ $

$( \xbm^{-}2)$ $Y \xcs \xbm^{-}(X) \xEd \xCQ $ $ \xcp $ $ \xbm^{+}(X \xcv
Y)= \xbm^{+}(Y)$

$( \xbm^{-}3)$ $Y \xcs \xbm^{-}(X) \xEd \xCQ $ $ \xcp $ $ \xbm^{-}(Y) \xcs
X= \xCQ $

$( \xbm^{-}4)$ $ \xbm^{+}(A) \xcc \xbm^{+}(B)$ $ \xcp $ $ \xbm^{-}(A) \xcc
\xbm^{-}(B)$

$( \xbm^{-}5)$ $ \xbm^{+}(X \xcv Y)= \xbm^{+}(X) \xcv \xbm^{+}(Y)$ $ \xcp
$ $ \xbm^{-}(X \xcv Y)= \xbm^{-}(X) \xcv \xbm^{-}(Y)$

\ecd

\bfa

$\hspace{0.01em}$

% (+++ Orig. No.:  Fact Booth-Fact-1 +++)

\label{Fact Booth-Fact-1}

$( \xbm^{-}1)$ and $( \xbm \xCQ ),$ $( \xbm \xcc )$ for $ \xbm^{+}$ imply

(1) $ \xbm^{+}(X) \xcs Y \xEd \xCQ $ $ \xcp $ $ \xbm^{+}(X) \xcs
\xbm^{-}(Y)= \xCQ $

(2) $X \xcs \xbm^{-}(X)= \xCQ.$

\efa

\subparagraph{
Proof
}

$\hspace{0.01em}$

% (+++ Orig.:  Proof +++)

(1) Let $ \xbm^{+}(X) \xcs \xbm^{-}(Y) \xEd \xCQ,$ then $X \xcs
\xbm^{-}(Y) \xEd \xCQ,$ so by $( \xbm^{-}1)$ $ \xbm^{+}(X) \xcs Y= \xCQ
.$

(2) Set $X:=Y,$ and use $( \xbm \xCQ ),$ $( \xbm \xcc ),$ $( \xbm^{-}1),$
(1).

$ \xcz $
\\[3ex]

\bp

$\hspace{0.01em}$

% (+++ Orig. No.:  Proposition Booth-Prop-1 +++)

\label{Proposition Booth-Prop-1}

$ \xbn:\{M( \xbf ): \xbf \xbe F( \xdl )\} \xcp \xdD_{ \xdl }$ is
representable by $<$ and $ \xej,$ where $<$ is a
smooth ranked relation, and $ \xej $ a subrelation of $<,$ and $
\xbm^{+}(X)$ is the usual set
of $<-$minimal elements of $X,$ and $ \xbm^{-}(X)$ $=$ $\{x \xce X: \xcE y
\xbe \xbm^{+}(X).(x \xej y)\},$ iff
the following conditions hold:
$( \xbm \xcc ),$ $( \xbm \xCQ ),$ $( \xbm =)$ for $ \xbm^{+},$ and $(
\xbm^{-}1)-( \xbm^{-}5)$ for $ \xbm^{+}$ and $ \xbm^{-}.$

\ep

\subparagraph{
Proof
}

$\hspace{0.01em}$

% (+++ Orig.:  Proof +++)

\paragraph{
Soundness
}

$\hspace{0.01em}$

% (+++ Orig.:  Soundness +++)

\label{Section Soundness}

The first three hold for smooth ranked structures, and the others are
easily
verified.

\paragraph{
Completeness
}

$\hspace{0.01em}$

% (+++ Orig.:  Completeness +++)

\label{Section Completeness}

We first show how to generate the ranked relation $<:$
\label{Booth-rel-generation}

There is a small problem.

The authors first thought that one may take any result for ranked
structures off
the shelf, plug in the other relation somehow (see the second half), and
that's
it. No, that $isn' t$ it: Suppose there is $x,$ and a sequence $x_{i}$
converging to $x$
in the usual topology. Thus, if $x \xbe M( \xbf ),$ then there will always
be some $x_{i}$ in
$M( \xbf ),$ too. Take now a ranked structure $ \xdz,$ where all the
$x_{i}$ are strictly
smaller than $x.$ Consider $ \xbm ( \xbf ),$ this will usually not contain
$x$ (avoid some
nasty things with definability), so in the usual construction $( \xec_{1}$
below),
$x$ will not be forced to be below any element $y,$ how high up $y>x$
might be.
However, there is $ \xbq $ separating $x$ and $y,$ e.g. $x \xcm \xCN \xbq
,$ $y \xcm \xbq,$ and if we take as
the second relation just the ranking again, $x \xbe \xbm^{-}( \xbq ),$ so
this becomes visible.

Consequently, considering $ \xbm^{-}$ may give strictly more information,
and we have to
put in a little more work. We just patch a proof for simple ranked
structures,
adding information obtained through $ \xbm^{-}.$

We follow closely the strategy of the proof of 3.10.11 in  \cite{Sch04}. We
will,
however, change notation at one point: the relation $R$ in  \cite{Sch04} is
called $ \xec $
here. The proof goes over several steps, which we will enumerate.

Note that by Fact \ref{Fact Mu-Base} (page \pageref{Fact Mu-Base}), taken
from  \cite{Sch04}, see also  \cite{GS08c}, $( \xbm \xFO ),$
$( \xbm \xcv ),$
$( \xbm \xcv ' ),$ $( \xbm =' )$ hold, as the prerequisites about the
domain are valid.

(1) To generate the ranked relation $<,$
we define two relations, $ \xec_{1}$ and $ \xec_{2},$ where $ \xec_{1}$ is
the usual one for
ranked structures, as defined in the proof of 3.10.11 of  \cite{Sch04},
$a \xec_{1}b$ iff $a \xbe \xbm^{+}(X),$ $b \xbe X,$ or $a=b,$ and
$a \xec_{2}b$ iff $a \xbe \xbm^{-}(X),$ $b \xbe X.$

Moreover, we set $a \xec b$ iff $a \xec_{1}b$ or $a \xec_{2}b.$

(2) Obviously, $ \xec $ is reflexive, we show that $ \xec $ is transitive
by looking at the
four different cases.

(2.1) In  \cite{Sch04}, it was shown that $a \xec_{1}b \xec_{1}c$ $
\xcp $ $a \xec_{1}c.$ For completeness' sake,
we repeat the argument:
Suppose $a \xec_{1}b,$ $b \xec_{1}c,$ let $a \xbe \xbm^{+}(A),$ $b \xbe
A,$ $b \xbe \xbm^{+}(B),$ $c \xbe B.$ We show
$a \xbe \xbm^{+}(A \xcv B).$ By $( \xbm \xFO )$ $a \xbe \xbm^{+}(A \xcv
B)$ or $b \xbe \xbm^{+}(A \xcv B).$ Suppose $b \xbe \xbm^{+}(A \xcv B),$
then $ \xbm^{+}(A \xcv B) \xcs A \xEd \xCQ,$ so by $( \xbm =)$ $
\xbm^{+}(A \xcv B) \xcs A= \xbm^{+}(A),$ so $a \xbe \xbm^{+}(A \xcv B).$

(2.2) Suppose $a \xec_{1}b \xec_{2}c,$ we show $a \xec_{1}c:$ Let $c \xbe
Y,$ $b \xbe \xbm^{-}(Y) \xcs X,$ $a \xbe \xbm^{+}(X).$
Consider $X \xcv Y.$ As $X \xcs \xbm^{-}(Y) \xEd \xCQ,$ by $( \xbm^{-}2)$
$ \xbm^{+}(X \xcv Y)= \xbm^{+}(X),$ so $a \xbe \xbm^{+}(X \xcv Y)$ and
$c \xbe X \xcv Y,$ so $a \xec_{1}c.$

(2.3) Suppose $a \xec_{2}b \xec_{2}c,$ we show $a \xec_{2}c:$ Let $c \xbe
Y,$ $b \xbe \xbm^{-}(Y) \xcs X,$ $a \xbe \xbm^{-}(X).$
Consider $X \xcv Y.$ As $X \xcs \xbm^{-}(Y) \xEd \xCQ,$ by $( \xbm^{-}2)$
$ \xbm^{+}(X \xcv Y)= \xbm^{+}(X),$ so by $( \xbm^{-}5)$
$ \xbm^{-}(X \xcv Y)= \xbm^{-}(X),$ so $a \xbe \xbm^{-}(X \xcv Y)$ and $c
\xbe X \xcv Y,$ so $a \xec_{2}c.$

(2.4) Suppose $a \xec_{2}b \xec_{1}c,$ we show $a \xec_{2}c:$ Let $c \xbe
Y,$ $b \xbe \xbm^{+}(Y) \xcs X,$ $a \xbe \xbm^{-}(X).$
Consider $X \xcv Y.$ As $ \xbm^{+}(Y) \xcs X \xEd \xCQ,$ $ \xbm^{+}(X)
\xcc \xbm^{+}(X \xcv Y).$ (Here is the argument:
By $( \xbm \xFO ),$ $ \xbm^{+}(X \xcv Y)= \xbm^{+}(X) \xFO \xbm^{+}(Y),$
so, if $ \xbm^{+}(X) \xcC \xbm^{+}(X \xcv Y),$ then
$ \xbm^{+}(X) \xcs \xbm^{+}(X \xcv Y)= \xCQ,$ so $ \xbm^{+}(X) \xcs (X
\xcv Y- \xbm^{+}(X \xcv Y)) \xEd \xCQ $ by $( \xbm \xCQ ),$ so by
$( \xbm \xcv ' )$ $ \xbm^{+}(X \xcv Y)= \xbm^{+}(Y).$ But if $ \xbm^{+}(Y)
\xcs X= \xbm^{+}(X \xcv Y) \xcs X \xEd \xCQ,$
$ \xbm^{+}(X)= \xbm^{+}(X \xcv Y) \xcs X$ by $( \xbm =),$ so $ \xbm^{+}(X)
\xcs \xbm^{+}(X \xcv Y) \xEd \xCQ,$ $contradiction.)$
So $ \xbm^{-}(X) \xcc \xbm^{-}(X \xcv Y)$ by $( \xbm^{-}4),$ so $c \xbe X
\xcv Y,$ $a \xbe \xbm^{-}(X \xcv Y),$ and $a \xec_{2}c.$

(3) We also see:

(3.1) $a \xbe \xbm^{+}(A),$ $b \xbe A- \xbm^{+}(A)$ $ \xcp $ $b \xeC a.$

(3.2) $a \xbe \xbm^{-}(A),$ $b \xbe A$ $ \xcp $ $b \xeC a.$

Proof of (3.1):

(a) $ \xCN (b \xec_{1}a)$ was shown in  \cite{Sch04}, we repeat
again the argument:
Suppose there is $B$ s.t. $b \xbe \xbm^{+}(B),$ $a \xbe B.$ Then by $(
\xbm \xcv )$ $ \xbm^{+}(A \xcv B) \xcs B= \xCQ,$
and by $( \xbm \xcv ' )$ $ \xbm^{+}(A \xcv B)= \xbm^{+}(A),$ but $a \xbe
\xbm^{+}(A) \xcs B,$ $contradiction.$

(b) Suppose there is $B$ s.t. $a \xbe B,$ $b \xbe \xbm^{-}(B).$ But $A
\xcs \xbm^{-}(B) \xEd \xCQ $ implies
$ \xbm^{+}(A) \xcs B= \xCQ $ by $( \xbm^{-}1).$

Proof of (3.2):

(a) Suppose $b \xec_{1}a,$ so there is $B$ s.t. $a \xbe B,$ $b \xbe
\xbm^{+}(B),$ so $B \xcs \xbm^{-}(A) \xEd \xCQ,$ so
$ \xbm^{+}(B) \xcs A= \xCQ $ by $( \xbm^{-}1).$
- -
(b) Suppose $b \xec_{2}a,$ so there is $B$ s.t. $a \xbe B,$ $b \xbe \xbm $
(B), so $B \xcs \xbm $ $(A) \xEd \xCQ,$ so
$ \xbm^{-}(B) \xcs A= \xCQ $ by $( \xbm^{-}3).$

(4) Let, by Fact \ref{Fact Abs-Rel-Ext} (page \pageref{Fact Abs-Rel-Ext}),
$S$ be a total, transitive, reflexive relation on
$U$ which extends $ \xec $ s.t. $ \xCf xSy,ySx$ $ \xcp $ $x \xec y$
(recall that $ \xec $ is transitive and
reflexive). But note that we loose ignorance, here.
Define $a<b$ iff $ \xCf aSb,$ but not $ \xCf bSa.$
If $a \xcT b$ (i.e. neither $a<b$ nor $b<a),$ then, by totality of $S,$ $
\xCf aSb$ and $ \xCf bSa.$
$<$ is ranked: If $c<a \xcT b,$ then by transitivity of $S$ $ \xCf cSb,$
but
if $ \xCf bSc,$ then again by transitivity of $S$ $ \xCf aSc.$ Similarly
for $c>a \xcT b.$

(5) It remains to show that $<$ represents $ \xbm $ and is $ \xdy
-$smooth:

Let $a \xbe A- \xbm^{+}(A).$ By $( \xbm \xCQ ),$ $ \xcE b \xbe
\xbm^{+}(A),$ so $b \xec_{1}a,$
but by case (3.1) above $a \xeC b,$ so bSa, but not aSb, so $b<a,$ so $a
\xbe A- \xbm_{<}(A).$
Let $a \xbe \xbm^{+}(A),$ then for all $a' \xbe A$ $a \xec a',$ so aSa',
so there is
no $a' \xbe A$ $a' <a,$ so $a \xbe \xbm_{<}(A).$
Finally, $ \xbm^{+}(A) \xEd \xCQ,$ all $x \xbe \xbm^{+}(A)$ are minimal
in A as we just saw, and for
$a \xbe A- \xbm^{+}(A)$ there is $b \xbe \xbm^{+}(A),$ $b \xec_{1}a,$ so
the structure is smooth.

The subrelation $ \xej $:
\label{Booth-subrelation}

Let $x \xbe \xbm^{-}(X),$ we look for $y \xbe \xbm^{+}(X)$ s.t. $x \xej y$
where $ \xej $ is the smaller,
additional relation. By the definition of the relation $ \xec_{2}$ above,
we know that
$ \xej \xcc \xec $ and by (3.2) above $ \xej \xcc <.$

Take an arbitrary enumeration of the propositional variables of $ \xdl,$
$p_{i}:i< \xbk.$
We will inductively decide for $p_{i}$ or $ \xCN p_{i}.$ $ \xbs $ etc.
will denote a finite
subsequence of the choices made so far, i.e. $ \xbs = \xCL p_{i_{0}}, \Xl
, \xCL p_{i_{n}}$ for some $n< \xbo.$
Given such $ \xbs,$ $M( \xbs ):=M( \xCL p_{i_{0}}) \xcs  \Xl  \xcs M(
\xCL p_{i_{n}}).$ $ \xbs + \xbs ' $ will be the union of two
such sequences, this is again one such sequence.

Take an arbitrary model $m$ for $ \xdl,$ i.e. a function $m:v( \xdl )
\xcp \{t,f\}.$ We will use
this model as a ``strategy'', which will tell us how to decide, if we have
some
choice.

We determine $y$ by an inductive process, essentially cutting away $
\xbm^{+}(X)$ around
$y.$ We choose $p_{i}$ or $ \xCN p_{i}$ preserving the following
conditions inductively:
For all finite sequences $ \xbs $ as above we have:

(1) $M( \xbs ) \xcs \xbm^{+}(X) \xEd \xCQ,$

(2) $x \xbe \xbm^{-}(X \xcs M( \xbs )).$

For didactic reasons, we do the case $p_{0}$ separately.

Consider $p_{0}.$ Either $M(p_{0}) \xcs \xbm^{+}(X) \xEd \xCQ,$ or $M(
\xCN p_{0}) \xcs \xbm^{+}(X) \xEd \xCQ,$ or both.
If e.g. $M(p_{0}) \xcs \xbm^{+}(X) \xEd \xCQ,$ but $M( \xCN p_{0}) \xcs
\xbm^{+}(X)= \xCQ,$ then we have no choice, and
we take $p_{0},$ in the opposite case, we take $ \xCN p_{0}.$ E.g. in the
first case,
$ \xbm^{+}(X \xcs M(p_{0}))= \xbm^{+}(X),$ so $x \xbe \xbm^{-}(X \xcs
M(p_{0}))$ by $( \xbm^{-}4).$
If both intersections are non-empty, then by $( \xbm^{-}5)$ $x \xbe
\xbm^{-}(X \xcs M(p_{0}))$ or
$x \xbe \xbm^{-}(X \xcs M( \xCN p_{0})),$ or both. Only in the last case,
we use our strategy to
decide whether to choose $p_{0}$ or $ \xCN p_{0}:$ if $m(p_{0})=t,$ we
choose $p_{0},$ if not, we choose
$ \xCN p_{0}.$

Obviously, (1) and (2) above are satisfied.

Suppose we have chosen $p_{i}$ or $ \xCN p_{i}$ for all $i< \xba,$ i.e.
defined a partial function
from $v( \xdl )$ to $\{t,f\},$ and the induction hypotheses (1) and (2)
hold.
Consider $p_{ \xba }.$ If there is no finite subsequence $ \xbs $ of the
choices done so far
s.t. $M( \xbs ) \xcs M(p_{ \xba }) \xcs \xbm^{+}(X)= \xCQ,$ then $p_{
\xba }$ is a candidate. Likewise for $ \xCN p_{ \xba }.$

One of $p_{ \xba }$ or $ \xCN p_{ \xba }$ is a candidate:

Suppose not, then there are $ \xbs $ and $ \xbs ' $ subsequences of the
choices
done so far, and $M( \xbs ) \xcs M(p_{ \xba }) \xcs \xbm^{+}(X)= \xCQ $
and
$M( \xbs ' ) \xcs M( \xCN p_{ \xba }) \xcs \xbm^{+}(X)= \xCQ.$ But then
$M( \xbs + \xbs ' ) \xcs \xbm^{+}(X)$ $=$ $M( \xbs ) \xcs M( \xbs ' ) \xcs
\xbm^{+}(X)$
$ \xcc $ $M( \xbs ) \xcs M(p_{ \xba }) \xcs \xbm^{+}(X)$ $ \xcv $ $M( \xbs
' ) \xcs M( \xCN p_{ \xba }) \xcs \xbm^{+}(X)$ $=$ $ \xCQ,$
contradicting (1) of the induction hypothesis.

So induction hypothesis (1) will hold again.

Recall that for each candidate and any $ \xbs $ by induction hypothesis
(1)
$M( \xbs ) \xcs M(p_{ \xba }) \xcs \xbm^{+}(X)$ $=$ $ \xbm^{+}(M( \xbs )
\xcs M(p_{ \xba }) \xcs X)$ by $( \xbm =' ),$ and also for $ \xbs \xcc
\xbs ' $
$ \xbm^{+}(M( \xbs ' ) \xcs M(p_{ \xba }) \xcs X)$ $ \xcc $ $ \xbm^{+}(M(
\xbs ) \xcs M(p_{ \xba }) \xcs X)$ by $( \xbm =' )$ and $M( \xbs ' ) \xcc
M( \xbs ),$
and thus by $( \xbm^{-}4)$ $ \xbm^{-}(M( \xbs ' ) \xcs M(p_{ \xba }) \xcs
X)$ $ \xcc $ $ \xbm^{-}(M( \xbs ) \xcs M(p_{ \xba }) \xcs X).$

If we have only one candidate left, say e.g. $p_{ \xba },$ then for each
sufficiently
big sequence $ \xbs $ $M( \xbs ) \xcs M( \xCN p_{ \xba }) \xcs
\xbm^{+}(X)= \xCQ,$ thus for such $ \xbs $
$ \xbm^{+}(M( \xbs ) \xcs M(p_{ \xba }) \xcs X)$ $=$ $M( \xbs ) \xcs M(p_{
\xba }) \xcs \xbm^{+}(X)$ $=$ $M( \xbs ) \xcs \xbm^{+}(X)$ $=$ $
\xbm^{+}(M( \xbs ) \xcs X),$
and thus by $( \xbm^{-}4)$ $ \xbm^{-}(M( \xbs ) \xcs M(p_{ \xba }) \xcs
X)$ $=$ $ \xbm^{-}(M( \xbs ) \xcs X),$ so $ \xCN p_{ \xba }$ plays
no really important role. In particular, induction hypothesis (2) holds
again.

Suppose now that we have two candidates, thus for $p_{ \xba }$ and $ \xCN
p_{ \xba }$ and each $ \xbs $
$M( \xbs ) \xcs M(p_{ \xba }) \xcs \xbm^{+}(X) \xEd \xCQ $ and $M( \xbs )
\xcs M( \xCN p_{ \xba }) \xcs \xbm^{+}(X) \xEd \xCQ.$

By the same kind of argument as above we see that either for $p_{ \xba }$
or for $ \xCN p_{ \xba },$
or for both, and for all $ \xbs $
$x \xbe \xbm^{-}(M( \xbs ) \xcs M(p_{ \xba }) \xcs X)$ or $x \xbe
\xbm^{-}(M( \xbs ) \xcs M( \xCN p_{ \xba }) \xcs X).$

If not, there are $ \xbs $ and $ \xbs ' $ and
$x \xce \xbm^{-}(M( \xbs ) \xcs M(p_{ \xba }) \xcs X)$ $ \xcd $ $
\xbm^{-}(M( \xbs + \xbs ' ) \xcs M(p_{ \xba }) \xcs X)$ and
$x \xce \xbm^{-}(M( \xbs ' ) \xcs M( \xCN p_{ \xba }) \xcs X)$ $ \xcd $ $
\xbm^{-}(M( \xbs + \xbs ' ) \xcs M( \xCN p_{ \xba }) \xcs X),$ but
$ \xbm^{-}(M( \xbs + \xbs ' ) \xcs X)$ $=$ $ \xbm^{-}(M( \xbs + \xbs ' )
\xcs M(p_{ \xba }) \xcs X)$ $ \xcv $ $ \xbm^{-}(M( \xbs + \xbs ' ) \xcs M(
\xCN p_{ \xba }) \xcs X),$
so $x \xce \xbm^{-}(M( \xbs + \xbs ' ) \xcs X),$ contradicting the
induction hypothesis (2).

If we can choose both, we let the strategy decide, as for $p_{0}.$

So induction hypotheses (1) and (2) will hold again.

This gives a complete description of some $y$ (relative to the strategy!),
and
we set $x \xej y.$ We have to show: for all $Y \xbe \xdy $ $x \xbe
\xbm^{-}(Y)$ $ \xcr $ $x \xbe \xbm_{ \xej }(Y)$ $: \xcr $
$ \xcE y \xbe \xbm^{+}(Y).x \xej y.$
`` $ \xcp $ '': As we will do above construction for all $Y,$ it suffices
to
show that $y \xbe \xbm^{+}(X).$
`` $ \xcq $ '': Conversely, if the $y$ constructed above is in
$ \xbm^{+}(Y),$ then $x$ has to be in $ \xbm^{-}(Y).$

If $y \xce \xbm^{+}(X),$ then $Th(y)$ is inconsistent with $Th(
\xbm^{+}(X)),$ as $ \xbm^{+}$ is
definability preserving, so by classical compactness there is a suitable
finite
sequence $ \xbs $ with $M( \xbs ) \xcs \xbm^{+}(X)= \xCQ,$ but this was
excluded by the induction
hypothesis (1). So $y \xbe \xbm^{+}(X).$

Suppose $y \xbe \xbm^{+}(Y),$ but $x \xce \xbm^{-}(Y).$ So $y \xbe
\xbm^{+}(Y)$ and $y \xbe \xbm^{+}(X),$ and $Y=M( \xbf )$ for
some $ \xbf,$ so there will be a suitable finite sequence $ \xbs $ s.t.
for all $ \xbs ' $ with
$ \xbs \xcc \xbs ' $ $M( \xbs ' ) \xcs X \xcc M( \xbf )=Y,$
and by our construction $x \xbe \xbm^{-}(M( \xbs ' ) \xcs X).$ As $y \xbe
\xbm^{+}(X) \xcs \xbm^{+}(Y) \xcs (M( \xbs ' ) \xcs X),$
$ \xbm^{+}(M( \xbs ' ) \xcs X) \xcc \xbm^{+}(Y),$ so by $( \xbm^{-}4)$ $
\xbm^{-}(M( \xbs ' ) \xcs X) \xcc \xbm^{-}(Y),$ so $x \xbe \xbm^{-}(Y),$
$contradiction.$

We do now this construction for all strategies. Obviously, this does not
modify our results.

This finishes the completeness proof. $ \xcz $
\\[3ex]

As we postulated definability preservation, there are no problems to
translate
the result into logic. (Note that $ \xbn $ was applied to formula defined
model sets,
but the resulting sets were perhaps theory defined model sets.)

Comment:

One might try a construction similar to the one for Counterfactual
Conditionals, see  \cite{SM94},
and try to patch
together several ranked structures, one for each $K$ on the left, to
obtain
a general distance, by repeating elements.

So we would have different ``copies'' of A, say $A_{i},$ more precisely of
its
elements, and the natural definition seems to be: $A* \xbf \xcl \xbq $ iff
for all $i$
$A_{i}* \xbf \xcl \xbq,$ so $A \xfA B= \xcV \{A_{i} \xfA B:i \xbe I\}.$

But this does not work: Take $A:=\{a,a',a'' \},$ $B:=\{b,b' \},$ with $A
\xfA B:=\{b,b' \},$
and $a \xfA B=a' \xfA B=a'' \xfA B=\{b\}.$ Then for all copies of the
singletons, the result
cannot be empty, but must be $\{b\}.$ But $A \xfA B$ can only be a
``partial'' union of
the $x \xfA B,$ $x \xbe A,$ so it must be $\{b\}$ for all copies of A,
contradiction.

(Alternative definitions with copies fail too, but no systematic
investigation
was done.)

% ******* BEGIN LATEX SOURCE FILE 10-is.tex *******
%
% Uebers. aus Karltex File: 10-is.m
%
%
\chapter{
An analysis of defeasible inheritance systems
}
\label{Chapter Inheritance}
\label{Section Is-In}
\index{Section Is-In}
\section{
Introduction
}
\subsection{
Terminology
}

``Inheritance'' will stand here for ``nonmonotonic or defeasible
inheritance''. We will use indiscriminately ``inheritance system'',
``inheritance
diagram'', ``inheritance network'', ``inheritance net''.

In this introduction, we first give the connection to reactive diagrams,
then
give the motivation, then
describe in very brief terms some problems of inheritance diagrams, and
mention the basic ideas of our analysis.
\subsection{
Inheritance and reactive diagrams
}

Inheritance sytems or diagrams have an intuitive appeal. They seem close
to
human reasoning, natural, and are also implemented (see  \cite{Mor98}). Yet,
they
are a more procedural approach to nonmonotonic reasoning, and, to the
authors'
knowledge, a conceptual analysis, leading to a formal semantics, as well
as a
comparison to more logic based formalisms like the systems $P$ and $R$ of
preferential systems are lacking. We attempt to reduce
the gap between the more procedural and the more analytical approaches in
this
particular case. This will also give indications how to modify the systems
$P$ and
$R$ to approach them more to actual human reasoning. Moreover, we
establish a link
to multi-valued logics and the logics of information sources (see e.g.
 \cite{ABK07} and forthcoming work of the same authors, and also
 \cite{BGH95}).

%% Introduction to inheritance paper, 326

An inheritance net is a directed graph with two types of connections
between nodes $x \to y$ and $x \not\to y$. Diagram
\ref{Diagram Nixon-Diamond}
(page \pageref{Diagram Nixon-Diamond})
is such an
example.  The meaning of $x \to y$ is that $x$ is also a $y$ and the
meaning of $x \not\to y$ is that $x$ is not a $y$.  We do not allow
the combinations
$x \not\to y \not \to z$ or
$x \not\to y \to z$
but we do allow $x \to y \to z$ and
$x \to y \not \to z$.

Given a complex diagram such as Diagram
\ref{Diagram Up-Down-Chaining}
(page \pageref{Diagram Up-Down-Chaining})
and two points say $z$
and $y$, the question we ask is to determine from the diagram whether
the diagram says that
\begin{enumerate}
\item $z$ is $y$
\item $z$ is not $y$
\item nothing to say.
\end{enumerate}

Since in Diagram
\ref{Diagram Up-Down-Chaining}
(page \pageref{Diagram Up-Down-Chaining})
there are paths to $y$ from $z$ either through
$x$ or through $v$, we need to have an algorithm to decide. Let ${\cal A}$
be such an algorithm.

We need ${\cal A}$ to decide
\begin{enumerate}
\item are there valid paths from $z$ to $y$
\item of the opposing paths (one which supports `$z$ is $y$' and one
which supports `$z$ is not $y$'), which one wins (usually winning makes use
of being more specific but there are other possible options).
\end{enumerate}

So for example, in Diagram
\ref{Diagram Up-Down-Chaining}
(page \pageref{Diagram Up-Down-Chaining}),
the connection $x \to v$ makes paths
though $x$ more specific than paths through $v$. The question is whether
we have a valid path from $z$ to $x$.

In the literature, as well as in this paper, there are algorithms
for deciding the valid paths and the relative specificity of paths.
These are complex inductive algorithms, which may need the help of
a computer for the case of the more complex diagrams.

It seems that for inheritance networks we cannot adopt a simple minded
approach and just try to `walk' on the graph from $z$ to $y$, and
depending on what happens during this `walk' decide whether $z$ is
$y$ or not.

To explain what we mean, suppose we give the network a different meaning,
that of fluid flow. $x \to y$ means there is an open pipe from $x$
to $y$ and $x \not \to y$ means there is a blocked pipe from $x$ to
$y$.

To the question `can fluid flow from $z$ to $y$
in Diagram
\ref{Diagram Up-Down-Chaining}'
(page \pageref{Diagram Up-Down-Chaining}),
there
is a simple answer:
\begin{quote}
Fluid can flow iff there is a path comprising of $\to$ only (without
any $\not\to$ among them).
\end{quote}

Similarly we can ask in the inheritance network something
like (*) below:
\begin{quote}
\begin{itemize}
\item [(*)] $z$ is (resp.\ is not) $y$ according to diagram $D$, iff
there is a path $\pi$ from $z$ to $y$ in $D$ such that some non-inductive
condition $\psi(\pi)$ holds for the path $\pi$.
\end{itemize}
\end{quote}
Can we offer the reader such a $\psi$?

If we do want to help the user to `walk' the graph and get an answer,
we can proceed as one of the following options:

\begin{itemize}
\item [Option 1.]  Add additional annotations to paths to obtain $D^*$
from $D$, so that a predicate $\psi$ can be defined on $D^*$ using
these annotations.  Of course these annotations will be computed using
the inductive algorithm in ${\cal A}$, i.e.\ we modify
${\cal A}$ to ${\cal A}^*$ which also executes the annotations.
\item [Option 2.] Find a transformation $\tau$ on diagrams $D$ to
transform $D$ to $D'=\tau(D)$, such that a predicate $\psi$ can be
found for $D'$. So we work on $D'$ instead of on $D$.
\end{itemize}
We require a compatibility condition on options 1 and 2:
\begin{itemize}
\item [(C1)] If we apply ${\cal A}^*$ to $D^*$ we get $D^*$ again.
\item [(C2)] $\tau(\tau(D)) = \tau(D)$.
\end{itemize}

We now present the tools we use for our annotations and transformation.
These are the reactive double arrows.

Consider the following Diagram
\ref{Diagram Dov-Is-1}
(page \pageref{Diagram Dov-Is-1}):

\vspace{10mm}

\begin{diagram}

\label{Diagram Dov-Is-1}
\index{Diagram Dov-Is-1}

\centering
\setlength{\unitlength}{0.00083333in}
{\renewcommand{\dashlinestretch}{30}

\begin{picture}(2630,4046)(0,-10)
\put(1289.500,1811.500){\arc{2059.430}{4.5293}{7.6709}}
\blacken\path(1217.372,2868.600)(1102.000,2824.000)(1224.788,2809.060)(1217.372,2868.600)
\path(1102,274)(1852,1324)
\blacken\path(1806.663,1208.915)(1852.000,1324.000)(1757.839,1243.789)(1806.663,1208.915)
\path(1102,274)(352,1324)
\blacken\path(446.161,1243.789)(352.000,1324.000)(397.337,1208.915)(446.161,1243.789)
\path(352,1324)(1102,2374)
\blacken\path(1056.663,2258.915)(1102.000,2374.000)(1007.839,2293.789)(1056.663,2258.915)
\path(1852,1324)(1102,2374)
\blacken\path(1196.161,2293.789)(1102.000,2374.000)(1147.337,2258.915)(1196.161,2293.789)
\path(1102,2374)(1102,3274)
\blacken\path(1132.000,3154.000)(1102.000,3274.000)(1072.000,3154.000)(1132.000,3154.000)
\path(1412,2854)(1257,2844)
\blacken\path(1374.820,2881.664)(1257.000,2844.000)(1378.683,2821.788)(1374.820,2881.664)

\path(2452,274)(1852,1324)
\blacken\path(1937.584,1234.695)(1852.000,1324.000)(1885.489,1204.927)(1937.584,1234.695)

\put(1950,1300){\xssc $c$}
\put(1252,2299){\xssc $d$}
\put(1177,3199){\xssc $e$}
\put(277,1249) {\xssc $b$}
\put(1102,49)  {\xssc $a$}
\put(2477,59)  {\xssc $a'$}

\put(20,3500) {{\rm\bf Reactive graph}}

\end{picture}
}

\end{diagram}

\vspace{4mm}

We want to walk from $a$ to $e$. If we go to $c$, a double arrow from
the arc $a \to c$ blocks the way from $d$ to $e$. So the only way
to go to $e$ is throgh $b$. If we start at $a'$ there is no such block.
It is the travelling through the arc $(a,c)$ that triggers the double
arrow $(a,c)\twoheadrightarrow (d,e)$.

We want to use $\twoheadrightarrow$ in ${\cal A}^*$ and in $\tau$.

So in Diagram
\ref{Diagram Up-Down-Chaining}
(page \pageref{Diagram Up-Down-Chaining})
the
path $u \to x \not\to y$ is winning over the
path $u \to v \to y$, because of the specificity arrow $x \to v$. However,
if we start at $z$ then the path $z \to u \to v \to y$ is valid because
of $z \not\to x$.  We can thus add the following double arrows to
the diagram. We get Diagram
\ref{Diagram Dov-Is-2}
(page \pageref{Diagram Dov-Is-2}).

\vspace{10mm}

\begin{diagram}

\label{Diagram Dov-Is-2}
\index{Diagram Dov-Is-2}

\centering
\setlength{\unitlength}{0.00083333in}
{\renewcommand{\dashlinestretch}{30}

\begin{picture}(3253,3820)(0,-10)
\put(-6848.000,6274.000){\arc{20135.789}{0.4303}{0.5839}}
\blacken\path(2278.077,1952.642)(2302.000,2074.000)(2223.782,1978.177)(2278.077,1952.642)
\put(903.476,2620.958){\arc{992.406}{4.9432}{6.4593}}
\blacken\path(1139.941,3090.381)(1017.000,3104.000)(1119.068,3034.128)(1139.941,3090.381)
\put(531.951,2405.800){\arc{1183.940}{5.2342}{7.0408}}
\blacken\path(941.882,2873.151)(827.000,2919.000)(906.791,2824.483)(941.882,2873.151)
\put(1490.621,1557.466){\arc{1660.763}{1.5028}{4.2591}}
\blacken\path(1038.516,2217.568)(1127.000,2304.000)(1008.252,2269.376)(1038.516,2217.568)
\path(2902,2524)(352,2524)
\blacken\path(472.000,2554.000)(352.000,2524.000)(472.000,2494.000)(472.000,2554.000)
\path(1552,1474)(352,2524)
\blacken\path(462.064,2467.557)(352.000,2524.000)(422.554,2422.402)(462.064,2467.557)
\path(1552,1474)(2902,2524)
\blacken\path(2825.696,2426.647)(2902.000,2524.000)(2788.860,2474.008)(2825.696,2426.647)
\path(2902,2524)(1552,3574)
\blacken\path(1665.140,3524.008)(1552.000,3574.000)(1628.304,3476.647)(1665.140,3524.008)
\path(352,2524)(1552,3574)
\blacken\path(1481.446,3472.402)(1552.000,3574.000)(1441.936,3517.557)(1481.446,3472.402)
\path(1552,274)(1552,1324)
\blacken\path(1582.000,1204.000)(1552.000,1324.000)(1522.000,1204.000)(1582.000,1204.000)
\path(1627,274)(2902,2524)
\blacken\path(2868.939,2404.807)(2902.000,2524.000)(2816.738,2434.388)(2868.939,2404.807)
\path(2157,1774)(2237,1939)
\blacken\path(2211.642,1817.934)(2237.000,1939.000)(2157.653,1844.111)(2211.642,1817.934)
\path(1262,2959)(1162,3049)
\blacken\path(1271.264,2991.023)(1162.000,3049.000)(1231.126,2946.425)(1271.264,2991.023)
\path(1037,2694)(957,2814)
\blacken\path(1048.526,2730.795)(957.000,2814.000)(998.603,2697.513)(1048.526,2730.795)
\path(867,2119)(1002,2219)
\blacken\path(923.430,2123.466)(1002.000,2219.000)(887.716,2171.679)(923.430,2123.466)
\path(2212,3184)(2097,3059)
\path(2162,1344)(2282,1279)

\put(1552,1324){\xssc $u$}
\put(1552,49)  {\xssc $z$}
\put(1552,3649){\xssc $y$}
\put(2977,2449){\xssc $x$}
\put(277,2449) {\xssc $v$}

\end{picture}

}

\end{diagram}

\vspace{4mm}

If we start from $u$, and go to $u \to v$, then $v \to y$ is cancelled.
Similary $u \to x \to v$ cancelled $v \to y$. So the only path is
$u\to x \not\to y$.

If we start from $z$ then $u\to x$ is cancelled and so is the cancellation
$(u,v)\twoheadrightarrow (v,y)$. hence the path $z \to u \to v \to
y$ is open.

We are not saying that $\xbt(Diagram \ref{Diagram Dov-Is-1} )=
Diagram \ref{Diagram Dov-Is-2}$ but something
effectively similar will be done by $\tau$.

We emphasize that the construction depends on the point of departure.
Consider Diagram \ref{Diagram Up-Down-Chaining} (page \pageref{Diagram
Up-Down-Chaining}).
Starting at $u,$ we will have to block the path uvy. Starting at $z,$ the
path
zuvy has to be free.
See Diagram \ref{Diagram U-D-reactive} (page \pageref{Diagram U-D-reactive}).

\vspace{10mm}

\begin{diagram}

\label{Diagram U-D-reactive}
\index{Diagram U-D-reactive}

\unitlength1.0mm
\begin{picture}(130,100)

\newsavebox{\ZWEIvierreac}
\savebox{\ZWEIvierreac}(140,110)[bl]
{

\put(0,95){{\rm\bf The problem of downward chaining - reactive}}

\put(43,27){\vector(1,1){24}}
\put(37,27){\vector(-1,1){24}}
\put(13,57){\vector(1,1){24}}
\put(67,57){\vector(-1,1){24}}

\put(53,67){\line(1,1){4}}

\put(67,54){\vector(-1,0){54}}

\put(40,7){\vector(0,1){14}}
\put(43,7){\line(3,5){24}}
\put(58,28.1){\line(-5,3){3.6}}

\put(24,42){\vector(0,1){24}}
\put(24,42){\vector(0,1){22.5}}

\put(39,3){$z$}
\put(39,23){$u$}
\put(9,53){$v$}
\put(69,53){$x$}
\put(39,83){$y$}

}

\put(0,0){\usebox{\ZWEIvierreac}}
\end{picture}

\end{diagram}

\vspace{4mm}

So we cannot just add a double arrow from $u \xcp v$ to $v \xcp y,$
blocking $v \xcp y,$ and leave it there when we start from $z.$ We will
have to erase
it when we change the origin.

At the same time, this shows an advantage over just erasing the arrow $v
\xcp y:$

When we change the starting point, we can erase simply all double arrows,
and do not have to remember the original diagram.

How do we construct the double arrows given some origin $x?$

First, if all possible paths are also valid, there is nothing to do.
(At the same time, this shows that applying the procedure twice will not
result in anything new.)

Second, remember that we have an upward chaining formalism. So if a
potential
path fails to be valid, it will do so at the end.

Third, suppose that we have two valid paths $ \xbs:x \xcp y$ and $ \xbt
:x \xcp y.$

If they are negative (and they are either both negative or positive, of
course),
then they cannot be continued. So if there is an arrow
$y \xcp z$ or $y \xcP z,$ we will block it by a double arrow from the
first arrow of $ \xbs $ and
from the first arrow of $ \xbt $ to $y \xcp z$ $(y \xcP z$ respectively).

If they are positive, and there is an arrow $y \xcp z$ or $y \xcP z,$ both
$ \xbs y \xcp z$ and $ \xbt y \xcp z$
are potential paths (the case $y \xcP z$ is analogue). One is valid iff
the other one
is, as $ \xbs $ and $ \xbt $ have the same endpoint, so preclusion, if
present, acts on
both the same way. If they are not valid, we block $y \xcp z$
by a double arrow from the first arrow of $ \xbs $ and
from the first arrow of $ \xbt $ to $y \xcp z.$

Of course, if there is only one such $ \xbs,$ we do the same, there is
just to
consider to see the justification.

We summarize:

Our algorithm switches impossible continuations off by making them
invisible,
there is just no more arrow to concatenate. As validity in inheritance
networks is not forward looking - validity of $ \xbs:x \xcp y$ does not
depend on what
is beyond $y$ - validity in the old and in the new network starting at $x$
are
the same. As we left only valid paths, applying the algorithm twice will
not give anything new.

We illustrate this by
considering Diagram \ref{Diagram InherUniv} (page \pageref{Diagram InherUniv}).

First, we add double arrows for starting point $c,$ see
Diagram \ref{Diagram I-U-reac-c} (page \pageref{Diagram I-U-reac-c}),

\vspace{10mm}

\begin{diagram}

\label{Diagram I-U-reac-c}
\index{Diagram I-U-reac-c}

\centering
\setlength{\unitlength}{1mm}
{\renewcommand{\dashlinestretch}{30}
\begin{picture}(150,150)(0,0)

\put(10,60){\circle*{1}}
\put(70,60){\circle*{1}}
\put(130,60){\circle*{1}}

\put(70,10){\circle*{1}}
\put(70,110){\circle*{1}}

\put(40,50){\circle*{1}}
\put(40,70){\circle*{1}}
\put(100,50){\circle*{1}}
\put(100,70){\circle*{1}}

\path(69,11)(12,58)
\path(13.5,55.4)(12,58)(14.8,57)
\path(71,11)(128,58)
\path(125.2,57)(128,58)(126.5,55.4)
\path(12,62)(69,109)
\path(66.2,108)(69,109)(67.5,106.4)
\path(128,62)(71,109)
\path(72.5,106.4)(71,109)(73.8,108)

\path(70,61)(70,108)
\path(68,86)(72,86)
\path(69,105.2)(70,108)(71,105.2)

\path(11,60.5)(39,69.5)
\path(14,60.4)(11,60.5)(13.4,62.3)
\path(11,59.5)(39,50.5)
\path(13.4,57.7)(11,59.5)(14,59.6)
\path(41,69.5)(69,60.5)
\path(43.4,67.7)(41,69.5)(44,69.6)
\path(41,50.5)(69,59.5)
\path(44,50.4)(41,50.5)(43.4,52.3)

\path(24,67)(26,63)

\path(71,60.5)(99,69.5)
\path(71,59.5)(99,50.5)
\path(101,69.5)(129,60.5)
\path(101,50.5)(129,59.5)
\path(84,67)(86,63)

\path(74,60.4)(71,60.5)(73.4,62.3)
\path(73.4,57.7)(71,59.5)(74,59.6)
\path(103.4,67.7)(101,69.5)(104,69.6)
\path(104,50.4)(101,50.5)(103.4,52.3)

\path(100,52)(100,68)
\path(99.2,65.9)(100,68)(100.8,65.9)

\put(69,5){{\xssc $x$}}
\put(69,115){{\xssc $y$}}

\put(7,60){{\xssc $a$}}
\put(69,57){{\xssc $b$}}
\put(132,60){{\xssc $c$}}

\put(39,47){{\xssc $d$}}
\put(99,47){{\xssc $e$}}
\put(39,72){{\xssc $f$}}
\put(99,72){{\xssc $g$}}

\path(109,68)(109,75)(89,75)(89,68)
\put(89,75){\vector(0,-1){7}}
\put(89,75){\vector(0,-1){6}}
\path(115,56)(115,76)(86,76)(86,67)
\put(86,76){\vector(0,-1){9}}
\put(86,76){\vector(0,-1){8}}
\path(116,56.5)(116,77)(72,77)
\put(116,77){\vector(-1,0){44}}
\put(116,77){\vector(-1,0){43}}
\path(117,57)(117,120)(22,120)(22,66)
\put(22,120){\vector(0,-1){54}}
\put(22,120){\vector(0,-1){53}}
\path(117,54)(117,3)(22,3)(22,54)
\put(22,3){\vector(0,1){51}}
\put(22,3){\vector(0,1){50}}

\end{picture}
}

\end{diagram}

\vspace{4mm}

and then for starting point $x,$ see
Diagram \ref{Diagram I-U-reac-x} (page \pageref{Diagram I-U-reac-x}).

\vspace{10mm}

\begin{diagram}

\label{Diagram I-U-reac-x}
\index{Diagram I-U-reac-x}

\centering
\setlength{\unitlength}{1mm}
{\renewcommand{\dashlinestretch}{30}
\begin{picture}(150,150)(0,0)

\put(10,60){\circle*{1}}
\put(70,60){\circle*{1}}
\put(130,60){\circle*{1}}

\put(70,10){\circle*{1}}
\put(70,110){\circle*{1}}

\put(40,50){\circle*{1}}
\put(40,70){\circle*{1}}
\put(100,50){\circle*{1}}
\put(100,70){\circle*{1}}

\path(69,11)(12,58)
\path(13.5,55.4)(12,58)(14.8,57)
\path(71,11)(128,58)
\path(125.2,57)(128,58)(126.5,55.4)
\path(12,62)(69,109)
\path(66.2,108)(69,109)(67.5,106.4)
\path(128,62)(71,109)
\path(72.5,106.4)(71,109)(73.8,108)

\path(70,61)(70,108)
\path(68,86)(72,86)
\path(69,105.2)(70,108)(71,105.2)

\path(11,60.5)(39,69.5)
\path(14,60.4)(11,60.5)(13.4,62.3)
\path(11,59.5)(39,50.5)
\path(13.4,57.7)(11,59.5)(14,59.6)
\path(41,69.5)(69,60.5)
\path(43.4,67.7)(41,69.5)(44,69.6)
\path(41,50.5)(69,59.5)
\path(44,50.4)(41,50.5)(43.4,52.3)

\path(24,67)(26,63)

\path(71,60.5)(99,69.5)
\path(71,59.5)(99,50.5)
\path(101,69.5)(129,60.5)
\path(101,50.5)(129,59.5)
\path(84,67)(86,63)

\path(74,60.4)(71,60.5)(73.4,62.3)
\path(73.4,57.7)(71,59.5)(74,59.6)
\path(103.4,67.7)(101,69.5)(104,69.6)
\path(104,50.4)(101,50.5)(103.4,52.3)

\path(100,52)(100,68)
\path(99.2,65.9)(100,68)(100.8,65.9)

\put(69,5){{\xssc $x$}}
\put(69,115){{\xssc $y$}}

\put(7,60){{\xssc $a$}}
\put(69,57){{\xssc $b$}}
\put(132,60){{\xssc $c$}}

\put(39,47){{\xssc $d$}}
\put(99,47){{\xssc $e$}}
\put(39,72){{\xssc $f$}}
\put(99,72){{\xssc $g$}}

\path(92,27)(140,27)(140,75)(89,75)(89,68)
\put(89,75){\vector(0,-1){7}}
\put(89,75){\vector(0,-1){6}}

\path(90,25)(142,25)(142,77)(72,77)
\put(116,77){\vector(-1,0){44}}
\put(116,77){\vector(-1,0){43}}

\path(88,23)(144,23)(144,120)(22,120)(22,66)
\put(22,120){\vector(0,-1){54}}
\put(22,120){\vector(0,-1){53}}

\path(80,17)(80,3)(22,3)(22,54)
\put(22,3){\vector(0,1){51}}
\put(22,3){\vector(0,1){50}}

\end{picture}
}

\end{diagram}

\vspace{4mm}

For more information on reactive diagrams, see also
 \cite{Gab08b}, and an earlier
version:  \cite{Gab04}.
\subsection{
Conceptual analysis
}

Inheritance diagrams are deceptively simple. Their conceptually
complicated
nature is seen by e.g. the fundamental difference between direct links and
valid paths, and the multitude of existing formalisms, upward vs. downward
chaining, intersection of extensions vs. direct scepticism, on-path vs.
off-path
preclusion (or pre-emption), split validity vs. total validity preclusion
etc.,
to name a few, see the discussion in Section \ref{Section IS-2.3} (page
\pageref{Section IS-2.3}).
Such a proliferation of
formalisms usually hints at deeper
problems on the conceptual side, i.e. that the underlying ideas are
ambigous,
and not sufficiently analysed. Therefore, any clarification and resulting
reduction of possible formalisms seems a priori to make progress. Such
clarification will involve conceptual decisions, which need not be shared
by
all, they can only be suggestions. Of course, a proof that such decisions
are
correct is impossible, and so is its contrary.

We will introduce into the analysis of inheritance systems a number of
concepts
not usually found in the field, like multiple truth values, access to
information, comparison of truth values, etc. We think that this
additional
conceptual burden pays off by a better comprehension and analysis of the
problems behind the surface of inheritance.

We will also see that some distinctions between inheritance formalisms go
far
beyond questions of inheritance, and concern general problems of treating
contradictory information - isolating some of these is another objective
of
this article.

The text is essentially self-contained, still some familiarity with the
basic
concepts of inheritance systems and nonmonotonic logics in general is
helpful.
For a presentation, the reader might look into  \cite{Sch97-2} and
 \cite{Sch04}.

The text is organized as follows. After an introduction to inheritance
theory,
connections with reactive diagrams
in Section \ref{Section Def-Reac} (page \pageref{Section Def-Reac}),
and big and small subsets and the systems $P$ and $R$
in Section \ref{Section Is-Inh-Intro} (page \pageref{Section Is-Inh-Intro}), we
turn to an informal description of the fundamental differences between
inheritance and the systems $P$ and $R$ in Section \ref{Section 4.2} (page
\pageref{Section 4.2}),
give an analysis of
inheritance systems in terms of information and information flow in
Section \ref{Section 4.3} (page \pageref{Section 4.3}), then
in terms of reasoning with prototypes in
Section \ref{Section 4.4} (page \pageref{Section 4.4}), and
conclude in Section \ref{Section Translation} (page \pageref{Section
Translation})
with a translation
of inheritance into (necessarily deeply modified)
coherent systems of big and small sets, respectively logical systems $P$
and $R.$
One of the main modifications will be to relativize the notions of small
and
big, which thus become less ``mathematically pure'' but perhaps closer to
actual
use in ``dirty'' common sense reasoning.
\section{
Introduction to nonmonotonic inheritance
}
\label{Section Is-Inh-Intro}
\index{Section Is-Inh-Intro}
%  2.1  Basic discussion
%  2.1  Basic discussion
% %
% ======================
\subsection{
Basic discussion
}
\label{Section IS-2.1}

We give here an informal discussion. The reader unfamiliar with
inheritance
systems should consult in parallel Definition \ref{Definition 2.3} (page
\pageref{Definition 2.3})  and
Definition \ref{Definition 2.4} (page \pageref{Definition 2.4}). As
there are many variants of the definitions, it seems
reasonable to discuss them before a formal introduction, which, otherwise,
would seem to pretend to be definite without being so.

\paragraph{
(Defeasible or nonmonotonic) inheritance networks or diagrams
}

$\hspace{0.01em}$

% (+++ Orig.:  (Defeasible or nonmonotonic) inheritance networks or diagrams
% +++)

\label{Section (Defeasible or nonmonotonic) inheritance networks or diagrams}

Nonmonotonic inheritance systems describe situations like
``normally, birds fly'', written $birds \xcp fly.$ Exceptions are
permitted, ``normally penguins $don' t$ fly'', $penguins \xcP fly.$

\bd

$\hspace{0.01em}$

% (+++ Orig. No.:  Definition IS-2.1 +++)

\label{Definition IS-2.1}

\label{Definition Inheritance-Net}

A nonmonotonic inheritance net is a finite DAG, directed, acyclic graph,
with two types of arrows or links, $ \xcp $ and $ \xcP,$ and labelled
nodes. We will use
$ \xbG $ etc. for such graphs, and $ \xbs $ etc. for paths - the latter to
be defined below.

\ed

Roughly (and to be made precise and modified below, we try to give here
just a
first intuition), $X \xcp Y$ means that
``normal'' elements of $X$ are in $Y,$ and $X \xcP Y$ means that ``normal''
elements of $X$ are
not in $Y.$ In a semi-quantitative set interpretation, we will read ``most''
for
``normal'', thus ``most elements of $X$ are in $Y'',$ ''most elements of $X$
are not in
$Y'',$ etc. These are by no means the only interpretations, as we will
see - we
will use these expressions for the moment just to help the reader's
intuition.
We should add immediately a word of warning: ``most'' is here not
necessarily,
but only by default, transitive, in the following sense. In the Tweety
diagram,
see Diagram \ref{Diagram Tweety} (page \pageref{Diagram Tweety})  below,
most penguins are birds, most birds fly, but it
is not the case that most penguins fly. This is the problem of transfer of
relative size which will be discussed extensively, especially
in Section \ref{Section Translation} (page \pageref{Section Translation}).

According to the set interpretation,
we will also use informally expressions like $X \xcs Y,$ $X-$Y, $ \xdC X$
- where $ \xdC $ stands
for set complement -, etc. But we
will also use nodes informally as formulas, like $X \xcu Y,$ $X \xcu \xCN
Y,$ $ \xCN X,$ etc.
All this will only be used here as an appeal to intuition.

Nodes at the beginning of an arrow can also stand for individuals,
so $Tweety \xcP fly$ means something like: ``normally, Tweety will not
fly''. As always
in nonmonotonic systems, exceptions are permitted, so the soft rules
``birds
fly'', ``penguins $don' t$ fly'', and (the hard rule) ``penguins are birds''
can coexist
in one diagram, penguins are then abnormal birds (with respect to flying).
The direct link $penguins \xcP fly$ will thus be accepted, or considered
valid, but
not the composite path $penguins \xcp birds \xcp fly,$ by specificity -
see below.
This is illustrated by Diagram \ref{Diagram Tweety} (page \pageref{Diagram
Tweety}),
where a stands for Tweety, $c$ for penguins,
$b$ for birds, $d$ for flying animals or objects.

(Remark: The arrows $a \xcp c,$ $a \xcp b,$ and $c \xcp b$ can also be
composite
paths - see below for the details.)
% Diagram 2.1
% Diagram 2.1
% -----------
% The Tweety diagram
% The Tweety diagram
% @@         %d
% @@         %d
% @@
% @@       E    F
% @@      E      f
% @@     E        F
% @@    E          F
% @@
% @@ %b caaaaaaaaaaa %c
% @@
% @@    F          E
% @@     F        E
% @@      F      E
% @@       F    E
% @@
% @@         %a

\vspace{10mm}

\begin{diagram}

\label{Diagram Tweety}
\index{Diagram Tweety}

\unitlength1.0mm
\begin{picture}(130,100)(0,0)

\newsavebox{\SECHSacht}
\savebox{\SECHSacht}(140,110)[bl]
{

\put(0,95){{\rm\bf The Tweety diagram}}

\put(43,27){\vector(1,1){24}}
\put(37,27){\vector(-1,1){24}}
\put(13,57){\vector(1,1){24}}
\put(67,57){\vector(-1,1){24}}

\put(53,67){\line(1,1){4}}

\put(67,54){\vector(-1,0){54}}

\put(39,23){$a$}
\put(9,53){$b$}
\put(69,53){$b$}
\put(39,83){$d$}

}

\put(0,0){\usebox{\SECHSacht}}
\end{picture}

\end{diagram}

\vspace{4mm}

(Of course, there is an analogous case for the opposite polarity, i.e.
when the arrow from $b$ to $d$ is negative, and the one from $c$ to $d$ is
positive.)

The main problem is to define in an intuitively acceptable way
a notion of valid path, i.e. concatenations of arrows satisfying certain
properties.

We will write $ \xbG \xcm \xbs,$ if $ \xbs $ is a valid path in the
network $ \xbG,$ and if
$x$ is the origin, and $y$ the endpoint of $ \xbs,$ and $ \xbs $ is
positive, we will write
$ \xbG \xcm xy,$ i.e. we will accept the conclusion that x's are $y' $s,
and analogously
$ \xbG \xcm x \ol{y}$ for negative paths. Note that we
will not accept any other conclusions, only those established by a valid
path, so many questions about conclusions have a trivial negative answer:
there
is obviously no path from $x$ to $y.$ E.g., there is no path from $b$ to
$c$ in
Diagram \ref{Diagram Tweety} (page \pageref{Diagram Tweety}). Likewise, there
are no disjunctions,
conjunctions etc. in our conclusions, and negation is present only in a
strong
form: ``it is not the case that x's are normally $y' $s'' is not a possible
conclusion, only ``x's are normally not $y' $s'' is one. Also, possible
contradictions are contained, there is no EFQ.

To simplify matters, we assume that for no two nodes $x,y \xbe \xbG $ $x
\xcp y$ and $x \xcP y$ are
both in $ \xbG,$ intuitively, that $ \xbG $ is free from (hard)
contradictions. This
restriction is inessential for our purposes. We admit, however, soft
contradictions, and preclusion, which allows us to solve some soft
contradictions -
as we already did in the penguins example. We will also assume that all
arrows
stand for rules with possibly exceptions, again, this restriction is not
important for our purposes. Moreover, in the abstract treatment, we will
assume that all nodes stand for (nonempty) sets, though this will not be
true
for all examples discussed.

This might be the place for a remark on absence of cycles. Suppose we also
have a positive arrow from $b$ to $c$ in Diagram \ref{Diagram Tweety} (page
\pageref{Diagram Tweety}).
Then, the concept of preclusion collapses,
as there are now equivalent arguments to accept $a \xcp b \xcp d$ and $a
\xcp c \xcP d.$ Thus,
if we do not want to introduce new complications, we cannot rely on
preclusion
to decide conflicts. It seems that this would change the whole outlook on
such diagrams. The interested reader will find more on the subject in
 \cite{Ant97},  \cite{Ant99},
 \cite{Ant05}.

Inheritance networks were introduced about 20 years ago (see e.g.  \cite{Tou84},
 \cite{Tou86},  \cite{THT87}), and exist
in a multitude of more or less differing formalisms, see
e.g.  \cite{Sch97-2} for a brief discussion. There still does not
seem
to exist a satisfying semantics for these networks. The authors' own
attempt
 \cite{Sch90} is an a posteriori semantics, which cannot explain or
criticise or
decide between the different formalisms. We will give here a conceptual
analysis, which provides also at least some building blocks for a
semantics,
and a translation into (a modified version of) the language of small and
big
subsets, familiar from preferential structures, see
Definition \ref{Definition Gen-Filter} (page \pageref{Definition Gen-Filter}).

We will now discuss the two fundamental situations of contradictions,
then give a detailed inductive definition of valid paths for a certain
formalism so the reader has firm ground under his feet, and then present
briefly
some alternative formalisms.

As in all of nonmonotonic reasoning, the interesting questions arise in
the
treatment of contradictions and exceptions. The
difference in quality of information is expressed by ``preclusion''
(or $'' pre-emption'' ).$ The basic diagram is the Tweety diagram,
see Diagram \ref{Diagram Tweety} (page \pageref{Diagram Tweety}).

Unresolved contradictions give either rise to a branching
into different extensions, which may roughly be seen as maximal consistent
subsets, or to mutual cancellation in directly sceptical approaches.
The basic diagram for the latter is the Nixon Diamond, see
Diagram \ref{Diagram Nixon-Diamond} (page \pageref{Diagram Nixon-Diamond}),
where
$a=Nixon,$ $b=Quaker,$ $c=Republican,$ $d=pacifist.$

In the directly sceptical approach, we will not accept any path from a to
$d$
as valid, as there is an unresolvable contradiction between the two
candidates.
% Diagram 2.2
% Diagram 2.2
% -----------
% The Nixon Diamond
% The Nixon Diamond
% @@         %d
% @@         %d
% @@
% @@       E    F
% @@      E      f
% @@     E        F
% @@    E          F
% @@
% @@ %b              %c
% @@
% @@    F          E
% @@     F        E
% @@      F      E
% @@       F    E
% @@
% @@         %a

\vspace{10mm}

\begin{diagram}

\label{Diagram Nixon-Diamond}
\index{Diagram Nixon-Diamond}

\unitlength1.0mm
\begin{picture}(130,100)

\newsavebox{\ZWEIzwei}
\savebox{\ZWEIzwei}(140,110)[bl]
{

\put(0,95){{\rm\bf The Nixon Diamond}}

\put(43,27){\vector(1,1){24}}
\put(37,27){\vector(-1,1){24}}
\put(13,57){\vector(1,1){24}}
\put(67,57){\vector(-1,1){24}}

\put(53,67){\line(1,1){4}}

\put(39,23){$a$}
\put(9,53){$b$}
\put(69,53){$c$}
\put(39,83){$d$}

}

\put(0,0){\usebox{\ZWEIzwei}}
\end{picture}

\end{diagram}

\vspace{4mm}

The extensions approach can be turned into an indirectly sceptical one,
by forming first all extensions, and then taking the intersection of
either the
sets of valid paths, or of valid conclusions, see [MS91] for a detailed
discussion.
See also Section \ref{Section Review} (page \pageref{Section Review})  for more
discussion on directly
vs. indirectly
sceptical approaches.

In more detail:

\paragraph{
Preclusion
}

$\hspace{0.01em}$

% (+++ Orig.:  Preclusion +++)

\label{Section Preclusion}

In the above example, our intuition tells us that it is not admissible
to conclude from the fact that penguins are birds, and that most birds
fly that most penguins fly. The horizontal arrow $c \xcp b$ together with
$c \xcP d$ bars
this conclusion, it expresses specificity. Consequently, we have
to define the conditions under which two potential paths neutralize each
other, and when one is victorious. The idea is as follows:
1) We want to be sceptical, in the sense that we do not believe every
potential path. We will not arbitrarily chose one either.
2) Our scepticism will be restricted, in the sense that we will often make
well defined choices for one path in the case of conflict:
a) If a compound potential path is in conflict with a direct link,
the direct link wins.
$b)$ Two conflicting paths of the same type neutralize each other, as in
the
Nixon Diamond, where neither potential path will be valid.
$c)$ More specific information will win over less specific one.

(It is essential in the Tweety diagram that the arrow $c \xcP d$ is a
direct link, so
it is in a
way stronger than compound paths.) The arrows $a \xcp b,$ $a \xcp c,$ $c
\xcp b$ can also
be composite paths: The path from $c$ to $b$ (read $c \xcc $  \Xl  $ \xcc
b,$ where $ \xcc $ stands
here for soft inclusion),
however, tells us, that the information coming from $c$ is more specific
(and thus considered more reliable), so the negative path from a to $d$
via
$c$ will win over the positive one via $b.$ The precise inductive
definition
will be given below. This concept is evidently independent of the lenght
of the
paths, $a \xFB \xcp c$ may be much longer than $a \xFB \xcp b,$ so this is
not shortest
path reasoning (which has some nasty drawbacks, discussed e.g. in
[HTT87]).

A final remark: Obviously, in some cases, it need not be specificity,
which
decides conflicts. Consider the case where Tweety is a bird, but a dead
animal.
Obviously, Tweety will not fly, here because the predicate ``dead'' is very
strong and overrules many normal properties. When we generalize this, we
might
have a hierarchy of causes, where one overrules the other, or the result
may be
undecided. For instance, a falling object might be attracted in a magnetic
field, but a gusty wind might prevent this, sometimes, with unpredictable
results. This is then additional information (strength of cause), and this
problem is not addressed directly in traditional inheritance networks, we
would have to introduce a subclass ``dead bird'' - and subclasses often have
properties of ``pseudo-causes'', being a penguin probably is not a ``cause''
for not flying, nor bird for flying, still, things change from class to
subclass for a reason.

Before we give a formalism based on these ideas, we refine them, adopt one
possibility (but indicate some modifications), and discuss alternatives
later.
%  2.2  Directly sceptical split validity upward chaining off-path inheritance
%  2.2  Directly sceptical split validity upward chaining off-path inheritance
% %
% ============================================================================
\subsection{
Directly sceptical split validity upward chaining off-path inheritance
}
\label{Section IS-2.2}

Our approach will be directly sceptical, i.e. unsolvable contradictions
result
in the absence of valid paths, it is upward chaining, and split-validity
for
preclusions (discussed below, in particular
in Section \ref{Section IS-2.3} (page \pageref{Section IS-2.3}) ). We will
indicate modifications to make it
extension based, as well as for total validity preclusion.
This approach is strongly inspired by classical work in the field
by Horty, Thomason, Touretzky, and others, and we claim no priority
whatever.
If it is new at all, it is a very minor modification of existing
formalisms.

Our conceptual ideas to be presented in detail
in Section \ref{Section 4.3} (page \pageref{Section 4.3})  make split
validity, off-path preclusion and upward chaining a natural choice.

For the reader's convenience, we give here a very short resume of these
ideas:
We consider only arrows as information, e.g. $a \xcp b$ will be considered
information $b$ valid at or for a. Valid composed positive paths will not
be
considered information in our sense. They will be seen as a way to obtain
information, so a valid path $ \xbs:x \Xl  \xcp a$ makes information $b$
accessible to $x,$
and, secondly, as a means of comparing information
strength, so a valid path $ \xbs:a \Xl. \xcp a' $ will make information
at a stronger than
information at $a'.$ Valid negative paths have no function, we will only
consider
the positive initial part as discussed above, and the negative end arrow
as
information, but never the whole path.

Choosing direct scepticism is a decision beyond the scope of this article,
and
we just make it. It is a general question how to treat contradictory and
absent
information, and if they are equivalent or not, see the remark in
Section \ref{Section 4.4} (page \pageref{Section 4.4}). (The fundamental
difference
between intersection of extensions and direct scepticism for defeasible
inheritance was shown in  \cite{Sch93}.)
See also Section \ref{Section Review} (page \pageref{Section Review})  for more
discussion.

We turn now to the announced variants as well as a finer distinction
within the
directly sceptical approach.
Again, see also Section \ref{Section Review} (page \pageref{Section Review}) 
for more discussion.

Our approach generates another problem, essentially that of the treatment
of a
mixture of contradictory and concordant information of multiple strengths
or
truth values. We bundle the decision of this problem with that for direct
scepticism into a ``plug-in'' decision, which will be used in three
approaches:
the conceptual ideas, the inheritance algorithm, and the choice of the
reference
class for subset size (and implicitly also for the treatment as a
prototype
theory). It is thus well encapsulated, and independent from the context.

These decisions (but, perhaps to a lesser degree, (1)) concern a wider
subject
than only inheritance
networks. Thus, it is not surprising that there are different formalisms
for
solving such networks, deciding one way or the other. But this multitude
is not
the fault of inheritance theory, it is only a symptom of a deeper
question. We
first give an overview for a clearer overall picture,
and discuss them in detail below, as they involve sometimes quite subtle
questions.

(1) Upward chaining against downward or double chaining.

(2.1) Off-path against on-path preclusion.

(2.2) Split validity preclusion against total validity preclusion.

(3) Direct scepticism against intersection of extensions.

(4) Treatment of mixed contradiction and preclusion situations,
no preclusion by paths of the same polarity.

(1) When we interpret arrows as causation (in the sense that $X \xcp Y$
expresses
that condition $X$ usually causes condition $Y$ to result), this can also
be seen as
a difference in reasoning from cause to effect vs. backward reasoning,
looking
for causes for an effect. (A word of warning: There is a
well-known article  \cite{SL89} from which a superficial reader
might conclude that
upward chaining is tractable, and downward chaining is not. A more careful
reading reveals that, on the negative side, the authors only show that
double
chaining is not tractable.) We will adopt upward chaining in all our
approaches.
See Section \ref{Section 4.4} (page \pageref{Section 4.4})  for more remarks.

(2.1) and (2.2) Both are consequences of our view - to be discussed below
in
Section 4.3 - to see valid paths also as
an absolute comparison of truth values, independent of reachability of
information. Thus, in Diagram \ref{Diagram Tweety} (page \pageref{Diagram
Tweety}), the
comparison between the truth values
``penguin'' and ``bird'' is absolute, and does not depend on the point of view
``Tweety'', as it can in total validity preclusion - if we continue to
view preclusion as a comparison of information strength (or truth value).
This question of absoluteness transcends obviously inheritance
networks. Our decision is, of course, again uniform for all our
approaches.

(3) This point, too, is much more general than the problems of
inheritance.
It is, among other things, a question of whether only the two possible
cases
(positive and negative) may hold, or whether there might be still other
possibilities. See Section \ref{Section 4.4} (page \pageref{Section 4.4}).

(4) This concerns the treatment of truth values in more complicated
situations, where we have a mixture of agreeing and contradictory
information.
Again, this problem reaches far beyond inheritance networks.

We will group (3) and (4) together in one general, ``plug-in''
decision, to be found in all approaches we discuss.

\bd

$\hspace{0.01em}$

% (+++ Orig. No.:  Definition IS-2.2 +++)

\label{Definition IS-2.2}

This is an informal definition of a plug-in decision:

We describe now more precisely a situation which we will meet in all
contexts
discussed, and
whose decision goes beyond our problem - thus, we have to adopt one or
several
alternatives, and translate them into the approaches we will discuss.
There will be one global decision, which is (and can be) adapted to the
different contexts.

Suppose we have information about $ \xbf $ and $ \xbq,$ where $ \xbf $
and $ \xbq $ are presumed to be
independent - in some adequate sense.

Suppose then that we have information sources $A_{i}:i \xbe I$ and
$B_{j}:j \xbe J,$ where
the $A_{i}$ speak about $ \xbf $ (they say $ \xbf $ or $ \xCN \xbf ),$ and
the $B_{j}$ speak about $ \xbq $ in the
same way. Suppose further that we have a partial, not necessarily
transitive (!),
ordering $<$ on the information sources $A_{i}$ and $B_{j}$ together.
$X<Y$ will say that $X$ is
better (intuition: more specific) than $Y.$ (The potential lack of
transitivity
is crucial, as valid paths do not always concatenate to valid paths - just
consider the Tweety diagram.)

We also assume that there are contradictions, i.e. some $A_{i}$ say $ \xbf
,$ some $ \xCN \xbf,$
likewise for the $B_{j}$ - otherwise, there are no problems in our
context.

We can now take several approaches, all taking contradictions and the
order $<$
into account.
 \xEI
 \xDH (P1) We use the global relation $<,$ and throw away all information
coming from
sources of minor quality, i.e. if there is $X$ such that $X<Y,$ then no
information
coming from $Y$ will be taken into account. Consequently, if $Y$ is the
only source
of information about $ \xbf,$ then we will have no information about $
\xbf.$ This seems
an overly radical approach, as one source might be better for $ \xbf,$
but not
necessarily for $ \xbq,$ too.

If we adopt this radical approach, we can continue as below, and can even
split
in analogue ways into (P1.1) and (P1.2), as we do below for (P2.1) and
(P2.2).

 \xDH (P2) We consider the information about $ \xbf $ separately from the
information about
$ \xbq.$ Thus, we consider for $ \xbf $ only the $A_{i},$ for $ \xbq $
only the $B_{j}.$ Take now e.g. $ \xbf $
and the $A_{i}.$ Again, there are (at least) two alternatives.
 \xEI
 \xDH (P2.1) We eliminate again all sources among the $A_{i}$ for which
there is a better
$A_{i' },$ irrespective of whether they agree on $ \xbf $ or not.
 \xEI
 \xDH (a) If the sources left are contradictory, we conclude nothing about
$ \xbf,$ and
accept for $ \xbf $ none of the sources. (This is a directly sceptical
approach of
treating unsolvable contradictions, following our general strategy.)

 \xDH (b) If the sources left agree for $ \xbf,$ i.e. all say $ \xbf,$
or all say $ \xCN \xbf,$ then we
conclude $ \xbf $ (or $ \xCN \xbf ),$ and accept for $ \xbf $ all the
remaining sources.
 \xEJ
 \xDH (P2.2) We eliminate again all sources among the $A_{i}$ for which
there is a better
$A_{i' }$, but only if $A_{i}$ and $A_{i' }$ have contradictory
information.
Thus, more sources may survive than in approach (P2.1).

We now continue as for (P2.1):
 \xEI
 \xDH (a) If the sources left are contradictory, we conclude nothing about
$ \xbf,$ and
accept for $ \xbf $ none of the sources.

 \xDH (b) If the sources left agree for $ \xbf,$ i.e. all say $ \xbf,$
or all say $ \xCN \xbf,$ then we
conclude $ \xbf $ (or $ \xCN \xbf ),$ and accept for $ \xbf $ all the
remaining sources.
 \xEJ
 \xEJ
 \xEJ
The difference between (P2.1) and (P2.2) is illustrated by the following
simple
example. Let $A<A' <A'',$ but $A \xEc A'' $ (recall that $<$ is not
necessarily transitive),
and $A \xcm \xbf,$ $A' \xcm \xCN \xbf,$ $A'' \xcm \xCN \xbf.$ Then
(P2.1) decides for $ \xbf $ $( \xCB $ is the only
survivor), (P2.2) does not decide, as $ \xCB $ and $ \xCB '' $ are
contradictory, and both
survive in (P2.2).

There are arguments for and against either solution: (P2.1) gives a
uniform
picture, more independent from $ \xbf,$ (P2.2) gives more weight to
independent
sources, it ``adds'' information sources, and thus gives potentially more
weight
to information from several sources. (P2.2) seems more in the tradition of
inheritance networks, so we will consider it in the further development.

The reader should note that our approach is quite far from a fixed point
approach in two ways: First, fixed point approaches seem more appropriate
for extensions-based approaches, as both try to collect a maximal set of
uncontradictory information. Second, we eliminate information when there
is better, contradicting information, even if the final result agrees with
the first. This, too, contradicts in spirit the fixed point approach.

\ed

After these preparations, we turn to a formal definition of validity of
paths.

\paragraph{
The definition of $ \xcm $ (i.e. of validity of paths)
}

$\hspace{0.01em}$

% (+++ Orig.:  The definition of  m (i.e. of validity of paths) +++)

\label{Section The definition of  m (i.e. of validity of paths)}

All definitions are relative to a fixed diagram $ \xbG.$
The notion of degree will be defined relative to all nodes of $ \xbG,$ as
we will
work with split validity preclusion, so the paths to consider may have
different
origins. For simplicity, we consider $ \xbG $ to be just a set of points
and
arrows, thus e.g. $x \xcp y \xbe \xbG $ and $x \xbe \xbG $ are defined,
when $x$ is a point in $ \xbG,$ and
$x \xcp y$ an arrow in $ \xbG.$ Recall that we have two types of arrows,
positive and
negative ones.

We first define generalized and potential paths, then the notion of
degree, and
finally validity of paths, written $ \xbG \xcm \xbs,$ if $ \xbs $ is a
path, as well as $ \xbG \xcm xy,$
if $ \xbG \xcm \xbs $ and $ \xbs:x \Xl. \xcp y.$

\bd

$\hspace{0.01em}$

% (+++ Orig. No.:  Definition 2.3 +++)

\label{Definition 2.3}

(1) Generalized paths:

A generalized path is an uninterrupted chain of positive or negative
arrows
pointing in the same direction, more precisely:

$x \xcp p \xbe \xbG $ $ \xcp $ $x \xcp p$ is a generalized path,

$x \xcP p \xbe \xbG $ $ \xcp $ $x \xcP p$ is a generalized path.

If $x \xFB \xcp p$ is a generalized path, and $p \xcp q \xbe \xbG $,
then
$x \xFB \xcp p \xcp q$ is a generalized path,

if $x \xFB \xcp p$ is a generalized path, and $p \xcP q \xbe \xbG $,
then
$x \xFB \xcp p \xcP q$ is a generalized path.

(2) Concatenation:

If $ \xbs $ and $ \xbt $ are two generalized paths, and the end point of $
\xbs $ is the same
as the starting point of $ \xbt,$ then $ \xbs \xDM \xbt $ is the
concatenation of $ \xbs $ and $ \xbt.$

(3) Potential paths (pp.):

A generalized path, which contains at most one negative arrow, and this at
the
end, is a potential path. If the last link is positive, it is a positive
potential path, if not, a negative one.

(4) Degree:

As already indicated, we shall define paths inductively. As we do not
admit cycles in our systems, the arrows define a well-founded relation
on the vertices. Instead of using this relation for
the induction, we shall first define the auxiliary notion of
degree, and do induction on the degree.
Given a node $x$ (the origin), we need a (partial) mapping $f$ from the
vertices to
natural numbers such that $p \xcp q$ or $p \xcP q$ $ \xbe $ $ \xbG $
implies $f(p)<f(q),$ and define (relative
to $x):$

Let $ \xbs $ be a generalized path from $x$ to $y,$ then $deg_{ \xbG,x}(
\xbs ):=deg_{ \xbG,x}(y):=$
the maximal length of any generalized path parallel to $ \xbs,$ i.e.
beginning in $x$ and ending in $y.$

\ed

\bd

$\hspace{0.01em}$

% (+++ Orig. No.:  Definition 2.4 +++)

\label{Definition 2.4}

\label{Definition ValidPaths}

Inductive definition of $ \xbG \xcm \xbs:$

Let $ \xbs $ be a potential path.
 \xEI
 \xDH Case $I$:

$ \xbs $ is a direct link in $ \xbG.$ Then $ \xbG \xcm \xbs $

(Recall that we have no hard contradictions in $ \xbG.)$

 \xDH Case II:

$ \xbs $ is a compound potential path, $deg_{ \xbG,a}( \xbs )=n,$ and $
\xbG \xcm \xbt $ is defined
for all $ \xbt $ with degree less than $n$ - whatever their origin and
endpoint.

 \xDH Case II.1:

Let $ \xbs $ be a positive pp. $x \xFB \xcp u \xcp y,$ let $ \xbs ':=x
\xFB \xcp u,$ so $ \xbs = \xbs ' \xDM u \xcp y$

Then, informally, $ \xbG \xcm \xbs $ iff
 \xEh
 \xDH (1) $ \xbs $ is a candidate by upward chaining,
 \xDH (2) $ \xbs $ is not precluded by more specific contradicting
information,
 \xDH (3) all potential contradictions are themselves precluded by
information
contradicting them.
 \xEj
Note that (2) and (3) are the translation of (P2.2)
in Definition \ref{Definition IS-2.2} (page \pageref{Definition IS-2.2}).

Formally, $ \xbG \xcm \xbs $ iff
 \xEh
 \xDH (1) $ \xbG \xcm \xbs ' $ and $u \xcp y \xbe \xbG.$

(The initial segment must be a path, as we have an upward chaining
approach.
This is decided by the induction hypothesis.)

 \xDH (2) There are no $v,$ $ \xbt,$ $ \xbt ' $ such that $v \xcP y \xbe
\xbG $ and $ \xbG \xcm \xbt:=x \xFB \xcp v$ and
$ \xbG \xcm \xbt ':=v \xFB \xcp u.$ $( \xbt $ may be the empty path, i.e.
$x=v.)$

$( \xbs $ itself is not precluded by split validity preclusion and a
contradictory
link. Note that $ \xbt \xDM v \xcP y$ need not be valid, it suffices
that it is a better candidate (by $ \xbt ' ).)$

 \xDH (3) all potentially conflicting paths are precluded by information
contradicting them:

For all $v$ and $ \xbt $ such that $v \xcP y \xbe \xbG $ and $ \xbG \xcm
\xbt:=x \xFB \xcp v$ (i.e. for all potentially
conflicting paths $ \xbt \xDM v \xcP y)$ there is $z$ such that $z \xcp y
\xbe \xbG $ and either

$z=x$

(the potentially conflicting pp. is itself precluded by a direct link,
which is
thus valid)

or

there are $ \xbG \xcm \xbr:=x \xFB \xcp z$ and $ \xbG \xcm \xbr ':=z
\xFB \xcp v$ for suitable $ \xbr $ and $ \xbr '.$
 \xEj
 \xDH Case II.2: The negative case, i.e.
$ \xbs $ a negative pp. $x \xFB \xcp u \xcP y,$ $ \xbs ':=x \xFB \xcp u,$
$ \xbs = \xbs ' \xDM u \xcP y$
is entirely symmetrical.
 \xEJ

\ed

\br

$\hspace{0.01em}$

% (+++ Orig. No.:  Remark 2.1 +++)

\label{Remark 2.1}

The following remarks all concern preclusion.

(1) Thus, in the case of preclusion, there is a valid path from $x$ to
$z,$
and $z$ is more specific than $v,$ so
$ \xbt \xDM v \xcP y$ is precluded. Again, $ \xbr \xDM z \xcp y$ need not
be a valid path, but it is
a better candidate than $ \xbt \xDM v \xcP y$ is, and as $ \xbt \xDM v
\xcP y$ is in simple contradiction,
this suffices.

(2) Our definition is stricter than many popular ones, in the
following sense: We require - according to our general picture to treat
only
direct links as information - that the preclusion ``hits'' the precluded
path
at the end, i.e. $v \xcP y \xbe \xbG,$ and $ \xbr ' $ hits $ \xbt \xDM v
\xcP y$ at $v.$ In other definitions,
it is possible that the preclusion hits at some $v',$ which is somewhere
on the
path $ \xbt,$ and not necessarily at its end. For instance, in the Tweety
Diagram,
see Diagram \ref{Diagram Tweety} (page \pageref{Diagram Tweety}), if
there were a node $b' $ between $b$ and $d,$ we will need the
path $c \xcp b \xcp b' $ to be valid, (obvious) validity of the arrow $c
\xcp b$ will not
suffice.

(3) If we allow $ \xbr $ to be the empty path, then the case $z=x$ is a
subcase of the
present one.

(4) Our conceptual analysis has led to a very important simplification of
the
definition of validity. If we adopt on-path preclusion, we have to
remember
all paths which led to the information source to be considered: In the
Tweety
diagram, we have to remember that there is an arrow $a \xcp b,$ it is not
sufficient
to note that we somehow came from a to $b$ by a valid path, as the path $a
\xcp c \xcp b \xcp d$
is precluded, but not the path $a \xcp b \xcp d.$ If we adopt total
validity preclusion,
see also Section \ref{Section Review} (page \pageref{Section Review})  for more
discussion, we
have to remember the valid path $a \xcp c \xcp b$ to see that it precludes
$a \xcp c \xcp d.$ If
we allow preclusion to ``hit'' below the last node, we also have to remember
the entire path which is precluded. Thus, in all those cases, whole paths
(which can be very long) have to be remembered, but NOT in our definition.

We only need to remember (consider the Tweety diagram):

(a) we want to know if $a \xcp b \xcp d$ is valid, so we have to remember
a, $b,$ $d.$
Note that the (valid) path from a to $b$ can be composed and very long.

(b) we look at possible preclusions, so we have to remember $a \xcp c \xcP
d,$ again
the (valid) path from a to $c$ can be very long.

(c) we have to remember that the path from $c$ to $b$ is valid (this was
decided
by induction before).

So in all cases (the last one is even simpler), we need only remember the
starting node, a (or $c),$ the last node of the valid paths, $b$ (or $c),$
and the
information $b \xcp d$ or $c \xcP d$ - i.e. the size of what has to be
recalled is
$ \xck 3.$ (Of course, there may be many possible preclusions, but in all
cases we have to look at a very limited situation, and not arbitrarily
long
paths.)

We take a fast look forward to Section 4.3, where we describe diagrams as
information and its transfer, and nodes also as truth values. In these
terms -
and the reader is asked to excuse the digression - we may note above point
(a) as $a \xch_{b}d$ - expressing that, seen from a, $d$ holds with truth
value $b,$
(b) as $a \xch_{c} \xCN d,$ (c) as $c \xch_{c}b$ - and this is all we need
to know.

$ \xcz $
\\[3ex]

\er

We indicate here some modifications of the definition without discussion,
which
is to be found below.

(1) For on-path preclusion only:
Modify condition (2) in Case II.1 to:
$(2' )$ There is no $v$ on the path $ \xbs $ (i.e. $ \xbs:x \xFB \xcp v
\xFB \xcp u)$ such that $v \xcP y \xbe \xbG.$

(2) For total validity preclusion:
Modify condition (2) in Case II.1 to:
$(2' )$ There are no $v,$ $ \xbt,$ $ \xbt ' $ such that $v \xcP y \xbe
\xbG $ and $ \xbt:=x \xFB \xcp v$ and
$ \xbt ':=v \xFB \xcp u$ such that $ \xbG \xcm \xbt \xDM \xbt '.$

(3) For extension based approaches:
Modify condition (3) in Case II.1 as follows:
$(3' )$ If there are conflicting paths,
which are not precluded themselves by contradictory information, then we
branch recursively (i.e. for all such situations) into two extensions,
one,
where the positive non-precluded paths are valid, one, where the negative
non-precluded paths are valid.

\bd

$\hspace{0.01em}$

% (+++ Orig. No.:  Definition 2.5 +++)

\label{Definition 2.5}

Finally, define $ \xbG \xcm xy$ iff there is $ \xbs:x \xcp y$ s.th. $
\xbG \xcm \xbs,$ likewise for $x \ol{y}$ and
$ \xbs:x \xFB \xcP y.$

\ed

Diagram \ref{Diagram Complicated-Case} (page \pageref{Diagram Complicated-Case})
 shows
the most complicated situation for the positive case.
% Diagram 2.3
% Diagram 2.3
% -----------
%                                 y
%                                 y
%    @@                          EA  F
%    @@                         E A   F
%    @@                        E  A    F
%    @@                       E   D     F
%    @@                      E    A      F
%    @@                     E     A       F
%                          u      v@@caa!!!@z
%    @@                     F     A       E
%    @@                      F    A      E
%    @@                       F   A     E
%    @@                        )  ``    (
%    @@                         ) ''   (
%    @@                          )"  (
%                                 x

\vspace{10mm}

\begin{diagram}

\label{Diagram Complicated-Case}
\index{Diagram Complicated-Case}

\unitlength1.0mm
\begin{picture}(130,100)

\newsavebox{\Preclusion}
\savebox{\Preclusion}(140,90)[bl]
{

\multiput(43,8)(1,1){5}{\circle*{.3}}
\put(48,13){\vector(1,1){17}}

\multiput(37,8)(-1,1){5}{\circle*{.3}}
\put(32,13){\vector(-1,1){17}}

\put(13,38){\vector(1,1){24}}
\put(67,38){\vector(-1,1){24}}

\put(40,37){\vector(0,1){23}}

\multiput(40,8)(0,1){5}{\circle*{.3}}
\put(40,13){\vector(0,1){17}}

\multiput(66,34)(-1,0){5}{\circle*{.3}}
\put(61,34){\vector(-1,0){17}}

\put(39,3){$x$}
\put(9,33){$u$}
\put(39,33){$v$}
\put(69,33){$z$}
\put(39,63){$y$}

\put(38,50){\line(1,0){3.7}}

\put(10,0) {{\rm\bf The complicated case}}

}

\put(0,0){\usebox{\Preclusion}}
\end{picture}

\end{diagram}

\vspace{4mm}

We have to show now that the above approach corresponds to the preceeding
discussion.

\bfa

$\hspace{0.01em}$

% (+++ Orig. No.:  Fact 2.2 +++)

\label{Fact 2.2}

The above definition and the informal one outlined
in Definition \ref{Definition IS-2.2} (page \pageref{Definition IS-2.2}) 
correspond,
when we consider valid positive paths as access to information and
comparison
of information strength as indicated at the beginning
of Section \ref{Section IS-2.2} (page \pageref{Section IS-2.2})  and
elaborated in Section \ref{Section 4.3} (page \pageref{Section 4.3}).

\efa

\subparagraph{
Proof
}

$\hspace{0.01em}$

% (+++ Orig.:  Proof +++)

As Definition \ref{Definition IS-2.2} (page \pageref{Definition IS-2.2})  is
informal, this cannot be a formal proof, but it
is obvious how to transform it into one.

We argue for the result, the argument for valid paths is similar.

Consider then case (P2.2) in Definition \ref{Definition IS-2.2} (page
\pageref{Definition IS-2.2}),
and start from some $x.$

\paragraph{
Case 1:
}

$\hspace{0.01em}$

% (+++ Orig.:  Case 1: +++)

\label{Section Case 1:}

Direct links, $x \xcp z$ or $x \xcP z.$

By comparison of strength via preclusion, as a direct link starts at $x,$
the
information $z$ or $ \xCN z$ is
stronger than all other accessible information. Thus, the link and the
information will be valid in both approaches. Note that we assumed $ \xbG
$ free from
hard contradictions.

\paragraph{
Case 2:
}

$\hspace{0.01em}$

% (+++ Orig.:  Case 2: +++)

\label{Section Case 2:}

Composite paths.

In both approaches, the initial segment has to be valid, as information
will
otherwise not be accessible. Also, in both approaches, information will
have
the form of direct links from the accessible source. Thus, condition (1)
in
Case II.1 corresponds to condition (1)
in Definition \ref{Definition IS-2.2} (page \pageref{Definition IS-2.2}).

In both approaches, information contradicted by a stronger source
(preclusion)
is discarded, as well as information which is contradicted by other, not
precluded sources, so (P2.2)
in Definition \ref{Definition IS-2.2} (page \pageref{Definition IS-2.2})  and
II.1
$(2)+(3)$ correspond. Note that variant (P2.1)
of Definition \ref{Definition IS-2.2} (page \pageref{Definition IS-2.2})  would
give a
different result - which we could, of course, also imitate in a modified
inheritance approach.

\paragraph{
Case 3:
}

$\hspace{0.01em}$

% (+++ Orig.:  Case 3: +++)

\label{Section Case 3:}

Other information.

Inheritance nets give no other information, as valid information is
deduced
only through valid paths
by Definition \ref{Definition 2.5} (page \pageref{Definition 2.5}). And we did
not add any other
information either in the approach
in Definition \ref{Definition IS-2.2} (page \pageref{Definition IS-2.2}). But
as is obvious in
Case 2, valid paths coincide in both cases.

Thus, both approaches are equivalent.

$ \xcz $
\\[3ex]
%  2.3  Review of other approaches and problems
%  2.3  Review of other approaches and problems
% #
% =============================================
\subsection{
Review of other approaches and problems
}
\label{Section IS-2.3}
\label{Section Review}

We now discuss shortly in more detail some of the differences between
various
major definitions of inheritance formalisms.

Diagram 6.8, $p.$ 179, in  \cite{Sch97-2} (which is probably due to
folklore of the
field) shows requiring downward chaining would be wrong. We repeat it
here,
see Diagram \ref{Diagram Up-Down-Chaining} (page \pageref{Diagram
Up-Down-Chaining}).
% Diagram 2.4
% Diagram 2.4
% -----------
% The problem of downward chaining:
% The problem of downward chaining:
% @@         %y
% @@         %y
% @@
% @@       E    F
% @@      E      f
% @@     E        F
% @@    E          F
% @@
% @@ %v caaaaaaaaaaa %x
% @@
% @@    F          E  A
% @@     F        E   A
% @@      F      E    A
% @@       F    E     A
% @@                  D
% @@         %u       A
% @@                  A
% @@          A       A
% @@          A       A
% @@                  A
% @@         %z aaaaaah

\vspace{10mm}

\begin{diagram}

\label{Diagram Up-Down-Chaining}
\index{Diagram Up-Down-Chaining}

\unitlength1.0mm
\begin{picture}(130,100)

\newsavebox{\ZWEIvier}
\savebox{\ZWEIvier}(140,110)[bl]
{

\put(0,95){{\rm\bf The problem of downward chaining}}

\put(43,27){\vector(1,1){24}}
\put(37,27){\vector(-1,1){24}}
\put(13,57){\vector(1,1){24}}
\put(67,57){\vector(-1,1){24}}

\put(53,67){\line(1,1){4}}

\put(67,54){\vector(-1,0){54}}

\put(40,7){\vector(0,1){14}}
\put(43,7){\line(3,5){24}}
\put(58,28.1){\line(-5,3){3.6}}

\put(39,3){$z$}
\put(39,23){$u$}
\put(9,53){$v$}
\put(69,53){$x$}
\put(39,83){$y$}

}

\put(0,0){\usebox{\ZWEIvier}}
\end{picture}

\end{diagram}

\vspace{4mm}

Preclusions valid above (here at $u)$ can be invalid at lower points (here
at $z),$
as part of the relevant information is no longer accessible (or becomes
accessible). We have $u \xcp x \xcP y$ valid, by downward chaining, any
valid path
$z \xcp u \Xl.y$ has to have a valid final segment u \Xl y, which can
only be $u \xcp x \xcP y,$
but intuition says that $z \xcp u \xcp v \xcp y$ should be valid. Downward
chaining prevents such changes, and thus seems inadequate, so we decide
for
upward chaining. (Already preclusion itself underlines upward chaining: In
the
Tweety diagram, we have to know that the path from bottom up to penguins
is
valid. So at least some initial subpaths have to be known - we need upward
chaining.) (The rejection of downward chaining seems at first sight to be
contrary to the intuitions carried by the word $'' inheritance''.)$
See also the remarks in Section \ref{Section 4.4} (page \pageref{Section 4.4}).

\paragraph{
Extension-based versus directly skeptical definitions
}

$\hspace{0.01em}$

% (+++ Orig.:  Extension-based versus directly skeptical definitions +++)

\label{Section Extension-based versus directly skeptical definitions}

As this distinction has already received detailed discussion in the
literature, we shall be very brief here.
An extension of a net is essentially a maximally consistent and in some
appropriate sense reasonable subset of all its potential paths. This can
of course be presented either as a liberal conception (focussing on
individual extensions) or as a skeptical one (focussing on their
intersection - or, the
intersection of their conclusion sets). The seminal presentation is that
of [Tou86], as refined by [San86]. The directly skeptical approach seeks
to
obtain a notion of skeptically accepted path and conclusion, but without
detouring through extensions. Its classic presentation is that of [HTT87].
Even while still searching for fully adequate definitions of either kind,
we may use the former approach as a useful ``control'' on the latter. For
if we can find an intuitively possible and reasonable extension supporting
a conclusion $x \ol{y},$ whilst a proposed definition for a directly
skeptical
notion of legitimate inference yields xy as a conclusion, then the
counterexemplary extension seems to call into question the adequacy
of the directly skeptical construction, more readily than inversely.
It has been shown in  \cite{Sch93} that the intersection of
extensions is
fundamentally different from the directly sceptical approach.
See also the remark in Section \ref{Section 4.4} (page \pageref{Section 4.4}).

From now on, all definitions considered shall be (at least) upward
chaining.

\paragraph{
On-path versus off-path preclusion
}

$\hspace{0.01em}$

% (+++ Orig.:  On-path versus off-path preclusion +++)

\label{Section On-path versus off-path preclusion}

This is a rather technical distinction, discussed in [THT87]. Briefly, a
path
$ \xbs $: $x \xcp  \Xl  \xcp y \xcp  \Xl  \xcp z$ and a direct link $y
\xcP u$ is an off-path preclusion of
$ \xbt $: $x \xcp  \Xl  \xcp z \xcp  \Xl  \xcp u$, but an on-path
preclusion only iff all nodes of $ \xbt $
between $x$ and $z$ lie on the path $ \xbs.$

For instance, in the Tweety diagram, the arrow $c \xcP d$ is an on-path
preclusion
of the path $a \xcp c \xcp b \xcp d,$ but the paths $a \xcp c$ and $c \xcp
b,$ together with $c \xcP d,$ is
an (split validity) off-path preclusion of the path $a \xcp b \xcp d.$

\paragraph{
Split validity versus total validity preclusion
}

$\hspace{0.01em}$

% (+++ Orig.:  Split validity versus total validity preclusion +++)

\label{Section Split validity versus total validity preclusion}

Consider again a preclusion $ \xbs:u \xcp  \Xl  \xcp x \xcp  \Xl  \xcp
v,$ and $x \xcP y$ of
$ \xbt:u \xcp  \Xl  \xcp v \xcp  \Xl  \xcp y.$ Most definitions demand
for the preclusion to be
effective - i.e. to prevent $ \xbt $ from being accepted - that the
total path $ \xbs $ is valid. Some ([GV89], [KK89], [KKW89a], [KKW89b])
content
themselves with the combinatorially simpler separate (split) validity of
the lower and upper parts of $ \xbs $: $ \xbs ':u \xcp  \Xl  \xcp x$ and
$ \xbs '':x \xcp  \Xl  \xcp v.$
In Diagram \ref{Diagram Split-Total-Preclusion} (page \pageref{Diagram
Split-Total-Preclusion}), taken
from  \cite{Sch97-2}, the path $x \xcp w \xcp v$ is valid,
so is $u \xcp x,$ but not the whole preclusion path $u \xcp x \xcp w \xcp
v.$
% Diagram 2.5
% Diagram 2.5
% -----------
% Split vs. total validity preclusion:
% Split vs. total validity preclusion:
% @@          %y
% @@          %y
% @@
% @@        E    F
% @@       E      f
% @@      E        F
% @@     E          F
% @@
% @@  %v caaa %w caaa %x
% @@
% @@     F     A    E
% @@      F    D   E
% @@       F   A  E
% @@        F  A E
% @@
% @@          %u
% @@

\vspace{10mm}

\begin{diagram}

\label{Diagram Split-Total-Preclusion}
\index{Diagram Split-Total-Preclusion}

\unitlength1.0mm
\begin{picture}(130,100)

\newsavebox{\SECHSneun}
\savebox{\SECHSneun}(140,90)[bl]
{

\put(0,75){{\rm\bf Split vs. total validity preclusion}}

\put(43,8){\vector(1,1){24}}
\put(37,8){\vector(-1,1){24}}
\put(13,38){\vector(1,1){24}}
\put(67,38){\vector(-1,1){24}}

\put(40,8){\vector(0,1){23}}
\put(66,34){\vector(-1,0){23}}
\put(36,34){\vector(-1,0){23}}

\put(39,3){$u$}
\put(9,33){$v$}
\put(39,33){$w$}
\put(69,33){$x$}
\put(39,63){$y$}

\put(38,20){\line(1,0){3.7}}
\put(52,49){\line(1,1){4}}

}

\put(0,0){\usebox{\SECHSneun}}
\end{picture}

\end{diagram}

\vspace{4mm}

Thus, split validity
preclusion will give here the definite result $u \ol{y}.$ With total
validity
preclusion, the diagram has essentially the form of a Nixon Diamond.
\section{
Defeasible inheritance and reactive diagrams
}
\label{Section Def-Reac}
\index{Section Def-Reac}

Before we discuss the relationship in detail, we first summarize our
algorithm.
\subsection{
Summary of our algorithm
}

We look for valid paths from $x$ to $y.$

(1) Direct arrows are valid paths.

(2) Consider the set $C$ of all direct predecessors of $y,$ i.e. all $c$
such that
there is a direct link from $c$ to $y.$

(2.1) Eliminate all those to which there is no valid positive path from
$x$
(found by induction), let the new set be $C' \xcc C.$

If the existence of a valid path has not yet been decided, we have to
wait.

(2.2) Eliminate from $C' $ all $c$ such that there is $c' \xbe C' $ and a
valid positive
path from $c' $ to $c$ (found by induction) - unless the arrows from $c$
and from $c' $ to
$y$ are of the same type. Let the new set be $C'' \xcc C' $ (this handles
preclusion).

If the existence of such valid paths has not yet been decided,
we have to wait.

(2.3) If the arrows from all elements of $C'' $ to $y$ have same polarity,
we have
a valid path from $x$ to $y,$ if not, we have an unresolved contradiction,
and
there is no such path.

Note that we were a bit sloppy here. It can be debated whether preclusion
by some $c' $ such that $c$ and $c' $ have the same type of arrow to $y$
should be
accepted. As we are basically only interested whether there is a valid
path,
but not in its details, this does not matter.
\subsection{
Overview
}

There are several ways to use reactive graphs to help us solve
inheritance diagrams - but also to go beyond them.

 \xEh

 \xDH We use them to record results found by adding suitable arrows to the
graph.

 \xDH We go deeper into the calculating mechanism, and do not only use new
arrows as memory, but also for calculation.

 \xDH We can put up ``sign posts'' to mark dead ends, in the following
sense:
If we memorize valid paths from $x$ to $y,$ then, anytime we are at a
branching
point $u$ coming from $x,$ trying to go to $y,$ and there is valid path
through an
arrow leaving $u,$ we can put up a sign post
saying ``no valid path from $x$ to $y$ through this arrow''.

Note that we have to state destination $y$ (of course), but also outset,
$x:$
There might be a valid path from $u$ to $y,$ which may be precluded or
contradicted
by some path coming from $x.$

 \xDH We can remember preclusion, in the following sense: If we found a
valid positive path from a to $b,$ and there are contradicting arrows from
a
and $b$ to $c,$ then we can create an arrow from a to the arrow from $b$
to $c.$
So, if, from $x,$ we can reach both a and $b,$ the arrow from $b$ to $c$
will be
incapacitated.

 \xEj

Before we discuss the first three possibilities in detail we shortly
discuss the
more general picture (in rough outline).

 \xEh

 \xDH Replacing labels by arrows and vice versa.

As we can switch arrows on and off, an arrow carries a binary value -
even without nay labels. So the idea is obvious:

If an arrow has 1 label with $n$ possible values, we can replace it
with $n$ parallel arrows (i.e. same source and destination), where we
switch
eaxctly one of them on - this is the label.

Conversely, we can replace $n$ parallel arrows without labels, where
exactly one
is active, by one arrow with $n$ labels.

We can also code labels of a node $x$ by an arrow $ \xba:x \xcp x,$ which
has the
same labels.

 \xDH Coding logical formulas and truth in a model

We take two arrows for each propositional variable, one stands for true,
the other for false. Negation blocks the positive arrow, enables the
negative one. Conjunction is solved by concatenation, disjunction by
parallel paths. If a variable occurs more than once, we make copies,
which are ``switched'' by the ``master arrow''.

 \xEj

We come back to the first three ways to treat inheritance by reactive
graphs,
and also mention a way to go beyond usual inheritance.
\subsection{
Compilation and memorization
}

When we take a look at the algorithm deciding which potential paths are
valid,
we see that, with one exception, we only need the results already
obtained,
i.e. whether there is a valid positive/negative path from a to $b,$
and not the actual paths themselves. (This is, of course, due to our
split validity definition of preclusion.) The exception is that preclusion
works ``in the upper part'' with direct links. But this is a local problem:
we only have to look at the direct predecessors of a node.

Consequently, we can do most of the calculation just once, in an induction
of increasing ``path distance'', memorize valid positive (the negative ones
cannot be used) paths with special arrows which we activated once their
validity is established, and work now as follows with the new arrows:

Suppose we want to know if there is a valid path from $x$ to $y.$ We look
backwards
at all predecessors $b$ of $y$ (using a simple backward pointer), and look
whether
there is a valid positive path from $x$ to $b,$ using the new arrows.
We then look at all arrows going from such $b$ to $y.$ If they all agree
(i.e. all
are positive, or all are negative), we need not look further, and have the
result. If they disagree, we have to look at possible comparisons by
specificity. For this, we see whether there are new arrows between the $b'
$s.
All such $b$ to which there is a new arrow from another $b$ are out of
consideration. If the remaining agree, we have a valid path (and activate
a
new arrow from $x$ to $y$ if the path is positive), if not, there is no
such
path. (Depending on the technical details of the induction, it can be
useful to note this also by activating a special arrow.)
\subsection{
Executing the algorithm
}

Consider any two points $x,$ $y.$

There can be no path, a positive potential path, a negative potential
path,
both potential paths, a valid positive path, a valid negative path
(from $x$ to $y).$

Once a valid path is found, we can forget potential paths, so we can code
the possibilities by $\{*,p+,p-,p+-,v+,v-\},$ in above order.

We saw that we can work either with labels, or with a multitude of
parallel
arrows, we choose the first possibility.

We create for each pair of nodes a new arrow, $(x,y),$ which we intialize
with
label $*.$

First, we look for potential paths.

If there is a direct link from $x$ to $y,$ we change the value $*$ of
$(x,y)$ directly
to $v+$ or $v-.$

If $(x,y' )$ has value $p+$ or $p+-$ or $v+,$ and there is a direct link
from $y' $ to $y,$
we change the value of $(x,y)$ from $*$ to $p+$ or $p-,$ depending on the
link from
$y' $ to $y,$ from $p+$ or p- to $p+-$ if adequate (we found both
possibilities),
and leave the value unchanged otherwise.

This determines all potential paths.

We make for each pair $x,y$ at $y$ a list of its predecessors, i.e. of all
$c$
s.t. there is a direct link from $c$ to $y.$ We do this for all $x,$ so we
can
work in parallel. A list is, of course, a concatenation of suitable
arrows.

Suppose we want to know if there is valid path from $x$ to $y.$

First, there might be a direct link, and we are done.

Next, we look at the list of predecessors of $y$ for $x.$

If one $c$ in the list has value $*,$ $p-,$ v- for $(x,c),$ it is
eliminated from
the list. If one $c$ has $p+$ or $p+-,$ we have to wait. We do this until
all $(x,c)$
have $*,$ $v+$ or $v-,$ so those remaining in the list will have $v+.$

We look at all pairs in the list. While at least one $(c,c' )$ has $p+$ or
$p+-,$
we have to wait. Finally, all pairs will have $*,$ $p-,$ $v+$ or $v-.$
Eliminate
all $c' $ s.t. there is $(c,c' )$ with value $v+$ - unless the arrows from
$c$ and $c' $
to $y$ are of the same type.

Finally, we look at the list of the remaining predecessors, if they all
have the same link to $y,$ we set $(x,y)$ to $v+$ or $v-,$ otherwise to
$*.$

All such operations can be done by suitable operations with arrows,
but it is very lenghty and tedious to write down the details.
\subsection{
Signposts
}

Putting up signposts requires memorizing all valid paths, as leaving
one valid path does not necessarily mean that there is no alternative
valid path. The easiest thing to do this is probably to put up a warning
post
everywhere, and collect the wrong ones going backwards through the valid
paths.

We illustrate this with Diagram \ref{Diagram InherUniv} (page \pageref{Diagram
InherUniv}) :

\vspace{10mm}

\begin{diagram}

\label{Diagram InherUniv}
\index{Diagram InherUniv}

\centering
\setlength{\unitlength}{1mm}
{\renewcommand{\dashlinestretch}{30}
\begin{picture}(150,150)(0,0)

\put(10,60){\circle*{1}}
\put(70,60){\circle*{1}}
\put(130,60){\circle*{1}}

\put(70,10){\circle*{1}}
\put(70,110){\circle*{1}}

\put(40,50){\circle*{1}}
\put(40,70){\circle*{1}}
\put(100,50){\circle*{1}}
\put(100,70){\circle*{1}}

\path(69,11)(12,58)
\path(13.5,55.4)(12,58)(14.8,57)
\path(71,11)(128,58)
\path(125.2,57)(128,58)(126.5,55.4)
\path(12,62)(69,109)
\path(66.2,108)(69,109)(67.5,106.4)
\path(128,62)(71,109)
\path(72.5,106.4)(71,109)(73.8,108)

\path(70,61)(70,108)
\path(68,86)(72,86)
\path(69,105.2)(70,108)(71,105.2)

\path(11,60.5)(39,69.5)
\path(14,60.4)(11,60.5)(13.4,62.3)
\path(11,59.5)(39,50.5)
\path(13.4,57.7)(11,59.5)(14,59.6)
\path(41,69.5)(69,60.5)
\path(43.4,67.7)(41,69.5)(44,69.6)
\path(41,50.5)(69,59.5)
\path(44,50.4)(41,50.5)(43.4,52.3)

\path(24,67)(26,63)

\path(71,60.5)(99,69.5)
\path(71,59.5)(99,50.5)
\path(101,69.5)(129,60.5)
\path(101,50.5)(129,59.5)
\path(84,67)(86,63)

\path(74,60.4)(71,60.5)(73.4,62.3)
\path(73.4,57.7)(71,59.5)(74,59.6)
\path(103.4,67.7)(101,69.5)(104,69.6)
\path(104,50.4)(101,50.5)(103.4,52.3)

\path(100,52)(100,68)
\path(99.2,65.9)(100,68)(100.8,65.9)

\put(69,5){{\xssc $x$}}
\put(69,115){{\xssc $y$}}

\put(7,60){{\xssc $a$}}
\put(69,57){{\xssc $b$}}
\put(132,60){{\xssc $c$}}

\put(39,47){{\xssc $d$}}
\put(99,47){{\xssc $e$}}
\put(39,72){{\xssc $f$}}
\put(99,72){{\xssc $g$}}

\put(110,65){{\xssc $STOP$}}
\put(110,54){{\xssc $STOP$}}

% \put(20,3) {{\rm\bf $\xda-$ ranked structure and accessibility}}

\end{picture}
}

\end{diagram}

\vspace{4mm}

There are the following potential paths from $x$ to $y:$
xcy, xceb-y, xcebday, xay.

The paths xc, xa, and xceb are valid. The latter: xce is valid. xceb is
in competition with xcg-b and xceg-b, but both are precluded by the arrow
(valid path) eg.

y's predecessors on those paths are $a,b,c.$

a and $b$ are not comparable, as there is no valid path from $b$ to a, as
bda is contradicted by bf-a. None is more specific than the other one.

$b$ and $c$ are comparable, as the path ceb is valid, since cg-b and ceg-b
are
precluded by the valid path (arrow) eg. So $c$ is more specific than $b.$

Thus, $b$ is out of consideration, and we are left with a and $c.$ They
agree, so
there is positive valid path from $x$ to $y,$ more precisely, one through
a,
one through $c.$

We have put STOP signs on the arrows ce and cg, as we cannot continue via
them
to $y.$
\subsection{
Beyond inheritance
}

But we can also go beyond usual inheritance networks.

Consider the following scenario:

 \xEI

 \xDH Museum airplanes usually will not fly, but usual airplanes will fly.

 \xDH Penguins $don' t$ fly, but birds do.

 \xDH Non-flying birds usually have fat wings.

 \xDH Non-flying aircraft usually are rusty.

 \xDH But it is not true that usually non-flying things are rusty and have
fat wings.

 \xEJ

We can model this with higher order arrows as follows:

 \xEI

 \xDH Penguins $ \xcp $ birds, museum airplanes $ \xcp $ airplanes,
birds $ \xcp $ fly, airplanes $ \xcp $ fly.

 \xDH Penguins $ \xcP $ fly, museum airplanes $ \xcP $ fly.

 \xDH Flying objects $ \xcP $ rusty, flying objects $ \xcP $ have fat
wings.

 \xDH We allow concatenation of two negative arrows:

Coming e.g. from penguins, we want to concatenate penguins $ \xcP $ fly
and
fly $ \xcP $ fat wings. Coming from museum aircraft, we want to
concatenate
museum airfcraft $ \xcP $ fly and fly $ \xcP $ rusty.

We can enable this as follows: we introduce a new arrow
$ \xba:(penguin$ $ \xcP $ fly) $ \xcp $ (fly $ \xcP $ fat wings), which,
when traversing
penguin $ \xcP $ fly enables the algorithm to concatenate with the arrow
it
points to, using the rule $'' -*-=+'',$ giving the result that penguins
usually have fat wings.

 \xEJ

See  \cite{Gab08c} for deeper discussion.
\section{
Interpretations
}
\label{Section Is-Interpr}
\index{Section Is-Interpr}
%  4.1  Introduction
%  4.1  Introduction
% %
% ==================
\subsection{
Introduction
}
\label{Section 4.1}

We will discuss in this Section three interpretations of inheritance nets.

First, we will indicate fundamental differences between inheritance and
the
systems $P$ and $R.$ They will be elaborated
in Section \ref{Section Translation} (page \pageref{Section Translation}),
where an
interpretation in terms of small sets will be tried nonetheless, and its
limitations explored.

Second, we will interpret inheritance nets as systems of information and
information flow.

Third, we will interpret inheritance nets as systems of prototypes.

Inheritance nets present many intuitively attractive properties, thus it
is
not surprising that we can interpret them in several ways. Similarly,
preferential structures can be used as a semantics of deontic and of
nonmonotonic logic, they express a common idea: choosing a subset of
models
by a binary relation. Thus, such an ambiguity need not be a sign for a
basic
flaw.
%  4.2  Informal comparison of inheritance with the systems P and R
%  4.2  Informal comparison of inheritance with the systems P and R
% %
% =================================================================
\subsection{
Informal comparison of inheritance with the systems P and R
}
\label{Section 4.2}

\paragraph{
The main issues
}

$\hspace{0.01em}$

% (+++ Orig.:  The main issues +++)

\label{Section The main issues}

In the authors' opinion, the following two properties of inheritance
diagrams
show the deepest difference to preferential and similar semantics, and the
first even to classical logic. They have to be taken seriously, as they
are at
the core of inheritance systems, are independent of the particular
formalism,
and show that there is a fundamental difference between the former and the
latter. Consequently, any attempt at translation will have to stretch one
or
both sides perhaps beyond the breaking point.

(1) Relevance,

(2) subideal situations, or relative normality

Both (and more) can be illustrated by the following
simple Diagram \ref{Diagram Information-Transfer} (page \pageref{Diagram
Information-Transfer})
(which also shows conflict resolution by specificity).

(1) Relevance: As there is no monotonous path whatever between $e$ and
$d,$
the question whether e's are d's or not, or vice versa, does not even
arise.
For the same reason, there is no question whether b's are $c' $s, or not.
(As
a matter of fact, we will see below
in Fact \ref{Fact 5.1} (page \pageref{Fact 5.1})  that b's are non-c's in
system $P$ - see Definition \ref{Definition Log-Cond-Ref-Size} (page
\pageref{Definition Log-Cond-Ref-Size}) ).
In upward chaining formalisms, as there is no valid positive path from a
to $d,$
there is no question either whether a's are f's or not.

Of course, in classical logic, all information is relevant to the rest, so
we
can say e.g. that e's are $d' $s, or e's are $non-d' $s, or some are $d'
$s, some are
not, but there is a connection. As preferential models are based on
classical
logic, the same argument applies to them.

(2) In our diagram, a's are $b' $s, but not ideal $b' $s, as they are not
$d' $s, the
more specific information from $c$ wins. But they are $e' $s, as ideal b's
are.
So they are not perfectly ideal $b' $s, but as ideal b's as possible.
Thus, we
have graded ideality, which does not exist in preferential and similar
structures. In those structures, if an element is an ideal element, it has
all
properties of such,
if one such property is lacking, it is not ideal, and we $can' t$ say
anything any
more. Here, however, we sacrifice as little normality as possible, it is
thus a minimal change formalism.

In comparison, questions of information transfer and strength of
information
seem lesser differences. Already systems $P$ and $R$
(see Definition \ref{Definition Log-Cond-Ref-Size} (page \pageref{Definition
Log-Cond-Ref-Size}) )
differ on information transfer. In both cases, transfer is
based on the same notion of smallness, which describes ideal situations.
But,
as said in Remark \ref{Remark Ref-Class-Short} (page \pageref{Remark
Ref-Class-Short}),
this is conceptually very different from the use of
smallness, describing normal situations. Thus,
it can be considered also on this level an independent question, and we
can
imagine systems based on absolutely ideal situations for normality, but
with a
totally different transfer mechanism.
% Diagram 4.1
% Diagram 4.1
% -----------
%              f
%              f
% @@           A
% @@           A
% @@           A
%        e     d
% @@     A    E F
% @@     A   E   f
% @@     A  E     F
%        b@@caaaaaa@c
% @@        F     E
% @@         F   E
% @@          F E
%              a

\vspace{10mm}

\begin{diagram}

\label{Diagram Information-Transfer}
\index{Diagram Information-Transfer}

\unitlength1.0mm
\begin{picture}(130,110)

\newsavebox{\Tweety}
\savebox{\Tweety}(130,110)[bl]
{

\put(57,18){\vector(1,1){24}}
\put(51,18){\vector(-1,1){24}}
\put(27,51){\vector(1,1){24}}
\put(81,51){\vector(-1,1){24}}

\put(67,61){\line(1,1){4}}

\put(81,47){\vector(-1,0){54}}

\put(24,51){\vector(0,1){22}}
\put(54,81){\vector(0,1){22}}

\put(53,16){$a$}
\put(23,46){$b$}
\put(83,46){$c$}
\put(53,76){$d$}
\put(23,76){$e$}
\put(53,106){$f$}

\put(10,0) {{\rm\bf Information transfer}}

}

\put(0,0){\usebox{\Tweety}}
\end{picture}

\end{diagram}

\vspace{4mm}

For these reasons, extending preferential and related semantics to cover
inheritance nets seems to stretch them to the breaking point, Thus, we
should
also look for other interpretations.
(The term ``interpretation'' is used here in a non-technical sense.) In
particular, it seems worth while to connect inheritance systems to other
problems, and see whether there are similarities there. This is what we do
now.
We come back to the systems $P$ and $R$
in Section \ref{Section Translation} (page \pageref{Section Translation}).

Note that Reiter Defaults behave much more like inheritance nets than like
preferential logics.
%  4.3  Inheritance as information transfer
%  4.3  Inheritance as information transfer
% %
% =========================================
\subsection{
Inheritance as information transfer
}
\label{Section 4.3}

An informal argument showing parallel ideas common to inheritance with an
upward chaining formalism and information transfer is as follows: First,
arrows represent certainly some kind of information, of the kind ``most a's
are
$b' $s'' or so.
(See Diagram \ref{Diagram Information-Transfer} (page \pageref{Diagram
Information-Transfer}).)
Second, to be able to use information,
e.g. ``d's are $f' $s'' at a, we have
to be able to connect from a to $d$ by a valid path, this information has
to be
made accessible to a, or, in other terms, a working information channel
from
a to $d$ has to be established. Third, specificity (when present) decides
conflicts (we take the split validity approach). This can be done
procedurally, or, perhaps simpler and certainly in a more transparent way,
by assigning a comparison of information strength to valid paths. Now,
information strength may also be called truth value (to use a term
familiar
in logic) and the natural entity at hand is the node itself - this is just
a
cheap formal trick without any conceptual meaning.

When we adopt this view, nodes,
arrows, and valid paths have multiple functions, and it may seem that we
overload the (deceptively) simple picture. But it is perhaps the charm and
the utility and naturalness of inheritance systems that they are not
``clean'',
and hide many complications under a simple surface, as human common sense
reasoning often does, too.

In a certain way, this is a poor man's interpretation, as it does not base
inheritance on another formalism, but gives only an intuitive reading.
Yet, it
gives a connection to other branches of reasoning, and is as such already
justified - in the authors' opinion. Moreover, our analysis makes a clear
distinction between arrows and composite valid paths. This distinction
is implicit in inheritance formalisms, we make it explicit through our
concepts.

But this interpretation is by no means the
only one, and can only be suggested as a possibility.

We will now first give the details, and then discuss our interpretation.

\paragraph{
(1) Information:
}

$\hspace{0.01em}$

% (+++ Orig.:  (1) Information: +++)

\label{Section (1) Information:}

Direct positive or negative arrows represent information, valid for their
source. Thus, in a set reading, if there is an arrow $A \xcp B$ in the
diagram, most
elements of $ \xCB $ will be in $B,$ in short: ``most $ \xCB ' $s are $B'
$s'' - and $A \xcP B$ will mean
that most $ \xCB ' $s are not $B' $s.

\paragraph{
(2) Information sources and flow:
}

$\hspace{0.01em}$

% (+++ Orig.:  (2) Information sources and flow: +++)

\label{Section (2) Information sources and flow:}

Nodes are information sources. If $A \xcp B$ is in the diagram, $ \xCB $
is the source
of the information ``most $ \xCB ' $s are $B' $s''.

A valid, composed or atomic positive path $ \xbs $ from $U$ to $ \xCB $
makes the
information of source $ \xCB $ accessible to $U.$ One can also say that $
\xCB ' $s
information becomes relevant to $U.$ Otherwise, information is considered
independent - only (valid) positive paths create the dependencies.

(If we want to conform to inheritance, we must not add trivialities like
``x's are $x' $s'', as this would require $x \xcp x$ in the corresponding
net, which, of
course, will not be there in an acyclic net.)

\paragraph{
(3) Information strength:
}

$\hspace{0.01em}$

% (+++ Orig.:  (3) Information strength: +++)

\label{Section (3) Information strength:}

A valid, composed or atomic positive path $ \xbs $ from $ \xCB ' $ to $
\xCB $ allows us to compare
the strength of information source $ \xCB ' $ with that of $ \xCB:$ $
\xCB ' $ is stronger than $ \xCB.$
(In the set reading, this comparison is the result of specificity: more
specific
information is considered more reliable.) If there is no such valid path,
we
cannot resolve contradictions between information from $ \xCB $ and $ \xCB
'.$
This interpretation results in split validity preclusion: the
comparison between information sources $ \xCB ' $ and $ \xCB $ is
absolute, and does NOT
depend on the $U$ from which both may be accessible - as can be the case
with
total validity preclusion. Of course, if desired, we can also adopt the
much more complicated idea of relative comparison.

Nodes are also truth values. They are the strength of the information
whose
source they are. This might seem an abuse of nodes, but we already have
them,
so why not use them?

\paragraph{
IS-Discussion:
}

$\hspace{0.01em}$

% (+++ Orig.:  IS-Discussion: +++)

\label{Section IS-Discussion:}

Considering direct arrows as information meets probably with little
objection.

The conclusion of a valid path (e.g. if $ \xbs:a \Xl  \xcp b$ is valid,
then its conclusion
is ``a's are $b' $s'') is certainly also information, but it has a status
different
from the information of a direct link, so we should distinguish it
clearly. At
least in upward chaining formalisms, using the path itself as some channel
through which information flows, and not the conclusion, seems more
natural.
The conclusion says little about the inner structure of the path, which is
very important in inheritance networks, e.g. for preclusion. When
calculating
validity of paths, we look at (sub- and other) paths, but not just their
results, and should also express this clearly.

Once we accept this picture of valid positive paths as information
channels, it
is natural to see their upper ends as information sources.

Our interpretation supports upward chaining, and vice versa, upward
chaining
supports our interpretation.

One of the central ideas of inheritance is preclusion, which, in the case
of
split validity preclusion, works by an absolute comparison between nodes.
Thus,
if we accept split validity preclusion, it is natural to see valid
positive
paths as
comparisons between information of different strengths. Conversely, if we
accept absolute comparison of information, we should also accept split
validity
preclusion - these interpretations support each other.

Whatever type of preclusion we accept, preclusion clearly compares
information
strength, and allows us to decide for the stronger one. We can see this
procedurally, or by giving different values to different information,
depending
on their sources, which we can call truth values to connect our picture to
other
areas of logic. It is then natural - as we have it already - to use the
source
node itself as truth value, with comparison via valid positive paths.

\paragraph{
Illustration:
}

$\hspace{0.01em}$

% (+++ Orig.:  Illustration: +++)

\label{Section Illustration:}

Thus, in a given node $U,$ information from $ \xCB $ is accessible iff
there is a
valid positive path from $U$ to $ \xCB,$ and if information from $ \xCB '
$ is also accessible,
and there is a valid positive path from $ \xCB ' $ to $ \xCB,$ then, in
case of conflict,
information from $ \xCB ' $ wins over that from $ \xCB,$ as $ \xCB ' $
has a better truth value.
In the Tweety diagram, see Diagram \ref{Diagram Tweety} (page \pageref{Diagram
Tweety}),
Tweety has access to penguins and birds, the horizontal link from penguin
to
bird compares the strengths, and the fly/not fly arrows are the
information
we are interested in.

Note that negative links and (valid) paths have much less function in our
picture than positive links and valid paths. In a way, this asymmetry is
not surprising, as there are no negative nodes (which would correspond to
something like the set complement or negation). To summarize:
A negative direct link can only be information. A positive direct
link is information at its source, but it can also be a comparison of
truth
values, or it can give access from its source to information at its end.
A valid positive, composed path can only be comparison of truth values, or
give
access to information, it is NOT information itself in the sense of direct
links. This distinction is very important, and corresponds to the
different
treatment of direct arrows and valid paths in inheritance, as it appears
e.g. in
the definition of preclusion. A valid negative composed path has no
function,
only its parts have.

We obtain automatically that direct information is
stronger than any other information: If $ \xCB $ has information $ \xbf,$
and there is a
valid path from $ \xCB $ to $B,$ making B's information accessible to $
\xCB,$ then this same
path also compares strength, and $ \xCB ' $s information is stronger than
B's
information.

Inheritance diagrams in this interpretation represent not only
reasoning with many truth values, but also reasoning $ \xCf about$ those
truth
values: their comparison is done by the same underlying mechanism.

We should perhaps note in this context a connection to an area currently
en
vogue: the problem of trust, especially in the context of web information.
We can see our truth values as the degrees of trust we put into
information
coming from this node. And, we not only use, but also reason about them.

\paragraph{
Further comments:
}

$\hspace{0.01em}$

% (+++ Orig.:  Further comments: +++)

\label{Section Further comments:}

Our reading also covers enriched diagrams, where arbitrary
information can be ``appended'' to a node.

An alternative way to see a source of information is to see it as a reason
to believe the information it gives. $U$ needs a reason to believe
something, i.e. a
valid path from $U$ to the source of the information, and also a reason to
disbelieve, i.e. if $U' $ is below $U,$ and $U$ believes and $U' $ does
NOT believe some
information of $ \xCB,$ then either $U' $ has
stronger information to the contrary, or there is not a valid path to $
\xCB $ any more
(and neither to any other possible source of this information).
$('' Reason'',$ a concept very important in this context, was introduced
by A.Bochman
into the discussion.)

The restriction that negative links can only be information applies to
traditional inheritance networks, and the
authors make no claim whatever that it should also hold for modified such
systems, or in still other contexts. One of the reasons why we do not have
``negative nodes'', and thus negated arrows also in the middle of paths
might be
the following (with $ \xdC $ complementation): If, for some $X,$ we also
have a node for
$ \xdC X,$ then we should have $X \xcP \xdC X$ and $ \xdC X \xcP X,$ thus
a cycle, and arrows from $Y$ to
$X$ should be accompanied by their opposite to $ \xdC X,$ etc.

We translate the analysis and decision
of Definition \ref{Definition IS-2.2} (page \pageref{Definition IS-2.2})  now
into the
picture of information sources, accessibility, and comparison via valid
paths.
This is straightforward:

(1) We have that information from $A_{i},$ $i \xbe I,$ about $B$ is
accessible from $U,$ i.e.
there are valid positive paths from $U$ to all $A_{i}.$ Some $A_{i}$ may
say $ \xCN B,$ some $B.$

(2) If information from $A_{i}$ is comparable with information from
$A_{j}$ (i.e. there
is a valid positive path from $A_{i}$ to $A_{j}$ or the other way around),
and $A_{i}$
contradicts $A_{j}$ with respect to $B,$ then the weaker information is
discarded.

(3) There remains a (nonempty, by lack of cycles) set of the $A_{i},$ such
that
for no such
$A_{i}$ there is $A_{j}$ with better contradictory information about $B.$
If the information
from this remaining set is contradictory, we accept none (and none of the
paths
either), if not, we accept the common conclusion and all these paths.

We continue now Remark \ref{Remark 2.1} (page \pageref{Remark 2.1}), (4),
and turn this into a formal system.

Fix a diagram $ \xbG,$ and do an induction as in Definition \ref{Definition
IS-2.2} (page \pageref{Definition IS-2.2}).

\bd

$\hspace{0.01em}$

% (+++ Orig. No.:  Definition 4.1 +++)

\label{Definition 4.1}

(1) We distinguish $a \xch b$ and $a \xch_{x}b,$ where the intuition of $a
\xch_{x}b$ is:
we know with strength $x$ that a's are $b' $s, and of $a \xch b$ that it
has been
decided taking all information into consideration that $a \xch b$ holds.

(We introduce this notation to point informally to our idea of
information strength, and beyond, to logical systems with
varying strengths of implications.)

(2) $a \xcp b$ implies $a \xch_{a}b,$ likewise $a \xcP b$ implies $a
\xch_{a} \xCN b.$

(3) $a \xch_{a}b$ implies $a \xch b,$ likewise $a \xch_{a} \xCN b$ implies
$a \xch \xCN b.$ This expresses
the fact that direct arrows are uncontested.

(4) $a \xch b$ and $b \xch_{b}c$ imply $a \xch_{b}c,$ likewise for $b
\xch_{b} \xCN c.$ This expresses
concatenation - but without deciding if it is accepted! Note that we
cannot make
$(a \xch b$ and $b \xch c$ imply $a \xch_{b}c)$ a rule, as this would make
concatenation of two
composed paths possible.

(5) We decide acceptance of composed paths as
in Definition \ref{Definition 2.3} (page \pageref{Definition 2.3}), where
preclusion uses accepted paths for deciding.

\ed

Note that we reason in this system not only with, but also about relative
strength of truth values, which are just nodes, this is then, of course,
used
in the acceptance condition, in preclusion more precisely.
%  4.4  Inheritance as reasoning with prototypes
%  4.4  Inheritance as reasoning with prototypes
% %
% ==============================================
\subsection{
Inheritance as reasoning with prototypes
}
\label{Section 4.4}

Some of the issues we discuss here apply also to the more general picture
of
information and its transfer. We present them here for motivational
reasons:
it seems easier to discuss them in the (somewhat!) more concrete setting
of prototypes than in the very general situation of information handling.
These issues will be indicated.

It seems natural to see information in inheritance networks as information
about prototypes. (We do not claim that our use of the word ``prototype''
has
more than a vague relation to the use in psychology. We do not try to
explain
the usefulness of prototypes either, one possibility is that there are
reasons
why birds fly, and why penguins $don' t,$ etc.) In the Tweety diagram, we
will thus say that prototypical birds will fly, prototypical penguins will
not
fly. More precisely, the property ``fly'' is part of the bird prototype, the
property $'' \xCN fly'' $ part of the penguin prototype. Thus, the
information is
given for some node, which defines its application or domain (bird or
penguin in
our example) - beyond this node, the property is not defined (unless
inherited,
of course). It might very well be that no element of the domain has ALL
the
properties of the prototype, every bird may be exceptional in some sense.
This
again shows that we are very far from the ideal picture of small and big
subsets as used in systems $P$ and $R.$ (This, of course, goes beyond the
problem
of prototypes.)

Of course, we will want to ``inherit'' properties of prototypes,
for instance
in Diagram \ref{Diagram Information-Transfer} (page \pageref{Diagram
Information-Transfer}),
a ``should'' inherit the property $e$ from $b,$ and
the property $ \xCN d$ from $c.$ Informally, we will argue as follows:
``Prototypical
a's have property $b,$ and prototypical b's have property $e,$ so it seems
reasonable to assume that prototypical a's also have property $e$ - unless
there
is better information to the contrary.'' A plausible solution is then to
use
upward chaining inheritance as described above to find all relevant
information, and then compose the prototype.

We discuss now three points whose importance goes beyond the treatment of
prototypes:

(1) Using upward chaining has an
additional intuitive appeal: We consider information
at a the best, so we begin with $b$ (and $c),$ and only then, tentatively,
add
information $e$ from $b.$ Thus, we begin with strongest information, and
add weaker
information successively - this seems good reasoning policy.

(2) In upward chaining, we also collect information
at the source (the end of the path), and do not use information which was
already filtered by going down - thus the information we collect has no
history, and we cannot encounter problems of iterated revision, which are
problems of history of change. (In downward chaining, we
only store the reasons why something holds, but not why something does not
hold,
so we cannot erase this negative information when the reason is not valid
any
more. This is an asymmetry apparently not much noted before.
Consider Diagram \ref{Diagram Up-Down-Chaining} (page \pageref{Diagram
Up-Down-Chaining}). Here, the
reason why $u$ does not accept $y$ as information, but $ \xCN y,$ is the
preclusion
via $x.$ But from $z,$ this preclusion is not valid any more, so the
reason
why $y$ was rejected is not valid any more, and $y$ can now be accepted.)

(3) We come back to the question of extensions vs. direct scepticism.
Consider the Nixon Diamond, Diagram \ref{Diagram Nixon-Diamond} (page
\pageref{Diagram Nixon-Diamond}).
Suppose Nixon were a subclass
of Republican and Quaker. Then the extensions approach reasons as follows:
Either the Nixon class prototype has the pacifist property, or the hawk
property, and we consider these two possibilities. But this is not
sufficient:
The Nixon class prototype might have neither property - they are normally
neither pacifists, nor hawks, but some are this, some are that. So the
conceptual basis for the extensions approach does not hold: ``Tertium non
datur''
just simply does not hold - as in Intuitionist Logic, where we may have
neither
a proof for $ \xbf,$ nor for $ \xCN \xbf.$

Once we fixed this decision, i.e. how to find the relevant information, we
can
still look upward or downward in the net and investigate the changes
between the
prototypes in going upward or downward,
as follows: E.g., in above example, we can look at the node a and its
prototype,
and then at the change going from a to $b,$ or, conversely, look at $b$
and its
prototype, and then at the
change going from $b$ to a. The problem of finding the information, and
this
dynamics of information change have to be clearly separated.

In both cases, we see the following:

(1) The language is kept small, and thus efficient.

For instance, when we go from a to $b,$ information about $c$ is lost, and
``c'' does
not figure any more in the language, but $f$ is added. When we go from $b$
to a,
$f$ is lost, and $c$ is won. In our simple picture, information is
independent,
and contradictions are always between two bits of information.

(2) Changes are kept small, and need a reason to be effective.
Contradictory, stronger information will override the old one, but
no other information, except in the following case:
making new information (in-) accessible will cause indirect changes,
i.e. information now made (in-) accessible via the new node.
This is similar to formalisms of causation: if a reason is not there
any more, its effects vanish, too.

It is perhaps more natural when going downward also to
consider ``subsets'', as follows:
Consider Diagram \ref{Diagram Information-Transfer} (page \pageref{Diagram
Information-Transfer}).
b's are $d' $s, and c's are $ \xCN d' $s, and c's are also $b' $s. So it
seems
plausible to go beyond the language of inheritance nets, and conclude that
b's which are not c's will be $d' $s, in short to consider $(b-c)' $s. It
is obvious
which such subsets to consider, and how to handle them: For instance,
loosely
speaking, in $b \xcs d$ $e$ will hold, in $b \xcs c \xcs d$ $ \xCN f$ will
hold, in $b \xcs d \xcs \xdC c$ $f$ will
hold, etc. This is just putting the bits of information together.

We turn to another consideration, which will also transcend the prototype
situation and we will (partly) use the intuition that nodes stand for
sets, and
arrows for (soft, i.e.possibly with exceptions) inclusion in a set or its
complement.

In this reading, specificity stands for soft, i.e. possibly with
exceptions, set
inclusion. So, if $b$ and $c$ are
visible from a, and there is a valid path from $c$ to $b$ (as
in Diagram \ref{Diagram Information-Transfer} (page \pageref{Diagram
Information-Transfer}) ),
then a is a subset both of $b$ and $c,$ and $c$ a subset of $b,$ so
$a \xcc c \xcc b$ (softly). But then a is closer to $c$ than a is to $b.$
Automatically,
a will be closest to itself. This results in a partial, and not
necessarily
transitive relation between these distances.

When we go now from $b$ to $c,$ we lose information $d$ and $f,$ win
information
$ \xCN d,$ but keep information $e.$ Thus, this is minimal change: we give
up (and
win) only the necessary information, but keep the rest. As our language is
very simple, we can use the Hamming distance between formula sets here.
(We
will make a remark on more general situations just below.)

When we look now again from a, we take the set-closest class (c), and use
the information of $c,$ which was won by minimal change (i.e. the Hamming
closest)
from information of $b.$ So we have the interplay of two distances, where
the set distance certainly is not symmetrical, as we need valid paths for
access and comparison. If there is no such valid path, it is reasonable to
make the distance infinite.

We make now the promised remark on more general situations: in richer
languages,
we cannot
count formulas to determine the Hamming distance between two situations
(i.e.
models or model sets), but have to take the difference in propositional
variables. Consider e.g. the language with two variables, $p$ and $q.$ The
models
(described by) $p \xcu q$ and $p \xcu \xCN q$ have distance 1, whereas $p
\xcu q$ and $ \xCN p \xcu \xCN q$
have distance 2. Note that this distance is NOT robust under re-definition
of the language. Let $p' $ stand for $(p \xcu q) \xco ( \xCN p \xcu \xCN
q),$ and $q' $ for $q.$ Of course,
$p' $ and $q' $ are equivalent descriptions of the same model set, as we
can define
all the old singletons also in the new language. Then
the situations $p \xcu q$ and $ \xCN p \xcu \xCN q$ have now distance 1,
as one corresponds to
$p' \xcu q',$ the other to $p' \xcu \xCN q'.$

There might be misunderstandings about the use of the word ``distance''
here. The
authors are fully aware that inheritance networks cannot be captured by
distance
semantics in the sense of preferential structures. But we do NOT think
here of
distances from one fixed ideal point, but of relativized distances: Every
prototype is the origin of measurements. E.g., the bird prototype is
defined
by ``flying, laying eggs, having feathers  \Xl.''. So we presume that all
birds
have these properties of the prototype, i.e. distance 0 from the
prototype.
When we see that penguins do not fly, we move as little as possible from
the
bird prototype, so we give up ``flying'', but not the rest. Thus, penguins
(better: the penguin prototype) will have distance 1 from the bird
prototype
(just one property has changed). So there is a new prototype for penguins,
and
considering penguins, we will not measure from the bird prototype, but
from the
penguin prototype, so the point of reference changes. This is exactly as
in
distance semantics for theory revision, introduced in  \cite{LMS01},
only the point
of reference is not the old theory $T,$ but the old prototype, and the
distance
is a very special one, counting properties assumed to be independent. (The
picture is a little bit more complicated, as the loss of one property
(flying)
may cause other modifications, but the simple picture suffices for this
informal argument.)

We conclude this Section with a remark on prototypes.

Realistic prototypical reasoning will probably neither always be upward,
nor
always be downward. A medical doctor will not begin with the patient's
most
specific category (name and birthday or so), nor will he begin with all he
knows about general objects. Therefore, it seems reasonable to investigate
upward and downward reasoning here.
\section{
Detailed translation of inheritance to modified systems of small sets
}

For background material on abstract size semantics, the reader is
referred to
Chapter \ref{Chapter Size} (page \pageref{Chapter Size}).
\label{Section Translation}
\label{Section Is-Trans}
\index{Section Is-Trans}
%  5.1  Normality
%  5.1  Normality
% %
% ===============
\subsection{
Normality
}
\label{Section 5.1}

As we saw already in Section \ref{Section 4.2} (page \pageref{Section 4.2}),
normality in inheritance (and Reiter defaults etc.) is relative, and as
much normality as possible is preserved. There is no set of absolute
normal
cases of $X,$ which we might denote $N(X),$ but only for $ \xbf $ a set
$N(X, \xbf ),$
elements of $X,$ which behave normally with respect to $ \xbf.$
Moreover, $N(X, \xbf )$ might be defined, but not $N(X, \xbq )$ for
different $ \xbf $ and $ \xbq.$
Normality in the sense of preferential structures is absolute: if $x$ is
not
in $N(X)$ $(=$ $ \xbm (X)$ in preferential reading), we do not know
anything beyond
classical logic.
This is the dark Swedes' problem: even dark Swedes should probably be
tall.
Inheritance systems are different: If birds usually lay eggs, then
penguins,
though abnormal with respect to flying, will still usually lay eggs.
Penguins are fly-abnormal birds, but will continue to be egg-normal birds
-
unless we have again information to the contrary.
So the absolute, simple $N(X)$ of preferential structures splits up into
many, by default independent, normalities.
This corresponds to intuition: There are no absolutely normal birds, each
one is particular in some sense, so $ \xcS \{N(X, \xbf ): \xbf \xbe \xdl
\}$ may well be empty, even
if each single $N(X, \xbf )$ is almost all birds.

What are the laws of relative normality?
$N(X, \xbf )$ and $N(X, \xbq )$ will be largely independent (except for
trivial situations,
where $ \xbf \xcr \xbq,$ $ \xbf $ is a tautology, etc.). $N(X, \xbf )$
might be defined, and $N(X, \xbq )$
not. Connections between the different normalities will be established
only
by valid paths.
Thus, if there is no arrow, or no path, between $X$ and $Y,$ then $N(X,Y)$
and
$N(Y,X)$ - where $X,Y$ are also properties - need not be defined. This
will get rid
of the unwanted connections found with absolute normalities, as
illustrated
by Fact \ref{Fact 5.1} (page \pageref{Fact 5.1}).

We interpret now ``normal'' by ``big set'', i.e. essentially $'' \xbf $ holds
normally in
$X'' $ iff ``there is a big subset of $X,$ where $ \xbf $ holds''. This
will, of course, be
modified.
%  5.2  Small sets
%  5.2  Small sets
% %
% ================
\subsection{
Small sets
}
\label{Section 5.2}

The main interest of this Section is perhaps to show the adaptations of
the
concept of small and big subsets necessary for a more ``real life''
situation,
where we have to relativize. The amount of changes illustrates the
problems and what can be done, but also perhaps what should not be done,
as the
concept is stretched too far.
For more background,
see Chapter \ref{Chapter Size} (page \pageref{Chapter Size}).

As said, the usual informal way of speaking about inheritance networks
(plus other considerations) motivates an interpretation by
sets and soft set inclusion - $A \xcp B$ means that ``most $ \xCB ' $s are
$B' $s''. Just as with
normality, the ``most'' will have to be relativized, i.e. there is a
$B-$normal part
of $ \xCB,$ and a $B-$abnormal one, and the first is $B-$bigger than the
second - where
``bigger'' is relative to $B,$ too. A further motivation for this set
interpretation
is the often evoked specificity argument for preclusion. Thus, we will now
translate our remarks about normality into the language of big and small
subsets.

Consider now the system $P$ (with Cumulativity),
see Definition \ref{Definition Log-Cond-Ref-Size} (page \pageref{Definition
Log-Cond-Ref-Size}).
Recall from Remark \ref{Remark Ref-Class-Short} (page \pageref{Remark
Ref-Class-Short})  that
small sets (see Definition \ref{Definition Gen-Filter} (page \pageref{Definition
Gen-Filter}) ) are
used in two
conceptually very distinct ways: $ \xba \xcn \xbb $ iff the set of $ \xba
\xcu \xCN \xbb -$cases is a small
subset (in the absolute sense, there is just one system of big subsets of
the
$ \xba -$cases) of the set of $ \xba -$cases. The second use is in
information transfer,
used in Cumulativity, or Cautious Monotony more precisely: if the set of
$ \xba \xcu \xCN \xbg -$cases is a small
subset of the set of $ \xba -$cases, then $ \xba \xcn \xbb $ carries over
to $ \xba \xcu \xbg:$ $ \xba \xcu \xbg \xcn \xbb.$
(See also the discussion in  \cite{Sch04}, page 86, after Definition
2.3.6.) It is
this transfer which we will consider here, and not things like AND, which
connect different $N(X, \xbf )$ for different $ \xbf.$

Before we go into details, we will show that e.g. the system $P$ is too
strong
to model inheritance systems, and that e.g. the system $R$ is to weak for
this
purpose. Thus, preferential systems are really quite different from
inheritance
systems.

\bfa

$\hspace{0.01em}$

% (+++ Orig. No.:  Fact 5.1 +++)

\label{Fact 5.1}

(a) System $P$ is too strong to capture inheritance.

(b) System $R$ is too weak to capture inheritance.

\efa

\subparagraph{
Proof
}

$\hspace{0.01em}$

% (+++ Orig.:  Proof +++)

(a) Consider the Tweety diagram,
Diagram \ref{Diagram Tweety} (page \pageref{Diagram Tweety}). $c \xcp b \xcp
d,$ $c \xcP d.$ There is no
arrow $b \xcP c,$ and we will see that $P$ forces one to be there. For
this, we take
the natural translation, i.e. $X \xcp Y$ will be $'' X \xcs Y$ is a big
subset of $X'',$ etc. We
show that $c \xcs b$ is a small subset of $b,$ which we write $c \xcs
b<b.$
$c \xcs b=(c \xcs b \xcs d) \xcv (c \xcs b \xcs \xdC d).$ $c \xcs b \xcs
\xdC d \xcc b \xcs \xdC d<b,$ the latter by $b \xcp d,$ thus
$c \xcs b \xcs \xdC d<b,$ essentially by Right Weakening. Set now $X:=c
\xcs b \xcs d.$ As $c \xcP d,$
$X:=c \xcs b \xcs d \xcc c \xcs d<c,$ and by the same reasoning as above
$X<c.$ It remains to show
$X<b.$ We use now $c \xcp b.$ As $c \xcs \xdC b<c,$ and $c \xcs X<c,$ by
Cumulativity $X=c \xcs X \xcs b<c \xcs b,$
so essentially by OR $X=c \xcs X \xcs b<b.$ Using the filter property, we
see that
$c \xcs b<b.$

(b) Second, even $R$ is too weak: In the diagram $X \xcp Y \xcp Z,$ we
want to conclude that
most of $X$ is in $Z,$ but as $X$ might also be a small subset of $Y,$ we
cannot
transfer the information ``most Y's are in $Z'' $ to $X.$

$ \xcz $
\\[3ex]

We have to distinguish direct information or arrows from inherited
information or valid paths. In the language of big and small sets, it is
easiest
to do this by two types of big subsets: big ones and very big ones. We
will
denote the first big, the second BIG. This corresponds to the distinction
between $a \xch b$ and $a \xch_{a}b$ in Definition \ref{Definition 4.1} (page
\pageref{Definition 4.1}).

We will have the implications $BIG \xcp big$ and $SMALL \xcp small,$ so we
have nested
systems. Such systems were discussed in  \cite{Sch95-1}, see also
 \cite{Sch97-2}. This distinction seems to be necessary to prevent
arbitrary
concatenation of valid paths to valid paths, which would lead to
contradictions.
Consider e.g. $a \xcp b \xcp c \xcp d,$ $a \xcp e \xcP d,$ $e \xcp c.$
Then concatenating $a \xcp b$ with $b \xcp c \xcp d,$
both valid, would lead to a simple contradiction with $a \xcp e \xcP d,$
and not to
preclusion, as it should be - see below.

For the situation $X \xcp Y \xcp Z,$ we will then conclude that:

If $Y \xcs Z$ is a $Z-$BIG subset of $Y$ and $X \xcs Y$ is a $Y-$big
subset of $X$ then $X \xcs Z$ is a
$Z-$big subset of $X.$ (We generalize already to the case where there is a
valid
path from $X$ to $Y.)$

We call this procedure information transfer.

$Y \xcp Z$ expresses the direct information in this context, so $Y \xcs Z$
has to be a $Z-$BIG
subset of $Y.$ $X \xcp Y$ can be direct information, but it is used here
as channel of
information flow, in particular it might be a composite valid path, so in
our
context, $X \xcs Y$ is a $Y-$big subset of $X.$ $X \xcs Z$ is a $Z-$big
subset of $X:$ this can
only be big, and not BIG, as we have a composite path.

The translation into big and small subsets and their modifications is now
quite
complicated: we seem to have to relativize, and we seem to need two types
of big and small. This casts, of course, a doubt on the enterprise of
translation. The future will tell if any of the ideas can be used in other
contexts.

We investigate this situation now in more detail, first without conflicts.

The way we cut the problem is not the only possible one. We were guided by
the
idea that we should stay close to usual argumentation about big and small
sets,
should proceed carefully, i.e. step by step, and should take a general
approach.

Note that we start without any $X-$big subsets defined, so $X$ is not even
a $X-$big
subset of itself.

(A) The simple case of two arrows, and no conflicts.

(In slight generalization:) If information $ \xbf $ is appended at $Y,$
and $Y$ is
accessible from $X$ (and there is no better information about $ \xbf $
available), $ \xbf $
will be valid at $X.$ For simplicity, suppose there is a direct positive
link from
$X$ to $Y,$ written sloppily $X \xcp Y \xcm \xbf.$
In the big subset reading, we will interpret this as: $Y \xcu \xbf $ is a
$ \xbf -$BIG
subset of $Y.$ It is important that this is now direct information, so we
have
``BIG'' and not ``big''. We read now $X \xcp Y$ also as: $X \xcs Y$ is an
$Y-$big subset of $X$ -
this is the channel, so just ``big''.

We want to conclude by transfer that $X \xcs \xbf $ is a $ \xbf -$big
subset of $X.$

We do this in two
steps: First, we conclude that $X \xcs Y \xcs \xbf $ is a $ \xbf -$big
subset of $X \xcs Y,$ and then, as
$X \xcs Y$ is an $Y-$big subset of $X,$ $X \xcs \xbf $ itself is a $ \xbf
-$big subset of $X.$ We do NOT
conclude that $(X-Y) \xcs \xbf $ is a $ \xbf -$big subset of $X-$Y, this
is very important, as we
want to preserve the reason of being $ \xbf -$big subsets - and this goes
via $Y!$
The transition from ``BIG'' to ``big'' should be at the first step, where we
conclude that $X \xcs Y \xcs \xbf $ is a $ \xbf -$big (and not $ \xbf
-$BIG) subset of $X \xcs Y,$ as it is
really here where things happen, i.e. transfer of information from $Y$ to
arbitrary subsets $X \xcs Y.$

We summarize the two steps in a slightly modified notation, corresponding
to the
diagram $X \xcp Y \xcp Z:$

(1) If $Y \xcs Z$ is a $Z-$BIG subset of $Y$ (by $Y \xcp Z),$ and $X \xcs
Y$ is a $Y-$big subset of $X$
(by $X \xcp Y),$ then $X \xcs Y \xcs Z$ is a $Z-$big subset of $X \xcs Y.$

(2) If $X \xcs Y \xcs Z$ is a $Z-$big subset of $X \xcs Y,$ and $X \xcs Y$
is a $Y-$big subset of $X$
(by $X \xcp Y)$ again, then $X \xcs Z$ is a $Z-$big subset of $X,$ so $X
\Xl  \xcp Z.$

Note that (1) is very different from Cumulativity or even Rational
Monotony, as
we do not say anything about $X$ in comparison to $Y:$ $X$ need not be any
big or
medium size subset of $Y.$

Seen as strict rules, this will not work, as it would result in
transitivity,
and thus
monotony: we have to admit exceptions, as there might just be a negative
arrow
$X \xcP Z$ in the diagram. We will discuss such situations below in (C),
where we
will modify our approach slightly, and obtain a clean analysis.

(Here and in what follows, we are very cautious, and relativize all
normalities. We could perhaps obtain our objective with a more daring
approach,
using absolute normality here and there. But this would be a purely
technical
trick (interesting in its own right), and we look here more for a
conceptual
analysis, and, as long as we do not find good conceptual reasons why to be
absolute here and not there, we will just be relative everywhere.)

We try now to give justifications for the two (defeasible) rules. They
will be
philosophical and can certainly be contested and/or improved.

For (1):

We look at $Y.$ By $X \xcp Y,$ Y's information is accessible at $X,$ so,
as $Z-$BIG is
defined for $Y,$ $Z-$big will be defined for $Y \xcs X.$ Moreover,
there is a priori nothing
which prevents $X$ from being independent from $Y,$ i.e. $Y \xcs X$ to
behave like $Y$
with respect to $Z$ - by default: of course, there could be a negative
arrow
$X \xcP Z,$ which would prevent this.
Thus, as $Y \xcs Z$ is a $Z-$BIG subset of $Y,$ $Y \xcs X \xcs Z$ should
be a $Z-$big subset of
$Y \xcs X.$ By the same argument (independence), we should also conclude
that
$(Y-X) \xcs Z$ is a $Z-$big subset of $Y-$X. The definition of $Z-$big for
$Y-X$ seems,
however, less clear.

To summarize, $Y \xcs X$ and $Y-X$ behave by default with respect to $Z$
as $Y$ does, i.e.
$Y \xcs X \xcs Z$ is a
$Z-$big subset of $Y \xcs X$ and $(Y-X) \xcs Z$ is a $Z-$big subset of
$Y-$X. The reasoning is
downward, from supersets to subsets, and symmetrical to $Y \xcs X$ and
$Y-$X. If the
default is violated, we need a reason for it.
This default is an assumption about the adequacy of the language. Things
do not
change wildly from one concept to another (or, better: from $Y$ to $Y \xcu
X),$ they
might change, but then we are told so - by a corresponding negative link
in the
case of diagrams.

For (2):

By $X \xcp Y,$ $X$ and $Y$ are related, and we assume that $X$ behaves as
$Y \xcs X$ does
with respect to $Z.$ This is upward reasoning, from subset to superset and
it is
NOT symmetrical: There is no reason to suppose that $X-Y$ behaves the same
way as
$X$ or $Y \xcs X$ do with respect to $Z,$ as the only reason for $Z$ we
have, $Y,$ does not
apply. Note that, putting relativity aside (which can also be considered
as being
big and small in various, per default independent dimensions) this is
close to
the reasoning with absolutely big and small sets: $X \xcs Y-(X \xcs Y \xcs
Z)$ is small in
$X \xcs Y,$ so a fortiori small in $X,$ and $X-(X \xcs Y)$ is small in
$X,$ so
$(X-(X \xcs Y)) \xcv (X \xcs Y-(X \xcs Y \xcs Z))$ is small in $X$ by the
filter property, so $X \xcs Y \xcs Z$
is big in $X,$ so a fortiori $X \xcs Z$ is big in $X.$

Thus, in summary, we conclude by default that,

(3) If $Y \xcs Z$ is a $Z-$BIG subset of $Y,$ and $X \xcs Y$ is a $Y-$big
subset of $X,$ then
$X \xcs Z$ is a $Z-$big subset of $X.$

(B) The case with longer valid paths, but without conflicts.

Treatment of longer paths: Suppose we have a valid composed path from $X$
to
$Y,$ $X \Xl  \xcp Y,$ and not any longer a direct link $X \xcp Y.$ By
induction, i.e. upward
chaining, we argue - using directly (3) - that $X \xcs Y$ is a $Y-$big
subset of $X,$ and
conclude by (3) again that $X \xcs Z$ is a $Z-$big subset of $X.$

(C) Treatment of multiple and perhaps conflicting information.

Consider Diagram \ref{Diagram Multiple} (page \pageref{Diagram Multiple}) :
% Diagram 5.1
% Diagram 5.1
% -----------
%                                    Z
%                                    Z
%    @@                            E AF
%    @@                           E  A F
%    @@                          E   A  F
%    @@                         E    D   F
%    @@                        E     A    F
%    @@                       E      A     F
%                          Y@@caaaaa@Y'     Y"
%    @@                       F      A
%    @@                        F     A
%    @@                         F    A
%    @@                          F   A
%    @@                           F  A
%    @@                            F A
%                          U@@caaaaa@X

\vspace{10mm}

\begin{diagram}

\label{Diagram Multiple}
\index{Diagram Multiple}

\unitlength1.0mm
\begin{picture}(130,100)

\newsavebox{\sets}
\savebox{\sets}(140,90)[bl]
{

\put(37,8){\vector(-1,1){22}}

\put(13,38){\vector(1,1){24}}
\put(67,38){\vector(-1,1){24}}

\put(40,37){\vector(0,1){23}}

\put(40,8){\vector(0,1){22}}

\put(36,34){\vector(-1,0){22}}

\put(36,4){\vector(-1,0){22}}

\put(39,3){$X$}
\put(9,33){$Y$}
\put(39,33){$Y'$}
\put(69,33){$Y''$}
\put(39,63){$Z$}
\put(9,3){$U$}

\put(38,50){\line(1,0){3.7}}

\put(10,90) {{\rm\bf Multiple and conflicting information}}

}

\put(0,0){\usebox{\sets}}
\end{picture}

\end{diagram}

\vspace{4mm}

We want to analyze the situation and argue that e.g. $X$ is mostly not in
$Z,$ etc.

First, all arguments about $X$ and $Z$ go via the $Y' $s. The arrows from
$X$ to
the $Y' $s, and from $Y' $ to $Y$ could also be valid paths. We look at
information
which concerns $Z$ (thus $U$ is not considered), and which is accessible
(thus $Y'' $
is not considered). We can
be slightly more general, and consider all possible combinations of
accessible
information, not only those used in the diagram by $X.$ Instead of arguing
on the level of $X,$ we will argue one level above, on the Y's and their
intersections, respecting specificity and unresolved conflicts.

(Note that in more general situations, with arbitrary information
appended, the
problem is more complicated, as we have to check which information is
relevant
for some $ \xbf $ - conclusions can be arrived at by complicated means,
just as in
ordinary logic. In such cases, it might be better first to look at all
accessible information for a fixed $X,$ then at the truth values and their
relation, and calculate closure of the remaining information.)

We then have (using the obvious language: ``most $ \xCB ' $s are $B' $s''
for $'' A \xcs B$ is a
big subset of $ \xCB '',$ and ``MOST $ \xCB ' $s are $B' $s'' for $'' A
\xcs B$ is a BIG subset of $ \xCB '' ):$

In $Y,$ $Y'',$ and $Y \xcs Y'',$ we have that MOST cases are in $Z.$
In $Y' $ and $Y \xcs Y',$ we have that MOST cases are not in $Z$ $(=$ are
in $ \xdC Z).$
In $Y' \xcs Y'' $ and $Y \xcs Y' \xcs Y'',$ we are UNDECIDED about $Z.$

Thus:

$Y \xcs Z$ will be a $Z-$BIG subset of $Y,$ $Y'' \xcs Z$ will be a $Z-$BIG
subset of $Y'',$
$Y \xcs Y'' \xcs Z$ will be a $Z-$BIG subset of $Y \xcs Y''.$

$Y' \xcs \xdC Z$ will be a $Z-$BIG subset of $Y',$ $Y \xcs Y' \xcs \xdC
Z$ will be a $Z-$BIG subset of
$Y \xcs Y'.$

$Y' \xcs Y'' \xcs Z$ will be a $Z-$MEDIUM subset of $Y' \xcs Y'',$ $Y
\xcs Y' \xcs Y'' \xcs Z$ will be
a $Z-$MEDIUM subset of $Y \xcs Y' \xcs Y''.$

This is just simple arithmetic of truth values, using specificity and
unresolved
conflicts, and the non-monotonicity is pushed into the fact that subsets
need
not preserve the properties of supersets.

In more complicated situations, we implement e.g. the general principle
(P2.2)
from Definition \ref{Definition IS-2.2} (page \pageref{Definition IS-2.2}), to
calculate the truth values. This
will use in our case specificity for conflict resolution, but it is an
abstract procedure, based on an arbitrary relation $<.$

This will result in the ``correct'' truth value for the intersections, i.e.
the one corresponding to the other approaches.

It remains to do two things: (C.1) We have to assure that $X$ ``sees'' the
correct
information, i.e. the correct intersection, and, (C.2), that $X$ ``sees''
the
accepted $Y' $s, i.e. those through which valid paths go, in order to
construct not
only the result, but also the correct paths.

(Note that by split validity preclusion, if there is valid path from $
\xCB $ through $B$
to $C,$ $ \xbs:A \xFB \xcp B,$ $B \xcp C,$ and $ \xbs ':A \xFB \xcp B$
is another valid path from $ \xCB $ to $B,$ then
$ \xbs ' \xDM B \xcp C$ will also be a valid path. Proof: If not, then $
\xbs ' \xDM B \xcp C$ is precluded,
but the same preclusion will also preclude $ \xbs \xDM B \xcp C$ by split
validity
preclusion, or it is contradicted, and a similar argument applies again.
This is the same argument as the one for the simplified definition of
preclusion
- see Remark \ref{Remark 2.1} (page \pageref{Remark 2.1}), (4).)

(C.1) Finding and inheriting the correct information:

$X$ has access to $Z-$information from $Y$ and $Y',$ so we
have to consider them. Most of $X$ is in $Y,$ most of $X$ is in $Y',$
i.e.
$X \xcs Y$ is a $Y-$big subset of $X,$ $X \xcs Y' $ is a $Y' -$big subset
of $X,$ so
$X \xcs Y \xcs Y' $ is a $Y \xcs Y' -$big subset of $X,$ thus most of $X$
is in $Y \xcs Y'.$

We thus have $Y,$ $Y',$ and $Y \xcs Y' $ as possible reference classes,
and use
specificity to choose $Y \xcs Y' $ as reference class. We do not know
anything e.g. about $Y \xcs Y' \xcs Y'',$ so this is not a possible
reference class.

Thus, we use specificity twice, on the $Y' s-$level (to decide that $Y
\xcs Y' $ is
mostly not in $Z),$ and on $X' s-$level (the choice of the reference
class), but
this is good policy, as, after all, much of nonmonotonicity is about
specificity.

We should emphasize that nonmonotonicity lies in the behaviour of the
subsets,
determined by truth values and comparisons thereof, and the choice of the
reference class by specificity. But both are straightforward now and local
procedures, using information already decided before. There is no
complicated
issue here like determining extensions etc.

We now use above argument, described in the simple case, but with more
detail,
speaking in particular about the most specific reference class for
information about $Z,$ $Y \xcs Y' $ in our example - this is used
essentially in
(1.4), where the ``real'' information transfer happens, and where we go
from BIG to big.

(1.1) By $X \xcp Y$ and $X \xcp Y' $ (and there are no other $Z-$relevant
information
sources), we have to consider $Y \xcs Y' $ as reference class.

(1.2) $X \xcs Y$ is a $Y-$big subset of $X$ (by $X \xcp Y)$
(it is even $Y-$BIG, but we are immediately more general to treat valid
paths),
$X \xcs Y' $ is a $Y' -$big subset of $X$ (by $X \xcp Y' ).$
So $X \xcs Y \xcs Y' $ is a $Y \xcs Y' -$big subset of $X.$

(1.3) $Y \xcs Z$ is a $Z-$BIG subset of $Y$ (by $Y \xcp Z),$
$Y' \xcs \xdC Z$ is a $Z-$BIG subset of $Y' $ (by $Y' \xcP Z),$
so by preclusion
$Y \xcs Y' \xcs \xdC Z$ is a $Z-$BIG subset of $Y \xcs Y'.$

(1.4) $Y \xcs Y' \xcs \xdC Z$ is a $Z-$BIG subset of $Y \xcs Y',$ and
$X \xcs Y \xcs Y' $ is a $Y \xcs Y' -$big subset of $X,$ so
$X \xcs Y \xcs Y' \xcs \xdC Z$ is a $Z-$big subset of $X \xcs Y \xcs Y'.$

This cannot be a strict rule without the reference class, as it would then
apply
to $Y \xcs Z,$ too, leading to a contradiction.

(2) If $X \xcs Y \xcs Y' \xcs \xdC Z$ is a $Z-$big subset of $X \xcs Y
\xcs Y',$ and
$X \xcs Y \xcs Y' $ is a $Y \xcs Y' -$big subset of $X,$ so
$X \xcs \xdC Z$ is a $Z-$big subset of $X.$

We make this now more formal.

We define for all nodes $X,$ $Y$ two sets: $B(X,Y),$ and $b(X,Y),$ where
$B(X,Y)$ is the set of $Y-$BIG subsets of $X,$ and
$b(X,Y)$ is the set of $Y-$big subsets of $X.$
(To distinguish undefined from medium/MEDIUM-size, we will also have to
define $M(X,Y)$ and $m(X,Y),$ but we omit this here for simplicity.)

The translations are then:

$(1.2' )$ $X \xcs Y \xbe b(X,Y)$ and $X \xcs Y' \xbe b(X,Y' )$ $ \xch $ $X
\xcs Y \xcs Y' \xbe b(X,Y \xcs Y' )$

$(1.3' )$ $Y \xcs Z \xbe B(Y,Z)$ and $Y' \xcs \xdC Z \xbe B(Y',Z)$ $ \xch
$ $Y \xcs Y' \xcs \xdC Z \xbe B(Y \xcs Y',Z)$ by preclusion

$(1.4' )$ $Y \xcs Y' \xcs \xdC Z \xbe B(Y \xcs Y',Z)$ and $X \xcs Y \xcs
Y' \xbe b(X,Y \xcs Y' )$ $ \xch $ $X \xcs Y \xcs Y' \xcs \xdC Z \xbe b(X
\xcs Y \xcs Y',Z)$
as $Y \xcs Y' $ is the most specific reference class

$(2' )$ $X \xcs Y \xcs Y' \xcs \xdC Z \xbe b(X \xcs Y \xcs Y',Z)$ and $X
\xcs Y \xcs Y' \xbe b(X,Y \xcs Y' )$ $ \xch $ $X \xcs \xdC Z \xbe b(X,Z).$

Finally:

$(3' )$ $A \xbe B(X,Y)$ $ \xcp $ $A \xbe b(X,Y)$ etc.

Note that we used, in addition to the set rules, preclusion, and the
correct
choice of the reference class.

(C.2) Finding the correct paths:

Idea:

(1) If we come to no conclusion, then no path is valid, this is trivial.

(2) If we have a conclusion:

(2.1) All contradictory paths are out: e.g. $Y \xcs Z$ will be $Z-$big,
but
$Y \xcs Y' \xcs \xdC Z$ will be $Z-$big. So there is no valid path via
$Y.$

(2.2) Thus, not all paths supporting the same conclusion are valid.

Consider the following Diagram \ref{Diagram Paths-Conclusions} (page
\pageref{Diagram Paths-Conclusions}) :
% Diagram 5.2
% Diagram 5.2
% -----------
%                                    Z
%                                    Z
%    @@                            E AF
%    @@                           E  A   F
%    @@                          E   A    F
%    @@                         E    D     F
%    @@                        E     A      F
%    @@                       E      A       F
%                          Y@@caaaaa@Y'@@caaaa@Y"
%    @@                       F      A       E
%    @@                        F     A      E
%    @@                         F    A     E
%    @@                          F   A    E
%    @@                           F  A   E
%    @@                            F AE
%                                    X

\vspace{10mm}

\begin{diagram}

\label{Diagram Paths-Conclusions}
\index{Diagram Paths-Conclusions}

\unitlength1.0mm
\begin{picture}(130,100)

\newsavebox{\setsx}
\savebox{\setsx}(140,90)[bl]
{

\put(37,8){\vector(-1,1){22}}
\put(41,8){\vector(1,1){22}}

\put(13,38){\vector(1,1){24}}
\put(67,38){\vector(-1,1){24}}

\put(40,37){\vector(0,1){23}}

\put(40,8){\vector(0,1){22}}

\put(36,34){\vector(-1,0){22}}
\put(66,34){\vector(-1,0){22}}

\put(39,3){$X$}
\put(9,33){$Y$}
\put(39,33){$Y'$}
\put(69,33){$Y''$}
\put(39,63){$Z$}

\put(38,50){\line(1,0){3.7}}

\put(10,90) {{\rm\bf Valid paths vs. valid conclusions}}

}

\put(0,0){\usebox{\setsx}}
\end{picture}

\end{diagram}

\vspace{4mm}

There might be a positive path
through $Y,$ a negative one through $Y',$ a positive one through $Y'' $
again, with
$Y'' \xcp Y' \xcp Y,$ so $Y$ will be out, and only $Y'' $ in. We can see
this, as there is a
subset, $\{Y,Y' \}$ which shows a change: $Y' \xcs Z$ is $Z-$BIG, $Y' \xcs
\xdC Z$ is $Z-$BIG,
$Y'' \xcs Z$ is $Z-$BIG, and $Y \xcs Y' \xcs \xdC Z$ is $Z-$BIG, and the
latter can only happen if
there is a preclusion between $Y' $ and $Y,$ where $Y$ looses. Thus, we
can see this
situation by looking only at the sets.

We show now equivalence with the inheritance formalism given
in Definition \ref{Definition 2.3} (page \pageref{Definition 2.3}).

\bfa

$\hspace{0.01em}$

% (+++ Orig. No.:  Fact 5.2 +++)

\label{Fact 5.2}

The above definition and the one outlined
in Definition \ref{Definition 2.3} (page \pageref{Definition 2.3})  correspond.

\efa

\subparagraph{
Proof
}

$\hspace{0.01em}$

% (+++ Orig.:  Proof +++)

By induction on the length of the deduction that $X \xcs Z$ (or $X \xcs
\xdC Z)$ is a $Z-$big
subset of $X.$ (Outline)

It is a corollary of the proof that we have to consider only subpaths and
information of all generalized paths between $X$ and $Z.$

Make all sets (i.e. one for every node) sufficiently different, i.e.
all sets and boolean combinations of sets differ by infinitely many
elements,
e.g. $A \xcs B \xcs C$ will have infinitely many less elements than $A
\xcs B,$ etc. (Infinite
is far too many, we just choose it by laziness to have room for the
$B(X,Y)$ and
the $b(X,Y).$

Put in $X \xcs Y \xbe B(X,Y)$ for all $X \xcp Y,$ and $X \xcs \xdC Y \xbe
B(X,Y)$ for all $X \xcP Y$ as base
theory.

$Length=1:$
Then big must be BIG, and, if $X \xcs Z$ is a $Z-$BIG subset of $X,$ then
$X \xcp Z,$ likewise
for $X \xcs \xdC Z.$

We stay close now to above Diagram \ref{Diagram Multiple} (page \pageref{Diagram
Multiple}),
so we argue for the negative case.

Suppose that we have deduced $X \xcs \xdC Z \xbe b(X,Z),$ we show that
there must
be a valid negative path from $X$ to $Z.$ (The other direction is easy.)

Suppose for simplicity that there is no negative link from $X$ to $Z$ -
otherwise
we are finished.

As we can distinguish intersections from elementary sets (by the starting
hypothesis about sizes), this can only be deduced using
$(2' ).$ So there must be some suitable $\{Y_{i}:i \xbe I\}$ and we must
have deduced
$X \xcs \xcS Y_{i} \xbe b(X, \xcS Y_{i}),$ the second hypothesis of $(2'
).$
If $I$ is a singleton, then we have the induction hypothesis, so there is
a valid
path from $X$ to $Y.$ So suppose $I$ is not a singleton. Then the
deduction of
$X \xcs \xcS Y_{i} \xbe b(X, \xcS Y_{i})$ can only be done by $(1.2' ),$
as this is the only rule having in the conclusion an elementary set on the
left
in $b(.,.),$ and a true intersection on the right. Going back along $(1.2'
),$
we find $X \xcs Y_{i} \xbe b(X,Y_{i}),$ and by the induction hypothesis,
there are valid paths
from $X$ to the $Y_{i}.$

The first hypothesis of $(2' ),$ $X \xcs \xcS Y_{i} \xcs \xdC Z \xbe b(X
\xcs \xcS Y_{i},Z)$ can be obtained by
$(1.3' )$ or $(1.4' ).$ If it was obtained by $(1.3' ),$ then $X$ is one
of the $Y_{i},$ but
then there is a direct link from $X$ to $Z$ (due to the ``B'', BIG). As a
direct link
always wins by specificity, the link must be negative, and we have a valid
negative path from $X$ to $Z.$ If it was obtained by $(1.4' ),$ then its
first
hypothesis $ \xcS Y_{i} \xcs \xdC Z \xbe B( \xcS Y_{i},Z)$ must have been
deduced, which can only be by
$(1.3' ),$ but the set of $Y_{i}$ there was chosen to take all $Y_{i}$
into account for which
there is a valid path from $X$ to $Y_{i}$ and arrows from the $Y_{i}$ to
$Z$ (the rule was
only present for the most specific reference class with respect to $X$ and
$Z!),$
and we are done by the definition of valid paths
in Section \ref{Section Is-Inh-Intro} (page \pageref{Section Is-Inh-Intro}).

$ \xcz $
\\[3ex]

We summarize our ingredients.

Inheritance was done essentially by (1) and (2) of (A) above
and its elaborations (1.i), (2) and $(1.i' ),$ $(2' ).$ It consisted of a
mixture of bold and careful (in comparison to systems $P$ and $R)$
manipulation of
big subsets. We had to be bolder than the systems $P$ and $R$ are, as we
have to
transfer information also to perhaps small subsets. We had to be more
careful,
as $P$ and $R$ would have introduced far more connections than are
present. We also
saw that we are forced to loose the paradise of absolute small and big
subsets,
and have to work with relative size.

We then have a plug-in decision what to do with contradictions. This is a
plug-in, as it is one (among many possible) solutions to the much more
general
question of how to deal with contradictory information, in the presence of
a
(partial, not necessarily transitive) relation which compares strength. At
the
same place of our procedure, we can plug in other solutions, so our
approach
is truly modular in this aspect. The comparing relation is defined by the
existence of valid paths, i.e. by specificity.

This decision is inherited downward using again the specificity criterion.

Perhaps the deepest part of the analysis can be described as follows:
Relevance is coded by positive arrows, and valid positive paths, and thus
is similar to Kripke structures for modality, where the arrows code
dependencies between the situations for the evaluation of the modal
quantifiers.
In our picture, information at A can become relevant only to node $B$ iff
there
is a valid positive path from $B$ to A. But, relevance (in this reading,
which is closely related to causation) is profoundly non-monotonic, and
any
purely monotonic treatment of relevance would be insufficient. This seems
to
correspond to intuition. Relevance is then expressed in our translation to
small and big sets formally by the possibility
of combining different small and big sets in information transfer. This
is, of
course, a special form of relevance, there might be other forms.

% ******* BEGIN LATEX SOURCE FILE 11-arg.tex *******
%
% Uebers. aus Karltex File: 11-arg.m
%
%
\chapter{
Argumentation
}
\section{Blocking paths}

We have two arrows $a \xcp b$ (positive) and $a \xcP b$ (negative). The
basic idea is that
one negative argument blocks all positive ones. (We will see later that
this
is not a restriction.)

This is a generalization over just ending the relation $ \xcp $ at a
point, as
some $x$ might be accessible from a, but not from $\{a,b\},$ as $b$
introduces a
blocking argument - thus it becomes non-monotonic. Of course, the blocking
might itself be blocked by a new $c$ in $\{a,b,c\},$ etc.

E.g. $a \xcp b \xcP c$ will block $a \xcp c.$

The question is which nodes are visible from a (set of) $node(s),$ we
denote
this set $ \ol{a}$ or $ \ol{A}$ (in the case of a set). In above example,
$b \xbe \ol{a},$ but $c \xce \ol{a}.$
Intuitively, we can see this as a relative horizon.

We will assume that there are no cycles, and that the networks are finite.

Suppose $a \xcP b$ is in a network, then we can look at $ \ol{\{a,b\}},$
and we then force $b$
to be valid, despite $a \xcP b.$ ``Outside forcing'' overrides negative
arrows.

\bfa

$\hspace{0.01em}$

% (+++ Orig. No.:  Fact 1: +++)

\label{Fact 1:}

(1) A version of Cumulativity holds:
(Cum) $A \xcc B \xcc \ol{A}$ $ \xch $ $ \ol{B}= \ol{A}$

(2) If $x \xbe \ol{A},$ but $x \xce \ol{A \xcv \{a\}},$ then $a \xcP x,$
or there is some new $a' \xbe \ol{A \xcv \{a\}}- \ol{A}.$

(3) $x \xce \ol{A},$ $x \xbe \ol{A \xcv \{a\}}$ $ \xcH $ $x \xbe \ol{a}$
(x might just have destroyed a counterargument).

(4) $x \xbe \ol{A},$ $x \xbe \ol{B}$ $ \xch $ $x \xbe \ol{A \xcv B}$?

\efa

\subparagraph{
Proof
}

$\hspace{0.01em}$

% (+++ Orig.:  Proof +++)

 \Xl.

$ \xcz $
\\[3ex]

There are problems of ambiguity.

\be

$\hspace{0.01em}$

% (+++ Orig. No.:  Example 1: +++)

\label{Example 1:}

Consider $a \xcp b,$ $c,$ and then add
(1) nothing
(2) $a \xcP c$
(3) $b \xcP c$
(4) $a \xcP c$ and $b \xcP c$
(5) $b \xcP c$ and $a \xcp c$
They cannot be distinguished - $b \xbe \ol{a}$ will always be the case,
but neither a
nor $b$ ``lead'' to anything else, as $a \xcp c$ is destroyed by $a \xcp b
\xcP c.$

\ee

Of course, giving more information, e.g. some $d \xcp b$ and $d \xcp c$
allows us to see
if $b \xcP c$ is present. Thus, giving some special nodes to ``read out''
information and
``feed in'' information can disambiguate.

A more serious problem is the following:

\be

$\hspace{0.01em}$

% (+++ Orig. No.:  Example 2: +++)

\label{Example 2:}

$b \xcp c \xcP x,$ $b \xcP d \xcp x,$ $a \xcP c,$ $a \xcp d$ are the base
arrows. Adding at least one of
$a \xcp x$ or $b \xcp x$ will give the same information - irrespective of
whether
we add just one or both. But the case without adding any of the two is
different. This is easily seen by examining the cases.

For $ \ol{a},$ $a \xcp x$ is not necessary, as $a \xcp d \xcp x$ is valid.
For $ \ol{b},$ $b \xcp c \xcP x$ will destroy $b \xcp x.$
For $ \ol{\{a,b\}}:$ $a \xcP c$ blocks $c,$ so $c \xcP x$ is not valid,
and $x$ is free. But $b \xcP d$ blocks
$d \xcp x,$ so $x$ is not accessible via $d,$ so we need at least one of
$a \xcp x$ or $b \xcp x.$

\ee

This example makes a representation proof difficult, as we have to branch
into
several possibilities. The following results show that one can probably
just
add both, when there is an ``or''.

\paragraph{
New results (5.4.08)
}

$\hspace{0.01em}$

% (+++ Orig.:   New results (5.4.08) +++)

\label{Section New results (5.4.08)}

We will write now $X+x$ for $X \xcv \{x\}$ etc.

\bfa

$\hspace{0.01em}$

% (+++ Orig. No.:  Fact 2: +++)

\label{Fact 2:}

(1) Let $a \xcp x,$ $b \xcp y$ be such that

(a) for all A either

(1.1) $ \ol{A}= \ol{A+a \xcp x}= \ol{A+b \xcp y}= \ol{A+a \xcp x+b \xcp
y}$ or

(1.2) $ \ol{A} \xEd \ol{A+a \xcp x}= \ol{A+b \xcp y}= \ol{A+a \xcp x+b
\xcp y}$

(b) there is at least one A with property (1.2).

Then $x=y.$

(2) Let $a \xcp x,$ $b \xcp y$ be such that

(a) for all A either

(1.1) $ \ol{A}= \ol{A+a \xcp x}= \ol{A+b \xcp y}= \ol{A+a \xcp x+b \xcp
y}$ or

$(1.2' )$ $ \ol{A}= \ol{A+a \xcp x}= \ol{A+b \xcp y} \xEd \ol{A+a \xcp x+b
\xcp y}$

(b) there is at least one A with property $(1.2' ).$

Then $x=y.$

(3) Situation (2) is impossible.

(4) Let $ \ol{A} \xEd \ol{A+a \xcp x}= \ol{A+b \xcp y}.$ Then $x=y.$

(5) Let $ \ol{A} \xEd \ol{A+a \xcp x}= \ol{A+b \xcp x}.$ Then $ \ol{A+a
\xcp x}= \ol{A+b \xcp x}= \ol{A+a \xcp x+b \xcp x}.$

Analogous properties hold for negative arrows:

(6) Let $ \ol{A} \xEd \ol{A+a \xcP x}= \ol{A+b \xcP y}.$ Then $x=y.$

(7) Let $ \ol{A} \xEd \ol{A+a \xcP x}= \ol{A+b \xcP x}.$ Then $ \ol{A+a
\xcP x}= \ol{A+b \xcP x}= \ol{A+a \xcP x+b \xcP x}.$

(8) Let $ \ol{A}= \ol{A+a \xcP x}= \ol{A+b \xcP y} \xEd \ol{A+a \xcP x+b
\xcP y},$ then $x=y.$

(9) $ \ol{A}= \ol{A+a \xcP x}= \ol{A+b \xcP x} \xEd \ol{A+a \xcP x+b \xcP
x}$ is impossible.

\efa

\subparagraph{
Proof
}

$\hspace{0.01em}$

% (+++ Orig.:  Proof +++)

(1)
Let A be with property (1.2). As $ \ol{A} \xEd \ol{A+a \xcp x},$ $a \xcp
x$ has an effect, this can only
be because $x \xce \ol{A},$ $x \xbe \ol{A+a \xcp x}.$ Analogously $y \xce
\ol{A},$ $y \xbe \ol{A+b \xcp y}.$ Thus $x,y \xce \ol{A},$
$x,y \xbe \ol{A+a \xcp x}= \ol{A+b \xcp y}= \ol{A+a \xcp x+b \xcp y}.$
Then $ \ol{A+a \xcp x}= \ol{A+x},$ and $y \xce \ol{A},$ $y \xbe \ol{A+x},$
thus $x$ is
before $y$ (not necessarily $x \xcp y$ or so, maybe $x \xcP z \xcP y$
etc.). Analogously
$y$ is before $x,$ as the diagram is free from cycles, $x=y.$

(2)
If $x,y \xbe \ol{A},$ then $ \ol{A}= \ol{A+a \xcp x+b \xcp y}.$ If $x \xbe
\ol{A},$ then $ \ol{A+b \xcp y}= \ol{A+a \xcp x+b \xcp y},$ as $x$ is
already
present, contradiction. Analogously for $y,$ so $x,y \xce \ol{A}.$
As $ \ol{A+a \xcp x} \xEd \ol{A+a \xcp x+b \xcp y},$ $y \xbe \ol{A+a \xcp
x+b \xcp y}$ (otherwise, $+b \xcp y$ would have no
consequences). Analogously for $x,$ thus $x,y \xbe \ol{A+a \xcp x+b \xcp
y}.$
As $x \xce \ol{A+a \xcp x},$ but $x \xbe \ol{A+a \xcp x+b \xcp y},$ $y$
has to be before $x,$ analogously $x$ before $y,$
so again $x=y.$

(3)
Thus $ \ol{A}= \ol{A+a \xcp x}= \ol{A+b \xcp x} \xEd \ol{A+a \xcp x+b \xcp
x}.$ But e.g. $ \ol{A+b \xcp x} \xEd \ol{A+a \xcp x+b \xcp x}$ is
impossible,
as the number of supporting arguments is unimportant.

(4)
The proof of (1) uses only the new prerequisites.

(5)
The number of arguments are unimportant. So $ \ol{A+a \xcp x}= \ol{A+b
\xcp x}= \ol{A+a \xcp x+b \xcp x}= \ol{A+x}.$

(6)
$a \xcP x$ has effect, so $x \xbe \ol{A},$ $x \xce \ol{A+a \xcP x},$
analogously for $y \xbe \ol{A},$ $y \xce \ol{A+b \xcP y}.$
Thus, again, $x$ is before $y,$ $y$ before $x,$ so $x=y.$

(7)
Again, the number of arguments is unimportant.

(8)
If $x \xce \ol{A},$ then $ \ol{A+b \xcP y}= \ol{A+a \xcp x+b \xcP y},$
contradiction. Thus, $x,y \xbe \ol{A}.$
As $ \ol{A+a \xcP x} \xEd \ol{A+a \xcP x+b \xcP y},$ $y \xce \ol{A+a \xcP
x+b \xcP y},$ otherwise, $b \xcP y$ would have no
consequences. Analogously for $x,$ so $x,y \xce \ol{A+a \xcP x+b \xcP y}.$
As $x \xbe \ol{A+a \xcP x},$ but
$x \xce \ol{A+a \xcP x+b \xcP y},$ $x$ has to be before $y,$ and vice
versa, so $x=y.$

(9)
$x \xbe \ol{A}= \ol{A+a \xcP x}= \ol{A+b \xcP x},$ so $a,b \xce \ol{A},$
and $x \xbe \ol{A+a \xcP x+b \xcP y},$ contradiction.

$ \xcz $
\\[3ex]

\bfa

$\hspace{0.01em}$

% (+++ Orig. No.:  Fact 3: +++)

\label{Fact 3:}

(1) $ \ol{A} \xEd \ol{A+a \xcp x}$ (meaning: same A, but added to graph $a
\xcp x)$ $ \xcp $
(1.1) $a \xbe \ol{A}$
(1.2) $x \xce \ol{A},$ $x \xbe \ol{A+a \xcp x}$
(1.3) all $z$ s.t. $z \xbe \ol{A}$ $ \xcr $ $z \xce \ol{A+a \xcp x}$ are
downstream from $x$
(1.4) $ \ol{A+a \xcp x}= \ol{A \xcv \{x\}}$ (meaning: graph unchanged, but
added $x$ to A)

(2) $ \ol{A} \xEd \ol{A+a \xcP x}$ $ \xcp $
(1.1) $a \xbe \ol{A}$
(1.2) $x \xbe \ol{A},$ $x \xce \ol{A+a \xcP x}$
(1.3) all $z$ s.t. $z \xbe \ol{A}$ $ \xcr $ $z \xce \ol{A+a \xcP x}$ are
downstream from $x$
(1.4) to have analogon, perhaps make it possible to force $x$ out? Easier
representation?

(3) If $x \xbe \ol{A},$ then $ \ol{A}= \ol{A+a \xcp x}.$ If $x,y \xbe
\ol{A},$ then $ \ol{A}= \ol{A+a \xcp x+b \xcp y}.$

\efa

New proof for (2) in above Fact 2:
As $ \ol{A+b \xcp y} \xEd \ol{A+b \xcp y+a \xcp x},$ so by (1.2) above $x
\xce \ol{A}= \ol{A+b \xcp y},$ and $x \xbe \ol{A+a \xcp x+b \xcp y}.$
Analogously, $y \xce \ol{A},$ $y \xbe \ol{A+a \xcp x+b \xcp y}.$ So $x,y
\xce \ol{A},$ $x,y \xbe \ol{A+a \xcp x+b \xcp y}.$
The rest of the argument holds by (1.3) above.

New proof for (3) in above Fact 2:
$x \xce \ol{A}= \ol{A+a \xcp x}$ $ \xcp $ $a \xce \ol{A}$ or ex. $c \xbe
\ol{A},$ $c \xcP x.$ $x \xbe \ol{A+a \xcp x+b \xcp y},$ so $b \xbe \ol{A+a
\xcp x},$ so $b \xbe \ol{A},$
as $b$ is upstream from $x,$ and there is no $c \xbe \ol{A},$ $c \xcP x.$
Analogously for $a \xbe \ol{A}.$
\chapter{
Acknowledgements
}

Torre, Herzig, div. journals for permission, David,
greasy spoon

% \include{1-0-do}

% ******* BEGIN LATEX SOURCE FILE 12-BIB.tex *******
%
% Uebers. aus Karltex File: 12-BIB.m
%
%

\end{document}